\documentclass[a4paper,11pt, oneside]{book}%aggiungi opzione `` openany `` per togliere le pagine bianche a fine capitolo; aggiungi opzione "oneside" per avere tutto sulla stessa pagina (include automaticamente openany).
\usepackage{geometry}
\usepackage{pdfpages} %Per inserire PDF esterni

\usepackage{changepage} %Per usare adjustwidth nella copertina

\date{\today}
\usepackage{authblk}

\usepackage[italian, english]{babel} %ATT. Viene preso in considerazione l'ultimo scritto! In questo caso, l'inglese
\usepackage[OT2, T1]{fontenc} %languages and accents
\usepackage[utf8]{inputenc} %languages and accents

\usepackage{latexsym}
\usepackage{amssymb,amsmath}%,mathrsfs}
\usepackage{amsthm,graphicx}
\usepackage{enumerate}
\usepackage{array}

\usepackage{xcolor} %per usare i colori
\usepackage{color} %per usare i colori

\definecolor{uniba}{RGB}{0, 64, 112}
\usepackage[colorlinks = true,
 linkcolor = uniba,
 urlcolor = uniba,
 citecolor = uniba,
 anchorcolor = uniba
]{hyperref}
\usepackage{bm}

\usepackage[hyperpageref]{backref} %Ti inserisce nella bibliografia le pagine dove è stato citato l'articolo.

\usepackage{esint}

\usepackage{centernot}

\usepackage{lipsum}

\usepackage[noadjust]{cite} %per avere referenze consecutive con trattini [5--10]

\usepackage{sectsty} % sets colour of Chapters

\chapterfont{\color{uniba}} 
\sectionfont{\color{uniba}}
\subsectionfont{\color{uniba}}
\subsubsectionfont{\color{uniba}}

\usepackage{titlesec}
\usepackage{microtype}
\usepackage{tikz}

\titleformat{\chapter}[display]
 {\normalfont\bfseries\color{uniba}}
 {\filleft%
 \begin{tikzpicture}
 \node[
 outer sep=0pt,
 text width=2.5cm,
 minimum height=3cm,
 fill=uniba,
 font=\color{white}\fontsize{80}{90}\selectfont,
 align=center
 ] (num) {\thechapter};
 \node[
 rotate=90,
 anchor=south,
 font=\color{black}\Large\normalfont
 ] at ([xshift=-5pt]num.west) {\textls[180]{\textsc{\chaptertitlename}}}; 
 \end{tikzpicture}%
 }
 {10pt}
 {\titlerule[2.5pt]\vskip3pt\titlerule\vskip4pt\LARGE\sffamily}

\usepackage[]{frontespizio}

\usepackage{emptypage} %per togliere le testatine dalle pagine vuote a fine capitolo

\usepackage{lmodern}

\usepackage{ifthen}
\usepackage{fancyhdr} %per gestire le testatine. VEDI http://profs.sci.univr.it/~gregorio/breveguida.pdf
\renewcommand{\sectionmark}[1]{\markright{\emph{\thesection.\ #1}}}
\newenvironment{abstract}{\cleardoublepage
\vspace*{3\baselineskip}\center
\huge{\textbf{\tu{Abstract}} \normalsize}\endcenter
\quote}{\endquote\clearpage}

\allowdisplaybreaks %per spezzare le equazioni su più pagine 
\def\claim#1.{\noindent {\bf #1.}}

\def\flushright#1{{\unskip\nobreak\hfil\penalty50\hskip2em\hbox{}\nobreak\hfil%
#1\parfillskip=0pt\finalhyphendemerits=0\par}}

\def\bull{\vrule height 1.8ex width 1.0ex depth .1ex }
\def\QED{\ifmmode\eqno\hbox{$\bull$}\else\flushright{\hbox{$\bull$}}\fi}

\newcommand{\parag}[1]{\left\{ \begin{aligned} #1 \end{aligned}\right.}

\newcommand{\abs}[1]{\left| #1 \right|} 
\newcommand{\norm}[1]{\left\Vert#1\right\Vert} 
\newcommand{\tnorm}[1]{{\left\vert\kern-0.25ex\left\vert\kern-0.25ex\left\vert #1 
		\right\vert\kern-0.25ex\right\vert\kern-0.25ex\right\vert}}
		
\newcommand{\mc}[1]{\mathcal{#1}}

\usepackage{scalerel}
\usepackage{tikz}
\usetikzlibrary{svg.path}
\definecolor{orcidlogocol}{HTML}{A6CE39}
\tikzset{
 orcidlogo/.pic={
 \fill[orcidlogocol] svg{M256,128c0,70.7-57.3,128-128,128C57.3,256,0,198.7,0,128C0,57.3,57.3,0,128,0C198.7,0,256,57.3,256,128z};
 \fill[white] svg{M86.3,186.2H70.9V79.1h15.4v48.4V186.2z}
 svg{M108.9,79.1h41.6c39.6,0,57,28.3,57,53.6c0,27.5-21.5,53.6-56.8,53.6h-41.8V79.1z M124.3,172.4h24.5c34.9,0,42.9-26.5,42.9-39.7c0-21.5-13.7-39.7-43.7-39.7h-23.7V172.4z}
 svg{M88.7,56.8c0,5.5-4.5,10.1-10.1,10.1c-5.6,0-10.1-4.6-10.1-10.1c0-5.6,4.5-10.1,10.1-10.1C84.2,46.7,88.7,51.3,88.7,56.8z};
 }
}
\newcommand\orcidicon[1]{\href{https://orcid.org/#1}{\mbox{\scalerel*{
\begin{tikzpicture}[yscale=-1,transform shape]
\pic{orcidlogo};
\end{tikzpicture}
}{|}}}}

\newtheorem{Theorem}{Theorem}[section]
\newtheorem{Proposition}[Theorem]{Proposition}
\newtheorem{Corollary}[Theorem]{Corollary}
\newtheorem{Lemma}[Theorem]{Lemma}
\newtheorem{Remark}[Theorem]{Remark}
\newtheorem{Definition}[Theorem]{Definition}
\newtheorem{Example}[Theorem]{Example}

\newcommand{\R}{\mathbb{R}}
\newcommand{\N}{\mathbb{N}}
\newcommand{\Z}{\mathbb{Z}}
\newcommand{\G}{{\mathbb G}}
\newcommand{\F}{\mathbb{F}}

\def\calD{{\cal D}}

\newcommand{\cupl}{{\rm cupl}}
\newcommand{\cat}{{\rm cat}}
\newcommand{\genus}{{\rm genus}} 

\def\supp{\mathop{\rm supp}}

\newcommand{\Tail}{{\rm Tail}}
\newcommand{\DG}{{\rm DG}}
\newcommand{\Ima}{{\rm Im}}
\newcommand{\Ker}{{\rm Ker}}
\newcommand{\dive}{{\rm div}}

\newcommand{\eps}{\varepsilon}
\def\epsilon{\varepsilon}
\def\eps{\varepsilon}

\newcommand{\wto}{\rightharpoonup}
\newcommand{\nequiv}{\not \equiv}
\newcommand{\sgn}{{\rm sgn}}
\def\meas{\mathop{\rm meas}\nolimits}

\def\half{\frac{1}{2}} %\frac{1}{2}}

\def\abs#1{|#1|}
\def\pabs#1{\left|{#1}\right|}
\def\norm#1{\|#1\|}
\def\intRN{\int_{\R^N}}
\def\RRint{\iint_{\R^N\times\R^N}}

\newcommand{\dist}{{\rm dist}}
\newcommand{\dista}{{\rm dist}}

\newcommand{\tu}[1]{\textcolor{uniba}{#1}}

\usepackage{biolinum}

\usepackage{minitoc} %\usepackage[tight,english]{minitoc} %genera mini indici all'inizio di ogni capitolo: \minitoc\mtcskip all'inizio di ogni capitolo

\begin{document}

\pagenumbering{Roman}

\newpage

%VERSIONE STAMPATA
%\newpage $\,$

%FRONTESPIZIO STANDARD UNIBA (INGLESE) - VERSIONE STAMPA 

\thispagestyle{empty}

\begin{figure*}%[!h]
\begin{center}
\includegraphics*[width=0.35\textwidth]{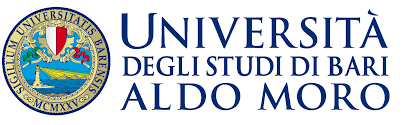}
\end{center}
\vspace{-2em}
\end{figure*}
\begin{center}
%\noindent
\rule{13.5cm}{0.4pt}%\\
%\end{center}
%
%
%\begin{center}

\bigskip
\bigskip
\bigskip

{\large Università degli Studi di Bari Aldo Moro
\\}

\medskip

{\large Department of Mathematics
\\}

\bigskip
\medskip

{\large Ph.D. in Computer Science and Mathematics\\
}

\medskip

{
%\vskip.3cm
\normalsize
XXXV Cycle\\}

\bigskip
\medskip

{%\large 
%Academic discipline %s list 
%%for italian university research and teaching 
Mathematical Analysis
MAT/05\\}

%\bigskip
%\medskip
%
%{\large Ph.D. Thesis
%\\}

\bigskip
\bigskip
\bigskip
\bigskip

%{\Large \textbf{LOCAL AND NONLOCAL ASPECTS\\}}
%{\Large \textbf{IN NONLINEAR VARIATIONAL PROBLEMS\\} }
%{\Large \textbf{NONLOCAL ELLIPTIC PDES\\}}
%{\Large \textbf{WITH GENERAL NONLINEARITIES\\} }
\tu{
{\huge \textbf{Nonlocal Elliptic PDEs\\}}
{\huge \textbf{with General Nonlinearities\\} }
}

%Nonlocal Elliptic PDEs with general nonlinearities

\end{center}

\bigskip
\bigskip
\bigskip
\bigskip
\bigskip

%\begin{minipage}{0.5\textwidth}
%Supervisor:\\
%\textbf{Prof. Silvia Cingolani}
%\end{minipage}
%\begin{minipage}{0.5\textwidth}
%%\begin{flushright}
%PhD Candidate:\\
%\textbf{%Dott. 
%Marco Gallo}
%%\end{flushright}
%\end{minipage}
%\vskip.5cm
%\begin{minipage}{0.5\textwidth}
%Co-Supervisor:\\
%\textbf{Prof. Denis Bonheure}
%\end{minipage}
%\vskip.9cm
%\begin{minipage}{0.5\textwidth}
%Coordinator:\\
%\textbf{Prof. Francesca Mazzia}
%\end{minipage}
%\begin{minipage}{0.5\textwidth}
%\end{minipage}

\begin{minipage}{0.5\textwidth}
%Ph.D. Candidate:\\
Ph. Doctor:\\
\textbf{%Dott. 
Marco Gallo}
\end{minipage}
\begin{minipage}{0.5\textwidth}
Supervisor:\\
\textbf{Prof. Silvia Cingolani}
\end{minipage}
\vskip0.5cm
\begin{minipage}{0.5\textwidth}
$\,$\\
$\,$
\end{minipage}
\begin{minipage}{0.5\textwidth}
Co-Supervisor:\\
\textbf{Prof. Denis Bonheure}
\end{minipage}
%%
%\vskip1.7cm
%%
%\begin{minipage}{0.5\textwidth}
%$\,$\\
%$\,$
%\end{minipage}
%%
%\begin{minipage}{0.5\textwidth}
%Ph.D. Coordinators:\\
%\textbf{Prof. Maria Francesca Costabile}\\
%\textbf{Prof. Francesca Mazzia}
%\end{minipage}
%
%\begin{minipage}{0.5\textwidth}
%\end{minipage}
%\begin{flushright}
%\end{flushright}

\bigskip
\bigskip
\bigskip
\bigskip
\bigskip
\bigskip
\bigskip

\begin{center}
\noindent\rule{8cm}{0.4pt}\\

\medskip

%Final Exam 
Ph.D. Thesis, March 2023\\
\noindent\rule{8cm}{0.4pt}
\end{center}

%%%%DA FLAVIA
%\include{copertinadedica} 
%%%%%

%FRONTESPIZIO CLASSICO

%\begin{figure*}[!h]
%\includegraphics*[width=1\textwidth]{logo}
%\end{figure*}
%\noindent\rule{13.5cm}{0.4pt}\\
%
%\begin{center}
%
%{\large Dottorato di Ricerca in Informatica e Matematica\\
%\vskip.2cm
%XXXV Ciclo\\}
%
%\bigskip
%
%{\large Settore Scientifico Disciplinare MAT/05\\}
%
%\bigskip
%\bigskip
%
%{\large Tesi di Dottorato\\}
%
%\bigskip
%\bigskip
%\bigskip
%
%{\Large \textbf{LOCAL AND NONLOCAL ASPECTS\\}}
%{\Large \textbf{IN NONLINEAR VARIATIONAL PROBLEMS\\} }
%\end{center}
%
%\bigskip
%\bigskip\bigskip
%
%\begin{minipage}{0.5\textwidth}
%Dottorando:\\
%\textbf{Dott. Marco Gallo}
%\end{minipage}
%\begin{minipage}{0.5\textwidth}
%%\begin{flushright}
%Coordinatore:\\
%\textbf{Prof.ssa Francesca Mazzia}
%%\end{flushright}
%\end{minipage}
%\vskip.5cm
%\begin{minipage}{0.5\textwidth}
%Supervisore:\\
%\textbf{Prof.ssa Silvia Cingolani}
%\end{minipage}
%\vskip.5cm
%\begin{minipage}{0.5\textwidth}
%Co-Supervisore:\\
%\textbf{Prof. Denis Bonheure}
%\end{minipage}
%\begin{minipage}{0.5\textwidth}
%\end{minipage}
%
%\bigskip
%\bigskip
%\bigskip
%\bigskip
%
%\begin{center}
%\noindent\rule{8cm}{0.4pt}\\
%Esame Finale 2022\\
%\noindent\rule{8cm}{0.4pt}
%
%\end{center}

%FRONTESPIZIO CAMBIATO

%\usepackage[firstpage=true]{background}
%\usepackage{lipsum}

\thispagestyle{empty}

\newgeometry{bottom=-6in}

\begin{adjustwidth}{-1in}{-1in}

\begin{tikzpicture}[remember picture,overlay]
%\node [opacity=1,scale=1] at (current page.center) {\includegraphics{webb7}};
\node [opacity=1,scale=1] at (current page.center) {\includegraphics{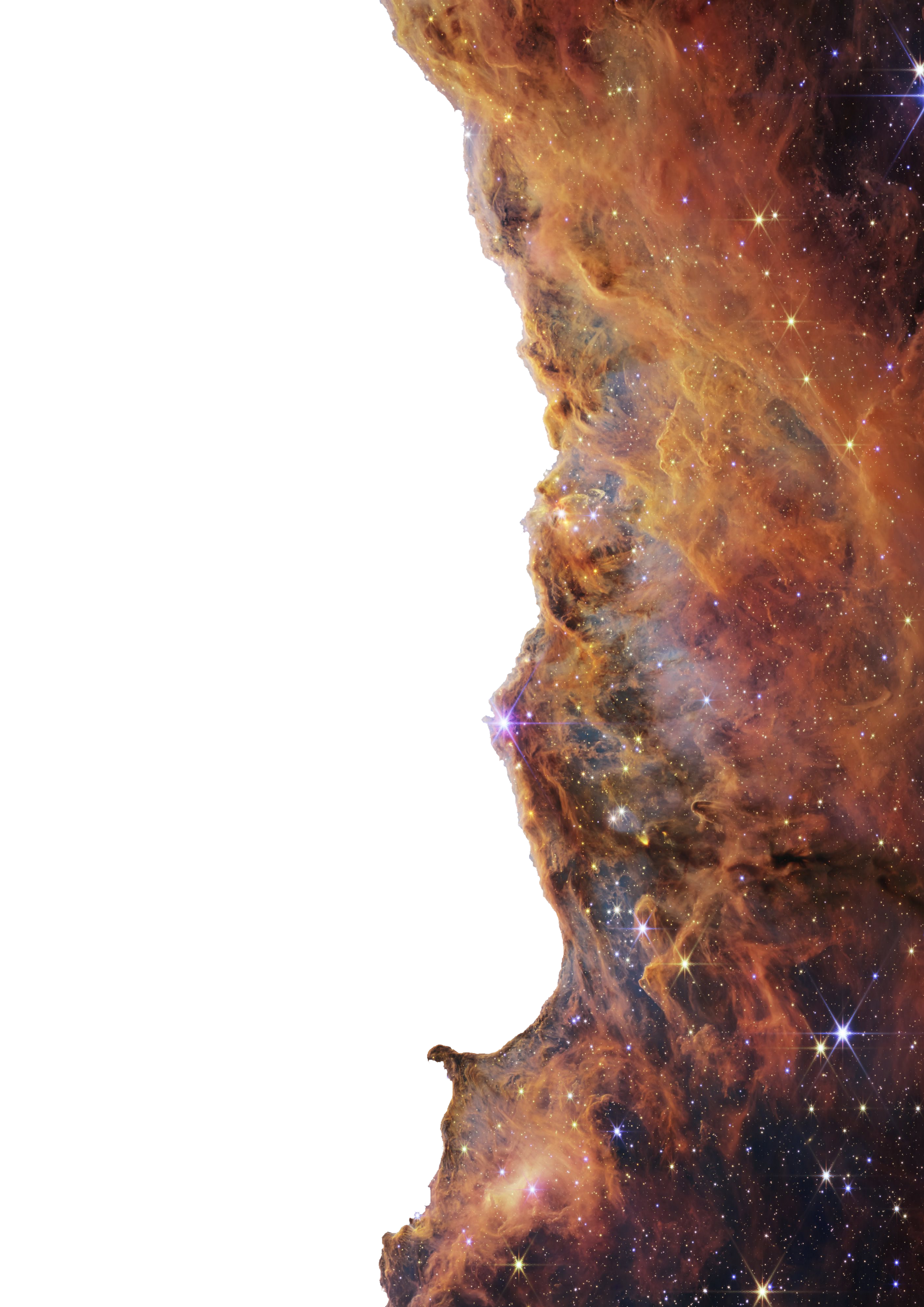}};
%\node [opacity=1,scale=1] at (current page.center) {\includegraphics{webb2(blue)}};
\end{tikzpicture}

\biolinum{

\begin{flushleft} 

\huge{
\textbf{Marco Gallo}
}

\vskip1cm

%\sffamily

\tu{
{\Huge \textbf{Nonlocal Elliptic PDEs\\}}
{\Huge \textbf{with General Nonlinearities\\} }
}
%Nonlocal elliptic PDEs with general nonlinearities

\vskip3cm

\textsf{
\Large{Supervised by}}
 \\
\textsf{
\LARGE{
Prof. Silvia Cingolani}
}
\\
\textsf{
\LARGE{
Prof. Denis Bonheure}
}

\vskip4cm

%\Large{Ph.D. Thesis, 2022}

\Large{Ph.D. Thesis} %\\2022}

\vskip8.5cm

\includegraphics*[width=0.42\textwidth]{logo_uniba}

\end{flushleft} 

}

\end{adjustwidth}

\restoregeometry

\newpage

\cleardoublepage
\thispagestyle{empty}

%VERSIONE STAMPATA (NON PIU')
%\newgeometry{paperwidth=17cm, paperheight=24cm, textheight=250mm, textwidth=100mm, top=-5mm, left=5mm, right=3mm}

\begin{abstract} 

\normalsize

\medskip

%The general spirit of this thesis is to investigate how the nonlocalities -- both in the operator and in the source -- comes into account in the study of different PDEs, and analyze these equations under almost optimal assumptions of the nonlinearity.

$\; $ 
In this thesis we investigate how the nonlocalities affect the study of different PDEs coming from physics, and we analyze these equations under almost optimal assumptions of the nonlinearity. In particular, we focus on the fractional Laplacian operator and on sources involving convolution with the Riesz potential, as well as on the interaction of the two, and we aim to do it through variational and topological methods.

$\; $ 
We examine both quantitative and qualitative aspects, proving multiplicity of solutions for nonlocal nonlinear problems with free or prescribed mass, showing regularity, positivity, symmetry and sharp asymptotic decay of ground states, and exploring the influence of the topology of a potential well in presence of concentration phenomena. On the nonlinearities we consider general assumptions which avoid monotonicity and homogeneity: this generality obstructs the use of classical variational tools and forces the implementation of new ideas.

$\; $ 
Throughout the thesis we develop some new tools: among them, a Lagrangian formulation modeled on Pohozaev mountains is used for the existence of normalized solutions, annuli-shaped multidimensional paths are built for genus-based multiplicity results, a fractional chain rule is proved to treat concave powers, and a fractional center of mass is defined to detect semiclassical standing waves. We believe that these tools could be used to face problems in different frameworks as well.

\end{abstract}

\frontmatter %Così compaiono i numeri romani nella prefazione

\newpage

%\cleardoublepage
%%\enlargethispage{1.5\baselineskip}
%\vspace*{-2\baselineskip}

%\addtocontents{toc}{\vspace{-50.0pt}}

\tableofcontents %Indice

\chapter*{Introduction
%\tb{da scrivere}%COMMENT NOW
}
\addcontentsline{toc}{chapter}{Introduction}
\markboth{Introduction}{Introduction}

Nonlinear phenomena pervade natural and social sciences, and lots of them are modeled by nonlinear equations: %Nonlinear Partial Differential Equations: 
there has been an enormous progress in the study of the structure and in the qualitative understanding of these equations in recent years, and many astonishing interrelations have been found. 
In this thesis we aim to contribute to these studies. 

In particular, the goal is to detect local and nonlocal effects in some nonlinear partial differential equations, having as a common feature a variational structure. Mathematically, \emph{nonlocality} is an intrinsic feature of integral operators and of associated energy functionals, which have the peculiarity -- contrary to the classical local ones -- of capturing long-range interactions or self-interactions.
%
%, i.e. events that happen far away, may that be in time or in space. 
%
In the context of functional variational principles and associated inequalities, nonlocal energy functionals are currently receiving great attention since %as
 they are closely related to problems %both 
in geometry, physics, engineering, biology, finance %social sciences 
and many others, manifesting both in the operator and in the source.
%in quantum mechanics, nuclear physics and also surprisingly in classical geometric optimization problems.
In this setting, classical PDE theory fails because of the presence of the nonlocality.

\bigskip

A first goal of our research is the study of some %Some objectives of our research are the 
generalized nonlinear Schr\"odinger equations (here the Planck's constant and the mass are normalized $\hbar=m=1$) %with dispersion
%[CAMBIA QUESTA $g$, e VEDI la $f$ dopo]
$$ %\begin{equation}\label{pseudo}
i \partial_t u = P(D)u- h(|u|)u, \quad x \in \R^N, \, t>0
$$ %\end{equation} 
where $P(D)$ denotes a %self-adjoint 
pseudo-differential operator with constant coefficients, defined by multiplication in Fourier spaces as $ \widehat{P(D)u}(\xi)= p(\xi) \widehat{u}(\xi)$, and $h \in C(\R_+)$.
%, under suitable assumptions on the multiplier $p(\xi)$.
In particular, %we will focus In particular, 
we are interested to the case %the symbol is chosen as 
$p(\xi)\equiv|\xi|^{2s}$, $s\in (0,1)$, and to %, and we will focus on 
the study %in the study %existence 
of standing waves solutions
$$u(t,x)= \linebreak e^{i \mu t} Q(x)$$
with some nontrivial profile $Q$, depending on the frequency $\mu>0$: %, where the symbol is chosen as $p(\xi)=|\xi|^{2s}$: % \in \R$. % and the velocity $v \in \R^N$. 
this leads to investigate the so called \emph{fractional nonlinear Schr\"odinger equation} (fNLS), % pseudo-differential equation
$$%\begin{equation*}\label{boost}
 (-\Delta)^s Q + \mu Q = h(|Q|)Q, \quad x \in \R^N
$$ 
where $P(D)\equiv (-\Delta)^s$ is known as \emph{fractional Laplacian}.
In 1948 Feynman \cite{Fey0} proposed indeed a new suggestive description of the %time 
evolution of the state of a non-relativistic quantum particle: 
according to Feynman, the wave function solution of the Schr\"odinger equation should be given by a %{\sl 
sum over all possible histories of the system, %}, 
that is by a heuristic %integral over the space of paths. The classical notion of a single, unique classical trajectory for a system is replaced by a functional 
integral over an infinity of quantum-mechanically possible trajectories. 
Following this % Feynman's %path integral 
approach, % to quantum mechanics, 
Laskin \cite{Las0,Las1,Las2,Las3} 
%generalized the path integral over Brownian motions (random motion seen in swirling gas molecules) to L\'evy flights (a mix of long trajectories and short, random movements found in turbulent fluids) and 
derived the fractional Schr\"odinger equation (fNLS): 
numerous applications of these %(fNLS) 
equations in the physical sciences could be mentioned, ranging from image reconstruction %the description of boson stars 
to water wave dynamics, passing through %from 
jump processes in probability theory with applications to financial mathematics. %, as well as applications to .

%\smallskip

In this thesis we are interested in detecting %implementing new deformation theorems to detect 
existence of one or more %standing wave 
solutions of % looking at 
(fNLS) equations, or more generally %\eqref{}, or more generally 
problems related to equations of the type
\begin{equation}\label{eq_introd_frac}
 (-\Delta)^s u + \mu u = g(u), \quad x \in \R^N, \tag{I.1}
\end{equation}
where $g \in C(\R)$, % and $\widehat{(-\Delta)^s u}(\xi)= |\xi|^{2s} \widehat{u}(\xi)$, 
%we are interested in detecting %implementing new deformation theorems to detect 
%existence of one or more %standing wave 
%solutions, %as well as multiplicity, 
%some issues related to the existence of boosted solutions for equation \eqref{pseudo} and their symmetry properties.
and in studying their qualitative %symmetry 
properties. We aim to do it by looking at solutions as critical points of suitable real-valued functionals, % defined on space of functions, 
as well as by exploiting methods coming from both algebraic topology and geometry. Here, the influence of an external potential $V=V(x)$ may be considered as well.

\medskip

%Second we aim to analyze nonlocal differential equations in isotropic and aniso\-tropic media.
%We are interested in analyzing 
Another target of this thesis is the analysis of the so-called \emph{Pekar nonlinear problem}, which describes a polaron -- namely a quantum electron in a polar crystal -- at rest. This problem was raised by Pekar \cite{Pek0} in 1954: 
the atoms of the crystal are displaced due to the electrostatic force induced by the charge of the electron and the resulting deformation is then felt by the electron itself. 
%Mathematically, this model belongs to the class of equations
%%equation has the following form
%$$
%-\Delta u + \mu u= (W \ast F(u)) F'(u), \quad x \in \R^N
%$$
%where $W$ is an even, self-interaction potential, %$N \geq 2$, 
%$\mu>0$ and $F \in C^1(\R)$. 
Afterwards, Choquard \cite{Cho0} (see also Lieb \cite{Lie1,Lie2} and Lions \cite{Lio1}) developed a similar theory to study steady states of the one-component plasma approximation in the Hartree-Fock theory; %,Lio1}; 
the same model was then also derived by Penrose in his discussion about the self-gravitational collapse of a quantum-mechanical wave function \cite{Pen1,Pen2,Pen3}, coupling together the Schr\"odinger equation with the Newton law.
Mathematically, these models belong to the class of equations
%equation has the following form
$$
-\Delta u + \mu u= (W \ast F(u)) F'(u), \quad x \in \R^N
$$
where $W$ is a radially symmetric %n even, 
%self-interaction 
potential, %$N \geq 2$, 
$\mu>0$ and $F \in C^1(\R)$. In particular, the abovementioned physics problems are set in the case $N=3$, $F$ power and $W(x)\equiv \frac{1}{4 \pi |x|}$ \emph{Newton potential}.

We address to study existence, multiplicity and qualitative results for these integro-differential equations, in the wider (model) class of \emph{Riesz potentials} $W(x)\equiv I_{\alpha}(x):= \frac{C_{N,\alpha}}{|x|^{N-\alpha}}$, with $\alpha \in (0,N)$ and $C_{N,\alpha}>0$ constant, that is
\begin{equation}\label{eq_introd_choq}
-\Delta u + \mu u= (I_{\alpha} \ast F(u)) F'(u), \quad x \in \R^N \tag{I.2}
\end{equation}
also known as \emph{Choquard-Hartree-Pekar equation}. 
%When $N=3$ and $\alpha=2$, $I_2(x)=\frac{1}{4 \pi |x|}$ is also known as \emph{Newton potential}.

\medskip

When dealing with the mathematical description of the gravitational collapse of %(massless) 
exotic %boson 
stars, \emph{double nonlocalities} arise naturally, both in the operator and in the source: this was observed already by Chandrasekhar \cite{Chd0} in 1931, and then developed by Lieb, Thirring and Yau \cite{Thir1,LT0, LiYa0,LiYa1}. % in the '80s. 
Other applications can be found for example in quantum chemistry and in the study of graphene. 
This is why part of the thesis will be devoted to the study of equations of the type
\begin{equation}\label{eq_introd_doubl}
 (-\Delta)^s u + \mu u = (I_{\alpha} \ast F(u)) F'(u), \quad x \in \R^N, \tag{I.3}
\end{equation}
%focusing especially on %the qualitative aspects related to this equations and 
highlighting especially how the two nonlocalities interact.

\bigskip

The approach of this thesis will be mainly of variational type: in the last thirty % 30 
years, the study of abstract variational methods and their applications to nonlinear differential equations have greatly developed. % in the last 30 years. 
In the past, variational methods have been applied to solve nonlinear differential equations, both ordinary and partial, taking advantage of a related functional with some specific features: among them we can find compactness properties (typically the Palais-Smale condition), natural constraints of Nehari type, use of integral identities (such as the Pohozaev identity), presence of a local operator, restriction to %n operator of local type or 
 bounded domains, and others.
In the subsequent years, the study of nonlinear differential equations arising in geometry, physics and applied mathematics has suggested developments in which at least one of the previous assumptions is not satisfied.

The substantial progress %in the study of nonlinear PDEs 
made in the %recent 
last years allows now to tackle equations with particular features, as %the 
nonlocal PDEs. 
The greatly increased interest in %Recently, the interest in 
nonlocal operators %has greatly increased, %because of the applications also to image reconstruction and to mathematical finance, 
has motivated a systematic study of the properties of the fractional Laplacian and pseudo-differential operators in general %. See for instance the papers
 \cite{CafSil1,CafSil2,CafSir,CabSir,Gar0,Sil0,SiV0,DpPV,FLS,FQT,ROS1,ROS2}; 
variational techniques have %also 
been employed also to obtain quantitative and qualitative results for elliptic PDEs with nonlocal nonlinearities
\cite{Len1, MZo0, MPT,MS0, CCS1, GV0, GMV, RV0, ClSa, Wet0, WX0,XW0}.

\smallskip

A key aspect in the study of partial differential equations %a key point 
consists also on the hypothesis assumed on the nonlinearity: considering very general ones allows to include different models coming from different frameworks. In 1983 Berestycki and Lions \cite{BL1,BL2} proposed a set of assumptions which relies, essentially, only on the growth of the nonlinearity in zero and at infinity: these assumptions may be considered, from a variational point of view, \emph{almost optimal}, and %type assumptions 
include for instance the most common power type functions $g(t) \sim t^p$, but also combined powers representing cooperation $g(t) \sim t^p+t^q$ and competition $g(t) \sim t^p-t^q$, as well as asymptotically linear saturable sources arising in nonlinear optics $g(t) \sim \frac{t^3}{1+t^2}$ and in the study of semiconductors $g(t) \sim t - \frac{t}{\sqrt{1+t^2}}$, and many others.
 The generality of these assumptions, which do not include regularity, homogeneity, Ambrosetti-Rabinowitz-type or monotonicity conditions, precludes the possibility of using classical tools of the variational analysis, such as minimization on Nehari manifolds and %Pohozaev 
fibering methods \cite{Neh1, Poh1, Wil0, KP0, BrZh0}, use of Pohozaev identities \cite{Poh0}, as well as boundedness of standard Palais-Smale sequences and classical Mountain Pass geometries \cite{AR0}.

Goal of this thesis is to investigate the abovementioned PDEs avoiding the use of these additional assumptions, examining especially how the geometry and the compactness of the problems can be tackled in this generality. 
In particular, we solve here also some problems which were left open in literature, and their resolution requires the implementation of new ideas.

%\smallskip

\bigskip

Studying equations \eqref{eq_introd_frac} and \eqref{eq_introd_choq}, % and \eqref{eq_introd_doubl}, 
%At this point (set \hbar=1) 
the research has been pursued essentially in two main directions: % the problems can be treated by two classical approaches: 
the first is to assign the frequency $\mu \in(0,+\infty)$, and let the mass (given by the $L^2$-norm of $u$) to be free. 
This \emph{unconstrained approach} has been extensively studied in the literature 
\cite{BL1,BL2,BGK,JT0,Med1, %local
ChWa0,BKS,Iko1,Iko2,MS2, MS3
}. 
A second approach is to prescribe the mass $\int_{\R^N} u^2 =m>0$ and let instead the frequency to be an unknown 
\cite{Shi0,HT0, LZ0, Lio1, BaLiLi }: 
this \emph{constrained approach} is also significantly meaningful in physics, for instance in quantum mechanics due to the normalization of probability.
% where information on the mass is of key importance. 

In this thesis we aim to find existence and multiplicity results for $L^2$-constrained problems, that is
$$\parag{ &(-\Delta)^s u + \mu u = g(u), \quad x \in \R^N, & \\ &\int_{\R^N} u^2 =m,&} \quad 
\parag{ &-\Delta u + \mu u= (I_{\alpha} \ast F(u)) F'(u), \quad x \in \R^N, & \\ &\int_{\R^N} u^2 =m.&} 
$$
When dealing with nonlocalities, the classical minimization approach on the $L^2$-sphere is rather involved, since the techniques require a delicate control on the tails of the functions; moreover, this approach is less suitable for the research of % especially when searching for 
multiple solutions. Here we propose instead a minimax approach, related to a \emph{Lagrangian formulation} of the problem, and modeled on suitable \emph{mountains on the product space}: we believe that this method may be applied to a wider class of equations. 
A posteriori, we show also that the found solution with minimal energy is indeed an $L^2$-minimum.
Even though the approach to the two problems is similar, % for the two problems, %Nevertheless, 
the study of the two abovementioned equations gives rise to different problems. 

A particular feature of the fractional Laplacian, indeed, is the lack of a regularizing effect: this fact does not allow to prove the well known Pohozaev identity, a quite useful tool in the framework of PDEs. This lack of regularity is here tackled by implementing a suitable modification of the Palais-Smale condition, that we call \emph{Palais-Smale-Pohozaev condition}, and a deformation argument around the set of critical points satisfying the Pohozaev identity. 
Here we face for the first time the problem of the existence of a normalized solution for a fractional framework, where the Pohozaev identity is no more ensured; moreover, we highlight that the multiplicity result presented is new even in the power setting $g(u)=|u|^{p-2}u$.
This is done in Chapter \ref{chap_fract_normal}.

In the case of Choquard %-Pekar 
nonlinearities, instead, a delicate issue is the research of multiple solutions: indeed, this is typically based on the construction of suitable multidimensional Mountain Pass paths. On the other hand, when the nonlinearity is not local, this is not obvious, and that is why we need to implement a delicate construction based on \emph{multidimensional annuli} which takes into account the interaction of far components. 
In particular, as a peculiar feature of the nonlocal setting we are allowed to consider \emph{odd} (and not only even) functions $F$, which make the energy functional symmetric as well: this possibility has not been developed in the common literature. Nevertheless, the case $F$ odd makes much more involved the control of far nonlocal contributions; here we include this case in our study.
% arises to some interesting phenomenon in the study of the interaction of different nonlocal contributions
Moreover, as a byproduct of this construction, we find existence of infinitely many solutions for the unconstrained Choquard problem \eqref{eq_introd_doubl}, solving a problem which was left open in the literature \cite{MS2} 
and extending to nonlocal nonlinearities the seminal paper by Berestycki and Lions \cite{BL2}. 
We do this in Chapter \ref{chap_choq_multi}.

\medskip

When studying \emph{fractional Choquard equations} \cite{DSS1} of the type \eqref{eq_introd_doubl}, the combination of the two nonlocalities and of the nonhomogeneous nonlinearity heavily influences the investigation of qualitative properties of the solutions. The lack of explicit computations, the absence of a proper chain rule and the singularities of the Fourier symbol and of the convolution kernel obstruct classical approaches in the study of boundedness, $L^1$-summability and regularity of solutions, as well as positivity and asymptotic decay of ground states. 
Again, also here we consider the possibility of $F$ to be odd in the study of some symmetry properties: 
all the abovementioned difficulties require new ideas and the implementation of more delicate arguments. 
Some of the cited results are, in addition, new even for the case $s=1$, improving some results in \cite{MS2}. 

 The nonlocal interaction of the fractional Laplacian and of the Choquard %-Pekar 
term gives rise moreover to new phenomena: for instance, when $F$ has a subquadratic growth in the origin, the asymptotic behaviour at infinity of the solutions seems to be connected to a new \emph{growth threshold}, differently from the local case $s=1$. All these properties are examined in Chapter \ref{chap_doubly}.

\medskip

Finally, another problem we aim to investigate is the \emph{concentration} of solutions in fractional nonlinear Schrödinger equations. Indeed, given an external potential $V=V(x)$, physicists are interested in studying the effect of this potential on the solutions of the equation
$$ \hbar^{2s} (-\Delta)^s u+V(x)u=g(u), \quad x\in\R^N$$
as long as the term $\hbar$ goes to zero, which somehow describes the passage from quantum to classical mechanics \cite{BrJe, SeSq}; this is why solutions of this equation for $\hbar>0$ small are also called \emph{semiclassical}. In particular it has been proved that, if a family of solutions has maxima which concentrate in a point, then that point is critical for $V$ \cite{Wan0, FMV}. This is the reason why a huge literature is focused on studying concentration on different types of critical points, both in a local framework \cite{BT1, CJT} and nonlocal \cite{CT0, FiSi0, AA0, Che0,Seo0}.
Our aim is to investigate concentration phenomena on local minima of $V$, in the framework of fractional equations: in this case, the spreading of the mass carried by the fractional Laplacian strongly opposes the research of solutions \emph{localized} in a prescribed domain of $\R^N$. Despite this obstruction, we find the existence of multiple solutions with this behaviour, whose number is related to some algebraic-topological information on the set of local minima of $V$. 

In order to achieve this, some careful analysis is needed: 
indeed, the possible degeneracy of the local minimum of $V$ does not allow to implement finite-dimensional reduction arguments, while the generality of the function $g$ hinders the possibility of working on natural constraints, such as Nehari manifolds. In order to study sets of local minima we combine perturbation and penalization arguments and implement delicate deformation theorems on some set of expected solutions. %, by exploitiong some algebraic tool to deduce multiplicity from the topology of $V$. 
In this discussion, we include a posteriori the case of a lost of compactness given by a Sobolev-critical growth of $g$, through the use of a truncation argument and suitable a priori estimates.

As already highlighted, the presence of a nonlocality makes the whole study much more involved: the lack of a proper Leibniz rule and of the preservation of the supports prevents the use of classical cut-off functions and standard penalization arguments. Moreover, a strong control on the tails of the functions is needed, especially when trying to localize their \emph{fractional center of mass}, and we do this by means of a suitable mixed fractional seminorm. 
This study is made in Chapter \ref{chap_concentr}.

\bigskip

The general spirit of this thesis is thus to investigate how the nonlocalities -- both in the operator and in the source -- comes into account in the study of different PDEs, and analyze these equations under almost optimal assumptions on the nonlinearity.

\bigskip

The thesis is organized as follows. 
In Chapter \ref{chap_prelim} we 
recall and revisit some known results in literature, furnishing the proofs whenever it was not possible to find a precise reference, 
% some facts about the fractional Laplacian and the Riesz potential: 
%Here we collect and revisit some known results in literature, furnishing some proofs whenever it was not possible to find a precise reference.
%some of the results are well known but a precise reference was not found, %while 
and we present some new results %are shown 
as well.
Chapter \ref{chap_fract_normal} is dedicated to the study of autonomous fractional equations: after having recalled what is known for the unconstrained problem, we focus on the study of the mass-constrained problem, obtaining both existence and multiplicity of solutions for general nonlinearities. % $g$.
Then, %similar ideas are applied 
in Chapter \ref{chap_choq_multi} we %to the 
research for multiple solutions to the Choquard %-Pekar 
problem: in this case, one of the main issues is the construction of suitable multidimensional paths, since the general and nonlocal nonlinearity heavily affects the geometry of the problem.
In Chapter \ref{chap_doubly} we move to study the case of doubly nonlocal equations, where we mainly focus on the qualitative properties of the solutions, investigating how the interaction of the two nonlocalities influences both the techniques and the results.
Finally, we face the fractional semiclassical problem in Chapter \ref{chap_concentr}, by studying how the nonlocality of the fractional operator comes into play while searching for multiple solutions concentrating to a local minimum of the potential.
Appendix \ref{chap_app_alg_top} is dedicated to a little survey % an introduction 
on the algebraic and topological tools used throughout the thesis.

\medskip

This thesis is mainly based on the papers \cite{CG0, CG1, CGT1, CGT2, CGT3, CGT4, CGT5, Gal0, Gal1}.

\mainmatter 

%%%%%%%%%%%%%%%%%%%%%%%%%%%%%%%%%%%%%%%%%%%%%%%%%%%%%%%%%%%%%%%%%%%%%%
%%%%%%%%%%%%%%%%%%%%%%%%%%%%%%%%%%%%%%%%%%%%%%%%%%%%%%%%%%%%%%%%%%%%%%
%%%%%%%%%%%%%%%%%%%%%%%%%%%%%%%%%%%%%%%%%%%%%%%%%%%%%%%%%%%%%%%%%%%%%%

\chapter{Some facts about nonlocalities}
\label{chap_prelim}

In this Chapter we introduce some preliminary results about the fractional Laplacian (Section \ref{sec_prelim_fract}) and the Riesz potential (Section \ref{sec_riesz_poten}), as well as some considerations about nonlinear functionals (Section \ref{sec_prelim_assump}). 
Here we collect and revisit some known results in literature, furnishing some proofs whenever it was not possible to find a precise reference.
%The recalled results are mainly known in literature, even if 
%%Most of them will be recalls of known results; 
%for some of them we where not able to find a precise reference, so some short proofs are provided. 
%Moreover, some new results are here presented (see for instance Section \ref{sec_chain_rule}).

Moreover, we present here some new results: in particular, in Section \ref{sec_chain_rule} we deal with a fractional chain rule in presence of concave compositions, by working with a viscosity formulation; this can be found in paper \cite{Gal1}.
In Section \ref{sec_regol_degiorgi} instead, we present an $L^{\infty}$-bound for non-positive solutions of fractional nonautonomous elliptic -- possibly critical -- equations, which adapts also to the %(superlinear) 
Choquard framework; this has been developed in papers \cite{Gal0,CGT3}.

%%%%%%%%%%%%%%%%%%%%%%%%%%%%%%%%%%%%%%%%%%%%%%%%%%%%%%%%%%%%%%%%%%%%%%
%%%%%%%%%%%%%%%%%%%%%%%%%%%%%%%%%%%%%%%%%%%%%%%%%%%%%%%%%%%%%%%%%%%%%%

\section{Notations}

We start by writing down some notations used throughout the thesis.
We write $ \R_+:=(0,+\infty)$ and
\begin{eqnarray*}
 && B_r(x_0):=B(x_0,r) :=\{ x\in\R^N \mid |x-x_0| <r\} \quad \hbox{ for $x_0 \in \R^N$ and $r>0$}, \\
 && D_N:=\{\xi\in\R^N \mid \abs\xi\leq 1\} \quad \hbox{ for $N \in \N^*$}, \\
 && A(R,h):=\big\{x\in\R^N \mid \abs x\in [R-h,R+h]\big\}, \quad \hbox{for $R>0$, $h>0$}
\end{eqnarray*}
for balls, disks, annuli; in particular, $B_r:=B_r(0)$, and $\chi(R,h;\cdot):= \chi_{A(R,h)}$. In addition, 
$$A_{\delta}:= \{x \in X \mid d(x,A)\leq \delta\}$$
denotes a neighborhood for any $A\subset (X,d)$ metric space. Sometimes we will write 
$$\complement(A):=A^c := X \setminus A $$
for $A \subset X$ to avoid cumbersome notation, if the ambient space is clear from the context. The function $P_i$ will denote, generally, the projection on the $i$-th component (in some product space).

We write 
\begin{eqnarray*}
 &&\norm u_r := \left(\intRN \abs{u}^r\, dx\right)^{1/r} \quad \hbox{for $r\in [1,\infty)$ and $u \in L^r(\R^N)$},\\
 &&\norm u_{\infty} := \textnormal{ess sup}_{\R^N} |u| \quad \hbox{for $u \in L^{\infty}(\R^N)$},
\end{eqnarray*}
the classical $L^p$-norm in the entire space, $p \in [1, +\infty]$; we will use also the following notation
$$\norm{f}_{\infty, \theta}:= \norm{f(\cdot )(1+|\cdot|^{\theta})}_{\infty}$$
for any $\theta>0$. By $\mc{F}(u)$ or $\widehat{u}$ we will denote, moreover, the Fourier transform of a function $u$, and by $u_{\pm}$ its positive and negative parts, $u=u_+-u_-$.

The function $\Gamma(\cdot)$ will denote the standard Gamma function, while ${}_2F_1(\cdot, \cdot, \cdot \; ;\, \cdot)$ will denote the Gauss hypergeometric function.

%\tor{
We write $\mc{S}$ for the Schwartz function space. 
%}
For any $k \in \N$ and $\sigma \in (0,1)$, we denote by $C_0(\R^N)$ the space of continuous functions decaying to zero at infinity, by $C_b^k(\R^N)$ (resp. $C^k_c(\R^N)$) the space of $k$ times differentiable functions with bounded (resp. compactly supported) and continuous $j$-derivative, $j=0, \dots, k$, by $C^{k,\sigma}(\R^N)$ the space of $k$ times differentiable functions with $\sigma$-H\"older continuous $k$-derivatives (on $\R^N$), where
$$[u]_{C^{0,\sigma}(A)}:= \sup_{\substack{x, y \in A \\ x\neq y}} \frac{\abs{u(x)-u(y)}}{\abs{x-y}^{\sigma}}$$
denotes the usual seminorm in H\"older spaces for $\sigma \in (0,1]$ and $A\subseteq \R^N$. %We allow $\sigma>1$ by simply writing $C^{\sigma}(\R^N)$.
By $C^{k,\sigma}_{loc}(\R^N)$ we consider functions whose $k$-derivatives are locally $\sigma$-H\"older continuous; if $\sigma=1$ we also write $Lip(\R^N):=C^{0,1}(\R^N)$ and similarly $Lip_{loc}(\R^N)$ and $Lip_{c}(\R^N)$.
More briefly we will sometimes write 
$$C^{\beta}(\R^N):= C^{[\beta], \beta-[\beta]}(\R^N)$$
for any $\beta>0$, observing that this notation throws out spaces $C^{k,1}(\R^N)$, usually subsituted by Zygmund spaces (see Remark \ref{rem_Zyg} below);
%\tr{così però non sto includendo $C^{k,1}(\R^N)$ (Lipschitziane)! Non va bene come notazione!}
%we convey
%$$C^{k,\sigma}(\R^N):=C^{k,1}(\R^N)$$
%\tr{anche questa la eliminerei come notazione!}
%for $\sigma > 1$, in order to avoid cumbersome notation; 
similarly $C^{\beta}_{loc}(\R^N)$. %for local spaces.

\begin{Remark}\label{rem_DuPl}
In \cite{DuP0} it is defined, for $\sigma \in (0,1]$, $u\in Lip(\sigma)$ if there exist $C> 0$ and $\delta>0$ such that, for each $x, y \in \R^N$,
$$0<|x-y|\leq \delta \implies \frac{|u(x)-u(y)|}{|x-y|^{\sigma}} \leq C.$$
We notice that
$$C^{0,\sigma}(\R^N) \subset Lip(\sigma) \subset C^{0,\sigma}_{loc}(\R^N)$$
and moreover
$$Lip(\sigma) \cap L^{\infty}(\R^N) \subset C^{0,\sigma}(\R^N);$$
indeed, %if $|u(x)|\leq M$ for each $x\in \R^N$, then, 
for each $x, y \in \R^N$,
%$$|x-y|> \delta \implies \frac{|u(x)-u(y)|}{|x-y|^{\sigma}} \leq \frac{2M}{\delta^{\sigma}}.$$
$$|x-y|> \delta \implies \frac{|u(x)-u(y)|}{|x-y|^{\sigma}} \leq \frac{2\norm{u}_{\infty}}{\delta^{\sigma}}.$$
\end{Remark}

\begin{Remark}\label{rem_Zyg}
To state some results it is useful to introduce also the \emph{Zygmund space} $\Lambda_1(\R^N)$ \cite[Section 6]{Sti0} as the space of the continuous functions $u$ such that
$$\sup_{x,h \in \R^N} \frac{|u(x+h)-2u(x) + u(x-h)|}{|h|} < \infty.$$
We notice that $u\in C^{0,\sigma}(\R^N)$ for $\sigma \in (0,1)$ if equivalently 
$$\sup_{x,h \in \R^N} \frac{|u(x+h)-2u(x) + u(x-h)|}{|h|^{\sigma}} < \infty,$$
but the same does not hold true for $\sigma=1$; indeed
$$C^{0,1}(\R^N) \subsetneq \Lambda_1(\R^N).$$
 We can further define $\Lambda_2(\R^N)$ as the space of functions in $C^1(\R^N)$ with partial derivatives in $\Lambda_1(\R^N)$; also in this case $C^{1,1}(\R^N) \subsetneq \Lambda_2(\R^N).$
The following relations hold true \cite[Propositions 5.5.8, 5.5.9 and 5.5.10]{Ste0}:
%$$C^{0,1}(\R^N) \cap L^{\infty}(\R^N) \subset \Lambda_1(\R^N) \cap L^{\infty}(\R^N) \subset C^{0,\sigma_2}(\R^N) \cap L^{\infty}(\R^N) \subset C^{0,\sigma_1}(\R^N) \cap L^{\infty}(\R^N) $$
%for each $0 < \sigma_1 < \sigma_2 < 1$.
%$$C^{1,1}(\R^N) \cap L^{\infty}(\R^N) \subset \Lambda_2(\R^N) \cap L^{\infty}(\R^N) \subset C^{1,\sigma_2}(\R^N) \cap L^{\infty}(\R^N) \subset C^{1,\sigma_1}(\R^N) \cap L^{\infty}(\R^N) $$
%for each $0 < \sigma_1 < \sigma_2 < 1$.
%$$C^{k,1}(\R^N) \cap L^{\infty}(\R^N) \subset \Lambda_{k+1}(\R^N) \cap L^{\infty}(\R^N) \subset C^{k,\sigma_2}(\R^N) \cap L^{\infty}(\R^N) \subset C^{k,\sigma_1}(\R^N) \cap L^{\infty}(\R^N) $$
$$C^{k,1} \cap L^{\infty} \subset \Lambda_{k+1} \cap L^{\infty} \subset C^{k,\sigma_2} \cap L^{\infty} \subset C^{k,\sigma_1} \cap L^{\infty} $$
for $k=0,1$ and each $0 < \sigma_1 < \sigma_2 < 1$.
\end{Remark}

Here we write $f\sim g$ as $x \to x_0 \in \overline{\R}$ if there exist constants $C_1$, $C_2>0$ independent of $x$ such that
 $$ C_1g(x) \leq f(x) \leq C_2g(x) \quad \hbox{for $x$ near $x_0$}, $$
while by $f \stackrel{.}\sim g$ as $x \to x_0$ we mean that 
$$\lim_{x \to x_0} \frac{f(x)}{g(x)}=1.$$

Moreover, by $\approx$ we will mean \emph{approximately equal to} (in a sense clear from the context) or \emph{isomorphic to}. Symbols $\lesssim$, $\simeq$ and $\gtrsim$ will mean less, equal or greater up to (positive) constants.

Finally, for every $A\subset B \subset \R^N$, we will write
$$A\prec \phi \prec B$$
to indicate a Urysohn-type regular function $\phi \in C^{\infty}_c(\R^N)$ such that
$$\phi_{|A}=1 \quad \hbox{ and } \quad \phi_{|\R^N \setminus B}=0.$$

We introduce the following terminology: if $\G$ is a group acting on a set $X$, we say that $A\subset X$ is \emph{invariant} under $\G$ if $g A = A$ for each $g \in \G$, while we say that a function $f:X\to Y$ ($Y$ another set) is \emph{invariant} under $\G$ if $f(g \cdot ) = f$ for each $g \in \G$; finally we say that $f: X \to X$ is \emph{equivariant} if $f(g \cdot ) = g \cdot f$ for each $g \in \G$.
When $\G=\Z_2 \equiv \{\pm 1\}$ acting on some vectorial space $X$, we have that $A\subset X$ invariant means symmetric with respect to the origin ($A=-A$), $f$ invariant means even ($f(-\cdot)=f$), $f$ equivariant means odd ($f(-\cdot )=-f$).

We highlight that, all throughout the thesis, we will actually assume $N\geq 2$ when dealing with the fractional framework $s \in (0,1)$ (despite the beginning of the preliminaries, where generally $N>2s$), and $N\geq 3$ in the local framework $s=1$.
Moreover the constants $C, C'$ appearing in inequalities may change from a passage to another; to avoid cumbersome notations, we will not stress the dependence of such constants, which will be based only on the fixed quantities in play. 

%%%%%%%%%%%%%%%%%%%%%%%%%%%%%%%%%%%%%%%%%%%%%%%%%%%%%%%%%%%%%%%%%%%%%%
%%%%%%%%%%%%%%%%%%%%%%%%%%%%%%%%%%%%%%%%%%%%%%%%%%%%%%%%%%%%%%%%%%%%%%

\section{The fractional Laplacian}
\label{sec_prelim_fract}

Let $s \in (0,1)$ and $N>2s$. For this Section we mainly refer to \cite{DnPV, Gar0}, together with \cite{AV0, BKS, FQT, Sil0}; other interesting references are \cite{Amb5, BRS, BuV, CabSir, DD0, SeV} (see also %the recent books 
\cite{ChLiMa}%, FRRO}
). For motivations and a physical introduction we refer to Sections \ref{sec_introd_frac_lap} and \ref{sec_boson_stars}.

Let the fractional Laplacian be defined by \cite{DnPV}
$$(-\Delta)^s u(x) := C_{N,s} \textnormal{PV} \int_{\R^N} \frac{u(x)-u(y)}{|x-y|^{N+2s}} dy$$
where 
$$C_{N,s}:=\frac{4^s \Gamma(\frac{N+2s}{2})}{\pi^{N/2} |\Gamma(-s)|}>0$$
is a normalization constant with $\Gamma$ the Gamma function, and the integral is in the Principal Value sense, that is
$$(-\Delta)^s u(x) = C_{N,s} \lim_{\eps \to 0^+} \int_{B_{\eps}^c(x)} \frac{u(x)-u(y)}{|x-y|^{N+2s}} dy;$$
notice that, when $s \in (0, \frac{1}{2})$, we actually do not need to employ the Principal Value formulation 
%(when $u$ is, for instance, a Schwartz function, see \cite[Remark 3.1]{DnPV}).
(when $u$ belongs, for instance, to $C^{0,\sigma}_{loc}(\R^N) \cap L^{\infty}(\R^N)$ for some $\sigma \in (2s,1]$ \cite[Remark 3.1]{DnPV}, see also the proof of Proposition \ref{prop_diff_represent} below).

A sufficient condition in order to have $(-\Delta)^s u$ well defined pointwise is given by \cite[Proposition 2.4]{Sil0} (see also \cite[Proposition 2.15]{Gar0} and \cite[Lemma 2.4]{BKS}).
\begin{Proposition}[Fractional well posedness]\label{prop_well_posed}
Let $x_0 \in \R^N$. Then, if
\begin{itemize}
\item $u \in L^p(\R^N)\cap C^{\gamma}(U)$ for some $p \in [1,+\infty]$, $\gamma >2s$ and $U$ open neighborhood of $x_0$,
\end{itemize}
then $(-\Delta)^s u(x_0)$ is well defined; in this case, actually, $(-\Delta)^s u \in C(U)$. In particular, $(-\Delta)^s u$ is everywhere well defined pointwise if
\begin{itemize}
\item $u \in L^p(\R^N) 
\cap C^{\gamma}_{loc}(\R^N)$ for some $p \in [1,+\infty]$ and $\gamma >2s$, 
\end{itemize}
and we have $(-\Delta)^s u \in C(\R^N)$. 
\end{Proposition}

Actually the condition $u \in L^p(\R^N)$ can be substituted by the more general condition
\begin{equation} \label{eq_spazio_gener_s}
\int_{\R^N} \frac{|u(x)|}{(1+|x|)^{N+2s}} < \infty.
\end{equation}

%%%

A different pointwise representation is given in the following Proposition \cite[Lemma 3.2]{DnPV} (see also \cite[Proposition 2.8]{Gar0}).

\begin{Proposition}
\label{prop_diff_represent}
%Let $s\in (0,1)$. 
Assume $u \in L^p(\R^N) 
\cap C^{\gamma}_{loc}(\R^N)$ for some $p \in [1,+\infty]$ and $\gamma >2s$, 
%\begin{equation} \label{eq_spazio_gener_s}
%\int_{\R^N} \frac{|u(x)|}{(1+|x|)^{N+2s}} dx< \infty;
%\end{equation}
%e.g., $u \in (L^1+ L^{\infty})(\R^N)$. 
%% $u \in L^p(\R^N)$ for some $p \in [1,+\infty]$. 
%Assume moreover $u \in C^{\gamma}_{loc}(\R^N)$ for some $\gamma >2s$.
%\begin{itemize}
%\item $u \in C^{\gamma}_{loc}(\R^N)$ for some $\gamma >2s$.
%\end{itemize}
%, if
%\begin{itemize}
%\item $u \in C^{\gamma}(U)$ for some $\gamma >2s$ and $U$ open neighborhood of $x_0$,
%\end{itemize}
%then $(-\Delta)^s u(x_0)$ is well defined; in this case, actually, $(-\Delta)^s u \in C(U)$. In particular, 
Then
$$
(-\Delta)^s u(x) =\frac{C_{N,s}}{2} \int_{\R^N}\frac{ 2u(x)-u(x+y)-u(x-y)}{|y|^{N+2s}} dy,
$$
%\begin{align*}
%(-\Delta)^s u(x) &= %C_{N,s} \textnormal{PV} \int_{\R^N} \frac{u(x)-u(y)}{|x-y|^{N+2s}} dy = 
%C_{N,s} \lim_{\eps \to 0} \int_{\R^N \setminus B_{\eps}(x)} \frac{u(x)-u(y)}{|x-y|^{N+2s}} dy \\
%&= \frac{C_{N,s}}{2} \int_{\R^N}\frac{ 2u(x)-u(x+y)-u(x-y)}{|y|^{N+2s}} dy
%\end{align*}
%for every $x \in \R^N$, and in this case we have $(-\Delta)^s u \in C(\R^N)$. Moreover the last 
and the integral is absolutely convergent.
\end{Proposition}

\claim Proof.
We check only the absolute convergence. Indeed, let $x \in \R^N$ and $R>2|x|+1$. Notice that, for $|y|\geq R$, we have, for $|y|\geq R$,
$$|x+y| \geq |y|-|x| \geq R - |x|>|x|+1 $$
and
$$|x+y| - |x| \geq \frac{|x+y|+1}{|x|+2} $$
thus
\begin{eqnarray*}
\lefteqn{\int_{B_R^c}\frac{ |2u(x)-u(x+y)-u(x-y)|}{|y|^{N+2s}} dy} \\
 &\leq& \int_{B_R^c}\frac{2 |u(x)|}{|y|^{N+2s}} dy + \int_{B_R^c}\frac{|u(x+y)|}{|y|^{N+2s}}dy + \int_{B_R^c}\frac{|u(x-y)|}{|y|^{N+2s}} dy \\
&\leq& 2 |u(x)| \int_{B_R^c}\frac{1}{|y|^{N+2s}} dy + 2\int_{B_{|x|+1}^c}\frac{|u(z)|}{(|z|-|x|)^{N+2s}}dz\\ %+ \int_{B_R^c}\frac{|u(x-y)|}{|y|^{N+2s}} dy
&\leq& C_R |u(x)| + 2(2+|x|)^{N+2s}\int_{B_{|x|+1}^c}\frac{|u(z)|}{(|z|+1)^{N+2s}}dz < \infty.
\end{eqnarray*}
Let now $s \in (0, \frac{1}{2})$. Then, being $u \in C^{0,\gamma}_{loc}(\R^N)$ for some $\gamma>2s$,
$$\int_{B_R}\frac{ |2u(x)-u(x+y)-u(x-y)|}{|y|^{N+2s}} dy \leq 2 C \int_{B_R}\frac{ 1}{|y|^{N+2s-\gamma}} dy < \infty;$$
notice that a similarly argument shows also that the integral in the definition of the fractional Laplacian %of the Proposition's claim 
does not need the Principal Value, being absolute convergent.

If instead $s \in [\frac{1}{2}, 1)$, then, being $u \in C^{1,\gamma}_{loc}(\R^N)$ for some $\gamma>2s-1$, for each $x,y \in \R^N$ there exists $\sigma=\sigma(x,y) \in (0,1)$ %, \tau=\tau(x,y)$ in $ (0,1)$ 
such that
%\begin{align*}
%\int_{B_R}\frac{ |2u(x)-u(x+y)-u(x-y)|}{|y|^{N+2s}} dy & = \int_{B_R}\frac{ |\nabla u(x+\tau y)\cdot y - \nabla u(x-\sigma y) \cdot y|}{|y|^{N+2s}} dy \\
%& \leq \int_{B_R}\frac{|\tau + \sigma|}{|y|^{N+2s-\gamma-1}} dy < \infty.
%\end{align*}
\begin{align*}
\int_{B_R}\frac{ |2u(x)-u(x+y)-u(x-y)|}{|y|^{N+2s}} dy & = \int_{B_R}\frac{ |\nabla u(x+\sigma y)\cdot y - \nabla u(x-\sigma y) \cdot y|}{|y|^{N+2s}} dy \\
& \leq \int_{B_R}\frac{2\sigma}{|y|^{N+2s-\gamma-1}} dy < \infty.
\end{align*}
Joining the pieces, we have the claim.
\QED

\bigskip

%%%

It is well known that the fractional Laplacian is a \emph{nonlocal} operator. This means, for example, that 
$$\supp\big((-\Delta)^s u\big) \not \subset \supp(u);$$
in particular, if $\psi$ is a cut-off function with support in some $A\subset \R^N$, we cannot localize $(-\Delta)^s (\psi u)$ inside $A$ as well. Notice that the fact that $(-\Delta)^s u$ is expressed through an integral does not directly implies that the operator is nonlocal (see, for instance, \cite[Section 2.1]{Aba0}); anyway we can see this considering, for example, a nonnegative $u\in C^2_c(\R^N)$ with $u\geq 1$ on $B_1(0)$, and a point $x\in \R^N$ far from the support: we thus have
%$$(-\Delta)^s u(x) %= \int_{\R^N} \frac{u(x)-u(y)}{|x-y|^{N+2s}} dy %= -\int_{\R^N} \frac{u(y)}{|x-y|^{N+2s}} dy 
%\leq -\int_{B_1(0)} \frac{u(y)}{|x-y|^{N+2s}} dy %\leq -\int_{B_1(0)} \frac{1}{|x-y|^{N+2s}} dy 
%\leq -\int_{B_1(0)} \frac{1}{(1+|x|)^{N+2s}} dy = - \frac{C}{(1+|x|)^{N+2s}} <0.$$
\begin{align*}
(-\Delta)^s u(x) %= \int_{\R^N} \frac{u(x)-u(y)}{|x-y|^{N+2s}} dy %= -\int_{\R^N} \frac{u(y)}{|x-y|^{N+2s}} dy 
&\leq -\int_{B_1(0)} \frac{u(y)}{|x-y|^{N+2s}} dy %\leq -\int_{B_1(0)} \frac{1}{|x-y|^{N+2s}} dy 
\leq -\int_{B_1(0)} \frac{1}{(1+|x|)^{N+2s}} dy \\
&= - \frac{C}{(1+|x|)^{N+2s}} <0.
\end{align*}
 Moreover, a proper Leibniz rule lacks in this framework, thus in general
$$(-\Delta)^{s/2}(\psi u) \neq (-\Delta)^{s/2} u \psi + (-\Delta)^{s/2} \psi u,$$
formula which instead holds when $(-\Delta)^{s/2}$ is substituted with the gradient $\nabla$ in the local framework $s=1$ (see Remark \ref{rem_grad_fract12}). In the fractional framework a correction term is needed \cite[Proposition 1.5]{BWZ}
$$(-\Delta)^{s/2}(\psi u) = (-\Delta)^{s/2} u \psi + (-\Delta)^{s/2} \psi u + C_{N,s} \int_{\R^N} \frac{\big(u(x)-u(y)\big)\big(\psi(x)-\psi(y)\big)}{|x-y|^{N+s}} dy$$
or different approaches, like error estimates \cite{GO0} or approximation arguments \cite[Lemma 2.6]{Sec0} must be employed.
All these issues create problems, for example, in concentration arguments (see Chapter \ref{chap_concentr}). A proper chain rule lacks as well, and we will make some comments in Section \ref{sec_chain_rule}.

The operator anyway enjoys some trivial but useful scaling properties
$$(-\Delta)^s(\lambda u) = \lambda (-\Delta)^s u, \quad (-\Delta)^s \big( u(\beta \cdot)\big) = |\beta|^{2s} \big((-\Delta)^s u\big)(\beta \cdot).$$
for any $\lambda, \beta \in \R$, as well as linearity.
%$$(-\Delta)^s(u+v) = (-\Delta)^s u + (-\Delta)^s v.$$

\medskip

We further have the following relation with the Fourier transform \cite[Proposition 3.3]{DnPV} (see also \cite[Proposition 2.8]{Gar0}) 
whenever $u \in \mc{S}$ % (e più in generale?)}
\begin{equation}\label{eq_fourier_def} 
(-\Delta)^s u = \mc{F}^{-1}(|\xi|^{2s} \mc{F}(u));
\end{equation}
this relation can be extended to the setting of Proposition \ref{prop_well_posed}, that is for functions $u\in C^{\gamma}_{loc}(\R^N) \cap L^p(\R^N)$ for $\gamma>2s$, see \cite[proof of Proposition 2.4]{Sil0} (see also \cite[Lemma 2.4]{BKS}).

When $u$ is not regular enough, relation \eqref{eq_fourier_def} might be taken as a definition, whenever for example $|\xi|^{2s} \mc{F}(u) \in L^2(\R^N)$; in this case the fractional Laplacian is defined up to a set of zero Lebesgue measure.
Notice moreover that \eqref{eq_fourier_def} could be interpreted more generally also in the sense of tempered distributions $\mc{S}'$.

%CCCOMMENT NOW
\begin{Remark}
We notice that relation \eqref{eq_fourier_def}, i.e.
$$(-\Delta)^s u (x) = \int_{\R^N} (|\xi|^2)^{s} (u, e^{i \xi \cdot})_2 \, e^{i \xi \cdot x} d\xi$$
for almost every $x \in \R^N$, can be interpreted in terms of the spectral theorem by considering the continuum of eigenvalues $\xi \in \R^N \mapsto \lambda_{\xi}:=|\xi|^2$ with eigenfunctions $e_{\xi}(x):=e^{i \xi \cdot x} \in L^{\infty}(\R^N)$, and applying the power function $h(t):=t^s$. This is indeed how the \emph{spectral fractional Laplacian} is defined on bounded sets (see \cite[Section 2.3]{AV0}).
\end{Remark}

\begin{Remark}
Actually there are several equivalent ways to define the fractional Laplacian \cite{Kwa0}. One of the most used is the Caffarelli-Silvestre \emph{$s$-harmonic extension}, where the fractional Laplacian in $\R^N$ is seen as the trace of a divergence-form operator (possibly singular) in $\R^{N+1}$ \cite{CafSil1}: this formulation is widely used in order to bring the computations from a nonlocal framework to a local framework. Anyway we stress that we will not make use of the $s$-harmonic extension in this thesis, by mean of working directly with integral quantities. This has the further advantage of possibly extending our results to other nonlocal frameworks where the harmonic extension is no more available, see e.g. \cite{FaSi0}. %, for instance the $p$-fractional Laplacian .
%COMMENT NOW
%\tr{for instance...? INSERISCI ESEMPIO} 
\end{Remark}

Relation \eqref{eq_fourier_def} shows, informally, that
$$(-\Delta)^s u \stackrel{s \to 0^+} \to u, \quad (-\Delta)^s u \stackrel{s \to 1^-} \to -\Delta u$$
which motivates the symbol with a fractional power of the Laplacian; see \cite[Theorems 3 and 4]{Sti0} for a precise statement (see also \cite[Proposition 4.4]{DnPV}).

Moreover, \eqref{eq_fourier_def} is suitable to extend the notion of fractional Laplacian to every $s>0$ \cite{Sil0, ROS3, BMP0,CFW} (see also \cite[Proposition 3.1]{AJS2}); see \cite{Aba0} for an overview on the topic (see in particular \cite{AJS1,AJS2,AJS3,AJS4} and \cite[Section 3.1]{Sam0} for a hypersingular integral definition, \cite{CDQ} for a recursive pointwise definition, \cite{GFR} for a harmonic-extension definition).

\medskip

Another feature of the fractional Laplacian is its \emph{polynomial decay}, that is, whenever $u$ is good enough (for example, Schwartz), then \cite[Proposition 2.9]{Gar0} (see also Remark \ref{rem_decay_of_frac})
\begin{equation}\label{eq_stima_classica_s}
|(-\Delta)^s u(x)| \leq \frac{C}{1+ |x|^{N+2s}} \quad x \in \R^N;
\end{equation}
%this in particular implies that $(-\Delta)^s u$ is bounded and summable.
generally, one can not expect a faster decay: 
this is the case, for example, of $u(x)=\frac{1}{(1+|x|^2)^{\frac{N-2s}{2}}}$ (see Section \ref{sec_hypergeo}, and also \cite[Lemma 8.6]{Gar0} and \cite[Proposition 2.12]{Sil0}). 
Even when $u$ is a Schwartz function, by \eqref{eq_fourier_def} we notice that $(-\Delta)^s u$ has generally not a fast decay, since $|\xi|^{2s}$ is not regular enough near zero when $s<1$; thus 
$$(-\Delta)^s \mc{S} \not \subset \mc{S} \quad \hbox{ for $s \in (0,1)$}.$$
On the other hand, one can show \cite[Lemma 8.1]{Gar0} that, for every $u \in \mc{S}$, one has $(-\Delta)^s u \in C^{\infty}(\R^N)$ with
$$|D^{\beta}\big((-\Delta)^s u\big)(x)| \leq \frac{C}{1+ |x|^{N+2s}} \quad x \in \R^N$$
for each multi-index $\beta$.
We find the asymptotic decay \eqref{eq_stima_classica_s} also in fundamental solutions of the operator $(-\Delta)^s + id$ (see Lemma \ref{lem_Bessel_kernel}) and actually it will be a key feature of the solutions of fractional PDEs (see Section \ref{sez_limit_eq}), at least when there is not a too strong nonlocal effect coming from the nonlinearity (see Section \ref{sec_sublin_case}).

%%%%%%%%%%%%%%%%%%%%%%%%%%%%%%%%%%%%%%%%%%%%%%%%%%%%%%%%%%%%%%%%%%%%%%

\subsection{Fractional Sobolev spaces}
\label{sec_prel_sobolev}

We introduce, for any $\Omega \subseteq \R^N$ and $s \in (0,1)$, the fractional Sobolev space
$$H^s(\Omega) := \left\{ u \in L^2(\Omega) \mid [u]_{H^s(\Omega)}^2:=\int_{\Omega} \int_{\Omega} \frac{|u(x)-u(y)|^2}{|x-y|^{N+2s}}dy < +\infty \right \},$$
endowed with
$$\norm{u}_{H^s(\Omega)}^2:= \norm{u}_{L^2(\Omega)}^2 + [u]_{H^s(\Omega)}^2.$$
The finite quantity $ [u]_{H^s(\Omega)}$ is said Gagliardo seminorm. We will denote the dual space by %$H^{-s}(\Omega) = (H^s)^*(\Omega) = 
$(H^s(\Omega))^*$. %\tr{UNIFORMA} 
%\begin{Remark}
%If $u\in H^s(\R^N)$ and $h: \R \to \R$ is a Lipschitz function with $h(0)=0$, then $h(u) \in H^s(\R^N)$. Indeed
%$$\norm{h(u)}_2^2 = \int_{\R^N}|h(u)-h(0)|^2 \, dx \leq \int_{\R^N} \norm{h'}_{\infty}^2 |u-0|^2 \, dx = \norm{h'}_{\infty}^2 \norm{u}_2^2$$
%and
%$$%\norm{(-\Delta)^{s/2} h(u)}_2^2 
%[h(u)]_{\R^N}\leq C(N,S) \int_{\R^{2N}} \frac{\norm{h'}_{\infty}^2 |u(x)-u(y)|^2}{|x-y|^{N+2s}} \, dx \, dy = \norm{h'}_{\infty}^2 [u]_{\R^N}
%%\norm{(-\Delta)^{s/2} u}_2^2
%.$$
%\end{Remark}

We recall, whenever $\Omega=\R^N$ or $\Omega$ has a Lipschitz and bounded boundary, the continuous embedding \cite[Theorem 6.7]{DnPV}
\begin{equation}\label{eq_embedd_cont_Hs}
H^s(\Omega) \hookrightarrow L^p(\Omega)
\end{equation}
for every $p \in [2, 2^*_s]$ with
$$2^*_s = \frac{2N}{N-2s}$$
the fractional Sobolev critical exponent, 
and, if $p\in [2,2^*_s)$, the compact embedding \cite[Corollary 7.2]{DnPV}
$$H^s(\R^N) \hookrightarrow \hookrightarrow L^p_{loc}(\R^N)$$
in the sense that for every $(u_n)_n$ bounded in $H^s(\R^N)$, and for every $A\subset \R^N$ bounded and regular enough (e.g. $\partial A$ Lipschitz), we have that $(u_n)_n$ restricted to $A$ admits a convergent subsequence in $L^p(A)$.

Moreover we set 
$$H^s_{loc}(\R^N):= \left\{ u: \R^N \to \R \mid u \in H^s(\Omega) \hbox{ for each $\Omega \subset \subset \R^N$}\right\}$$
and, for any $\Omega \subset \R^N$, \cite[Section 4.3.2]{Tri0} 
\begin{align*}
X^s_0(\Omega) := & \left\{ w \in H^s(\R^N) \mid w=0 \; \hbox{on $\Omega^c$}\right\} \\
 = & \left\{ w \in H^s(\R^N) \mid \supp(w) \subset \overline{\Omega}\right\}.
\end{align*}

\begin{Remark}
The following density result holds in $\R^N$ \cite[Proposition 4.27]{DD0} (see also \cite[Proposition 4.11]{DD0}):
$$H^s(\R^N) = \overline{C^{\infty}_c(\R^N)}^{\norm{\cdot}_{H^s(\R^N)}}.$$
Assume now $\Omega$, with $\partial \Omega$ compact, to be a Lipschitz domain \cite[Definition 3.28]{McL0}.
Then \cite[Theorem 3.29]{McL0}
%When $\Omega$ is a continuous domain (for instance $\partial \Omega$ compact), then
$$X^s_0(\Omega) = \overline{C^{\infty}_c(\Omega)}^{\norm{\cdot}_{H^s(\R^N)}}.$$
If moreover %$\partial \Omega$ is Lipschitz and 
$s \neq \frac{1}{2}$, then \cite[Theorem 3.33]{McL0}
$$X^s_0(\Omega) = \overline{C^{\infty}_c(\Omega)}^{\norm{\cdot}_{H^s(\Omega)}}.$$
See also \cite[Theorem 1 in Section 4.3.2]{Tri0} for more results on these spaces.
\end{Remark}

In the case $\Omega = \R^N$ we also have the following relation \cite[Proposition 3.4%3.6
]{DnPV} 
$$ [u]_{H^s(\R^N)}^2 = \tfrac{2}{C_{N,s}} \norm{|\xi|^s \widehat{u}}_2^2; $$
by interpreting the fractional Laplacian through the Fourier transform definition \eqref{eq_fourier_def} we may also write
\begin{equation}\label{eq_semin_gagl}
 [u]_{H^s(\R^N)}^2 = \tfrac{2}{C_{N,s}} \norm{(-\Delta)^{s/2}u}_2^2. 
\end{equation}
%$$\norm{(-\Delta)^{s/2} u}_{L^2(\R^N)}^2 = \tfrac{1}{2} C_{N,s} [u]_{H^s(\R^N)}^2$$
%\tr{di sicuro in $H^s$ vale $[u]_{H^s(\R^N)}^2 = \norm{|\xi|^{s} \widehat{u}}_2^2$ (a meno di costanti); il membro di destre eguaglia il laplaciano in senso di Fourier. Ma per eguagliare il Laplaciano in senso integrale, ho bisogno di $u\in \mc{S}$... CAPISCI}
Moreover, by polarization
\begin{align}
 \int_{\R^N} (-\Delta)^{s/2} u (-\Delta)^{s/2} v dx &= \int_{\R^N} |\xi|^{2s} \widehat{u} \widehat{v} d\xi \notag \\
&= \tfrac{1}{2} C_{N,s} \int_{\R^N} \int_{\R^N} \frac{\big(u(x)-u(y)\big)\big(v(x)-v(u)\big)}{|x-y|^{N+2s}} dx dy \label{eq_sobolev_polaraz}
\end{align}
%\begin{equation}\label{eq_sobolev_polaraz}
%\int_{\R^N} \int_{\R^N} \frac{\big(u(x)-u(y)\big)\big(v(x)-v(y)\big)}{|x-y|^{N+2s}} dx \, dy = \tfrac{1}{2} C_{N,s} \int_{\R^N} (-\Delta)^{s/2} u \, (-\Delta)^{s/2} v \,dx
%\end{equation}
for every $u, \, v \in H^s(\R^N)$. Relation \eqref{eq_semin_gagl} leads also to an equivalent definition for the fractional Sobolev space
\begin{align*}
H^s(\R^N) &= \left\{ u \in L^2(\R^N) \mid |\xi|^{s} \widehat{u} \in L^2(\R^N) \right \} \\
 &= \left\{ u \in L^2(\R^N) \mid (-\Delta)^{s/2} u \in L^2(\R^N) \right \} 
\end{align*}
endowed with
\begin{align*}
\norm{u}_{H^s}^2 &\equiv \norm{u}_2^2 + \norm{|\xi|^{s} \widehat{u}}_2^2 \\
&= \norm{u}_2^2 + \norm{(-\Delta)^{s/2}u}_2^2.
\end{align*}

Together with $H^s(\R^N) \hookrightarrow L^2(\R^N) \cap L^{2^*_s}(\R^N)$
we have the following embedding of the homogeneous fractional space $\dot{H}^s(\R^N) \hookrightarrow L^{2^*_s}(\R^N)$ \cite[Theorem 6.5]{DnPV} (see also \cite{CFW}), where 
$$\dot{H}^s(\R^N):= \big\{ u \hbox{ measurable} \mid (-\Delta)^{s/2} u \in L^2(\R^N)\big\};$$
here the fractional Laplacian is intended in the sense of tempered distributions. That is, for some optimal constant $\mc{S}>0$,
\begin{equation}\label{eq_embd_homog} %\label{eq_embd_homog}
\norm{u}_{2^*_s} \leq \mc{S}^{-1/2} \norm{(-\Delta)^{s/2}u}_2.
\end{equation}

Moreover, we recall the fractional version of the Gagliardo-Nirenberg inequality \cite{Par1} (see also \cite{BGMMV})
\begin{equation}\label{GN}
\| u \|_r \leq C \|(-\Delta)^{s/2} u \|^{\beta}_2 \ \|u \|_2^{1-\beta}	
\end{equation}
for $u \in H^s(\R^N)$, 
$r \in [2,2^*_s]$ and $\beta$ satisfying 
$$\frac{1}{r} =\frac{\beta}{2^{^*}_s} + \frac{1-\beta}{2}.$$

\subsubsection{Extension to \boldmath{$p\in [1,\infty]$} and \boldmath{$s>0$}}

%Relation \eqref{eq_fourier_def} leads to the following formulation via Fourier transform 
Consider now again the relation
\begin{align}
H^s(\R^N) &= \left\{ u \in L^2(\R^N) \mid \mc{F}^{-1}\big(|\xi|^{s} \widehat{u}\big) \in L^2(\R^N) \right \} \notag \\
 &= \left\{ u \in L^2(\R^N) \mid \mc{F}^{-1}\big((1+|\xi|^{s}) \widehat{u}\big) \in L^2(\R^N) \right \}. \label{eq_def_sob_Hs}
\end{align}
%endowed with
%$$\norm{u}_{H^s}^2 \equiv\norm{u}_2^2 + \norm{|\xi|^{s} \widehat{u}}_2^2.$$
This last expression is suitable for defining the fractional Sobolev space $W^{s,p}(\R^N)$ also for $s\geq 1$ and $p\geq 1$, by \cite{FQT}
\begin{equation}\label{eq_def_sob_Ws1}
W^{s,p}(\R^N) := \big\{ u \in L^p(\R^N) \mid \mc{F}^{-1}\big((1+|\xi|^{s})\widehat{u}\big) \in L^p(\R^N)\big\}.
\end{equation}
It has been proved in \cite[Theorem 3.1]{FQT} that this definition coincide with the following
\begin{equation*}\label{eq_def_sob_Ws2}
\overline{W}^{s,p}(\R^N) := \left\{ u \in L^p(\R^N) \mid \mc{F}^{-1}\left((1+|\xi|^2)^{s/2}\widehat{u}\right) \in L^p(\R^N)\right\}
\end{equation*}
that is
$$W^{s,p}(\R^N) \equiv \overline{W}^{s,p}(\R^N).$$

\begin{Remark}
\label{rem_Bessel_kernels}
We want to highlight that the last equality is actually not trivial a priori. Indeed, we can rewrite the spaces as 
%\tr{(controlla formalmente le inversioni e le convoluzioni di Fourier in spazi $L^p$ e di distribuzioni temperate)} %CCOMMENT NOW
$$W^{s,p}(\R^N) \equiv \big\{ u \in L^p(\R^N) \mid u=\mc{K}_{2s}*g \hbox{ for some $g \in L^p(\R^N)$}\big\},$$
$$\overline{W}^{s,p}(\R^N) \equiv \big\{ u \in L^p(\R^N) \mid u=\mc{G}_{2s}*g \hbox{ for some $g \in L^p(\R^N)$}\big\}$$
where
$$\mc{K}_{2s}:= \mc{F}^{-1}\left(\frac{1}{1+|\xi|^{2s}}\right)$$
is the \emph{Bessel kernel}, and
$$\mc{G}_{2s}:= \mc{F}^{-1}\left( \frac{1}{(1+|\xi|^2)^s}\right)$$
is the \emph{pseudorelativistic kernel}.
%are both called \emph{Bessel kernels}, 
%\tr{(Trova nomi diversi!)} %COMMENT NOW
The two functions are the fundamental solutions, respectively, of
$$(-\Delta)^s \mc{K}_{2s} + \mc{K}_{2s} = \delta_0, \quad (-\Delta + id)^s \mc{G}_{2s}=\delta_0$$
in $\R^N$, where $\delta_0$ is the Dirac delta; the operator $(-\Delta + id)^s$ is also called \emph{pseudorelativistic} operator (see Section \ref{sec_boson_stars}). Even if 
$$\xi \mapsto \frac{1}{1+|\xi|^{2s}} \quad \hbox{and} \quad \xi \mapsto \frac{1}{(1+|\xi|^2)^s}$$
have same behaviour in zero and at infinity and same summability, the fact that the two functions have different regularity (the first is nonregular in the origin, the second is analytic) brings $\mc{K}_{2s}$ and $\mc{G}_{2s}$ to be quite different kernels: for instance, $\mc{K}_{2s}$ has a polynomial decay at infinity (of order $\frac{1}{|x|^{N+2s}}$, see Lemma \ref{lem_Bessel_kernel}), while $\mc{G}_{2s}$ decays exponentially \cite[equation (1.2.15)]{AH0}. 
These properties influence the qualitative behaviours of %may be inherited by 
the solutions of the linear equations 
%\tr{(spiega meglio)} %CCOMMENT NOW
$$(-\Delta)^s u + u = g, \quad (-\Delta + id)^s u =g$$
in $\R^N$, given by $u=\mc{K}_{2s}*g $ and $u=\mc{G}_{2s}*g$ respectively (see e.g. Lemma \ref{lem_Bessel_kernel}). % and \cite{GalSch}).
 Because of these representation formulas, we also write
$$\mc{K}_{2s}* \equiv \big((-\Delta)^s + id)^{-1}, \quad \mc{G}_{2s}* \equiv (-\Delta + id)^{-s}.$$
These considerations also show that the pseudorelativistic operator is quite different from the fractional Laplacian by giving more regularity and decay to solutions, but without enjoying the same scaling properties; its study is an interesting line of research for the future.
%does not enjoy the same scaling properties of $(-\Delta)^s+id$: in this thesis we will not deal with the pseudorelativistic operator, but understanding if our results can be adapted to this different framework is a good of research to be investigated in the future}.
\end{Remark}

\begin{Remark}
Notice that in \eqref{eq_def_sob_Hs} and \eqref{eq_def_sob_Ws1} the request $u \in L^p(\R^N)$ is actually superfluous. This is the same for $\overline{W}^{s,p}(\R^N)$ as well, since by \cite[Theorem 1.2.4]{AH0} we have (if $N>ps$)
%\begin{equation}\label{eq_Bessel_guadagna}
$$
 \mc{F}^{-1}\left((1+|\xi|^2)^{s/2}\widehat{u}\right) \in L^p(\R^N) \implies u \in L^q(\R^N) \hbox{ for each $q \in [p, \frac{pN}{N-ps}]$};
$$
%\end{equation}
this result is in accordance to the continuous embeddings \eqref{eq_embedd_cont_Hs} stated before for $p=2$. 
In particular the previous embedding is \emph{continuous}, which means that (for $q=p$)
$$\norm{u}_p \leq C \norm{\mc{F}^{-1}\big((1+|\xi|^2)^{s/2}\widehat{u}\big)}_p;$$
this relation can be rephrased by saying that 
%\tr{vedi meglio l'inversione} %COMMENT NOW
$$\norm{\big(-\Delta + id)^{-s} u}_p = \norm{\mc{G}_{2s}*u}_p \leq C \norm{u}_p$$
and this can be obtained directly by Young's inequality with $ C= \norm{\mc{G}_{2s}}_1 $ (indeed $\mc{G}_{2s} \in L^1(\R^N)$, see \cite[equation (1.2.12)]{AH0}).
A similar argument holds for $ ((-\Delta)^s + id)^{-1}$, since 
$$\norm{\big((-\Delta)^s + id)^{-1} u}_p = \norm{\mc{K}_{2s}*u}_p \leq \norm{\mc{K}_{2s}}_1 \norm{u}_p$$
being $\mc{K}_{2s} \in L^1(\R^N)$ (see Lemma \ref{lem_Bessel_kernel}), thus
\begin{equation}\label{eq_operat_contin_K}
\big((-\Delta)^s + id\big)^{-1} :L^p(\R^N) \to L^p(\R^N)
\end{equation}
is a continuous operator for every $p \in (1,+\infty)$ 
%(\tr{vedi meglio la restrizione}) %COMMENT NOW
and 
%\tr{vedi meglio l'inversione (se $u\in L^2(\R^N)$ dovrebbe bastare)} %COMMENT NOW
%L'inversione vale nel senso delle distribuzioni temperate. Più nello specifico, se $u\in L^2$ (come sarà quando applichiamo questo risultato), l'inversione è nel senso classico delle funzioni. 
\begin{equation}\label{eq_Bessel_guadagna} 
\norm{u}_p \leq \norm{\mc{K}_{2s}}_1 \norm{\mc{F}^{-1}\big(1+|\xi|^{2s}\big)\widehat{u}}_p.
\end{equation}
\end{Remark}

\medskip

We observe, by \eqref{eq_fourier_def}, that if $u\in W^{2s,p}(\R^N)$ for some $p$, then $(-\Delta)^s u$ is well defined pointwise up to a set of zero Lebesgue measure. 
\\Moreover, by \cite[Theorem 3.2]{FQT} we obtain the following embedding, 
% for every $s \in (0,2]$ and $p \in [1,+\infty)$,
%$$W^{s,p}(\R^N) \hookrightarrow C^{\gamma}(\R^N) 	\quad \hbox{for $\gamma \in (0, s - \frac{N}{p}]$}; $$
%in particular this implies, 
for every $s \in (0,1)$, 
\begin{equation}\label{eq_immers_Holder}
%H^{2s}(\R^N) \cap W^{2s,\infty}(\R^N) \hookrightarrow C^{\gamma}(\R^N)\quad \hbox{for $\gamma \in (0, 2s)$}.
H^{2s}(\R^N) \cap W^{2s,\infty}(\R^N) \hookrightarrow %C^{\gamma}(\R^N) = 
\parag{ 
&C^{0,\gamma}(\R^N) \; \hbox{for $\gamma \in (0, 2s)$} & \quad \hbox{ if $2s\leq 1$}, \\
&C^{1,\gamma-1}(\R^N) \; \hbox{for $\gamma \in (0, 2s)$} & \quad \hbox{ if $2s> 1$}.
}
%Not enough per la definizione puntuale integrale di $(-\Delta)^s$.
\end{equation}

\begin{Remark}\label{rem_emb_e_spezz}
We observe that, if $s\geq s'$ and $p \in (1,\infty)$, then \cite[equation (9) in Section 2.3.3]{Tri0}
$$W^{s,p}(\R^N) \hookrightarrow W^{s',p}(\R^N);$$
this is easily seen for $p=2$: 
indeed, by the fact that $|\xi|^{2s'} \leq 1+ |\xi|^{2s}$ we have 
%for $\frac{1}{p} + \frac{1}{p'}=1$ ($p \in [1,2]$)
%$$%\norm{u}_{W^{s',p}(\R^N)}^p = 
%\int_{\R^N} |(1 + |\xi|^{2s'}) \widehat{u}|^{p'} \leq \int_{\R^N} |(2 + |\xi|^{2s}) \widehat{u}|^{p'} \leq 2^{p'} \int_{\R^N} |(1 + |\xi|^{2s}) \widehat{u}|^{p'}.$$% \leq 2 \norm{u}_{W^{s,p}(\R^N)}^p$$
$$%\norm{u}_{W^{s',p}(\R^N)}^p = 
\int_{\R^N} |(1 + |\xi|^{2s'}) \widehat{u}|^{2} \leq \int_{\R^N} |(2 + |\xi|^{2s}) \widehat{u}|^{2} \leq 4 \int_{\R^N} |(1 + |\xi|^{2s}) \widehat{u}|^{2}.$$% \leq 2 \norm{u}_{W^{s,p}(\R^N)}^p$$
In particular, %whenever $2s\geq 1$, we have
$$H^{2s}(\R^N) \hookrightarrow H^1(\R^N) \quad \hbox{for $2s\geq 1$}.$$
Moreover, for every $s>0$, since %$2([s]-s) \leq 2[s]\leq 2s$, which means 
$H^{2s}(\R^N) \hookrightarrow H^{2[s]}(\R^N) \hookrightarrow H^{2([s]-s)}(\R^N)$, we notice that, for $u \in H^{2s}(\R^N)$,
%since We further notice that, for $u \in H^{2s}(\R^N)$ we have
%\begin{align*}
%(-\Delta)^{\alpha/2} u &= \mc{F}^{-1}\left( |\xi|^{\alpha} \widehat{u}\right) \\
%& = \mc{F}^{-1}\left( |\xi|^{[\alpha]-\alpha} |\xi|^{[\alpha]} \widehat{u}\right) \\
%&=\mc{F}^{-1}\left( |\xi|^{[\alpha]-\alpha} \mc{F}\Big(\mc{F}^{-1} \big(|\xi|^{[\alpha]} \widehat{u}\big)\Big)\right) \\
%&=(-\Delta)^{\frac{[\alpha]-\alpha}{2}} (-\Delta)^{\frac{[\alpha]}{2}} u
%\end{align*}
\begin{align*}
(-\Delta)^{s} u %&= \mc{F}^{-1}\left( |\xi|^{2s} \widehat{u}\right) \\
& = \mc{F}^{-1}\left( |\xi|^{2([s]-s)} |\xi|^{2[s]} \widehat{u}\right) \\
&=\mc{F}^{-1}\left( |\xi|^{2([s]-s)} \mc{F}\Big(\mc{F}^{-1} \big(|\xi|^{2[s]} \widehat{u}\big)\Big)\right) \\
&=(-\Delta)^{[s]-s}\left( (-\Delta)^{[s]} u\right)
\end{align*}
and similarly
$$ (-\Delta)^{s} u = (-\Delta)^{[s]}\left( (-\Delta)^{[s]-s} u\right).$$
See also \cite[Proposition 2.1]{BMP0}, \cite[Remark 3.2]{AJS2}, \cite{CDQ} and \cite[Theorems 1.2 and 1.8 and Corollary 1]{AJS1}. 
\end{Remark}

\begin{Remark}\label{rem_sobol_altern}
%Vedi anche [Notes on Sobolev Spaces - A. Visintin, 2017-2018]
By exploiting the Gagliardo seminorm one can define a fractional Sobolev space, for $p \in [1,\infty)$ and $s \in (0,1)$, by
$$\widetilde{W}^{s,p}(\R^N):= \left\{u \in L^p(\R^N) \mid \int_{\R^N} \int_{\R^N} \frac{|u(x)-u(y)|^p}{|x-y|^{N+ps}}dy < +\infty \right \};$$
this is a possible good choice \cite%[Section 1.2, Remark 3.5]
{DnPV},
but generally it does not coincide with $W^{s,p}(\R^N)$ for $p\neq 2$ \cite[Remark 3.5]{DnPV}. See also \cite[Remark 4 in Section 2.3.3]{Tri0} and \cite[Remark 6 in Section 2.1.1]{RS0}. %Section 2.1.2, Remark 2.1.3, Theorem 2.4.2
%$$\widetilde{W}^{s,p}(\R^N) = W^{s,p}(\R^N) \iff p=2.$$
The space $W^{s,p}(\R^N)$ is also known as \emph{Triebel-Lizorkin} space, or \emph{Bessel-potential space}, o \emph{Liouville space}, while $\widetilde{W}^{s,p}(\R^N)$ is also known as \emph{Besov space} or \emph{Slobodecki\u{\i} space}. %Slobodeckij space
\end{Remark}

%%%%%%%%%%%%%%%%%%%%%%%%
\subsubsection{Radially symmetric functions}

In order to gain some compactness on the entire space, we consider also the subspace of radially symmetric functions 
%$$H^s_r(\R^N) = \big\{ u \in H^s(\R^N) \mid u(x)= u(|x|) \big\};$$
$$H^s_r(\R^N) = \big\{ u \in H^s(\R^N) \mid \exists \, v: \R_+ \to \R \; \hbox{ s.t. } \; u(x)= v(|x|) \big\};$$
to avoid cumbersom notation, we will alway write $u(x) \equiv u(|x|)$. 
We notice that the fractional Laplacian inherits the radial symmetry of the function (this is immediate by use of the Fourier transform \eqref{eq_fourier_def}, see also \cite[Lemma 2.7]{Gar0}); % \tr{qui vale per $C^2 \cap L^{\infty}$; in general?}); 
anyway, it has not an easy representation in radial coordinates (see \cite{FV0} and \cite[Lemma 7.1]{Gar0}) based on Gaussian hypergeometric functions (see Section \ref{sec_hypergeo}), and this obstructs, for example, ODE's methods for resolution of PDEs.

We recall that, whenever $N\geq 2$, Lions proved the compact embedding \cite{Lio2} (see also \cite[Proposition 1.7.1]{Caz0} and \cite{EbSc0}) 
\begin{equation}\label{eq_immer_compatt}
H^s_r(\R^N) \hookrightarrow \hookrightarrow L^p(\R^N)
\end{equation}
for every $p \in (2, 2^*_s)$; however, as shown in \cite{CO0} for general $s\in (0, \frac{1}{2})$, a result in the spirit of Radial Lemma by Strauss \cite{Stra0}
$$|u(x)|^2 \lesssim \frac{1}{|x|^{N-2s}} \norm{(-\Delta)^{s/2}u}_2^2, \quad x \in \R^N\setminus\{0\}$$
is not available in the fractional framework $H^s_r(\R^N)$. We highlight that the embedding is not compact for $q=2^*_s$ even on bounded subsets of $\R^N$.
Sometimes we will write $\norm{\cdot}_{H^s_r}:=\norm{\cdot}_{H^s}$.
%$$\|u\|^2_{H^s_r}= \int_{\R^N} \big(|(-\Delta)^{s/2} u|^2 + u^2\big) \, dx.$$

\begin{Remark}
We observe that
$$H^s_r(\R^N)=\textnormal{Fix}(\textnormal{O}(N))=\{ u \in H^s(\R^N) \mid \tau(Q, u)=u \; \hbox{ for each $Q \in \textnormal{O}(N)$}\},$$
where $\textnormal{O}(N)$ is the orthogonal group of rotation matrices and the isometric action is given by
%$$\tau: (Q, \mu, u) \in \textnormal{O}(N)\times (\R_+ \times H^s(\R^N)) \mapsto (\mu, u(Q \cdot)) \in \R_+ \times H^s(\R^N);$$
$$\tau: (Q, u) \in \textnormal{O}(N)\times H^s(\R^N) \mapsto u(Q \cdot) \in H^s(\R^N);$$
working with a variational formulation, we will often work with $\textnormal{O}(N)$-invariant functionals: by the Principle of Symmetric Criticality of Palais \cite{Pal0} we will obtain that every critical point on $H^s_r(\R^N)$ is actually a critical point on the whole $H^s(\R^N)$, which justifies our restriction onto the radial setting.
\end{Remark}

\begin{Remark}\label{rem_grad_fract12}
Notice that, when $s=1$, we have %\tr{(up to constants)}
%Risulta solo $\mc{F}(\partial_i u) = i |\xi|\widehat{u}$, che con il modulo sparisce; oppure, con una divera definizione della trasformata, comparirebbe $\mc{F}(\partial_i u) = 2 \pi i |\xi|\widehat{u}$.
\begin{align*}
\norm{(-\Delta)^{1/2}u}_2^2 &= \int_{\R^N} \pabs{ \mc{F}^{-1} (|\xi| \widehat{u})}^2 = \int_{\R^N}\pabs{|\xi|\widehat{u}}^2 = \sum_i \int_{\R^N} \pabs{|\xi_i| \widehat{u}}^2 \\
&= \sum_i \int_{\R^N} \pabs{ \mc{F}(\partial_i u) }^2 = \sum_i \int_{\R^N} \pabs{\partial_i u}^2 = \int_{\R^N} |\nabla u|^2 \\
&= \norm{\nabla u}_2^2
\end{align*}
and this justifies, for example, the use of $(-\Delta)^{s/2}$ in the weak formulation of PDEs (see Definition \ref{def_weak_solut}). We highlight, anyway, the nontriviality of the relation, since $(-\Delta)^{1/2}$ is a nonlocal operator, while $\nabla$ is a local operator (see also \cite[Section 6]{Gar0}).

When $s=1$ thus we will actually consider the classical Sobolev space $H^1(\R^N)$ endowed with
$$\norm{u}_{H^1} := \left(\intRN \big(\abs{\nabla u}^2+u^2\big)\, dx\right)^{1/2} \quad \hbox{ for $u \in H^1(\R^N)$}$$
and its subspace
$$ H^1_r(\R^N):= \{ u \in H^1(\R^N) \mid u \hbox{ radially symmetric}\}.$$
%
%Moreover, by the fact that $|\xi|^2 \leq 1+ |\xi|^{2t}$ for every $t>1$, we see that
%$$\norm{u}_{H^{2s}(\R^N)}^2 = \int_{\R^N} |1 + |\xi||^{2}| |\widehat{u}|^2 \leq \int_{\R^N} |2 + |\xi|^{4s}||^2 |\widehat{u}|^2 \leq 4 \norm{u}_{H^1(\R^N)}^2$$
%whenever $2s\geq 1$, and thus
%$$H^{2s}(\R^N) \hookrightarrow H^1(\R^N) \quad \hbox{if $2s\geq 1$}.$$
%
\end{Remark}

%%%%%%%%%%%%%%%
\subsubsection{Tail-controlling mixed norms}
\label{sec_mixed_norm}

In order to handle the long range interaction of the fractional norms, we will make use of the following \emph{mixed} Gagliardo seminorm 
$$[u]_{A_1, A_2}^2:= \int_{A_1} \int_{A_2} \frac{|u(x)-u(y)|^2}{|x-y|^{N+2s}} dx \, dy, \quad [u]_A:=[u]_{A,A}$$
for any $A_1, A_2,A \subset \R^N$ and $u \in H^s(\R^N)$;
by using that $\varphi_u(x,y):= \frac{|u(x)-u(y)|}{|x-y|^{N/2+s}}$ satisfies $\varphi_{u+v} \leq \varphi_u + \varphi_v$ and $[u]_{A_1, A_2} = \norm{\varphi_u}_{L^2(A_1 \times A_2)}$, we have that $[u]_{A_1, A_2}$ is actually a seminorm. 
This seminorm has been introduced in \cite{CG0}, although after the publication the authors discovered that similar tools were implemented in different frameworks \cite{FKV, JW0, CDMV}.
% were anyway other authors made use of similar tools, 

For any $u \in H^s(\R^N)$ and $A \subset \R^N$ it will be useful to work also with the following norms:
$$%\begin{equation*} %\label{eq_def_Htilde}
\norm{u}_{A}^2:= \norm{u}_{L^2(A)}^2 + [u]_{A,\R^N}^2
$$%\end{equation*}
and
$$%\begin{equation*}%\label{eq_def_triplanorm}
\tnorm{u}_{A}:=\norm{u}_{L^{p+1}(A)} + \norm{u}_{A},
$$%\end{equation*}
for some suitable $p \in (2, 2^*_s)$.
%where $p$ is introduced in assumption (f1.3). 
We highlight that $\norm{u}_{\R^N}=\norm{u}_{H^s(\R^N)}$, but generally $\norm{u}_A \geq \norm{u}_{H^s(A)}$ for $A \neq \R^N$. 
By $H^s(A) \hookrightarrow L^{p+1}(A)$ the norms $\norm{\cdot}_A$ and $\tnorm{\cdot}_A$ are equivalent: on the other hand, the constant such that $\tnorm{u}_A \leq C_A \norm{u}_A$ depends on $A$, thus not useful for $\eps$-dependent sets $A=A(\eps)$ (see Chapter \ref{chap_concentr}). This is why we will make direct use also of $\tnorm{\cdot}_A$. 

\smallskip

Regarding $\eps$-dependent norms, we will use also
$$\norm{u}_{H_{V, \eps}^s(\R^N)}^2:= \norm{(-\Delta)^{s/2} u}_2^2 + \int_{\R^N} V(\eps x) u^2 dx$$
which is an equivalent norm on $H^s(\R^N)$ whenever $V \in L^{\infty}(\R^N)$ with $V \geq V_0>0$; 
%, thanks to the positivity and the boundedness of $V$ and \eqref{eq_equiv_norm_laplac}; 
the space $H^s_{\eps}(\R^N)$ is defined straightforwardly.

%%%%%%%%%%%%%%%%%%%%%%%%%%%%%%%%%%%%%%%%%%%%%%%%%%%%%%%%%%%%%%%%%%%%%%

\subsection{Some computations: hypergeometric Gaussian functions}
\label{sec_hypergeo}

In order to implement some comparison argument (see Section \ref{sec_frac_aux}), we search for a function which behaves like $\sim\frac{1}{|x|^{\beta}}$, $\beta >0$, and which lies in $H^s(\R^N)$; in order to handle the presence of a pole in the origin when $\beta \geq N$, we 
%it would be useful to work with the fractional Laplacian of the function $\frac{1}{|x|^{\beta}}$, $\beta >0$. 
%On the other hand, $(-\Delta)^s \frac{1}{|\,\cdot\,|^{\beta}}$ is well defined only for $\beta \in (0,N)$, because of the presence of a pole in the origin. 
%Thus we 
make the following choice, by considering, for any $\beta>0$, 
$$h_{\beta}(x):= \frac{1}{(1+|x|^2)^{\frac{\beta}{2}}};$$
notice that, when $\beta=N+2s$, this function is related to the extremals of the fractional Sobolev inequality \cite{Lie2,CFW} and to the solutions of the zero mass critical fractional Choquard equation \cite{Le0} (see also Proposition \ref{prop_HLS} below). 
Chosen $h_{\beta}$ in this way, we have \cite[Table 1 page 168]{Kwa1} (see also \cite[Sections 4 and 6]{FV0})
\begin{equation}\label{eq_frac_hyper}
(-\Delta)^s h_{\beta} (x) = C_{\beta, N, s} \, {}_2F_1\left(\frac{N}{2}+s, \frac{\beta}{2} + s, \frac{N}{2}; -|x|^2\right)
\end{equation}
where 
$$C_{\beta, N, s}:= 2^{2s} \frac{\Gamma\big(\frac{N}{2}+s\big) \Gamma\big(\frac{\beta}{2} + s\big)}{\Gamma\big(\frac{N}{2}\big) \Gamma\big(\frac{\beta}{2}\big)}>0 $$ 
and ${}_2F_1$ 
denotes the Gauss hypergeometric function (see also \cite[Corollary 2]{DKK}, observed that $h_{\beta}(x)= {}_2F_1(\frac{N}{2}, \frac{\beta}{2}, \frac{N}{2}, -|x|^2)$). 
Notice that we will be interested in 
$$\beta \in (0, N+2s].$$

%VEDI ANCHE GAROFALO PAG 26

%Recall that, by \eqref{eq_frac_hyper}, 
%\begin{equation}\label{eq_frac_hyper_2}
%(-\Delta)^s h_{\beta} (x) = C_{\beta, N, s} \, {}_2F_1\left(\frac{N}{2}+s, \frac{\beta}{2} + s, \frac{N}{2}; -|x|^2\right)
%\end{equation}
%where $h_{\beta}(x)= \frac{1}{(1+|x|^2)^{\frac{\beta}{2}}}$, and $C_{\beta, N, s}= 2^{2s} \frac{\Gamma\big(\frac{N}{2}+s\big) \Gamma\big(\frac{\beta}{2} + s\big)}{\Gamma\big(\frac{N}{2}\big) \Gamma\big(\frac{\beta}{2}\big)}>0 $.

The asymptotic behaviour at infinity of the hypergeometric function appearing in \eqref{eq_frac_hyper} can be found in \cite[pages 559-560]{AbSt0} (see also \cite[pages 78-79, 88]{AAR} and \cite[page 161]{WG0}).
Recall that the Gamma function $\Gamma(z)$ is well defined whenever $z \in \R \setminus (-\N)$ and $|\Gamma(z)|\to +\infty$ as $z$ approaches $-\N$ (so that the \emph{reciprocal Gamma function} is well defined on $-\N$ and equals zero). 
Moreover, recall the symmetry property ${}_2F_1(a,b,c;x)={}_2F_1(b,a,c;x)$ and the fact that ${}_2F_1(0,b,c;x) =1$ and ${}_2F_1(-1,b,c;x) = 1 - \frac{b}{c} z$.
\begin{Lemma}[\cite{AbSt0}]
Consider ${}_2F_1(a,b,c; \cdot)$. For the sake of simplicity, assume a priori that $a,b,c>0$ and
$$a-c \in \R_+ \setminus \N,$$
$$a-b \in \Z \iff a-b \in \N,$$ 
$$b-c \in \N \iff b-c \in \{0,1\};$$
in particular $a-b$ and $b-c$ do not lie in $\Z$ at the same time. 
We have the following asymptotic estimates as $x\to -\infty$.
\begin{itemize}
\item If $a-b \notin \Z$ and $b-c \notin \N$, then
$${}_2F_1(a,b,c; x) \stackrel{.}\sim \frac{\Gamma(c) \Gamma(b-a)}{\Gamma(c-a)\Gamma(b)} \frac{1}{(-x)^{a}} + \frac{\Gamma(c) \Gamma(a-b)}{\Gamma(c-b)\Gamma(a)} \frac{1}{(-x)^{b}} ;$$
\item If $b=c$ (and $a-b \notin \Z$), then
$${}_2F_1(a,b,b;x) =\frac{1}{(1-x)^{a}}; $$ 
\item If $b=c+1$ (and $a-b \notin \Z$), then
%$${}_2F_1(a,b,b-1;x) = - \frac{\Gamma(b-1)\Gamma(b-a)}{\Gamma(b-a-1)\Gamma(b)}\frac{x}{(1-x)^{a+1}} + \frac{1}{(1-x)^{a+1}} \stackrel{.}\sim \frac{\Gamma(b-1)\Gamma(b-a)}{\Gamma(b-a-1)\Gamma(b)} 
% \frac{1}{(-x)^{a}};$$
\begin{align*}
{}_2F_1(a,b,b-1;x) &= - \frac{\Gamma(b-1)\Gamma(b-a)}{\Gamma(b-a-1)\Gamma(b)}\frac{x}{(1-x)^{a+1}} + \frac{1}{(1-x)^{a+1}}\\
& \stackrel{.}\sim \frac{\Gamma(b-1)\Gamma(b-a)}{\Gamma(b-a-1)\Gamma(b)} 
 \frac{1}{(-x)^{a}};
\end{align*}
\item If $a=b$ (and $b-c \notin \N$), then
$${}_2F_1(a,a,c;x) \stackrel{.}\sim \frac{\Gamma(c)}{\Gamma(a)\Gamma(c-a)} \frac{\log(-x)}{(-x)^a} + \frac{C_1}{(-x)^a} \stackrel{.}\sim \frac{\Gamma(c)}{\Gamma(a)\Gamma(c-a)} \frac{\log(-x)}{(-x)^a}; $$
\item If $a-b \in \N^*$ (and $b-c \notin \N$), then
$${}_2F_1(a,b,c;x) \stackrel{.}\sim \frac{\Gamma(c) \Gamma(a-b)}{\Gamma(c-b)\Gamma(a)} \frac{1}{(-x)^{b}} + C_2 \frac{\log(-x)}{(-x)^a} + \frac{C_3}{(-x)^a} \stackrel{.}\sim \frac{\Gamma(c) \Gamma(a-b)}{\Gamma(c-b)\Gamma(a)} \frac{1}{(-x)^{b}}. $$
\end{itemize}
Here $C_i$, $i=1,2,3$, are some strictly positive constants.
\end{Lemma}

Notice that $a= \frac{N}{2} + s$, $b=\frac{\beta}{2} + s$, $c= \frac{N}{2}$ satisfy the assumptions of the previous Lemma, whenever $s \in (0,1)$ and $\beta \in (0, N+2s]$. 
Thus, exploiting the representation of $(-\Delta)^s h_{\beta}$ given in \eqref{eq_frac_hyper} and the results on Gauss hypergeometric functions, we come up with the following lemma.

\begin{Lemma}\label{lem_calcol_2potenz}
Let $\beta \in (0,N+2s]$. 
Then $(-\Delta)^s h_{\beta}(x)$ is well defined for every $x \neq 0$. Moreover, we have the following asymptotic behaviours:
\begin{itemize}
\item if $\beta\in (N, N+2s]$, then
$$(-\Delta)^s h_{\beta}(x) \stackrel{.}\sim C'_{\beta, N, s} \frac{1}{|x|^{N+2s}} \quad \hbox{ as $|x|\to +\infty$}$$
where $C'_{\beta,N,s}:=2^{2s} \frac{\Gamma\big(\frac{N}{2}+s\big) \Gamma\big(\frac{\beta}{2}-\frac{N}{2}\big)}{\Gamma\big(\frac{\beta}{2}\big) \Gamma\big(-s\big)}<0.$ 
This in particular includes the case $\beta = N-2s+2$ (possible if $s > \frac{1}{2}$), with $C_{N-2s+2,N,s}'= - 2^{2s+1}\frac{s }{N-2s}<0$. 
Notice moreover that $C'_{N+2s,N,s}=2^{2s} \frac{\Gamma(s)}{\Gamma(-s)} \to 0$ as $s \to 1^{-}$.

\item if $\beta=N$, then 
$$(-\Delta)^s h_N(x) \stackrel{.}\sim C'_{N,N,s} \frac{\log(|x|)}{|x|^{N+2s}} \quad \hbox{ as $|x|\to +\infty$}$$
where 
$C'_{N,N,s}:=
2^{2s+1} \frac{\Gamma\big(\frac{N}{2}+s\big)}{\Gamma\big(\frac{N}{2}\big) \Gamma\big(-
s\big)} <0.$

\item if $\beta\in (N-2s, N)$, then
$$(-\Delta)^s h_{\beta}(x) \stackrel{.}\sim C_{\beta,N,s}' \frac{1}{|x|^{\beta + 2s}} \quad \hbox{as $|x|\to +\infty$}$$
where 
$C_{\beta,N,s}':= 2^{2s} \frac{ \Gamma\big(\frac{\beta}{2} +s\big)\Gamma\big(\frac{N}{2}-\frac{\beta}{2}\big)}{\Gamma\big(\frac{\beta}{2}\big)\Gamma\big(\frac{N}{2}-\frac{\beta}{2}-s\big) }<0.$

\item if $\beta=N-2s$, then 
\begin{align*}
(-\Delta)^s h_{N-2s}(x) &= C_{N-2s,N,s}' h_{N+2s}%\frac{N}{2}+s}
(x) \quad \hbox{for $x \in \R^N\setminus \{0\}$}\\
& \stackrel{.}\sim C_{N-2s,N,s}' \frac{1}{|x|^{N + 2s}} \quad \hbox{as $|x|\to +\infty$}
\end{align*} 
where 
$C_{N-2s,N,s}':= 2^{2s} \frac{ \Gamma\big(\frac{N}{2} +s\big)}{\Gamma\big(\frac{N}{2}-s\big) } >0.$

\item if $\beta\in (0, N-2s)$, then 
$$(-\Delta)^s h_\beta(x) \stackrel{.}\sim C_{\beta,N,s}' \frac{1}{|x|^{\beta + 2s}} \quad \hbox{as $|x|\to +\infty$}$$
where 
$C_{\beta,N,s}':= 2^{2s} \frac{ \Gamma\big(\frac{\beta}{2} +s\big)\Gamma\big(\frac{N}{2}-\frac{\beta}{2}\big)}{\Gamma\big(\frac{\beta}{2}\big)\Gamma\big(\frac{N}{2}-\frac{\beta}{2}-s\big) }>0.$ 
This in particular includes the case $\beta = N - 2k$ with $k = 1, \dots, [\frac{N}{2}]$.
\end{itemize}
\end{Lemma}

\begin{Remark}
Notice that, for $\beta \in \{N-2s\} \cup [N, N+2s]$, % +\infty)$, 
the asymptotic behaviour of $\abs{(-\Delta)^s h_\beta(x)}$ does not depend on $\beta$; on the other hand, the sign and the precise constant depend on $\beta$.

In the case $\beta \in (0, N)\setminus \{N-2s\}$, we may use $x\mapsto \frac{1}{|x|^{\beta}}$, whose fractional Laplacian has a close (simple) representation:
$$\left((-\Delta)^s \frac{1}{|\cdot|^\beta} \right)(x) = C_{\beta,N,s} \frac{1}{|x|^{\beta+2s}},$$
see \cite[Table 1 and Theorem 3.1]{Kwa1} (see also \cite[Lemma 4.1]{Fal1}, \cite[Appendix 1, page 798]{Vaz0} and \cite[Lemma A.2]{BrMoSq}). 
In particular
$$(-\Delta)^s h_\beta(x) \stackrel{.}\sim \left((-\Delta)^s \frac{1}{|\cdot|^\beta} \right)(x) \quad \hbox{as $|x|\to +\infty$}.$$
On the other hand, if $\beta=N-2s$, we obtain, far from the origin, $(-\Delta)^s \frac{1}{|\,\cdot\,|^\beta} \equiv 0$ (recall that the Riesz potential $\frac{1}{|\cdot|^{N-2s}} \equiv I_{2s}$ is a fundamental solution, see Proposition \ref{prop_inversa_dx}); thus, in particular, the two functions have different asymptotic behaviours. 
This is the same reason why, for $h_{\beta}$, we have a discontinuity on the behaviour at infinity around $\beta=N-2s$.

Finally we highlight that, when $\beta=N+2s$, we may use the function found in Lemma \ref{lem_esist_sol_part}.
\end{Remark}

%%%%%%%%%%%%%%%%%%%%%%%%%%%%%%%%%%%%%%%%%%%%%%%%%%%%%%%%%%%%%%%%%%%%%%

\subsection{Definitions of solutions: weak, viscosity, strong, classical}

In the majority of the thesis we will work with the notion of weak solutions, by exploiting a variational formulation. Anyway, sometimes we will need to exploit different formulations, in particular strong, classical and viscosity formulations; that is why we recall them here for the sake of clarity.

%\tr{In what follows, we will not focus on the assumptions on $g$, assuming that they are minimal in order to give sense to the considered quantities.}

%\tor{
%\begin{Definition}[Strong and classical solution]
%Let $\Omega \subseteq \R^N$ and $g: \Omega \times \R \to \R$ be \tr{measurable}. We say that $u$ is a \emph{strong} solution to
%$$(-\Delta)^s u = g(x,u) \quad \hbox{in $\Omega$}$$
%if $u$ and $(-\Delta)^s u$ are almost everywhere defined (e.g., and $u \in H^{2s}(\Omega)$) and $u$ satisfies the relation for almost every $x \in \Omega$. 
%%
%\\
%We say instead that $u$ is a \emph{classical} solution if $u$ and $(-\Delta)^s u$ are continuous (e.g., $u \in L^p(\R^N) \cap C^{\gamma}_{loc}(\R^N)$ for some $p \in [1,+\infty]$ and $\gamma >2s$) and the relation is satisfied pointwise everywhere on $\Omega$.
%\end{Definition}
%}

\begin{Definition}[Strong and classical solution]
Let $\Omega \subseteq \R^N$ and $g: \Omega %\times \R 
\to \R$. % be \tr{defined}. 
We say that $u$ is a \emph{strong} solution to
$$(-\Delta)^s u = g(x) \quad \hbox{in $\Omega$}$$
if $u$ and $(-\Delta)^s u$ are almost everywhere defined (e.g. $u \in H^{2s}(\Omega)$) and $u$ satisfies the relation for almost every $x \in \Omega$. 

We say instead that $u$ is a \emph{classical} solution if $u$ and $(-\Delta)^s u$ are continuous (e.g. $u \in L^p(\R^N) \cap C^{\gamma}_{loc}(\R^N)$ for some $p \in [1,+\infty]$ and $\gamma >2s$) and the relation is satisfied pointwise everywhere on $\Omega$.
\end{Definition}

%\begin{Definition}[Weak solution]\label{def_weak_solut}
%Let $\Omega \subseteq \R^N$ and $g: \Omega \times \R \to \R$ be \tr{measurable}. We say that $u\in H^s(\Omega)$ is a \emph{weak} subsolution [supersolution] of 
%$$(-\Delta)^s u = g(x,u) \quad \hbox{ in $\Omega$}$$
%if
%\begin{align}\label{eq_weak_sol}
%&\int_{\R^N} (-\Delta)^{s/2} u (-\Delta)^{s/2} \varphi dx + \mu \int_{\R^N} u \varphi dx \leq \int_{\R^N} g(x,u) \varphi dx
%\\ &\qquad \qquad
%\Big[\int_{\R^N} (-\Delta)^{s/2} u (-\Delta)^{s/2} \varphi dx + \mu \int_{\R^N} u \varphi dx \geq \, \int_{\R^N} g(x,u) \varphi dx \Big] \notag
%\end{align}
%is well defined and holds for each positive $\varphi \in X^s_0(\Omega)$. 
%We say that $u$ is a weak solution if it is both a subsolution and a supersolution, i.e. if it satisfies the equality in \eqref{eq_weak_sol} for every $\varphi \in X^s_0(\Omega)$. 
%Notice that, when $\Omega=\R^N$, we have $X^s_0(\R^N) \equiv H^s(\R^N)$.
%\end{Definition}

\begin{Definition}[Weak solution]\label{def_weak_solut}
Let $\Omega \subseteq \R^N$ and $g: \Omega %\times \R 
\to \R$ be measurable. We say that $u\in H^s(\Omega)$ is a \emph{weak} subsolution [supersolution] of 
$$(-\Delta)^s u = g(x) \quad \hbox{ in $\Omega$}$$
if
\begin{align}\label{eq_weak_sol}
&\int_{\R^N} (-\Delta)^{s/2} u (-\Delta)^{s/2} \varphi dx %+ \mu \int_{\R^N} u \varphi dx 
\leq \int_{\R^N} g(x) \varphi dx
\\ &\qquad \qquad
\Big[\int_{\R^N} (-\Delta)^{s/2} u (-\Delta)^{s/2} \varphi dx %+ \mu \int_{\R^N} u \varphi dx 
\geq \, \int_{\R^N} g(x) \varphi dx \Big] \notag
\end{align}
is well defined (finite) and holds for each positive $\varphi \in X^s_0(\Omega)$. 
We say that $u$ is a weak solution if it is both a subsolution and a supersolution, i.e. if it satisfies the equality in \eqref{eq_weak_sol} for every $\varphi \in X^s_0(\Omega)$. 
Notice that, when $\Omega=\R^N$, we have $X^s_0(\R^N) \equiv H^s(\R^N)$.
\end{Definition}

\begin{Remark}
By \eqref{eq_sobolev_polaraz} and \eqref{eq_fourier_def} we may interpret the left-hand side of \eqref{eq_weak_sol} as 
%$$ \tfrac{1}{2} C_{N,s} \int_{\R^N} \int_{\R^N} \frac{\big(u(x)-u(y)\big)\big(\varphi(x)-\varphi(u)\big)}{|x-y|^{N+2s}} dx dy = \int_{\R^N} (-\Delta)^{s/2} u (-\Delta)^{s/2} \varphi dx =\int_{\R^N} |\xi|^{2s} \widehat{u} \widehat{\varphi} d\xi.$$
\begin{align*}
 \int_{\R^N} (-\Delta)^{s/2} u (-\Delta)^{s/2} \varphi dx &\equiv \tfrac{1}{2} C_{N,s} \int_{\R^N} \int_{\R^N} \frac{\big(u(x)-u(y)\big)\big(\varphi(x)-\varphi(y)\big)}{|x-y|^{N+2s}} dx dy \\
&\equiv \int_{\R^N} |\xi|^{2s} \widehat{u} \widehat{\varphi} d\xi.
\end{align*}
Moreover we see that the definition of weak solution is justified by % actually makes sense, since 
the following \emph{integration by parts} rule %holds
$$\int_{\R^N} (-\Delta)^s u v = \int_{\R^N} |\xi|^{2s} \widehat{u} \widehat{v} = \int_{\R^N} |\xi|^{s} \widehat{u} |\xi|^s \widehat{v} = \int_{\R^N} (-\Delta)^{s/2} u (-\Delta)^{s/2} v$$
which holds whenever $u \in H^{2s}(\R^N)$ and $v \in H^s(\R^N)$. In particular, if both $u,v \in H^{2s}(\R^N)$ we have (see also \cite[Lemma 5.4]{Gar0})
$$\int_{\R^N} (-\Delta)^s u v = \int_{\R^N} u (-\Delta)^s v.$$
\end{Remark}

\begin{Remark}
If $\varphi \in C^{\infty}_c(\R^N)$ and $u \in L^{\infty}(\R^N)$, we can show the relation
$$\int_{\R^N} u (-\Delta)^s \varphi = \tfrac{1}{2} C_{N,s} \int_{\R^N} \int_{\R^N} \frac{\big(u(x)-u(y)\big)\big(\varphi(x)-\varphi(y)\big)}{|x-y|^{N+2s}} dx dy$$
also by exploiting the pointwise definition of the fractional Laplacian. Indeed (assume for simplicity $s \in (0,\frac{1}{2})$ to avoid the technicality of the Principal Value) we have
\begin{equation}\label{eq_dim_int_ppar}
\int_{\R^N} u(x) (-\Delta)^s \varphi(x) dx = C_{N,s} \int_{\R^N} \left( \int_{\R^N} \frac{u(x)\big(\varphi(x)-\varphi(y)\big)}{|x-y|^{N+2s}} dy \right) dx.
\end{equation}
First, we rewrite \eqref{eq_dim_int_ppar} by applying Fubini-Tonelli theorem, possible because
$$ \int_{\R^{2N}} \pabs{ \frac{u(x)\big(\varphi(x)-\varphi(y)\big)}{|x-y|^{N+2s}}} \leq \norm{u}_{\infty} \int_{\R^{2N}} \frac{|\varphi(x)-\varphi(y)|}{|x-y|^{N+2s}} < \infty;$$
notice that we are actually using that $\varphi \in \widetilde{W}^{2s,1}(\R^N)$ (see Remark \ref{rem_sobol_altern}). 
Thus
$$\int_{\R^N} u(x) (-\Delta)^s \varphi(x) dx = C_{N,s} \int_{\R^N} \left( \int_{\R^N} \frac{u(x)\big(\varphi(x)-\varphi(y)\big)}{|x-y|^{N+2s}} dx \right) dy.$$
Secondly, we rewrite \eqref{eq_dim_int_ppar} by simply renaming the variables, that is
%\begin{align*}
%\int_{\R^N} u(x) (-\Delta)^s \varphi(x) dx &= \int_{\R^N} \left( \int_{\R^N} \frac{u(y)\big(\varphi(y)-\varphi(x)\big)}{|y-x|^{N+2s}} dx \right) dy \\
%&= \int_{\R^N} \left( \int_{\R^N} \frac{\big(-u(y)\big)\big(\varphi(x)-\varphi(y)\big)}{|x-y|^{N+2s}} dx \right) dy.
%\end{align*}
$$\int_{\R^N} u(x) (-\Delta)^s \varphi(x) dx = C_{N,s} \int_{\R^N} \left( \int_{\R^N} \frac{\big(-u(y)\big)\big(\varphi(x)-\varphi(y)\big)}{|x-y|^{N+2s}} dx \right) dy.$$
By summing the two expressions obtained, we get the claim.
%$$2\int_{\R^N} u(x) (-\Delta)^s \varphi(x) dx = C_{N ,s} \int_{\R^N} \left( \int_{\R^N} \frac{\big(u(x)-u(y)\big)\big(\varphi(x)-\varphi(y)\big)}{|x-y|^{N+2s}} dx \right) dy$$
%which is the claim.
\end{Remark}

For the following definition, see e.g. \cite[page 136]{SeV} or \cite[Definition 2.1]{CFQ}. 
%\begin{Definition}[Viscosity solution] 
%Let $\Omega \subseteq \R^N$ and $g: \Omega \times \R \to \R$ be \tr{some function}. We say that $u\in C(\R^N)$ is a \emph{viscosity} subsolution [supersolution] of
%$$(-\Delta)^s u = g(x,u) \quad \hbox{ in $\Omega$}$$
% if, for any $x_0 \in \Omega$, every $U\subset \Omega$ open neighborhood of $x_0$, and every $\phi \in C^2(U)$ such that
%$$\phi(x_0)=u(x_0), \quad \phi \geq u \; [\phi \leq u] \; \hbox{ in $U$}$$
%set
%$$v:= \phi \chi_U + u \chi_{U^c}$$ 
%we have
%\begin{align}\label{eq_visc_sol}
%&(-\Delta)^s v(x_0) \leq g(x_0, u(x_0)) \\ &\qquad 
% \big[(-\Delta)^s v(x_0) \geq g(x_0, u(x_0))\big]. \notag
%\end{align}
%We say that $u$ is a viscosity solution if it is both a viscosity subsolution and a viscosity supersolution.
%\end{Definition}

\begin{Definition}[Viscosity solution] 
Let $\Omega \subseteq \R^N$ and $g: \Omega %\times \R 
\to \R$. % be \tr{some function}. 
We say that $u\in C(\R^N)$ is a \emph{viscosity} subsolution [supersolution] of
$$(-\Delta)^s u = g(x) \quad \hbox{ in $\Omega$}$$
 if, for any $x_0 \in \Omega$, every $U\subset \Omega$ open neighborhood of $x_0$, and every $\phi \in C^2(U)$ such that
$$\phi(x_0)=u(x_0), \quad \phi \geq u \; [\phi \leq u] \; \hbox{ in $U$}$$
set
$$v:= \phi \chi_U + u \chi_{U^c}$$ 
we have
%\begin{align}\label{eq_visc_sol}
%&(-\Delta)^s v(x_0) \leq g(x_0) \\ &\qquad 
% \big[(-\Delta)^s v(x_0) \geq g(x_0)\big]. \notag
%\end{align}
\begin{equation*}\label{eq_visc_sol}
(-\Delta)^s v(x_0) \leq g(x_0) \quad \big[(-\Delta)^s v(x_0) \geq g(x_0)\big].
\end{equation*}

We say that $u$ is a viscosity solution if it is both a viscosity subsolution and a viscosity supersolution.
\end{Definition}

We observe that, generally, the function $v$ appearing in the definition of viscosity solution might be discontinuous. 
More generally, this definition involves lower and upper semicontinuity of $u$ (see for instance \cite[Definition 2.2]{CafSil2}). 
Furthermore, one can easily check that every %(continuous) 
classical solution is a viscosity solution, that the sum of two viscosity solutions is still a viscosity solution (with source the sum of the sources), and that the notion of viscosity solution is conserved on subdomains $\Omega' \subset \Omega$.

We refer to \cite[Remark 2.11]{ROS1} and \cite[Theorem 1]{SeV} for some discussions on the relation between classical, weak and viscosity solutions on bounded domains.

\smallskip

When dealing with equations with nonlinearities of the type $h=h(x,u)$, $h: \Omega \times \R \to \R$, we interpret the equation by saying that $u$ is a classic/strong/weak/viscosity solution if $u$ satisfies the equation with nonlinearity $g(x):=h(x,u(x))$. The same interpretation will be given in the case of nonlocal nonlinearities (see Section \ref{sec_riesz_poten}).

\begin{Remark}
In this preliminary Chapter, we will sometimes mention \emph{distributional} solutions of equations of the type $(-\Delta)^s u = T$, with $T$ distribution on some $\Omega$. By this, we mean that $u\in L^1_{loc}(\R^N)$ satisfies \eqref{eq_spazio_gener_s} and
$$\int_{\Omega} u (-\Delta)^{s/2} \varphi = T(\varphi)$$
for every $\varphi \in C^{\infty}_c(\Omega)$. The extra condition required on $u$ (differently form the usual definition of distributional solution) is due to the fact that $(-\Delta)^{s/2} \varphi$ has generally not compact support; here we use thus \eqref{eq_stima_classica_s} to well define the integral.
\end{Remark}

%%%%%%%%%%%%%%%%%%%%%%%%%%%%%%%%%%%%%%%%%%%%%%%%%%%%%%%%%%%%%%%%%%%%%%%%%%

\subsection{A concave Chain rule}
\label{sec_chain_rule}

We already pointed out how the fractional Laplacian does not satisfy a proper Lebiniz formula. The same conclusion is actually true looking at chain rule formulas. A first result is given by the following lemma.

\begin{Lemma}
\label{lem_Lipschitz_Sobolev}
Let $\Omega \subseteq \R^N$. If $u\in H^s(\Omega)$ and $h: \R \to \R$ is a Lipschitz function with $h(0)=0$, then $h(u) \in H^s(\Omega)$. 
%Notice that a general chain rule is not available for the fractional Laplacian. See also Lemma \ref{lem_chain_cordoba}.
\end{Lemma}

\claim Proof.
The proof is straightforward. Indeed
$$\norm{h(u)}_{L^2(\Omega)}^2 = \int_{\Omega}|h(u)-h(0)|^2 \, dx \leq \int_{\Omega} \norm{h'}_{\infty}^2 |u-0|^2 \, dx = \norm{h'}_{\infty}^2 \norm{u}_{L^2(\Omega)}^2$$
and
$$%\norm{(-\Delta)^{s/2} h(u)}_2^2 
[h(u)]_{H^s(\Omega)}\leq C_{N,s} \int_{\Omega} \int_{\Omega} \frac{\norm{h'}_{\infty}^2 |u(x)-u(y)|^2}{|x-y|^{N+2s}} \, dx \, dy = \norm{h'}_{\infty}^2 [u]_{H^s(\Omega)}
%\norm{(-\Delta)^{s/2} u}_2^2
.
\QED
$$

\bigskip

% (see also Remark \ref{rem_Lipschitz_Sobolev}):
We look now to proper pointwise chain rules. What one can prove is that, whenever $\varphi$ is convex (and Lipschitz), then the following inequality holds (see \cite[Theorem 1.1]{CafSir}, \cite[Theorem 19.1]{Gar0})
$$(-\Delta)^s \varphi(u) \leq \varphi'(u) (-\Delta)^s u$$
in the weak sense. 
One may expect the inverse inequality when handling concave functions: and this is actually what we need in the study of the asymptotic behaviour of ground state in doubly nonlocal equations (see Section \ref{sec_est_bel_vis}).

On the other hand, %We need first some preliminary results, in order to deal with the fractional Laplacian of the concave power of a function: 
since we do not know if $u^{\theta} \notin H^s(\R^N)$ when $u \in H^s(\R^N)$ and $\theta \in (0,1)$, the weak formulation seems not to be the right choise; pointwise formulation seems not good as well, since $(-\Delta)^s u^{\theta}$ might be not well defined, even by assuming $u$ regular.
%Notice that knowing a priori that $u$ is continuous, radially symmetric and decreasing seems of no use. 
The idea is thus to take advantage of a viscosity formulation.

We prove hence the following inequality in the case of concave (not globally Lipschitz) function, in the framework of viscosity solutions. Notice that we do not require $u$ to be in $L^2(\R^N)$.

\begin{Lemma}[C\'ordoba-C\'ordoba chain rule inequality]\label{lem_chain_cordoba}
Let $\varphi: I \to \R$ be a concave function, $I\subseteq \R$ interval, such that $\varphi \in C^1(I)$. % \cap Lip(I)$. 
Let $u: \R^N \to I$. %Then
\begin{itemize}
%\item If $(-\Delta)^{s/2} u \in L^2(\R^N)$, then $(-\Delta)^{s/2} \varphi(u)\in L^2(\R^N)$ and
%$$\norm{(-\Delta)^{s/2}\varphi(u)}_2 \leq \norm{\varphi'}_{L^{\infty}(I)} \norm{(-\Delta)^{s/2} u}_2. $$ 
\item Let $\Omega \subset \R^N$, and assume $\varphi \in Lip(u(\Omega))$. Then
%$$\left(\int_{\Omega} \int_{\Omega} \frac{|\varphi(u)(x)-\varphi(u)(y)|^2}{|x-y|^{N+2s}} dx dy \right)^{1/2} \leq \norm{\varphi'}_{L^{\infty}(u(\Omega))} \left(\int_{\Omega} \int_{\Omega} \frac{|u(x)-u(y)|^2}{|x-y|^{N+2s}} dx dy \right)^{1/2},$$
$$[\varphi(u)]_{H^s(\Omega)} \leq \norm{\varphi'}_{L^{\infty}(u(\Omega))} [u]_{H^s(\Omega)}.$$
%thus if $u \in H^s(\Omega)$, then $\varphi(u) \in H^s(\Omega)$. 
In particular, if $\varphi \in Lip(I)$ and $(-\Delta)^{s/2} u \in L^2(\R^N)$, then $(-\Delta)^{s/2} \varphi(u)\in L^2(\R^N)$ and
$$\norm{(-\Delta)^{s/2}\varphi(u)}_2 \leq \norm{\varphi'}_{L^{\infty}(I)} \norm{(-\Delta)^{s/2} u}_2. $$
\item If $u$ is defined pointwise, then
$$(-\Delta)^s (\varphi(u))(x) \geq \varphi'(u(x)) (-\Delta)^s u(x)$$
for every $x \in \R^N$ such that $(-\Delta)^s (\varphi(u))(x)$ and $(-\Delta)^s u(x)$ are well defined.
\item Assume in addition $\varphi$ invertible, increasing, with $\varphi^{-1} \in C^2$ increasing. 
If $u$ is a continuous viscosity supersolution of
$$(-\Delta)^s u \geq g \quad \hbox{ in $\Omega$}$$
for some function $g$ and $\Omega \subseteq \R^N$, then $\varphi(u)$ is a viscosity supersolution of
$$(-\Delta)^s (\varphi(u)) \geq \varphi'(u) g \quad \hbox{ in $\Omega$}.$$
\end{itemize}
\end{Lemma}

\claim Proof.
The first claim is a direct consequence of the Lipschitz continuity
$$ \int_{\Omega}\int_{\Omega} \frac{|\varphi(u(x))-\varphi(u(y))|^2}{|x-y|^{N+2s}} dx dy \leq \norm{\varphi'}_{L^{\infty}(u(\Omega))}^2 \int_{\Omega} \int_{\Omega} \frac{|u(x)-u(y)|^2}{|x-y|^{N+2s}} dx dy.$$
%First, by the Lipschitz continuity of $\varphi$ we obtain
%\begin{align*}
%\norm{(-\Delta)^{s/2}\varphi(u)}_2^2 &= C_{N,s} \int_{\R^N}\int_{\R^N} \frac{|\varphi(u(x))-\varphi(u(y))|^2}{|x-y|^{N+2s}} dx dy \\
%& \leq C_{N,s} \int_{\R^N} \int_{\R^N} \norm{\varphi'}_{L^{\infty}(I)}^2 \frac{|u(x)-u(y)|^2}{|x-y|^{N+2s}} dx dy \\
%&= \norm{\varphi'}_{L^{\infty}(I)}^2 \norm{(-\Delta)^{s/2} u}_2^2.
%\end{align*}
%Moreover, 
Secondly, by the concavity of $\varphi$, for each $t,r \in I$ we have
$$\varphi(t) - \varphi(r) \geq \varphi'(t) (t-r)$$
thus
\begin{align*}
(-\Delta)^s(\varphi(u))(x) &= C_{N,s} \int_{\R^N} \frac{\varphi(u(x))-\varphi(u(y))}{|x-y|^{N+2s}} dy \\
&\geq C_{N,s} \int_{\R^N} \frac{\varphi'(u(x)) \big(u(x)-u(y)\big)}{|x-y|^{N+2s}} dy = \varphi'(u(x)) (-\Delta)^s u(x).
\end{align*}

We move to the third part. Let $x_0 \in U \subset \Omega$ and $\phi \in C^2(U)$ be such that $\phi(x_0)= \varphi(u(x_0))$ and $\phi\leq \varphi(u)$ in $U$, and set $v:= \phi \chi_{U} + \varphi(u) \chi_{U^c}$. Let now
$$\psi:= \varphi^{-1} \circ \phi, \quad w:=\varphi^{-1} \circ v = \psi \chi_{U} + u \chi_{U^c}.$$
By the assumptions on $\varphi^{-1}$ we have $\psi \in C^2(U)$, $\psi(x_0)=u(x_0)$ and $\psi \leq u$ in $U$. Thus
$$(-\Delta)^s w(x_0) \geq g(x_0).$$
On the other hand, $w = \psi \in C^2$ on $U$ and $\varphi(w) = \phi \in C^2$ on $U$, hence both the functions are regular enough in a neighborhood of $x_0$ to state that both the fractional Laplacians are well defined (see Proposition \ref{prop_well_posed}). 
 Thus we may apply the previous point and obtain
$$(-\Delta)^s (\varphi(w))(x_0) \geq \varphi'(w(x_0)) (-\Delta)^s w(x_0).$$
Since $w(x_0)=u(x_0)$, $\varphi(w)=v$ and $\varphi'$ is positive, we obtain, by joining the two previous inequalities
$$(-\Delta)^s v (x_0) \geq \varphi'(u(x_0)) g(x_0)$$
which is the claim. 
This concludes the proof. 
\QED

\bigskip

%With a similar argument, 
As a corollary, we obtain the following result. 
%\\ \tr{scrivi i dettagli qui nella tesi?}%CCCOMMENT NOW, li scrivo se me li chiedono
\begin{Corollary}\label{corol_concav_u2}
Let $\theta \in (0,1)$, and let $u \in C(\R^N)$ be strictly positive. 
We have the following results.
\begin{itemize}
%\item If $(-\Delta)^{s/2} u \in L^2(\R^N)$, then
%$$\left(\int_{\Omega} \int_{\Omega} \frac{|u^{\theta}(x)-u^{\theta}(y)|^2}{|x-y|^{N+2s}} dx dy \right)^{1/2} \leq \frac{\theta}{\min_{\Omega} u^{1-\theta} } \norm{(-\Delta)^{s/2} u}_2$$
%for each $\Omega \subset \subset \R^N$. 
\item We have
%$$\left(\int_{\Omega} \int_{\Omega} \frac{|u^{\theta}(x)-u^{\theta}(y)|^2}{|x-y|^{N+2s}} dx dy \right)^{1/2} \leq \frac{\theta}{\min_{\Omega} u^{1-\theta} } \tr{\left(\int_{\Omega} \int_{\Omega} \frac{|u(x)-u(y)|^2}{|x-y|^{N+2s}} dx dy \right)^{1/2}}
%%\norm{(-\Delta)^{s/2} u}_2
%$$
$$[u^{\theta}]_{H^s(\Omega)} \leq \frac{\theta}{\min_{\Omega} u^{1-\theta} } [u]_{H^s(\Omega)}
%\norm{(-\Delta)^{s/2} u}_2
$$
In particular, if $u \in H^s_{loc}(\R^N)$, %L^2_{loc}(\R^N)$, 
then $ u^{\theta} \in H^s_{loc}(\R^N)$.%
\footnote{Indeed, if $u\in L^2_{loc}(\R^N)$, then $u^{\theta} \in L^2_{loc}(\R^N)$ can be deduced by the inverse H\"older inequality: $\int_{\Omega} u^2 = \int_{\Omega} u^2 \cdot 1 \geq \int_{\Omega} u^{\frac{2}{p}} \int_{\Omega} 1^{-\frac{1}{p-1}} = \int_{\Omega} u^{2 \theta} \cdot m(\Omega)$, if $p := \frac{1}{\theta}>1$ and $\Omega$ is bounded (with positive measure).}
As a consequence, if $u \in H^s(\R^N)$, then
$$[u^{\theta}]_{H^s(\Omega)} \leq \frac{\theta}{\min_{\Omega} u^{1-\theta} }\norm{(-\Delta)^{s/2} u}_2.
$$
\item If $(-\Delta)^s u$ is well defined pointwise, then
$$(-\Delta)^s u^{\theta}(x) \geq \frac{\theta}{(u(x))^{1-\theta}} (-\Delta)^s u(x)$$
for every $x \in \R^N$ such that $(-\Delta)^s u^{\theta}(x)$ is well defined. 
\item If $u$ is a viscosity supersolution of
$$(-\Delta)^s u \geq g \quad \hbox{ in $\Omega$}$$
for some function $g$ and $\Omega \subseteq \R^N$, then $u^{\theta}$ is a viscosity supersolution of
$$(-\Delta)^s u^{\theta} \geq \frac{\theta}{u^{1-\theta}} g \quad \hbox{ in $\Omega$}.$$
\end{itemize}
\end{Corollary}

%%%%%%%%%%%%%%%%%%%%%%%%%%%%%%%%%%%%%%%%%%%%%%%%%%%%%%%%%%%%%%%%%%%%%%

\subsection{Regularity: tail functions and De Giorgi classes}
\label{sec_regol_degiorgi}

%\subsubsection{Boundedness of signed solutions} 

%We observe that a direct proof of the boundedness for generally signed solutions, but assuming also (f5).
%, can be found in Appendix \ref{sec_bound_sign}. 
%
%In order to achieve the boundedness of general signed solution, we ask in addition that $f$ satisfies (f5). 
%
We gain now some $L^{\infty}$-bound for sign-changing solutions, in a fractional, possibly critical, framework.
We adapt some arguments from the papers \cite{Gal0, CGT3}. 
This result will be then implemented in the study of sign-changing solutions for doubly nonlocal equations (Theorem \ref{th_INT_regular}), and in the study of uniform bounds for semiclassical critical problems (Proposition \ref{prop_upperbound}). 
Notice that we avoid the use of the Caffarelli-Silvestre $s$-harmonic extension, and this allows to extend our proof to different frameworks where this tool is not available. 
%, giving here the details for the reader's convenience. 
%\tor{sfoltisci da una delle due parti}

\begin{Proposition}\label{prop_prop_u_Linf}
Let $u\in H^s(\R^N)$ be a weak subsolution of
$$(-\Delta)^s u \leq g(x,u) \quad \hbox{in $\mathbb{R}^N$}$$
with
$$|g(x,t)| \leq C\big(|t| 
+ |t|^{2^*_s-1}\big) \quad \hbox{ for all $x \in \R^N$, $t \in \R$} $$
for some uniform $C>0$.
Then $u \in L^{\infty}(\R^N)$.
\end{Proposition}

%In particular this apply to \eqref{eq_introduction} with
%\tr{
%$$g(x,t) := \big(I_{\alpha}*F(u)\big)(x) f(t) - \mu u,$$
%}
%whenever $u$ is a fixed solution and (f5) holds (together with (f1)--(f2)), thanks to Proposition \ref{prop_conv_C0}.

%\medskip

\claim Proof.
We already know that $ u \in L^2(\R^N) \cap L^{2^*_s}(\R^N)$. Let us introduce $\gamma >1$, to be fixed, and an arbitrary $T>0$, and set a \emph{$\gamma$-linear (positive) truncation at $T$}
$$h(t)\equiv h_{T, \gamma}(t):= \parag{ &0& \quad \hbox{ if $t \leq 0$}, \\ &t^{\gamma}& \quad \hbox{if $t \in (0, T]$}, \\ &\gamma T^{\gamma-1} t - (\gamma -1) T^{\gamma}& \quad \hbox{ if $t > T$}.}$$
We have that $h\in C^1(\R)\cap W^{1,\infty}(\R)$, it is positive (increasing and convex), zero on the negative halfline, and by direct computations it satisfies the following properties
\begin{equation}\label{eq_dim_h_1}
0 \leq h(t)\leq |t|^{\gamma}, \quad t \in \R,
\end{equation}
\begin{equation}\label{eq_dim_h_2}
0\leq t h'(t) \leq \gamma h(t), \quad t \in \R,
\end{equation}
\begin{equation}\label{eq_dim_h_3}
\lim_{T\to +\infty} h_{T, \gamma}(t) = t^{\gamma}, \quad t \geq 0.
\end{equation}
The goal is to estimate $\norm{h( u)}_{2^*_s}$ and give thus a bound of $ u$ in $L^{2^*_s \gamma}(\R^N)$, where $2^*_s \gamma > 2^*_s$. In order to handle the weak formulation of the notion of solution we introduce 
$$\tilde{h}(t):= \int_0^t (h'(r))^2 \, dr, \quad t \in \R$$
and observe that $\tilde{h}\in C^1(\R) \cap W^{1, \infty}(\R)$ is positive, increasing, convex and zero on the negative halfline. In particular
\begin{equation}\label{eq_dim_htild_0}
\tilde{h}'(t) = (h'(t))^2, \quad t \in \R
\end{equation}
by definition and
\begin{equation}\label{eq_dim_htild_1}
\tilde{h}(t)-\tilde{h}(r) \leq \tilde{h}'(t) (t-r), \quad t, \,r \in \R
\end{equation}
by convexity, and we gain also the Lipschitz continuity
\begin{equation*}\label{eq_dim_htild_abs_1}
|\tilde{h}(t)-\tilde{h}(r)| \leq \norm{\tilde{h}'}_{\infty} |t-r|, \quad t, \,r \in \R.
\end{equation*}
 Combining the definition of $\tilde{h}$, \eqref{eq_dim_h_2} and \eqref{eq_dim_h_1} we obtain
\begin{equation}\label{eq_dim_htild_2}
0 \leq \tilde{h}(t) \leq \norm{h'}_{\infty} |t|^{\gamma}, \quad t \in \R.
\end{equation}
Finally, by a direct application of Jensen inequality we gain
\begin{equation}\label{eq_dim_h_htild}
|h(t)-h(r)|^2 \leq \big( \tilde{h}(t) - \tilde{h}(r)\big) (t-r) , \quad t,\, r \in \R.
\end{equation}
We observe that $\tilde{h}( u) \in H^s(\R^N)$ since $\tilde{h}$ is Lipschitz continuous and $\tilde{h}(0)=0$ (see Lemma \ref{lem_Lipschitz_Sobolev});
moreover, since $2^*_s$ is the best summability exponent, if we assume 
\begin{equation}\label{eq_dim_stima_gamma}
1<\gamma \leq \frac{2^*_s}{2}
\end{equation}
by \eqref{eq_dim_htild_2} we obtain also 
$$\tilde{h}( u) \leq \norm{h'}_{\infty} |u|^{\gamma} \in L^2(\R^N).$$

We use now the embedding \eqref{eq_embd_homog} and combine \eqref{eq_semin_gagl}, \eqref{eq_dim_h_htild} and \eqref{eq_sobolev_polaraz} %the polarized version of \eqref{eq_semin_gagl} 
to obtain
\begin{align*}
\norm{h( u)}_{2^*_s}^2 &\leq \mc{S}^{-1} \norm{(-\Delta)^{s/2} h( u)}_2^2 \\
&= (C'(N,s))^{-1} \mc{S}^{-1} \int_{\R^{2N}} \frac{|h( u(x))-h( u(y))|^2}{|x-y|^{N+2s}} \, dx \, dy \\
&\leq (C'(N,s))^{-1} \mc{S}^{-1} \int_{\R^{2N}} \frac{\big(\tilde{h}( u(x))- \tilde{h}( u(y))\big)\big( u(x)- u(y)\big)}{|x-y|^{N+2s}} \, dx \, dy \\
&= \mc{S}^{-1} \int_{\R^N} (-\Delta)^{s/2} u \, (-\Delta)^{s/2} \tilde{h}( u) \, dx.
\end{align*}
Since $\tilde{h}( u) \in H^s(\R^N)$ we can choose it as a test function in the equation and gain
$$\norm{h( u)}_{2^*_s}^2 \leq \mc{S}^{-1} \int_{\R^N} g(x,u) 
 \tilde{h}( u) \, dx.$$
By the assumptions on $g$ 
and the positivity of $\tilde{h}( u)$ we get
\begin{equation*}\label{eq_posit_to_negat}
\norm{h(u)}_{2^*_s}^2 \leq \mc{S}^{-1} \int_{\R^N} |g(x,u)| \tilde{h}( u) \, dx \leq C \mc{S}^{-1} \int_{\R^N}\big(|u|+ |u|^{2^*_s-1} \big)\tilde{h}( u) \, dx.
\end{equation*}
Since $h(u)$ and $\tilde{h}(u)$ are zero when $u$ is negative, we obtain
$$\norm{h(u_+)}_{2^*_s}^2 \leq C \mc{S}^{-1} \int_{\R^N}\big(u_++ u_+^{2^*_s-1} \big)\tilde{h}( u_+) \, dx.$$
Now we use \eqref{eq_dim_htild_1} (with $r=0$), \eqref{eq_dim_htild_0}, and \eqref{eq_dim_h_2}
\begin{eqnarray}
\lefteqn{ \norm{h(u_+)}_{2^*_s}^2 \leq C \mc{S}^{-1} \int_{\R^N} \big(u_++ u_+^{2^*_s-1} \big) u_+ \tilde{h}'( u_+) \, dx} \nonumber \\
&\leq& C \mc{S}^{-1} \int_{\R^N} \big(u_++ u_+^{2^*_s-1} \big) u_+ (h'(u_+))^2 \, dx \leq \gamma^2 C \mc{S}^{-1} \int_{\R^N} \big(1+ u_+^{2^*_s-2} \big) (h(u_+))^2 \, dx \nonumber \\
&\leq& \gamma^2 C \mc{S}^{-1} \int_{\R^N} (h(u_+))^2 \, dx + \gamma^2 C \mc{S}^{-1} \int_{\R^N} u_+^{2^*_s-2} (h(u_+))^2 \, dx. \label{eq_dim_stima_buona}
\end{eqnarray}
Let now $R>0$ to be fixed; splitting the second piece of the right-hand side of \eqref{eq_dim_stima_buona} and by using the H\"older inequality we gain
\begin{align*}
\int_{\R^N} u_+^{2^*_s-2} (h(u_+))^2 \, dx &= \int_{ u \leq R} u_+^{2^*_s-2} (h( u_+))^2 \, dx + \int_{ u>R} u_+^{2^*_s-2} (h( u_+))^2 \, dx \\
&\leq R^{2^*_s-2}\norm{h( u_+)}_2^2 +\left( \int_{ u>R} u^{2^*_s} \, dx \right)^{\frac{2^*_s -2}{2^*_s}}\norm{h( u_+)}_{2^*_s}^2 .
\end{align*}
Since $u \in L^{2^*_s}(\R^N)$, 
we can find a sufficiently large $R=R(\gamma, m_0, \mc{S}^{-1}) $ such that
$$\left( \int_{ u>R} u^{2^*_s} \, dx \right)^{\frac{2^*_s -2}{2^*_s}} < \frac{1}{2} \frac{1}{\gamma^{2} C \mc{S}^{-1}}.$$
Thus, plugging this information into \eqref{eq_dim_stima_buona}, and absorbing the second piece on the right-hand side into the left-hand side, we obtain by \eqref{eq_dim_h_1}
$$\norm{h( u_+)}_{2^*_s}^2 \leq 2\gamma^2 C \mc{S}^{-1} (1+R^{2^*_s-2})\norm{h( u_+)}_2^2 
\leq 2\gamma^2 C \mc{S}^{-1}(1+R^{2^*_s-2})\norm{ u_+}_{2 \gamma }^{2 \gamma}.$$
Recalled that $h=h_{T, \gamma}$, by \eqref{eq_dim_h_3} and Fatou's Lemma we have
\begin{align*}
\norm{ u_+}_{2^*_s \gamma}^{2 \gamma } &= \left( \int_{\R^N} \liminf_{T \to +\infty} h_{T, \gamma}^{2^*_s}( u_+) \, dx \right)^{\frac{2}{2^*_s}} 
\leq \left( \liminf_{T \to +\infty} \int_{\R^N} h_{T, \gamma}^{2^*_s}( u_+) \, dx \right)^{\frac{2}{2^*_s}} \\
&\leq 2\gamma^2 C \mc{S}^{-1} (1+R^{2^*_s-2})\norm{ u_+}_{2 \gamma }^{2 \gamma}.
\end{align*}
By our choice \eqref{eq_dim_stima_gamma} of $\gamma$ we gain that $ u_+ \in L^{2^*_s \gamma}(\R^N)$, which was the claim. By an iteration argument, with 
$$\gamma_0 := \frac{1}{2} 2^*_s, \quad \gamma_i:= \frac{1}{2}2^*_s \gamma_{i-1}, \quad \gamma_i \to +\infty,$$
 we obtain $u_+\in L^r(\R^N)$ for each $r \in [2, +\infty)$. In order to achieve $u_+ \in L^{\infty}(\R^N)$ we need to be careful on the bound on the $L^r$-norms.

Knowing that $u_+$ lies in every Lebesgue space for $r<\infty$ we can implement a more precise iteration argument, where we drop the dependence of the constant on $R$. 
We exploit once more \eqref{eq_dim_stima_buona}. Applying again Fatou's Lemma to \eqref{eq_dim_stima_buona} and using \eqref{eq_dim_h_1} we obtain
\begin{equation}\label{eq_nuovo_bootstrap}
\norm{ u_+}_{2^*_s \gamma}^{2 \gamma } \leq \gamma^2 C \mc{S}^{-1}\int_{\R^N}\big(u_+^{2 \gamma} + u_+^{2^*_s-2+2\gamma}\big)\, dx.
\end{equation}
Focusing on the second term on the right-hand side, exploiting first the generalized H\"older inequality with
$$\frac{1}{N/s} + \frac{1}{2} + \frac{1}{2^*_s} = 1,$$
possible since $u_+^{2^*_s-2} \in L^{\frac{N}{s}}(\R^N)$ because $(2^*_s-2) \frac{N}{s}= \frac{4N}{N-2s}\geq 2$, 
and the generalized Young's inequality then, we get
\begin{eqnarray*}
\lefteqn{\int_{\R^N} u_+^{2^*_s-2+2\gamma}\, dx = \int_{\R^N} u_+^{2^*_s-2} u_+^{\gamma} u_+^{\gamma} \, dx \leq \norm{u_+^{2^*_s-2}}_{\frac{N}{s}} \, \norm{u_+^{\gamma}}_2 \, \norm{u_+^{\gamma}}_{2^*_s}} \\
&\leq& \norm{u_+^{2^*_s-2}}_{\frac{N}{s}} \Big( \frac{1}{2 \eps} \norm{u_+^{\gamma}}_2^2+ \frac{\eps}{2} \norm{u_+^{\gamma}}_{2^*_s}^2 \Big) = \norm{u_+}_{\frac{4N}{N-2s}}^{2^*_s-2} \Big( \frac{1}{2 \eps} \norm{u_+}_{2 \gamma}^{2 \gamma}+ \frac{\eps}{2} \norm{u_+}_{2^*_s\gamma}^{2 \gamma} \Big) .
\end{eqnarray*}
Plugging this into \eqref{eq_nuovo_bootstrap}, set $a:= \norm{u_+}_{\frac{4N}{N-2s}}^{2^*_s-2} $, choosing $\eps = \frac{1}{a \gamma^2 C \mc{S}^{-1}}$ and bringing the $L^{2^*_s \gamma}$-norm on the left hand side, we gain
$$ \norm{ u_+}_{2^*_s \gamma}^{2 \gamma } \leq 2 \gamma^2 C \mc{S}^{-1} \big( 1 + \tfrac{1}{2}a^2\gamma^{2} C \mc{S}^{-1} \big) \norm{u_+}_{2\gamma}^{2\gamma} \leq C' \gamma^{4} \norm{u_+}_{2\gamma}^{2\gamma}$$
for some $\gamma$-independent $C'>0$. Choosing $2\gamma_i := 2^*_s \gamma_{i-1}$ we have
$$ \norm{ u_+}_{2^*_s \gamma_i} \leq \big(C' \gamma_i^{4}\big)^{\frac{1}{2\gamma_i}} \norm{u_+}_{2^*_s\gamma_{i-1}}$$
and thus
%$$ \norm{ u_+}_{2^*_s \gamma_i} 
%\leq \prod_{j=0}^i \big(C' \gamma_j^{4}\big)^{\frac{1}{2\gamma_j}} \norm{u_+}_{2^*_s\gamma_0} 
%= e^{\sum_{j=0}^i \frac{\log\big(C' \gamma_j^{4}\big)}{2 \gamma_j}} \norm{u_+}_{2^*_s\gamma_0} 
%=e^{\sum_{j=0}^i \frac{\log\big(C' \left(\frac{2^*_s}{2}\right)^{4j}\gamma_0^{4}\big)}{2 \left(\frac{2^*_s}{2}\right)^j \gamma_0}} \norm{u_+}_{2^*_s\gamma_0} 
% $$
\begin{align*}
 \norm{ u_+}_{2^*_s \gamma_i} 
&\leq \prod_{j=0}^i \big(C' \gamma_j^{4}\big)^{\frac{1}{2\gamma_j}} \norm{u_+}_{2^*_s\gamma_0} 
= e^{\sum_{j=0}^i \frac{\log\left(C' \gamma_j^{4}\right)}{2 \gamma_j}} \norm{u_+}_{2^*_s\gamma_0} \\
&=e^{\sum_{j=0}^i \frac{\log\big(C' \left(\frac{2^*_s}{2}\right)^{4j}\gamma_0^{4}\big)}{2 \left(\frac{2^*_s}{2}\right)^j \gamma_0}} \norm{u_+}_{2^*_s\gamma_0} 
\end{align*}
 and finally, sending $i\to +\infty$ (recall that $\norm{\cdot}_p\to \norm{\cdot}_{\infty}$ as $p\to+\infty$),
 $$ \norm{ u_+}_{\infty}
\leq e^{\sum_{j=0}^{\infty}\frac{\log\big(C' \left(\frac{2^*_s}{2}\right)^{4j}\gamma_0^{4}\big)}{2 \left(\frac{2^*_s}{2}\right)^j \gamma_0}} \norm{u_+}_{2^*_s\gamma_0} 
 $$
where the constant is finite. Thus $u_+ \in L^{\infty}(\R^N)$.

To deal with $u_-$ we consider
$$k(t)\equiv k_{T, \gamma}(t):=h_{T, \gamma}(-t), \quad \tilde{k}(t):= \int_t^0 (k'(r))^2 \, dr = \tilde{h}(-t)$$
and choose $\tilde{k}(u)$ as test function. With the same passages as before we obtain 
$$\norm{k( u)}_{2^*_s}^2 \leq -\mc{S}^{-1} \int_{\R^N} g(x,u) \tilde{k}( u) \, dx$$
and thus
$$\norm{k(u)}_{2^*_s}^2 \leq \mc{S}^{-1}\int_{\R^N} |g(x,u)| \tilde{k}( u) \, dx \leq C \mc{S}^{-1} \int_{\R^N}\big(|u|+ |u|^{2^*_s-1} \big)\tilde{k}( u) \, dx$$
which implies
$$\norm{k(-u_-)}_{2^*_s}^2 \leq C \mc{S}^{-1} \int_{\R^N}\big(|-u_-| +|-u_-|^{2^*_s-1} \big)\tilde{k}(-u_-) \, dx$$
and hence
$$\norm{h(u_-)}_{2^*_s}^2 \leq C \mc{S}^{-1} \int_{\R^N}\big(|u_-| +|u_-|^{2^*_s-1} \big)\tilde{h}(u_-) \, dx;$$
we then proceed as before to gain $u_- \in L^{\infty}(\R^N)$. This concludes the proof.
\QED

\bigskip

Once obtained that $u\in L^{\infty}(\R^N)$, we can improve the regularity. The following result can be found in \cite[Theorem 15]{Sti0} (see also \cite[Propositions 2.8 and 2.9]{Sil0}); see Remark \ref{rem_Zyg} for the definition of $\Lambda_1$ and $\Lambda_2$.
\begin{Proposition}\label{prop_reg_fraz}
Let $s\in (0,1)$ and $u\in %H^{2s}(\R^N) \cap 
L^{\infty}(\R^N)$ be a strong solution of
$$(-\Delta)^s u = g \quad \hbox{ in $\R^N$}.$$
%Then
\begin{itemize}
\item[i)] If $g \in L^{\infty}(\R^N)$, then 
$$u \in \parag{
&C^{0,\gamma}(\R^N) \; \hbox{ for $\gamma < 2s$} & \quad \hbox{if $2s \in (0,1) $}, \\
&\Lambda_1(\R^N), \hbox{ thus } C^{0,\gamma}(\R^N) \; \hbox{ for $\gamma < 1$} & \quad \hbox{if $2s = 1$}, \\
&C^{1,\gamma-1}(\R^N) \; \hbox{ for $\gamma < 2s$} & \quad \hbox{if $2s \in (1,2)$}.
}
$$
%\item[ii)] if $g \in C^{0,\sigma}_{loc}(\R^N)$ for some $\sigma \in (0,1]%)
%$, then 
%%$u \in C^{\gamma}_{loc}(\R^N)$ for every $\gamma \leq \sigma +2s$ \tr{with $\gamma <2$;} 
%%$$u \in \parag{
%%&C^{0,\sigma+2s}_{loc}(\R^N) & \quad \hbox{if $\sigma+2s \leq 1$}, \\
%%&C^{1,\sigma+2s-1}_{loc}(\R^N) & \quad \hbox{if $\sigma+2s \in(1,2]$}, \\
%%&\tr{C^{1,1}_{loc}(\R^N)} & \tr{\quad \hbox{if $\sigma+2s >2$}};
%%}
%%$$
%$$u \in \parag{
%&C^{0,\sigma+2s}_{loc}(\R^N) & \quad \hbox{if $\sigma+2s \leq 1$}, \\
%&C^{1,\tr{\min\{\sigma+2s-1,1\}}}_{loc}(\R^N) & \quad \hbox{if $\sigma+2s >1$}, \\
%}
%$$
\item[ii)] If $g \in C^{0,\sigma}(\R^N)$ for some $\sigma \in (0,1]$, %)
%$ \tr{with $\sigma+2s\leq2$}, 
then
% $u \in C^{\sigma+2s}(\R^N)$.
%$$u \in \parag{
%&C^{0,\sigma+2s}(\R^N) & \quad \hbox{if $\sigma+2s \leq 1$}, \\
%&C^{1,\sigma+2s-1}(\R^N) & \quad \hbox{if $\sigma+2s \in(1,2]$};
%}
%$$
$$u \in \parag{
&C^{0,\sigma+2s}(\R^N) & \quad \hbox{if $\sigma+2s \in (0,1)$}, \\
&\Lambda_1, \hbox{ thus } C^{0,\gamma}(\R^N) \; \hbox{ for $\gamma < 1$} & \quad \hbox{if $\sigma+2s = 1$}, \\
&C^{1,\sigma+2s-1}(\R^N) & \quad \hbox{if $\sigma+2s \in(1,2)$}, \\
&\Lambda_2, \hbox{ thus } C^{1,\gamma}(\R^N) \; \hbox{ for $\gamma < 1$} & \quad \hbox{if $\sigma+2s = 2$}, \\
&C^{2,\sigma+2s-2}(\R^N) & \quad \hbox{if $\sigma+2s \in(2,3)$};
}
$$
the previous relations holds also if we substitute global spaces with local spaces.
%\tr{Posso spingermi fino a $C^2$? Dalla dimostrazione di Silvestre sembra di sì, ma sembra strano (e non è nel suo enunciato..)}
%\\\tr{Che succede se $\sigma + 2s$ è $\geq 2$? Da chiarire, anche nell'articolo!}
\end{itemize}
%\tr{If $s\geq 1$ then for every $\gamma<1$ we have $u \in C^{1,\gamma}(\R^N)$ in case of $i)$ and $u \in C^{2,\gamma}(\R^N)$ in case of $ii)$. [E' vero? Vedi meglio]}
\end{Proposition}
Notice that the conclusion in $i)$ was partially contained in the embedding \eqref{eq_immers_Holder}.

\begin{Remark}
We see that regularity theory of Proposition \ref{prop_reg_fraz} extends to $s\geq 1$. Indeed, by Remark \ref{rem_emb_e_spezz}, if $u \in H^{2s}(\R^N) \cap W^{2[s], \infty}(\R^N)$ %non posso dire $W^{2s,\infty}(\R^N)$ perché non so se vale l'immersione per $p=\infty$!
then
$$(-\Delta)^s u = g \implies (-\Delta)^{[s]-s} \big((-\Delta)^{[s]} u\big) = g$$
%if $(-\Delta)^{[s]}u \in L^{\infty}(\R^N)$, then 
and all the regularity results apply to $(-\Delta)^{[s]} u$. At this point it is sufficient to apply regularity theory for polyharmonic operators \cite[Section 3.20]{Wlo0}. %to get the claims.
See also \cite[Theorem 1.2]{ROS3} and \cite[Theorem 3.7]{AJS4}.
\end{Remark}

\medskip

We want now to investigate in a more detailed way the regularity of solutions. 
Set first
\begin{equation}\label{eq_def_tail}
\Tail(u; x_0, R):=(1-s) R^{2s} \int_{\R^N \setminus B_R(x_0)} \frac{|u(x)|}{|x-x_0|^{N+2s}}dx
\end{equation}
the \emph{tail function} of $u \in H^s(\R^N)$, centered in $x_0 \in \R^N$ with radius $R>0$, introduced in \cite{DKP1, DKP2}. 
We recall properties of the fractional De Giorgi class stated in \cite{Coz0}, to which we refer for a complete introduction on the topic; we focus only on the linear case. %, where the exponent of the space is given by $2$.

By \cite[Paragraph 6.1]{Coz0} we have the following definition.

\begin{Definition}
Let $A \subset \R^N$ be open, $\zeta \geq 0$, $H\geq 1$, $k_0 \in \R$, $\mu \in (0, 2s/N]$, $\lambda \geq 0$ and $R_0 \in (0, +\infty]$. We say that $u$ belongs to the \emph{fractional De Giorgi class} $DG_+^{s,2}(A, \zeta, H, k_0, \mu , \lambda, R_0)$ if and only if
\begin{eqnarray*}
\lefteqn{[(u-k)_+]^{2}_{B_r(x_0)}+ \int_{B_r(x_0)} (u(x)-k)_+ \Bigg( \int_{B_{2R_0}(x)} \frac{(u(y)-k)_{-}}{|x-y|^{N+2s}} dy\Bigg) dx } \\
&\leq & \frac{H}{1-s} \bigg(\Big(R^{\lambda} \zeta^{2} + \frac{|k|^{2}}{R^{N \mu }}\Big) \big| \supp((u-k)_+) \cap B_R(x_0) \big |^{1- \frac{2s}{N} + \mu } +\\
&&+ \frac{R^{2(1-s)}}{(R-r)^{2}} \norm{(u-k)_+}_{L^{2}(B_R(x_0))}^{2} + \\
&& + \frac{R^{N }}{(R-r)^{N+2s}} \norm{(u-k)_+}_{L^1(B_R(x_0))} \Tail((u-k)_+; x_0, r) \bigg)
\end{eqnarray*}
for any $x_0 \in A$, $0<r<R< \min\{R_0, d(x_0, \partial A)\}$ and $k\geq k_0$. 
%tHere $(u-k)_{\pm}$ denote the positive and negative parts of the function $u-k$.
\end{Definition}

We see now how this class of functions is related to the PDE setting. By a careful analysis of the proof of \cite[Proposition 8.5]{Coz0} we obtain the following result.

\begin{Theorem}\label{thm_degiorgi_class}
Let $N\geq 2$ and let $u\in H^s(\R^N)$ be a weak subsolution of
\begin{equation*}\label{eq_eq_cozzi}
(-\Delta)^s u \leq g(x,u), \quad x \in \R^N
\end{equation*}
where $g: \R^N \times \R \to \R$ satisfies, for a.$\,$e. $x \in \R^N$ and every $t \in \R$,
$$\abs{g(x,t)} \leq d_1 + d_2 |t|^{q-1}$$
for some $q \in (2, 2^*_s)$. Then there exist $\alpha=\alpha(N,s,q)>0$, $C=C(N,s,q,d_2)>0$ and $H=H(N, s, q, d_2)\geq 1$ such that, for each $x_0\in \R^N$ and each $R_0$ verifying
$$0<R_0 \leq C(N,s,q,d_2) \min \left \{1, \norm{u}_{L^{2^*_s}(\R^N)}^{-\alpha(N,s,q)}\right\},$$
it results that
$$u \in \DG_+^{s,2} \Big( B_{R_0}(x_0), d_1, H, 0, 1-\frac{q}{2^*_s}, 2s, R_0\Big).$$
\end{Theorem}

As shown in \cite[Proposition 6.1 and Theorem 8.2]{Coz0}, the belonging to a De Giorgi class implies useful $L^{\infty}_{loc}$ and $C^{0,\sigma}_{loc}$ estimates. %, which will be implemented in some regularity pro

\medskip

For other regularity results we refer to \cite{JLX, CabSir, ROS1,BKS, FLS, Amb5}.

%%%%%%%%%%%%%%%%%%%%%%%%%%%%%%%%%%%%%%%%%%%%%%%%%%%%%%%%%%%%%%%%%%%%%%

\subsection{Existence theorems and comparison principles}

We collect here some results regarding existence and comparison principles.
%which are already known in literature, even if the author was not able to find a precise reference.

As a consequence of the Riesz representation theorem, we start by recalling the situation for linear equations in $\R^N$ \cite[page 1241, Theorem 3.3 and Lemma 4.1]{FQT} (see also \cite[Lemma C.1]{FLS}).
\begin{Lemma}[Representation in $\R^N$]\label{lem_Bessel_kernel}
Consider the equation in the weak sense
$$(-\Delta)^s u + \lambda u = g \quad \hbox{in $\R^N$}$$
where $\lambda>0$ and $g \in L^2(\R^N)$. % \tr{scrivi}$.
Then $u$ is given by
$$u = \mc{K}_{2s,\lambda} * g$$
where $\mc{K}_{2s,\lambda}$ is the \emph{Bessel Kernel} (see Remark \ref{rem_Bessel_kernels})
$$\mc{K}_{2s,\lambda}:= \mc{F}^{-1}\left(\frac{1}{\lambda +|\xi|^{2s}}\right).$$
Moreover
\begin{itemize}
\item $\mc{K}_{2s,\lambda}$ is non-negative, radially symmetric and decreasing,
\item $ \frac{C_1}{|x|^{N+2s}} \leq \mc{K}_{2s,\lambda}(x) \leq \frac{C_2}{|x|^{N+2s}}$ for $|x| \geq 1$ and some $C_1, C_2>0$, while $|\mc{K}(x)|\leq \frac{C_3}{|x|^{N-2s}}$ for $|x|\leq 1$ and some $C_3>0$,
\item $\mc{K}_{2s,\lambda} \in L^q(\R^N)$ for every $q \in [1, 1 + \tfrac{2s}{N-2s})$,
\item $\mc{K}_{2s,\lambda}$ solves $(-\Delta)^s \mc{K}_{2s,\lambda} + \lambda \mc{K}_{2s,\lambda} = \delta_0$ (in a distributional sense), where $\delta_0$ is the standard Dirac delta.
\end{itemize}
\end{Lemma}

We notice that the fundamental solution of $(-\Delta)^s u = \delta_0$, instead, is given (up to constants) by $I_{2s}:= \mc{F}^{-1}\left(\frac{1}{|\xi|^{2s}}\right)= \frac{1}{|x|^{N-2s}}$, which lies in $L^q_{loc}(\R^N)$ for every $q< \frac{N}{N-2s}$ but in no $L^p(\R^N)$. This \emph{Riesz potential} will be better studied in Section \ref{sec_riesz_poten}.

\medskip

The Bessel kernel allows also to find suitable comparison function with no restriction on the boundary; the result can be found in \cite[Lemma A.2]{CG0} (see also \cite[Lemmas 4.2 and 4.3]{FQT}). 

\begin{Lemma}[Comparison function]
\label{lem_esist_sol_part}
Let $b>0$. Then there exists a strictly positive continuous function $W_b\in H^s(\R^N)$ such that, for some positive constants $C'_b, C_b''$, it verifies 
$$(-\Delta)^s W_b + \frac{b}{2} W _b = 0,\quad x \in \R^N \setminus B_{r_b}$$
pointwise, with $r_b:= \left(\frac{2}{b}\right)^{1/2s}$, and
\begin{equation}\label{eq_stima_fun_confr}
\frac{C'_b}{|x|^{N+2s}}<W_b(x)< \frac{C_b''}{|x|^{N+2s}}, \quad \textit{ for $|x|>2 r_b$}.
\end{equation}
The constants $r_b, C'_b, C_b''$ remain bounded by letting $b$ vary in a compact set far from zero.
\end{Lemma}

\claim Proof.
Let $B_{1/2}\prec \varphi \prec B_1$, and define $\tilde{W}:= \mc{K}_{2s} * \varphi$, where $\mc{K}_{2s}$ is the Bessel potential. Arguing as in \cite{FQT} (see also \cite[Theorem 1.3]{BKS}) we obtain
$$(-\Delta)^s \tilde{W} + \tilde{W} = \varphi, \quad x \in \R^N$$
and
$$\frac{C'}{|x|^{N+2s}}<\tilde{W}(x)\leq \frac{C''}{|x|^{N+2s}} \quad \hbox{for $|x|\geq 2$}.$$
By scaling $W:=\tilde{W}(r_b \cdot)$ we reach the claim.
\QED

\bigskip

We give now an existence result (see also \cite[Corollary 1.15]{ShSp}).

\begin{Lemma}[Existence for weak solutions]\label{lem_lax_milgr}
Let $\Omega\subset \R^N$ be of class $C^{0,1}$ with bounded boundary, $\lambda>0$, $\psi \in H^{s}(\Omega^c)$, 
and $g \in L^q(\Omega)$, for some $q \in [\frac{2N}{N+2s}, 2]$.
Then there exists a (unique) function $v \in H^s(\R^N)$ such that 
$$\parag{ &(- \Delta)^s v + \lambda v = g& \quad \hbox{in $\Omega$}, \\ &v = \psi& \quad \hbox{on $\Omega^c$,}}$$
in the weak sense, which in particular means %that is 
$v \in X^s_0(\Omega)+\psi$. %and 
%$$\int_{\R^N} (-\Delta)^{s/2} v (-\Delta)^{s/2} \varphi dx + \lambda \int_{\R^N} v \varphi dx = \int_{\R^N} g \varphi dx$$
%for every $\varphi \in X^s_0(\Omega)$. 
If moreover $g \in L^{q}_{loc}(\R^N)$ for some $q>\frac{N}{2s}$, then $v\in L^{\infty}_{loc}(\R^N)$. If instead $g \in C^{0, \sigma}_{loc}(\R^N)$ for some $\sigma \in (0,1]$, then $v \in C^{2s+\sigma}_{loc}(\R^N)$.
\end{Lemma}

\begin{Remark}
The result is still valid in a whatever $\Omega^c$ \emph{extension domain} (see \cite{DnPV}). 
\end{Remark}

\claim Proof.
By \cite[Theorem 5.4]{DnPV} we know that there exists $\tilde{\psi} \in H^s(\R^N)$ such that $\tilde{\psi}_{|\Omega^c}\equiv \psi$.
The problem is thus equivalent to
$$\parag{ &(- \Delta)^s v + \lambda v = g& \quad \hbox{in $\Omega$}, \\ &v = \tilde{\psi}& \quad \hbox{on $\Omega^c$}.}$$
Consider $u= v-\tilde{\psi}$ and rewrite the weak formulation as
$$\int_{\R^N} (-\Delta)^{s/2} u (-\Delta)^{s/2} \varphi + \lambda \int_{\R^N} u \varphi = \int_{\R^N} (g-\lambda \psi) \varphi - \int_{\R^N} (-\Delta)^{s/2} \tilde{\psi} (-\Delta)^{s/2} \varphi .$$
It is easy to see that the left-hand side is a bilinear, continuous coercive map on $X^s_0(\Omega)$, while
$$\varphi \in X^s_0(\Omega) \mapsto \int_{\R^N} (g-\lambda \psi) \varphi - \int_{\R^N} (-\Delta)^{s/2} \tilde{\psi} (-\Delta)^{s/2} \varphi $$
belongs to the dual space $(X^s_0)^*(\Omega)$. By Lax-Milgram theorem, we obtain a solution $u \in X^s_0(\Omega)$, which implies that $v:=u+\tilde{\psi}$ is the desired function.

Finally, the regularity results are a consequence of De Giorgi-Nash-Moser estimates \cite[Proposition 2.6]{JLX} and Schauder estimates \cite[Theorem 2.11]{JLX}. 
\QED

\bigskip

The following existence result can be found in \cite[Lemma 2.2 and Remark 4.1]{CFQ} for bounded domains, and in \cite[Theorem A.1]{SoV} for the homogeneous case $\psi\equiv 0$.

\begin{Lemma}[Existence for viscosity solutions]\label{lem_lax_3milgr}
Let $\Omega \subset \R^N$ be a $C^2$-domain, $\lambda>0$, $\psi \in L^{\infty}(\Omega^c) \cap C(\Omega^c)$, 
and $g \in L^{\infty}(\Omega) \cap C(\overline{\Omega})$.
Then there exists a function $v \in C(\R^N) \cap L^{\infty}(\R^N)$
such that 
$$\parag{ &(- \Delta)^s v + \lambda v = g& \quad \hbox{in $\Omega$}, \\ &v = \psi& \quad \hbox{on $\Omega^c$,}}$$
in the viscosity sense. If $g \in C^{\sigma}_{loc}(\Omega)$ for some $\sigma \in (0,1)$, then $v\in C^{\gamma}_{loc}(\Omega)$, for some $\gamma>2s$ is a pointwise solution.
If $\psi \equiv 0$, we further have $v \in C^s(\R^N) \cap C^{\gamma}_{loc}(\Omega)$, for some $\gamma>\max\{1,2s\}$ and $\frac{w}{(\dist(\cdot,\partial \Omega))^s} \in C^{0,\theta}(\overline{\Omega})$ for some $\theta \in (0,1)$.
\end{Lemma}

\claim Proof.
First notice that, by extension, we may assume $g \in L^{\infty}(\R^N) \cap C(\R^N)$. % \cap C^{\sigma}_{loc}(\Omega)$.
Since $\Omega$ is a $C^2$-domain, $g\in C(\R^N)$, $\psi\in C(\Omega^c) \cap L^{\infty}(\Omega^c)$, by \cite[Theorem 4]{BCI} with $b\equiv c \equiv 0$, we obtain the existence of a (unique) viscosity solution $v \in C(\R^N)$, satisfying the boundary condition pointwise (see also \cite[page 615]{CafSil2}). 
Since the cited theorem is a corollary of \cite[Theorem 1]{BCI}, with $F(x,u,p,X,l)\equiv F(x,u,l)=l + \lambda u - g(x)$, $l=\mc{I}[u]\equiv (-\Delta)^s u$, $d \mu_x(z) = \frac{dz}{|z|^{N+\alpha}}$, one can notice, looking carefully at the proof, that the found solution is actually bounded (see also \cite[Corollary 4]{SeV}). 
Thus $v$ is a bounded viscosity solution. 

By \cite[Theorem 2.6]{QX0}, since $(-\Delta)^s v =-\lambda v + g \in L^{\infty}(\Omega)$ 
with $v\in C(\overline{\Omega})$, we have $v \in C^{\gamma_1}_{loc}(\R^N)$ for some $\gamma_1>0$. 
Since $\psi \in L^{\infty}(\Omega^c)$ and $g-\lambda v \in C^{\min\{\sigma,\gamma_1\}}_{loc}(\Omega)$, by \cite[Theorem 2.5]{QX0} we have that $v\in C^{\gamma}_{loc}(\Omega)$ for some $\gamma>2s$; thus $(-\Delta)^s v$ is pointwise defined (actually H\"older continuous). 
As observed in \cite[Remark 2.3]{QX0}, we conclude that $v$ is a pointwise solution.
\QED

\bigskip

We write down now the following two maximum principles (for unbounded domains). %again 
%for the reader's convenience. 
See \cite[Lemma A.1]{CG0} for the first 
(see also \cite[Lemma 6]{SeV} and \cite{Jar0}).

\begin{Lemma}[Maximum Principle (weak)]\label{lem_comp_prin}
Let $\Omega \subset \R^N$, %(possibly unbounded), 
$\lambda>0$, and let $u\in H^s(\R^N)$ be a weak subsolution of
$$(-\Delta)^s u + \lambda u \leq 0 \quad \hbox{in $\Omega$}.$$
Assume moreover that
$$u(x)\leq 0 \quad \hbox{for a.$\,$e. $x \in \Omega^c$}.$$
Then
\begin{equation}\label{eq_dis_comp_prin}
u(x)\leq 0 \quad \hbox{for a.$\,$e. $x\in\R^N$}.
\end{equation}
\end{Lemma}

\claim Proof.
By the assumption we have $u^+=0$ on $\Omega^c$, thus $u^+\in X^s_0(\Omega)$ is a suitable test function (see Lemma \ref{lem_dis_modul_1}) and we obtain, using $u=u^+-u^-$ and $u^+ u^-\equiv 0$,
\begin{align*}
0 &\geq \int_{\R^N} |(-\Delta)^{s/2} u^+|^2 \, dx + \lambda \int_{\R^N} |u^+|^2 \, dx - \int_{\R^N} (-\Delta)^{s/2} u^- (-\Delta)^{s/2} u^+ \, dx \\
&= \norm{(-\Delta)^{s/2}u^+}_2^2 + \lambda \norm{u^+}_2^2+C\int_{\R^{2N}} \frac{u^-(x)u^+(y)+u^-(y)u^+(x)}{|x-y|^{N+2s}} \, dx \, dy \\
&\geq \norm{(-\Delta)^{s/2}u^+}_2^2 + \lambda\norm{u^+}_2^2
\end{align*}
which implies $u^+=0$ on $\R^N$.
\QED

\medskip

\begin{Remark}
We point out that if $u$ is assumed continuous, then \eqref{eq_dis_comp_prin} is actually pointwise. Moreover, the constant $\lambda>0$ may be substituted by a more general $V(x)>0$ which gives sense to the integrals. 
\end{Remark}

\begin{Lemma}[Maximum Principle (viscosity)]\label{lem_comp_3prin}
Let $\Omega \subset \R^N$ be open, %(possibly unbounded), 
$\lambda>0$, and let $u$ be a viscosity, continuous subsolution of
$$(-\Delta)^s u + \lambda u \leq 0 \quad \hbox{in $\Omega$}$$
such that
$$\lim_{|x|\to +\infty} u(x) \leq 0.$$
%$$u(x) \to 0 \quad \hbox{as $|x|\to +\infty$}.$$
Assume moreover that
$$u(x)\leq 0 \quad \textit{on $\Omega^c$}.$$
Then
\begin{equation}\label{eq_dis_comp_2prin}
u(x)\leq 0 \quad \textit{on $\R^N$}.
\end{equation}
The result applies, in particular, to pointwise solutions.
\end{Lemma}

\claim Proof.
We first observe that $u\in L^{\infty}(\R^N)$ and set $M:=\sup_{x \in \R^N} u(x)$. By contradiction, assume $M>0$. Let $(x_n)_n$ be a maximizing sequence, i.e. $u(x_n) \to M$ as $n\to +\infty$; we can assume that $x_n \in \Omega$. 
We observe that $(x_n)_n$ is bounded (up to a subsequence) since, if not, we would have $|x_n| \to +\infty$ and thus $\lim_n u(x_n) \leq 0$, %$u(x_n) \to 0$, 
which is an absurd. 
Thus $x_n \to x_0 \in \overline{\Omega}$, and by continuity $u(x_0)=M>0$; since $u(x)\leq 0$ on $\overline{\Omega^c} \supset \partial \Omega$, we have $x_0 \in \Omega$. In particular, $x_0$ is a point of maximum for $u$. 

We can thus choose a whatever $ U\subset \Omega $ neighborhood of $x_0$ and set $\phi \equiv u(x_0)$ as contact function in the definition of viscosity solution: indeed $\phi \in C^2(U)$, $\phi(x_0)=u(x_0)$ and $\phi \geq u$ in $U$. 
Hence, set $v:= \phi \chi_{U} + u \chi_{U^c}$ we have
\begin{align*}
0& \geq (-\Delta)^s v(x_0) + \lambda v(x_0) = C_{N,s} \int_{\R^N} \frac{u(x_0) - v(y)}{|x_0-y|^{N+2s}} dy + \lambda u(x_0) \\
&= C_{N,s} \int_{U^c} \frac{M - u(y)}{|x_0-y|^{N+2s}} dy + \lambda M >0,
\end{align*}
which is a contradiction. This concludes the proof.
\QED

%%%%%%%%%%%%%%%%%%%%%%%%%%%%%%%%%%%%%%%%%%%%%%%%%%%%%%%%%%%%%%%%%%%%%%
%%%%%%%%%%%%%%%%%%%%%%%%%%%%%%%%%%%%%%%%%%%%%%%%%%%%%%%%%%%%%%%%%%%%%%

\section{The Riesz potential}
\label{sec_riesz_poten}

Let $\alpha \in (0,N)$. We recall here some results on the Riesz kernel \cite[Appendix]{MS3}
\begin{equation}\label{eq_def_Riesz}
I_{\alpha}(x) := \frac{C_{N,\alpha}}{|x|^{N-\alpha}}
\end{equation}
where 
$$C_{N,\alpha}:=\frac{\Gamma(\frac{N-\alpha}{2})}{2^{\alpha} \pi^{N/2} \Gamma(\frac{\alpha}{2})}>0$$
is a normalization constant. For motivations and a physical introduction we refer to Sections \ref{sec_intro_choquard} and \ref{sec_boson_stars}.

We are interested in studying the behaviour of the convolution
$$I_{\alpha}*g$$
for some $g$. We will use the following notation, whenever well defined for some $g$ and $h$:
$$\mc{D}_{\alpha}(g,h):= \int_{\R^N} (I_{\alpha}*g)h = \int_{\R^N} \int_{\R^N} \frac{g(x) h(y)}{|x-y|^{N-\alpha}} dx dy.$$
We start observing that the operator enjoys a trivial but useful scaling property
$$\mc{D}_{\alpha}\big(g(\theta \cdot), h(\theta \cdot)\big) = |\theta|^{-(N+\alpha)} \mc{D}_{\alpha}(g,h).$$
for any $\theta \in \R$.

%%%%%%%%%%%%%%%%%%
\subsubsection{Well posedness}

The following theorem ensures the well posedness of the Riesz potential: see \cite[Theorem 4.3]{LiLo}, \cite[pages 61-62]{Lan0} and \cite[Section 4.2]{Miz0} for a proof. % (see also \cite[Section 4.2]{Miz0}). 
%See also \cite[Appendix A.2]{MS3}, Proposition \ref{prop_conv_C0} and Remark \ref{rem_conv_welldef} for some further investigation.
% integrals involved in the variational formulation of PDEs.
\begin{Proposition}[Hardy-Littlewood-Sobolev inequality]\label{prop_HLS}
%We have the following.
Let $\alpha \in (0,N)$.
\begin{itemize}
\item Let $g$ be a measurable function. Then $I_{\alpha}*g$ is finite almost everywhere if and only if
\begin{equation}\label{eq_buona_def_riesz}
\int_{\R^N} \frac{|g(x)|}{(1+|x|)^{N-\alpha}} <\infty.
\end{equation}
In particular, $I_{\alpha}*g$ is well defined if %A condition that ensures \eqref{eq_buona_def_riesz} is 
$g \in L^1_{loc}(\R^N) \cap L^r(B_R^c)$ for some $R\geq 0$ and some $r \in [1, \frac{N}{\alpha})$.
Moreover, if \eqref{eq_buona_def_riesz} does not hold, then $I_{\alpha}*|g| \equiv \infty$. 
\item 
%Let $r,h \in (1,+\infty)$ be such that $\frac{1}{r}-\frac{1}{h}=\frac{\alpha}{N}$. Then the map
%$$g \in L^r(\R^N) \mapsto I_{\alpha}*g \in L^h(\R^N)$$
%is continuous. 
%Equivalently, 
%Let $r \in [1, \frac{N}{\alpha})$. Then the map
Let $r \in (1, \frac{N}{\alpha})$. Then, for some $C=C(N,\alpha,r)>0$ we have
$$\norm{I_{\alpha}*g}_{\frac{Nr}{N-\alpha r}} \leq C \norm{g}_r$$
for all $g \in L^r(\R^N)$, 
thus the map
$$g \in L^r(\R^N) \mapsto I_{\alpha}*g \in L^{\frac{Nr}{N-\alpha r}}(\R^N)$$
is continuous.
 In particular, since the operator is linear,
$$g_n \wto g \hbox{ in $L^r(\R^N)$} \implies I_{\alpha}*g_n \wto I_{\alpha}*g \hbox{ in $L^{\frac{Nr}{N-\alpha r}}(\R^N)$}.$$
\item Let $r,t \in (1,+\infty)$ be such that $\frac{1}{r}+\frac{1}{t}= \frac{N+\alpha}{N}$. Then there exists a constant $C=C(N,\alpha,r,t)>0$ such that
%$$\pabs{\int_{\R^N} (I_{\alpha}*g)h dx} 
$$\pabs{\mc{D}_{\alpha}(g,h)}\leq C \norm{g}_r \norm{h}_t$$
for all $g \in L^r(\R^N)$ and $h \in L^t(\R^N)$. Thus the bilinear map
$$(g,h) \in L^r(\R^N)\times L^t(\R^N) \mapsto \mc{D}_{\alpha}(g,h) \in \R$$
is continuous. 
If $r=t=\frac{2N}{N+\alpha}$, then equality is reached in the previous inequality if and only if $g \equiv h$ (up to multiplicative constants), and $g(x)=(1+|x|^2)^{-\frac{N+\alpha}{2}}$ (up to translations and rescaling).
\end{itemize}
In the limiting case $g \in L^{\frac{N}{\alpha}}(\R^N)$ (i.e. $\frac{Nr}{N-\alpha r} \to \infty$) we have that $I_{\alpha}*g$ is a BMO function (see \cite[Appendix A.2]{MS3} and references therein). Anyway we have
\begin{itemize}
\item If $g \in L^{\frac{N}{\alpha} -\eps}(\R^N)\cap L^{\frac{N}{\alpha} +\eps}(\R^N)$ for some $\eps>0$, then $I_{\alpha}*g \in C_0(\R^N) \subset L^{\infty}(\R^N)$.
\end{itemize}
\end{Proposition}
%
%\tor{
%We remark that $I_{\alpha}*g$ is well defined (almost everywhere pointwise) for $g \in L^r(\R^N)$ with $r \in [1, \frac{N}{\alpha})$.
%}

\claim Proof.
We show only the last claim, i.e. \cite[Lemma 4.5(ii)]{LoMa08}; we argue as in \cite[Proposition 4.5]{CGT3} (see also Proposition \ref{prop_conv_C0} and Remark \ref{rem_conv_welldef}). 
Recall theat, by Young's Theorem, if two functions belong to two Lebesgue spaces with conjugate (finite) indexes, then their convolution belong to $C_0(\R^N)$. 
% $g*h \in C_0(\R^N)$. 
We first split
$$I_{\alpha}*g = (I_{\alpha}\chi_{B_1})*g + (I_{\alpha}\chi_{B_1^c})*g$$
where
$$ I_{\alpha}\chi_{B_1} \in L^{r_1}(\R^N), \quad \hbox{ for $r_1 \in [1, \frac{N}{N-\alpha})$},$$
$$ I_{\alpha}\chi_{B_1^c} \in L^{r_2}(\R^N), \quad \hbox{ for $r_2 \in (\frac{N}{N-\alpha}, \infty]$}.$$
We need that $g \in L^{q_1}(\R^N)\cap L^{q_2}(\R^N)$ for some $q_i$ satisfying
$$\frac{1}{q_i} + \frac{1}{r_i} = 1, \quad i=1,2$$
that is
$$\frac{q_1}{q_1-1} \in \left[1, \frac{N}{N-\alpha}\right), \quad \frac{q_2}{q_2-1}\in \left(\frac{N}{N-\alpha}, \infty\right]$$
or equivalently $q_2 < \frac{N}{\alpha} < q_1$. Thus we have the claim.
\QED

\bigskip

%We highlight that, generally, in the limiting case $g \in L^{\frac{N}{\alpha}}(\R^N)$ (i.e. $\frac{Nr}{N-\alpha r} \to \infty$) we have that $I_{\alpha}*g$ is a BMO function (see \cite[Appendix A.2]{MS3} and references therein); anyway, if we assume that $g \in L^{\frac{N}{\alpha} -\eps}(\R^N)\cap L^{\frac{N}{\alpha} +\eps}(\R^N)$ for some $\eps>0$, then $I_{\alpha}*g \in C_0(\R^N) \subset L^{\infty}(\R^N)$ (see Remark \ref{rem_conv_welldef}).

We emphasize the similarity of condition \eqref{eq_buona_def_riesz} and condition \eqref{eq_spazio_gener_s}, when formally $\alpha = -2s$.

%\medskip

\subsubsection{Positivity}

We observe the following: if $g \in \mc{S}$ \cite[Lemma 5.1.2]{Ste0} or if $\alpha \in (0, \frac{N}{2})$ and $g \in L^{\frac{2N}{N+2\alpha}}(\R^N)$ \cite[Corollary 5.10]{LiLo} then we have 
$$\mc{D}_{\alpha}(g,g)=\int_{\R^N} \big(I_{\alpha}*g\big) g dx = \int_{\R^N} \widehat{I_{\alpha}*g} \widehat{g} d\xi= \int_{\R^N} \widehat{I_{\alpha}} |\widehat{g}|^2 d\xi = \int_{\R^N} \frac{|\widehat{g}|^2}{|\xi|^{\alpha}} d\xi \geq 0$$
(see also \cite[Lemma 4.5(v)]{LoMa08}, \cite[Lemma 2.7]{Buc0}, \cite[Section 1.1]{Lan0}, \cite[Theorem 2.8 in Chapter 2]{Miz0}, \cite[Sections 2.1.1 and 2.3.3]{Sam0} and \cite[Theorem 5.9]{LiLo}).
This shows that
$$g \mapsto \mc{D}_{\alpha}(g,g)$$ %\int_{\R^N} \big(I_{\alpha}*g\big) g $$
is a \emph{positive functional} (i.e. its sign does not depend on the sign of $g$). 
A more general result can be adapted from $\alpha=2$ \cite[Theorem 9.8]{LiLo} to a generic $\alpha \in (0,N)$ as follows.
\begin{Proposition}[\cite{LiLo}]\label{prop_positivity_Riesz}
Let $g:\R^N \to \R$ measurable be such that
%$$\int_{\R^N} \big(I_{\alpha}*|g|\big) |g| 
$$\mc{D}_{\alpha}(|g|,|g|)< \infty.$$
Then
%$$\int_{\R^N} \big(I_{\alpha}*g\big) g 
$$\mc{D}_{\alpha}(g,g)\geq 0 $$
and the above quantity is zero if and only if $g\equiv 0$ almost everywhere.
In particular the following representation holds %\tr{questa vale per $\alpha=2$!}
%$$\int_{\R^N} \big(I_{\alpha}*g\big) g = \int_0^{+\infty} \frac{1}{t^{\alpha+1}} \int_{\R^N} \abs{ \big(t^N h(t\cdot)\big)* g }^2 dx dt \geq 0 $$
%$$\int_{\R^N} \big(I_{\alpha}*g\big) g 
$$\mc{D}_{\alpha}(g,g)= \int_0^{+\infty} t^{2N-\alpha-1} \int_{\R^N} \abs{ h(t\cdot)* g }^2 dx dt \geq 0 $$
for a whatever nonnegative, radially symmetric $h \in C^{\infty}_c(\R^N)$ normalized in such a way that $\int_0^{+\infty} t^{N-\alpha-1} (h*h)(t) dt = C_{N,\alpha}$. %\frac{1}{2}$.
\end{Proposition}

\subsubsection{Decay}

We investigate now the decay of $I_{\alpha}*g$: indeed, if $g \in L^1_{loc}(\R^N)$, $g\geq 0$ and $g>0$ on some ball, then
$$(I_{\alpha} * g)(x) \geq I_{\alpha}(2x) \int_{B_{2|x|}(x)} g \gtrsim I_{\alpha}(x) \simeq \frac{1}{|x|^{N-\alpha}} \quad \hbox{for $|x|\gg0$}$$
which shows a polynomial bound from below on the Riesz potential, whatever the decay of $g$ is (even with compact support). 
%In particular we see that if $\alpha in [\frac{N}{2},N)$, then $I_{\alpha}*g$ might not lie $L^1$.
Moreover, if $g\geq 0$ %\in L^1_{loc}(\R^N)$ 
has at least a polynomial decay
$$ g(x) \lesssim \frac{1}{|x|^{\theta}} \quad \hbox{ as $|x|\to +\infty$}$$
with $\theta>\alpha$, then the following estimates from above hold \cite[Lemma A.1]{MS1} (see also \cite[Lemma 2.1]{GMM} and \cite[Lemma 4.6]{GM0})
\begin{equation}\label{eq_decay_Riesz}
(I_{\alpha}*g)(x) \lesssim \parag{
&\frac{1}{|x|^{\theta-\alpha}} & \quad \hbox{ if $\theta \in (\alpha,N)$}, \\
&\frac{\log(x)}{|x|^{N-\alpha}} & \quad \hbox{ if $\theta =N$}, \\
&\frac{1}{|x|^{N-\alpha}} & \quad \hbox{ if $\theta \in (N, +\infty)$}.
}
\end{equation}
In particular, if $\theta>N$ the decay of $I_{\alpha}*g$ is exactly the same of $I_{\alpha}$, as stated 
%; if $g \in L^1(\R^N)$ we further have
%$$I_{\alpha}*g \sim \norm{g}_1 I_{\alpha} \quad \hbox{ as $|x|\to +\infty$};$$
%this last estimate is furhter investigated 
in the following result.

\begin{Lemma}[\cite{MS0}]\label{lem_stima_1_Riesz}
Let $g \in L^{\infty}(\R^N)$ be continuous and $\theta>N$ be such that
$$\sup_{x \in \R^N} |g(x)| |x|^{\theta} < +\infty.$$
Then there exists $C=C(N,\alpha)>0$ such that
%$$\pabs{\int_{\R^N} \frac{g(y)}{|x-y|^{N-\alpha}} dy - \frac{1}{|x|^{N-\alpha}} \int_{\R^N} g(y) dy} \leq \frac{C \norm{g}_{\infty,\theta}}{|x|^{N-\alpha}} \left( \frac{1}{1+|x|} + \frac{1}{1+ |x|^{\theta-N}}\right)$$
$$\pabs{(I_{\alpha}*g)(x) - % \norm{g}_1 
I_{\alpha}(x) \int_{\R^N} g(y) dy} \leq \frac{C \norm{g}_{\infty,\theta}}{|x|^{N-\alpha}} \left( \frac{1}{1+|x|} + \frac{1}{1+ |x|^{\theta-N}}\right)$$
for each $x \in \R^N$, $x \neq 0$, where we recall that $\norm{g}_{\infty, \theta}= \norm{g(\cdot )(1+|\cdot|^{\theta})}_{\infty}$.
\end{Lemma}

\claim Proof.
See \cite[Lemma 6.2]{MS0}. See also \cite[Lemma C.3]{FLS}.
\QED

\bigskip

The rigidity of the previous result in particular highlights that it is not possible to implement a \emph{bootstrap-type} argument in order to show fine results on the decay of a solutions. See Section \ref{sec_est_bel_vis}.

%%%%%%%%%%%%%%%%%%%%%%%%%%%%%%%%%%%%%%
\subsection{The Riesz potential as the inverse of the fractional Laplacian}
\label{sec_riesz_vs_frac}

%\tr{When $\alpha \in (0,2)$} 
Since roughly 
$$(-\Delta)^{\alpha/2} I_{\alpha} = \mc{F}^{-1}\left( |\xi|^{\alpha} \widehat{I_{\alpha}}(\xi)\right) = \mc{F}^{-1}\left( |\xi|^{\alpha} \frac{1}{|\xi|^{N-(N-\alpha)}}\right) = \mc{F}^{-1}(1) = \delta_0 $$
then the Riesz kernel can be seen as the fundamental (distributional) solution for the fractional Laplacian \cite[Theorem 5.10]{AJS3} (see also \cite[Theorem 8.4]{Gar0} and \cite[Theorem 2.3]{Buc0} for the case $\alpha \in (0,2)$) %for a formal proof when $\alpha \in (0,2)$)
\begin{equation}\label{eq_fond_delta}
(-\Delta)^{\alpha/2} I_{\alpha} = \delta_0 \quad \hbox{in $\R^N$};
\end{equation}
thus the Riesz potential %actually 
generates the solutions of fractional equations in $\R^N$, 
%$$g \mapsto \phi = I_{\alpha}*g$$
%actually gives the solutions to the equations
%$$(-\Delta)^{\alpha/2} \phi = g \quad \hbox{in $\R^N$}.$$
that is
$$ \phi = I_{\alpha}*g \iff (-\Delta)^{\alpha/2} \phi = g \quad \hbox{in $\R^N$}.$$
Therefore we may roughly say that (\cite[Section 5.1]{Ste0}, \cite[equation (2.7)]{Gar0}, \cite[equation (2.3)]{Sil0} and \cite[equation (1.2.7)]{AH0})
$$I_{\alpha} * \equiv (-\Delta)^{-\alpha/2}.$$
More precisely 
$$I_{\alpha}*((-\Delta)^{\alpha/2} v) = v = (-\Delta)^{\alpha/2}(I_{\alpha}*v) \quad \hbox{ for every $v \in C^{\infty}_c(\R^N)$};$$
indeed, when the fractional Laplacian is defined through hypersingular integrals, the first equality can be found in \cite[proof of Theorem 2.9 in Chapter 2 and Section 4.5]{Miz0} for Schwartz functions, while the second equality for $L^p$ functions in \cite[Theorems 3.22 and 3.24]{Sam0}: anyway the hypersingular definition coincides with the Fourier transform one at least on $C^{\infty}_c$ functions (see \cite[Lemma 3.1]{Sam0} and \cite[Theorem 1.8]{AJS1}). See also \cite[equation (1.1.12')]{Lan0}. For the case $\alpha \in (0,2)$ see also \cite[Theorem 2.8 and Corollary 2.9]{Buc0} and \cite[Theorem 6]{Sti0}, while for $\alpha=2$ see also \cite[Lemma 4.5(iii)]{LoMa08}.
% for $\alpha \in (0,2)$ with $g \in L^{\infty}(\R^N)$ with compact support.
%Osservazione: potrebbe sembrare ovvio, ma non lo è, poiché la formula di Plancharel su $\int_{\R^N} (I_{\alpha}*g) (-\Delta)^s \varphi$ non la posso applicare (per ottenere $\int_{\R^N} \frac{1}{|\xi|^{\alpha}} \widehat{g} |\xi|^{2s} \widehat{\varphi}$) poiché $(-\Delta)^s \varphi$ non è in generale una funzione Schwartz (vedi tutti i riferimenti nella sezione Positività).

Let us state this relation in a more general framework. % in a rigorous way. %: see and

\begin{Proposition}\label{prop_inversa_dx}
Let $\alpha \in (0, N)$ (i.e., set $s:=\frac{\alpha}{2}$, we ask $N>2s$). 
\begin{itemize}
\item[i)] Assume %$g \in L^1(\R^N)$ with compact support
$g \in L^p(\R^N)$ for some $p \in [1, \frac{N}{\alpha})$. %\tr{serve supporto compatto?}
Then
$$(-\Delta)^{\alpha/2} \big(I_{\alpha}*g\big) = g \quad \hbox{in $\R^N$}$$
in the strong sense; notice, in particular, that the fractional Laplacian of $\phi=I_{\alpha}*g$ is well defined pointwise (i.e., finite) almost everywhere. 
Moreover, if %$g$ is a measurable function being continuous at $x$ (or more generally, $x$ is a Lebesgue point for $g$), %il che significa che vale quasi ovunque per qualsiasi funzione in $L^1_{loc}(\R^N)$!
 $x$ is a Lebesgue point for $g$ (e.g., $g$ is continuous at $x$), 
then the previous relation holds at $x$. 
\item[ii)] If $g \in L^p(\R^N) \cap X$ for some $p \in [1, \frac{N}{\alpha})$ and some function space $X$, then $(-\Delta)^{\alpha} \big(I_{\alpha}*g\big) \in L^p(\R^N) \cap X$; 
in particular if $g \in L^p(\R^N)\cap L^{\frac{N}{N-\alpha p}}(\R^N)$, then
$$I_{\alpha}*g \in W^{\alpha, \frac{N}{N-\alpha p}}(\R^N).$$
\item[iii)] %Assume $\alpha \in (0,2)$ and 
Let $g \in L^p(\R^N)$ for some $p \in [1, \frac{N}{\alpha})$. 
Then $\phi=I_{\alpha}*g$ is the only (distributional) solution to
$$(-\Delta)^{\alpha/2} \phi = g \quad \hbox{in $\R^N$}$$
belonging to $L^{\frac{Np}{N-\alpha p}}(\R^N)$.
\item[iv)] %Assume $\alpha \in (0,2)$ and 
Let $\phi \in W^{\alpha, \frac{N}{N-\alpha p}}(\R^N)$ for some $p\in [1, \frac{N}{\alpha})$; assume moreover that $(-\Delta)^{\alpha/2} \phi \in L^p(\R^N)$. 
Then
$$ I_{\alpha}* \big((-\Delta)^{\alpha/2} \phi \big) = \phi \quad \hbox{in $\R^N$}$$
in the strong sense.
\end{itemize}
\end{Proposition}

\medskip

\claim Proof.
Point $i)$ is stated in \cite[page 22, Definition 2.5 and Proposition 7.1]{Kwa0} (see also \cite[Corollary 5.16]{AJS3} for compactly supported $g\in L^1(\R^N)$ and \cite[Corollary 5.10]{LiLo} for $\alpha \in (0,\frac{N}{2})$ and $g \in L^{\frac{2N}{N+2\alpha}}(\R^N)$); see instead \cite[Theorems 3.22 and 3.24%, and Lemma 3.1
]{Sam0}, and \cite[Theorem 5.1 and Remark 5.1 in Chapter 4]{Miz0} for a hypersingular approach. 
 % for a proper proof of $\implies$ in the case $\alpha \in (0, \frac{N}{2})$ and $g \in L^{\frac{2N}{N+\alpha}}(\R^N)$, and notice that $\frac{2N}{N+\alpha} <\frac{N}{\alpha}$.);
Point $ii)$ is a direct consequence.

To show $iii)$, by linearity it is sufficient to prove the statement for $g=0$; this can be done as in \cite[Theorems 1.3 and Theorem 3.1]{CDL}. See also \cite[Corollary 1.4]{Fal0}, \cite[Corollary 1.3]{FaWe0} and \cite[Theorem 1.5]{DOV} for $\alpha \in (0,2)$, \cite[Theorem 5.17]{AJS3} for $\alpha \notin 2\N$ and \cite{Hui0} for $\alpha \in 2\N$.

We give some details only on $iv)$ (see also \cite{CFW}). %, even if a similar result is already known in the literature. 
Indeed, consider
$$(-\Delta)^{\alpha/2} \phi = 0 \quad \hbox{in $\R^N$};$$
by $iii)$ we know that the only solution $\phi \in W^{\alpha, \frac{Np}{N-\alpha p}}(\R^N)$ is the null function. Thus the kernel of the linear operator
$$(-\Delta)^{\alpha/2} : W^{\alpha, \frac{Np}{N-\alpha p}}(\R^N) \to L^{\frac{Np}{N-\alpha p}}(\R^N)$$
is null, and hence the operator injective. In particular, considered the homogeneous space
$$\dot{W}^{\alpha,p}(\R^N):=\left \{ u \hbox{ measurable} \mid (-\Delta)^{\alpha/2} u \in L^p(\R^N)\right\}$$
we have that
$$(-\Delta)^{\alpha/2} : \dot{W}^{\alpha,p}(\R^N) \cap W^{\alpha, \frac{Np}{N-\alpha p}}(\R^N) \to L^p(\R^N) \cap L^{\frac{Np}{N-\alpha p}}(\R^N) $$
is injective, and thus admits a left inverse. On the other hand, by $ii)$ we have
$$I_{\alpha}* : L^p(\R^N) \cap L^{\frac{Np}{N-\alpha p}}(\R^N) \to \dot{W}^{\alpha,p}(\R^N) \cap W^{\alpha, \frac{Np}{N-\alpha p}}(\R^N)$$
and moreover, by $i)$, it is a right inverse for $(-\Delta)^{\alpha/2}$. Therefore the left and right inverse must coincide, which means that $I_{\alpha}*$ is a right inverse for $(-\Delta)^{\alpha/2}$. 
This concludes the proof.
%if $N\geq 2\alpha$, then $0 \in L^{\frac{2N}{N+2\alpha}}(\R^N)$ and by $ii)$ there exists a unique (distributional) solution $\phi\in L^2(\R^N)$, which is $\phi \equiv 0$, while if $N\leq 2\alpha$ then $0 \in L^1(\R^N)$ and again by $i)$ there exists a unique (distributional) solution $\phi \in L^{\frac{N}{N-\alpha}}(\R^N)$, $\phi \equiv 0$. Thus we see that, generally, the kernel of $(-\Delta)^{\alpha/2}$ on $H^{\alpha}(\R^N)\cap L^1(\R^N) \subset L^2(\R^N) \cap L^{\frac{N}{N-\alpha}}(\R^N)$ is null (in the strong sense). Therefore the linear operator
%$$L: H^{\alpha}(\R^N)\cap L^1(\R^N) \to L^2(\R^N), \quad L:=(-\Delta)^{\alpha/2}$$
%in injective, and 
%%since the null function belong to $L^{q}(\R^N)$ for $q \in [\frac{N}{N-\alpha}, \frac{2N}{N-4 \alpha}]$, with $\frac{2N}{N+2\alpha \in [1, \frac{N}{\alpha})$, we know by $i)$ that there exists a unique solution in $L^2(\R^N)$
%%Indeed, by $i)$ we know that 
%%$$L: H^{\alpha}(\R) \to L^2(\R), \quad L:=(-\Delta)^{\alpha/2}$$
%%is injective, since the kernel is null. 
%thus it admists a left inverse. On the other hand, by $i)$, we know that $I_{\alpha}$ is a right inverse over $L^p(\R^N)$.
%%
%We need two spaces $X$ and $Y$ such that
%$$ X \subset H^{\alpha}\cap L^1 \subset L^2 \cap L^{\frac{N}{N-\alpha}} \stackrel{ (-\Delta)^{\frac{\alpha}{2}}} \to Y \subset L^1 \cap L^{\frac{N}{\alpha})} \stackrel{I_{\alpha}} \to X.$$
%\tr{Schwartz non va bene, perché il Laplaciano non lo conserva}
\QED

\medskip

\begin{Remark}
By the previous proof, we see that, if $p\in [1, \frac{N}{\alpha})$, then
$$I_{\alpha} * \equiv (-\Delta)^{-\alpha/2}.$$
when looked on the spaces $\dot{W}^{\alpha,p}(\R^N) \cap W^{\alpha, \frac{Np}{N-\alpha p}}(\R^N)$ and $ L^p(\R^N) \cap L^{\frac{Np}{N-\alpha p}}(\R^N)$.
\end{Remark}

%%%%%%%%%%%%%%%%
\subsubsection{Regularity}

As already seen by point $ii)$ of Proposition \ref{prop_inversa_dx}, the Riesz potential has a regularizing effect. We give more details in the following result.
\begin{Proposition}\label{prop_regol_riesz}
Let $\alpha \in (0,N)$ and $p \in [1, \frac{N}{\alpha})$, and let $g \in L^p(\R^N)$.
\begin{itemize}
%\item if $g \in L^p(\R^N) \cap X$ for some function space $X$, then $(-\Delta)^{\alpha} \big(I_{\alpha}*g\big) \in L^p(\R^N) \cap X$; 
%in particular if $g \in L^p(\R^N)\cap L^{\frac{N}{N-\alpha p}}(\R^N)$, then
%$$I_{\alpha}*g \in W^{\alpha, \frac{N}{N-\alpha p}}(\R^N);$$
%
%\item if $\alpha \in (1, \frac{N}{2}]$, $N\geq 2 \alpha$, and $g \in L^{\frac{N}{2\alpha}}(\R^N) \cap L^2(\R^N)$, then
%$$I_{\alpha}*g \in H^{\alpha}(\R^N);$$
%
%\item if $\alpha \in (1, 2)$, $N\geq 2 \alpha$, and $g \in L^{\frac{N}{2\alpha}}(\R^N) \cap L^2(\R^N) \cap L^{\infty}(\R^N)$, 
%
%
%\item[i)] Assume $g \in L^q(\R^N)$ for some $q\in (\frac{N}{\alpha},\infty)$ with $\alpha- \frac{N}{q} \in (0,1)$.
%%$\frac{N}{q} < \alpha < 1+\frac{N}{q}$. 
%Then $I_{\alpha}*g \in Lip(\alpha-\frac{N}{q})$. 
%\\In particular, if we assume a priori also $I_{\alpha}*g \in L^{\infty}(\R^N)$, then $I_{\alpha}*g \in C^{\alpha-\frac{N}{q}}(\R^N)$.
\item[i)] Assume $g \in L^q(\R^N)$ for some $q\in (\frac{N}{\alpha},\infty)$ with $\alpha- \frac{N}{q} \in (0,1)$. %, satysfing \eqref{eq_buona_def_riesz}.
%$\frac{N}{q} < \alpha < 1+\frac{N}{q}$. 
Then 
%$I_{\alpha}*g \in Lip(\alpha-\frac{N}{q})$. 
%\\In particular, if we assume a priori also $I_{\alpha}*g \in L^{\infty}(\R^N)$, then 
$I_{\alpha}*g \in C^{\alpha-\frac{N}{q}}(\R^N)$.
\item[ii)] Assume $\alpha \in (0,2)$ and 
$g \in L^{\infty}(\R^N)$, 
and we assume a priori that $I_{\alpha}*g \in L^{\infty}(\R^N)$.
Then
$$I_{\alpha}*g \in \parag{
&C^{0,\gamma}(\R^N) \; \hbox{ for $\gamma < \alpha$} & \quad \hbox{if $\alpha \in (0,1) $}, \\
&\Lambda_1(\R^N), \hbox{ thus } C^{0,\gamma}(\R^N) \; \hbox{ for $\gamma < 1$} & \quad \hbox{if $\alpha = 1$}, \\
%&C^{1,\gamma-1}(\R^N) \; \hbox{ for $\gamma < \min\{2,\alpha\}$} & \quad \hbox{if $\alpha >1$}.
&C^{1,\gamma-1}(\R^N) \; \hbox{ for $\gamma < \alpha$} & \quad \hbox{if $\alpha \in (0,2)$}.
}
$$
%
%\item if $\alpha \in (1, 2)$, $N\geq 2 \alpha$, and $g \in L^{\frac{N}{2\alpha}}(\R^N) \cap L^2(\R^N) \cap C^{0,\sigma}(\R^N)$ for some $\sigma \in (0,1]$, 
%
\item[iii)] Assume $g \in C^{0,\sigma}(\R^N) \cap L^{\infty}(\R^N)$ for some $\sigma \in (0,1]$, and $g\geq 0$. Then $I_{\alpha}*g \in C^{0,\gamma}(\R^N)$ for each $\gamma < (1-\frac{\alpha}{N} p) \sigma$.\footnote{Actually, if in addition $\int_{\R^N} g^q < \infty$ for some $q\in (0,1)$, then we can take $\gamma < (1-\frac{\alpha}{N} q) \sigma$. In particular, if $q=\frac{N}{N+\alpha}$, then $\gamma<\frac{N}{N+\alpha} \sigma$.}
%\frac{N-\alpha p}{N} \sigma$.
%Holder (o anche uniformemente continua) e L^p implicano limitata, ad esempio (adattandola al caso multidimensionale): https://mathproblems123.wordpress.com/2009/10/01/barbalats-lemma/
%\item[iii)] 
\item[iv)] Assume $\alpha \in (0,1)$ and $g \in C^{0,\sigma}(\R^N)$ % Lip(\sigma)$ 
for some $\sigma \in (0,1)$ such that $\sigma+ \alpha \in (0,1)$. 
Then $I_{\alpha}*g \in Lip(\sigma + \alpha)$. 
\\In particular, if we assume a priori also $I_{\alpha}*g \in L^{\infty}(\R^N)$, then $I_{\alpha}*g \in C^{\sigma + \alpha}(\R^N)$.
%
%\item[iv)]
\item[v)] Assume $\alpha \in (0,2)$ and % Assume %$\alpha \in (0,2)$ and 
$g \in C^{0,\sigma}(\R^N)$ for some $\sigma \in (0,1]$, 
and we assume a priori that $I_{\alpha}*g \in L^{\infty}(\R^N)$. 
Then
$$I_{\alpha}*g \in \parag{
&C^{0,\sigma+\alpha}(\R^N) & \quad \hbox{if $\sigma+\alpha \in (0,1)$}, \\
&\Lambda_1, \hbox{ thus } C^{0,\gamma}(\R^N) \; \hbox{ for $\gamma < 1$} & \quad \hbox{if $\sigma+\alpha = 1$}, \\
&C^{1,\sigma+\alpha-1}(\R^N) & \quad \hbox{if $\sigma+\alpha \in(1,2)$}, \\
&\Lambda_2, \hbox{ thus } C^{1,\gamma}(\R^N) \; \hbox{ for $\gamma < 1$} & \quad \hbox{if $\sigma+\alpha = 2$}, \\
%&C^{2,\gamma-2}(\R^N) \; \hbox{ for $\gamma \leq \sigma + \alpha$, $\gamma <3$} & \quad \hbox{if $\sigma+\alpha >2$};
&C^{2,\sigma+\alpha-2}(\R^N) & \quad \hbox{if $\sigma+\alpha \in (2,3)$};
}
$$
the previous relations holds also if we substitute global H\"older spaces with local spaces.
\end{itemize}
\end{Proposition}

\claim Proof.
Point $i)$ can be found in \cite[Theorem 2.2 in Section 4.2]{Miz0} (see also \cite[Theorem 2]{DuP0} and \cite[Theorem 1.6]{ROS1}); point %$iii)$ 
$iv)$ can be found in \cite[Theorem 1]{DuP0} and Remark \ref{rem_DuPl}.
Points $ii)$ and %$iv)$ 
$v)$ are consequences of Proposition \ref{prop_reg_fraz} and Proposition \ref{prop_inversa_dx}.

We are left to prove $iii)$. Indeed, let 
$r:= \frac{\sigma}{\gamma}>\frac{N}{N-\alpha p}>1$. We can find thus $\theta \in (1, \frac{N}{\alpha})$ such that $(1-\frac{1}{r}) \theta \geq p $. We thus have, for $x, y,z \in \R^N$, exploiting $|a^r-b^r| \lesssim |a-b| |a^{r-1}-b^{r-1}|$
\begin{align*}
|g(x-z) - g(y-z)| & \lesssim |(g(x-z))^{\frac{1}{r}} - (g(y-z))^{\frac{1}{r}}||(g(x-z))^{\frac{r-1}{r}} - (g(y-z))^{\frac{r-1}{r}}| \\
& \lesssim |g(x-z) - g(y-z)|^{\frac{1}{r}}|(g(x-z))^{\frac{r-1}{r}} - (g(y-z))^{\frac{r-1}{r}}| \\
& \lesssim |x-y|^{\frac{\sigma}{r}}|(g(x-z))^{\frac{r-1}{r}} - (g(y-z))^{\frac{r-1}{r}}|;
\end{align*}
as a consequence %(noticed that $g \in L^{\infty}(\R^N)$)
\begin{align*}
\abs{(I_{\alpha}*g)(x) - (I_{\alpha}*g)(y)} \lesssim & \, |x-y|^{\gamma} \int_{\R^N} \frac{|(g(x-z))^{\frac{r-1}{r}} - (g(y-z))^{\frac{r-1}{r}}|}{|y|^{N-\alpha}} dz \\
 \leq & \, |x-y|^{\gamma} \Big(  \int_{B_1(0)} \frac{|g(x-z)|^{\frac{r-1}{r}} +|g(x-z)|^{\frac{r-1}{r}}}{|y|^{N-\alpha}} dz + \\
& +  \int_{B_1^c(0)} \frac{|g(x-z)|^{\frac{r-1}{r}} +|g(x-z)|^{\frac{r-1}{r}}}{|y|^{N-\alpha}} dz \Big) \\
 \leq & \, |x-y|^{\gamma} \Big( 2 \norm{g}_{\infty}^{\frac{r-1}{r}} \int_{B_1(0)} \frac{1}{|y|^{N-\alpha}} dz +  2 \norm{g}_{\frac{r-1}{r} \theta}^{\frac{r-1}{r}}  \int_{B_1^c(0)} \frac{1}{|y|^{(N-\alpha)\theta'}} dz \Big)% \\
%\leq & C |x-y|^{\gamma}
\end{align*}
which gives the claim, being $(N-\alpha)\theta'>N$ and $\frac{r-1}{r} \theta \in [p,\infty)$. % $g \in L^{\frac{r-1}{r} \theta}(\R^N)$.
\QED

%\smallskip

%%%%%%%%%%%%%%%%%%
\subsubsection{Limiting cases}

We wonder what happens when $\alpha \to 0$ or $\alpha \to N$. In the first case, the Riesz potential collapses into a local operator (as one may expect from the representation $I_{\alpha} *g\equiv \mc{F}^{-1}(|\xi|^{-\alpha} \widehat{g})$), that is
\begin{equation}\label{eq_alpha_to_0}
I_{\alpha} * g \stackrel{\alpha \to 0^+} \to \delta_0* g=g ;
\end{equation}
 in the second case, as one may imagine looking at the Poisson equation \eqref{eq_fond_delta} with $\alpha=N$ 
%$$(-\Delta)^{N/2} u=g \quad \hbox{in $\R^N$}$$
(for example, in the planar case $N=2$ with the classical Laplacian), the Riesz kernel converges (up to constants) to a logarithm kernel 
\begin{equation}\label{eq_alpha_to_N}
I_{\alpha} * g \stackrel{\alpha \to N^-} \to \log\big(\tfrac{1}{|\cdot|}\big) *g 
\end{equation}
whenever computed on a function with zero mean $\int_{\R^N} g=0$. 
%Differently
%$$I_{\alpha} * \big(g - \int_{\R^N} g\big) \stackrel{\alpha \to N^-} \to \log\big(\tfrac{1}{|\cdot|}\big) * \big(g - \int_{\R^N} g\big) ;$$
See \cite[pages 46 and 50]{Lan0} for precise statements.

%VEDI anche [Comi-Stefani2021, Proposition 4.3]♥

%\bigskip

%%%%%%%%%%%%%%%%%%%%%%%%%%%%%%%%%%%%
\subsubsection{Definitions of solutions}

The definitions of weak and viscosity solutions apply, mutatis mutandis, to nonlocal equations of the type
\begin{equation}\label{eq_main_prel}
(-\Delta)^s u + \mu u = (I_{\alpha}*F(u))f(u) \quad \hbox{on $\R^N$}
\end{equation}
where we ask $u$ to satisfy the equation in the classic/strong/weak/viscosity sense with nonlinearity $g(x):= (I_{\alpha}*F(u))(x)f(u(x))$.
When dealing with weak solution, we need the term to be summable (see Remark \ref{rem_buona_posit_f}); while, when dealing with classical and viscosity solutions, we need $I_{\alpha}*F(u)$ to be well defined pointwise (see Remark \ref{rem_conv_welldef}).
%by substituting \eqref{eq_weak_sol} with 
%$$\int_{\R^N} (-\Delta)^{s/2} u (-\Delta)^{s/2} \varphi dx + \mu \int_{\R^N} u \varphi dx \leq \int_{\R^N} \big(I_{\alpha} * F(u)\big) f(u) \varphi dx,$$
%whenever we assume that the right-hand side is well defined, 
%%where we implicitly assume \hyperref[(f1)]{\textnormal{(f1)}}-\hyperref[(f2)]{\textnormal{(f2)}} to give sense to the integrals, 
%and by substituting \eqref{eq_visc_sol} with
%$$(-\Delta)^s v(x_0) \leq \big(I_{\alpha}*F(u)\big)(x_0) f(u(x_0)) ;$$
%in this last case, we need some assumptions on $f$ and $u$ to have $I_{\alpha}*F(u)$ well defined pointwise, see Remark \ref{rem_conv_welldef}.
%Similarly for strong a classical solutions.

We notice that, under the assumptions of invertibility of the fractional Poisson equation (see Proposition \ref{prop_inversa_dx}), equation \eqref{eq_main_prel} can be rewritten as a fractional Schr\"odinger-Newton system
\begin{equation}\label{eq_SN_syst}
\parag
{(-\Delta)^s u + \mu u = \phi f(u) \quad \hbox{ in $\R^N$}, \\
(-\Delta)^{\alpha/2} \phi = F(u) \quad \hbox{ in $\R^N$}.
}
\end{equation}

%%%%%%%%%%%%%%%%%%%%%%%%%%%%%%%%%%%%%%%%%%%%%%%%%%%%%%%%%%%%%%%%%%%%%%
%%%%%%%%%%%%%%%%%%%%%%%%%%%%%%%%%%%%%%%%%%%%%%%%%%%%%%%%%%%%%%%%%%%%%%

\section{Some manipulations: absolute value and polarization}
\label{sec_absol_pol}

If one considers a function $u\in H^1(\R^N)$ and its absolute value, it is easy to see that
$$|\nabla |u|| = |\nabla u|.$$
%for every $u \in H^1(\R^N)$. 
Actually, the equality is not the case of the fractional Laplacian, generally. 
We show thus how the fractional Laplacian, and the Riesz potential, behave with respect to the absolute value. 
\begin{Lemma}\label{lem_dis_modul_1}
Let $s\in (0,1)$ and $\alpha \in (0,N)$.
\begin{itemize}
\item Let $u \in H^s(\R^N)$. Then $|u| \in H^s(\R^N)$ and
$$\norm{(-\Delta)^{s/2} |u|}_2 \leq \norm{(-\Delta)^{s/2} u}_2.$$
As a consequence, if $u=u_+ - u_-$, then $u_{\pm} = \frac{|u|\pm u}{2} \in H^s(\R^N)$.
\item Let $F: \R \to \R$ continuous and $u:\R^N \to \R$ measurable be such that $\mc{D}_{\alpha}\big(F(|u|), F(|u|)\big)< \infty$.
\\If $F$ is even, then
$$\mc{D}_{\alpha}\big(F(|u|), F(|u|)\big) = \mc{D}_{\alpha}\big(F(u), F(u)\big);$$
if $F$ is odd and has constant sign on $(0,+\infty)$, then
$$\mc{D}_{\alpha}\big(F(|u|), F(|u|)\big) \geq \mc{D}_{\alpha}\big(F(u), F(u)\big).$$
\end{itemize}
%\tr{
%\begin{itemize}
%\item $F$ is even, or
%\item $F$ is odd, and $F$ has constant sign on % $f$ has constant sign on 
%$(0, +\infty)$,
%\end{itemize}
%}
%% and that $f$ is even or odd, 
%then
%$$\mc{D}(|u|)\geq \mc{D}(u);$$
%if $F$ is even, equality holds. As a consequence
%$$\mc{J}_{\mu}(|u|)\leq \mc{J}_{\mu}(u), \quad \mc{P}_{\mu}(|u|)\leq\mc{P}_{\mu}(u).$$
\end{Lemma}

\claim Proof.
By \eqref{eq_semin_gagl} we have
\begin{align*} % Package amsmath Error: \begin{split} won't work here.
\norm{(-\Delta)^{s/2} |u|}_2^2 &= C_{N,s} \int_{\R^{2N}} \frac{\big(|u(x)|-|u(y)|\big)^2}{|x-y|^{N+2s}} \, dx dy \\
&= C_{N,s} \int_{\R^{2N}} \frac{|u|^2(x) + |u|^2(y) - 2|u|(x)|u|(y)}{|x-y|^{N+2s}} \, dx dy \\
&\leq C_{N,s} \int_{\R^{2N}} \frac{u^2(x) + u^2(y) - 2u(x)u(y)}{|x-y|^{N+2s}} \, dx dy \\
&= C_{N,s} \int_{\R^{2N}} \frac{\big(u(x)-u(y)\big)^2}{|x-y|^{N+2s}} \, dx dy = \norm{(-\Delta)^{s/2}u}_2^2,
\end{align*}
thus the first claim.
Focus on the second claim: if $F$ is odd and with constant sign on $(0,+\infty)$, then, set for brevity $A^{\pm}:=\{\pm u >0\}$, 
\begin{eqnarray*} % Package amsmath Error: \begin{split} won't work here.
\lefteqn{\mc{D}_{\alpha}\big(F(|u|), F(|u|)\big)} \\ %=& \int_{\R^{2N}} I_{\alpha}(x-y) F(|u(x)|) F(|u(y)|) \, dx dy \\
&=& \int_{A^+ \times A^+} I_{\alpha}(x-y) F(u(x)) F(u(y)) -\int_{A^- \times A^+} I_{\alpha}(x-y) F(u(x)) F(u(y)) - \\
&& - \int_{A^+ \times A^-} I_{\alpha}(x-y) F(u(x)) F(u(y)) + \int_{A^- \times A^-} I_{\alpha}(x-y) F(u(x)) F(u(y)) \\
&\geq& \int_{A^+ \times A^+} I_{\alpha}(x-y) F(u(x)) F(u(y)) + \int_{A^- \times A^+} I_{\alpha}(x-y) F(u(x)) F(u(y)) + \\
&&+ \int_{A^+ \times A^-} I_{\alpha}(x-y) F(u(x)) F(u(y)) + \int_{A^- \times A^-} I_{\alpha}(x-y) F(u(x)) F(u(y)) \\
&=& %\int_{\R^{2N}} I_{\alpha}(x-y) F(u(x)) F(u(y)) \, dx dy = 
\; \mc{D}_{\alpha}\big(F(u), F(u)\big),
\end{eqnarray*}
which concludes the proof, observing that equality holds if $F$ is instead even.
\QED

\bigskip

\begin{Remark}\label{rem_conv_deb_assol} 
We highlight that, in Sobolev spaces, the absolute value conserves the weak convergences. 
Indeed, assume $u_k \wto u$ in $H^s(\R^N)$. Since $u_k$ is bounded and 
%then by using the fact that $-|u_k(x)||u_k(y)|\leq -u_k(x)u_k(y)$ we obtain 
$\norm{|u_k|}_{H^s(\R^N)}\leq \norm{u_k}_{H^s(\R^N)}$ we have that $\abs{u_k}$ is bounded too. Therefore, $\abs{u_k}\wto v$ in $H^s(\R^N)$ up to a subsequence. 
As a consequence, up to a subsequence, $u_k \to u$ and $|u_k| \to v$ almost everywhere, which means that $\abs{u}=v$ almost everywhere. This means that %, up to a subsequence, 
$\abs{u_k}\wto \abs{u}$ in $H^s(\R^N)$.

Notice that, in general, for weak convergences in $L^p$-spaces the implication is not true \cite[Section 5]{Za0}.
\end{Remark}

We turn now to the study of symmetries. We exploit the tool of the polarization, useful in the presence of the Riesz potential.
%%
%\\
Let 
$$\mc{H}:= \big\{ H\subset \R^N \hbox{closed half-space}, \; 0 \in H\big\}.$$
For any $H \in \mc{H}$ let $\sigma_H$ be the reflection with respect to $\partial H$. The \emph{polarization} (or \emph{two-points symmetrization}) of a function $u:\R^N \to \R$ is defined as
$$u^H(x):=\parag{ 
&\max\{ u(x), u(\sigma_H(x))\big\} & \quad \hbox{if $x \in H$}, \\
&\min\{ u(x), u(\sigma_H(x))\big\} & \quad \hbox{if $x \notin H$}.
}
$$
For example, if $u=\chi_{\Omega}$, with $\Omega \subset \R^N$ crossing $\partial H$, then
$$(\chi_{\Omega})^H(x)=\parag{ 
&\chi_{\Omega \cup \sigma_H(\Omega)}(x) & \quad \hbox{if $x \in H$}, \\
&\chi_{\Omega \cap \sigma_H(\Omega)}(x) & \quad \hbox{if $x \notin H$},
}
$$
which roughly means that $u^H$ brings mass from $H^c$ to $H$.
%%
%\\
One can see \cite{Bur0}, \cite[Section 8.3]{Wil1} and \cite{VaS0} and references therein for an introduction on the topic and some relations with the symmetric decreasing rearrangement.

Clearly we have
$$u^H \equiv u \iff u \geq u \circ \sigma_H \quad \hbox{on $H$},$$
$$u^H \equiv u \circ \sigma_H \iff u \leq u \circ \sigma_H \quad \hbox{on $H$}$$
which means, roughly, that there is more mass of $u$ on $H$ than on $H^c$. We expect that, if $u$ coincide with $u^H$ for all the hyperplanes, then some symmetry must hold. 
This is actually stated in the following result \cite[Lemma 5.4]{MS0} (see also \cite[Proposition 3.15]{VSW0} and \cite[Lemma 6.3]{BS0}).
\begin{Proposition}
Let $u \in L^p(\R^N)$, for some $p \in [1,+\infty)$, be nonnegative. Then $u$ is radially symmetric %up to a translation 
if and only if for every $H \in \mc{H}$ it results that $u^H=u$, while $u$ is radially symmetric up to a translation if and only if for every $H\in \mc{H}$ it results that $u^H=u$ or $u^H = u\circ \sigma_H$.
\end{Proposition}
We state now a proposition which shows both how the Riesz potential behaves with respect to polarization, and why this tool is particularly effective in this framework \cite[Lemma 5.3]{MS0}.
\begin{Proposition}
Let $\alpha \in (0,N)$ and $H\in \mc{H}$, and let $g \in L^{\frac{2N}{N+\alpha}}(\R^N)$ be nonnegative. Then
$$\mc{D}_{\alpha}(g^H, g^H) \geq \mc{D}_{\alpha}(g,g)$$
and equality holds if and only if $u^H \equiv u$ or $u^H \equiv u \circ \sigma_H$.
\end{Proposition}
We investigate now how the fractional Laplacian behaves with respect to polarization, see \cite[equation (2.14)]{BBDG} and \cite[Lemma 5.3]{BS0} (see also \cite[page 4818%Theorem 3 (corollary of)
]{Bec0}).
\begin{Proposition}\label{prop_prel_polar_lapl}
Let $s \in (0,1)$ and $H \in \mc{H}$, and let $ u \in H^s(\R^N)$. Then $u^H \in H^s(\R^N)$ and
$$\norm{(-\Delta)^{s/2} u^H}_2 \leq \norm{(-\Delta)^{s/2} u}_2.$$
When $s=1$, the equality holds. % (see ).
\end{Proposition}
%\noindent Actually, 
%
%\\
Finally, it is easy to verify that \cite[Proposition 8.3.7]{Wil1}, for every $p \in [1,+\infty)$,
$$\norm{u^H}_p = \norm{u}_p$$
and that, % a direct consequence of the definition that, 
if $F:\R \to \R$ is nondecreasing, then
\begin{equation}\label{eq_prel_pol_F_incr}
F(u^H)=(F(u))^H.
\end{equation}

We refer to \cite{LS0} for other interesting results about manipulations of nonlocal quantities.

%\begin{Proposition}
%Let $F:\R \to \R$ be nondecreasing on $(0,+\infty)$ and $u:\R^n \to \R$ be nonnegative. Then
%$$F(u^H)=(F(u))^H$$. % \quad \hbox{whenever $v \geq 0$}.$$
%\end{Proposition}

%Let $u^H$ be the polarization of $u$ with respect to a closed half-space $H\subset \R^N$. Observe, , that we have
%
%
%Assume $f$ has constant sign on $(0, +\infty)$ in addition to \textnormal{(f1)-(f2)}. 
%Assume moreover that $f \geq 0$ on $(0,+\infty)$ (again, if we substitute $f$ with $-f$ the Hartree-type terms are conserved). Observe that $F$ is nondecreasing on $(0,+\infty)$; this means that

%%%%%%%%%%%%%%%%%%%%%%%%%%%%%%%%%%%%%%%%%%%%%%%%%%%%%%%
%%%%%%%%%%%%%%%%%%%%%%%%%%%%%%%%%%%%%%%%%%%%%%%%%%%%%%%

\section{Berestycki-Lions type assumptions: some convergences} % and comments}
\label{sec_prelim_assump}

The assumptions considered throughout the thesis are in the spirit of the ones proposed by Berestycki and Lions \cite{BL1,BL2}, adapted then to the fractional framework in \cite{ChWa0, BKS} and to the Choquard-Hartree-Pekar framework by Moroz and Van Schaftingen \cite{MS2}. These assumptions cover different models which arise in applications, see Examples \ref{ex_esempi_nonlin}.

In the case of the unconstrained problem (frequency fixed, mass free), as shown in the abovementioned papers (see also \cite{JT0,MS0,DSS1}), these assumptions are somehow \emph{almost optimal}, in the sense that when the nonlinearity collapses to a power, the growth condition are optimal for the existence of a (sufficiently regular) variational solution. See also \cite{Stra0, CMP, ZhZo1} for the case of combination of powers and %anyway 
\cite{Med2} for some further generalizations to the so called infinity-mass regime. 

We highlight that no pointwise condition of Ambrosetti-Rabinowitz type, nor of monotonicity type, is assumed, and this lack of additional assumptions obstructs some classical arguments related both to compactness and geometry of the problems.

In the $L^2$-constrained case (frequency free, mass prescribed), different qualitative phenomena are related to sub and super $L^2$-critical cases: for instance, the sub or super $L^2$-criticality of the exponent influences the boundedness of the functional %$\mc{J}$ on $\mc{S}_m$, 
on the $L^2$-sphere, as well as %and also 
the lifespan and the stability of the solutions in some related equations %the time-dependent Choquard equation 
(see \cite{CaLi0}).
In this thesis we restrict our analysis to the $L^2$-subcritical regime: we aim to extend our results to the $L^2$-critical and supercritical regime in the future.

\medskip

In this Section, for the sake of clarity, we list all the assumptions on the nonlinearities that will come into play in the following Chapters, both in the fractional framework and in the Choquard framework; %, in this Chapter or following ones. 
we let here $s \in (0,1]$ and $\alpha \in (0,N)$.
In particular, we show the role of the subcriticality growth in the convergence of nonlinear functionals. 

\smallskip

We highlight that the labeling here introduced will be changed throughout different Chapters, in order to avoid cumbersome notations.

\subsection{Local nonlinearities} %Assumptions on the nonlinearity}

%\subsubsection{Local nonlinearities}

%In the rest of thesis we will use the following notations:
%$$2^{\#}:=2$$
%$$2^*_s : =2 + \frac{4s}{N-2s} = \frac{2N}{N-2s}$$
%$$2^{m}_s:=2+ \frac{4s}{N} = \frac{2N+4s}{N}$$ 
%$$2^{\#}_{\alpha}:= 1 + \frac{\alpha}{N}=\frac{N+\alpha}{N}$$
%$$2^{*}_{\alpha,s}:= 1 + \frac{\alpha+2s}{N-2s}= \frac{N+\alpha}{N-2s}$$
%$$2^{m}_{\alpha,s}:= 1 + \frac{\alpha+2s}{N}= \frac{N+\alpha+2s}{N}$$
%$2^{*}_{\alpha,1} \equiv 2^{*}_{\alpha} = 1 + \frac{\alpha+2s}{N-2}= \frac{N+\alpha}{N-2}$, $2^{m}_{\alpha,1} \equiv 2^{m}_{\alpha} = 1 + \frac{\alpha+2}{N}= \frac{N+\alpha+2}{N}$
%\bigskip

%Let $s \in (0,1]$. 
For local nonlinearities of the type $g(u)$, $G(t)= \int_0^t g(\tau) d\tau$, we introduce the following notations:
\begin{itemize}
\item \emph{Lower critical exponent:} $2^{\#}:=2$,
\item \emph{Upper critical exponent:} $2^*_s : =2 + \frac{4s}{N-2s} = \frac{2N}{N-2s} \in (2, +\infty)$,
\item \emph{$L^2$-critical exponent:} $2^{m}_s:=2+ \frac{4s}{N} = \frac{2N+4s}{N} \in (2, 2+\frac{4}{N})$, %\min\{2+\frac{4}{N},2^*_s\})$,
\end{itemize}
and notice that
$$2=2^{\#} < 2^{m}_s < 2^*_s < +\infty.$$

Moreover we introduce the following set of assumptions:
\begin{itemize}
\item[(h0)] \label{(h0)}
 \emph{Continuity:} $g \in C(\R)$,
\item[(h0')] \label{(h0')}
\emph{Pohozaev regularity:} $s \in (\frac{1}{2},1)$ or $g \in C^{\sigma}_{loc}(\R)$ for some $\sigma > 1-2s$,
%$g \in C(\R)$ if $s \in (\frac{1}{2},1)$, $g \in C^{\sigma}_{loc}(\R)$ for some $\sigma > 1-2s$ if $s\in (0, \frac{1}{2}]$,
\item[(h1)] \label{(h1)}
\emph{Nontriviality (frequency free):} there exists $t_0>0$ such that $G(t_0)>0$,
\item[(h1')] \label{(h1')}
\emph{Nontriviality (frequency $\mu>0$ fixed):} there exists $t_0=t_0(\mu)>0$ such that $G(t_0) \geq \frac{\mu}{2} t_0^2$,

\smallskip

%\hyperref[(g1)]{\textnormal{(g1)}}

\item[(h2)] \label{(h2)}
\emph{Supercriticality in $0$:} $\lim_{t \to 0} \frac{g(t)}{t} =0$,
\item[(h2*)] \label{(h2*)}
\emph{$L^2$-subcriticality in $0$:} 
%$\lim_{t \to 0} \frac{g(t)}{|t|^{2^{m}_s-1}} =+\infty$ \tr{basta $\lim_{t \to 0} \frac{G(t)}{|t|^{2^{m}_s}} =+\infty$?},
$\lim_{t \to 0} \frac{G(t)}{|t|^{2^{m}_s}} =+\infty$, %\tr{basta $\lim_{t \to 0} \frac{G(t)}{|t|^{2^{m}_s}} =+\infty$?},

\smallskip

\item[(h3)] \label{(h3)}
\emph{Subcriticality at $\infty$:} $\lim_{|t| \to +\infty} \frac{g(t)}{|t|^{2^*_s-2}t} =0$,
\item[(h3')] \label{(h3')}
\emph{Subcriticality (strict) at $\infty$:} $\lim_{|t| \to +\infty} \frac{g(t)}{|t|^{p-2}t} =0$ for some $p \in (2, 2^*_s)$,
\item[(h3'')] \label{(h3'')}
\emph{$L^2$-subcriticality at $\infty$:} $\lim_{|t| \to +\infty} \frac{g(t)}{|t|^{2^{m}_s-2}t} =0$, 
%\\ \tr{basta $\limsup_{|t| \to +\infty} \frac{G(t)}{|t|^{2^{m}_s}} =\limsup_{|t| \to +\infty} \frac{g(t)}{|t|^{2^{m}_s-1}} =0$?} %CCOMMENT NOW
%\tr{NON basta $\lim_{|t| \to +\infty} \frac{G(t)}{|t|^{2^{m}_s}} =0$} %COMMENT NOW
\item[(h3*)] \label{(h3*)}
\emph{Criticality %(not pure)
 at $\infty$:} $\lim_{|t| \to +\infty} \frac{g(t)}{|t|^{2^*_s-2}t} =a\neq 0$; if $a>0$ we also assume $g(t) \geq a t^{2^*_s-1} + Ct^{p-1}$ 
for some $C>0$ and $p \in (\max\{2^*_s-2s,2\}, 2^*_s)$ and every $t>0$,

\smallskip

\item[(h4)] \label{(h4)}
\emph{Symmetry:} $g$ odd,
\item[(h5)] \label{(h5)}
\emph{Negative-cut (for positivity):} $g\equiv 0$ on $(-\infty,0]$.
\end{itemize}

% \hyperref[(h2*)]{\textnormal{($h2^*$)}}
 
Notice that
$$ \hyperref[(h2*)]{\textnormal{(h2*)}}
\vee 
\hyperref[(h3*)]{\textnormal{(h3*)}}
\implies 
\hyperref[(h1')]{\textnormal{(h1')}}
\implies 
\hyperref[(h1)]{\textnormal{(h1)}}, 
\quad 
\hyperref[(h3'')]{\textnormal{(h3'')}}
\implies 
\hyperref[(h3')]{\textnormal{(h3')}} 
\implies 
\hyperref[(h3)]{\textnormal{(h3)}}.$$

\begin{Example}\label{ex_esempi_nonlin}
These general assumptions include different models arising in applications. For examples, they cover
%Some nonlinear models satisfying ... are given by 
pure powers $g(t)=|t|^{q-2}t$, with $q \in (2, 2^*_s)$ (or $q \in (2, 2^m_s)$), and combined powers like $g(t)=|t|^{q-2} t + |t|^{r-2} t$ (cooperation models) and $g(t)=|t|^{q-2} t - |t|^{r-2} t$ (competion models). %, with $2\leq r < q < 2^*_s$ (or $2^m_s$). 
Other physical models can be found for example in asymptotically linear functions 
$$g(t)= \frac{t^3}{1+ t^2}, \quad G(t)= \frac{1}{2} \left(t^2 - \log(1+t^2)\right),$$
which arise in the saturation effect in nonlinear optics for photorefractive media \cite{DLQWZPH,WCCYZH,MMP,HLRZ,RAAY}, 
%Gatz, S., Herrmann, J.: Propagation of optical beams and the properties of two-dimensional spatial solitons in media with a local saturable nonlinear refractive index. J. Opt. Soc. Am. B 14, 1795–1806 (1997)
%Kelley, P.L.: Self-focusing of optical beams. Phys. Rev. Lett. 15, 1005 (1965)
%Marburger, J.H., Dawesg, E.: Dynamical formation of a small-scale filament. Phys. Rev. Lett. 21(8), 556–558 (1968)
%Merhasin, I.M., Malomed, B.A., Senthilnathan, K., Nakkeeran, K.,Wai, P.K.A., Chow, K.W.: Solitons in Bragg gratings with saturable nonlinearities. J. Opt. Soc. Am. B 24, 1458–1468 (2007)
%Petrovi´c, M.S., Beli´c, M.R., Denz, C., Kivshar, Y.S.: Counterpropagating optical beams and solitons. Laser Photonics Rev. 5, 214–233 (2011)
or also
$$g(t)= \left( 1- \frac{1}{\sqrt{1+t^2}}\right)t, \quad G(t) =\frac{1}{2} \left( t^2 - 2 \sqrt{1+t^2} +2\right) $$
of square-root type, which describes narrow-gap semiconductors \cite{PPB, SBB}.
%Petroskia, M.M., Petrovica, M.S., Belica, M.R.: Quasi-stable propagation of vortices and soliton clusters in previous saturable Kerr media with square-root nonlinearity. Opt. Commun. 279, 196202 (2007)
%Skarka, V., Berezhiani, V.I., Boucher, V.: Self-trapping of electromagnetic pulses in narrow-gap semiconductors. Phys. Rev. B 61, 10201 (2000)
\end{Example}

\begin{Remark}
We trivially observe that assigning a condition on $g$ is generally stronger than assigning a similar condition on $G$. 
Indeed, by De l'Hôpital theorem, % dalla generalizzazione del teorema, non c'è bisogno di assumere: if $G(t) \to +\infty$ as $t \to +\infty$, then
$$\lim_{|t|\to +\infty} \frac{g(t)}{|t|^{q-2}t}= l \in \overline{\R} \implies \lim_{|t|\to +\infty} \frac{G(t)}{|t|^q}= l,$$
or more generally
$$\liminf_{|t|\to +\infty} \frac{g(t)}{|t|^{q-2}t} \leq \liminf_{|t|\to +\infty} \frac{G(t)}{|t|^q} \leq \limsup_{|t|\to +\infty} \frac{G(t)}{|t|^q} \leq \limsup_{|t|\to +\infty} \frac{g(t)}{|t|^{q-2}t}.$$
The viceversa is generally not true: consider for example 
%$G(t)=t^{q-1}(t+\cos(t))$ 
%$G(t)=\frac{t^q \cos(t)}{\log(t+2)}$ (for $t\geq 0$) 
%$G(t)=t^{q-\frac{1}{2}} \cos(t)$
%which verifies
%$$\lim_{t \to+\infty} \frac{G(t)}{t^q}=0, \quad \limsup_{t \to+\infty} \frac{g(t)}{t^{q-1}}=+\infty, \quad \liminf_{t \to + \infty} \frac{g(t)}{t^{q-1}} = - \infty.$$
$G(t):=t^q \left(\int_0^{t}\frac{\sin(\tau)}{\tau} - \frac{\pi}{2}\right)$
which verifies
$$\lim_{t \to+\infty} \frac{G(t)}{t^q}=0, \quad \liminf_{t \to+\infty} \frac{g(t)}{t^{q-1}}=-1, \quad \limsup_{t \to + \infty} \frac{g(t)}{t^{q-1}} = 1;$$
%as $t\to+\infty$.
notice that the $\limsup$ is finite (consider $G(t)=t^{q-\frac{1}{2}} \cos(t)$ for an infinite $\limsup$). 
On the other hand, if one assume a priori that $\lim_{|t|\to +\infty} \frac{g(t)}{|t|^{q-2}t}$ exists, then the viceversa holds true. %\tr{E se assumo che il limsup sia finito?}

Moreover, if $\delta \in (0,1)$, by choosing $\eps \in (0,1-\delta)$ and setting $G(t)=t^{q-\eps}\cos(t)$ we see that 
$$\lim_{|t|\to +\infty} \frac{G(t)}{|t|^q} =0 \quad \hbox{ but } \quad \lim_{|t|\to +\infty} \frac{g(t)}{|t|^{q+\delta-1}} \neq 0;$$
in particular, since generally $2^*_s - 2^m_s \in (0,+\infty)$, we have
$$\lim_{|t|\to +\infty} \frac{G(t)}{|t|^{2^m_s}} =0 \centernot \implies \lim_{|t|\to +\infty} \frac{g(t)}{|t|^{2^*_s-1}} =0.$$

Similar considerations can be done for $t\to 0$ (consider $G(t):=t^q \left(\int_0^{1/t}\frac{\sin(\tau)}{\tau} - \frac{\pi}{2}\right)$ or $G(t)=t^{q+\eps} \cos\big( \frac{1}{t}\big)$).
%\\ Moreover
%$$\lim_{t\to 0} \frac{g(t)}{|t|^{q-2}t}= l \in \overline{\R} \implies \lim_{t\to 0} \frac{G(t)}{|t|^q}= l,$$
%or more generally
%$$\liminf_{t\to 0} \frac{g(t)}{|t|^{q-2}t} \leq \liminf_{t\to 0} \frac{G(t)}{|t|^q} \leq \limsup_{t\to 0} \frac{G(t)}{|t|^q} \leq \limsup_{t\to 0} \frac{g(t)}{|t|^{q-2}t}.$$
%The viceversa is generally not true: consider for example 
%%$G(t)=t^{q+\frac{1}{2}} \cos\big( \frac{1}{t}\big)$
%%which verifies
%%$$\lim_{t \to 0} \frac{G(t)}{t^q}=0, \quad \limsup_{t \to 0} \frac{g(t)}{t^{q-1}}=+\infty, \quad \liminf_{t \to 0} \frac{g(t)}{t^{q-1}} = - \infty.$$
%$G(t)=t^q \int_0^{1/t}\frac{\sin(\tau)}{\tau}$
%which verifies
%$$\lim_{t \to 0} \frac{G(t)}{t^q}=0, \quad \liminf_{t \to 0} \frac{g(t)}{t^{q-1}}=-1, \quad \liminf_{t \to 0} \frac{g(t)}{t^{q-1}} =1.$$
%On the other hand, if one assume a priori that $\lim_{t\to 0} \frac{g(t)}{|t|^{q-2}t}$ exists, then the viceversa holds true. 
%% \tr{E se assumo che il limsup sia finito?}
%\\ \tr{IMPORTANTE: Nel caso nonlocale impongo su $F$, ma credo serva imporre anche su $f$. Nel caso locale invece impongo su $g$ (se il caso nonlocale è effettivamente giusto, potrei imporre anche qui su $G$?)}
\end{Remark}

\begin{Remark}\label{rem_somma_lebesg}
Generally, when $u \in H^s(\R^N)$, $g(u)$ %and $G(u)$ 
will not lie on a precise Lebesgue space, but on a summation of spaces. To handle these quantities we remark 
%We observe 
that the following properties are equivalent \cite[Proposition 2.3]{BPR}, for any $p, q \in (1,+\infty)$:
\begin{itemize}
\item $g \in L^p(\R^N) + L^q(\R^N)$,
\item $|g| \in L^p(\R^N) + L^q(\R^N)$,
\item $|g| \leq h$ for some $h \in L^p(\R^N) + L^q(\R^N)$.
\end{itemize}
\end{Remark}

\begin{Remark}\label{rem_buona_posit_g}
%Since of key importance in the good definition of the functionals, as well as in bootstrap argument in the rest of the thesis, w
We write here in which spaces lie the considered quantities. 
Let $u\in H^s(\R^N) \subset L^2(\R^N) \cap L^{2^*_s}(\R^N)$. By assuming
\begin{equation}\label{eq_g_wellposed}
\limsup_{t\to 0} \frac{|g(t)|}{|t|} <\infty, \quad \limsup_{|t| \to +\infty} \frac{|g(t)|}{|t|^{2^*_s}} < \infty
\end{equation}
(for instance given by \hyperref[(h2)]{\textnormal{(h2)}} and \hyperref[(h3)]{\textnormal{(h3)}})
we have (see Remark \ref{rem_somma_lebesg})
\begin{align*}
g(u) &\in L^{2}(\R^N)\cap L^{\frac{2N}{N-2s} }(\R^N)+ L^{2 \frac{N-2s}{N+2s}}\cap L^{\frac{2N}{N+2s}}(\R^N) \\
&\subset L^{2}(\R^N) + L^{\frac{2N}{N+2s}}(\R^N),\\
G(u) &\in L^{1}(\R^N)\cap L^{\frac{N}{N-2s}}(\R^N) + L^{2 \frac{N-2s}{N}}(\R^N)\cap L^{1}(\R^N) \\
&\subset L^{1}(\R^N).
\end{align*}
%Thus, seen as an operator, $G: L^2(\R^N) \cap L^{2^*_s}(\R^N) \to L^1(\R^N)$.
If $\varphi \in H^s(\R^N) \subset L^2(\R^N) \cap L^{2^*_s}(\R^N)$ is a test function, we notice that the found summability is enough to have
$\int_{\R^N} g(u) \varphi \, dx$ 
well defined. 
\end{Remark}

We state now the convergence properties of the nonlinear functionals, in the case of a subcritical growth \cite{ChWa0} (see also \cite[Theorem 2 and Corollary 2]{MT1}). %; we give some details for the reader's convenience. %we have the following
\begin{Proposition}
\label{prop_converg_generiche_loc}
 Assume \hyperref[(h0)]{\textnormal{(h0)}} and \eqref{eq_g_wellposed}.
\begin{itemize}
\item Let $u_n \wto u$ in $H^s(\R^N)$. Then 
%$$g(u_n) \wto g(u) \quad \hbox{in $ L^{\frac{2N}{N+2s}}(\R^N)$};$$
%in particular, 
for any $\varphi \in H^s(\R^N)$ we have
$$\int_{\R^N} g(u_n) \varphi \to \int_{\R^N} g(u) \varphi. $$
%Notice that $H^s(\R^N) \not \subset L^{\frac{2N}{N+2s}}(\R^N)$.
\item Assume in addition \hyperref[(h2)]{\textnormal{(h2)}} and \hyperref[(h3)]{\textnormal{(h3)}}. Let $u_n \wto u$ in $H^s_r(\R^N)$. Then 
$$\int_{\R^N} \abs{G(u_n) - G(u)} \to 0, \quad \int_{\R^N}\abs{g(u_n) u_n - g(u) u} \to 0$$
as well as $ \int_{\R^N}\abs{g(u_n) v - g(u) v} \to 0$ for each $v\in H^s(\R^N)$.
\item Assume in addition \hyperref[(h2)]{\textnormal{(h2)}} and \hyperref[(h3)]{\textnormal{(h3)}}. Let $u_n \wto u$ in $H^s(\Omega)$ with $\Omega \subset \R^N$ bounded. Then 
$$\int_{\Omega} \abs{G(u_n) - G(u)} \to 0, \quad \int_{\Omega}\abs{g(u_n) u_n - g(u) u} \to 0$$
as well as $ \int_{\R^N}\abs{g(u_n) v - g(u) v} \to 0$ for each $v\in H^s(\Omega)$.
%\item %(o serve (g3')?). 
%Assume in addtion (g2) and (g3).
%Let $u_n \wto u$ in $H^s(\R^N)$. Then, for each $v \in C^{\infty}_c(\R^N)$ (or more generally $v\in H^s(\R^N)$ with compact support) \tr{vale per generico $H^s(\R^N)$? Vedi meglio (serve $\sim$ Nemitskii continuo weak-weak)},
%%$$\int_{\R^N} \abs{g(u_n) v - g(u) v} \to 0.$$
%%$$\int_{\R^N} \big(g(u_n) v - g(u) v\big) \to 0.$$
%$$\int_{\R^N} g(u_n) v \to \int_{\R^N} g(u) v.$$
\end{itemize}
\end{Proposition}

\medskip

\claim Proof.
We prove the first claim. 
Let $\varphi \in C^{\infty}_c(\R^N)$, and
%we have (by the $L^p$-dominated convergence theorem) $f(u_n) \varphi \to f(u) \varphi$ in $L^{\frac{2N}{N+\alpha}}(\R^N)$.
let $\Omega:= \supp(\varphi)$. Since $u_n \to u$ in $L^r(\Omega)$ for each $r \in [2, 2^*_s)$, we have (by the $L^r$-dominated convergence theorem) $g(u_n) \to g(u)$ in $L^r(\Omega)$ for each $r \in [1, \frac{2N}{N+2s})$. For a whatever of such $r$, let $q$ be its conjugate; since $\varphi \in L^q(\R^N)$ for such $q$, we have $g(u_n) \varphi \to g(u) \varphi$ in $L^{1}(\Omega)$. Thus
$$ \int_{\R^N} g(u_n) \varphi \to \int_{\R^N} g(u) \varphi \quad \forall \varphi \in C^{\infty}_c(\R^N) .$$
We want to extend the relation to $H^s(\R^N)$. 
Indeed, observe first that, for $\varphi \in H^s(\R^N)$, %we easily see that $g(u_n)$ is bounded in $L^{\frac{2N}{N+2s}}(\R^N)$. Thus, if 
\begin{align*}
\pabs{ \int_{\R^N} g(u_n) \varphi} &\lesssim \int_{\R^N} \big( |u_n| + |u_n|^{2^*_s-1}\big) \pabs{\varphi} \\
&\leq \norm{u_n}_2 \norm{\varphi}_2 + \norm{u_n}_{2^*_s}^{\frac{N+2s}{N-2s}} \norm{\varphi}_{2^*_s} 
\leq C \norm{\varphi}_{H^s}
\end{align*}
uniform in $n \in \N$, since $u_n$ are equibounded in $L^2(\R^N) \cap L^{2^*_s}(\R^N)$. Let now $\varphi_{\eps} \in C^{\infty}_c(\R^N)$ approximating a fixed $\varphi$ in $H^s(\R^N)$. Then
%$$\int_{\R^N} g(u_n) \varphi - \int_{\R^N} g \varphi= \int_{\R^N} g(u_n) (\varphi-\varphi_{\eps}) + \int_{\R^N} (g(u_n) -g)\varphi_{\eps} + \int_{\R^N} g (\varphi_{\eps}-\varphi)$$
$$\int_{\R^N} g(u_n) \varphi - \int_{\R^N} g(u) \varphi= \int_{\R^N} g(u_n) (\varphi-\varphi_{\eps}) + \int_{\R^N} (g(u_n) -g(u))\varphi_{\eps} + \int_{\R^N} g(u) (\varphi_{\eps}-\varphi);$$
thus the first and the third quantities are small in $\eps$ (uniformly in $n$), and the second is small for $n=n(\eps)\gg 0$. Hence we have the first claim.

\smallskip

The second and the third claims are a consequence of \cite[Lemma 2.4]{ChWa0}. 
We exhibit here an easier proof of the second point, by assuming the stronger condition \hyperref[(h3')]{\textnormal{(h3')}}.

Recall that $H^s_r(\R^N)$ is compactly embedded in $L^p(\R^N)$, being $p\in (2, 2^*_s)$ introduced in \textnormal{(h3')}. Then by standard argument one has, up to a subsequence, that 
\begin{itemize}
\item $u_n \to u$ almost everywhere,
\item $u_n\to u$ strongly in $L^{p}(\R^N)$, with $|u_n|, |u| \leq w \in L^{p}(\R^N)$.
%\item $u_n\to u$ strongly in $L^2_{loc}$, with $|u_n|, |u| \leq w \in L^2_{loc}$. ((non serve))
\end{itemize}
By the assumption %\footnote{Indeed, let $\delta>0$. By $(g2)$, for $|x|\leq \omega$ small enough we have $|g(x)|\leq \delta |x|$; by $(g3)$, for $|x|\geq M$ big enough we have $|g(x)| \leq \delta |x|^p$. On the other hand, by $(g1)$ on the interval $[\omega, M]$ $g$ is bounded by a constant $C$, i.e. $|g(x)|\leq C = \frac{C}{\omega}\omega\leq \frac{C}{\omega}|x|$. By setting $C_{\delta}:= \max\{\delta,\frac{C}{\omega}\}$, we have the claim.},
 there exists an $M$ such that
$$|g(t)t|\leq \parag{C_{\delta} |t|^2 \quad \hbox{if $|t|\leq M$}, \\ \delta |t|^{p} % |t|^{p+1} 
\quad \hbox{if $|t|\geq M$}.}$$
Fixed a whatever $R>0$, set
$$M_n:=\{|u_n|\leq M\}\cap B_R(0),$$
we have
%\begin{align*}
%|g(u_n)u_n| &= |g(u_n)u_n| \chi_{M_n} + |g(u_n)u_n| \chi_{\R^N \setminus M_n} \\
%&\leq C_{\delta} |u_n|^2 \chi_{M_n} + \delta |u_n|^{p+1} \chi_{\R^N\setminus M_n} \\
%%&\leq& C_{\delta} |u_n|^2 \chi_{B_T(0)} + \delta |u_n|^{p+1} \\
%%&\leq& C_{\delta}|w|^2 \chi_{B_T(0)} + \delta |w|^{p+1} \in L^1
%&\leq C_{\delta} M^2 \chi_{M_n} + \delta |u_n|^{p+1} \\
%&\leq C_{\delta} M^2 \chi_{B_R(0)} + \delta |w|^{p+1} \in L^1(\R^N)
%\end{align*}
\begin{align*}
|g(u_n)u_n| &= |g(u_n)u_n| \chi_{M_n} + |g(u_n)u_n| \chi_{\R^N \setminus M_n}
\leq C_{\delta} |u_n|^2 \chi_{M_n} + \delta |u_n|^{p}% |u_n|^{p+1}
 \chi_{\R^N\setminus M_n} \\
%&\leq& C_{\delta} |u_n|^2 \chi_{B_T(0)} + \delta |u_n|^{p+1} \\
%&\leq& C_{\delta}|w|^2 \chi_{B_T(0)} + \delta |w|^{p+1} \in L^1
&\leq C_{\delta} M^2 \chi_{M_n} + \delta |u_n|^{p}  %|u_n|^{p+1}
\leq C_{\delta} M^2 \chi_{B_R(0)} + \delta |w|^{p}  %|w|^{p+1}
\in L^1(\R^N)
\end{align*}
and similarly for $G(u_n)$ 
%((al terzo/quarto rigo, si può più facilmente stimare: $ |u_n|^2 \chi_{M_n} \leq M^2 \chi_{M_n} \leq M^2 \chi_{B_T(0)} \in L^1$)).
%and similarly for 
and $|g(u_n)v|$. 
Moreover, since $g$ is continuous, we have $g(u_n)\to g(u)$ almost everywhere. By dominated convergence theorem, we obtain the claim.
%$$\int g(u_n)u_n \to \int g(u)u \quad \textnormal{ and } \quad \int g(u_n)u \to \int g(u)u.$$ 
\QED

\bigskip

%\section{Assumptions on the nonlinearity}
\subsection{Nonlocal nonlinearities}

%Let $s\in (0,1]$. 
For nonlocal nonlinearities of the type $\big(I_{\alpha}*F(u)\big) f(u)$, % with $\alpha \in (0,N)$, 
$F(t)= \int_0^t f(\tau) d\tau$, we introduce the following notations:
\begin{itemize}
\item \emph{Lower critical exponent:} $2^{\#}_{\alpha}:= 1 + \frac{\alpha}{N}=\frac{N+\alpha}{N} \in (1,2)$,
\item \emph{Upper critical exponent:} $2^{*}_{\alpha,s}:= 1 + \frac{\alpha+2s}{N-2s}= \frac{N+\alpha}{N-2s} \in (1, +\infty)$,
\item \emph{$L^2$-critical exponent:} $2^{m}_{\alpha,s}:= 1 + \frac{\alpha+2s}{N}= \frac{N+\alpha+2s}{N} \in (1, 2 + \frac{2}{N})$, %\min\{2 + \frac{2}{N}, 2^*_{\alpha,s}\})$,
\end{itemize}
and notice that
$$1<2^{\#}_{\alpha} < 2^{m}_{\alpha,s} < 2^{*}_{\alpha,s} < +\infty;$$
if $s=1$, if there is no ambiguity from the framework, we write $ 2^{*}_{\alpha} \equiv 2^{*}_{\alpha,1} = 1 + \frac{\alpha+2s}{N-2}= \frac{N+\alpha}{N-2}$ and $2^{m}_{\alpha} \equiv 2^{m}_{\alpha,1} = 1 + \frac{\alpha+2}{N}= \frac{N+\alpha+2}{N}$.

\begin{Remark}
We observe that, defining the Riesz potential by $x \mapsto \frac{A_{N,\beta}}{|x|^{\beta}}$, as some authors do, we have that the critical exponents %$\frac{N+\alpha}{N}<\frac{N+\alpha+2s}{N}<\frac{N+\alpha}{N-2s}$ 
become %respectively 
$\frac{2N-\beta}{N}<\frac{2N-\beta+2s}{N}<\frac{2N-\beta}{N-2s}$.
\end{Remark}

We introduce the following set of assumptions: % (valid for $s\in (0,1]$):
\begin{itemize}
\item[(H0)] \label{(H0)}
\emph{Continuity:} $f \in C(\R)$ (i.e. $F \in C^1(\R)$),
%\item Regularity (for Pohozaev): $s \in (\frac{1}{2},1)$ or $g \in C^{\sigma}_{loc}(\R)$ for some $\sigma > 1-2s$,
\item[(H0')] \label{(H0')}
\emph{Additional regularity:} $f \in C^{\sigma}_{loc}(\R)$ (i.e. $F\in C^{1,\sigma}_{loc}(\R)$) for some $\sigma \in (0,1]$,
\item[(H1)] \label{(H1)}
\emph{Nontriviality:} $F\nequiv 0$, i.e. there exists $t_0\in\R^*$ such that $F(t_0)\neq 0$,

\smallskip

\item[(H2)] \label{(H2)}
\emph{Well posedness:} $\limsup_{t \to 0} \frac{|f(t)|}{|t|^{2^{\#}_{\alpha}-1}} <\infty$, $\limsup_{|t| \to +\infty} \frac{|f(t)|}{|t|^{2^*_{\alpha,s}-1}}<\infty$, or equivalently $|t f(t)| \leq C\big( |t|^{2^{\#}_{\alpha}} + |t|^{2^*_{\alpha,s}}\big)$ for some $C<0$,
\item[(H2')] \label{(H2')}
\emph{$L^2$-well posedness:} $\limsup_{t \to 0} \frac{|f(t)|}{|t|^{2^{\#}_{\alpha}-1}} <\infty$, $\limsup_{|t| \to +\infty} \frac{|f(t)|}{|t|^{2^{m}_{\alpha,s}-1}}<\infty$, or equivalently 
$|t f(t)| \leq C\big( |t|^{2^{\#}_{\alpha}} + |t|^{2^{m}_{s,\alpha}}\big)$ for some $C<0$, 
\smallskip

\item[(H3)] \label{(H3)}
\emph{Supercriticality in $0$:} 
%$\lim_{t \to 0} \frac{F(t)}{|t|^{2^{\#}_{\alpha}}} =0$ \tr{serve $\lim_{t \to 0} \frac{f(t)}{|t|^{2^{\#}_{\alpha}-1}} =0$?},
%$\lim_{t \to 0} \frac{\tr{f}(t)}{|t|^{2^{\#}_{\alpha}}} =0$,
$\lim_{t \to 0} \frac{F(t)}{|t|^{2^{\#}_{\alpha}}} =0$,
\item[(H3')] \label{(H3')}
\emph{(Super)linerarity in $0$:} $\limsup_{t \to 0} \frac{|f(t)|}{|t|} <\infty$,
%$\limsup_{t \to 0} \frac{|F(t)|}{|t|^2} <\infty$,

\item[(H3*)] \label{(H3*)}
\emph{$L^2$-subcriticality in $0$:} $\lim_{t \to 0} \frac{|F(t)|}{|t|^{2^{m}_{\alpha,s}}} =+\infty$, 
\item[(H3*')] \label{(H3*')}
\emph{Sublinearity in $0$:} $\lim_{t \to 0} \frac{|f(t)|}{|t|} = +\infty$,
%$\lim_{t \to 0} \frac{|F(t)|}{|t|^2} = +\infty$,

\smallskip

\item[(H4)] \label{(H4)}
\emph{Subcriticality at $\infty$:} 
%$\lim_{|t| \to +\infty} \frac{F(t)}{|t|^{2^*_{\alpha,s}}} =0$ \tr{serve $\lim_{|t| \to +\infty} \frac{f(t)}{|t|^{2^*_{\alpha,s}-1}} =0$?},
%$\lim_{|t| \to +\infty} \frac{\tr{f}(t)}{|t|^{2^*_{t\alpha,s}}} =0$, %\tr{serve $\lim_{|t| \to +\infty} \frac{f(t)}{|t|^{2^*_{\alpha,s}-1}} =0$?},
$\lim_{|t| \to +\infty} \frac{F(t)}{|t|^{2^*_{\alpha,s}}} =0$, %\tr{serve $\lim_{|t| \to +\infty} \frac{f(t)}{|t|^{2^*_{\alpha,s}-1}} =0$?},
%\item Strict subcriticality at $\infty$: $\lim_{|t| \to +\infty} \frac{g(t)}{t^{q}} =0$ for some $q \in (1, 2^*_s-1)$,
\item[(H4')] \label{(H4')}
\emph{$L^2$-subcriticality at $\infty$:} $\lim_{|t| \to +\infty} \frac{F(t)}{|t|^{2^{m}_{\alpha,s}}} =0$, 
%\tr{basta questo, purché combinato con (f2')}
%\item Criticality at $\infty$: $\infty$: $\lim_{|t| \to +\infty} \frac{g(t)}{t^{2^*_s-1}} =a>0$, together with
%$f(t) \geq a t^{2^*_s-1} + Ct^{p-1}$
%with $C>0$ and $p \in (\max\{2^*_a-2s,2\}, 2^*_s)$,

\smallskip

\item[(H5)] \label{(H5)}
\emph{Symmetry:} $f$ is odd or even,
\item[(H6)] \label{(H6)}
\emph{Sign:} $f$ has constant sign on $(0,+\infty)$.
\end{itemize}

Notice that
$$ 
\hyperref[(H0')]{\textnormal{(H0')}}
\implies
\hyperref[(H0)]{\textnormal{(H0)}}, 
\quad 
\hyperref[(H2')]{\textnormal{(H2')}} 
\implies 
\hyperref[(H2)]{\textnormal{(H2)}}, 
\quad 
\hyperref[(H3')]{\textnormal{(H3')}} 
\implies 
\hyperref[(H3)]{\textnormal{(H3)}}, 
\quad 
%(f4') \implies (f4), 
\hyperref[(H2')]{\textnormal{(H2')}}
\implies 
\hyperref[(H4)]{\textnormal{(H4)}}, 
$$
$$
\hyperref[(H3*)]{\textnormal{(H3*)}}
 \vee 
 \hyperref[(H3*')]{\textnormal{(H3*')}}
 \implies \hyperref[(H1)]{\textnormal{(H1)}},
 \quad 
\hyperref[(H3)]{\textnormal{(H3)}}
\wedge 
\hyperref[(H4)]{\textnormal{(H4)}}
\implies 
\hyperref[(H2)]{\textnormal{(H2)}}, % \quad 
%(f3)\wedge (f4')\implies (f2'), 
$$
while generally \hyperref[(H3*)]{\textnormal{(H3*)}} and \hyperref[(H3*')]{\textnormal{(H3*')}} are not related (since $2$ and $2^m_{\alpha, s}$ are not so).

When searching for multiple normalized solutions in Choquard equations, in addition to \hyperref[(H3*)]{\textnormal{(H3*)}} and \hyperref[(H5)]{\textnormal{(H5)}} we will ask the following technical assumption (see also Remark \ref{rem_extra_cond}):
\begin{itemize}
\item[(H7)] \label{(H7)}
\emph{Almost monotonicity:} if $F$ is odd, then $F$ has a constant sign in $(0,\delta_0]$ and 
$$\sup_{t \in (0,\delta_0], \, h \in [0,1]} \pabs{\frac{F(th)}{F(t)}} < \infty$$ 
for some $\delta_0>0$ (e.g., $|F|$ is non-decreasing in $[0,\delta_0]$).
\end{itemize}

\begin{Remark}\label{rem_buona_posit_f}
%Since of key importance in the good definition of the functionals, as well as in bootstrap argument in the rest of the thesis, w
We write here in which spaces lie the considered quantities. 
Let $u\in H^s(\R^N) \subset L^2(\R^N) \cap L^{2^*_s}(\R^N)$. By \hyperref[(H2)]{\textnormal{(H2)}} we have (see Remark \ref{rem_somma_lebesg})
\begin{align*}
f(u) &\in L^{\frac{2N}{\alpha}}(\R^N)\cap L^{\frac{N}{\alpha} \frac{2N}{N-2s}}(\R^N)+ L^{2 \frac{N-2s}{\alpha+2s}}\cap L^{\frac{2N}{\alpha+2s}}(\R^N) \\
&\subset L^{\frac{2N}{\alpha}}(\R^N) + L^{\frac{2N}{\alpha+2s}}(\R^N),\\
F(u), f(u)u &\in L^{\frac{2N}{N+\alpha}}(\R^N)\cap L^{\frac{N}{N+\alpha} \frac{2N}{N-2s}}(\R^N) + L^{2 \frac{N-2s}{N+\alpha}}(\R^N)\cap L^{\frac{2N}{N+\alpha}}(\R^N) \\
&\subset L^{\frac{2N}{N+\alpha}}(\R^N).
\end{align*}
%Thus, seen as an operator, $F: L^2(\R^N) \cap L^{2^*_s}(\R^N) \to L^{\frac{2N}{N+\alpha}}(\R^N)$..
Thus by the Hardy-Littlewood-Sobolev inequality we obtain
\begin{align*}
I_{\alpha}*F(u) &\in L^{\frac{2N}{N-\alpha}}(\R^N) \cap L^{\frac{2N^2}{N^2-(\alpha+2s)N-2s\alpha}}(\R^N) + L^{\frac{2N(N-2s)}{N^2-\alpha N+4s \alpha}}(\R^N) \cap L^{\frac{2N}{N-\alpha}}(\R^N) \\
&\subset L^{\frac{2N}{N-\alpha}}(\R^N).
\end{align*}
Finally, by the H\"older inequality, we have
$$(I_{\alpha}*F(u))F(u) \in L^1(\R^N)$$
and
\begin{align*}
(I_{\alpha}*F(u))f(u) &\in L^2(\R^N)\cap L^{\frac{2N^2}{N^2 - 2s\alpha}}(\R^N) + L^{\frac{2N(N-2s)}{N^2+2 \alpha s}}(\R^N)\cap L^{\frac{2N}{N+2s}}(\R^N)\\
& \subset L^2(\R^N) + L^{\frac{2N}{N+2s}}(\R^N);
\end{align*}
we observe that $(I_{\alpha}*F(u))f(u)$ does not lie in $L^2(\R^N)$, generally. On the other hand, if $\varphi \in H^s(\R^N) \subset L^2(\R^N) \cap L^{2^*_s}(\R^N)$ is a test function, we notice that the found summability of $(I_{\alpha}*F(u))f(u)$ is enough to have
$\int_{\R^N} (I_{\alpha}*F(u))f(u) \varphi \, dx$ 
well defined, since $f(u) \varphi \in L^{\frac{2N}{N+\alpha}}(\R^N)$. 
\end{Remark}

\begin{Remark}\label{rem_conv_welldef}
By Propsition \ref{prop_HLS} %Arguing as in Proposition \ref{prop_conv_C0} 
we see that $I_{\alpha}*F(u) \in C_0(\R^N)$ (and thus it is well defined pointwise) if $F(u)$ lies in $L^{\frac{N}{\alpha}-\eps}(\R^N) \cap L^{\frac{N}{\alpha}+\eps}(\R^N)$ for some $\eps>0$. 
In particular, if $u\in L^1(\R^N) \cap L^{\infty}(\R^N)$, it is sufficient to assume that $F$ grows at most polynomially (and at least superlinearly) in zero and at infinity.
Moreover, assuming %\hyperref[(f1)]{\textnormal{(f1)}}-\hyperref[(f2)]{\textnormal{(f2)}} 
\hyperref[(H0)]{\textnormal{(H0)}} and \hyperref[(H2)]{\textnormal{(H2)}} on $f$, we need to assume that $u \in L^{\frac{N+\alpha}{\alpha} - \eps}(\R^N) \cap L^{\frac{N}{\alpha}\frac{N+\alpha}{N-2s} + \eps}(\R^N)$ for some $\eps>0$; in particular, the convolution is pointwise well defined if $u \in L^2(\R^N) \cap L^{\frac{N}{\alpha} \frac{2N}{N-2s}}(\R^N)$.
\end{Remark}

We state now the convergences for the nonlinear Choquard terms in the case of a subcritical growth (see also \cite[pages 6565 and 6577]{MS2}, \cite[page 11]{BaLiLi} and \cite[page 353]{Amb5}).
\begin{Proposition}
\label{prop_converg_generiche_nonloc}
 Assume \hyperref[(H0)]{\textnormal{(H0)}} and \hyperref[(H2)]{\textnormal{(H2)}}.
\begin{itemize}
\item Let $u_n \wto u$ in $H^s(\R^N)$. Then 
%$$\big(I_{\alpha}*F(u_n)\big) F(u_n) \wto \big(I_{\alpha}*F(u)\big) F(u) \quad \hbox{in $ L^{\frac{2N}{N+2s}}(\R^N)$};$$
%in particular, 
for any $\varphi \in H^s(\R^N)$ we have
$$\int_{\R^N} \big(I_{\alpha}*F(u_n)\big) f(u_n) \varphi \to \int_{\R^N} \big(I_{\alpha}*F(u)\big) f(u) \varphi. $$
%Notice that $H^s(\R^N) \not \subset L^{\frac{2N}{N+2s}}(\R^N)$.
\item Assume in addition \hyperref[(H3)]{\textnormal{(H3)}} and \hyperref[(H4)]{\textnormal{(H4)}}. % (actually we need (f4')). 
Let $u_n \wto u$ in $H^s_r(\R^N)$. Then 
$$\int_{\R^N} \big(I_{\alpha}*F(u_n)\big) F(u_n) \to \int_{\R^N} \big(I_{\alpha}*F(u)\big) F(u) $$
and %\tr{IMPORTANTE: qui serve la condizione di sottocriticalità in $f$! Non basta quella su $F$ (?)}
$$\int_{\R^N} \big(I_{\alpha}*F(u_n)\big) f(u_n) u_n \to \int_{\R^N} \big(I_{\alpha}*F(u)\big) f(u) u. $$
%as well as $\int_{\R^N} \big(I_{\alpha}*F(u_n)\big) f(u_n) v \to \int_{\R^N} \big(I_{\alpha}*F(u)\big) f(u) v$ for each $v\in H^s(\R^N)$.
\end{itemize}
%\begin{itemize}
%\item Let $u_n \wto u$ in $H^s_r(\R^N)$. Then 
%$$\int_{\R^N} \abs{G(u_n) - G(u)} \to 0, \quad \int_{\R^N}\abs{g(u_n) u_n - g(u) u} \to 0$$
%as well as $ \int_{\R^N}\abs{g(u_n) v - g(u) v} \to 0$ for each $v\in H^s(\R^N)$.
%\item Let $u_n \wto u$ in $H^s(\Omega)$ with $\Omega \subset \R^N$ bounded. Then 
%$$\int_{\Omega} \abs{G(u_n) - G(u)} \to 0, \quad \int_{\Omega}\abs{g(u_n) u_n - g(u) u} \to 0$$
%as well as $ \int_{\R^N}\abs{g(u_n) v - g(u) v} \to 0$ for each $v\in H^s(\Omega)$.
%\item %(o serve (g3')?). 
%Let $u_n \wto u$ in $H^s(\R^N)$. Then, for each $v \in C^{\infty}_c(\R^N)$ (or more generally $v\in H^s(\R^N)$ with compact support) \tr{vale per generico $H^s(\R^N)$? Vedi meglio (serve $\sim$ Nemitskii continuo weak-weak)},
%%$$\int_{\R^N} \abs{g(u_n) v - g(u) v} \to 0.$$
%%$$\int_{\R^N} \big(g(u_n) v - g(u) v\big) \to 0.$$
%$$\int_{\R^N} g(u_n) v \to \int_{\R^N} g(u) v.$$
%\end{itemize}
\end{Proposition}

\claim Proof.
%Assume \hyperref[(H2)]{\textnormal{(H2)}}. 
Let $u_n \wto u$ in $H^s(\R^N)$, then $u_n$ is bounded in $L^2(\R^N) \cap L^{2^*_s}(\R^N)$. % and (up to a subsequence) $u_n \to u$ a.$\,$e. pointwise.
 By Remark \ref{rem_buona_posit_f} we have $F(u_n)$ bounded in $L^{\frac{2N}{N+\alpha}}(\R^N)$. %, with $\frac{2N}{N+\alpha}>1$.
 %, thus (up to a subsequence) $F(u_n) \wto v$ in $L^{\frac{2N}{N+\alpha}}(\R^N)$. 
%moreover $F(u_n) \to F(u)$ a.$\,$e. pointwise. By \cite[Lemma 1]{MT1} we have $v=F(u)$, thus $F(u_n) \wto F(u)$ in $L^{\frac{2N}{N+\alpha}}(\R^N)$.
Moreover we can assume $u_n \to u$ in $L^p_{loc}(\R^N)$ for $p \in [1, 2^*_s)$, and thus $F(u_n) \to F(u)$ in $L^q_{loc}(\R^N)$ for $q \in [1, \frac{2N}{N+\alpha})$. This two information on $F(u_n)$ imply $F(u_n) \wto F(u)$ in $L^{\frac{2N}{N+\alpha}}(\R^N)$%
\footnote{We argue in this way. First, fix $\Omega \subset \R^N$ bounded and $q \in [1, \frac{2N}{N+\alpha})$, so that $L^{\frac{2N}{N+\alpha}}(\Omega) \subset L^q(\Omega)$ for every $q \in [1, \frac{2N}{N+\alpha})$. Since $F(u_n)$ is bounded in $L^{\frac{2N}{N+\alpha}}(\R^N) \subset L^{\frac{2N}{N+\alpha}}(\Omega)$ and $\frac{2N}{N+\alpha}>1$, it converges to some $v \in L^{\frac{2N}{N+\alpha}}(\Omega) \subset L^q(\Omega)$; on the other hand $F(u_n) \wto F(u) \in L^q(\R^N)\subset L^q(\Omega)$, thus by uniqueness $v=F(u_n)$. Let now $\varphi$ be in the dual $L^{\frac{N-\alpha}{2N}}(\R^N)$, and consider $\varphi_k \in C^{\infty}_c(\R^N)$ approximating $\varphi$. Thus
$$\int_{\R^N} \big(F(u_n)-F(u)\big) \varphi \leq \norm{F(u_n)-F(u)}_{\frac{2N}{N+\alpha}} \norm{\varphi-\varphi_k}_{\frac{N-\alpha}{2N}} + \int_{\supp(\varphi_k)} \big( F(u_n)-F(u)\big) \varphi_k;$$
exploiting that $F(u_n)$ is bounded, the first piece is small for $k$ large (uniform in $n$), while the second is small (fixed this $k$), for $n$ large, by the previous argument with $\Omega=\supp(\varphi_k)$.}%
.%
\footnote{We can deduce the implication also in this way: since $u_n \to u$ a.$\,$e. pointwise and $F$ is continuous, then $F(u_n) \to F(u)$ a.$\,$e. pointwise; moreover, being bounded, then $F(u_n) \wto v$ in $L^{\frac{2N}{N+\alpha}}(\R^N)$ for some $v$, where $\frac{2N}{N+\alpha}>1$; hence by \cite[Lemma 1]{MT1} we have $v=F(u)$.}
%Moreover we can $F(u_n) \to F(u)$ a.$\,$e. pointwise. By \cite[Lemma 1]{MT1} we have $v=F(u)$, thus $F(u_n) \wto F(u)$ in $L^{\frac{2N}{N+\alpha}}(\R^N)$. 
%%
%Let $u_n \wto u$ in $H^s(\R^N)$, then (up to a subsequence) \tr{$u_n\to u$ in $L^p_{loc}(\R^N)$ for $p \in [1, 2^*_s)$}, and a.$\,$e. pointwise.
%By assumptions $F(u_n)$ is bounded in $L^{\frac{2N}{N+\alpha}}(\R^N)$ and \tr{$F(u_n) \to F(u)$ in $L^q_{loc}(\R^N)$ for $q \in [1, \frac{2N}{N+\alpha})$.} 
%\tr{Thus $F(u_n) \wto F(u)$ in $L^{\frac{2N}{N+\alpha}}(\R^N)$.}
%
By some standard topological argument, the convergence holds for the whole sequence. 
Moreover, by Proposition \ref{prop_HLS} we gain 
$$I_{\alpha}*F(u_n) \wto I_{\alpha}*F(u) \quad \hbox{in $L^{\frac{2N}{N-\alpha}}(\R^N)$}.$$

%In a similar way, we obtain, for any $\varphi \in H^s(\R^N)$, $f(u_n) \varphi \wto f(u) \varphi$ in $L^{\frac{2N}{N+\alpha}}(\R^N)$ and a.$\,$e. pointwise. 
%In addition, if 
Let now $\varphi \in C^{\infty}_c(\R^N)$, 
%we have (by the $L^p$-dominated convergence theorem) $f(u_n) \varphi \to f(u) \varphi$ in $L^{\frac{2N}{N+\alpha}}(\R^N)$.
and set $\Omega:= \supp(\varphi)$. Since $u_n \to u$ in $L^p(\Omega)$ for each $p \in [2, 2^*_s)$, we have (by the $L^p$-dominated convergence theorem) $f(u_n) \to f(u)$ in $L^p(\Omega)$ for each $p \in [1, \frac{2N}{\alpha+2s})$. Let $p \in (\frac{2N}{N+\alpha}, \frac{2N}{\alpha+2s})$ be whatever and let $q $ be such that $\frac{1}{p}+\frac{1}{q} = \frac{N+\alpha}{2N}$; since $\varphi \in L^q(\R^N)$ for such $q$, we have $f(u_n) \varphi \to f(u) \varphi$ in $L^{\frac{2N}{N+\alpha}}(\Omega)$. Thus
$$ \int_{\R^N} \big(I_{\alpha}*F(u_n)\big) f(u_n) \varphi \to \int_{\R^N} \big(I_{\alpha}*F(u)\big) f(u) \varphi \quad \forall \varphi \in C^{\infty}_c(\R^N) .$$
To extend the relation to $\varphi \in H^s(\R^N)$ we argue as in Proposition \ref{prop_converg_generiche_loc}, after having observed that
\begin{align*}
\pabs{ \int_{\R^N} \big(I_{\alpha}*F(u_n)\big) f(u_n) \varphi} &\lesssim \norm{F(u_n)}_{\frac{2N}{N+\alpha}} \norm{f(u_n) \varphi}_{\frac{2N}{N+\alpha}} \\
&\lesssim \norm{|u_n|^{2^{\#}_{\alpha}-1}\varphi}_{\frac{2N}{N+\alpha}} + \norm{|u_n|^{2^*_{\alpha,s}-1} \varphi}_{\frac{2N}{N+\alpha}} 
\\
&\leq \norm{|u_n|^{\frac{\alpha}{N}}}_{\frac{2N}{\alpha}} \norm{\varphi}_{2} + \norm{|u_n|^{\frac{\alpha+2s}{N-2s}}}_{\frac{2N}{\alpha+2s}} \norm{\varphi}_{2^*_s} 
\\
&\leq \norm{u_n}_2^{2^{\#}_{\alpha}-1} \norm{\varphi}_{H^s} + \norm{u_n}_{2^*_s}^{2^*_{\alpha,s}-1} \norm{\varphi}_{H^s} 
 \lesssim \norm{\varphi}_{H^s}.
\end{align*}
%. Indeed, we easily see that $g_n:=\big(I_{\alpha}*F(u_n)\big) f(u_n)$ is bounded in $L^{\frac{2N}{N+2s}}(\R^N)$. Thus, if $\varphi_{\eps} \in C^{\infty}_c(\R^N)$ approximate $\varphi$, we have
%$$\int_{\R^N} g_n \varphi - \int_{\R^N} g= \int_{\R^N}g_n (\varphi-\varphi_{\eps}) + \int_{\R^N} (g_n-g)\varphi_{\eps} + \int_{\R^N} g (\varphi_{\eps}-\varphi)$$
%where the first and the third quantities are small in $\eps$ (uniformly in $n$), and the second is small for $n=n(\eps)\gg 0$. Thus we have the first claim.

Assume now \hyperref[(H3)]{\textnormal{(H3)}} and \hyperref[(H4)]{\textnormal{(H4)}}. Let $G(t):= (F(t))^{\frac{N+\alpha}{2N}}$. By the assumptions we have
$$\lim_{t \to 0} \frac{G(t)}{|t|^2} = \lim_{t \to 0} \left( \frac{F(t)}{|t|^{2^{\#}_{\alpha}}}\right)^{\frac{N+\alpha}{2N}}=0, \quad \quad \lim_{t \to \infty} \frac{G(t)}{|t|^{2^*_s}} = \lim_{t \to 0} \left( \frac{F(t)}{|t|^{2^{*}_{\alpha, s}}}\right)^{\frac{N+\alpha}{2N}}=0.$$
Thus, by Proposition \ref{prop_converg_generiche_loc} we gain $G(u_n) \to G(u)$ in $L^1(\R^N)$, which means $F(u_n) \to F(u)$ in $L^{\frac{2N}{N+\alpha}}(\R^N)$. In particular, by Proposition \ref{prop_HLS} we obtain 
$$I_{\alpha}*F(u_n) \to I_{\alpha}*F(u) \quad \hbox{in $L^{\frac{2N}{N-\alpha}}(\R^N)$}.$$
Thus we get the first claim. Moreover, arguing as before we get $f(u_n)u_n \wto f(u) u$ in $L^{\frac{2N}{N+\alpha}}(\R^N)$, and this concludes the proof.
%\tr{Similarly we argue for $G(t):=f(t)t$. %(VEDI).
%}
%\\ \tr{Can we say that $f(u_n)u_n \wto f(u) u$ and exploit the product strong-weak?}
\QED

\medskip

\begin{Remark}
When $\alpha \to 0$, by \eqref{eq_alpha_to_0}, we know that, under suitable assumptions,
$$(I_{\alpha}*F(u)) f(u) \stackrel{\alpha \to 0} \to F(u) f(u)=:g(u);$$
notice that (by integration by parts) $G(u)=\frac{1}{2} F^2(u)$. 
This relation is coherent with the definitions of the critical exponents of the local and nonlocal frameworks; indeed:
%$$2^{\#}_0 + (2^{\#}_0-1) %= 2 
%= 2^{\#}-1 ,$$
%$$2^{*}_{0,s} + (2^{*}_{0,s}-1) %= 1+\frac{4s}{N-2s} 
%= 2^{*}_s-1, $$
%$$2^{m}_{0,s} + (2^{m}_{0,s}-1) %= 1+\frac{4s}{N} 
%= 2^{m}_s-1. $$
$$2^{\#}_0 + (2^{\#}_0-1) %= 2 
= 2^{\#}-1 ,\quad
2^{*}_{0,s} + (2^{*}_{0,s}-1) %= 1+\frac{4s}{N-2s} 
= 2^{*}_s-1, \quad
2^{m}_{0,s} + (2^{m}_{0,s}-1) %= 1+\frac{4s}{N} 
= 2^{m}_s-1. $$
This correspondence lacks when comparing the nontriviality assumptions $F(t_0)\neq 0$ and $G(t_0) \geq \frac{\mu}{2} t_0^2$: this is due to the fact that, for any $\alpha \neq 0$, the pieces $\mu u$ and $(I_{\alpha}*F(u))f(u)$ scales differently.
Moreover, we see that while the subcriticality assumptions for the local problem are made for $g$, for the nonlocal problem are made for $F$, since essentially the product $Ff$ automatically becomes subcritical if $F$ is so.
\end{Remark}

%%%%%%%%%%%%%%%%%%%%%%%%%%%%%%%%%%%%%%%%%%%%%%%%%%%%%%%
%%%%%%%%%%%%%%%%%%%%%%%%%%%%%%%%%%%%%%%%%%%%%%%%%%%%%%%
%%%%%%%%%%%%%%%%%%%%%%%%%%%%%%%%%%%%%%%%%%%%%%%%%%%%%%%

%\chapter{Fractional Schr\"odinger equations: prescribed and free mass problems} %prescribed frequency problems} 
%VERSIONE STAMPATA (NON PIU')
%\chapter{\textls[-10]{Fractional Schr\"odinger equations: prescribed and free mass problems}} 
\chapter{Fractional Schr\"odinger equations: prescribed and free mass problems} 

%Normalized solutions for fractional equations with general nonlinearities}
\label{chap_fract_normal}
%Existence and multiplicity: nonlocal operators and nonlocal sources} %: prescribed mass solutions} 

%In this Chapter we aim to study existence and multiplicity of solutions, focusing especially on the case of normalized ($L^2$-constrained) solutions.
%In the first Section we will deal with a fractional problem (with local nonlinearities): here some new ideas (about the compactness and the geometry) will be introduced, based on a Lagragian approach.
%In the second Section we will treat a Choquard problem (with local Laplacian): here we implement the previous ideas, but the main focus will be made on existence of multiple solutions, since when the nonlinearity is nonlocal the construction of a multidimensional path is not obvious.

%\section{Normalized solutions for fractional equations with general nonlinearities}

%\bigskip

In this Chapter we study the following fractional Schr\"odinger equation 
$$	(-\Delta)^{s} u + \mu u =g(u) \quad \text{in $\mathbb{R}^N$},$$
where $N\geq 2$, $s\in (0,1)$, $u \in H^s(\mathbb{R}^N)$, $\mu>0$ is a frequency and $g \in C(\mathbb{R}, \mathbb{R})$ satisfies Berestycki-Lions type conditions. First, we recall some known facts about the \emph{unconstrained problem}, i.e. when $\mu$ is fixed, which has been investigated in \cite{BKS,ChWa0}. 
Then we study the \emph{constrained} problem
\begin{equation*}
%(P_m) \quad	
\parag{
	(-\Delta)^{s} u + \mu u &=g(u) & \; \text{in $\mathbb{R}^N$}, \cr
	\int_{\mathbb{R}^N} u^2 dx &= m, & \cr
	}
\end{equation*}
%where $N\geq 2$, $s\in (0,1)$, $m>0$, $u \in H^s_r(\mathbb{R}^N)$, $\mu$ is an unknown Lagrange multiplier and $g \in C(\mathbb{R}, \mathbb{R})$ satisfies Berestycki-Lions type conditions.
where $m>0$ is a prescribed mass, $u \in H^s_r(\mathbb{R}^N)$ and $\mu$ is a Lagrange multiplier, part of the unknowns.
Using a Lagrangian formulation, % of the problem $(P_m)$, 
we prove the existence of a weak solution with prescribed mass when $g$ has an $L^2$-subcritical growth. 
The approach relies on the construction of a minimax structure, by means of a \emph{Pohozaev mountain} in a product space and some deformation arguments under a weaker version of the Palais-Smale condition. 
A multiplicity result of infinitely many normalized solutions is also obtained if $g$ is odd, and this is new even for $g$ power.

\medskip

The present Chapter is mainly based on the paper \cite{CGT1} (see also \cite{CGT2}).

%%%%%%%%%%%%%%%%%%%%%%%%%%%%%%%%%%%%%%%%%%%%%%%%%%%%%%%
%\setcounter{equation}{0} %FOR ARXIV
\section{The fractional Schr\"odinger equation: a long-range interaction 
%\tb{da distribuire con l'intro iniziale}%COMMENT NOW
}
\label{sec_introd_frac_lap}

In 1948, following a suggestion by %P.A.M. 
Dirac, %R.P. 
Feynman \cite{Fey0} proposed a new suggestive description of the time evolution of the state of a non-relativistic quantum particle. 
According to Feynman, the wave function solution of the Schr\"odinger equation should be given by a heuristic integral over the space of paths: 
the classical notion of a single, unique classical trajectory for a system is replaced by a functional integral over an infinity of quantum-mechanically possible trajectories. 
Following Feynman's path integral approach to quantum mechanics, Laskin \cite{Las0,Las1,Las2,Las3} generalized the path integral over Brownian motions (random motion seen in swirling gas molecules) to L\'evy flights (a mix of long trajectories and short, random movements found in turbulent fluids) and derived the \emph{fractional nonlinear Schr\"odinger} ((fNLS) for short) equation 
\begin{equation}\label{eq_fNLS_1}
	i \hbar \partial_t \psi = \hbar^{2s} (-\Delta)^s \psi + V(x) \psi- g(\psi), \quad (t,x) \in (0,+\infty) \times \R^N
%i \partial_t \psi = (-\Delta)^s \psi - g(\psi), \quad (t,x) \in \R \times \R^N ,
\end{equation}
where %$\psi=\psi(t,x)$ is a complex wave, 
$s \in (0,1)$, $N>2s$, the symbol $(-\Delta)^{s}$ denotes the fractional power of the Laplace operator (defined via Fourier transform on the spatial variable), $\hbar$ designates the usual Planck constant, $V$ is a real potential and $g$ is a Gauge invariant nonlinearity, i.e. $g(e^{i \theta} \rho) = e^{i \theta} g(\rho)$ for any $\rho$, $\theta \in \R$. 
The complex wave function $\psi(x,t)$ represents the quantum mechanical probability amplitude for a given unit mass particle to have position $x$ at time $t$, under the confinement due to the potential $V,$ and $|\psi|^2$ is the corresponding probability density.

Fractional integrals and derivatives in the calculation methods have been used for the explanation of physical phenomena which do not comply with the laws of classical statistical physics, for instance in modeling Bose-Einstein condensates. 
It is known that Bose-Einstein condensation, theoretically discovered in 1924 and observed experimentally with alkali metals %(rubidium and sodium atoms) 
in 1995, represents a topical subject due to the explanation of quantum effects seen on a macroscopic scale, transmission of matter and the behaviour of superconductivity and superfluids. 
In this respect, not only experimental studies are important but theoretical studies too, which lead to the analysis of class of (fNLS) equations (also known as fractional Gross-Pitaievskii equations). 
Numerical simulations show existence of standing waves solutions, having a soliton behaviour and bound states \cite{DZ0,ZHSWW}, including mass conservation, energy conservation and dispersion relation, in which the fractional order exponent influences the shape of the state.

In 2015 a first optical realization of the fractional Schr\"odinger equation, based on transverse light dynamics in aspherical optical cavities, was achieved by Longhi \cite{Lon0}; 
subsequently, the propagation dynamics of wave packets were reported in Kerr nonlinearities, with constant or double-barrier potential. 
Numerical results showed the existence of solitons for (fNLS) equations where the L\'evy index $s$ and the saturation parameter can significantly affect the stability of these solitons \cite{KSM,WCCYZH,LiMaMi,YL0}.
Numerous other applications of the (fNLS) equation arise in the physical sciences, ranging from models of boson stars (see Section \ref{sec_boson_stars}) to geo-hydrology \cite{Ata0}, from charge transport in biopolymers, like DNA \cite{KLS} to anomalous diffusion phenomena \cite{BuV,Vaz1,MK1}, from water wave dynamics \cite{IP0} to jump processes in probability theory with applications to financial mathematics (see also \cite{DnPV} and the references therein). 
%[RIPETIZIONE]
Applications for wide ranges of $s$ appears, for example, also in the dynamics of populations \cite{CDV}: here small values of $s\approx 0$ or large values of $s\approx 1$ better model specific behaviours, according to the environments.
We refer also to %\cite{MK1,Vaz1} for a discussion on recent developments in the description of anomalous diffusion via fractional dynamics and to
 \cite{BJMR,BER} for some recent applications %of fractional operators to different frameworks (
to the analysis of the amount of bromsulphthalein in the human liver and to the study of thermostat systems, and others.

\medskip

From a mathematical view point, when searching for standing waves to \eqref{eq_fNLS_1}, i.e. factorized solutions 
$$\psi (t,x) = e^{i \mu t } u(x), \quad \mu >0,$$
 two possible directions can be pursued. 
A first possibility is to study \eqref{eq_fNLS_1} with a prescribed frequency $\mu$ and free mass. 
This approach, which we call the \emph{unconstrained} problem, has been deeply developed: 
the literature concerning the local version of the unconstrained problem starts from the seminal papers of Berestycki and Lions \cite{BL1,BL2} (see also \cite{JT0,Med1,BeCaNi,BN}) and it is so large that we do not even make an attempt to summarize it.
Some fundamental contributions for the fractional case $s\in (0,1)$ instead can be found in \cite{CafSil1,CabSir,FLS}; in particular, the existence and qualitative properties of the solutions for more general classes of fractional NLS equations with local source were studied in \cite{FQT,ChWa0,BKS,Iko1,Iko2,Amb5}.

A second approach is to prescribe the mass of $u$, thus conserved by $\psi$ in time
$$\int_{\R^N} |\psi(x,t)|^2 \, dx= m, \quad \forall \, t \in (0,+\infty)$$
and let the frequency $\mu$ to be free, becoming an unknown. 
This second approach is of considerable significance in physics, not only for the quantum probability normalization and the information on the mass itself, but also because the mass may also have specific meaning, such as the power supply in nonlinear optics, or the total number of atoms in Bose-Einstein condensation. 
Moreover, it can give better insights into the dynamical properties, such as the orbital stability or instability of solutions of \eqref{eq_fNLS_1} (see \cite{CaLi0}).

 In the local framework ($s=1$) the seminal contribution to the study of \emph{constrained} problems is due to Stuart \cite{Stu2}, Cazenave and Lions \cite{CaLi0}; 
see \cite{Jea0,BV0,BJS,Shi0,HT0,BM0,MeSc0,Sch0,BCDN,BCGJ,BFDST} for more recent contributions in the local case.

In the fractional case, the existence of a mass-constrained solution was, instead, recently considered in \cite{Fen0, Yan0, Din0} for pure powers and in \cite{LZ0} for combined powers. 
It remains an open problem anyway to derive analytically the existence of infinitely many bound states with higher energy, including mass conservation.

\medskip

The present Chapter is dedicated to the study of standing waves solutions of \eqref{eq_fNLS_1} (when $V=const$ and we fix $\hbar=1$) with prescribed mass, by means of a new variational method. 
Namely, we are interested to seek for radially symmetric solutions of the fractional problem
\begin{equation} \label{problem_frac}
%\tag{$P_m$} We need the package amsmath, which does not work well in Nonlinearity (as the journal itself says).
%(P_m) 
%\quad 
\parag{
&(-\Delta)^{s} u + \mu u = g(u) \quad \hbox{in $\R^N$,} & \cr
&\int_{\R^N} u^2 dx = m,& 
}
\end{equation}
where $N\geq 2$, $s\in (0,1)$, $m>0$ and $\mu$ is a Lagrange multiplier. 
We assume that the function $g$ satisfies the following Berestycki-Lions type conditions:
\vskip2pt
\begin{itemize}
	\item[(g1)] \label{(g1)}
$g : \R \to \R$ continuous and $\lim_{t \to 0} \frac{g(t)}{t}=0$,
	\item[(g2)] \label{(g2)}
$\lim_{|t| \to \infty} \frac{g(t)}{|t|^p} =0$ where $p = 2^m_s=1 + \frac{4s}{N}$ (see also Remark \ref{rem_weak_g2}),
	\item[(g3)] \label{(g3)}
there exists $t_0>0$ such that $G(t_0) >0$,
\end{itemize}
\vskip 2pt
where $G(t)= \int_0^t g(\tau) d\tau$. 
We recall that %the exponent $p = 1 + \frac{4s}{N}$ appears as a $L^2$-critical exponent for the nonlinear fractional equations with $L^2$-constraint and thus assumption 
 \hyperref[(g2)]{\textnormal{(g2)}} 
means that $g$ has an $L^2$-subcritical growth.

The solutions to \eqref{problem_frac} can be characterized as critical points of the $C^1$-functional $\mc{L} : H^s_r(\R^N) \to \R$
$$ \mathcal{L}(u) :=\half \int_{\R^N} |(-\Delta)^{s/2} u|^2 - \int_{\R^N} G(u) $$
constrained on the sphere 
$$ \mathcal{S}_m := \big\{ u \in H^s_r(\R^N) \mid \|u\|^2_2= m \big\}; $$
here we consider thus, as in \cite{HT0}, % (see also \cite{CT1}), 
a Lagrangian formulation of the problem \eqref{problem_frac}. 
In order to avoid technical issues with the boundary of $\R_+$ (see Section \ref{sec_ground_state} for a different approach), %For technical reasons 
we write 
$$\mu \equiv e^\lambda$$
with $\lambda \in \R$ and define the $C^1$-functional $\mc{I}^m : \R \times H^s_r(\R^N) \to \R$ by setting 
%\begin{equation}\label{functlag}
$$
\mc{I}^m(\lambda, u) :=\half \int_{\R^N} |(-\Delta)^{s/2} u|^2 - \int_{\R^N} G(u)\, + 
\frac{e^\lambda}{2} \bigl( \|u\|_2^2 -m \bigr).
$$%\end{equation}
We seek for critical points $(\lambda,u) \in \R \times H^s_r(\R^N)$ of 
$\mc{I}^m$, namely weak solutions of $\partial_u \mc{I}^m(\lambda, u)=0$ and $\partial_\lambda \mc{I}^m(\lambda, u)=0$ or equivalently 
$$
\quad \parag{
&\int_{\R^N} \left((-\Delta)^{s/2} u \ (-\Delta)^{s/2} \phi + e^\lambda u \phi \right)= \int_{\R^N} g(u) \phi, \quad \forall \phi \in H_r^s(\R^N), & \cr 
&\int_{\R^N} u^2 dx = m. &
}
$$
We implement a minimax approach to detect normalized solutions in the nonlocal framework using a Pohozaev type function.
More precisely, inspired by the Pohozaev (or Pohozaev-Derrick) identity \cite{Poh0} 
\begin{equation}\label{eq_Pohozaev_fractional}
\frac{N-2s}{2} \int_{\R^N} | (-\Delta)^{s/2}u|^2 + N \int_{\R^N} \left( \frac{\mu}{2} u^2 - G(u)\right) =0,
\end{equation}
for any $s \in (0,1)$ we introduce the Pohozaev function 
$\mc{P}:\R \times H^s_r(\R^N) \to \R$ by setting 
$$
\mc{P}(\lambda, u) :=
\frac{N-2s}{2} \int_{\R^N} |(-\Delta)^{s/2} u|^2 + N \int_{\R^N} \left( \frac{e^\lambda}{2} u^2 - G(u)\right) 
$$
and the Pohozaev set
$$ \Omega :=\big\{(\lambda,u) \in \R \times H^s_r(\R^N) \mid \mc{P}(\lambda,u)>0\big\} \cup\big\{(\lambda,0) \mid \lambda \in \R \big\}. $$
We note that, for each $\lambda\in\R$, the set
$\{u\in H_r^s(\R^N) \mid \mc{P}(\lambda,u)>0\}\cup\{ 0\}$ is a neighborhood of $u=0$, and thus 
$$ \partial \Omega =\big\{(\lambda,u) \in \R \times H^s_r(\R^N) \mid \mc{P}(\lambda,u)=0, \ u \neq 0 \big\}. $$
Therefore $(\lambda, u) \in \partial \Omega$ if and only if $u \neq 0$ and $u$ satisfies the Pohozaev identity. However we emphasize that under assumptions \hyperref[(g1)]{\textnormal{(g1)}}--\hyperref[(g3)]{\textnormal{(g3)}}, 
if $u \in H^s(\R^N)$ solves $\partial_u \mc{I}^m(\lambda, \cdot)=0$ with $\lambda \in \R$ fixed, then $\mc{P}(\lambda, u)=0$ when $s \in (\frac{1}{2},1)$. A similar result for $s \in (0,\frac{1}{2}]$ is not available since the weak solutions are not proved to be $C^1$, in general (see Section \ref{sec_frac_unconstr}). 

In spite of this lack of regularity, which is a special feature of the nonlocal framework, we recognize a Mountain Pass structure \cite{AR0} for the functional $\mc{I}^m$, where the mountain is given by the subset $\partial\Omega$. We refer to it as the \emph{Pohozaev mountain}. 
This approach can be useful to deal with different problems in other contexts.

Inspired by \cite{HT0,IT0}, we need to use a new variant of the Palais-Smale condition which takes into account the Pohozaev identity, and we establish some deformation theorems which enable us to perform our minimax arguments in the \emph{product space} $\R \times H^s_r(\R^N)$. 

As a byproduct, our solutions satisfy the Pohozaev identity, even if we assume that $f$ is only a continuous function (see Corollary \ref{coroll_esist_Pm}).
We also note that solutions with the Pohozaev identity are essential, in the following sense: our deformation argument shows that only critical points with the Pohozaev identity contribute to the topology; 
that is, solutions without the Pohozaev identity are deformable with a suitable deformation flow and have no topological relevance.

\medskip

Firstly we prove the following existence results for \eqref{problem_frac}. 

\begin{Theorem}\label{S:1.1_frac}
Suppose $N\geq 2$ and \hyperref[(g1)]{\textnormal{(g1)}}--\hyperref[(g3)]{\textnormal{(g3)}}.
Then there exists $m_0 \geq 0$ such that for any $m>m_0$, the problem \eqref{problem_frac} has a solution, satisfying the Pohozaev identity \eqref{eq_Pohozaev_fractional}.
\end{Theorem}	
	
\begin{Theorem}\label{S:1.12_frac}
		Suppose $N\geq 2$, \hyperref[(g1)]{\textnormal{(g1)}}--\hyperref[(g3)]{\textnormal{(g3)}} and 
		\begin{itemize}
%		\item[\textnormal{(g4)}] $\lim_{t \to 0} \frac{g(t)}{|t|^{\frac{4s}{N}} t} = + \infty.$
		\item[\textnormal{(g4)}] \label{(g4)}
$\lim_{t \to 0} \frac{G(t)}{|t|^{p+1}} = + \infty$, where $p = 2^m_s=1 + \frac{4s}{N}$.
		\end{itemize}
		Then for any $m >0$, the problem \eqref{problem_frac} has a solution, satisfying the Pohozaev identity \eqref{eq_Pohozaev_fractional}.
\end{Theorem}
	
We highlight that the found solution is actually a minimum for $\mc{L}$ constrained to the sphere (see Proposition \ref{minimizing}), which furnishes a strong indication to its stability properties.
The techniques employed in \cite{Shi0} for the local case $s=1$, to get directly the existence of a minimum for $\mc{L}$, are not easily adaptable to the fractional framework, because of the need of a control on the tails in the Brezis-Lieb lemma and in the Concentration-Compactness techniques. 
Anyway, our method not only gets around these difficulties, but moreover it is also suitable to get multiple solutions.

%\label{pag_comm_minimum}

Indeed, if we also suppose the oddness of $g$, namely
\begin{itemize}
	\item[(g5)] \label{(g5)}
$g(-t)= - g (t)$ for all $t \in \R$,
\end{itemize}
we have $\mc{I}^m(\lambda,-u)= \mc{I}^m(\lambda,u)$ for all $(\lambda,u) \in \R \times H^s_r(\R^N)$ and we can establish the existence of infinitely many $L^2$-constrained standing waves solutions for the (fNLS) equation. 

We prove the following multiplicity result.

\begin{Theorem}\label{S:1.13}	
Suppose $N\geq 2$ and \hyperref[(g1)]{\textnormal{(g1)}}--\hyperref[(g3)]{\textnormal{(g3)}} and \hyperref[(g5)]{\textnormal{(g5)}}. 
Then we have:
\begin{itemize}
	\item[\textnormal{(i)}]
	 For any $k \in \N$ there exists $m_k \geq 0$ such that for each $m > m_k$, the problem \eqref{problem_frac} has at least $k$ nontrivial, distinct pairs of solutions, satisfying the Pohozaev identity \eqref{eq_Pohozaev_fractional}. 
	\item[\textnormal{(ii)}]
	In addition assume \hyperref[(g4)]{\textnormal{(g4)}}. For any $m >0$ the problem \eqref{problem_frac} has countably many solutions $(u_n)_{n }$ (satisfying the Pohozaev identity \eqref{eq_Pohozaev_fractional}), which verify 
	$$\mathcal{L}(u_n) <0 \quad \hbox{for all $ n \in \N$},$$ 
	$$\mathcal{L}(u_n) \to 0 \quad \hbox{as $ n \to + \infty$}.$$
\end{itemize}
\end{Theorem}

We remark that our subcritical multiplicity result seems new even in the case of the pure power $g(t)=|t|^{q-2}t$ and in the non-monotone case of competing powers $g(t)= |t|^{q-2}t - |t|^{r-2}t$, and it has a physical relevance since it describes the existence of multiple bound states with arbitrary high energies (see e.g. \cite{DZ0}).
We stress that the analytical solutions for fractional differential equations are still limited, while there is a large amount of numerical methods in discretizing the fractional differential operators. 
In Theorem \ref{S:1.13} we furnish an analytical rigorous approach to detect infinitely many symmetric solitons, which can be applied to the computation of ground and excited states to (fNLS) equations.
% arising from Bose-Einstein condensation theory or nonlinear optics phenomena with saturation.
%
%
%
% 
\begin{Remark}\label{rem_weak_g2}
%CONTROLLA: \\%CCCOMMENT NOW
We highilight that \hyperref[(g2)]{\textnormal{(g2)}} can be weakened, with no changes in the proofs, by asking, for some $q \in (2^m_s-1, 2^*_s-1)$
$$\lim_{|t| \to + \infty} \frac{g(t)}{|t|^{q}} =0 \quad \hbox{ and } \quad \limsup_{|t| \to +\infty} \frac{G(t)}{|t|^{2^{m}_s}} =\limsup_{|t| \to +\infty} \frac{g(t)}{|t|^{2^{m}_s-1}} =0.$$
See also Remark \ref{rem_caso_subcr_nonstr} %and \cite{GalSch} 
for some additional discussions.
\end{Remark}

\begin{Remark}
We highlight that we assume \emph{a priori} the positivity of the Lagrange multiplier $\mu$ in \eqref{problem_frac}. As a matter of fact, this condition seems to be quite natural: %indeed, 
if $u$ is a ground state on the sphere $\int_{\R^N} u^2 \, dx=m$ and its energy is negative, then \emph{a posteriori} the corresponding Lagrange multiplier $\mu$ is strictly positive (see Proposition \ref{prop_equiv_gs}). 
In addition, from a physical perspective, in the study of standing waves the multiplier $\mu$ describes the frequency of the particle, and thus it is positive; moreover, this prescribed sign is characteristic also of chemical potentials in the description of ideal gases, see \cite{LSSY,PS0}.
\end{Remark}

\medskip

The Chapter is organized as follows. 
%\tr{aggiusta label sezioni} %COMMENT NOW
In Section \ref{sec_frac_unconstr}, we establish some preliminaries related to the unconstrained problem. 
In Section \ref{sec_lagrange_frac} we give the Lagrangian formulation of the problem \eqref{problem_frac} and a description of the geometry of a %an auxiliary 
functional in a product space. 
Section \ref{section:4} concerns with the Palais-Smale-Pohozaev ((PSP) for short) condition and Section \ref{section:5} is devoted to the construction of the deformation argument under this (PSP) condition. 
Section \ref{section:6} deals with our minimax procedure to detect the normalized solutions by means of the Pohozaev mountain. 
Finally in Section \ref{sec_frac_multiple_sol} we derive the multiplicity result of infinitely many normalized solutions when $g$ is odd.

%%%%%%%%%%%%%%%%%%%%%%%%%%%%%%%%%%%%%%%%%%%%%%%%%%%%%%%

%\setcounter{equation}{0} %FOR ARXIV
\section{The unconstrained problem} 
% in $H^s(\R^N)$}
%\label{section:2}
\label{sec_frac_unconstr}

In this Section we consider the unconstrained fractional equation
%\begin{equation} \label{problemxx}
%\parag{
%&(-\Delta)^{s} u + \mu u =	f(u)& \quad \hbox{in $\R^N$,} \cr 
%& u \in H^s(\R^N),& \cr
%& u>0,&
%}
%\end{equation}
\begin{equation} \label{problemxx}
(-\Delta)^{s} u + \mu u =	g(u) \quad \hbox{in $\R^N$,} 
\end{equation}
where $s \in (0,1)$, $N \geq 2$, $u \in H^s(\R^N)$, % (strictly positive), 
$\mu >0$ is fixed and $g$ satisfies \hyperref[(g1)]{\textnormal{(g1)}} together with the following assumptions
\vskip2pt
\begin{itemize}
%	\item[(g1)] \label{(g1')}
%$g : \R \to \R$ continuous %, $g(t)= 0$ for $t \leq 0$ 
%and $\lim_{t \to 0} \frac{g(t)}{t}=0;$
	\item[(g2')] \label{(g2')}
 $\limsup_{|t| \to \infty} \frac{g(t)}{|t|^q} =0$ where $ q \in (1, 2^*_s -1)$, where $2^*_s = \frac{2N}{N-2s}$; %1+ \frac{4s}{N-2s})$;
	\item[(g3')] \label{(g3')}
 there exists $t_0 > 0$ such that $G(t_0) >\frac{\mu}{2}t_0^2$, where $G(t)= \int_0^t g(\tau) d\tau$. 
\end{itemize}
\vskip 2pt

Under the assumptions \hyperref[(g1)]{\textnormal{(g1)}}-\hyperref[(g2')]{\textnormal{(g2')}}, it is standard to show that any weak solution of \eqref{problemxx} is a critical point of the $C^1$-functional $\mc{J}_{\mu}: H^s(\R^N) \to \R$ defined by
$$
\mc{J}_{\mu} (u):=\half \int_{\R^N} | (-\Delta)^{s/2} u|^2 \, dx + \frac{\mu}{2} \int_{\R^N} u^2 \, dx - \int_{\R^N} G(u) \, dx. 
$$

In the celebrated paper \cite{BL1}, for the local case $s=1$, Berestycki and Lions proved the existence of a classical solution to \eqref{problemxx}, which is radially symmetric and has an exponentially decay, under the assumption \hyperref[(g1)]{\textnormal{(g1)}}-\hyperref[(g2')]{\textnormal{(g2')}}-\hyperref[(g3')]{\textnormal{(g3')}}; 
these conditions are almost optimal for the existence of \eqref{problemxx}. 
The found solution is of least energy among all nontrivial solutions, and a Mountain Pass (MP for short) solution as shown by Jeanjean and Tanaka \cite{JT0}. 
%indeed, when $s=1$, in proved that the last energy solution is indeed a Mountain Pass (MP for short) solution. 
Successively Byeon, Jeanjean, Maris \cite{BJM0} showed that every least energy solution of \eqref{problemxx} has constant sign and is radially symmetric (and decreasing) up to translations.

For the nonlocal case $s \in (0,1)$, 
%the equivalence between MP weak solution and least energy solution is still a partially open problem (see Section \ref{sec_proper_pohozaev_gs}).
we begin to recall that in the recent paper \cite{BKS}, Byeon, Kwon and Seok established the following results (see also \cite{ChWa0}).
\begin{Proposition}[Regularity]\label{Regularity}
	Suppose \hyperref[(g1)]{\textnormal{(g1)}}-\hyperref[(g2')]{\textnormal{(g2')}}. Let $u \in H^s(\R^N)$ be a weak solution 
	of the fractional equation \eqref{problemxx}.
	Then $u \in C^1(\R^N)$ if one of the following assumptions holds:
	\begin{itemize}
		\item[\textnormal{(i)}]	$s \in (1/2,1)$;
	\item[\textnormal{(ii)}]
	$s \in (0,1/2]$ and $g \in C^{0,\sigma}_{loc}(\R)$ for some $\sigma \in (1-2s,1)$.
	\end{itemize}	
\end{Proposition}

\begin{Proposition}[Fractional Pohozaev identity]\label{Pohozaev-prop}
	Suppose \hyperref[(g1)]{\textnormal{(g1)}}-\hyperref[(g2')]{\textnormal{(g2')}} %-\hyperref[(g3')]{\textnormal{(g3')}} 
and
	\begin{itemize}
		\item[\textnormal{(g4')}] \label{(g4')} 
if $s \in (0, 1/2]$, $g \in C^{0, \sigma}_{loc} (\R)$ for some $\sigma \in (1-2s,1)$.
	\end{itemize}
	\vskip 2pt
	Then every weak solution $u \in H^s(\R^N)$ of the fractional equation of \eqref{problemxx} satisfies the Pohozaev identity \eqref{eq_Pohozaev_fractional}, % with $e^{\lambda}=\mu$, 
%(or Pohozaev-Derrick identity)
%\begin{equation}\label{eq_Pohozaev_fractional}
%\frac{N-2s}{2} \int_{\R^N} | (-\Delta)^{s/2}u|^2 + N \int_{\R^N} \left( \frac{\mu}{2} u^2 - G(u)\right) =0.
%\end{equation}
%%	 \int_{\R^N} | (-\Delta)^{s/2}u|^2 \, dx + 2^*_s \int_{\R^N} \left( \frac{\mu}{2} u^2 - G(u)\right) dx =0$$
which can be rewritten as
$$%	\begin{equation*}%\label{Pohozaev}
%	\frac{N-2s}{2} \int_{\R^N} | (-\Delta)^{s/2}u|^2 \, dx + N \int_{\R^N} \left( \frac{\mu}{2} u^2 - F(u)\right) dx =0.
%	 \int_{\R^N} | (-\Delta)^{s/2}u|^2 \, dx + 2^*_s \int_{\R^N} \left( \frac{\mu}{2} u^2 - G(u)\right) dx =0.
\frac{1}{2^*_s} \int_{\R^N} | (-\Delta)^{s/2}u|^2 + \frac{\mu}{2^{\#}} \int_{\R^N} u^2 - \int_{\R^N}G(u) =0,
$$%	\end{equation*}
where $2^*_s = \frac{2N}{N-2s}$ and $2^{\#}=2$ are the upper and lower critical exponents. 
\end{Proposition}
Roughly, we see that the Pohozaev identity essentially means $\frac{d}{d\theta} \mc{J}_{\mu}(u(\cdot/e^{\theta}))_{|\theta=0}=0$, thus it is strictly related to the scaling invariance of the problem (which will be exploited through an augmented functional, see \eqref{ugual}). Anyway the fact that $u$ is a critical point for $\mc{J}_{\mu}'$ does not imply directly this relation, since $\frac{d}{d\theta} \mc{J}_{\mu}(u(\cdot/e^{\theta}))_{|\theta=0} = (\mc{J}_{\mu}'(u), \nabla u \cdot x)=0$ requires some restriction on $\mc{J}_{\mu}$ and $u$.

%We remark that 
Indeed, the $C^1$-regularity of the weak solution seems crucial for proving formally a Pohozaev type identity.
Under \hyperref[(g1)]{\textnormal{(g1)}}-\hyperref[(g2')]{\textnormal{(g2')}}-\hyperref[(g3')]{\textnormal{(g3')}} we know \cite{BKS} that each weak solution of \eqref{problemxx} belongs to $H^s(\R^N) \cap C^\beta (\R^N)$ with $\beta \in (0, 2s)$ and thus it is not known if the Pohozaev identity holds when $s \in (0,1/2]$, without additional regularity assumptions on the nonlinearity $g$.

\smallskip

 In \cite%[Theorem 1.2]
{BKS}, the authors further investigated the existence of MP weak solutions of \eqref{problemxx}. 
 We recall that a weak solution $u$ is said of MP type if 
\begin{equation}\label{MPlevel}
\mc{J}_{\mu}(u) = a(\mu), %C_{mp},
\end{equation}
where
$$
%C_{mp} 
a(\mu):= \inf_{\gamma\in\Gamma}\max_{t\in [0,1]} \mathcal{J}_{\mu}(\gamma(t))
$$
and
\begin{equation}\label{classpath}
\Gamma_{\mu} := \big\{ \gamma(t)\in C\big([0,1], H^s_r(\R^N)\big) \mid \gamma(0)=0, \, \mc{J}_{\mu}(\gamma(1))<0\big\}.
\end{equation}
As for $s=1$, the functional $\mc{J}_{\mu}$ does not satisfies the Palais-Smale condition at level $a(\mu)$ under the assumptions \hyperref[(g1)]{\textnormal{(g1)}}-\hyperref[(g2')]{\textnormal{(g2')}}-\hyperref[(g3')]{\textnormal{(g3')}}, thus one can not directly apply the MP theorem. For the local case $s=1$, any weak solution is $C^1$ and it satisfies the Pohozaev identity, so that one can reduce the search of MP solutions to that of minimizers on the Pohozaev type constraint.
For the fractional case, this approach seems to work for $s \in (1/2, 1)$, while requires additional regularity on the nonlinearity if $s \in (0,1/2]$.
 
Conversely, in \cite{BKS} the authors established that every minimizer of $\mc{J}_{\mu}$ on the Pohozaev %type 
constraint corresponds to a MP weak solution and derived some radially symmetric properties of the minimizer using a fractional version of the Polya-Szego inequality.
Namely they introduce the Pohozaev functional $\mathcal{P}: H^s_r(\R^N) \to \R$ by setting 
$$ \mathcal{P}_{\mu}(u) := \frac{N-2s}{2} \int_{\R^N} |(-\Delta)^{s/2}u|^2 + N\int_{\R^N} \left(\frac{\mu}{2} u^2 - G(u)\right) $$
%$$ \mathcal{P}(u) := \tor{ \frac{1}{2^*_s} \int_{\R^N} |(-\Delta)^{s/2}u|^2 + \frac{\mu}{2} \int_{\R^N} u^2 - \int_{\R^N}G(u),} $$
and 
$${P}_{\mu} :=\big\{ u \in H^s_r(\R^N) \setminus \{ 0 \} \, | \, \mathcal{P}_{\mu}(u)=0 \big \},$$
$$ %C_{po} 
p(\mu):= \min_{u \in{P}_{\mu}} \mathcal{J}_{\mu}(u).$$
In \cite[Theorem 1.2]{BKS} they established the following result.
\begin{Theorem}\label{Byeonminimizers}
Assume \hyperref[(g1)]{\textnormal{(g1)}}-\hyperref[(g2')]{\textnormal{(g2')}}-\hyperref[(g3')]{\textnormal{(g3')}}. Let $s \in (0,1)$ and $\mu >0$. Then 
	\begin{itemize}
		\item[\textnormal{(i)}]
			there exists a minimizer of $\mathcal{J}_{\mu}$ subject to ${P}_{\mu}$;
			\item[\textnormal{(ii)}]
			every minimizer of $\mathcal{J}_{\mu}$ subject to ${P}_{\mu}$ is a MP weak solution of 
		\eqref{problemxx};
%		\item[\textnormal{(iii)}]	
%			every positive minimizer of $\mc{J}_{\mu}$ subject to ${P}_{\mu}$ is radially symmetric up to a translation.
\end{itemize}
\end{Theorem}
From Theorem \ref{Byeonminimizers} it follows that 	
$$a(\mu)=p(\mu).$$
%$$ C_{mp} =C_{po}. $$
While the equivalence between Mountain Pass solutions and least energy solutions is shown for $s\in(1/2,1)$, it is yet an open problem for $s\in(0,1/2]$ under the assumptions \hyperref[(g1)]{\textnormal{(g1)}}-\hyperref[(g2')]{\textnormal{(g2')}}-\hyperref[(g3')]{\textnormal{(g3')}}. 
In \cite{BKS}, this equivalence is established under the same regularity assumption of Proposition 2.1, namely $g\in C^{0,\sigma}(\R^N)$ for some $\sigma\in (1-2s,1)$; see Section \ref{sec_proper_pohozaev_gs} for more comments on this relation.
In the following Sections, in contrast, we will show that, under $L^2$-constraint, least energy solutions have Mountain Pass characterization.
See Proposition \ref{prop_equiv_gs}. % Proposition \ref{minimizing}. 
%\label{page_minim_2}

\begin{Remark}\label{rem_byeon_addit}
We highlight that in \cite{BKS} they define $a(\mu)$ and $p(\mu)$ on the whole space $H^s(\R^N)$, by additionally assuming 
$$g(t) \equiv 0 \quad \hbox{for $t\leq 0$};$$
indeed, thanks to this assumptions, they can pass from a generic minimization sequence to a positive one, and thus to a radially symmetric one (and exploit then compactness). With this additional assumption they also show that 
\begin{itemize}
\item[\textnormal{(iii)}]	
			every minimizer of $\mc{J}_{\mu}$ subject to ${P}_{\mu}$ is positive and radially symmetric up to a translation.
\end{itemize}
Without this assumption we notice that their arguments show that every \emph{positive} minimizer of $\mc{J}_{\mu}$ subject to ${P}_{\mu}$ is radially symmetric up to a translation.

On the other hand, without this additional assumption but by assuming \hyperref[(g4')]{\textnormal{(g4')}}, one may argue as follows: by a result similar to Proposition \ref{prop_esist_ground} (see \cite[Theorem 3]{BL1} and \cite[Lemma 1]{BJM0} for details) the ground state problem can be seen as a minimization problem; thus we can apply \cite[Theorem 4.1]{LoMa08} to deduce the radial symmetry (up to a translation) of any minimizer. See also \cite{LMS}.
%every minimizer is positive and radially symmetric up to a translation this is actually needed to ensure that the found solutions are positive. 
%In that case point $(iii)$ of Theorem \ref{Byeonminimizers} holds for every minimizer.
\end{Remark}

\begin{Remark}\label{rem_esist_nonstrett_unc}
The existence of a Pohozaev minimum (Mountain Pass solution) when 
$$\limsup_{|t| \to \infty} \frac{g(t)}{|t|^{2^*_s-1}} =0$$
substitutes \hyperref[(g2')]{\textnormal{(g2')}} can be found in \cite{ChWa0}, where %: in this case, on the other hand, 
it is assumed that $g \in C^1(\R)$ (or, more specifically, it is sufficient that Pohozaev holds for every solution). A result involving this assumption together with $g \in C(\R)$ seems to lack in literature, even though the proof by \cite{BKS} can be easily adapted. % (see also \cite{GalSch}). 
%\tr{CONTROLLARE:} \\ %CCCOMMENT NOW
 Anyway, we can obtain this result as a byproduct of our argument, similarly to Section \ref{sec_unconstrained_choquard}. 
See instead \cite{Iko2} for the existence of infinitely many solutions. 
\end{Remark}
%\tr{Il risultato di esistenza di ground state si può estendere al caso frazionario sottocritico (non strettamente)? Non riesco a trovare una referenza (Byeon è sottocritico stretto, Chang chiedere $C^1$ per usare Pohozaev, Ambrosio chiede anche lui regolarità, Ikoma chiede a massa nulla (?), ...)}

Some further properties of this autonomous equation will be invetigated in Section \ref{sez_limit_eq} and in Section \ref{sec_recalls}.

%%%%%%%%%%%%%%%%%%%%%%%%%%%%%%%%%%%%%%%%%%%%%%%%%%%%%%%

%\setcounter{equation}{0} %FOR ARXIV
 \section{Lagrangian formulation and Pohozaev geometry} 
% of $\mc{I}^m(\lambda,u)$} 
%\label{section:3}
 \label{sec_lagrange_frac}

We come back to the constrained case; from now on in this Chapter we briefly denote 
$$ p:=2^m_{s}-1= 1+ \frac{4s}{N}.$$
We consider the Lagrangian formulation of the problem \eqref{problem_frac} in the space of radially symmetric functions $H_r^s(\R^N)$.
Namely, we seek for critical points of the functional $\mc{I}^m: \R \times H^s_r(\R^N) \to \R$ % defined in \eqref{functlag}, i.e.
\begin{equation}\label{equation}
 \mc{I}^m(\lambda,u):=\half \int_{\R^N} |(-\Delta)^{s/2} u|^2 - \int_{\R^N} G(u)\, + 
 \frac{e^\lambda}{2} \big(\|u\|_2^2 -m\big).
 \end{equation}
Under the assumption \hyperref[(g1)]{\textnormal{(g1)}}--\hyperref[(g3)]{\textnormal{(g3)}}, it is standard to prove that $\mc{I}^m$ is $C^1$ in the product space $\R \times H^s_r(\R^N)$. 
It is immediate to recognize that for any $m >0$
$$\mc{I}^m(\lambda,u)=\mc{J}(\lambda,u) -{\frac{e^\lambda}{2}}m$$
where $\mc{J}:\R \times H^s_r(\R^N) \to \R$ is the $C^1$-functional defined by $\mc{J}(\lambda, u) :=\mc{J}_{e^{\lambda}}(u)$, i.e.
$$
\mc{J}(\lambda, u) :=\half \int_{\R^N} |(-\Delta)^{s/2} u|^2 - \int_{\R^N} G(u)\, + 
\frac{e^\lambda}{2} \int_{\R^N} u^2.
$$
For a fixed $\lambda \in \R$, $u$ is critical point of $\mc{J}(\lambda, \cdot)$ means that $u\in H^s_r(\R^N)$ solves, in the weak sense,
\begin{equation}\label{problemxxx}
%\parag{
%&(-\Delta)^{s} u + e^\lambda u = g(u)& \qquad \hbox{in $\R^N$,} \\ 
%& u \in H^s_r(\R^N).&
%}
(-\Delta)^{s} u + e^\lambda u = g(u) \qquad \hbox{in $\R^N$}. 
\end{equation}
Inspired by the Pohozaev identity \eqref{eq_Pohozaev_fractional}, for any $s \in (0,1)$ we also introduce the Pohozaev functional $\mc{P}:\R \times H^s_r(\R^N) \to \R$ by setting 
%\tr{cambia la definizione?}
$$ 
\mc{P}(\lambda, u) :=
\frac{N-2s}{2} \int_{\R^N} |(-\Delta)^{s/2} u|^2 + N \int_{\R^N} \left( \frac{e^\lambda}{2} u^2 - G(u)\right).
$$
%$$ 
%\mc{P}(\lambda, u) :=
%\tor{\frac{1}{2^*_s} \int_{\R^N} |(-\Delta)^{s/2} u|^2 - \int_{\R^N} G(u)+ \frac{e^\lambda}{2} \int_{\R^N} u^2.}
%$$
By Proposition \ref{Pohozaev-prop}, it follows that for any $\lambda \in \R$, if $u \in H^s_r(\R^N)$ solves \eqref{problemxxx}, then $\mc{P}(\lambda, u)=0$ when $s \in (\frac{1}{2},1)$.
A similar result for $s \in (0,\frac{1}{2}]$ is not known under \hyperref[(g1)]{\textnormal{(g1)}}--\hyperref[(g3)]{\textnormal{(g3)}}.

\smallskip

We introduce now the Pohozaev set
$$ \Omega :=\big\{(\lambda,u) \in \R \times H^s_r(\R^N) \mid \mc{P}(\lambda,u)>0\big\} \cup\big\{(\lambda,0) \mid \lambda \in \R \big\}.$$
Since $\int_{\R^N} G(u)=o(\norm u_{H^s}^2)$ as $u \to 0$ we have the following.
\begin{Lemma}
We have
	\begin{equation}\label{eq_interior_omega}
\{ (\lambda,0) \mid \lambda\in \R\} \subset int(\Omega).
	\end{equation}
\end{Lemma}

\claim Proof.
For any fixed $\delta>0$ there exists a suitable $C_{\delta}>0$ such that
$$G(t)\leq \delta |t|^2 + C_{\delta} |t|^{p+1},$$
where $p+1 < 2^*_s$. Thus
\begin{align*}
	0&= \frac{1}{N} \mathcal{P}(\lambda,u) \geq 
%\tor{ \mathcal{P}(\lambda,u_{\lambda}) } 
	\frac{1}{2^*_s} %\frac{N-2s}{2N}
	\norm{(-\Delta)^{s/2}u}_2^2 +
	\left(\frac{e^{\lambda}}{2} -\delta\right) \norm{u}_2^2 - C_{\delta}\norm{u}_{p+1}^{p+1} \\
	&\gtrsim \norm{u}_{H^s}^2 - \norm{u}_{H^s}^{p+1}>0
\end{align*}
for $\delta$ small and $ \norm{u}_{H^s}$ small, $u\neq 0$.
\QED
%La positività non è uniforme in $\lambda$: per $\lambda\to -\infty$ serve $\delta\to 0$, e quindi $C_{\delta}$ arbitrarimaente grande, quindi $ \norm{u}_{H^s}$ arbitrariamente piccolo. Quindi non posso dire che $\partial \Omega$ è uniformemente lontano da $u=0$: anzi, si avvicina sempre di più, quindi l'inf è nullo!

\bigskip

%We note that for each $\lambda\in\R$, $\mc{P}(\lambda,u)>0$ in a small neighborhood of $u=0$ except $0$, thus
%\begin{equation}\label{eq_interior_omega}
%\big\{(\lambda,0) \mid \lambda \in \R \big\} \subset int(\Omega)
%\end{equation}
This last result implies
%which implies
$$ \partial \Omega =\big\{(\lambda,u) \in \R \times H^s_r(\R^N) \mid \mc{P}(\lambda,u)=0, \ u \neq 0 \big\}; $$
we call this set the \emph{Pohozaev mountain}. % for $\mc{J}$. 
We remark that $(\lambda, u) \in \partial \Omega$ if and only if $u \neq 0$ and $u$ satisfies the Pohozaev identity $\mc{P}(\lambda, u)=0$.

Contrary to assumption \hyperref[(g3')]{\textnormal{(g3')}}, the arbitrariness of the frequency $\mu$ and the corresponding assumption \hyperref[(g3)]{\textnormal{(g3)}} lead to different interactions between the pieces $\mu u$ and $g(u)$, which have to be taken into account; these interactions are described by the quantity
%Set
\begin{equation}\label{muzero}
\mu_0 := 2 \sup_{t \in \R, t \neq 0} \frac{G(t)}{t^2};
\end{equation}
we deduce $\mu_0 \in (0, +\infty]$ under the assumptions \hyperref[(g1)]{\textnormal{(g1)}}--\hyperref[(g3)]{\textnormal{(g3)}}. 
%In what follows, 
We also denote 
\begin{equation}\label{lambdazero}
\lambda_0 := \log (\mu_0), \quad \hbox{if} \ \mu_0 \in (0, \infty),
\end{equation}
otherwise $\lambda_0 := + \infty$. 
Analysing the two cases ($\lambda_0 \in \R$ and $\lambda_0=+\infty$) will be of key importance in the study of the geometry of the problem.

Taking into account that $2^m_s < 2^*_s$, %$1+ \frac{4s}{N} < 1 + \frac{4s}{N-2s},$ 
we deduce by $(i)$ Theorem \ref{Byeonminimizers} that for any $\lambda \in (- \infty, \lambda_0)$ the functional 
$$ u \in H^s_r(\R^N)\mapsto \mc{J}(\lambda,u) \in \R $$
has a minimizer $u_\lambda$ subject to 
$$(\partial \Omega)_{\lambda} :=
\big\{ u \in H_r^s (\R^N) \setminus \{ 0 \} \mid \mc{P}(\lambda,u)=0 \big\}, $$
namely
\begin{equation}\label{minimumpohozaev}
\mc{J}(\lambda,u_\lambda)= \min_{u \in (\partial \Omega)_{\lambda}} \mc{J}(\lambda,u).
\end{equation}
Furthermore by $(ii)$ of Theorem \ref{Byeonminimizers} such $u_\lambda$ is a Mountain Pass critical point of $\mathcal{J}(\lambda, \cdot)$ at level $a(\lambda)$, i.e.
$$ \mc{J}(\lambda,u_\lambda)= a(\lambda) $$
where 
\begin{equation}\label{MPlevelx}
a(\lambda) := \inf_{\gamma\in\Gamma(\lambda)}\max_{t\in [0,1]} \mc{J}(\lambda,\gamma(t))
 \end{equation}
and
\begin{equation}\label{classpathx}
\Gamma(\lambda) := \big\{ \gamma \in C\big([0,1], H^s_r(\R^N)\big) \mid \gamma(0)=0, \, \mc{J}(\lambda,\gamma(1))<0\big\}.
\end{equation}
We notice that $\lambda \in (-\infty,\lambda_0)\mapsto a(\lambda)\in\R$ is strictly monotone increasing: this can be shown, for example, by relying on the fact that $a(\lambda)$ coincides with the Pohozaev minimum and exploiting some scaling argument.\footnote{ % on $\R$.
%Let $\mu<\nu$, and $v$ be a $\nu$-Pohozaev minimum (i.e. $J_{\nu}(v)=a(\nu)$). Let rescale $v$ such that it belongs to the $\mu$-Pohozaev, i.e. $u:=v(\lambda \cdot)$: straightforwards computation shows the explicit value of $\lambda$ and that (since $\mu<\nu$) $\lambda <1$. Thus $a(\mu) \leq J_{\mu}(u) < J_{\nu}(v) = a(\nu)$. The strict inequalities comes from the Pohozaev identities, the scaling and from $\lambda<1$.}
Let $\lambda_1, \lambda_2 \in \R$ %<\lambda_2$, 
and $v$ be a $\lambda_2$-Pohozaev minimum (i.e. $\mc{J}(\lambda_2,v)=a(\lambda_2)$ and $\mc{P}(\lambda_2,v)=0$). Let rescale $v$ in such a way it belongs to the $\lambda_1$-Pohozaev set, i.e. $u:=v( \cdot/\theta)$ with $\mc{P}(\lambda_1,u)=0$, for some explicit $\theta = \left(1 + \frac{2^*_s}{2} (\lambda_2-\lambda_1) \frac{\norm{v}_2^2}{\norm{(-\Delta)^{s/2}v}_2^2}\right)^{-\frac{1}{2s}}$. Thus, by the Pohozaev identities,
$$a(\lambda_1) \leq \mc{J}(\lambda_1,u) = \frac{s}{N} \norm{(-\Delta)^{s/2}u}_2^2 = \theta^{N-2s}\frac{s}{N}\norm{(-\Delta)^{s/2}v}_2^2 = \theta^{N-2s} \mc{J}(\lambda_2,v)= \theta^{N-2s} a(\lambda_2).$$
%the equality among the functionals is due to the Pohozaev identities. 
If $\lambda_1<\lambda_2$ then $\theta <1$ and thus the have claim $a(\lambda_1)< a(\lambda_2)$. %strict inequalities comes from the Pohozaev identities, the scaling and from $\theta<1$.

As a further result, since $\theta \to 1$ as $\lambda_1$ and $\lambda_2$ approach, we obtain also $ a(\lambda_1) \leq \liminf_{\lambda_2 \to \lambda_1} a(\lambda_2)$ and $\limsup_{\lambda_1 \to \lambda_2} a(\lambda_1) \leq a(\lambda_2)$. Swapping the role of $\lambda_1$ and $\lambda_2$ 
%We are using here that one liminf exists to say that the limif of the product is the product of the liminf (since the inequality of the general case is not in the right direction).
%Since in this last passages we did not use the monotonicity relation among $\lambda_1$ and $\lambda_2$, 
actually we obtain the (extra) continuity property: $\lim_{\lambda \to \lambda_0} a(\lambda) = a(\lambda_0)$.
}

\begin{Lemma}\label{lem_buona_def}
Let $\lambda \in \R$. 
Then the following statements are equivalent:
\begin{itemize}
\item[\textnormal{(a)}] $\lambda < \lambda_0$.
\item[\textnormal{(b)}]
There exists a $t_0=t_0(\lambda)>0$ such that
$$ G(t_0) > \frac{e^{\lambda}}{2} t_0^2.$$
\item[\textnormal{(c)}]
There exists $u\in H^s_r(\R^N)\setminus \{0\}$ such that $\mc{P}(\lambda, u)=0$; in particular $(\partial \Omega)_{\lambda}\neq \emptyset$.
\item[\textnormal{(d)}]
$\Gamma(\lambda)\neq \emptyset$, and thus $a(\lambda)$ is well defined.
\end{itemize}
As further consequence, we see that $\partial \Omega \neq \emptyset$. Finally, $a(\lambda)>0$.
\end{Lemma}

\claim Proof.
(a) $\iff$ (b).
This is a straightforward consequence of the definition of $\lambda_0$.

(b) $\implies$ (c)
Let $u\in H^s_r(\R^N)$ to be fixed. 
We have, for $t>0$,
$$\mc{P}(\lambda, u(\cdot/t))= \frac{N-2s}{2} t^{N-2s} \int_{\R^N} |(-\Delta)^{s/2}u|^2 - N t^{N} \int_{\R^N} \left(G(u) -\frac{ e^{\lambda}}{2} u^2 \right).$$
%$$\mc{P}(\lambda, u(\cdot/t))= \tor{\frac{1}{2^*_s} \int_{\R^N} |(-\Delta)^{s/2}u|^2 - t^{N} \int_{\R^N} \left(G(u) -\frac{ e^{\lambda}}{2} u^2 \right)}.$$
We notice that $\mc{P}(\lambda, u(\cdot/t))>0$ for small $t>0$. In order to get a $\bar{t}$ such that $\mc{P}(\lambda, u(\cdot/\bar{t}))=0$ we need the quantity
$$ \int_{\R^N} \left(G(u) - \frac{e^{\lambda}}{2} u^2\right)$$
to be positive. 
For any $R>0$ we choose a smooth $u=u_R\in C^{\infty}_c$ such that $u_R=t_0$ in $B_R(0)$ and $u_R=0$ out of $B_{R+\frac{1}{R^N}}(0)$, $0\leq u_R\leq t_0$. 
We set
$$C:=\sup_{t\in[0,t_0]} \left|G(t)-\tfrac{e^{\lambda}}{2} |t|^2\right|<+\infty.$$
Then
%\begin{eqnarray*}
%\lefteqn{ \int_{\R^N} \left(G(u_R) - \frac{e^{\lambda}}{2} u_R^2 \right)} \\
%&=& \int_{B_{R+\frac{1}{R^N}}\setminus B_R}\left(G(u_R) - \frac{e^{\lambda}}{2} u_R^2 \right) +\int_{B_R} \left(G(u_R) - \frac{e^{\lambda}}{2} u_R^2 \right) \\
%&\geq& -C \abs{B_{R+\frac{1}{R^N}} \setminus B_R} + |B_R| \left(G(t_0) - \frac{e^{\lambda}}{2} |t_0|^2\right) \to +\infty
%\end{eqnarray*}
\begin{align*}
 \int_{\R^N} \left(G(u_R) - \frac{e^{\lambda}}{2} u_R^2 \right)
&= \int_{B_{R+\frac{1}{R^N}}\setminus B_R}\left(G(u_R) - \frac{e^{\lambda}}{2} u_R^2 \right) +\int_{B_R} \left(G(u_R) - \frac{e^{\lambda}}{2} u_R^2 \right) \\
&\geq -C \abs{B_{R+\frac{1}{R^N}} \setminus B_R} + |B_R| \left(G(t_0) - \frac{e^{\lambda}}{2} |t_0|^2\right) \to +\infty
\end{align*}
and in particular it is positive for a sufficiently large $R$.

(c) $\implies$ (d).
Let $u\in H^s_r(\R^N) $, $u\nequiv 0$ such that $\mc{P}(\lambda,u)=0$. 
We define $\gamma(t):=u(\cdot/t)$ for $t\neq 0$ and $\gamma(0)= 0$, so that $\gamma:[0, \infty)\to H^s_r(\R^N)$ is continuous.
We have
$$\mc{J}(\lambda, \gamma(t))= \frac{1}{2} t^{N-2s} \int_{\R^N} |(-\Delta)^{s/2}u|^2 - t^N \int_{\R^N} \left(G(u) - \frac{e^{\lambda}}{2} u^2\right).$$
Noting $\int_{\R^N}\left(G(u)-\frac{e^\lambda}2 u^2\right)>0$ by $\mc{P}(\lambda,u)=0$, we have
$\mc{J}(\lambda,\gamma(t))\to -\infty$ as $t\to\infty$ and thus $\Gamma(\lambda)\not=\emptyset$.

(d) $\implies$ (b).
If $\gamma\in \Gamma(\lambda)$, then $\mathcal{J}(\lambda,\gamma(1))<0$, thus
$$\int_{\R^N} \left(G(\gamma(1)) - \frac{e^{\lambda}}{2} \gamma(1)^2\right)>0,$$
which implies that there exists an $x_0\in \R^N$ such that
$$G(\gamma(1)(x_0)) - \frac{e^{\lambda}}{2} \gamma(1)^2(x_0)>0.$$
The claim comes by setting $t_0:=\gamma(1)(x_0)$.

Finally, by Theorem \ref{Byeonminimizers}, there exists a Pohozaev minimum $u_{\lambda}$ which is also a Mountain Pass solution, %\tr{potrei evitare di usare l'esistenza per il problema unconstrained, così come abbiamo fatto per il frazionario-Choquard}, 
thus $\mathcal{J}(\lambda,u_{\lambda})=a(\lambda)$, $D_u \mathcal{J}(\lambda, u_{\lambda})=0$ and $\mathcal{P}(\lambda,u_{\lambda})=0$, which imply
$$ a(\lambda)=\frac{s}{N} \| (-\Delta)^{s/2} u_\lambda \|^2_2>0.
\QED
$$
%see also Remark \ref{rem_altern_proof_>0} %\ref{...} for an approach 
%for a proof which does not require an existence theorem for the unconstrained problem.
%\QED

\medskip

\begin{Remark} 
Assume $\lambda_0 <+\infty$. 
We observe that, in this case, for $\lambda\geq \lambda_0$ we have $\mc{P}(\lambda, u)\geq 0$ and $\mc{J}(\lambda, u)\geq 0$ for each $u$, both strictly positive for $u\nequiv 0$. 
This means that $[\lambda_0,+\infty)\times H^s_r(\R^N) \subset \Omega$.
\end{Remark}

In the next result, we consider the case $\lambda_0 \in \R$ and we investigate the behaviour of $a(\lambda)$ as $\lambda$ approach $\lambda_0$.

\begin{Proposition}\label{S:2.33}
	Assume \hyperref[(g1)]{\textnormal{(g1)}}--\hyperref[(g3)]{\textnormal{(g3)}} and $\lambda_0 \in \R$. We have 
	\smallskip
	\begin{itemize}
		\item[\textnormal{(a)}] if $(\lambda,u) \in \partial \Omega$ for some $u \in H^s_r(\R^N)$, then $\lambda < \lambda_0$.
		\item[\textnormal{(b)}] $\lim_{\lambda \to \lambda_0^-} a(\lambda)=+\infty$.
	\end{itemize}
\end{Proposition}

\claim Proof.
Let $(\lambda, u) \in \partial \Omega$, namely $\mathcal{P}(\lambda,u)=0$ and $u \neq 0$. 
This implies that for some $x_0 \in\R^N$ 
$$ G(u(x_0)) - \frac{e^\lambda}{2} u(x_0)^2 >0 $$
and thus $\lambda < \lambda_0$ and (a) holds.

Now we show point (b). Let $\lambda<\lambda_0$; by contradiction, since by Lemma \ref{lem_buona_def} $a(\lambda)$ is increasing and strictly positive, we assume that $a(\lambda)\to c \in (0, +\infty)$ as $\lambda \to \lambda_0^-$, from which we deduce that $\norm{(-\Delta)^{s/2}u_{\lambda}}_2$ is bounded. Moreover, for any fixed $\delta>0$ there exists a suitable $C_{\delta}>0$ such that
$$G(t)\leq \delta |t|^2 + C_{\delta} |t|^{p+1},$$
where we recall that $p=1+\frac{4s}{N}$.
%Basta su $G$. Va bene anche critico.
 
Thus we have by the fractional Gagliardo-Nirenberg inequality \eqref{GN} and the fact that $\norm{(-\Delta)^{s/2} u_{\lambda}}_2$ is bounded,
\begin{align*}
	0&= \frac{1}{N} \mathcal{P}(\lambda,u_{\lambda}) \geq 
%\tor{ \mathcal{P}(\lambda,u_{\lambda}) } 
	\frac{1}{2^*_s} %\frac{N-2s}{2N}
	\norm{(-\Delta)^{s/2}u_{\lambda}}_2^2 +
	\left(\frac{e^{\lambda}}{2} -\delta\right) \norm{u_{\lambda}}_2^2 - C_{\delta}\norm{u_{\lambda}}_{p+1}^{p+1} \\
	&\geq \frac{1}{2^*_s} \norm{(-\Delta)^{s/2}u_{\lambda}}_2^2 + \left(\frac{e^{\lambda}}{2} -\delta\right) \norm{u_{\lambda}}_2^2 
 - C'C_{\delta}
	 \norm{(-\Delta)^{s/2} u_{\lambda}}_2^2 \norm{u_{\lambda}}_2^{p-1} \\
	&\geq \left(\frac{e^{\lambda}}{2} -\delta\right) \norm{u_{\lambda}}_2^2 - C'' C_{\delta} \norm{u_{\lambda}}_2^{\frac{4s}{N}}
\end{align*}
for some $C'$, $C''>0$. 
By choosing $\delta < \frac{e^{\lambda}}{2}$, since $\frac{4s}{N}<2$, also $\norm{u_{\lambda}}_2$ must be bounded, which means that $(u_\lambda)_{\lambda<\lambda_0}$ is bounded in $H^s_r(\R^N)$. 
Hence, up to a subsequence, $u_{\lambda}\rightharpoonup u_0$ in $H^s_r(\R^N)$. By the immersion \eqref{eq_immer_compatt} %Lemma \ref{compact} 
and taking into account that $\partial_u\mathcal{J}(\lambda, u_{\lambda})=0$, we deduce that $u_\lambda \to u_0$ strongly in $H^s_r(\R^N)$ with $\mathcal{J}(\lambda_0, u_0)=c$, $\partial_u\mathcal{J}(\lambda_0, u_0)=0$, $\mathcal{P}(\lambda_0, u_0)=0$.
Since $c>0$, we have $u_0 \neq 0$. By $\mathcal{P}(\lambda_0, u_0) =0$, we conclude
$$ G(u_0(x)) - \frac{e^{\lambda_0}}{2} u_0(x)^2 >0 $$
for some $x \in \R^N$, which contradicts the definition of $\lambda_0$. 
\QED

\bigskip

We consider now the case $\lambda_0= + \infty$ and we investigate the behaviour of $a(\lambda)$ for $\lambda$ large.
\begin{Proposition}\label{S:lim}
Assume that $\lambda_0= + \infty$. Then
$$ \lim_{\lambda\to +\infty} \frac{a(\lambda)}{e^{\lambda}} = +\infty.$$
\end{Proposition}

\claim Proof.
%By (g1)-(g2) we have that for any $\delta>0$ there exists $C_\delta>0$ such that for all $t\in\R$
%\tor{
%\begin{eqnarray}
%\abs{g(t)} &\leq& \delta\abs t^p + C_\delta\abs t, \label{maggiorazione-0} \\
%\abs{G(t)} &\leq& 
%\frac{\delta}{p+1}|t|^{p+1} + \frac{C_\delta}{2} \abs t^2, \label{maggiorazione} 
%\end{eqnarray}
By \hyperref[(g1)]{\textnormal{(g1)}}-\hyperref[(g2)]{\textnormal{(g2)}} we have that for any $\delta>0$ there exists $C_\delta>0$ such that for all $t\in\R$
\begin{equation}
%\abs{G(t)} \leq
G(t) \leq
\frac{\delta}{p+1}|t|^{p+1} + \frac{C_\delta}{2} \abs t^2. \label{maggiorazione} 
\end{equation}
%where $p = 1 + \frac{4s}{N}$.} 
%Serve sottocritico!
We also denote by $b(\delta)$ the MP value of $\mc{H}_{\delta} : H_r^s(\R^N) \to \R$ defined by
$$\mc{H}_{\delta}(v) := \half \norm{(-\Delta)^{s/2}v}_2^2 + \half \norm{v}_2^2 -\frac{\delta}{p+1}\norm{v}_{p+1}^{p+1}.$$
It is easy to see that\footnote{
Indeed, by scaling, $\mc{H}_{\delta}(\delta^{-\frac{1}{p-1}}\cdot) = \delta^{-\frac{2}{p-1}} \mc{H}_1$, which implies $\Gamma_{\delta}= \delta^{-\frac{1}{p-1}} \Gamma_1$; here $\Gamma_{\delta}$ is the set of paths related to $\mc{H}_{\delta}$. Using these two relations one obtains $b(\delta)=\delta^{-\frac{2}{p-1}} b(1) \to +\infty$.}
$$b(\delta) \to +\infty \quad \hbox{as}\ \delta\to 0^+.$$
For $v\in H^s_r(\R^N)\setminus\{ 0\}$, we set 
$$u_\theta :=\theta^{N/2}v(\theta \cdot),$$
and for simplicity we write $\mu \equiv e^\lambda$ and $\mc{J}(\mu,\cdot) =\mc{J}(\lambda,\cdot)$. 
By \eqref{maggiorazione}, we pass to evaluate
$$\mc{J}(\mu,u_\theta) \geq \theta^{2s} \left(\half \norm{(-\Delta)^{s/2}v}_2^2 
 + \half (\mu - C_\delta) \theta^{-2s} \norm{v}_2^2 
 - \frac{\delta}{p+1} \norm{v}_{p+1}^{p+1} \right).$$
Setting $\theta := (\mu - C_\delta)^{1/{2s}}$ for $\mu =\mu_{\delta}> C_\delta$, we have
\begin{equation}\label{key1}
\mc{J}(\mu, u_{(\mu - C_\delta)^{1/{2s}}} ) \geq (\mu -C_\delta) \mc{H}_{\delta} (v)
\end{equation}
and hence
\begin{equation}\label{key2}
\frac{\mc{J}(\mu, u_{(\mu - C_\delta)^{1/{2s}}})}{\mu} \geq \frac{\mu - C_\delta}{\mu}
\mc{H}_\delta (v).
\end{equation}
Thus we have %\tr{RICONTROLLA ($\frac{1}{\mu}\to 0$ quando $\delta \to 0$...)}
\begin{equation}\label{key3}
 \frac{a(\mu)}{\mu} \geq \frac{\mu-C_\delta}{\mu} b(\delta).
\end{equation}
Choosing $\mu= \mu_{\delta}>2 C_{\delta}$ we obtain
$$ \frac{a(\mu)}{\mu} \geq \frac{b(\delta)}{2};$$
since $\delta>0$ is arbitrary, we derive
$$ \lim_{\mu\to +\infty} \frac{a(\mu)}{\mu} = +\infty. 
\QED
$$

\bigskip

Finally we investigate the behaviour of $a(\lambda)$ for $\lambda \to -\infty$, under some more restrictive assumption in the origin.

\begin{Proposition}\label{out}
	Assume \hyperref[(g4)]{\textnormal{(g4)}} in addition to \hyperref[(g1)]{\textnormal{(g1)}}--\hyperref[(g3)]{\textnormal{(g3)}}. 
	Then
	\begin{equation}\label{third}
	\lim_{\lambda\to -\infty} \frac{a(\lambda)}{e^{\lambda}} = 0.
	\end{equation}
\end{Proposition}
\claim Proof.
We fix $u\in H^s_r(\R^N) \cap L^\infty(\R^N)$ with $\norm u_\infty=1$. Recalled $p = 1+ \frac{4s}{N}$, by \hyperref[(g4)]{\textnormal{(g4)}} there exists $M_r>0$
such that for all $r \in (0,1]$ 
$$G(ru(x)) \geq M_r %\frac{M_r}{p+1}
r^{p+1} \abs{u(x)}^{p+1}, \ \quad \forall x\in\R^N$$
with 
$$M_r\to+\infty \quad \hbox{as}\ r\to 0.$$
%Basta su $G$
%"Basta" $\liminf$ (il che significa che serve il limite): $M_r:= \inf_{t \in (0,r]} \frac{G(t)}{|t|^{p+1}}$, da cui la disequazione per $t=r|u(x)|$, e inltre $\lim_{r\to 0} M_r = \lim_{r \to 0} \inf_{t \in (0,r]} \frac{G(t)}{|t|^{p+1}} = \liminf_{r \to 0}\frac{G(t)}{|t|^{p+1}} = +\infty$.
We write again $\mu \equiv e^{\lambda}$ for the sake of simplicity. Therefore for $t>0$ we have 
\begin{eqnarray*}
\lefteqn{ \mc{J}(\mu,r u(\cdot/t)) \leq \half r^2 t^{N-2s} \|(-\Delta)^{s/2}u \|_2^2 +\frac{\mu}{2}r^2t^N \norm u_2^2
	- M_r %\frac{M_r}{p+1}
 r^{p+1} t^{N} \|u\|_{p+1}^{p+1} } \\ 
&=& r^2 \mu^{-\frac{N-2s}{2s}} \left( \half t^{N-2s} \mu^{\frac{N-2s}{2s}} 
	\|(-\Delta)^{s/2}u \|_2^2 +\frac{1}{2} \mu^{\frac{N}{2s}} t^N \norm u_2^2
	- M_r %\frac{M_r}{p+1}
 r^{\frac{4s}{N}} \mu^{\frac{N-2s}{2s}} t^{N} \|u\|_{p+1}^{p+1}\right) 	\\ 
&=& r^2 \mu^{-\frac{N-2s}{2s}} \left( \frac{1}{2}\tau^{N-2s} 
	\|(-\Delta)^{s/2}u \|_2^2 +\frac{1}{2} \tau^{N} \norm u_2^2
	- M_r %\frac{M_r}{p+1}
 r^{\frac{4s}{N}} \mu^{-1} \tau^{N} \|u\|_{p+1}^{p+1} \right)
	\end{eqnarray*}
after setting $\tau := \mu^{\frac{1}{2s}}t$. Moreover choosing $r:=\mu^{\frac{N}{4s}}$ we infer 	
\begin{eqnarray*}
\lefteqn{ \mc{J}\left(\mu,\mu^{\frac{N}{4s}} u(\cdot/(\mu^{-1/(2s)}\tau))\right) }
\\ &\leq& \mu\left( \half\tau^{N-2s} 	\|(-\Delta)^{s/2}u \|_2^2 + \half \tau^N \norm u_2^2
-M_{\mu^{N/(4s)}}%\frac{M_{\mu^{N/(4s)}}}{p+1} 
\tau^N \|u\|_{p+1}^{p+1}\right). 
\end{eqnarray*}
For $\mu\in (0,1)$, the map
$$ \tau \in (0,\infty)\mapsto \mu^{\frac{N}{4s}} u(\cdot/\mu^{-1/(2s)}\tau)\in H^s_r(\R^N)$$
can be regarded as a path in $\Gamma(\mu)$ after rescaling. 
Thus
$$ \frac{a(\mu)}{\mu} \leq \max_{\tau\in [0,\infty)} 
\left( \half 
\|(-\Delta)^{s/2}u \|_2^2 \tau^{N-2s} + \half \norm u_2^2\tau^N - M_{\mu^{N/(4s)}}%\frac{M_{\mu^{N/(4s)}}}{p+1} 
 \|u\|_{p+1}^{p+1}\tau^{N}\right). $$
Since $M_{\mu^{N/(4s)}} \to \infty$ as $\mu\to 0$, we derive the conclusion. 
\QED

\medskip

\begin{Proposition}\label{S:2.2}
	Assume \hyperref[(g1)]{\textnormal{(g1)}}--\hyperref[(g3)]{\textnormal{(g3)}}. 
	Then we have 
\begin{itemize}
\item[\textnormal{(a)}] $\mc{J}(\lambda,u)\geq 0$ for all $(\lambda,u)\in \Omega$;
\item[\textnormal{(b)}] $\mc{J}(\lambda,u)\geq a(\lambda)>0$ for all $(\lambda,u)\in \partial\Omega$.
\end{itemize}
\end{Proposition}

\claim Proof.
We notice that for all $(\lambda,u)\in \Omega$ 
$$ \mc{J}(\lambda,u) \geq \mc{J}(\lambda,u) - \frac{1}{N} \mc{P}(\lambda,u) 
= \frac{s}{N} \norm{(-\Delta)^{s/2} u}_2^2 \geq 0 $$
and thus (a) follows.

The proposition (b) follows from the fact that every minimizer of $\mc{J}(\lambda,\cdot)$ subject to $(\partial \Omega)_{\lambda}$ is a mountain pass weak solution of \eqref{problemxx} at level $a(\lambda)$ (see Theorem \ref{Byeonminimizers}).
\QED

\bigskip

%\begin{Remark}\label{rem_altern_proof_>0}
%In order to show that $a(\lambda)>0$, without exploiting the existence result for the unconstrained problem, we argue as follows (see also \cite{CGT2}). 
%Let $\gamma \in \Gamma(\lambda)$; by definition of $\Gamma(\lambda)$ and by Proposition \ref{S:2.2} (a) there exists $t^*$ such that $\gamma(t^*) \in \partial \Omega$ and $\gamma(t^*) \neq 0$, thus $\mc{P}(\lambda, \gamma(t^*))=0$. This means that
%$$\mc{J}(\lambda, \gamma(t^*)) =\frac{\alpha +2s}{2(N + \alpha)} \norm{(-\Delta)^{s/2} \gamma(t^*)}_2^2+ \frac{\alpha \mu}{2(N + \alpha)} \norm{\gamma(t^*)}_2^2 \simeq \norm{\gamma(t^*)}_{H^s}^2 $$
%thus
%$$a(\lambda) \gtrsim \inf_{u \in (\partial \Omega)_{\lambda}} \norm{u}_{H^s}^2.$$
%Since, by \eqref{eq_interior_omega}, $ (\partial \Omega)_{\lambda}$ is far from the line $(\lambda,0)$, we obtain that the right-hand side is strictly positive, which is the claim.
%\end{Remark}

We are ready to show that for any $m >0$ the functional $\mc{I}^m$ is bounded from below on the Pohozaev set $\partial \Omega$.
\begin{Proposition}\label{negativo}
Assume \hyperref[(g1)]{\textnormal{(g1)}}--\hyperref[(g3)]{\textnormal{(g3)}}.
For any $m>0$, we set
$$B_m:=\inf_{\lambda < \lambda_0} \left(a(\lambda) -\frac{e^\lambda}{2}m \right)$$
and
$$E_m:= \inf_{(\lambda,u)\in\partial\Omega}\mc{I}^m(\lambda,u).$$
Then
\begin{equation}\label{limitinf}
E_m \geq B_m >-\infty.
\end{equation}
\end{Proposition}
\claim Proof.
Let $m>0$. 
If $(\lambda, u) \in \partial \Omega$, by (b) of Proposition \ref{S:2.2} we have
$$\mc{I}^m(\lambda, u) = \mathcal{J}(\lambda, u) - \frac{e^{\lambda}}{2} m \geq a(\lambda) - \frac{e^{\lambda}}{2} m;$$
since, by (a) of Proposition \ref{S:2.33} it results that $\lambda < \lambda_0$, we have, passing to the infimum,
$$E_m \geq B_m.$$
We distinguish now two cases.
\smallskip
Firstly we assume $\lambda_0 \in \R$. From (b) of Proposition \ref{S:2.33}
we have $a(\lambda) \to + \infty$ as $\lambda \to \lambda_0^-$, and thus we conclude 
$$\inf_{\lambda < \lambda_0} \left(a(\lambda) -\frac{e^\lambda}{2}m \right)> -\infty.$$
Secondly, we suppose that $\lambda_0= +\infty$. We have
$$a(\lambda) -\frac{e^\lambda}{2}m = e^{\lambda} \left(\frac{a(\lambda)}{e^{\lambda}}- \frac{m}{2}\right)$$
and thus, by Proposition \ref{S:lim}
$$\inf_{\lambda \in \R} \left(a(\lambda) -\frac{e^\lambda}{2}m \right)> -\infty.
\QED
$$

%%%%%%%%%%%%%%%%%%%%%%%%%%%%%%%%%%%%%%%%%%%%%%%%%%%%%%%

%\setcounter{equation}{0} %FOR ARXIV
\section{Compactness by scaling} %Palais-Smale-Pohozaev condition and augmented functional}
\label{section:4}

Firstly we introduce the notations:
\begin{align*}
K_b :=&\, \big\{ (\lambda,u)\in \R\times H^s_r(\R^N) \mid \mc{I}^m(\lambda, u)=b,\,
 \partial_\lambda\mc{I}^m(\lambda, u)=0,\, \partial_u\mc{I}^m(\lambda, u)=0 \big\}, \\
K^{PSP}_b :=& \, \big\{ (\lambda,u)\in K_b \mid \mc{P}(\lambda,u)=0 \big\}.
\end{align*}
Clearly, we have $K^{PSP}_b \subset K_b.$
We note that for the definition of $K^{PSP}_b$ we do not need additional regularity about $g$.

Under the assumptions \hyperref[(g1)]{\textnormal{(g1)}}--\hyperref[(g3)]{\textnormal{(g3)}}, it seems difficult to verify the standard Palais-Smale condition for the functional $\mc{I}^m$.
Therefore we cannot recognize that the set $K_b$ is compact.

Inspired \cite{HT0,IT0}, we introduce the Palais-Smale-Pohozaev (shortly (PSP)) condition, which is a weaker compactness condition than the standard Palais-Smale one. Such (PSP) condition takes into account the scaling properties of $\mc{I}^m$ through the Pohozaev functional $\mc{P}$. 
Using this new condition we will show that $K^{PSP}_b$ is compact when $b<0$.

%%%%%%%%%%%%%%%%%%%%%%%%%%%%%%%%%%%%%

\subsection{A limiting Pohozaev identity} %The Palais-Smale-Pohozaev condition}

We give the definition of (PSP) condition in the radial setting.

\smallskip
\begin{Definition}\label{PSPcondition2}
	For $b \in \R$, we say that $\mc{I}^m$ satisfies the Palais-Smale-Pohozaev condition at level $b$ (shortly the $(PSP)_b$ condition), if for any sequence $(\lambda_n, u_n)_n \subset \R \times H^s_r(\R^N)$ such that 
	\begin{equation}\label{prima}
	\mc{I}^m (\lambda_n, u_n) \to b,
	\end{equation}
	\begin{equation}\label{seconda}	
\partial_{\lambda} 	\mc{I}^m(\lambda_n, u_n) \to 0, 
	\end{equation}
	\begin{equation}\label{seconda3}
\norm{\partial_u 	\mc{I}^m(\lambda_n, u_n)}_{(H^s_r(\R^N))^*} \to 0,
\end{equation}
	\begin{equation}\label{terza4.14}
	\mc{P}(\lambda_n, u_n) \to 0,
	\end{equation}
	it happens that $(\lambda_n, u_n)_n$ has a strongly convergent subsequence in $\R \times H^s_r(\R^N)$.
\end{Definition}

We will show the following result.

\begin{Proposition}\label{PSP}
	Assume 
	\hyperref[(g1)]{\textnormal{(g1)}}--\hyperref[(g3)]{\textnormal{(g3)}}. Let $b <0$. Then $\mc{I}^m$ satisfies the $(PSP)_b$ condition on $\R\times H_r^s(\R^N)$.
\end{Proposition}

\claim Proof.
Let $b <0$ and suppose that $(\lambda_n, u_n) \subset \R \times H^s_r(\R^N)$ satisfies \eqref{prima}--\eqref{terza4.14}. 
We will show that $(\lambda_n, u_n)$ has a strongly convergent subsequence in several steps.

\smallskip

\noindent
\textbf{Step 1:} \emph{$\lambda_n$ is bounded from below.}
\\
Indeed
\begin{align*}
	\frac{m}{2} e^{\lambda_n} &= \frac{1}{N} \mc{P}(\lambda_n, u_n) - \mc{I}^m(\lambda_n, u_n) + \frac{s}{N}\norm{(-\Delta)^{s/2}u_n }_2^2 \\
	&\geq \frac{1}{N} \mc{P}(\lambda_n, u_n) -\mc{I}^m(\lambda_n, u_n)
\end{align*}
hence
$$\frac{m}{2} \liminf_n e^{\lambda_n} \geq 0 -b >0,$$
which implies (since $m>0$) that $\lambda_n$ is bounded from below. 

\smallskip

\noindent
\textbf{Step 2:} \emph{$\norm{u_n}_2^2\to m$.} 
\\
Indeed, we have
$$\partial_{\lambda} \mc{I}^m(\lambda_n, u_n) = \frac{e^{\lambda_n}}{2} \left(\norm{u_n}_2^2 -m\right)\to 0,$$
which implies the claim by Step 1.

\smallskip

\noindent
\textbf{Step 3:} \emph{$\norm{(-\Delta)^{s/2}u_n}_2^2$ and $\lambda_n$ are bounded (from above) as $n\to +\infty$.} 
\\
Indeed, by \hyperref[(g1)]{\textnormal{(g1)}}-\hyperref[(g2)]{\textnormal{(g2)}} we have that for any $\delta>0$ there exists $C_\delta>0$ such that for all $t\in\R$
\begin{equation}
%\abs{g(t)} \leq \delta\abs t^p + C_\delta\abs t. \label{maggiorazione-0}
g(t) \leq \delta\abs t^p + C_\delta\abs t. \label{maggiorazione-0}
\end{equation}
 %Serve su $g$. Serve sottocritico!
By \eqref{maggiorazione-0} and the fractional Gagliardo-Nirenberg inequality \eqref{GN} we have 
\begin{align*}
%	\abs{\partial_{u}\mc{I}^m(\lambda_n, u_n)u_n} &\geq \norm{(-\Delta)^{s/2}u_n}_2^2 + e^{\lambda_n} \norm{u_n}_2^2 - \intRN |g(u_n)u_n| \\
	\partial_{u}\mc{I}^m(\lambda_n, u_n)u_n &= \norm{(-\Delta)^{s/2}u_n}_2^2 + e^{\lambda_n} \norm{u_n}_2^2 - \intRN g(u_n)u_n \\
	&\geq \norm{(-\Delta)^{s/2}u_n}_2^2 + \left(e^{\lambda_n} - C_{\delta}\right) \norm{u_n}_2^2 - \delta \norm{u_n}_{p+1}^{p+1} \\
	&\geq \norm{(-\Delta)^{s/2}u_n}_2^2 + \left(e^{\lambda_n} - C_{\delta}\right) \norm{u_n}_2^2 - \delta C \norm{(-\Delta)^{s/2}u_n}_2^2 \norm{u_n}_2^{p-1};
\end{align*}
moreover
\begin{align*}
	\abs{\partial_{u} \mc{I}^m(\lambda_n, u_n)u_n} &\leq \norm{\partial_{u}\mc{I}^m(\lambda_n, u_n)}_{(H^s_r(\R^N))^*} \norm{u_n}_{H^s_r(\R^N)}\\
	&= \norm{\partial_{u}\mc{I}^m(\lambda_n, u_n)}_{(H^s_r(\R^N))^*} \sqrt{\norm{(-\Delta)^{s/2} u_n}_2^2 + \norm{u_n}_2^2}.
\end{align*}
Set $\varepsilon_n := \norm{\partial_{u}\mc{I}^m(\lambda_n, u_n)}_{(H^s_r(\R^N))^*}$ and (by Step 2) $\norm{u_n}_2^2= m + o(1)$, we finally have, joining the previous two inequalities, that
\begin{eqnarray*}
\lefteqn{\norm{(-\Delta)^{s/2} u_n}_2^2 \left(1-\delta C(m+o(1))^{\frac{p-1}{2}}\right) + \left(e^{\lambda_n} - C_{\delta}\right) ( m + o(1))}\\ 
&&\qquad \qquad \qquad \qquad \qquad %\qquad 
\leq \varepsilon_n \sqrt{\norm{(-\Delta)^{s/2} u_n}_2^2 + m + o(1)}.
\end{eqnarray*}
Choosing $\delta>0$ small so that $\delta Cm^{\frac{p-1}{2}}<1$, we obtain the claim.

\smallskip

\noindent
\textbf{Step 4:} \emph{Conclusion.} 
\\
By Steps 1-3, we have that $(\lambda_n, u_n)$ is bounded in $\R\times H^s_r(\R^N)$. 
Hence, up to a subsequence, $\lambda_n \to \lambda$ and $u_n \rightharpoonup u$ in $H^s_r(\R^N)$. 
Therefore, we obtain (see Proposition \ref{prop_converg_generiche_loc})
$$\intRN g(u_n)u_n \to \intRN g(u)u \quad \hbox{ and } \quad \intRN g(u_n)u \to \intRN g(u)u.$$ 
%Serve Sobolev sottocritico per $g$ (quindi va bene $L^2$-(sotto)critico per $g$, ma non per $G$)
Again by the assumption $\partial_{u}\mc{I}^m(\lambda_n, u_n)\to 0$ we get
\begin{align}
	0&=\lim_n \partial_{u}\mc{I}^m(\lambda_n, u_n)u \notag \\ 
	&= \lim_n \left(\intRN (-\Delta)^{s/2}u_n (-\Delta)^{s/2}u +e^{\lambda_n} \intRN u_n u - \intRN g(u_n)u\right) \notag \\
	&= \norm{(-\Delta)^{s/2} u}_2^2 + e^{\lambda} \norm{u}_2^2 - \intRN g(u)u. \label{eq_dim_PSP1}
\end{align}
Since $\partial_{u}\mc{I}^m (\lambda_n, u_n)\to 0$ and $u_n \rightharpoonup u$, we have $\partial_{u}\mc{I}^m (\lambda_n, u_n)u_n\to 0$; thus
\begin{align}
	0&=\lim_n \partial_{u}\mc{I}^m(\lambda_n, u_n)u_n \notag \\
	 &= \lim_n \left(\norm{(-\Delta)^{s/2}u_n}_2^2 +e^{\lambda_n} \norm{u_n}_2^2- \intRN g(u_n)u_n\right) \notag \\
	&= \lim_n\left(\norm{(-\Delta)^{s/2} u_n}_2^2 + e^{\lambda_n} \norm{u_n}_2^2\right) - \intRN g(u)u \label{eq_dim_PSP2}
\end{align}
and hence, joining \eqref{eq_dim_PSP1} and \eqref{eq_dim_PSP2},
$$\norm{(-\Delta)^{s/2} u_n}_2^2 + e^{\lambda_n} \norm{u_n}_2^2 \to \norm{(-\Delta)^{s/2} u}_2^2 + e^{\lambda} \norm{u}_2^2, $$
which easily implies (since $e^{\lambda_n}\to e^{\lambda}$ and $\norm{u_n}_2^2$ is bounded)
$$\norm{u_n}_{\lambda}^2 \to \norm{u}_{\lambda}^2, $$
where $\norm{\cdot}_{\lambda}^2 := \norm{(-\Delta)^{s/2} \cdot}_2 + e^{\lambda} \norm{\cdot}_2^2$ is an equivalent norm on $H^s_r(\R^N)$. 
This, together with $u_n \rightharpoonup u$ in $H^s_r(\R^N)$ %and the fact that $H^s_r(\R^N)$ is a Hilbert space, 
gives $u_n\to u$ strongly in $H^s_r(\R^N)$. 
\QED

\medskip
%\bigskip

%\newpage %€

\begin{Corollary}\label{PSP-cor}
	Assume \hyperref[(g1)]{\textnormal{(g1)}}--\hyperref[(g3)]{\textnormal{(g3)}}. 
	Let $b \in \R$, $b <0$. Then $K^{PSP}_b \cap (\R \times \{0\}) = \emptyset$ and $K^{PSP}_b$ is compact.
\end{Corollary}
\claim Proof.
Since 	$\partial_\lambda \mc{I}^m(\lambda, 0) = - \frac{e^\lambda}{2m} \neq 0$, we have $K^{PSP}_b \cap (\R \times \{0\}) = \emptyset$.
Proposition \ref{PSP} implies that $K^{PSP}_b$ is compact. 
\QED

\medskip

\begin{Remark}
	We emphasize that the $(PSP)_b$ condition does not hold at level $b=0$. Indeed we can consider the unbounded sequence $(\lambda_j, 0)$ with $\lambda_j \to - \infty$ such that
	$$ \mc{I}^m(\lambda_j, 0) = \partial_\lambda \mc{I}^m(\lambda_j, 0) = - \frac{e^{\lambda_j}}{2} m \to 0$$ 
	and 
	$$ \partial_u \mc{I}^m(\lambda_j, 0) = 0, \quad \mc{P}(\lambda_j, 0) =0. $$ 
\end{Remark}

%%%%%%%%%%%%%%%%%%%%%%%%%%%%%%%%%%%%%

\subsection{A functional in an augmented space} %An augmented functional}
\label{sec_frac_augmented}

Following \cite{Jea0,HIT,HT0} we introduce the augmented functional $\mc{H}^m: \R \times \R \times H^s_r(\R^N) \to \R$
\begin{equation}\label{ugual}
\mc{H}^m(\theta, \lambda, u) := \mc{I}^m(\lambda, u(e^{-\theta}\cdot)) 
\end{equation}
for $(\theta, \lambda,u ) \in \R \times \R \times H^s_r(\R^N)$. 
By the scaling properties of $\mc{I}^m$ we can recognize that
\begin{equation}\label{eq:22}
\mc{H}^m(\theta, \lambda, u)= 
\frac{e^{(N-2s)\theta}}{2}
\int_{\R^N} |(-\Delta)^{s/2}u|^2 - e^{N \theta}
\int_{\R^N} G(u)\, + 
\frac{e^\lambda}{2} \bigl( e^{N \theta}\|u\|_2^2 -m \bigr) 
\end{equation}
for all $(\theta, \lambda,u ) \in \R \times \R \times H^s_r(\R^N).$

Moreover, by standard calculations we have the following proposition.
\begin{Proposition}
	For all $(\theta, \lambda,u ) \in \R \times \R \times H^s_r(\R^N)$, $h \in H^s_r(\R^N)$, $\beta \in \R$, we have
	\begin{itemize}
	\item[\textnormal{(i)}]
		$\partial_\theta \mc{H}^m(\theta, \lambda, u) = \mc{P}(\lambda, u(\cdot /e^\theta)),	$		
	\item[\textnormal{(ii)}]			
		$\partial_\lambda \mc{H}^m(\theta, \lambda, u) = \partial_\lambda \mc{I}^m(\lambda, u(\cdot / e^\theta)), $
	\item[\textnormal{(iii)}]		 
		$\partial_u \mc{H}^m(\theta, \lambda, u) h = \partial_u \mc{I}^m(\lambda, u(\cdot/ e^\theta)) h(\cdot/e^\theta), $		
	\item[\textnormal{(iv)}]
		$ \mc{H}^m(\theta + \beta, \lambda, u(e^\beta \cdot)) = \mc{H}^m(\theta, \lambda, u). $		
	\end{itemize}
\end{Proposition}

Now we define a metric on the Hilbert manifold
$$M := \R \times \R \times H^s_r(\R^N)$$
by setting 
\begin{align*}
{\|(\alpha, \nu, h)\|}_{(\theta,\lambda, u)}^2 :=& \pabs{\left(\alpha, \nu,\norm{h(e^{-\theta} \cdot)}_{H^s_r(\R^N)}\right)}^2 \\
=& \, \alpha^2 + \nu^2 + e^{N\theta} \norm{h}_2^2+ e^{(N-2s)\theta} \norm{(-\Delta)^{s/2} h}_2^2 
\end{align*}
for any $(\alpha, \nu, h) \in T_{(\theta,\lambda,u)} M \equiv \R \times \R \times H^s_r(\R^N)$.
We also denote the dual norm on 
$T^*_{(\theta,\lambda,u)}M$ by $\|\cdot \|_{(\theta,\lambda, u), *}$. 
We notice that 
${\|(\cdot, \cdot, \cdot)\|}_{(\theta,\lambda, u)}^2$ depends only on $\theta$ and we can write ${\|(\cdot, \cdot, \cdot )\|}_{(\theta,\cdot, \cdot)}^2$.
Moreover for any $(\alpha, \nu, h) \in T_{(\theta,\lambda,u)}M$ and $\beta \in \R$ we have
\begin{equation}\label{eq_shift_norma}
{\|(\alpha, \nu, h(e^\beta x))\|}_{(\theta + \beta,\cdot, \cdot)}^2=
{\|(\alpha, \nu, h)\|}_{(\theta,\cdot, \cdot)}^2.
\end{equation}
Furthermore we define the standard distance between two points as the infimum of length of curves connecting the two points, namely
$$ \dist_M\big((\theta_0, \lambda_0, h_0), (\theta_1, \lambda_1, h_1)\big) := 
\inf_{\gamma \in \mathcal{G}} \int_0^1 \|\dot \gamma(t)\|_{\gamma(t)} dt $$
where 
$$\mathcal{G} :=\left\{\gamma \in C^1([0,1],M) \, \middle | \, \gamma(0)= (\theta_0, \lambda_0, h_0), \gamma(1)= (\theta_1, \lambda_1, h_1) \right\}.$$

Observe that, if $\sigma$ is a path connecting $(\alpha_0, \nu_0, h_0)$ and $(\alpha_1, \nu_1, h_1)$, then by \eqref{eq_shift_norma} $\tilde{\sigma}(t):=(\sigma_1(t)+\beta, \sigma_2(t), (\sigma_3(t))(e^{\beta}\cdot))$ is a path connecting $(\alpha_0 +\beta, \nu_0, h_0(e^{\beta}\cdot))$ and $(\alpha_1+\beta, \nu_1, h_1(e^{\beta}\cdot))$ with same length, and hence
\begin{equation}\label{eq_shift_distance}
\dist_M\big((\alpha_0, \nu_0, h_0), (\alpha_1, \nu_1, h_1)\big) = \dist_M\big((\alpha_0 +\beta, \nu_0, h_0(e^{\beta}\cdot)), (\alpha_1+\beta, \nu_1, h_1(e^{\beta}\cdot))\big).
\end{equation}

Denote now for simplicity $D:=(\partial_\theta,\partial_\lambda,\partial_u)$ the gradient with respect to all the variables; a direct computation shows that
$$D\mc{H}^m(\theta, \lambda, u)(\alpha,\nu,h) 
= \mc{P} (\lambda, u(e^{-\theta} \cdot))\alpha +\partial_{\lambda} \mc{I}^m(\lambda, u(e^{-\theta}\cdot))\nu
+\partial_u \mc{I}^m(\lambda, u(e^{-\theta}\cdot))h(e^{-\theta} \cdot)$$
and thus we obtain
\begin{eqnarray*}
\lefteqn{\|{D\mc{H}^m(\theta, \lambda, u)\|}_{(\theta, \lambda, u),*}^2 }\\
	&=& \left|\left(\mathcal{P}(\lambda, u(e^{-\theta} \cdot)), \partial_{\lambda}\mc{I}^m(\lambda, u(e^{-\theta}\cdot)), \norm{\partial_u \mc{I}^m(\lambda, u(e^{-\theta}\cdot))}_{(H^s_r(\R^N))^*}\right)\right|^2 \\
	&=& \abs{\mathcal{P}(\lambda, u(e^{-\theta} \cdot))}^2 + \abs{\partial_{\lambda}\mc{I}^m(\lambda, u(e^{-\theta}\cdot))}^2 
		+ \norm{\partial_u \mc{I}^m(\lambda, u(e^{-\theta}\cdot))}_{(H^s_r(\R^N))^*}^2 .
\end{eqnarray*}
Now defined
$$\tilde{K}_b :=\big \{ (\theta, \lambda, u) \in M \mid \mc{H}^m(\theta, \lambda, u)=b,\, D \mc{H}^m(\theta, \lambda,u)=0\big\}$$
the set of critical points at level $b$ of $\mc{H}^m$, we deduce
\begin{equation}\label{eq_confronto_K}
\tilde{K}_b = \big\{(\theta, \lambda, u(e^{\theta} \cdot)) \mid (\lambda, u)\in K^{PSP}_b, \; \theta \in \R\big\}.
\end{equation}

\begin{Proposition}\label{PSPtilde}
	Assume \hyperref[(g1)]{\textnormal{(g1)}}--\hyperref[(g3)]{\textnormal{(g3)}}. 
	Let $b \in \R$, $b <0$. Then the functional $\mc{H}^m$ satisfies the following Palais-Smale type condition $(\widetilde{PSP})_b$: 
	%That is, 
for each sequence $(\theta_n, \lambda_n, u_n)_n$ such that
	$$\mc{H}^m(\theta_n, \lambda_n, u_n) \to b,$$
	$$\norm{D \mc{H}^m(\theta_n, \lambda_n, u_n)}_{(\theta_n, \lambda_n, u_n),*} \to 0,$$
	we have, up to a subsequence,
	$$\dist_M((\theta_n, \lambda_n, u_n), \tilde{K}_b)\to 0.$$
\end{Proposition}
We note that $(\widetilde{PSP})_b$ condition is different from the standard Palais-Smale condition and it ensures the compactness of the sequence $(\theta_n,\lambda_n,u_n)_n$ after a suitable scaling. 
By \eqref{eq_confronto_K} we also highlight that, if $\tilde K_b\not=\emptyset$, then $\tilde K_b$ is not compact. % (see \eqref{eq_confronto_K}).

\medskip

\claim Proof.	
Let $(\theta_n, \lambda_n, u_n)_n$ be as in $(\widetilde{PSP})_b$. Then set $\tilde{u}_n := u_n(e^{-\theta_n} \cdot)$ we have
$$ \mc{P}(\lambda_n, \tilde{u}_n)\to 0,$$
$$\partial_{\lambda}\mc{I}^m(\lambda_n, \tilde{u}_n) \to 0,$$
$$ \norm{\partial_u \mc{I}^m(\lambda_n,\tilde{u}_n)}_{(H^s_r(\R^N))^*}\to 0,$$
and thus by Proposition \ref{PSP} the sequence $(\lambda_n, \tilde{u}_n)$ is convergent (up to subsequences) to a $(\lambda, \tilde{u})\in K^{PSP}_b$. 
Observe that, for each $n$, set $v_n : =\tilde{u}(e^{\theta_n}\cdot)$, we have $(\theta_n, \lambda, v_n)\in \tilde{K}_b$. 
Therefore by \eqref{eq_shift_distance}
\begin{align*}
	\dist_M((\theta_n, \lambda_n, u_n), \tilde{K}_b) &\leq \dist_M((\theta_n, \lambda_n, u_n), (\theta_n, \lambda, v_n)) \\
	&= \dist_M((0,\lambda_n, \tilde{u}_n), (0, \lambda, \tilde{u})) \\
	&\leq \sqrt{|\lambda_n -\lambda|^2 + \norm{\tilde{u}_n-\tilde{u}}_{H^s_r(\R^N)}^2} \to 0,
\end{align*}
which reaches the claim.
\QED

\bigskip

\noindent
{\bf Notation.}
We use the following notation: for $\tilde A\subset M$ and $\rho>0$ we set
	$$ \tilde N_\rho(\tilde A) := \{ (\theta,\lambda,u)\in M \mid \dist_M((\theta,\lambda,u),\tilde A) <\rho\}, $$
while for $A\subset\R\times H_r^s(\R^N)$ and $R>0$ we set
	$$ N_R(A) := \{(\lambda,u)\in \R\times H_r^s(\R^N) \mid d((\lambda,u),A)<R\}, $$
where
	$$ d\big((\lambda,u),(\lambda',u')\big) := (|\lambda-\lambda'|^2 +\norm{u-u'}_{H_r^s}^2)^{1/2}. $$
We also write for $a<b$ %$-\infty<a<b<\infty$
\begin{align*}
% \mc{I}^b 
[\mc{I}^m\leq b]&:= \{ (\lambda,u)\in\R\times H_r^s(\R^N) \mid \mc{I}(\lambda,u) \leq b\}, \\
% \mc{I}^b_a 
[a\leq \mc{I}^m\leq b] &:= \{ (\lambda,u)\in\R\times H_r^s(\R^N) \mid a\leq \mc{I}(\lambda,u) \leq b\}, \\
 [\mc{H}^m\leq b]_M%\mc{H}^b
 &:= \{ (\theta,\lambda,u)\in M \mid \mc{H}(\theta,\lambda,u) \leq b\}, \\
[a\leq \mc{H}^m\leq b]_M%\mc{H}^b_a 
&:= \{ (\theta,\lambda,u)\in M \mid a\leq \mc{H}(\theta,\lambda,u) \leq b\}.
\end{align*}
Using these notations, as a corollary to Proposition \ref{PSPtilde}, we have
\begin{Corollary}\label{cor_PStilde}
For any $\rho>0$ there exists a $\delta_\rho>0$ such that
		\begin{equation}\label{eq_curv_pend_gen}
		\forall \, (\theta,\lambda,u) \in [b-\delta_{\rho} \leq \mc{H}^m \leq b+\delta_{\rho}]_M %\mc{H}^{b+\delta_{\rho}}_{b-\delta_{\rho}}
\setminus \tilde N_\rho(\tilde K_b)
 \; : \; \norm{D\mc{H}(\theta,\lambda,u)}_{(\theta,\lambda,u),*}> \delta_{\rho}.
		\end{equation}
Here, if $\tilde{K}_b=\emptyset$, we regard $\tilde N_\rho(\tilde K_b)=\emptyset$.
\end{Corollary}

%%%%%%%%%%%%%%%%%%%%%%%%%%%%%%%%%%%%%%%%%%%%%%%%%%%%%%%

%\setcounter{equation}{0} %FOR ARXIV
\section{A deformation flow by projections}
\label{section:5}

Exploiting an idea in \cite{HT0} (see also \cite{IT0}), we aim to prove the following Deformation Theorem in the fractional framework.

\begin{Theorem}\label{defarg}
Let $b<0$, and assume $K_b^{PSP}= \emptyset$. 
Let $\bar{\varepsilon}>0$, then there exist $\varepsilon \in (0,\bar{\varepsilon})$ and $\eta: [0,1]\times (\R\times H^s_r(\R^N))\to \R\times H^s_r(\R^N)$ continuous such that
\begin{enumerate}
\item $\eta(0, \cdot,\cdot)=id_{\R\times H^s_r(\R^N)}$;
\item $\eta$ fixes $%\mc{I}^{b-\bar{\varepsilon}}
[\mc{I}^m\leq b-\bar{\eps}]$, that is, $\eta(t, \cdot,\cdot)=id_{%\mc{I}^{b-\bar{\varepsilon}}
[\mc{I}^m\leq b-\bar{\eps}]}$
for all $t\in [0,1]$;
\item $\mc{I}^m$ is non-increasing along $\eta$, and in particular $\mc{I}^m(\eta(t,\cdot, \cdot))\leq \mc{I}^m(\cdot, \cdot)$ for all $t \in [0,1]$;
\item $\eta(1, [\mc{I}^m\leq b+\eps]%\mc{I}^{b+\varepsilon}
)\subseteq [\mc{I}^m\leq b-\eps]%\mc{I}^{b-\varepsilon}
$.
\end{enumerate}

\end{Theorem}
We omit the proof of the Theorem since it will be very similar to the one made in the case of multiplicity (see Theorem \ref{thm_def_gen}). 
We remark that this deformation flow is not $C^1$ and it does not satisfy the two properties of the standard deformation flows, in general \cite[Remark 3.2]{HT0}:
\begin{itemize}
	\item[(1)] $\eta(s+t, \lambda, u)= \eta(t, \eta(s,\lambda, u))$ with $s+t \in [0,1], (\lambda,u) \in \R \times H^s_r(\R^N)$;	
	\item[(2)] for $t \in [0,1]$, the map $(\lambda,u) \mapsto \eta(t,\lambda, u)$ is a homeomorphism.
\end{itemize}
This is due to the fact that this deformation will be built through a projection of another deformation, built for the augmented functional $\mc{H}^m$.
%We refer to .

We also stress that the deformation argument in Theorem \ref{defarg} works for $K^{PSP}_b$ but not for $K_b$ and thus, if $K^{PSP}_b=\emptyset$, then we have the statement (4) in Theorem \ref{defarg} even if $K_b\not=\emptyset$. 
By classical arguments, we derive the following existence theorem (see also the proof of Corollary \ref{coroll_esist_Pm}).

\begin{Corollary}[Existence]\label{dedu}
	Let $\bar b <0$ be a MP minimax value for $\mc{I}^m$. 
	Then $K^{PSP}_{\bar b} \not=\emptyset$, that is, $\mc{I}^m$ has a critical point $(\bar \lambda, \bar u)$ satisfying the Pohozaev identity, namely $\mc{P}(\bar \lambda, \bar u)=0$.
\end{Corollary}

%%%%%%%%%%%%%%%%%%%%%%%%%%%%%%%%%%%%%%%%%%%%%%%%%%%%%%%

%\setcounter{equation}{0} %FOR ARXIV
\section{Minimax critical points in the product space}
\label{section:6}

For any $m >0$, let $B_m$ and $E_m$ be the constants defined in Proposition \ref{negativo}, namely
$$ B_m = \inf_{\lambda<\lambda_0} \left(a(\lambda) -\frac{e^{\lambda}}{2}m \right), \quad E_m=\inf_{(\lambda,u)\in\partial\Omega}\mc{I}^m(\lambda,u).$$
As a minimax class for $\mc{I}^m$, we define the paths going from $\Omega$ to $\Omega^c$, such that the energy of the ending points is below the minimal energy on the mountain $\partial \Omega$:
\begin{eqnarray*}
	\Gamma^m := \big\{\xi\in C\big([0,1], \R \times H^s_r(\R^N)\big) & \mid 
	& \xi(0) \in \R \times \{ 0\},\ \mc{I}^m(\xi(0)) \leq B_m -1,\\ &&
	\xi(1) \not \in \Omega,		
	\ \mc{I}^m(\xi(1)) \leq B_m -1\big\}. 
\end{eqnarray*}

We have the following result.

\begin{Proposition}\label{tom} 
Assume \hyperref[(g1)]{\textnormal{(g1)}}--\hyperref[(g3)]{\textnormal{(g3)}}.
	\begin{itemize}
	\item[\textnormal{(i)}]
				 For any $m>0$, we have $\Gamma^m \neq \emptyset$.
	\item[\textnormal{(ii)}]	
		For sufficiently large $m>0$ there exists $\xi\in \Gamma^m$ such that
		\begin{equation}\label{terza}
		\max_{t\in [0,1]} \mc{I}^m(\xi(t)) <0.
		\end{equation}
	\item[\textnormal{(iii)}]	
		Assume \hyperref[(g4)]{\textnormal{(g4)}}. 
		Then for any $m>0$ there exists $\xi\in \Gamma^m$ with the property \eqref{terza}.
	\end{itemize}
\end{Proposition}

\claim Proof.
Let $\lambda_0\in (-\infty,\infty]$ be defined in \eqref{lambdazero}. 
For any $\lambda<\lambda_0$ we show there exists a path $\psi_\lambda\in\Gamma^m$ such that
\begin{equation} \label{aa}
	\max_{t\in [0,1]} \mc{I}^m(\psi_\lambda(t)) \leq a(\lambda)-\frac{e^{\lambda}}{2}m.
\end{equation}
Let $u_\lambda$ be a MP solution of $\partial_u\mc{J}(\lambda,u)=0$ (by Theorem \ref{Byeonminimizers}). 
Set $\zeta_\lambda(t):=u_\lambda(\cdot/t)$ for $t>0$ and $\zeta_\lambda(0):=0$ and note that, since $u_{\lambda}$ satisfies the Pohozaev identity, we have $\mc{I}^m (\lambda,\zeta_\lambda(t)) \to-\infty$ and $\mc{P}(\lambda,\zeta_\lambda(t)) \to-\infty$ as $t\to+\infty$. 
We can find $\gamma_\lambda:=\zeta_\lambda(L\cdot)$ for $L\gg 1$ satisfying
$$ a(\lambda) = \max_{t\in [0,1]} \mc{J}(\lambda,\gamma_\lambda(t)),$$
$$ \mc{I}^m(\lambda,\gamma_\lambda(1)) \leq B_m-1, \quad (\lambda,\gamma_{\lambda}(1))\notin \Omega.$$
%Indeed, by Pohozaev, $\mc{I}^m(\lambda, \gamma_{\lambda}(t))=\norm{(-\Delta)^{s/2} u_{\lambda}}_2^2 \left( \frac{1}{2} (tL)^{N-2s} - \frac{1}{2^*_s} (tL)^N\right)$ which is maximized for $t=1/L$; this corresponds to have $\mc{I}^m(\lambda, u_{\lambda})=a(\lambda)$.
We also note that $t\mapsto \mc{I}^m(t,0)=-\frac{e^t}{2}m$ is decreasing and tending to $-\infty$
as $t\to+\infty$. Thus, joining $\gamma_\lambda$ and $t\mapsto (\lambda+Lt,0);\, [0,1]\to\R\times H_r^s(\R^N)$
for $L\gg 1$, we find a path $\psi_\lambda\in\Gamma^m$, defined as
$$	\psi_{\lambda}(t) := \parag{
		(\lambda+L(1-2t), \,0) & \quad \hbox{ if $t \in [0,1/2]$,} \\ 
		\left(\lambda, \, \gamma_{\lambda}(2t-1) \right) & \quad \hbox{ if $t \in (1/2,1]$}
		}
$$
 with \eqref{aa}.
Thus in particular we have (i).

Next we deal with (ii) and (iii). By \eqref{aa}, we have that (ii) follows easily; (iii) also follows from Proposition \ref{out}.
\QED

\bigskip

We notice that each path in $\Gamma^m$ passes through $\partial \Omega$, thus the minimax value
\begin{equation}\label{crti}
b_m := \inf_{\xi\in\Gamma^m}\max_{t\in [0,1]} \mc{I}^m(\xi(t))
\end{equation} 
verifies $b_m\geq E_m$ and hence by Proposition \ref{negativo} it is well defined and finite. 
Since the Palais-Smale-Pohozaev condition holds on $(-\infty,0)$, it is important to estimate $b_m$. 
We have the following result.
\begin{Proposition}\label{blue}
	 Assume \hyperref[(g1)]{\textnormal{(g1)}}--\hyperref[(g3)]{\textnormal{(g3)}}.
We have
\begin{equation}\label{mor}
 b_m \leq a(\lambda)-\frac{e^{\lambda}}{2}m 
 %&= e^\lambda\left(\frac{a(\lambda)}{e^{\lambda}}-\frac{m}{2}\right)
 \quad \hbox{for all}\ \lambda<\lambda_0.
\end{equation}
Moreover
	\begin{itemize}
		\item[\textnormal{(i)}]
		%There exists $m_0>0$ such that
Setting 
$$ m_0 := 2 \inf_{\lambda<\lambda_0}\frac{a(\lambda)}{e^{\lambda}}\geq 0, $$
we have
			$$ b_m<0 \quad \hbox{for $m>m_0$}. $$
		\item[\textnormal{(ii)}] Assume \hyperref[(g4)]{\textnormal{(g4)}} in addition, then $m_0=0$, that is, 
			$$ b_m <0 \quad \hbox{for all $m>0$}.$$
	\item[\textnormal{(iii)}]	
		We have $b_m=E_m=B_m$.
\item[\textnormal{(iv)}]
$\limsup_{m\to+\infty} \frac{b_m}{m} \leq -\frac{e^{\lambda_0}}{2}$. If $\lambda_0=+\infty$, then $\lim_{m\to+\infty} \frac{b_m}{m} = -\infty$ (see \cite{CGT2}).
%We have $\lim_{m\to+\infty} \frac{b_m}{m} = -\infty$ (see \cite{CGT2}); in particular, 
%$$b_m \to - \infty \quad \hbox{ as $m \to +\infty$}.$$
	\end{itemize}
\end{Proposition}

\claim Proof.
By \eqref{aa} we have \eqref{mor}, and thus %, we have
$$ b_m \leq e^\lambda\left(\frac{a(\lambda)}{e^{\lambda}}-\frac{m}{2}\right)
 \quad \hbox{for all}\ \lambda<\lambda_0.
 $$
%\begin{align}\label{mor}
% b_m &\leq a(\lambda)-\frac{e^{\lambda}}{2}m \nonumber \\
% &= e^\lambda\left(\frac{a(\lambda)}{e^{\lambda}}-\frac{m}{2}\right)
% \quad \hbox{for all}\ \lambda<\lambda_0.
%\end{align}
By definition of $m_0$, 
%Setting
%$$ m_0 := 2 \inf_{\lambda<\lambda_0}\frac{a(\lambda)}{e^{\lambda}}\geq 0, $$
we have $b_m<0$ for $m>m_0$. 
Thus we have (i). 
By Proposition \ref{out}, we have $m_0=0$ under the assumption \hyperref[(g4)]{\textnormal{(g4)}} and thus we have (ii).

Furthermore, from \eqref{mor} it follows $b_m \leq B_m$. As already observed $b_m \geq E_m \geq B_m$, from which we deduce (iii).

Finally for any $\lambda \in \R$ we have, again by \eqref{aa}, 
$$\limsup_{m\to+\infty} \frac{b_m}{m} \leq 
\lim_{m\to+\infty} \left(\frac{a(\lambda)}{m} - \frac{e^{\lambda}}{2}\right) =- \frac{e^{\lambda}}{2}.
$$
Since $\lambda$ is arbitrary, we get (iv).
\QED

\bigskip

By Proposition \ref{blue} and Corollary \ref{dedu} we conclude that the level $b_m$, defined in \eqref{crti}, is a critical value of $\mc{I}^m$ in the product space $\R \times H^s_r(\R^N)$ and thus Theorem \ref{S:1.1_frac} and Theorem \ref{S:1.12_frac} hold.

\begin{Corollary}\label{coroll_esist_Pm}
Let $m>m_0$.
 Then there exists a solution of problem \eqref{problem_frac} which satisfies the Pohozaev identity \eqref{eq_Pohozaev_fractional}.
If moreover \hyperref[(g4)]{\textnormal{(g4)}} holds, then there exists a solution of \eqref{problem_frac} for each $m>0$.
\end{Corollary}

\claim Proof.
Let $\bar{\varepsilon}\in (0, 1)$. By Theorem \ref{defarg}, in correspondence to $b_m<0$, there exists $\varepsilon\in (0,\bar{\varepsilon})$ and $\eta$ satisfying $1)-4)$. By definition of inf, there exists $\gamma \in \Gamma^m$ such that
$$\max_{t\in [0,1]} \mc{I}^m(\gamma(t)) < b_m + \varepsilon,$$
that is 
\begin{equation}\label{d1}
\gamma([0,1])\subseteq [\mc{I}^m\leq b_m+\eps]% \mc{I}^{b_m+\varepsilon}
.
\end{equation}
Set
$$\tilde{\gamma}(t):=\eta(1,\gamma(t)),$$
we show that $\tilde{\gamma}\in \Gamma^m$. Indeed for $i\in\{0,1\}$, since $\mc{I}^m(\gamma(i))\leq B_m - 1 \leq b_m-\bar{\varepsilon}$, 
Theorem \ref{defarg} implies that $\tilde{\gamma}(i)=\eta(1,\gamma(i))=\gamma(i)\in [\mc{I}^m \leq b_m + \bar{\eps}]%\mc{I}^{b_m-\bar{\varepsilon}}
$, and thus $\tilde{\gamma}(0)=\gamma(0)\in\R\times\{ 0\}$, $\tilde{\gamma}(1)=\gamma(1)\not\in\Omega$.
Therefore 
\begin{equation}\label{d2}
b_m \leq \max_{t\in [0,1]} \mc{I}^m(\tilde{\gamma}(t)).
\end{equation}
By contradiction, assume $K_{b_m}^{PSP} = \emptyset$. By the properties of $\eta$ and \eqref{d1} we obtain that $\tilde\gamma([0,1])=\eta(1,\gamma([0,1]))\subseteq [\mc{I}^m\leq b_m - \eps]%\mc{I}^{b_m-\varepsilon}
$, that is
$$\max_{t\in [0,1]} \mc{I}^m(\eta(1,\gamma(t)))\leq b_m-\varepsilon.$$
This is in contradiction with \eqref{d2}, and we conclude the proof.
\QED
\medskip

\begin{Remark}
We observe that, by Proposition \ref{blue} (iii), %we obtain that 
the found Mountain Pass solution $(\overline{\mu}, \overline{u})$ at level $b_m$ is a Pohozaev minimum on the product space $\R \times H^s_r(\R^N)$. 
This additionally implies that the found solution is a Pohozaev minimum for the unconstrained case, once fixed $\overline{\mu}$; see also Remark \ref{rem_constrain}.
\end{Remark}

%%%%%%%%%%%%%%%%%%%%%%%%%%%%%%%%%%%%%%%%%%%%%%%%%%%%%%%

%\setcounter{equation}{0} %FOR ARXIV
\section{Multiple normalized solutions}
%\label{section:7}
\label{sec_frac_multiple_sol}

We focus now on the existence of multiple solutions. In the whole Section we assume, in addition, \hyperref[(g5)]{\textnormal{(g5)}}.

%%%%%%%%%%%%%%%%%%%%%%%%%%%%%%%%%%%%%

\subsection{Symmetric deformation theorems}

In what follows we will use the following terminology. Consider the action $\sigma$ of $\G:=\Z_2$ on the last components of $\R\times H^s_r(\R^N)$ and $M =\R \times \R \times H^s_r(\R^N)$, that is
$$\sigma:\,(\pm1,\lambda, u) \in \G\times (\R\times H^s_r(\R^N)) \mapsto (\lambda, \pm u) \in \R\times H^s_r(\R^N),$$
$$\sigma:\,(\pm1,\theta, \lambda, u) \in \G\times M \mapsto (\theta, \lambda, \pm u) \in M.$$
We notice that $\mc{I}^m$ and $\mc{H}^m$ are invariant under this action (i.e. they are even in $u$), as well as the set $\Omega$ (i.e. it is symmetric with respect the axis $\R$). In particular this means that, if $u$ is a solution, then $-u$ is a solution as well. 
We highlight instead that the function $\eta=(\eta_1, \eta_2):\R\times H^s_r(\R^N) \to \R\times H^s_r(\R^N)$ (resp. $\tilde{\eta}=(\tilde{\eta}_0,\tilde{\eta}_1,\tilde{\eta}_2):M \to M$) is equivariant if $\eta_1$ is even in $u$ and $\eta_2$ is odd in $u$ (resp. $\tilde{\eta}_0$ and $\tilde{\eta}_1$ are even and $\tilde{\eta}_2$ is odd).
%\tbl{
%We recall that a function $f$ is said to be invariant under the action $g \cdot x$ if $f(g \cdot x) = f(x)$, while it is equivariant if $f(g \cdot x) = g \cdot f(x)$.
%We want to prove the following.
%}

%\enlargethispage{2\baselineskip}%€

\begin{Theorem}\label{thm_def_gen}
Let $b<0$, and let $\mc{O}$ be a neighborhood of $K_b^{PSP}$. 
Then for each $\bar{\varepsilon}>0$ there exist $\varepsilon \in (0,\bar{\varepsilon})$ and $\eta: [0,1]\times (\R\times H^s_r(\R^N))\to (\R\times H^s_r(\R^N))$ continuous such that
\begin{enumerate}
\item $\eta(0, \cdot,\cdot)=id_{\R\times H^s_r(\R^N)}$;
\item $\eta$ fixes $[\mc{I}^m\leq b-\bar{\eps}]%\mc{I}^{b-\bar{\varepsilon}}
$, that is, $\eta(t, \cdot,\cdot)=id_{[\mc{I}^m\leq b-\bar{\eps}]%\mc{I}^{b-\bar{\varepsilon}}
}$ for all $t \in [0,1]$;
\item $\mc{I}^m$ is non-increasing along $\eta$, and in particular $\mc{I}^m(\eta(t,\cdot, \cdot))\leq \mc{I}^m(\cdot, \cdot)$ for all $t \in [0,1]$;
\item if $K_b^{PSP}= \emptyset$, then $\eta(1, [\mc{I}^m\leq b+ \eps]%\mc{I}^{b+\varepsilon}
)\subseteq [\mc{I}^m\leq b-\eps]%\mc{I}^{b-\varepsilon}
$;
\item if $K_b^{PSP}\not=\emptyset$, then
$$\eta(1,[\mc{I}^m\leq b+\eps]%\mc{I}^{b+\varepsilon}
\setminus \mc{O}) \subseteq [\mc{I}^m\leq b-\eps]% \mc{I}^{b-\varepsilon}
$$
and
$$\eta(1, [\mc{I}^m\leq b+\eps]%\mc{I}^{b+\varepsilon}
) \subseteq [\mc{I}^m\leq b-\eps] %\mc{I}^{b-\varepsilon}
\cup \mc{O};$$
\item $\eta(t, \cdot,\cdot)$ is $\G$-equivariant, in the sense mentioned before.
\end{enumerate}
\end{Theorem}

To prove this, we work first on the functional $\mc{H}^m$, for which we obtained the $(\widetilde{PSP})$ condition.

\begin{Theorem}\label{thm_def_H}
Let $b<0$, $\rho>0$ and write $\tilde{\mc{O}}:=\tilde{N}_{\rho}(\tilde{K}_b)$. 
Then for each $\bar{\varepsilon}>0$ there exist $\varepsilon \in (0,\bar{\varepsilon})$ and $\tilde{\eta}: [0,1]\times M\to M$ continuous such that
\begin{enumerate}
\item $\tilde{\eta}(0, \cdot,\cdot)=id_M$;
\item $\tilde{\eta}$ fixes $[\mc{H}^m\leq b-\bar{\eps}]_M%\mc{H}^{b-\bar{\varepsilon}}
$, that is $\tilde{\eta}(t, \cdot,\cdot)=id_{[\mc{H}^m\leq b-\bar{\eps}]_M%\mc{H}^{b-\bar{\varepsilon}}
}$ for all $t\in [0,1]$;
\item $\mc{H}^m$ is non-increasing along $\tilde{\eta}$, and in particular $\mc{H}^m(\tilde{\eta}(t,\cdot,\cdot, \cdot))\leq \mc{H}^m(\cdot, \cdot, \cdot)$ for all $t\in [0,1]$; 
\item if $\tilde{K}_b= \emptyset$, then $\tilde{\eta}(1, [\mc{H}^m\leq b+\eps]_M%\mc{H}^{b+\varepsilon}
)\subseteq [\mc{H}^m\leq b-\eps]_M%\mc{H}^{b-\varepsilon}
$;
\item if $\tilde{K}_b \neq \emptyset$, then
$$\tilde{\eta}(1,[\mc{H}^m\leq b+\eps]_M%\mc{H}^{b+\varepsilon}
\setminus \tilde{\mc{O}}) \subseteq [\mc{H}^m \leq b- \eps]_M%\mc{H}^{b-\varepsilon}
$$
and
$$\tilde{\eta}(1,\mc{H}^{b+\varepsilon}) \subseteq \mc{H}^{b-\varepsilon}\cup \tilde{\mc{O}};$$
\item $\tilde{\eta}(t, \cdot,\cdot)$ is $\G$-equivariant, in the sense mentioned before.
\end{enumerate}
\end{Theorem}

We postpone the proof of Theorem \ref{thm_def_H} for $\mc{H}^m$ and see now how to use it to deduce the one for $\mc{I}^m$. Introduce first the following notation:
$$\pi: M \to \R\times H^s_r(\R^N), \; \pi(\theta, \lambda, u) := (\lambda, u(e^{-\theta}\cdot)),$$
$$\iota: \R\times H^s_r(\R^N) \to M, \; \iota(\lambda, u) := (0,\lambda, u),$$
which are a kind of rescaling projection and immersion. Observe that
$$\pi \circ \iota = id_{\R\times H^s_r(\R^N)}, \quad \hbox{(while $\iota \circ \pi \neq id_M$),}$$
$$\mc{H}^m\circ \iota=\mc{I}^m, \quad \mc{I}^m\circ \pi=\mc{H}^m,$$
$$\pi(\tilde{K}_b) = K_b^{PSP}.$$
For $\tilde{\eta}$ obtained in Theorem \ref{thm_def_H}, define "$\eta=\pi \circ \tilde{\eta} \circ \iota$" up to the time; more precisely
\begin{equation}\label{eq_def_flow}
\eta(t,\lambda, u):=\pi(\tilde{\eta}(t,\iota(\lambda, u))).
\end{equation}
It is now a straightforward computation showing that $\eta$ satisfies the requests of Theorem \ref{thm_def_gen}. 
A delicate issue, anyway, is to show the intuitive fact that neighborhoods of $\tilde{K}_b$ are brought to neighborhoods of $K_b^{PSP}$. More precisely we have the following result.

\begin{Lemma}\label{lem_intorni}
Assume that $K_b^{PSP}$ is compact (for instance, $b<0$). Let $\rho>0$, then there exists $R(\rho)>0$ such that, set $\tilde{\mc{O}}:=\tilde{N}_{\rho}(\tilde{K}_b)$ and $\mc{O}:=N_{R(\rho)}(K_b^{PSP})$, we have
$$\pi(\tilde{\mc{O}}) \subset \mc{O},$$
i.e.
$$\dist_M((\theta, \lambda, u), \tilde{K}_b) \leq \rho \implies d((\lambda, u(e^{-\theta}\cdot)), K_b^{PSP}) \leq R(\rho).$$
In particular, for $\theta =0$ we have
\begin{equation}\label{eq_theta0}
\dist_M((0, \lambda, u), \tilde{K}_b) \leq \rho \implies d((\lambda, u), K_b^{PSP}) \leq R(\rho),
\end{equation}
that is
$$\iota\big(\complement \mc{O}\big) \subseteq \complement \tilde{\mc{O}}$$
where $\complement$ denotes the complement of the set. Moreover
$$\lim_{\rho\to 0} R(\rho)=0.$$
\end{Lemma} 

\claim Proof.
We observe that is sufficient to prove \eqref{eq_theta0} since by \eqref{eq_shift_distance}
$$\dist_M((\theta,\lambda, u), \tilde{K}_b) = \dist_M((0,\lambda, u(e^{-\theta}\cdot)), \tilde{K}_b).$$
Let $\varepsilon>0$. By definition of $\dist_M((0,\lambda,u), \tilde K_b)$ there exists a $\sigma=\sigma(t)$, $\sigma=(\theta, \lambda, u)$, such that $\sigma(0)=(0,\lambda, u)$, $\sigma(1)\in \tilde{K}_b$ and
\begin{equation}\label{eq_dim_def_1}
\int_0^1 \norm{\dot{\sigma}(t)}_{\sigma(t)} dt \leq \rho + \varepsilon.
\end{equation}
By \eqref{eq_confronto_K} we have $(\lambda(1), u(1)(e^{-\theta(1)} \cdot))\in K_b^{PSP}$ and thus
\begin{eqnarray*}
\lefteqn{ \dist((\lambda, u), K_b^{PSP})}
\\ &\leq& \norm{(\lambda, u) - (\lambda(1), u(1)(e^{-\theta(1)} \cdot))}_{\R\times H^s_r(\R^N)}\\
&\leq & \norm{(\lambda, u) - (\lambda(1), u(1))}_{\R\times H^s_r(\R^N)}+ \norm{(\lambda(1), u(1)) - (\lambda(1), u(1)(e^{-\theta(1)} \cdot))}_{\R\times H^s_r(\R^N)}\\
&= & \norm{(\lambda(0), u(0)) - (\lambda(1), u(1))}_{\R\times H^s_r(\R^N)} + \norm{u(1) - u(1)(e^{-\theta(1)} \cdot)}_{H^s_r(\R^N)}\\
&=& I + II.
\end{eqnarray*}
Focus on $I$. We have, by the fundamental theorem of calculus and H\"older inequality,
\begin{align*}
I&=\norm{(\lambda(0), u(0)) - (\lambda(1), u(1))}_{\R\times H^s_r(\R^N)} \leq \int_0^1 \left(\abs{\dot{\lambda}(t)}^2 + \norm{\dot{u}(t)}_{H^s_r(\R^N)}^2\right)^{1/2} dt \\
&= \int_0^1 \left(\abs{\dot{\lambda}(t)}^2 + \norm{\dot{u}(t)}_{2}^2 + \norm{(-\Delta)^{s/2}\dot{u}(t)}_{2}^2\right)^{1/2} dt. 
\end{align*}
In order to use \eqref{eq_dim_def_1} it must appear the norm associated to $M$, which we recall is
$$ \norm{\dot{\sigma}(t)}_{\sigma(t)}^2= \dot{\theta}(t)^2 + \dot{\lambda}(t)^2 + e^{N\theta(t)} \norm{\dot{u}(t)}_2^2 + e^{(N-2s)\theta(t)} \norm{(-\Delta)^{s/2} \dot{u}(t)}_2^2.$$
Since we do not know the sign of $N\theta(t)$, 
we need an estimate on $\theta(t)$ and a corrective factor. Indeed, recalled that $\theta(0)=0$, we have
$$|\theta(t)| = |\theta(t)-\theta(0)| \leq \int_0^1 |\dot{\theta}(t)| dt \leq \int_0^1 \norm{\dot{\sigma}(t)}_{\sigma(t)} dt \leq \rho + \varepsilon.$$
Thus $\theta(t) \geq -(\rho + \varepsilon) \geq -\frac{N}{N-2s} (\rho + \varepsilon)$ which imply
$$e^{N(\rho+\varepsilon)} \geq 1, \quad e^{N(\rho+\varepsilon)} e^{N\theta(t)}\geq 1, \quad e^{N(\rho+\varepsilon)} e^{(N-2s)\theta(t)}\geq 1 $$
and hence we obtain
\begin{align*}
I& \leq e^{\frac{N(\rho+\varepsilon)}{2}} \int_0^1 \left(\abs{\dot{\lambda}(t)}^2 + e^{N\theta(t)} \norm{\dot{u}(t)}_{2}^2 + e^{(N-2s)\theta(t)} \norm{(-\Delta)^{s/2}\dot{u}(t)}_{2}^2\right)^{1/2} dt \\
& \leq e^{\frac{N(\rho+\varepsilon)}{2}} \int_0^1 \left(\abs{\dot{\theta}(t)^2} + \abs{\dot{\lambda}(t)}^2 + e^{N\theta(t)} \norm{\dot{u}(t)}_{2}^2 + e^{(N-2s)\theta(t)} \norm{(-\Delta)^{s/2}\dot{u}(t)}_{2}^2\right)^{1/2} dt \\
&= e^{\frac{N(\rho+\varepsilon)}{2}} \int_0^1 \norm{\dot{\sigma}(t)}_{\sigma(t)} dt 
\leq e^{\frac{N(\rho+\varepsilon)}{2}} (\rho + \varepsilon) \stackrel{\varepsilon \to 0} \to e^{\frac{N\rho}{2}} \rho.
\end{align*}
Focus now on $II$. Set $\bar{\omega}:= u(1)(e^{-\theta(1)}\cdot)$ we have $\bar{\omega} \in P_2(K_b^{PSP})$ (where $P_2$ is the projection on the second component) with $|\theta(1)|\leq \rho + \varepsilon$, and thus
\begin{align*}
II &= \norm{u(1) - u(1)(e^{-\theta(1)} \cdot)}_{H^s_r(\R^N)} = \norm{\bar{\omega}(e^{\theta(1)}\cdot) - \bar{\omega}}_{H^s_r(\R^N)}\\
&\leq \sup \left \{ \norm{\omega(e^{\alpha}\cdot)- \omega}_{H^s_r(\R^N)} \mid |\alpha| \leq \rho + \varepsilon, \; \omega \in P_2(K_b^{PSP}) \right \}.
\end{align*}
Since $P_2(K_b^{PSP})$ is compact, it is simple to show that, as $\varepsilon \to 0$,
$$II\leq \sup \left \{ \norm{\omega(e^{\alpha}\cdot)- \omega}_{H^s_r(\R^N)} \mid |\alpha| \leq \rho, \; \omega \in P_2(K_b^{PSP}) \right \}.$$ 
Summing up, we have
\begin{align*}
\dist((\lambda, u), K_b^{PSP})&\leq e^{\frac{N\rho}{2}} \rho + \sup \left \{ \norm{\omega(e^{\alpha}\cdot)- \omega}_{H^s_r(\R^N)} \mid \abs{\alpha} \leq \rho, \; \omega \in P_2(K_b^{PSP}) \right \} \\
&=: R(\rho)<\infty.
\end{align*}
Here we have
$$\lim_{\rho\to 0} R(\rho) =0,$$
which concludes the proof.
\QED

\bigskip

We are now ready to show that $\eta$ satisfies the desired properties.

\medskip

\claim Proof of Theorem \ref{thm_def_gen}.
Let $\mc{O}$ be a neighborhood of $K_b^{PSP}$, and choose $R$ such that $N_R(K_b^{PSP})\subset \mc{O}$. 
By Lemma \ref{lem_intorni} choose $\rho\ll 1$ satisfying $R(\rho)<R$ and thus $N_{R(\rho)}(K_b^{PSP}) \subset \mc{O}$. 
Consequently, by Theorem \ref{thm_def_H}, there exists a deformation $\tilde{\eta}$ corresponding to the neighborhood $\tilde{\mc{O}}:=\tilde{N}_{\rho}(\tilde{K}_b)$. 
We thus define $\eta$ by \eqref{eq_def_flow} and prove the properties. 
Start observing that
%$$(\lambda, u)\in [\mc{I}^m \leq b \pm \delta]%\mc{I}^{b\pm \delta} 
%\implies b\pm \delta > \mc{I}^m(\lambda, u) = \mc{H}^m(\iota(\lambda, u)) \implies \iota(\lambda, u)\in [\mc{H}^m\leq b \pm \delta]_M%\mc{H}^{b\pm \delta}
%,$$
\begin{align*}
(\lambda, u)\in [\mc{I}^m \leq b \pm \delta]%\mc{I}^{b\pm \delta} 
& \implies b\pm \delta > \mc{I}^m(\lambda, u) = \mc{H}^m(\iota(\lambda, u)) \\
&\implies \iota(\lambda, u)\in [\mc{H}^m\leq b \pm \delta]_M%\mc{H}^{b\pm \delta}
,
\end{align*}
i.e. $\iota( [\mc{I}^m \leq b \pm \delta]%\mc{I}^{b\pm \delta}
) \subset [\mc{H}^m\leq b \pm \delta]_M%\mc{H}^{b\pm \delta}
$; similarly, $\pi([\mc{H}^m \leq b \pm \delta]_M%\mc{H}^{b\pm \delta}
)\subset [\mc{I}^m\leq b \pm \delta]%\mc{I}^{b\pm \delta}
$.
\begin{enumerate}
\item $\eta(0,\lambda, u) = \pi(\tilde{\eta}(0,\iota(\lambda,u))) = \pi(\iota(\lambda,u)) = (\lambda, u)$.
\item If $(\lambda,u)\in [\mc{I}^m \leq b- \bar{\eps}]%\mc{I}^{b-\bar\varepsilon}
$, then $\iota(\lambda,u)\in [\mc{H}^m \leq b-\bar{\eps}]_M%\mc{H}^{b-\bar\varepsilon}
$.
Thus $\eta(t,\lambda, u) = \pi(\tilde{\eta}(t,\iota(\lambda,u))) = \pi(\iota(\lambda, u))= (\lambda, u)$.
\item $\mc{I}^m(\eta(t,\lambda, u)) = \mc{I}^m(\pi(\tilde{\eta}(t,\iota(\lambda,u)))) = \mc{H}^m(\tilde{\eta}(t,\iota(\lambda,u))) \leq \mc{H}^m(\iota(\lambda,u)) = \mc{I}^m(\lambda, u)$.
\item If $K_b^{PSP}=\emptyset$, then $\tilde{K}_b=\emptyset$. 
Thus for $(\lambda,u)\in [\mc{I}^m\leq b+\eps]%\mc{I}^{b+\varepsilon}
$, we have
$\mc{I}^m(\eta(1,\lambda, u)) = \mc{I}^m(\pi(\tilde{\eta}(1,\iota(\lambda,u)))) = \mc{H}^m(\tilde{\eta}(1,\iota(\lambda,u))) \leq b-\varepsilon$.
\item We have, by previous arguments and Lemma \ref{lem_intorni}, that $\iota([\mc{I}^m\leq b+\eps]%\mc{I}^{b+\varepsilon}
\setminus \mc{O}) 
= \iota([\mc{I}^m\leq b+\eps]%\mc{I}^{b+\varepsilon}
\cap \complement \mc{O}) 
\subseteq \iota([\mc{I}^m\leq b+\eps]%\mc{I}^{b+\varepsilon}
) \cap \iota(\complement \mc{O}) 
\subseteq [\mc{H}^m\leq b+\eps]_M%\mc{H}^{b+\varepsilon}
 \cap (\complement \tilde{\mc{O}})
= [\mc{H}^m\leq b+ \eps]_M%\mc{H}^{b+\varepsilon}
\setminus \tilde{\mc{O}}$ and thus
% $$\eta(1,[\mc{I}^m\leq b+\eps]%\mc{I}^{b+\varepsilon}
%\setminus \mc{O})
% = \pi(\tilde{\eta}(1,\iota([\mc{I}^m\leq b+\eps]%\mc{I}^{b+\varepsilon}
%\setminus \mc{O}))) 
% \subset \pi(\tilde{\eta}(1,\mc{H}^{b+\varepsilon}\setminus \tilde{\mc{O}})) 
% \subset\pi(\mc{H}^{b-\varepsilon}) \subset [\mc{I}^m\leq b-\eps]%\mc{I}^{b-\varepsilon}
%.$$ 
%\begin{align*}
%\eta(1,[\mc{I}^m\leq b+\eps]%\mc{I}^{b+\varepsilon}
%\setminus \mc{O})
% &= \pi(\tilde{\eta}(1,\iota([\mc{I}^m\leq b+\eps]%\mc{I}^{b+\varepsilon}
%\setminus \mc{O}))) 
% \subset \pi(\tilde{\eta}(1,[\mc{H}^m\leq b+\eps]_M%\mc{H}^{b+\varepsilon}
%\setminus \tilde{\mc{O}})) \\
%& \subset\pi([\mc{H}^m\leq b-\eps]_M%\mc{H}^{b-\varepsilon}
%) \subset [\mc{I}^m\leq b-\eps]%\mc{I}^{b-\varepsilon}
%.
%\end{align*}
\begin{eqnarray*}
\lefteqn{\eta(1,[\mc{I}^m\leq b+\eps]%\mc{I}^{b+\varepsilon}
\setminus \mc{O})} \\
 &&= \pi(\tilde{\eta}(1,\iota([\mc{I}^m\leq b+\eps]%\mc{I}^{b+\varepsilon}
\setminus \mc{O}))) 
 \subset \pi(\tilde{\eta}(1,[\mc{H}^m\leq b+\eps]_M%\mc{H}^{b+\varepsilon}
\setminus \tilde{\mc{O}})) \\
&& \subset\pi([\mc{H}^m\leq b-\eps]_M%\mc{H}^{b-\varepsilon}
) \subset [\mc{I}^m\leq b-\eps]%\mc{I}^{b-\varepsilon}
.
\end{eqnarray*}
The other inclusion is similar and easier.
\item 
We write $\tilde \eta(t,\theta,\lambda,u)=\big(\tilde \eta_0(t,\theta,\lambda,u), \tilde \eta_1(t,\theta,\lambda,u), \tilde \eta_2(t,\theta,\lambda,u)\big)$. Then by definition
$$\big(\eta_1(t,\lambda, u), \eta_2(t,\lambda, u)\big) = \Big(\tilde{\eta}_1(t,0,\lambda, u), \tilde{\eta}_2\big(t,0,\lambda, u(e^{-\tilde{\eta}_0(t,0,\lambda,u)}\cdot)\big)\Big)$$
thus by the property 6 of Theorem \ref{thm_def_H},
\begin{align*}
\big(\eta_1(t,\lambda, -u), \eta_2(t,\lambda, -u)\big) &= \Big(\tilde{\eta}_1(t,0,\lambda, -u), \tilde{\eta}_2\big(t,0,\lambda, -u(e^{-\tilde{\eta}_0(t,0,\lambda,-u)}\cdot)\big)\Big) \\
&= \Big(\tilde{\eta}_1(t,0,\lambda, u), -\tilde{\eta}_2\big(t,0,\lambda, u(e^{-\tilde{\eta}_0(t,0,\lambda,u)}\cdot)\big)\Big) \\
&=\big(\eta_1(t,\lambda, u), -\eta_2(t,\lambda, u)\big).
\end{align*}
\end{enumerate}
The theorem is hence proved.
\QED

\bigskip

Now we are ready to prove the main theorem for $\mc{H}^m$.

\medskip

\claim Proof of Theorem \ref{thm_def_H}.
To avoid cumbersome notation, we write $\xi=(\theta,\lambda, u)\in M$. Set 
$$M':=\{D\mc{H}^m(\xi)\neq 0\}.$$
It is known \cite{AmMa} that there exists a pseudo-gradient on the Hilbert manifold $M$ associated to $\mc{H}^m$, namely a locally Lipschitz vector field $\mc{V}:M'\to TM$ such that
\begin{itemize}
\item[(a)] $\norm{\mc{V}(\xi)}_{\xi} \leq 2 \norm{D\mc{H}^m(\xi)}_{\xi,*}$,
\item[(b)] $D\mc{H}^m(\xi)\cdot \mc{V}(\xi) \geq \norm{D\mc{H}^m(\xi)}_{\xi,*}^2$;
\end{itemize}
in particular,
\begin{equation}\label{eq_pseudograd}
\frac{1}{2} \norm{\mc{V}(\xi)}_{\xi} \leq \norm{D\mc{H}^m(\xi)}_{\xi,*} \leq \norm{\mc{V}(\xi)}_{\xi}.
\end{equation}
Moreover, we can ask, in the construction of the pseudo-gradient, that $\mc{V}$ is $\G$-equivariant, since $\mc{H}^m$ is $\G$-invariant.
Namely, set $\mc{V}=(\mc{V}_0,\mc{V}_1,\mc{V}_2)$, then $\mc{V}_0$ and $\mc{V}_1$ are even in $u$, while $\mc{V}_2$ is odd in $u$.

By Corollary \ref{cor_PStilde}, there exists $\delta=\delta_{\frac{\rho}{3}}>0$ such that
	\begin{equation}\label{eq_cor_PSP_delta}
	\forall \, \xi \in [b-\delta\leq \mc{H}^m\leq b+\delta]_M%\mc{H}^{b+\delta}_{b-\delta} 
\;\; s.t. \;\; \dist_M(\xi, \tilde{K}_b) > \frac{\rho}{3} \; : \; \norm{D\mc{H}^m(\xi)}_{\xi,*}> \delta .
	\end{equation}
We assume 
 \begin{equation}\label{choice-epsilon}
 \varepsilon < \min\Big\{ \frac{1}{2}\bar{\varepsilon}, \frac{1}{4} \delta, \frac16\rho\delta\Big\}.
 \end{equation}
Set the following
$$A:=[b-\eps \leq \mc{H}^m \leq b+ \eps]_M%\mc{H}^{b+\varepsilon}_{b-\varepsilon}
, \quad B:=[b-2\eps\leq\mc{H}^m\leq b+2\eps]_M%\mc{H}^{b+2\varepsilon}_{b-2\varepsilon}
$$
and choose a locally Lipschitz function $g\in C(M,[0,1])$ such that
$$g=1 \; \hbox{ on $A$}, \quad g=0 \; \hbox{ on $\complement B$},$$
for instance $g(\xi):=\frac{d(\xi,\complement B)}{d(\xi,\complement B)+d(\xi,A)}$. \\
When $\tilde K_b\not=\emptyset$, we choose a locally Lipschitz function $\tilde g\in C(M,[0,1])$ satisfying
$$ \tilde g=0 \ \hbox{ on }\, \tilde N_{\frac{\rho}{3}}(\tilde K_b), \quad	\tilde g=1 \ \hbox{ on }\, \complement \tilde N_{\frac{2}{3}\rho}(\tilde K_b). $$
When $\tilde K_b=\emptyset$, we set $\tilde g\equiv 1$. 
Moreover we introduce, for any $r\geq0$,
$$b(r):=\parag{\frac{1}{r} && \hbox{ if $r\geq 1$} \\ 1 && \hbox{ if $0\leq r < 1$}.}$$
Finally define, for $\xi \in M$,
$$W(\xi):=-g(\xi) \tilde g(\xi)b\left(\norm{\mc{V}(\xi)}_{\xi}\right)\mc{V}(\xi)$$
and, fixed $\xi \in M$, consider the Cauchy problem
$$\parag{ &\tilde{\eta}' = W(\tilde{\eta}),& \\ &\tilde{\eta}(0)=\xi.&}$$
We have that $W$ is well defined on $M$ and
$$\norm{W(\xi)}_{\xi}\leq \norm{\mc{V}(\xi)}_{\xi} \, b\left(\norm{\mc{V}(\xi)}_{\xi}\right) \leq 1, $$
where we have used that $|g|$, $|\tilde g|$ $\leq 1$.
Therefore we have the global existence of a flow $\tilde{\eta}=\tilde{\eta}(t,\xi)$; we are interested in $\tilde{\eta}$ restricted to $[0,1]$.
We now verify the desired properties.
\begin{itemize}
\item[1)] $\tilde{\eta}(0,\xi)=\xi$ by construction of the flow.
\item[2)] If $\xi \in [\mc{H}^m \leq b-\bar{\eps}]_M%\mc{H}^{b-\bar{\varepsilon}}
$, then $g(\xi)=0$, and thus $W(\xi)=0$. 
This means that $\tilde{\eta}(t,\xi)\equiv \xi$ is an equilibrium solution. 
Since $W\in Lip_{loc}(M)$ we have uniqueness of the solution, hence actually $\tilde{\eta}(t,\xi)\equiv \xi$.
\item[3)] We have
\begin{eqnarray*}
\lefteqn{\frac{d}{dt} \mc{H}^m(\tilde{\eta}(t,\xi)) 
= D\mc{H}^m(\tilde{\eta}(t,\xi)) \tilde{\eta}'(t,\xi)} \\
&=& -D\mc{H}^m(\tilde{\eta}(t,\xi)) \mc{V}(\tilde{\eta}(t,\xi)) g(\tilde{\eta}(t,\xi)) \tilde g(\tilde{\eta}(t,\xi))b\left(\norm{\mc{V}(\tilde{\eta}(t,\xi))}_{\tilde{\eta}(t,\xi)}\right)\\
&\leq& -\norm{D\mc{H}^m(\tilde{\eta}(t,\xi))}_{\tilde{\eta}(t,\xi),*}^2 g(\tilde{\eta}(t,\xi)) \tilde g(\tilde{\eta}(t,\xi)) b\left(\norm{\mc{V}(\tilde{\eta}(t,\xi))}_{\tilde{\eta}(t,\xi)}\right) \\
&\leq& 0
\end{eqnarray*}
that is the claim; we have used that $g$, $\tilde g$, $b$ are positive and the property $(b)$.

\item[4)]
We assume here $\tilde{K}_b= \emptyset$.
By using the fundamental theorem of calculus and previous arguments, we obtain
\begin{eqnarray*}
\lefteqn{\mc{H}^m(\tilde{\eta}(1,\xi)) - \mc{H}^m(\tilde{\eta}(0,\xi)) = \int_0^1 \frac{d}{ds} \mc{H}^m(\tilde{\eta}(s,\xi)) ds} \\
&=& -\int_0^1 D\mc{H}^m(\tilde{\eta}(s,\xi)) \mc{V}(\tilde{\eta}(s,\xi)) g(\tilde{\eta}(s,\xi)) b\left(\norm{\mc{V}(\tilde{\eta}(s,\xi))}_{\tilde{\eta}(s,\xi)}\right) ds \\
&\leq& -\int_0^1 \norm{D\mc{H}^m(\tilde{\eta}(s,\xi))}_{\tilde{\eta}(s,\xi),*}^2 g(\tilde{\eta}(s,\xi)) b\left(\norm{\mc{V}(\tilde{\eta}(s,\xi))}_{\tilde{\eta}(s,\xi)}\right) ds.
\end{eqnarray*}
Let now $\xi \in [\mc{H}^m\leq b+\eps]_M%\mc{H}^{b+\varepsilon}
$. This means, by point $3)$, that for $s\in [0,1]$
$$\mc{H}^m(\tilde{\eta}(s,\xi))\leq \mc{H}^m(\tilde{\eta}(0,\xi)) = \mc{H}^m(\xi) \leq b + \varepsilon,$$
thus $\tilde{\eta}(s,\xi)\in [\mc{H}^m\leq b+\eps]_M%\mc{H}^{b+\varepsilon}
$ and
$$ \mc{H}^m(\tilde{\eta}(1,\xi)) \leq b + \varepsilon - \int_0^1 \norm{D\mc{H}^m(\tilde{\eta}(s,\xi))}_{\tilde{\eta}(s,\xi),*}^2 g(\tilde{\eta}(s,\xi)) b\left(\norm{\mc{V}(\tilde{\eta}(s,\xi))}_{\tilde{\eta}(s,\xi)}\right) ds. $$
Assume now by contradiction that $\mc{H}^m(\tilde{\eta}(1,\xi))> b- \varepsilon$, which implies (again by point $3)$) $\mc{H}^m(\tilde{\eta}(s,\xi))> b- \varepsilon$, for all $s\in [0,1]$. 
Thus for all $s\in [0,1]$ we have $\tilde{\eta}(s,\xi) \in [b-\eps \leq \mc{H}^m \leq b+\eps]_M%\mc{H}^{b+\varepsilon}_{b-\varepsilon}
$ and in particular, since $\varepsilon < \frac{1}{2} \bar{\varepsilon}$, that $g(\tilde{\eta}(s,\xi))=1$; hence
$$ \mc{H}^m(\tilde{\eta}(1,\xi)) \leq b + \varepsilon - \int_0^1 \norm{D\mc{H}^m(\tilde{\eta}(s,\xi))}_{\tilde{\eta}(s,\xi),*}^2 b\left(\norm{\mc{V}(\tilde{\eta}(s,\xi))}_{\tilde{\eta}(s,\xi)}\right) ds. $$
By \eqref{eq_pseudograd}, by the fact that $\tilde{\eta}(s,\xi) \in [b-\eps \leq \mc{H}^m \leq b+\eps]_M%\mc{H}^{b+\varepsilon}_{b-\varepsilon}
\subset [b-\delta \leq \mc{H}^m \leq b+\delta]_M%\mc{H}^{b+\delta}_{b-\delta}
$ and by \eqref{eq_curv_pend_gen}, we have
\begin{equation}\label{eq_dim_stimaV}
\norm{\mc{V}(\tilde{\eta}(s,\xi))}_{\tilde{\eta}(s,\xi)} \geq \norm{D\mc{H}^m(\tilde{\eta}(s,\xi))}_{\tilde{\eta}(s,\xi),*} \geq \delta \geq 4 \varepsilon;
\end{equation}
in particular, 
$$ b\left(\norm{\mc{V}(\tilde{\eta}(s,\xi))}_{\tilde{\eta}(s,\xi)}\right) =\frac{1}{\norm{\mc{V}(\tilde{\eta}(s,\xi))}_{\tilde{\eta}(s,\xi)}}.$$
Thus, exploiting again \eqref{eq_pseudograd} and \eqref{eq_dim_stimaV} we obtain
\begin{align*}
\mc{H}^m(\tilde{\eta}(1,\xi)) &\leq b + \varepsilon - \frac{1}{2} \int_0^1 \norm{D\mc{H}^m(\tilde{\eta}(s,\xi))}_{\tilde{\eta}(s,\xi),*} ds\\
 &\leq b + \varepsilon - 2\int_0^1 \varepsilon ds = b - \varepsilon,
\end{align*}
which is an absurd.
\item[5)] We assume now $\tilde{K}_b \neq \emptyset$.
Let now $\xi \in [\mc{H}^m\leq b+\eps]_M%\mc{H}^{b+\varepsilon}
\setminus \tilde{\mc{O}}$. 
Assume again by contradiction that $\mc{H}^m(\tilde{\eta}(1,\xi))>b-\varepsilon$, which implies again $\tilde{\eta}(s,\xi)\in [b-\eps \leq \mc{H}^m \leq b+\eps]_M%\mc{H}^{b+\varepsilon}_{b-\varepsilon}
$. We distinguish two cases.
\\ \textbf{Case 1:} $\tilde{\eta}(t,\xi)\notin \tilde{N}_{\frac{2}{3}\rho}(\tilde{K}_b)$ for all $t\in [0,1]$. 
In this case we proceed as in the proof of 4). Indeed since $\varepsilon<\delta_{\frac{\rho}{3}}$, we are in the assumptions of \eqref{eq_curv_pend_gen} and thus
$$\norm{D\mc{H}^m(\tilde{\eta}(s,\xi))}_{\tilde{\eta}(s,\xi),*}> \delta>4 \varepsilon.$$
We argue %can do exactly the same passages 
as before and conclude.
\\ \textbf{Case 2:} $\tilde{\eta}(t^*,\xi)\in \tilde{N}_{\frac{2}{3}\rho}(\tilde{K}_b)$ for some $t^*\in [0,1]$. 
In this case $\varepsilon$ has to be better specified. We make a finer argument by choosing suitable $[\alpha,\beta]\subset [0,1]$ and observing that
\begin{align*}
 \mc{H}^m(\tilde{\eta}(1,\xi)) &\leq \mc{H}^m(\tilde{\eta}(\beta,\xi)) = \mc{H}^m(\tilde{\eta}(\alpha,\xi)) + \int_{\alpha}^{\beta} \frac{d}{ds} \mc{H}^m(\tilde{\eta}(s,\xi)) ds \\
&\leq \mc{H}^m(\tilde{\eta}(0,\xi)) + \int_{\alpha}^{\beta} \frac{d}{ds} \mc{H}^m(\tilde{\eta}(s,\xi)) ds \\
&\leq b+\varepsilon + \int_{\alpha}^{\beta} \frac{d}{ds} \mc{H}^m(\tilde{\eta}(s,\xi)) ds.
\end{align*}
Noting that $\tilde{\eta}(0,\xi)=\xi\notin \tilde{\mc{O}}=\tilde{N}_{\rho}(\tilde{K}_b)$ and $\tilde{\eta}(t^*,\xi)\in \tilde{N}_{\frac{2}{3}\rho}(\tilde{K}_b)$, we can find $\alpha$ and $\beta$ such that
$$\tilde{\eta}(\alpha)\in \partial \tilde{N}_{\rho}(\tilde{K}_b), \quad \tilde{\eta}(\beta)\in \partial \tilde{N}_{\frac{2}{3}\rho}(\tilde{K}_b),$$
and
$$\tilde{\eta}(s)\in \tilde{N}_{\rho}(\tilde{K}_b) \setminus \tilde{N}_{\frac{2}{3}\rho}(\tilde{K}_b) \quad \forall \, s\in (\alpha,\beta).$$
Hence we obtain by \eqref{eq_cor_PSP_delta}
$$\mc{H}^m(\tilde{\eta}(1,\xi)) \leq b+\varepsilon - \delta(\beta-\alpha).$$
We need an estimate from below of $\beta-\alpha$, which is obtained by observing that $\tilde{\eta}(\cdot,\xi)$ is a path connecting $\tilde{\eta}(\alpha, \xi)$ and $\tilde{\eta}(\beta, \xi)$, thus (recall that $1\geq \norm{W(\xi)}_{\xi}$)
\begin{align*}
\beta - \alpha &= \int_{\alpha}^{\beta} dt \geq \int_{\alpha}^{\beta} \norm{W(\tilde{\eta}(t,\xi))}_{\tilde{\eta}(t,\xi)} dt \\
&= \int_{\alpha}^{\beta} \norm{\tilde{\eta}'(t,\xi)}_{\tilde{\eta}(t,\xi)} dt \geq \dist_M(\tilde{\eta}(\alpha, \xi),\tilde{\eta}(\beta,\xi)) \\
&\geq\dist_M\left(\tilde{N}_{\rho}(\tilde{K}_b), \tilde{N}_{\frac{2}{3}\rho}(\tilde{K}_b)\right) \geq \frac{1}{3}\rho.
\end{align*}
Finally
$$\mc{H}^m(\tilde{\eta}(1,\xi)) \leq b+\varepsilon - \tfrac{1}{3} \rho \delta\leq b-\varepsilon$$
by our choice \eqref{choice-epsilon} of $\varepsilon$.

As regards the second inclusion, we argue in a similar way. Let $\xi \in [\mc{H}^m\leq b+\eps]_M%\mc{H}^{b+\varepsilon}
$. Case 1 can be done verbatim. 
In Case 2, if $\tilde{\eta}(1,\xi)\in \tilde{\mc{O}}$ we are done; if not, then we repeat the argument but with the path built thanks to $\tilde{\eta}(1,\xi)\notin \tilde{N}_{\rho}(\tilde{K}_b)$ and $\tilde{\eta}(t^*,\xi)\in \tilde{N}_{\frac{2}{3}\rho}(\tilde{K}_b)$.
\item[6)] Notice that, written $W=(W_0,W_1,W_2)$, we have that $W_0$ and $W_1$ are even in $u$ while $W_2$ is odd in $u$, since $\mc{V}$ is so and $g$, $b(\norm{D\mc{H}^m(\cdot)}_{\cdot,*})$ are even in $u$. 
Thus, by uniqueness of the solution, we have that $\tilde{\eta}$ satisfies the required symmetry properties.
\end{itemize}
The proof is thus concluded.
\QED

%%%%%%%%%%%%%%%%%%%%%%%%%%%%%%%%%%%%%

\subsection{Minimax values} % \boldmath{$a_j(\lambda)$}}

\subsubsection{Minimax values \boldmath{$a_j(\lambda)$}}

We write for $j\in\N$, $D_j:=\{\xi\in\R^j \mid \abs\xi\leq 1\}$ and we introduce the set of paths
$$\Gamma_j(\lambda) :=\big\{\gamma\in C(D_j, H^s_r(\R^N)) \, \mid \, \gamma \hbox{ odd}, \, \mc{J}(\lambda, \gamma(\xi))<0 \; \forall \xi \in \partial D_j \big\}$$
and
$$a_j(\lambda) :=\inf_{\gamma \in \Gamma_j(\lambda)} \sup_{\xi \in D_j} \mc{J}(\lambda, \gamma(\xi)).$$
By an odd extension from $[0,1]$ to $[-1,1]=D_1$, we may regard $\Gamma_1(\lambda)\equiv\Gamma(\lambda)$ and $a_1(\lambda)\equiv a(\lambda)$. 
Thus these quantities can be seen as %are a good 
generalizations. As for $j=1$, we prove the following properties.

\bigskip %PER FORMATTAZIONE PAGINA

\begin{Proposition}
Let $\lambda_0\in \R \cup \{+\infty\}$ be %the number 
given in %\eqref{muzero}--
\eqref{lambdazero}, $\lambda < \lambda_0$ and $j \in \N$.
\begin{enumerate}
\item $\Gamma_j(\lambda)\neq \emptyset$, thus $a_j(\lambda)$ is well defined. Moreover, it is increasing with respect to $\lambda$;
\item $a_j(\lambda)\leq a_{j+1}(\lambda)$;
\item $a_j(\lambda)> 0$; 
\item $\lim_{\lambda\to \lambda_0^-} \frac{a_j(\lambda)}{e^{\lambda}}=+\infty$;
\item if \hyperref[(g4)]{\textnormal{(g4)}} holds, then $\lim_{\lambda \to -\infty} \frac{a_j(\lambda)}{e^{\lambda}}=0$.
\end{enumerate}
\end{Proposition}

%\newpage%€

\claim Proof.
The proofs are quite the same of Propositions \ref{lem_buona_def}--\ref{out}. We point out just some slight differences.
\begin{enumerate}
\item For $\lambda<\lambda_0$, there exists $t_0>0$ such that 
	$$ G(t_0)-\frac{e^{\lambda}}{2}t_0^2 >0. $$
As in \cite{BL2}, we find that there exists a continuous odd map $\tilde{\gamma}:\,\partial D_j\to H_r^1(\R^N)\hookrightarrow H^s_r(\R^N)$ with $\mc{J}(\lambda,\tilde{\gamma}(\xi))<0$. 
Extending $\tilde{\gamma}$ onto $D_j$ we find $\Gamma_j(\lambda)\not=\emptyset$.
\item Since $D_j\subset D_{j+1}$, we observe $\gamma_{|D_j}\in\Gamma_j(\lambda)$ 
for $\gamma\in\Gamma_{j+1}(\lambda)$. Thus we regard $\Gamma_{j+1}(\lambda)\subset\Gamma_j(\lambda)$ 
and obtain 2).
\item Clear by $a_1(\lambda)=a(\lambda)>0$ and point $2)$. 
\item Again by $\lim_{\lambda\to \lambda_0^-} \frac{a(\lambda)}{e^{\lambda}}=+\infty$ and point $2)$.
\item We consider the path $\tilde{\gamma}:\partial D_j\to H^s_r(\R^N)$ obtained in 1) and introduce a path
$$ \xi \in D_j\mapsto \mu^{N/4} \tilde{\gamma}\left(\frac{\xi}{\abs{\xi}}\right)\big(\cdot/\mu^{-\frac{1}{2s}}\abs\xi\big) \in H^s_r(\R^N).
$$
Arguing as in Proposition \ref{out}, we have 5).
\QED
\end{enumerate}

%%%%%%%%%%%%%%%%%%%%%%%%%%%%%%%%%%%%%

\subsubsection{Minimax values \boldmath{$b^m_j$}}

We set
%\begin{eqnarray*}
%\Gamma_j^m:=\{\Theta\in C(D_j, \R\times H^s_r(\R^N)) &\mid&
% \hbox{$\Theta$ is $\G$-equivariant;} \\
%&& \mc{I}^m(\Theta(0)) \leq B_m-1; \\
%&& \Theta(\xi)\notin \Omega, \ \mc{I}^m(\Theta(\xi))\leq B_m-1 \ \hbox{for all
%}\ \xi \in \partial D_j\}
%\end{eqnarray*}
\begin{align*}
\Gamma_j^m:=\{\Theta\in C(D_j, \R\times H^s_r(\R^N)) \mid \; & 
 \hbox{$\Theta$ is $\G$-equivariant;} \\
& \mc{I}^m(\Theta(0)) \leq B_m-1; \\
& \Theta(\xi)\notin \Omega, \ \mc{I}^m(\Theta(\xi))\leq B_m-1 \ \hbox{for all
}\ \xi \in \partial D_j\}
\end{align*}
and
$$b_j^m:= \inf_{\Theta \in \Gamma_j^m} \sup_{\xi \in D_j} \mc{I}^m(\Theta(\xi)).$$
We notice that for $j=1$ we obtain $\Gamma_1^m \equiv \Gamma^m$ (up to an even/odd extension from $[0,1]$ to $[-1,1]=D_1$) and $b_1^m \equiv b_m$. 
So $\Gamma_j^m$ is a natural extension to build multiple solutions.

As in the case of $\Gamma^m$ and $b_m$, we want to prove that $\Gamma_j^m \neq \emptyset$ and that, for a fixed $k\in \N$, there exists an $m_k \gg 0$ (possibly equal to $0$) such that, if $m>m_k$, then $ b_j^m <0$ for $j=1 \dots k$.

%\newpage%€

\begin{Proposition}\label{prop_mk}
$\,$
\begin{itemize}
\item[\textnormal{(i)}] For any $\lambda <\lambda_0$, $m>0$, $j \in \N$, we have $\Gamma_j^m\neq \emptyset$ and $b_j^m\leq a_j(\lambda) - e^{\lambda} \frac{m}{2}$.
\item[\textnormal{(ii)}] For any $k\in\N$, set 
\begin{equation}\label{defmk}
	m_k := 2\inf_{\lambda<\lambda_0}\frac{a_k(\lambda)}{e^{\lambda}}\geq 0
\end{equation}
we have, % there exists $m_k\geq 0$ such that 
for any $m>m_k$
	$$ b_j^m < 0 \quad \hbox{for}\ j=1,2,\dots, k. $$
\item[\textnormal{(iii)}] $m_k=0$ for all $k\in\N$ if \hyperref[(g4)]{\textnormal{(g4)}} holds. That is,
	$$ b_j^m < 0 \quad \hbox{for all}\ j\in\N. $$
\end{itemize}
\end{Proposition}

\claim Proof. For (i), the proof is similar to Proposition \ref{tom}. We just need to set, for
$\zeta_{\lambda}\in \Gamma_j(\lambda)$, 
	$$ \psi_{\lambda}(\xi) := 
	\parag{
		(\lambda+L(1-2|\xi|), \,0) && \hbox{ if $|\xi| \in [0,1/2]$,} \\ 
		\left(\lambda, \, \zeta_{\lambda}\left(\frac{\xi}{\abs{\xi}}(2\abs\xi-1)\right)\right) && \hbox{ if $|\xi| \in (1/2,1]$}
		} $$	
and we come up again to the same proof.

For (ii), (iii), we come up with a proof similar to Proposition \ref{blue}, %, we set 
%\begin{equation}\label{defmk}
%	m_k := 2\inf_{\lambda<\lambda_0}\frac{a_k(\lambda)}{e^{\lambda}}\geq 0
%\end{equation}
observing in addition that $m_k\leq m_{k+1}$ since $a_k(\lambda)$ are increasing in $k$.
\QED

\bigskip

By Proposition \ref{PSP} and Theorem \ref{defarg} % and \ref{...}, 
$\mc{I}^m$ satisfies the $(PSP)_b$ condition for $b<0$ and the deformation lemma holds. 
Let $m_k\geq 0$ be a number given in Proposition \ref{prop_mk}. For $m>m_k$ we can see that $b_j^m<0$ for $j=1,2,\dots, k$ are critical values of $\mc{I}^m$.
If $b_j^m$ are different, we directly have multiplicity of solutions. To deal with the case $b_j^m=b_{j'}^m$ for some $j\not=j'$, we need another family of minimax methods, which exploits the topological information hidden in this equality. % we consider in the following Section.

%%%%%%%%%%%%%%%%%%%%%%%%%%%%%%%%%%%%%

\subsubsection{Minimax values \boldmath{$c_j^m$}}

Let us define minimax families $\Lambda_j^m$ which allow to find multiple solutions. We use an idea from \cite{Rab0}.
In what follows, we denote by $\genus(A)$ the genus of closed symmetric sets $A$ with $0\not\in A$ (see Appendix \ref{sec_app_genus}). 

Define, for each $j\in \N$, 
%\begin{eqnarray*}
%\Lambda_j^m:=\{A=\Theta(\overline{D_{j+l}\setminus Y}) &\mid & l\geq 0, \; \Theta\in \Gamma_{j+l}^m, \\
%&& Y\subseteq D_{j+l}\setminus \{0\} \; \hbox{ is closed, symmetric in $0$ and $\genus(Y)\leq l$} \}
%\end{eqnarray*}
\begin{eqnarray*}
\Lambda_j^m:=\{A=\Theta(\overline{D_{j+l}\setminus Y}) &\mid & l\geq 0, \; \Theta\in \Gamma_{j+l}^m, \\
&& Y\subseteq D_{j+l}\setminus \{0\} \; \hbox{ is closed, symmetric in $0$,} \\
&& \qquad \qquad \qquad \quad \; \hbox{and $\genus(Y)\leq l$} \}
\end{eqnarray*}
and
$$c_j^m:= \inf_{A\in \Lambda_j^m} \sup_{A} \mc{I}^m.$$

In the following lemma, we observe that $\Lambda_j^m$ includes, in some way, $\Gamma^m_j$ and that it inherits the property that the paths intersect $\partial \Omega$.

\begin{Lemma}\label{lem7.6}
$\,$
\begin{itemize}
\item[\textnormal{(i)}] $\Lambda_j^m \neq \emptyset$;
\item[\textnormal{(ii)}] $c_j^m\leq b_j^m$;
\item[\textnormal{(iii)}]
for any $A\in \Lambda_j^m$, we have $A\cap \partial \Omega \neq \emptyset$. As a consequence, we obtain
$$b_m=B_m= E_m\leq c_j^m.$$
\end{itemize}
\end{Lemma}

\claim Proof.
Indeed, we see that, by choosing $l=0$ and $Y=\emptyset$ we have
$$\{A=\Theta(D_j) \, \mid \, \Theta \in \Gamma_j^m\} \subset \Lambda_j^m$$
from which %easily 
come the first two claims.

Focus on the third claim. Let $A=\Theta(\overline{D_{j+l}\setminus Y})$ and set $U:=\Theta^{-1}(\Omega)$. 
By the symmetry in $(\lambda, u)$ of $\Theta$ and the symmetry in $u$ of $\Omega$ we have that $U$ is symmetric. 
Moreover, since $\Theta(0)\in \Omega$, we have that $U\subset D_{j+l} \subset \R^{j+l}$ is a symmetric neighborhood of the origin. 
By Proposition \ref{prop_genus_g} we have %that the genus of $\partial U$ is maximum, that is
\begin{equation}\label{eq_genus}
\genus(\partial U) = j+l.
\end{equation}
Observe in addition the following chain of inclusions
$$\overline{\partial U \setminus Y} = \overline{(\partial U \cap D_{j+l})\setminus Y} = \overline{(D_{j+l} \setminus Y) \cap \partial U} \subseteq \overline{D_{j+l} \setminus Y} \cap \overline{\partial U} = \overline{D_{j+l}\setminus Y} \cap \partial U $$
thus
$$\Theta\left(\overline{\partial U \setminus Y}\right) \subseteq \Theta\left( \overline{D_{j+l}\setminus Y} \cap \partial U\right) \subseteq \Theta\left( \overline{D_{j+l}\setminus Y}\right) \cap \Theta \left(\partial U\right) = A\cap \Theta \left(\partial U\right).$$
%\tr{Since $\Theta$ is continuous, we know that $\tilde{\partial} U = \tilde{\partial} \left( \Theta^{-1}(\Omega)\right) \subseteq \Theta^{-1}(\partial \Omega)$, and thus $\Theta(\tilde{\partial} U) \subset \partial \Omega$, where $\tilde{\partial}$ is the boundary with respect to the domain of $\Theta$, i.e. $D_{j+l}$; but $\tilde{\partial} U$ is generally smaller than $\partial U$ (the boundary made with respect to the whole space $\R^{j+l}$), which is the one appearing in \eqref{eq_genus}. Thus we need to be more precise.} 
Assume for the moment that it holds
\begin{equation}\label{eq_aggiunta_corr}
\Theta(\partial U) \subset \partial \Omega.
\end{equation}%but in this last expression we obtain
Then by the previous computation we have
$$\Theta\left(\overline{\partial U \setminus Y}\right) \subseteq A \cap \partial \Omega. $$
Thus, to reach the claim, we need %it is sufficient 
to show that $\overline{\partial U \setminus Y}\neq \emptyset$. 
But is an immediate consequence of \eqref{eq_genus} and Proposition \ref{prop_genus_g} that %the property of the genus
$$ \genus(\overline{\partial U \setminus Y}) \geq \genus(\partial U) - \genus(Y) \geq (j+l)-l = j \geq 1$$
which directly excludes the possibility that $\overline{\partial U \setminus Y}$ is empty. 

Focus now on \eqref{eq_aggiunta_corr}; we first observe that, by continuity, we have 
$$\partial_{D_{j+l}} U = \partial_{D_{j+l}} \big(\Theta^{-1}(\Omega)\big) \subset \Theta^{-1}(\partial \Omega),$$
 where $\partial_{D_{j+l}}$ is the boundary with respect to the topology restricted to $D_{j+l}$, but this is not enough, since $\partial_{D_{j+l}} U$ is generally smaller than $\partial U$ (the boundary made with respect to the whole space $\R^{j+l}$), which is the one appearing in \eqref{eq_genus}. 
Let thus $\xi \in \partial U$; we need to show that $\Theta(\xi) \in \partial \Omega$. By definition of $U$, we have $\Theta(\xi) \in %\Theta(\partial(\Theta^{-1}(\Omega))) \subset 
\Theta(\overline{\Theta^{-1}(\Omega)}) %\subset \overline{\Theta(\Theta^{-1}(\Omega))} 
\subset \overline{\Omega}$; assume by contradiction $\Theta(\xi) \in \Omega$. We first observe that $\xi \notin \partial D_{j+l}$, by definition of $\Theta \in \Gamma^m_{j+l}$, thus $\xi$ is in the interior of $D_{j+l}$. We then can find a neighborhood $N_1$ of $\xi$ (with respect to $\R^{j+l}$) contained in $D_{j+l}$, and a neighborhood $M_2$ of $\Theta(\xi)$ contained in $\Omega$; set $N:= N_1 \cap \Theta^{-1}(M_2)$, we have that $N$ is a neighborhood of $\xi$ (with respect to $\R^{j+l}$) contained in $U$, which implies that $\xi$ is in the interior of $U$, absurd.
This concludes the proof of the first part.

We prove now the consequence. Indeed, for each $A\in \Lambda_j^m$ we have
$$E_m = \inf_{\partial \Omega} \mc{I}^m \leq \inf_{\partial \Omega \cap A} \mc{I}^m \leq \sup_{\partial \Omega \cap A} \mc{I}^m \leq \sup_{A} \mc{I}^m$$
and thus the claim passing to the infimum over $\Lambda_j^m$.
\QED

\bigskip

Let us now show the main properties of $\Lambda_j^m$ and $c_j^m$,
%We point out that these classical properties are the only ones 
which will actually be the only ones used in the multiplicity result. % proof of the existence of multiple solutions.

\begin{Proposition} \label{prop_genus}
Let $j\in \N$.
\begin{enumerate}
\item $\Lambda_j^m \neq \emptyset$;
\item $\Lambda_{j+1}^m\subseteq \Lambda_j^m$, and thus $c_j^m\leq c_{j+1}^m$;
\item let $A\in \Lambda_j^m$ and $Z\subset \R \times H^s_r(\R^N)$ be $\G$-invariant, closed, 
and such that $0 \notin \overline{P_2(Z)}$ and $\genus(\overline{P_2(Z)})\leq i$. Then $\overline{A\setminus Z} \in \Lambda_{j-i}^m$.
\end{enumerate}
Fix now $k\in \N$, and let $m>m_k$, where $m_k$ has been introduced in %Proposition \ref{prop_mk}, i.e., in 
\eqref{defmk}.
Then
\begin{itemize}
\item[4.] $c_j^m<0$ and $\mc{I}^m$ satisfies $(PSP)_{c_j^m}$;
\item[5.] if $A\in \Lambda_j^m$ and $\eta$ is a deformation as in Theorem \ref{thm_def_gen} for $b=c_j^m$, then $\eta(1,A)\in \Lambda_j^m$.
\end{itemize}
\end{Proposition}

%\vspace*{-4\baselineskip} %€
%\enlargethispage{1.8\baselineskip}

\claim Proof.
Properties $1)$ and $4)$ have already been shown in the Lemma \ref{lem7.6}, while property $2)$ is a consequence of the definition. 
Let us see properties $3)$ and $5)$.
\begin{itemize}
\item[3)]
Let $A=\Theta(\overline{D_{j+l}\setminus Y}) \in \Lambda_j^m$ and let $Z$ be $\G$-invariant, closed and such that $0\notin \overline{P_2(Z)}$ and $\genus(\overline{P_2(Z)})\leq i$. 
Assume it holds
\begin{align}\label{eq_Gamma_g}
\overline{A\setminus Z} &= \Theta((\overline{D_{j+l}\setminus Y)\setminus \Theta^{-1}(Z)}) \\
&= \Theta(\overline{D_{(j-i)+(l+i)} \setminus (Y \cup \Theta^{-1}(Z))}); \notag
\end{align}
if $\genus(Y \cup \Theta^{-1}(Z))\leq l+i$ we have the claim. But this is a direct consequence of the assumptions and Proposition \ref{prop_genus_g}, since
\begin{align*}
\genus(Y \cup \Theta^{-1}(Z))&\leq \genus(Y)+\genus(\Theta^{-1}(Z)) \\
&\leq l + \genus(\overline{h(\Theta^{-1}(Z))}) \\
&= l + \genus(\overline{P_2(Z)}) \leq l+i
\end{align*}
where we have set $h:=P_2 \circ \Theta$, which is an odd map and thus admissible for the genus.

Turn now to \eqref{eq_Gamma_g}. Set $B:=D_{j+l}\setminus Y$ and $W:=\Theta^{-1}(Z)$ we have to prove
$$\overline{\Theta(\overline{B})\setminus \Theta(W)} = \Theta(\overline{B\setminus W}).$$
We have
$$\overline{\Theta(\overline{B})\setminus \Theta(W)} \subseteq
 \overline{\Theta(\overline{B}\setminus W)} \stackrel{(i)} \subseteq \overline{\Theta(\overline{B\setminus W})} \stackrel{(ii)}= \Theta(\overline{B\setminus W}) $$
and
$$\Theta(\overline{B\setminus W}) \stackrel{(iii)} \subseteq \overline{\Theta(B\setminus W)} \stackrel{(iv)} = \overline{\Theta(B)\setminus \Theta(W)} \subseteq \overline{\Theta(\overline{B})\setminus \Theta(W) }$$
where
\begin{itemize}
\item[\textnormal{(i)}] is due to the fact that $W$ is closed;
\item[\textnormal{(ii)}] $\overline{B\setminus W}\subseteq D_{j+l}$ is compact, thus $\Theta(\overline{B\setminus W})$ is closed;
\item[\textnormal{(iii)}] derives from the continuity of $\Theta$;
\item[\textnormal{(iv)}] is due to the fact that $W$ is a preimage.
\end{itemize}
\item[5)]
Consider $0<\bar{\varepsilon}<1$, $b=c_j^m\geq B_m$ and $\eta$ as in the deformation lemma, and fix $A=\Theta(\overline{D_{j+l}\setminus Y}) \in \Lambda_j^{m}$ with $\Theta \in \Gamma_{j+l}^m$. 
To show that $\eta(1,A)\in \Lambda_j^m$ and conclude the proof, it is sufficient to show that $\tilde{\Theta}:=\eta(1, \Theta)\in \Gamma_{j+l}^m$ as well.
\begin{itemize}
\item[$\bullet$] $\tilde{\Theta}(-\xi) = \eta(1,\Theta(-\xi)) = \eta(1,\Theta_1(-\xi), \Theta_2(-\xi)) = \eta(1,\Theta_1(\xi), -\Theta_2(\xi))$
and thus
\begin{eqnarray*}
\lefteqn{\big(\tilde{\Theta}_1(-\xi), \tilde{\Theta}_2(-\xi)\big) = \Big(\eta_1\big(1,\Theta_1(\xi), -\Theta_2(\xi)\big), \eta_2\big(1,\Theta_1(\xi), -\Theta_2(\xi)\big)\Big)} \\
&=& \Big(\eta_1\big(1,\Theta_1(\xi), \Theta_2(\xi)\big), -\eta_2\big(1,\Theta_1(\xi), \Theta_2(\xi)\big)\Big) 
= \big(\tilde{\Theta}_1(\xi), -\tilde{\Theta}_2(\xi)\big) 
\end{eqnarray*}
which shows that $\tilde{\Theta}_1$ is even and $\tilde{\Theta}_2$ is odd.
\item[$\bullet$] By Lemma \ref{lem7.6}, for $\xi=0$ and $\xi \in \partial D_{j+l}$ we have $\mc{I}^m(\Theta(\xi)) \leq B_m-1 = E_m-1 \leq c_j^m-\bar{\varepsilon}$, thus $\Theta(\xi)\in [\mc{I}^m\leq c_j^m-\bar{\eps}]%\mc{I}^{c_j^m-\bar{\varepsilon}}
$. 
Therefore $\tilde{\Theta}(\xi)=\eta(1,\Theta(\xi))=\Theta(\xi)$ for $\xi=0$ and $\xi \in \partial D_{j+l}$, and the same properties are satisfied.
\QED
\end{itemize}
\end{itemize}

%%%%%%%%%%%%%%%%%%%%%%%%%%%%%%%%%%%%%

\subsection{Multiplicity theorem}

Fix $k\in \N^*$, and let $\Lambda_j^m$ and $c_j^m$ be given in the previous Section for $j=1\dots k$.
Exploiting the properties given in Proposition \ref{prop_genus}, we can find multiple solutions.

%\enlargethispage{1\baselineskip}%€

\begin{Theorem}\label{theopassmult}
Fix $k\in \N^*$, and assume $m>m_k$. We have that
$$c_1^m \leq c_2^m \leq \dots \leq c_k^m<0$$
are critical values of $\mc{I}^m$. Moreover
\begin{itemize}
\item[\textnormal{(i)}]
 if, for some $q\geq 1$,
$$c_j^m < c_{j+1}^m < \dots < c_{j+q}^m$$
then we have $q+1$ different nonzero critical values, and thus $q+1$ different (pairs of) nontrivial solutions of \eqref{problem_frac};
\item[\textnormal{(ii)}]
if instead, for some $q\geq 1$,
\begin{equation}\label{eq7.6}
	c_j^m = c_{j+1}^m = \dots = c_{j+q}^m \equiv b
\end{equation}
then 
\begin{equation}\label{eq_genus_teo}
\genus(P_2(K_b^{PSP}))\geq q+1
\end{equation}
and thus (by Proposition \ref{prop_genus_g}) $\# P_2(K_b^{PSP})=+\infty$, which means that we have infinite different solutions of \eqref{problem_frac}. %
\end{itemize}%
Summing up, we have at least $k$ different (pairs of) solutions of \eqref{problem_frac} which satisfy the Pohozaev identity \eqref{eq_Pohozaev_fractional}.
\end{Theorem}

\claim Proof.
It is sufficient to show only the property \eqref{eq_genus_teo} on the genus: indeed by choosing $q=0$ we have that, for each $j$, $\#(K_{c_j^m}^{PSP})\geq 1$ and thus $c_j^m$ is a nontrivial critical value.

By the $(PSP)_b$ we have that $K_b^{PSP}$ is compact, thus $P_2(K_b^{PSP})$ is compact; moreover it is symmetric with respect to $0$ and does not contain $0$ (see Corollary \ref{PSP-cor}).

By Proposition \ref{prop_genus_g} %the property of the genus 
we can find a (closed, symmetric with respect to origin, not containing the zero) neighborhood $N$ of $P_2(K_b^{PSP})$ which preserves the genus, i.e. $\genus(N)=\genus(P_2(K_b^{PSP}))$. 
We can easily think $N$ as a projection of a neighborhood $Z$ of $K_b^{PSP}$ (i.e. $N=P_2(Z)$) satisfying the properties of Proposition \ref{prop_genus}. 

By Theorem \ref{thm_def_gen}, there exist a sufficiently small $\varepsilon$ and an $\eta$ such that $\eta([\mc{I}^m\leq b+\eps]%\mc{I}^{b+\varepsilon}
\setminus Z) \subseteq [\mc{I}^m\leq b-\eps]%\mc{I}^{b-\varepsilon}
$. 
Corresponding to $\varepsilon$, by definition of $c_j^m$, there exists an $A\in \Lambda_{j+q}^m$ such that $\sup_{A}\mc{I}^m<b+\varepsilon$, that is $A\subseteq [\mc{I}^m\leq b+\eps]%\mc{I}^{b+\varepsilon}
$. Thus, being $\eta(1,\cdot)$ continuous
%$$ \eta(1, \overline{A\setminus Z}) \subseteq \eta(1, \overline{[\mc{I}^m\leq b+\eps]%\mc{I}^{b+\varepsilon}
%\setminus Z}) \stackrel{\eta(1,\cdot) \hbox{ continuous}} \subseteq \overline{\eta(1,[\mc{I}^m\leq b+\eps]%\mc{I}^{b+\varepsilon}
%\setminus Z)} \subseteq \overline{[\mc{I}^m\leq b-\eps]%\mc{I}^{b-\varepsilon}
%} = [\mc{I}^m\leq b-\eps]%\mc{I}^{b-\varepsilon}
%,$$
\begin{align*}
\eta(1, \overline{A\setminus Z}) &\subseteq \eta(1, \overline{[\mc{I}^m\leq b+\eps]%\mc{I}^{b+\varepsilon}
\setminus Z}) 
%\stackrel{\eta(1,\cdot) \hbox{ continuous}} 
\subseteq \overline{\eta(1,[\mc{I}^m\leq b+\eps]%\mc{I}^{b+\varepsilon}
\setminus Z)} \\
&\subseteq \overline{[\mc{I}^m\leq b-\eps]%\mc{I}^{b-\varepsilon}
} 
= [\mc{I}^m\leq b-\eps]%\mc{I}^{b-\varepsilon}
,
\end{align*}
and hence
\begin{equation}\label{eq_th_molt}
\sup_{\eta(1,\overline{A\setminus Z})} \mc{I}^m \leq b-\varepsilon.
\end{equation}
\sloppy %va a capo ed evita l'overflow inline
On the other hand, assume by contradiction that $\genus(P_2(K_b^{PSP}))\leq q$, i.e. 
 $\genus(P_2(Z))\leq q$. We use now the properties on $c_j^m$ and $\Lambda_j^m$.

Replacing $j$ with $j+q$ and $i$ with $q$ and applying Proposition \ref{prop_genus}, we have $\overline{A\setminus Z}\in \Lambda_j^m$; by property 5) of Proposition \ref{prop_genus} we obtain $\eta(1,\overline{A\setminus Z}) \in \Lambda_j^m$, which implies (by definition of $c_j^m$) 
$$\sup_{\eta(1,\overline{A\setminus Z})} \mc{I}^m \geq c_j^m=b.$$
This is a contradiction with \eqref{eq_th_molt}, and thus concludes the proof.
\QED

\bigskip
\claim Proof of Theorem \ref{S:1.13}.
As consequence of Theorem \ref{theopassmult}, we derive (i). 
We pass to prove (ii). Under condition \hyperref[(g4)]{\textnormal{(g4)}}, we have $m_k=0$ for all $k \in \N$. 
Thus for any $j \in \N$, $c_j^m$ is a critical value of $\mc{I}^m$ and $c_j^m \leq b_j^m <0$. Since $c_j^m$ is an increasing sequence, we have $c_j^m \to \bar c \leq 0$ as $j \to \infty$. We need to show that $\bar c=0$.

By contradiction we assume $\bar c <0$. Then $K^{PSP}_{\bar c}$ is compact and $K^{PSP}_{\bar c} \cap (\R \times \{0\}) =\emptyset$. It follows that $q= \genus(P_2(K^{PSP}_{\bar c})) < \infty$.
Arguing as in the proof of Theorem \ref{theopassmult}, let $\delta>0$ be such that $q= \genus(P_2(N_\delta(K^{PSP}_{\bar c}))) < \infty$.
By Theorem \ref{thm_def_gen}, there exist $\varepsilon\in (0,1)$ small and $\eta: [0,1] \times \R \times H^s_r(\R^N) \to \R \times H^s_r(\R^N)$ satisfying
\begin{equation}\label{ax1}
\eta(1,[\mc{I}^m\leq \bar{c}+\eps]%\mc{I}^{\bar c+ \varepsilon}
\setminus N_\delta(K^{PSP}_{\bar c})) \subseteq [\mc{I}^m\leq \bar{c}-\eps]%\mc{I}^{ \bar c-\varepsilon}
\end{equation}
and
\begin{equation}\label{ax2}
\eta(t,\lambda, u) =(\lambda,u) \quad \hbox{if $\, \mc{I}^m(\lambda, u) \leq B_m -1$}. 
\end{equation} 
We can choose $j \in \N$ sufficiently large such that $c_j^m > \bar c - \varepsilon$ and take $B \in \Lambda_{j+q}^m$ such that $B \subset [\mc{I}^m\leq \bar{c}+\eps]%\mc{I}^{\bar c +\varepsilon}
$.
Then we have 
$$ \overline{B \setminus N_\delta(K^{PSP}_{\bar c})} \in \Lambda_j^m.$$
From equations \eqref{ax1}, \eqref{ax2} we derive $c_j^m \leq \bar c - \varepsilon$, which gives a contradiction.
\QED

\medskip

\begin{Remark}
We observe that, even if the problem is invariant under translations, the found solutions are not translations of a same solution since they all are radially symmetric. Moreover, assuming \hyperref[(g4)]{\textnormal{(g4)}}, since $0>c_j^m \to 0$ we easily find a sequence of solutions with distinct energy levels.
\end{Remark}

%%%%%%%%%%%%%%%%%%%%%%%%%%%%%%%%%%%%%%%%%%%%%%%%%

\section{$L^2$-minimum}
\label{sec_L2_min_prima}

In Theorems \ref{S:1.1_frac} and \ref{S:1.12_frac} we find a solution via mountain pass minimax methods. 
We remark that this solution is characterized as minimizer of the functional $\mc{L}$ on $\mc{S}_m$, where $\mc{L}:\, H_r^s(\R^N)\to\R$ is defined by
	$$	\mc{L}(u):=\half\norm{(-\Delta)^{s/2}u}_2^2 -\int_{\R^N} G(u)
	$$
and $\mc{S}_m$ is the $L^2$-sphere in $H_r^S(\R^N)$, i.e. 
$$\mc{S}_m:=\{u\in H_r^s(\R^N) \mid \norm u_2^2=m\}.$$
Set
	$$	\kappa_m:=\inf_{u\in\mc{S}_m} \mc{L}(u).
	$$
%we pass to recognize that the Mountain Pass solution found in Theorems \ref{S:1.1_frac} and \ref{S:1.12_frac} is a ground state solution, namely, a minimizer of ${\mathcal L}$ on the sphere.

\begin{Proposition}\label{prop_equiv_gs}
Assume \hyperref[(g1)]{\textnormal{(g1)}}--\hyperref[(g3)]{\textnormal{(g3)}}, and let $m \geq m_0$, where $m_0$ is introduced in Proposition \ref{blue}. 
 We have that the following statements hold.
\begin{itemize}
%\item[(i)] Every Mountain Pass solution at level $b_m$ is a Pohozaev minimum on the product space, that is
%$$ b_m = B_c$$
%where $B_c$ is defined in Corollary \ref{...}.
\item[(i)] The Mountain Pass level and the ground state level coincide, i.e.
\begin{equation}\label{eq_coinc_livel}
\kappa_m = b_m.
\end{equation}
In particular, thanks to Corollary \ref{coroll_esist_Pm}, %Theorem \ref{...}, 
there exists a ground state of $\mc{L}_{|\mc{S}_m}$.
\item[(ii)] Every ground state of $\mc{L}_{|\mc{S}_m}$ satisfies the Pohozaev identity \eqref{eq_Pohozaev_fractional} with $\mu$ the associated Lagrange multiplier. Thanks to \eqref{eq_coinc_livel}, the same conclusion holds for every Mountain Pass solution at level $b_m$.
\item[(iii)] Every ground state of $\mc{L}_{|\mc{S}_m}$ has a positive associated Lagrange multiplier. This means that every ground state of $\mc{L}_{|\mc{S}_m}$ is a solution of problem \eqref{problem_frac}. %\textnormal{($P_m$)}. %\eqref{...}.
\end{itemize}
Moreover, if \hyperref[(g4)]{\textnormal{(g4)}} holds, then $m_0=0$.
\end{Proposition}

\claim Proof. 
%(i) From \eqref{eq_stima_beta_c} and Proposition \ref{...} we have
%$$b_m \leq \inf_{\lambda \in \R} \left(a(\lambda) - \frac{e^{\lambda}}{2} c\right) \leq B_m (?).$$
%
%On the other hand, each path of $\Lambda_m$ passes through $\partial \Omega$ by definition, and thus
%$$\beta_c \geq B_m (?)$$
%which gives the claim $(i)$.
%
(i) 
Let $u_*$ be the Mountain Pass solution obtained in Corollary \ref{coroll_esist_Pm}, %Theorems \ref{...} and \ref{...}, 
which verifies $\norm{u_*}_2^2=m$. Thus,
\begin{equation}\label{u_star}
\kappa_m \leq {\mathcal L}(u_*) =b_m <0.
\end{equation}
In particular, by \eqref{u_star} we can find a minimizing sequence $(u_n)_n \subset \mc{S}_m$ for $\kappa_m$ satisfying ${\mathcal L}(u_n)<0$, and thus we can set
$$e^{\lambda_n} := \frac{2}{N m} \left( s \norm{(-\Delta)^{s/2} u_n}_2^2 - N \mc{L}(u_n)\right)>0
$$
so that ${\mathcal P}(\lambda_n, u_n)=0$, i.e., $(\lambda_n, u_n) \in \partial \Omega$. 
At this point Proposition \ref{blue} %Corollary \ref{...} and $...$ 
implies
$$\kappa_m + o(1) = {\mathcal L}(u_n) = {\mathcal I}^m(\lambda_n, u_n) \geq E_m %B_m \tr{(?)}
=b_m.$$
Passing to the limit, together with \eqref{u_star}, we have \eqref{eq_coinc_livel}.

\smallskip

(ii)
Let $u_0$ be a minimizer of ${\mathcal L}$ on $\mc{S}_m$. Corresponding to $u_0$, there exists a Lagrange multiplier $\mu_0\in\R$ such that
	$$
	(-\Delta)^{s/2} u_0 +\mu_0u_0 =g(u_0),
	$$
and thus, in particular,
	\begin{equation}\label{R1}
	\norm{(-\Delta)^{s/2} u_0}_2^2 + \mu_0\norm{u_0}_2^2 -\int_{\R^N} g(u_0)u_0\, dx = 0.
	\end{equation}
We show first that $u_0$ satisfies the Pohozaev identity. In fact, we consider the $\R$-action 
$\Phi:\, \R\times \mc{S}_m\to \mc{S}_m$ defined by
	\begin{equation}\label{Phi}
	(\Phi_\theta v)(x) := e^{\frac{N}{2}\theta}v(e^\theta x),
	\end{equation}
since $\norm{\Phi_\theta v}_2^2 = \norm{v}_2^2$. Then we have
	$$	{\mathcal L}(\Phi_\theta u_0) = \half e^{2 s \theta}	\norm{(-\Delta)^{s/2} u_0}_2^2	
- e^{-N\theta} \int_{\R^N} G\big(e^{\frac{N}{2} \theta} u_0\big).
		%-\half e^{-(N+\alpha)\theta}\int_{\R^N} (I_\alpha*F(e^{\frac{N}{2}\theta}u_0))F(e^{\frac{N}{2}\theta}u_0)\,dx.
	$$
Since $u_0$ is a minimizer, we have $\frac{d}{d \theta}\bigr|_{\theta=0}{\mathcal L}(\Phi_\theta u_0)=0$, that is,
%\begin{adjustwidth}{-4.6cm}{0cm} 
%\begin{adjustwidth}{-0.17cm}{0cm} 
	\begin{equation}\label{R2}
s \norm{(-\Delta)^{s/2} u_0}_2^2 
+ N \int_{\R^N} G(u_0) -\frac{N}{2} \int_{\R^N} g(u_0)u_0\, dx =0.
%+\frac{N+\alpha}{2}\int_{\R^N} (I_\alpha*F(u_0))F(u_0)\, dx -\frac{N}{2} \int_{\R^N} (I_\alpha*F(u_0))f(u_0)u_0\, dx =0.
	\end{equation}
%\end{adjustwidth} 
%
From \eqref{R1} and \eqref{R2}, the Pohozaev identity follows %mdpi2: please check intended meaning has been retained. %I confirm
	\begin{equation}\label{glom}
	\frac{N-2s}{2} 
	\norm{(-\Delta)^{s/2} u_0}_2^2
	 + \frac{N}{2} \mu_0 \|u_0\|_2^2 - N \int_{\R^N} G(u_0)=0.
%\frac{N + \alpha}{2} \, {\mathcal D}(u_0).
	\end{equation}

\smallskip
(iii)
Finally, from \eqref{u_star} we have ${\mathcal L}(u_0)=\kappa_m <0$, that is
	\begin{equation}\label{R3}
	\half \norm{(-\Delta)^{s/2} u_0}_2^2 - \int_{\R^N} G(u_0)=\kappa_m<0,
	\end{equation}
which joined to \eqref{glom} gives $\mu_0>0$. This concludes the proof.
\QED

\begin{Remark}
By \cite[Theorem 4.1]{LoMa08}, we have that actually every $L^2$-minimum is radially symmetric (up to a translation). Thus $\kappa_m$ coincide with the infimum made on the $L^2$-ball of the whole space $H^s(\R^N)$. 
\end{Remark}

%%%%%%%%%%%%%%%%%%%%%%%%%%%%%%%%%%%%%%%%%%%%%%%%%

\section{Relation between constrained and unconstrained problems}

Let $0<\mu<\mu_0$ and $m>0$. By joining the results of Proposition \ref{blue} and Proposition \ref{minimizing}, 
%by slightly changing the definition of $\kappa_m$, 
we proved the following relation.
\begin{equation}\label{eq_relaz_cons_uncons}
\kappa(m)=\inf_{\mu \in (0,\mu_0)} \big(p(\mu) -\mu \, m \big)
\end{equation}
where we slightly changed the definition of $L^2$-minimum
$$\kappa(m) := \inf_{
\substack{ u\in H^s_r(\R^N) \\ 
\half \norm{u}_2^2=m} }
\left( \half\norm{(-\Delta)^{s/2}u}_2^2 -\int_{\R^N} G(u) \right) $$
and of Pohozaev minimum
$$p(\mu) := \inf_{
\substack{u \in H^s_r(\R^N)\setminus\{0\} \\ 
\norm{(-\Delta)^{s/2}u}_2^2 + 2^*_s \left( \frac{\mu}{2}\norm{u}_2^2 - \int_{\R^N}G(u)\right) =0} } 
\left( \half \norm{(-\Delta)^{s/2}u}_2^2- \int_{\R^N} G(u)\, + \frac{\mu}{2} \|u\|_2^2 \right);
$$
we recall that, when $s\in (\frac{1}{2},1)$ or $g \in C^{\sigma}_{loc}(\R)$ for some $\sigma >1-2s$, $p(\mu)$ is actually a ground state level.

The relation between the unconstrained and the constrained problem is an old-fashioned problem, which has been deeply investigated in a recent paper by Jeanjean and Lu \cite{JL0} in the case $s=1$.
We see that equation \eqref{eq_relaz_cons_uncons} gives an interesting relation between the two energy levels: this relation may be also reformulated by saying that
\begin{equation}\label{eq_trasf_legendre}
\kappa(m) = - p^*(m)
\end{equation}
where $p^*$ is the \emph{Legendre transform} of $a$. A relation of this type, but in a different framework, has been also obtained by Dovetta, Serra and Tilli in a very recent paper \cite{DST}. Here, relying on the convexity of the energy functions (due to the polynomial shape of $g$), they exploit \eqref{eq_trasf_legendre} in order to achieve interesting results.

We believe thus that \eqref{eq_trasf_legendre} could give more insights in the study of the relation between these two problems.

\chapter{Choquard-Hartree-Pekar %-Hartree 
equations: 
multiplicity of solutions}
%Multiple solutions for Choquard equations with general nonlinearities}
\label{chap_choq_multi}

In this Chapter we study the following nonlinear \emph{Choquard-Hartree-Pekar equation}
$$
	- \Delta u + \mu u =(I_\alpha*F(u)) F'(u) \quad \hbox{in}\ \mathbb{R}^N, 
$$
where $N\geq 3$, $\alpha\in (0,N)$,
$I_\alpha$ % (x) := \frac{\Gamma(\frac{N-\alpha}{2})}{\Gamma(\frac{\alpha}{2}) \pi^{N/2} 2^\alpha } \frac{1}{|x|^{N- \alpha}}$, $x \in \mathbb{R}^N \setminus\{0\}$ 
is the Riesz potential, and $F$ is an almost optimal subcritical nonlinearity. 
The goal is to prove existence of infinitely many solutions $u \in H^1_r(\mathbb{R}^N)$, by assuming $F$ odd or even.

We analyze the two cases: $\mu$ is a fixed positive constant or $\mu$ is unknown and the $L^2$-norm of the solution is prescribed, i.e. $\int_{\mathbb{R}^N} u^2 =m>0$.
Since the presence of the nonlocality prevents to apply the classical approach introduced by Berestycki and Lions in \cite{BL2}, we implement a new construction of multidimensional odd paths, %where some estimates for the Riesz potential play an essential role, 
and we find a nonlocal counterpart of their multiplicity result. 
In particular we extend the existence result in \cite{MS2}, due to Moroz and Van Schaftingen.

\medskip

This Chapter is mainly based on the paper \cite{CGT4}.

%%%%%%%%%%%%%%%%%%%%%%%%%%%%%%%%%%%%%%%%%%%%%%%%%%%%%%%

\section{Convolution with Riesz potential: a self-interaction
%The Hartree equation: a self-interaction 
%\tb{spostare qualcosa all'inizio} %COMMENT NOW
}
\label{sec_intro_choquard}

Given a nonlinearity $F \in C^1(\R,\R)$ and set $f:=F'$, we are interested to seek for multiple solutions $u \in H^1_r(\R^N)$ of the nonlocal equation
\begin{equation}\label{eq_Choquard_genericaF}
%\label{eq:1.1}
- \Delta u + \mu u =(I_\alpha*F(u)) f(u) \quad \hbox{in}\ \R^N, 
\end{equation}
where $N\geq 3$ and $\alpha\in (0,N)$. 
%and $I_\alpha:\, \R^N\setminus\{ 0\}\to \R$ is the Riesz potential
%defined by
% $$I_\alpha(x) := \frac{\Gamma(\frac{N-\alpha}{2})}{\Gamma(\frac{\alpha}{2}) \pi^{N/2} 2^\alpha } \frac{1}{|x|^{N- \alpha}};$$
%here $H^1_r(\R^N)$ denotes the space of radially symmetric Sobolev functions.
 %
%\\ 
In literature the semilinear equation \eqref{eq_Choquard_genericaF} with nonlocal source has several physical motivations and it is usually called nonlinear \emph{Choquard} (or \emph{Hartree}, or \emph{Pekar}) \emph{equation}.

In 1954 the equation \eqref{eq_Choquard_genericaF} with $N=3$, $\alpha=2$ and $F(s)=\half |s|^2$, that is 
\begin{equation}\label{eq:1.2}
-\Delta u + \mu u= \left(\frac{1}{4 \pi|x|} * |u|^2\right) u \quad \text{in $\R^3$},
\end{equation}
 was elaborated by Pekar in \cite{Pek0} (see also \cite{LR0}) to describe the quantum theory of a polaron at rest, through the use of the Newton potential %; in this case 
$\frac{1}{4 \pi|x|} $. % is also known as Newton potential. 
The idea of the convolution as a feature of interaction of a body with itself was exploited also by other authors: in 1976 it was arisen in the work \cite{Lie1} suggested by Choquard \cite{Cho0} on the modeling of an electron trapped in its own hole, in a certain approximation to Hartree-Fock theory of one-component plasma (see also \cite{FroL0,FTY,Stu1}). 
In 1996 the same equation was derived by Penrose in his discussion on the self-gravitational collapse of a quantum mechanical wave-function \cite{Pen1,Pen2,Pen3,MPT} (see also \cite{TM0,Tod0, FTY}) and in that context it is referred as \emph{Schr\"odinger-Newton system} (see \eqref{eq_SN_syst}). % and \eqref{eq_SN_sys}). 
See also Section \ref{sec_intro_choquard} for a derivation concerning exotic stars.

If $u$ is a solution of \eqref{eq:1.2}, then we notice that the wave function 
$$\psi(x,t) = e^{i \mu t} u(x), \quad (x,t) \in \R^3 \times [0, +\infty)$$
is a solitary wave of the time-dependent \emph{Hartree equation} \cite{Har0}
\begin{equation} \label{eq:1.3}
	i \psi_t = - \Delta \psi - \left(\frac{1}{4\pi |x|} * |\psi|^2\right) \psi \quad \text{in $\R^3 \times (0, +\infty) $};
\end{equation}
thus \eqref{eq:1.2} represents the stationary nonlinear Hartree equation.

\medskip

%\\
%We aim to analyze the two cases, already presented in Chapter \ref{chap_fract_normal}: $\mu$ is a fixed positive constant or $\mu$ is unknown and the mass of the solution, described by its $L^2$-norm, is prescribed.
As already pointed out in Chapter \ref{chap_fract_normal}, the study of standing waves of \eqref{eq:1.3} has been pursed in two main directions, which opened two different challenging research fields.%

A first topic regards the search for solutions of \eqref{eq:1.2} with a prescribed frequency $\mu$ and free mass, the so-called \emph{unconstrained} problem. 
The second line of investigation of the problem \eqref{eq:1.3} consists of prescribing the mass $m >0$ of $u$, thus conserved by $\psi$ in time
%\begin{equation*}\label{eq:1.4}
$$
\int_{\R^3} |\psi(x,t)|^2 \, dx= m \quad \forall \, t \in [0, +\infty),
$$
%\end{equation*}
and letting the frequency $\mu$ to be free. Such problem is usually said \emph{constrained}.

\smallskip

For the unconstrained problem, the first investigations for existence and symmetry of the solutions to \eqref{eq:1.2} go back to the works of Lieb \cite{Lie2} and Menzala \cite{Men0}, and also to \cite{Stu1,CSV,MPT} by means of ordinary differential equations techniques. 
We mention also the recent papers by Lenzmann \cite{Len1} and by Winter and Wei \cite{WeWi0} about the nondegeneracy of the unique radial solution of \eqref{eq:1.2}.

Variational methods were also employed to derive existence and qualitative results of standing wave solutions for more generic values of $\alpha \in (0,N)$ and of power type nonlinearities $F(t)= \frac{1}{p} |t|^p$: 
in particular Moroz and Van Schaftingen \cite{MS0} (see also \cite{MS3}) considered the special model
\begin{equation}\label{eq:1.5}
		- \Delta u + \mu u =(I_\alpha \ast |u|^p) |u|^{p-2} u \quad \hbox{in}\ \R^N, 
\end{equation}
and they proved that \eqref{eq:1.5} has solutions if 
\begin{equation}\label{eq:1.6}
2^{\#}_{\alpha}=\frac{N+\alpha}{N} < p < \frac{N+\alpha}{N-2} = 2^*_{\alpha}.
\end{equation}
When dealing with variational (and regular) solutions, they proved that range \eqref{eq:1.6} is optimal. 
Moreover in \cite{MS0} they showed that all positive ground states of \eqref{eq:1.5} are radially symmetric and monotone decreasing about some point and derived the decay asymptotics at infinity of such ground states (see \cite{CCS1} for $p \geq 2$, and also \cite{MZo0}).
Furthermore, in \cite{GV0,GMV,RV0} the authors study, for some values of $p$ and $\alpha$, least energy nodal solutions, odd with respect to a hyperplane; see also \cite{CCS1,ClSa,Wet0,WX0,XW0} for other results on sign-changing solutions with various symmetries and saddle type solutions.

Recently in \cite{MS2} Moroz and Van Schaftingen considered the problem \eqref{eq_Choquard_genericaF} when $F$ is a Berestycki-Lions type function under the following general assumptions:
 \begin{itemize}
 	\item[(F1)] \label{(F1)}
$F \in C^1(\R, \R)$;
 	\item[(F2)] \label{(F2)}
there exists $C >0$ such that, for every $s \in \R$, 
$$|s f(s)| \leq C \big(|s|^{2^{\#}_{\alpha}} + |s|^{2^*_{\alpha}}\big);$$
 	\item[(F3)] \label{(F3)}
 	$$\lim_{s \to 0} \frac{F(s)}{|s|^{2^{\#}_{\alpha}}} =0, \quad
 	\lim_{ s \to + \infty} \frac{F(s)}{|s|^{2^*_{\alpha}}} =0;$$
% 	\item[(F3)] 
% 	$$\lim_{s \to 0} \frac{s\tr{f}(s)}{|s|^{\frac{N+ \alpha}{N}}} =0, \quad
% 	\lim_{ s \to + \infty} \frac{s\tr{f}(s)}{|s|^{\frac{N+ \alpha}{N-2}}} =0;$$
 	\item[(F4)] \label{(F4)}
$F(s)\not\equiv 0$, that is, there exists $s_0 \in \R$, $s_0 \neq 0$ 
 such that $F(s_0) \neq0$.
 \end{itemize}
 
In particular they prove the following theorem (see \cite[Theorems 1 and 4]{MS2}).
%\tor{Cita teorema}
\begin{Theorem}[\cite{MS2}]\label{thm_MVS_cita1}
We have the following results.
\begin{itemize}
\item Assume \hyperref[(F1)]{\textnormal{(F1)}}--\hyperref[(F4)]{\textnormal{(F4)}}. Then there exists a ground state solution $u \in H^1(\R^N)$. Moreover $u\in W^{2,q}_{loc}(\R^N)$ for each $q \geq 1$ (in particular, $u$ is H\"older continuous);
\item Assume \hyperref[(F1)]{\textnormal{(F1)}}-\hyperref[(F2)]{\textnormal{(F2)}}, $f$ odd and with constant sign on $(0,+\infty)$. Then every ground state has strict constant sign (strictly positive or negative) and it is radially symmetric with respect to some point in $\R^N$.
\end{itemize}
\end{Theorem}

The qualitative result contained in Theorem \ref{thm_MVS_cita1} will be extended in this thesis to the case $f$ even, see Theorem \ref{th_INT_LOC_positiv}.
The existence of an infinite number of standing wave solutions to \eqref{eq:1.2} was instead faced by Lions in \cite{Lio1} (see also \cite{ClSa}); here the homogeneity of the source plays a crucial role in order to work on finite dimensional subspaces. Similar ideas have been applied in \cite{Ack0,QRT} in presence of more general sources satisfying Ambrosetti-Rabinowitz type conditions. 
We remark that all these multiplicity results deal with odd power nonlinearities $f$. 

To our knowledge it is still an open problem the existence of infinitely many radially symmetric solutions for the nonlinear Choquard equation \eqref{eq_Choquard_genericaF} under the optimal assumptions \hyperref[(F1)]{\textnormal{(F1)}}--\hyperref[(F4)]{\textnormal{(F4)}} and symmetric conditions on the nonlocal source term $(I_\alpha*F(u))f(u)$, and this is the aim of this Chapter. 
We note that this nonlinear term is odd both if $f$ is even or odd.

\medskip

Existence of a solution for the nonlinear Choquard equation \eqref{eq:1.5} under mass constraint has been obtained by Ye \cite{Ye1}; see also \cite{LiYe0} for odd powers-sum type functions.
More recently, Cingolani and Tanaka in \cite{CT1} obtained existence of a solution $u\in H^1_r(\R^N)$ to
\begin{equation}\label{eq:1.7}
	\left \{
	\begin{aligned}
		- \Delta u & + \mu u =(I_\alpha*F(u))f(u) \quad \hbox{in}\ \R^N, \\ 
		&\int_{\R^N} u^2 dx = m,
% \\
%		& u \in H^1_r(\R^N), 
	\end{aligned}
	\right. 
\end{equation}
assuming that $F$ satisfies \hyperref[(F1)]{\textnormal{(F1)}}, \hyperref[(F4)]{\textnormal{(F4)}} and it is $L^2$-subcritical, namely
\begin{itemize}
	\item[(CF2)] \label{(CF2)}
there exists $C >0$ such that, set $2^m_{\alpha}= \frac{N+ \alpha+2}{N}$, for every $s \in \R$, $$|s f(s)| \leq C \big(|s|^{2^{\#}_{\alpha}} + |s|^{2^m_{\alpha}}\big);$$
	\item[(CF3)] \label{(CF3)}
	$$\lim_{s \to 0} \frac{F(s)}{|s|^{2^{\#}_{\alpha}}} =0, \quad
	%$$\lim_{s \to 0} \frac{s\tr{f}(s)}{|s|^{\frac{N+ \alpha}{N}}} =0, \quad
	\lim_{ s \to + \infty} \frac{F(s)}{|s|^{2^m_{\alpha}}} =0.$$
\end{itemize}
The existence result in \cite{CT1} relies on %a new approach, based on a 
a Lagrangian formulation of the problem, in the spirit of Chapter \ref{chap_fract_normal}. 

Multiplicity of radial standing wave solutions to \eqref{eq:1.3} with prescribed $L^2$-norm has been instead faced again by Lions in \cite{Lio1} (see also \cite{CJ0} for the planar logarithmic Choquard equation); 
as regards instead the case of general nonlinearities $f$, recently Bartsch et al. \cite{BaLiLi} obtained the existence of infinitely many solutions of \eqref{eq:1.7} by assuming that $f$ is an odd function which satisfies monotonicity and Ambrosetti-Rabinowitz conditions. 
We highlight that the restriction on odd functions is not just a matter of symmetry of the functional, but it is related also to some sign restriction on the function $f$. 
The authors in \cite{BaLiLi} rely on mountain pass and Concentration-Compactness arguments, together with the use of a stretched functional, i.e. a functional in an augmented space which takes into consideration scaling properties and the Pohozaev identity.

It remains open the challenging problem of the existence of infinitely many solutions for the constrained nonlinear Choquard equation \eqref{eq:1.7} under optimal assumptions on the nonlinearity $f$, when monotonicity and Ambrosetti-Rabinowitz type conditions do not hold or $f$ is not odd.

\medskip

In the present Chapter we will give an affirmative answer to both the \emph{unconstrained} and \emph{constrained} problems when $F$ satisfies the general Berestycki-Lions type assumptions \hyperref[(F1)]{\textnormal{(F1)}}--\hyperref[(F4)]{\textnormal{(F4)}} and \hyperref[(F1)]{\textnormal{(F1)}}-\hyperref[(CF2)]{\textnormal{(CF2)}}-\hyperref[(CF3)]{\textnormal{(CF3)}}-\hyperref[(F4)]{\textnormal{(F4)}} respectively, together with the symmetric condition
\begin{itemize}
	\item[(F5)] \label{(F5)}
$F$ is odd or even.
\end{itemize}

We begin to notice that despite \cite{CT1}, where existence is investigated, to gain multiplicity the symmetry of the function $F$ plays a crucial role.
In particular, we assume $F$ to be \emph{odd} or \emph{even}, which guarantees the evenness of the energy functional associated to \eqref{eq_Choquard_genericaF}. 
We emphasize that the possibility to assume both the symmetries on $F$ is a particular feature of the nonlocal source: 
indeed, in the source-local case \cite{BL2,HT0} (see also Chapter \ref{chap_fract_normal}), the nonlinear term is usually assumed odd in order to get the symmetry of the functional. 
We mention the recent paper \cite{DMP} where the existence of a single nonradial solution to (\ref{eq_Choquard_genericaF}) is obtained under the condition \hyperref[(F5)]{\textnormal{(F5)}}.

\medskip

We start to analyze the constrained case, which appears, as usual, more delicate. 
By virtue of \cite{Pal0}, radially symmetric solutions to \eqref{eq:1.7} can be characterized as critical points of the $C^1$-functional $\mc{L}: H^1_r(\R^N) \to \R$
$$ \mathcal{L}(u) := \half \int_{\R^N} |\nabla u|^2\, dx -\half \int_{\R^N} (I_\alpha*F(u))F(u)\, dx,$$
constrained on the sphere 
$$ \mathcal{S}_m := \left \{ u \in H^1_r(\R^N) \mid \int_{\R^N} u^2 \, dx= m \right\}.$$
A possible approach to problem \eqref{eq:1.7} is to minimize $\mc{L}$ on the sphere $\mc{S}_m$, whenever the functional is here bounded. 
Nevertheless, in the spirit of Chapter \ref{chap_fract_normal}, for the general class of nonlinearities related to \cite{BL1,MS2}, considered in this thesis, we introduce a Lagrangian formulation of the nonlocal problem \eqref{eq:1.7}, extending a new approach introduced by Hirata and Tanaka \cite{HT0} for the local case. 
We highlight again the advantage of this method, that can be suitably adapted to derive multiplicity results of normalized solutions in several different frameworks. 

We recall here briefly the ideas of Chapter \ref{chap_fract_normal}. Writing $\R_+:=(0, +\infty)$, a solution $(\mu,u)\in \R_+ \times H^1_r(\R^N)$ of \eqref{eq:1.7} corresponds to a critical point of the functional $\mc{I}^m: \R_+ \times H^1_r(\R^N)\to \R$ defined by 
$$ \mc{I}^m(\mu, u):=\half \int_{\R^N} |\nabla u|^2\, dx -\half \int_{\R^N} (I_\alpha*F(u))F(u)\, dx+ \frac{\mu}{2} \left(\int_{\R^N} u^2 \, dx -m \right). $$
We seek for critical points $(\mu,u) \in \R_+ \times H^1_r(\R^N)$ of $\mc{I}^m$, namely weak solutions of $\partial_u \mc{I}^m(\mu, u)=0$ and $\partial_{\mu} \mc{I}^m(\mu, u)=0$.

Inspired by the Pohozaev identity, we introduce the Pohozaev functional $\mc{P}:\R_+ \times H^1_r(\R^N) \to \R$ by setting 
\begin{equation*}
\mc{P}(\mu, u) := \frac{N-2}{2} \int_{\R^N} |\nabla u|^2 \, dx + N \frac{\mu}{2} \int_{\R^N} u^2 \, dx -\frac{N+ \alpha}{2} \int_{\R^N} (I_\alpha*F(u))F(u) \, dx
\end{equation*}
and the \emph{Pohozaev set}
$$ \Omega :=\big\{(\mu,u) \in \R_+ \times H^1_r(\R^N) \mid \mc{P}(\mu,u)>0\big\} \cup\big\{(\mu,0) \mid \mu \in \R_+ \big\}. $$
We note that $\{(\mu,0) \mid \mu\in\R_+\}\subset {\it int}(\Omega)$ and thus
$$ \partial \Omega =\big\{(\mu,u) \in \R_+ \times H^1_r(\R^N) \mid \mc{P}(\mu,u)=0, \ u \nequiv 0 \big\}, $$
where the interior and the boundary are taken with respect to the topology of $\R_+ \times H^1_r(\R^N)$. Therefore $(\mu, u) \in \partial \Omega$ if and only if $u \nequiv 0$ satisfies the Pohozaev identity $\mc{P}(\mu, u)=0$. We recognize a Mountain Pass structure for the functional $\mc{I}^m$ in $\R_+\times H_r^1(\R^N)$, where the mountain is given by $\partial\Omega$.
We call $\partial\Omega$ a \emph{Pohozaev mountain} for $\mc{I}^m$. We emphasize that under assumptions \hyperref[(F1)]{\textnormal{(F1)}}-\hyperref[(F2)]{\textnormal{(F2)}}, if $u \in H^1_r(\R^N)$ solves $\partial_u \mc{I}^m(\mu, u)=0$ with $\mu \in \R_+$ fixed, then $\mc{P}(\mu, u)=0$.

Using a variant of the Palais-Smale condition \cite{HT0,IT0}, which takes into account the Pohozaev identity, we will prove a deformation theorem which enables us to apply minimax arguments in the product space $\R_+ \times H^1_r(\R^N)$. 
We will prove the existence of multiple $L^2$-normalized solutions detecting minimax structures in such product space. 

\medskip

We state our main results.
\begin{Theorem}\label{S:1.1_choq}	
Suppose $N\geq 3$, $\alpha \in (0, N)$ and \hyperref[(F1)]{\textnormal{(F1)}}-\hyperref[(CF2)]{\textnormal{(CF2)}}-\hyperref[(CF3)]{\textnormal{(CF3)}}-\hyperref[(F4)]{\textnormal{(F4)}}-\hyperref[(F5)]{\textnormal{(F5)}}.
\begin{itemize}
\item[(i)] For any $k \in \N$ there exists $m_k \geq 0$ such that for every $m > m_k$, the problem \eqref{eq:1.7} has at least $k$ pairs of nontrivial, distinct, radially symmetric solutions. 
\item[(ii)] Assume in addition an $L^2$-subcritical growth also at zero, i.e.
	\begin{itemize}
\item[\textnormal{(CF4)}] \label{(CF4)} 
		$$	
		\lim_{s \to 0} \frac{|F(s)|}{|s|^{2^m_{\alpha}}} = + \infty;
		$$
	additionally, if $F$ is odd, assume that there exists $\delta_0>0$ such that $F$ has a constant sign in $(0,\delta_0]$ and
% \tr{ci vuole il valore assoluto?}
		\begin{equation}\label{eq_cond_more_gen}
		\sup_{s \in (0, \delta_0],\, h\in [0,1]} \frac{F(s h)}{F(s)}< +\infty;
		\end{equation}
for example, this is satisfied if %that 
$\abs{F(s)}$ is assumed non-decreasing in $[0,\delta_0]$. % for some $\delta_0>0$.
	\end{itemize}
Then $m_k=0$ for each $k \in \N$, that is for any $m >0$ the problem \eqref{eq:1.7} has countably many pairs of solutions $(\mu_n, u_n)_{n}$ satisfying $\mathcal{L}(u_n) <0$, $n \in \N$. Moreover we have 
$$\mathcal{L}(u_n) \to 0 \quad \hbox{as $ n \to + \infty$}.$$
\end{itemize}
\end{Theorem}

\begin{Remark}
\label{rem_extra_cond}
%	We point out that, for $F$ odd, the monotonicity near the origin in \textnormal{(CF4)} can be slightly weakened, with no change in the proof, with the following condition:
%	\begin{itemize}
%		\item[] For some $\delta_0>0$, $F$ has a constant sign in $(0,\delta_0]$ and \tr{ci vuole il valore assoluto?}
%		\begin{equation}\label{eq_cond_more_gen}
%		\sup_{s \in (0, \delta_0],\, h\in [0,1]} \frac{F(s h)}{F(s)}=:M < +\infty;
%		\end{equation}
%	\end{itemize}
We comment condition \eqref{eq_cond_more_gen}. Set 
$$ M:=\sup_{s \in (0, \delta_0],\, h\in [0,1]} \frac{F(s h)}{F(s)}< +\infty$$
we have, when $\abs{F(s)}$ is non-decreasing, $M=1$. As a nontrivial example one can consider $\beta \in (2^{\#}_{\alpha},2^m_{\alpha})$ and $F$ oscillating near zero between $|s|^{\beta}$ and $2|s|^{\beta}$, so that $M\leq 2$; for instance the odd extension of
	$$F(s):= s^{\beta} \big(2+ \sin(\tfrac{1}{s})\big) \quad \hbox {as $s \to 0^+$}.$$
	If instead $F$ oscillates (not strictly) between $|s|^{\beta_1}$ and $|s|^{\beta_2}$, with $2^{\#}_{\alpha} < \beta_1 < \beta_2 < 2^m_{\alpha}$, then $M=+\infty$; thus for instance the odd extension of
	$$F(s):= s^{\beta_1} \big(1+ \sin(\tfrac{1}{s})\big) + s^{\beta_2} \big(1- \sin(\tfrac{1}{s})\big) \quad \hbox {as $s \to 0^+$}$$
	is not covered by \eqref{eq_cond_more_gen}.
\end{Remark}

\begin{Remark}\label{R:2.2}
We observe that, by substituting $F$ with $-F$, there is no loss of generality in assuming 
$$F(s_0)>0 \quad \hbox{for some $s_0\neq0$}$$
in \hyperref[(F4)]{\textnormal{(F4)}} ($s_0$ can be chosen positive if, for example, \hyperref[(F5)]{\textnormal{(F5)}} holds) and 
$$\lim_{s \to 0^+} \frac{F(s)}{|s|^{2^m_{\alpha}}} = + \infty$$
in \hyperref[(CF4)]{\textnormal{(CF4)}}. %, together with $F$ non-decreasing when it is odd. 
Thus, for the remaining part of the Chapter, we assume this positivity on the right-hand side of zero.
\end{Remark}

\smallskip

A key point of the argument is the construction of multidimensional odd paths. 
When $f$ satisfies some Ambrosetti-Rabinowitz condition (i.e., $F$ can be estimated from below by an homogeneous function $|t|^p$), the construction of such a path classically relies on the equivalence of the $H^1$-norm and the $L^p$-norm on finite dimensional subspaces of $H^1(\R^N)$. When such condition is no longer available, a finer construction is needed: in the celebrated paper \cite{BL2} Berestycki and Lions build this path for a local problem by exploiting an inductive process based on piecewise affine functions.

In our nonlocal case, in order to prove the existence of multiple solutions for $m\gg 0$ (point ($i)$ of Theorem \ref{S:1.1_choq}), unlike the elaborated approach of \cite{BL2} we can obtain the existence of a multidimensional odd path by exploiting the positivity of the Riesz potential functional. 
A similar approach can be implemented to gain existence of infinitely many solutions for any $m>0$ when $F$ is even (first part of point $(ii)$ of Theorem \ref{S:1.1_choq}), since in this case $F$ can be assumed positive in a neighborhood of the origin. %(up to substituting $F$ with $-F$). 
See anyway Remark \ref{rem_general_kernels} below.

A quite delicate issue, instead, comes up when $F$ is odd. 
%, which appears delicate in the case of nonlocal nonlinearities, especially when $f$ is even. 
Differently from \cite{HT0} and Chapter \ref{chap_fract_normal}, the classical argument given by \cite{BL2} cannot be applied directly in the context of nonlinear Choquard equations because of the presence of a nonlocal source, and we need to implement a new approach to gain the existence of an admissible odd path. 

To this aim we proceed by finding suitable annuli: using characteristic functions corresponding to the annuli, we construct our multidimensional odd paths. Here interactions between these characteristic functions produced by the Riesz potential play a crucial role, in particular the index $\alpha$ is related to the strength of interaction and the case $\alpha\in (0,1]$ reveals to be more delicate.
To this aim we use sharp estimates for the Riesz potential obtained by Thim \cite{Thim0}. % in an essential way.

\medskip

As a further byproduct of the previous approach %of the construction of multidimensional odd paths 
we gain the existence of infinitely many solutions for the unconstrained problem. 
More precisely, defined the $C^1$-functional $\mc{J}_{\mu}: H^1_r(\R^N) \to \R$ by setting
$$
\mc{J}_{\mu} (u):=\half \int_{\R^N} |\nabla u|^2 \, dx + \frac{\mu}{2} \int_{\R^N} u^2 \, dx 
-\half \int_{\R^N} (I_\alpha*F(u))F(u)\, dx,
$$
we establish the following result.
\begin{Theorem}\label{S:1.3}	
Suppose $N\geq 3$, $\alpha \in (0, N)$ and $\mu>0$ fixed. Assume that \hyperref[(F1)]{\textnormal{(F1)}}--\hyperref[(F5)]{\textnormal{(F5)}} hold. Then there exist countably many radial solutions $(u_n)_{n}$ of the nonlinear Choquard equation \eqref{eq_Choquard_genericaF}. Moreover we have 
$$\mathcal{J}_{\mu}(u_n) \to +\infty \quad \hbox{as $ n \to + \infty$}.$$
\end{Theorem}

Our multiplicity result is the counterpart of what done in \cite{BL2} for the local case with odd nonlinearities and extend the existence result in \cite{MS2} due to Moroz and Van Schaftingen. 

%\tor{
%\begin{Remark}\label{rem_general_kernels}
%We highlight that the construction of a multidimensional odd path in the case of $F$ odd (point $(ii)$ of Theorem \ref{S:1.1_choq}) may be adapted to other different frameworks also in the case of $F$ even, or even for the unconstrained problem: indeed, when considering different kernels $K=K(x,y)$ which do not makes the functional
%$$g \mapsto \int_{\R^N} \int_{\R^N} K(x,y) g(x) g(y)$$
% positive (for example, $K(x,y)$ sign-changing kernels; see also Section \ref{sec_concent_sphere}), the easy approach based on Proposition \ref{prop_positivity_Riesz} do not work anymore. This is an interesting line of research for the future.
%\end{Remark}
%}

\begin{Remark}\label{rem_general_kernels}
We highlight that the easier approach for building a multidimensional path, based on the positivity of the Riesz kernel (Proposition \ref{prop_positivity_Riesz}), cannot generally be applied to more generally frameworks (also if $F$ is even); for examples, when dealing with kernels $K=K(x,y)$ which do not makes the functional
$$g \mapsto \int_{\R^N} \int_{\R^N} K(x,y) g(x) g(y)$$
 positive (for example, $K(x,y)$ sign-changing). %; see \eqref{eq_kernel_signchang} and \eqref{eq_logar_kernel}). % also Section \ref{sec_concent_sphere}). 
In this case, the approach here developed, based on suitable annuli, might instead be adapted. 
This is an interesting line of research for the future.
\end{Remark}

\smallskip

%\smallskip

The Chapter is organized as follows. 
%\tr{aggiusta label sezioni} %COMMENT NOW
In Section \ref{section:4..1} we focus on the construction of multidimensional paths, by dealing first with an easier version based on the positivity of the Riesz potential, and then a refined version based on some suitable annuli and essential interaction estimates for non-local terms. 
Section \ref{sec_choquard_asympt_MP} is then dedicated to the study of the asymptotic behaviour of the mountain pass values, according to variable values of $\mu$. 
Afterwards, in Section \ref{section:4.3}, we detect a mountain pass structure, built on the Pohozaev mountain, for the constrained case, and in Section \ref{sec_PSP_choquard} we derive a Palais-Smale-Pohozaev condition. 
In Section \ref{section:3.2} we introduce an augmented functional which will be used to gain a deformation lemma, and we further study suitable minimax values defined through the tool of the genus which allows to prove the main Theorem \ref{S:1.1_choq}.
Finally in Section \ref{sec_unconstrained_choquard} we deal with the unconstrained case by proving Theorem \ref{S:1.3}.

%%%%%%%%%%%%%%%%%%%%%%%%%%%%%%%%%%%%%%%%%%%%%

%\setcounter{equation}{0} %FOR Arxiv
\section{%Construction of 
Multidimensional annuli-shaped %odd 
paths: even and odd nonlinearities}
\label{section:4..1}

%\subsubsection{Construction of multidimensional odd paths}

In this Chapter we briefly denote by $q$ the lower-critical exponent $2^{\#}_{\alpha}$ and by $p$ the $L^2$-critical exponent $2^m_{\alpha}$, i.e.
$$q:=2^{\#}_{\alpha}=\frac{N+\alpha}{N}, \quad p:=2^m_{\alpha}= \frac{N+ \alpha+2}{N}.$$

Again, to avoid problems with the boundary of $\R_+$, we write from now on (see Section \ref{sec_ground_state} for a different approach)
$$\mu \equiv e^{\lambda}\in (0, +\infty), \quad \lambda \in \R.$$
We also set
$$	\mc{D}(u) := \mc{D}_{\alpha}(F(u),F(u))=\int_{\R^N} (I_\alpha*F(u))F(u)\, dx.	$$
Using Proposition \ref{prop_HLS} and \hyperref[(F1)]{\textnormal{(F1)}}-\hyperref[(F2)]{\textnormal{(F2)}}, we notice that $\mc{D}$ is continuous on $L^2(\R^N) \cap L^{2^*}(\R^N)$, where $2^*= \frac{2N}{N-2}$ is the Sobolev critical exponent, and thus continuous on $H^1_r(\R^N)$; notice that if we assume \hyperref[(CF2)]{\textnormal{(CF2)}}, then $\mc{D}$ is continuous also on $L^2(\R^N) \cap L^{2 + \frac{4}{N+\alpha}}(\R^N)$.

To deal with the unconstrained problem, we further define the $C^1$-functional 
$\mc{J}:\R \times H^1_r(\R^N) \to \R$ by setting %\tor{scrivi $\mc{J}_{\lambda}$}
\begin{equation}\label{eq:2.3}
\mc{J} (\lambda, u):=\half \norm{\nabla u}_2^2 -\half \mc{D}(u) + \frac{e^{\lambda}}{2} \|u\|_2^2, \quad (\lambda, u) \in \R\times H^1_r(\R^N).
\end{equation}
For a fixed $\lambda \in \R$, $u \in H^1_r(\R^N)$ is critical point of $\mc{J}(\lambda, \cdot)$ if and only if $u$ solves (weakly)
%\begin{equation}
%\left \{
%\begin{aligned}
%- \Delta u & + e^{\lambda} u =(I_\alpha*F(u))f(u) \quad \hbox{in}\ \R^N, \\ 
%& u \in H_r^1(\R^N);
%\end{aligned}
%\right. \label{eq:2.4}
%\end{equation}
\begin{equation} \label{eq:2.4}
- \Delta u + e^{\lambda} u =(I_\alpha*F(u))f(u) \quad \hbox{in}\ \R^N.
\end{equation}
%here by \emph{solution} we will always mean weak solution.

%\medskip

In this Section we study the geometry of
$$	u \in H_r^1(\R^N)\mapsto \mc{J}(\lambda,u) \in \R,	$$
for a fixed $\lambda\in\R$. We introduce a sequence of minimax values $a_n(\lambda)$, $n \in \N^*$: %=1,2,\dots$. 
these values play important roles to find multiple solutions for the constrained problem (Theorem \ref{S:1.1_choq}) as well as for the unconstrained problem (Theorem \ref{S:1.3}).

For $n\in\N^*$ and $\lambda \in \R$ we introduce the set of paths
$$\Gamma_n(\lambda):=\big\{\gamma\in C(D_n, H^1_r(\R^N)) \mid \gamma \hbox{ odd}, \, \mc{J}(\lambda, \gamma_{|\partial D_n})<0 \big\}$$
and the minimax values
$$a_n(\lambda):=\inf_{\gamma \in \Gamma_n(\lambda)} \sup_{\xi \in D_n} \mc{J}(\lambda, \gamma(\xi)).$$

For $n\geq 2$ the nonemptiness of $\Gamma_n(\lambda)$ has to be checked; for $n=1$ we refer to \cite[claim 1 of Proposition 2.1]{MS2}. 
Classically, in the local framework this fact was proved in \cite{BL2} by constructing inductively piecewise affine paths. This construction does not fit the nonlocality interaction given by the Choquard term, thus we need another approach. 

\begin{Proposition}\label{S:4.1}
Assume \hyperref[(F1)]{\textnormal{(F1)}}--\hyperref[(F4)]{\textnormal{(F4)}} and %\tr{and (F5) (it is sufficient 
$F(\pm s_0) \neq 0$.
Let $n \in \N^*$ and $\lambda \in \R$. 
Then $\Gamma_n(\lambda)\neq \emptyset$, thus $a_n(\lambda)$ is well defined. Moreover, $a_n(\lambda)>0$ and it is increasing with respect to $\lambda$ and $n$.
\end{Proposition}

\claim Proof.
Start observing that the polyhedron
 $$ \Sigma:=\left\{ t=(t_1,\dots,t_n) \mid \max_{i=1,\dots, n}\abs{t_i}=1\right\}
 $$
is homeomorphic to $\partial D_n$ (we passed from the $L^2$ to the $L^{\infty}$ norm). 
Let us fix $e_1, \dots, e_n \in C^{\infty}_c(\R^N)$, each of them between $0$ and $1$, %qui questa condizione non serve,
radially symmetric, equal to one in some annulus $A_i$, and such that their supports are mutually disjoint. Then set $\gamma: \Sigma \to H^1_r(\R^N)$ by
\begin{equation}\label{eq_cammino_easy}
\gamma(t)(x):=s_0 \sum_{i=1}^n t_i e_i(x)
\end{equation}
for every $t = (t_1, \dots, t_n) \in \Sigma$ and $x \in \R^N$. %, where $s_0$ appears in \eqref{...}. 
The map $\gamma$ is clearly odd and continuous. Moreover every $t \in \Sigma$ has at least a nontrivial component $|t_i|=1$, thus we have $F(\gamma(t)(x)) = F(s_0 t_i e_i(x)) = F( \pm s_0) \neq 0$ on $A_i$, hence $F(\gamma(t)) \nequiv 0$. By Proposition \ref{prop_positivity_Riesz} we have
$$\mc{D}(\gamma(t)) >0 \quad \hbox{ for each $t \in \Sigma$}.$$
Since
$$\mc{D}\circ \gamma : \Sigma \to \R$$
is continuous, and $\Sigma$ is compact, we obtain
$$ \min_{t \in \Sigma} \mc{D}(\gamma(t)) =: C>0,$$
i.e. 
$$
\mc{D}(\gamma(t)) \geq C >0 \quad \hbox{ for each $t \in \Sigma$}.
$$
Set moreover $M:= \max_{\Sigma} \norm{\gamma}_{H^1}^2 \in \R.$
By scaling, we obtain
\begin{align*}
\mc{J}_{\lambda}(\gamma(t)(\cdot/\theta))
&= \frac{\theta^{N-2}}{2} \norm{\nabla \gamma(t)}_2^2 + \frac{\theta^N e^{\lambda} }{2} \norm{\gamma(t)}_2^2 - \frac{\theta^{N+\alpha}}{2} \mc{D}(\gamma(t)) \\
&\leq \frac{\theta^{N-2}}{2} M + \frac{\theta^N e^{\lambda} }{2} M - \frac{\theta^{N+\alpha}}{2} C < 0
\end{align*}
for some $\theta=\theta^* \gg 0$. Thus we consider $\tilde{\gamma}:=\gamma(\cdot)(\cdot/\theta^*) : \partial D_n \to H^1_r(\R^N)$. Finally we extend $\tilde{\gamma}$ to $D_n$ by
$$\tilde{\gamma}(\xi):= |\xi| \gamma\left(\frac{\xi}{|\xi|}\right)$$
for every $\xi \in D_n \setminus \{0\}$, and $\tilde{\gamma}(0):=0$. Therefore $\tilde{\gamma} \in \Gamma_n(\lambda) \neq \emptyset$.

What remains to prove is the monotonicity and positivity of $a_n(\lambda)$.
Since $D_n\subset D_{n+1}$, we may regard for $\gamma\in \Gamma_{n+1}(\lambda)$,
	$$	\gamma_{|{D_n}}\in \Gamma_n(\lambda).	$$
Thus we have $a_n(\lambda)\leq a_{n+1}(\lambda)$. Since $\mc{J}(\lambda,u)$ is monotone
in $\lambda$, we also have the monotonicity with respect to $\lambda$.

The positivity of $a_1(\lambda)$ is essentially obtained in \cite{MS2} {(see also \cite{CT1})}. Thus
$$ a_n(\lambda) \geq a_1(\lambda)>0.			 
\QED
$$

\bigskip

In the proof of Proposition \ref{S:4.1} we hardly relied on the positivity of the Riesz potential functional given in Proposition \ref{prop_positivity_Riesz}, to obtain the existence of path $\gamma: D_n \to H^1_r(\R^N)$ and a $C>0$ such that
\begin{equation}\label{eq_stima_path_C}
\mc{D}(\gamma(\xi)) \geq C >0 \quad \hbox{ for each $\xi \in \partial D_n$}.
\end{equation}
Notice moreover that this $\gamma$ satisfies $\gamma(\theta \xi)= \theta\gamma( \xi)$ for any $\xi \in D_n$ and $\theta \in \R$. Anyway, no good information on the constant $C$ appearing in \eqref{eq_stima_path_C} are given by this result.

A useful estimate in order to get infinitely many solutions for any $m>0$, when \hyperref[(CF4)]{\textnormal{(CF4)}} holds, is the one which relates $\mc{D}(\theta\gamma)$ to $F(\theta s_0)$ (see Lemma \ref{S:4.7} and Section \ref{sec_choquard_asympt_MP}), that is
\begin{equation}\label{eq_stima_path_F}
\mc{D}(\theta\gamma( \xi)) \geq C (F(\pm \theta s_0))^2 \quad \hbox{ for each $\xi \in \partial D_n$ and $\theta \in [0,1]$}
\end{equation}
for some uniform $C>0$. When $F$ is positive in a neighborhood of the origin (which is the case of $F$ even and \hyperref[(CF4)]{\textnormal{(CF4)}}), then one can build a suitable $\gamma$ which satisfies \eqref{eq_stima_path_F}.
\begin{Proposition}\label{prop_F_posit_intorn}
Assume \hyperref[(F1)]{\textnormal{(F1)}}--\hyperref[(F4)]{\textnormal{(F4)}}. % \tr{and (F5)}. 
Assume moreover that $F$ is positive in some $[-s_0, s_0]$, $F(\pm s_0) \neq 0$. Let $n \in \N^*$ and $\lambda \in \R$. Then the path $\gamma \in \Gamma_n(\lambda)$ defined in \eqref{eq_cammino_easy} satisfies \eqref{eq_stima_path_F}.
\end{Proposition}
\claim Proof.
%Let us fix $e_1, \dots, e_n \in C^{\infty}_c(\R^N)$, each of them between $0$ and $1$, radially symmetric, equal to one in some ball $B_i$, and such that their supports are mutually disjoint. Then set $\gamma: \Sigma \to H^1_r(\R^N)$ by
%$$\gamma(\xi)(x):= s_0 \sum_{i=1}^n \sgn(\xi_i) e_i(x)$$
%for every $\xi = (\xi_1, \dots, \xi_n) \in \partial D_n$ and $x \in \R^N$. The map $\gamma$ is clearly odd and continuous. Moreover
Assume the notation of the proof of Proposition \ref{S:4.1}. For each $t\in \Sigma$, there exists $|t_{k}|=1$, thus, by exploiting that $\theta s_0 t_i e_i \in [-s_0, s_0]$ we obtain 
%\begin{align*}
%\mc{D}(\gamma(\theta t)) &= \int_{\R^N} \int_{\R^N} I_{\alpha}(x-y) F\left( \theta s_0\sum_{i=1}^n \sgn(t_i) e_i (x)\right) F\left(\theta s_0 \sum_{j=1}^n \sgn(t_j) e_j(x)\right) \\
%& = \sum_{i=1}^n \sum_{j=1}^n \int_{\R^N} \int_{\R^N}I_{\alpha}(x-y) F\big(\theta s_0 \sgn(t_i) e_i (x)\big) F\big( \theta s_0 \sgn(t_j) e_j(x)\big) \\
%& \geq \int_{\R^N} \int_{\R^N} I_{\alpha}(x-y) F\big(\theta s_0 \sgn(t_k) e_k (x)\big) F\big( \theta s_0 \sgn(t_k) e_k(y)\big) \\
%& = \int_{\R^N} \int_{\R^N} I_{\alpha}(x-y) F\big( \pm \theta s_0 e_k (x)\big) F\big( \pm \theta s_0 e_k(y)\big) \\
%& \geq \int_{B_k} \int_{B_k} I_{\alpha}(x-y) F( \pm \theta s_0 ) F( \pm \theta s_0 ) \\
%\end{align*}
%where we exploited that $\theta s_0 t_i e_i \in [-s_0, s_0]$
\begin{align*}
\mc{D}(\theta\gamma( t)) &= \int_{\R^N} \int_{\R^N} I_{\alpha}(x-y) F\left( \theta s_0\sum_{i=1}^n t_i e_i (x)\right) F\left(\theta s_0 \sum_{j=1}^n t_j e_j(x)\right) \\
& = \sum_{i=1}^n \sum_{j=1}^n \int_{\R^N} \int_{\R^N}I_{\alpha}(x-y) F\big(\theta s_0 t_i e_i (x)\big) F\big( \theta s_0 t_j e_j(x)\big) \\
& \geq \int_{\R^N} \int_{\R^N} I_{\alpha}(x-y) F\big(\theta s_0 t_k e_k (x)\big) F\big( \theta s_0 t_k e_k(y)\big) \\
% & = \int_{\R^N} \int_{\R^N} I_{\alpha}(x-y) F\big( \pm \theta s_0 e_k (x)\big) F\big( \pm \theta s_0 e_k(y)\big) \\
& \geq % \int_{B_k} \int_{B_k} I_{\alpha}(x-y) F( \pm \theta s_0 ) F( \pm \theta s_0 ) = 
\big( F( \pm \theta s_0 )\big)^2 \int_{A_k} \int_{A_k} I_{\alpha}(x-y) 
\end{align*}
which is the claim.
\QED

\bigskip

When $F$ is odd (and thus it cannot be positive around the origin) it seems not an easy task to build a $\gamma \in \Gamma_n(\lambda)$ satisfying \eqref{eq_stima_path_F}; indeed some estimate from below on 
$$\int_{\R^N} \int_{\R^N} I_{\alpha}(x-y) \frac{F\big(\theta \gamma(\xi)(x)\big) F\big(\theta \gamma(\xi)(y)\big)}{\big(F(\theta s_0)\big)^2}$$
uniform for $\theta \to 0$ seems required; this is related to quotients of the type $\frac{F(s h)}{F(s)}$ with $s \in (0,s_0]$ and $h \in [0,1]$. This is essentially the meaning of condition \eqref{eq_cond_more_gen}.

\bigskip

To deal with this case we need a deep understanding of the Riesz potential on radial functions. We thus give now a different construction for a $\gamma \in \Gamma_n(\lambda)$: this procedure might be investigated also for more general Choquard-type equations, where different kernels (possibly sign-changing) appear.

We start by recalling a result contained in \cite[Theorem 1]{Thim0} (see also \cite[Lemma 6.3]{MMV} and references therein).

\begin{Theorem}[\cite{Thim0}]\label{S:4.2}
Let $\alpha \in (0, N)$ and $u\in L^1(\R^N) \cap L^{\infty}(\R^N) % H_r^1(\R^N)
$ be radial. Then $I_{\alpha}*u$ is radial and
\begin{equation}
 (I_\alpha * u)(r)=r^\alpha \int_0^\infty F_\alpha\left(\frac{r}{\rho}\right)\left(\frac{\rho}{r}\right)^\alpha
 u(\rho)\, \frac{d\rho}{\rho}, \label{eq:4.1}
\end{equation}
where $F_\alpha$ is positive and it satisfies, for some constants $C_{N,0}$, $C_{N,\infty}$, $C_{N,\alpha}>0$, 
$$F_{\alpha}(s) \to C_{N,0}\ \hbox{ as $s \to 0$}, \quad \frac{F_{\alpha}(s)}{s^{\alpha-N}} \to C_{N,\infty}\ \hbox{ as $s \to +\infty$}$$
and 
\begin{equation}
\frac{F_\alpha(s)}{G_\alpha(s)} \to 1 \quad \hbox{as}\ s\to 1, 
 \label{eq:4.2}
\end{equation}
with
\begin{equation}
G_\alpha(s):=
\begin{cases}
 C_{N,\alpha} &\hbox{if $\alpha\in (1,N)$}, \cr
 {C_{N,\alpha}}
 \abs{\log\abs{s-1}} &\hbox{if $\alpha=1$}, \cr
 {C_{N,\alpha}}
 \abs{s-1}^{\alpha-1} &\hbox{if $\alpha\in (0,1)$}. \cr
\end{cases}
 \label{eq:4.3}
\end{equation}
\end{Theorem}

For a proof of Proposition \ref{S:4.1}, we prepare some notation and some estimates. We introduce the annuli
$$A(R,h):=\big\{x\in\R^N \mid \abs x\in [R-h,R+h]\big\}, \quad \chi(R,h;\cdot):= \chi_{A(R,h)}$$
% \begin{align*}
% &A(R,h):=\big\{x\in\R^N \mid \abs x\in [R-h,R+h]\big\}, \cr
% &\chi(R,h;x):=
%\begin{cases}
% 1 &\hbox{for $x\in A(R,h)$}, \cr
% 0 &\hbox{otherwise},
%\end{cases}
%\end{align*}
for any $R\gg h >0$.
We have the following key estimates.

\medskip

\begin{Lemma}\label{S:4.3}
It results as $h\to 0$
 $$ \RRint I_\alpha(x-y)\chi(1,h;x)\chi(1,h;y)\,dxdy
 \sim
\begin{cases}
 h^2 &\hbox{if $\alpha\in(1,N)$}, \cr
 h^2\abs{\log h} &\hbox{if $\alpha=1$}, \cr
 h^{1+\alpha} &\hbox{if $\alpha\in (0,1)$}. \cr
\end{cases}
 $$
\end{Lemma}
%
%Here we recall that $f\sim g$ if there exist constants $C_1$, $C_2>0$ independent of $h$ such that
% $$ C_1g(h) \leq f(h) \leq C_2g(h) \quad \hbox{for small $h$}. $$

%We postpone a proof of Lemma \ref{S:4.3} and give it in Section \ref{section:4.6}.

\medskip

\claim Proof. % of Lemma \ref{S:4.3}.
We apply Theorem \ref{S:4.2} to $u(|x|)=\chi(1,h;\abs x)$. In particular, by \eqref{eq:4.1} we have
%\begin{align*}
% S_h &:= \RRint I_\alpha(x-y)u(x)u(y)\,dxdy \cr
% &= C\int_0^\infty (I_\alpha* u)(r) u(r) r^{N-1}\, dr \cr
% &= C\int_0^\infty \int_0^\infty F_\alpha\Big(\frac{r}{\rho}\Big)
% \rho^{\alpha-1}r^{N-1} u(\rho)u(r)\, d\rho dr \cr
% &= C \iint_{[1-h,1+h]^2} F_\alpha\Big(\frac{r}{\rho}\Big)\rho^{\alpha-1}r^{N-1}\, d\rho dr. 
%\end{align*}
\begin{align*}
 S_h &:= \RRint I_\alpha(x-y)u(x)u(y)\,dxdy 
= C\int_0^\infty (I_\alpha* u)(r) u(r) r^{N-1}\, dr \cr
 &= C\int_0^\infty \int_0^\infty F_\alpha\left(\frac{r}{\rho}\right)
 \rho^{\alpha-1}r^{N-1} u(\rho)u(r)\, d\rho dr 
= C \iint_{[1-h,1+h]^2} F_\alpha\left(\frac{r}{\rho}\right)\rho^{\alpha-1}r^{N-1}. %\, d\rho dr. 
\end{align*}
First we note that
 $$ \sup_{\rho,r\in [1-h,1+h]} \pabs{\frac{r}{\rho}-1} \to 0 \quad \hbox{as}\ h\to 0. $$
%Viene usata quest'ultima? CAPISCI \tr{
We consider the following three cases separately:
 $$ \hbox{(i)} \ \alpha\in (1,N), \quad \hbox{(ii)}\ \alpha=1, \quad
 \hbox{(iii)}\ \alpha\in (0,1). $$

\noindent
(i) When $\alpha\in (1,N)$ we may assume $F(\frac{r}{\rho})\sim C_{N,\alpha}>0$.
Thus
$$
 S_h \sim 
 \iint_{[1-h,1+h]^2}\rho^{\alpha-1}r^{N-1}\,d\rho dr
 \sim h^2.
$$

\noindent
(ii) When $\alpha=1$
\begin{align*}
 F_{\alpha}\left(\frac{r}{\rho}\right)
 &\sim G_1\left(\frac{r}{\rho}\right) = C_{N,1}
 \Big |\log\Big|\frac{r}{\rho}-1\Big|\Big| \cr
 & \sim
\big| \log\abs{r-\rho}-\log\rho \big| 
 =
 -\log\abs{r-\rho}+\log\rho. 
\end{align*}
Thus
 $$ S_h \sim 
 \iint_{[1-h,1+h]^2} (-\log\abs{r-\rho}+\log\rho) r^{N-1}\, d\rho dr.
 $$
Set
 \begin{align*}
 A_h &:= \big\{ (\rho,r)|\, \abs{\rho-r}\leq \tfrac{1}{2} h,\, \abs{r-1}\leq \tfrac{1}{2} h\big\}, \\
 B_h &:= \big\{ (\rho,r)|\, \abs{\rho-r}\leq 2 h,\, \abs{r-1}\leq h\big\},
 \end{align*}
we have
 $$ A_h \subset [1-h,1+h]^2 \subset B_h. $$
Hence for some $C$, $C'>0$
 \begin{equation}\label{eq:4.24} 
C \iint_{A_h} (-\log\abs{r-\rho}+\log\rho) r^{N-1}\, d\rho dr
		\leq S_h \leq C'
 \iint_{B_h} (-\log\abs{r-\rho}+\log\rho) r^{N-1}\, d\rho dr.
 \end{equation}
We compute
 \begin{align*}
 &\iint_{B_h} (-\log\abs{r-\rho}+\log\rho) r^{N-1}\, d\rho dr \\
 &\leq \iint_{B_h} (-\log\abs{r-\rho}+\log(1+h)) (1+h)^{N-1}\, d\rho dr \\
 &= \iint_{[-2h,2h]\times[1-h,1+h]} (-\log\abs{\tau}+\log(1+h)) (1+h)^{N-1}\, d\tau dr \\
 &= 4h(1+h)^{N-1} \int_0^{2h} (-\log\tau)\,d\tau + 8h^2 (1+h)^{N-1}\log(1+h)\\
 &= 4h(1+h)^{N-1} \big(-2h\log(2h)+2h\big) + 8h^2(1+h)^{N-1}\log(1+h) \\
 &\leq C''
 h^2\abs{\log h} \quad \hbox{as $h\to 0$}. 
 \end{align*}
Similarly we have
 $$	\iint_{A_h} (-\log\abs{r-\rho}+\log\rho)r^{N-1}\, d\rho dr
		\geq C''' h^2\abs{\log h}, $$
from which we obtain
$$ S_h \sim h^2\abs{\log h} \quad \hbox{as}\ h \to 0.$$

\medskip

\noindent
(iii) When $\alpha\in (0,1)$
 $$ F_\alpha\left(\frac{r}{\rho}\right) \sim G_\alpha\left(\frac{r}{\rho}\right)
 =C_{N,\alpha} \pabs{\frac{r}{\rho}-1}^{\alpha-1}.
 $$
Thus 
$$
 S_h\sim \iint_{[1-h,1+h]^2} \pabs{\frac{r}{\rho}-1}^{\alpha-1}
 \rho^{\alpha-1} r^{N-1}\, d\rho dr 
 = \iint_{[1-h,1+h]^2} \pabs{r-\rho}^{\alpha-1}r^{N-1}\, d\rho dr.
$$
Since
 $$ C \iint_{A_h}\abs{r-\rho}^{\alpha-1}(1-h)^{N-1}\,d\rho dr \leq S_h
 \leq C' \iint_{B_h} \abs{r-\rho}^{\alpha-1}(1+h)^{N-1}\,d\rho dr, $$
we have as in \eqref{eq:4.24}
	$$	S_h \sim h^{1+\alpha} \quad \hbox{as}\ h\to 0.	$$
This completes the proof. 
\QED

\bigskip

We show how to use it to build a continuous odd map in $L^2(\R^N) \cap L^{2^*}(\R^N)$. By a regularization argument, we will obtain a map in $\Gamma_n(\lambda)$.

By scaling, we have
\begin{align*}
\lefteqn{\RRint I_\alpha(x-y)\chi(R,h;x)\chi(R,h;y)\, dxdy}\cr
 &= R^{N+\alpha}\RRint 
 I_\alpha(x-y)\chi\Big(1,\frac{h}{R};x\Big)
 \chi\Big(1,\frac{h}{R};y\Big)\, dxdy\cr
 &\sim 
\begin{cases}
 R^{N+\alpha}(\frac{h}{R})^2 &\hbox{if $\alpha\in (1,N)$}, \cr
 R^{N+1} (\frac{h}{R})^2\abs{\log\frac{h}{R}} &\hbox{if $\alpha=1$}, \cr
 R^{N+\alpha}(\frac{h}{R})^{1+\alpha} &\hbox{if $\alpha\in (0,1)$}.
 \end{cases}
\end{align*}
For $R\geq 2$, we set the thickness of the annuli as
 $$ h_R :=
\begin{cases}
 R^{-\frac{N-2+\alpha}{2}} &\hbox{if $\alpha\in (1,N)$}, \cr 
 R^{-\frac{N-1}{2}}(\log R)^{-1/2} &\hbox{if $\alpha=1$}, \cr
 R^{-\frac{N-1}{1+\alpha}} &\hbox{if $\alpha\in (0,1)$,}
\end{cases}
 $$
%so that
%\begin{equation}
% \RRint I_\alpha(x-y)\chi(R,h_R;x)\chi(R,h_R;y)\,dxdy
% \in [C_{01}, C_{02}] \quad \hbox{for}\ R\geq 2, \label{eq:4.4}
%\end{equation}
%where $C_{01}$, $C_{02}>0$ are independent of $R\geq 2$.
%We check \eqref{eq:4.4} only for $\alpha=1$. 
%We have
%\begin{eqnarray*}
%\lefteqn{\quad\RRint I_\alpha(x-y)\chi\Big(1,\frac{h_R}{R};x\Big)\chi\Big(1,\frac{h_R}{R};y\Big)\, dxdy }\cr
% &&\sim R^{N+1}\Big(\frac{h_R}{R}\Big)^2 \Big|\log\Big(\frac{h_R}{R}\Big)\Big| \cr
% &&= R^{N+1} \left(\frac{R^{-\frac{N-1}{2}}\abs{\log R}^{-1/2}}{R}\right)^2
% \pabs{\log\left(\frac{R^{-\frac{N-1}{2}}\abs{\log R}^{-1/2}}{R}\right)} \cr
% &&= (\log R)^{-1} \pabs{ \log\Big(R^{-\frac{N+1}{2}}(\log R)^{-1/2}\Big)} \cr
% &&= (\log R)^{-1} \left( \frac{N+1}{2}\log R +\half\log (\log R) \right) \cr
% &&\to \frac{N+1}{2} \quad \text{as}\ R\to \infty. 
%\end{eqnarray*}
so that a uniform bound is gained.
\begin{Corollary}
We have
\begin{equation}
 \RRint I_\alpha(x-y)\chi(R,h_R;x)\chi(R,h_R;y)\,dxdy
 \in [C_{01}, C_{02}] \quad \hbox{for}\ R\geq 2, \label{eq:4.4}
\end{equation}
where $C_{01}$, $C_{02}>0$ are independent of $R\geq 2$.
\end{Corollary}

\claim Proof.
We check \eqref{eq:4.4} only for $\alpha=1$. 
We have
\begin{eqnarray*}
\lefteqn{\quad\RRint I_\alpha(x-y)\chi\Big(1,\frac{h_R}{R};x\Big)\chi\Big(1,\frac{h_R}{R};y\Big)\, dxdy }\cr
 &&\sim R^{N+1}\Big(\frac{h_R}{R}\Big)^2 \Big|\log\Big(\frac{h_R}{R}\Big)\Big| \cr
 &&= R^{N+1} \left(\frac{R^{-\frac{N-1}{2}}\abs{\log R}^{-1/2}}{R}\right)^2
 \pabs{\log\left(\frac{R^{-\frac{N-1}{2}}\abs{\log R}^{-1/2}}{R}\right)} \cr
 &&= (\log R)^{-1} \pabs{ \log\Big(R^{-\frac{N+1}{2}}(\log R)^{-1/2}\Big)} \cr
 &&= (\log R)^{-1} \left( \frac{N+1}{2}\log R +\half\log (\log R) \right) \cr
 &&\to \frac{N+1}{2} \quad \text{as}\ R\to \infty. 
% \vspace{-1em}
\end{eqnarray*}
whic shows the claim.
 \QED

\bigskip

Next we compute the interaction effect between $\chi(R^i,h_{R^i};\cdot)$ and $\chi(R^j,h_{R^j};\cdot)$ with $i, j \in \N$, $i\neq j$ and $R\gg 1$.

\begin{Lemma} \label{S:4.4}
For $i<j$ we have
 $$ \RRint I_\alpha(x-y)\chi(R^i,h_{R^i};x)\chi(R^j,h_{R^j};y)\, dxdy 
 \to 0 \quad \hbox{as} \ R\to \infty. $$
\end{Lemma}

%We postpone the proof of Lemma \ref{S:4.4} and we will give it in Section \ref{section:4.6}.

\claim Proof. % of Lemma \ref{S:4.4}.
Since $\supp\chi(R,h_R; \cdot)= A(R,h_R)$ %\big\{x\in \R^N \mid \abs x\in [R-h_R,R+h_R]\big\}$ 
we get
\begin{align*}
 \dist\big(\supp \chi(R^i,h_{R^i};\cdot), 
 \supp \chi(R^j,h_{R^j};\cdot)\big) 
 &= (R^j-h_{R^j})-(R^i+h_{R^i}) \cr
 &= R^j -O(R^i). 
\end{align*}
Thus
\begin{align*}
I_R:=& \; \RRint I_\alpha(x-y)\chi(R^i,h_{R^i};x)\chi(R^j,h_{R^j};y)\, dxdy \\
 \leq& \; C(R^j+{O(R^i)})^{-(N-\alpha)} \norm{\chi(R^i,h_{R^i};\cdot)}_1
 \norm{\chi(R^j,h_{R^j};\cdot)}_1. 
\end{align*}
Here
 $$ \norm{\chi(R,h_R;\cdot)}_1 = \meas(A(R,h_R)) \sim CR^{N-1}h_R
 $$
hence
 \begin{align*}
I_R 
 &\le q C(R^j-O(R^i))^{-(N-\alpha)} R^{(N-1)i} h_{R^i} R^{(N-1)j} h_{R^j} \\
 &\leq 
 C' R^{(\alpha-1)j+(N-1)i} h_{R^i}h_{R^j}.
 \end{align*}
When $\alpha\in (1,N)$, we have by the definition of $h_R$ 
\begin{align*}
I_R
 &\leq C R^{(\alpha-1)j+(N-1)i} R^{-\frac12(N-2+\alpha)(i+j)} \\
 &= C'R^{-\frac12 (N-\alpha)(j-i)} 
 \to 0 \quad \text{as}\ R\to\infty; 
\end{align*}
when $\alpha=1$, we obtain
\begin{align*}
I_R
 &\leq C'R^{(N-1)i}R^{-\frac12(N-1)(i+j)}(\log R^i)^{-\frac12}(\log R^j)^{-\frac12} \\
 &=C' R^{-\frac12(N-1)(j-i)} (ij)^{-\frac12} (\log R)^{-1} 
 \to 0 \quad \text{as}\ R\to\infty;
\end{align*}
when $\alpha\in (0,1)$,
 \begin{align*}
I_R
 &\leq C'R^{(\alpha-1)j+(N-1)i} R^{{-}\frac{N-1}{
 	{1}+\alpha}(i+j)} \\
 &= C' R^{-\frac1{1+\alpha}( (N-
 	{\alpha^2})j
 	{-}\alpha(N-1)i)}
 \to 0 \quad \text{as}\ R\to\infty.
 \end{align*}
 This concludes the proof.
 \QED

\bigskip

We have now the tools to build a refined path $\gamma \in \Gamma_n(\lambda)$.

\medskip

\claim Proof of Proposition \ref{S:4.1} (refined).
%For $s_0>0$ with $F(s_0)>0$, which is given in \hyperref[(F4)]{\textnormal{(F4)}} or \hyperref[(CF4)]{\textnormal{(CF4)}} and Remark \ref{R:2.2}, 
We construct now a path $\gamma\in \Gamma_n(\lambda)$; this path will moreover satisfy % such that
	$$	\max_{\xi\in D_n, \, x\in\R^N} |\gamma(\xi)(x)|\leq s_0.	$$
%Sembra che di questa relazione non parliamo più! Non la dimostriamo, non la usiamo; comunque è vera... 

\smallskip

\noindent 
\textbf{Step 1:} \emph{ Construction of an odd path in $L^r$.}
%
%\smallskip
\\
For $n\geq 2$ %(for $n=1$ the construction is simpler), 
we consider again the polyhedron
 $$ \Sigma=\big\{ t=(t_1,\dots,t_n) \mid \max_{i=1,\dots, n}\abs{t_i}=1\big\}.
 $$
%and we recall that $\Sigma$ is homeomorphic to $\partial D_n$. 
For a large $R\gg 1$, which we will choose later, we define 
 $$ \gamma_R(t)(x):=\sum_{i=1}^n \hbox{sgn}(t_i) \chi\big(R^i, |t_i| h_{R^i};x\big):\,
	\Sigma\to L^r(\R^N) $$
 where $r \in [1, +\infty]$. 
Here we regard $\chi(R^i,0;x)\equiv 0$, and we notice that $\gamma_R(t)$ is radial for each $t \in \Sigma$.
%For $s_0>0$ with $F(s_0)>0$, which is given in \hyperref[(F4)]{\textnormal{(F4)}} or \hyperref[(CF4)]{\textnormal{(CF4)}}, 
Considered $s_0$, we have
\begin{align*}
 \calD(s_0\gamma_R(t)) &= \sum_{i,j} 
 F(\hbox{sgn}(t_i)s_0)F(\hbox{sgn}(t_j)s_0) \cdot \cr
 &\quad\cdot \RRint I_\alpha(x-y)
 \chi(R^i,|t_i|h_{R^i}; x)
 \chi(R^i,|t_j|h_{R^i}; y)\, dxdy.
\end{align*}
We note that
\begin{itemize}
\item[(i)] For any $t=(t_1,\dots,t_n)\in\Sigma$, there exists at least one $t_k$ such that $\abs{t_k}=1$.
\item[(ii)] By Lemma \ref{S:4.3},
 $$ (F(\pm s_0))^2 \RRint I_\alpha(x-y)
 \chi(R^k,h_{R^k};x)\chi(R^k,h_{R^k};y)
 \,dxdy \geq C_0.$$
\item[(iii)] By (i) and (ii),
 $$ \sum_{i=1}^n (F(\pm s_0))^2 \RRint I_\alpha(x-y)
 \chi(R^i,h_{R^i};x)\chi(R^i,h_{R^i};y)
 \,dxdy \geq C_0. $$
\item[(iv)] If $i\not=j$, by Lemma \ref{S:4.4},
 $$ \RRint I_\alpha(x-y)\chi(R^i,h_{R^i};x)
 \chi(R^j,h_{R^j};x)\,dxdy \to 0
 \quad \hbox{as}\ R\to\infty. $$
\end{itemize}
%
%\smallskip
%
%\noindent
By (i)--(iv), we have for sufficiently large $R\gg 1$,
\begin{equation}
 \calD(s_0\gamma_R(t)) \geq C > 0 %0 
\quad \hbox{for all}\ t\in\Sigma. \label{eq:4.6}
\end{equation}
In what follows we fix $R\gg 1$ so that \eqref{eq:4.6} holds.

\smallskip

\noindent
\textbf{Step 2:} \emph{Construction of an odd path in $H_r^1$.}
%
%\smallskip
\\
For $0\leq h\ll R$ and $\epsilon>0$, we set
	$$	\chi_\epsilon(R,h;x):=
\begin{cases}
			1 &\text{if} \ x\in A(R,h), \\
			1-\frac 1\epsilon \dista(x, A(R,h)), &\text{if}\ \dista(x,A(R,h))\in (0,\epsilon),\\
			0 &\text{otherwise}.
 \end{cases}
	$$
Here we regard
	$	A(R,0)=\{ x\in\R^N|\, |x|=R\}.	$
We note that
	\begin{align*}
	&\chi_\epsilon(R,h;\cdot) \in H_r^1(\R^N) \ \text{for}\ \epsilon>0, \\
 	&\chi_\epsilon(R,h;\cdot)\to \chi(R,h;\cdot) \ \text{in}\ L^r(\R^N) \ \text{as}\ \epsilon\to 0 \ \text{for all}\ r\in [1,\infty), \\
	&\supp\chi_\epsilon(R^i, h_{R^i};\cdot) \cap \supp\chi_\epsilon(R^j, h_{R^j};\cdot) =\emptyset \ \text{for}\ i\not=j \ \text{for $\epsilon$ small}.
	\end{align*}
We set
\begin{equation}\label{eq_def_gamma_epsR}
\gamma_{\epsilon,R}(t):=\sum_{i=1}^n \sgn(t_i)\chi_\epsilon(R^i, |t_i| h_{R^i}; \cdot), \quad t \in \Sigma, %:\ \Sigma\to H_r^1(\R^N),
\end{equation}
%for $\epsilon>0$, 
$\gamma_{\epsilon,R}:\Sigma\to H_r^1(\R^N)$, continuous. By \eqref{eq:4.6} 
{and the continuity of $\mc{D}$ on $L^2(\R^N) \cap L^{2^*}(\R^N)$,} we have for $\epsilon>0$ small
%Bisogna scrivere più dettagli su questa parte!!
	$$	\mc{D}(s_0\gamma_{\epsilon,R}(t)) \geq C > 0 %0 
 \quad \text{for all}\ t\in\Sigma.	$$
Since
 $$ \mc{J}(\lambda, u(\cdot/\theta)) 
 = \half \theta^{N-2}\norm{\nabla u}_2^2 +\frac{e^{\lambda}}{2}
 \theta^N\norm u_2^2 -\half\theta^{N+\alpha}\calD(u), $$
we have for large $\theta\gg 1$
%bisogna giustificare perché $\theta$ non dipende da $t$! 
%poiché t (o direttamente u=gamma(t)) viaggia in un compatto, e la scelta di theta da t (o u) è continua?
 $$ \mc{J}(\lambda, s_0\gamma_{\epsilon,R}(t)(\cdot/\theta))<0 \quad
 \hbox{for all}\ t\in\Sigma \approx \partial D_n. $$
Considering $D_n=\{ st \mid s\in [0,1], \, t\in\Sigma\}$ and extending $s_0\gamma_{\eps,{R}}(t)(\cdot/\theta)$ to $D_n$
by
	$$	\widetilde \gamma(st):=s	{s_0}\gamma_{\epsilon,R}(t)(\cdot/\theta),	$$
finally we obtain a path $\widetilde\gamma\in \Gamma_n(\lambda)$. 
\QED

\medskip

\begin{Remark}
Even without assuming the positivity of $F$ (see Proposition \ref{prop_F_posit_intorn}), we notice that the construction of an odd map in $L^r$ gets easier when $F$ is an even function. 
Indeed there is no negative contribution given by the mixed interactions. We give only an outline of the proof, highlighting that in this case we do not need to use the fine Theorem \ref{S:4.2} given by \cite{Thim0}.

Define for every $i=1,\dots n$ and $s\in [-1,1]$ the annuli
$$A_i(s):= \big\{ x \in \R^N \mid |x| \in [2ni-|s|, 2ni+|s|]\big\}.$$
For every $t=(t_1, \dots, t_n) \in \Sigma$ we have that $A_1(t_1), \dots, A_n(t_n)$ are disjoint. Moreover, if $t_i=0$, then $\meas(A_i(t_i))=0$. 
Thus we define a continuous, odd map by
	$$	\gamma(t)(x):=\sum_{i=1}^n \sgn(t_i) \chi_{A_i(t_i)}(x): \, \Sigma\to L^2(\R^N)\cap L^{2^*}(\R^N). $$
Since $F$ is even, we obtain
\begin{align*}
 &\mc{D}(s_0\gamma(t)) \\
 &= \sum_{i,j} \iint_{A_i(t_i)\times A_j(t_j)} 
	I_{\alpha}(x-y) F(s_0\sgn(t_i)\chi_{A_i(t_i)}(x)) F(s_0\sgn(t_j)\chi_{A_j(t_j)}(y)) \,dxdy\\
 &= \big(F(s_0)\big)^2 \sum_{i,j} \iint_{A_i(t_i)\times A_j(t_j)} I_{\alpha}(x-y) \,dxdy \big(F(s_0)\big)^2\geq C>0, 
\end{align*}
where $C$ does not depend on the specific $t$. The regularization to a $H^1_r$-path can be done as in the general case (or by mollification), as well as the extension to $D_n$.

We highlight that this construction can be adapted also to the local case, and thus it gives a simplified construction of a multidimensional path in the setting of Berestycki and Lions \cite{BL2}.
\end{Remark}

%We start noticing that, by (CF4) and Remark \ref{R:2.2}, for some $\delta_0>0$ 
%	$$	F(s)>0 \quad \hbox{for}\ s\in (0,\delta_0],	$$
%which implies
%\begin{itemize}
%\item[(i)] when $F$ is even, $F(s)>0$ for all $s\in [-\delta_0,\delta_0]\setminus\{ 0\}$;
%\item[(ii)] when $F$ is odd, $F(s)<0$ for all $s\in [-\delta_0,0)$.
%\end{itemize}
%By (CF4), we also note that there exists $L_s>0$ with $L_s\to\infty$ as
%$s\to 0^+$ such that
%	\begin{equation}\label{eq:4.8}	
%		F(\sigma) \tr{\geq} %\leq
% L_s \sigma^p \quad \hbox{for all}\ \sigma\in [0,s].
%	\end{equation}

We are ready now to show that $\gamma_{R,\epsilon}:\, \Sigma\to H_r^1(\R^N)$, defined in \eqref{eq_def_gamma_epsR}, has the desired property \eqref{eq_stima_path_F}.

\begin{Lemma} \label{S:4.7}
Assume \hyperref[(F1)]{\textnormal{(F1)}}--\hyperref[(F5)]{\textnormal{(F5)}}, and $F>0$ in some $(0,\delta_0]$.
If $F$ is odd, additionally assume \eqref{eq_cond_more_gen}. 
%Assume $F$ odd or even, and $F>0$ in some $(0,\delta_0]$. 
Then there exists a constant $A>0$ independent of $s\in (0,\delta_0]$ and $t \in \Sigma$ such that
	$$	\mc{D}(s\gamma_{R,\epsilon}(t)) \geq \frac{1}{2}\big(F(s)\big)^2(A+o_{\eps}(1)) \quad
		\hbox{as}\ \epsilon\to 0;	$$
here $o_{\eps}(1)$ is a quantity which goes to $0$ as $\epsilon\to 0$ uniformly in $t\in\Sigma$ and $s\in (0,\delta_0]$.
\end{Lemma}

\medskip

\claim Proof. 
We prove Lemma \ref{S:4.7} in two steps.

\smallskip

\noindent
\textbf{Step 1:} For $t\in\Sigma$, set
	$$	a_{ij}(t):= \RRint I_\alpha(x-y)\chi(R^i,\abs{t_i}h_{R^i};x)
			\chi(R^j,\abs{t_j}h_{R^j};y)\, dxdy.
	$$
Then for sufficiently large $R>0$, we have
	\begin{equation}\label{eq:4.9}	
	A := \inf_{t\in\Sigma} \left(\sum_{i=1}^n a_{ii}(t)-\sum_{i\not=j}^n a_{ij}(t)\right)>0.
	\end{equation}
This fact follows from \eqref{eq:4.4} and Lemma \ref{S:4.4}. We fix $R\gg 1$ so that \eqref{eq:4.9} holds.

\smallskip

\noindent
\textbf{Step 2:} \emph{$\mc{D}(s\gamma_{R,\epsilon}(t)) \geq\half F(s)^2 A$ as $\epsilon\to 0$.} 
%
%\smallskip
%
%\noindent
\\
We note that for $\epsilon>0$ small
	$$	\supp\chi_\epsilon(R^i,\abs{t_i}h_{R^i}; \cdot) 
		\cap \supp\chi_\epsilon(R^j,\abs{t_j}h_{R^j};\cdot)
		=\emptyset \quad \hbox{for}\ i\not=j.	$$
Thus we have
	\begin{align}
	&\mc{D}(s\gamma_{R,\epsilon}(t)) \nonumber \\
	&= \sum_{i,j} \RRint I_\alpha(x-y)
		F(s\,\sgn(t_i)\chi_\epsilon(R^i,\abs{t_i}h_{R^i};x))
		F(s\,\sgn(t_j)\chi_\epsilon(R^j,\abs{t_j}h_{R^j};y)) %\, dxdy \nonumber 
\\
	&=: \sum_{i,j} B_{ij}(s,t).	\label{eq:4.10}
	\end{align}
We consider cases $i=j$ and $i\not=j$ separately.

First we focus on the case $i=j$. For both even and odd $F$ we have
	\begin{align}
	&B_{ii}(s,t) \nonumber \\
	&= \RRint I_\alpha(x-y)
 		F(s\,\sgn(t_i)\chi_\epsilon(R^i,\abs{t_i}h_{R^i};x))
		F(s\,\sgn(t_i)\chi_\epsilon(R^j,\abs{t_i}h_{R^i};y)) %\, dxdy 
\nonumber \\
	&= \RRint I_\alpha(x-y)
 		F(s\chi_\epsilon(R^i,\abs{t_i}h_{R^i};x))
		F(s\chi_\epsilon(R^j,\abs{t_i}h_{R^i};y)) %\, dxdy 
\nonumber \\
	&\geq \RRint I_\alpha(x-y)
 		F(s\chi(R^i,\abs{t_i}h_{R^i};x))
		F(s\chi(R^j,\abs{t_i}h_{R^i};y)) %\, dxdy 
\nonumber \\
	&=(F(s))^2 a_{ii}(t),		\label{eq:4.11}
	\end{align}
where we used the positivity of $F$ and the monotonicity of the integral. 
Next we consider the case $i\not=j$ for even $F$.
Since $F(s)\geq 0$ for $s\in [-\delta_0,\delta_0]$ we obtain
	\begin{equation}\label{eq:4.12}
	B_{ij}(s,t) \geq 0 \quad \hbox{for all}\ t\in \Sigma.
	\end{equation}
Finally we consider the case $i\not=j$ for odd $F$.
Since $\abs{F(s)}=F(\abs s)$ for $s\in [-\delta_0,\delta_0]$
	\begin{align}
	&B_{ij}(s,t) \nonumber \\
	&= \RRint I_\alpha(x-y)
 		F(s{\,\sgn(t_i)}\chi_\epsilon(R^i,\abs{t_i}h_{R^i};x))
		F(s{\,\sgn(t_j)}\chi_\epsilon(R^j,\abs{t_j}h_{R^j};y)) %\, dxdy \nonumber 
\\
	&\geq - \RRint I_\alpha(x-y)
 		F(s\chi_\epsilon(R^i,\abs{t_i}h_{R^i};x))
		F(s\chi_\epsilon(R^j,\abs{t_j}h_{R^j};y)) %\, dxdy
. \label{eq:4.13}
	\end{align}
Setting 
$$C_i(t,\epsilon):=\big\{x \mid \dista(x,A(R^i,\abs{t_i}h_{R^i}))\in (0,\epsilon)\big\}$$
we have
	\begin{align*}
	&\chi_\epsilon(R^i,\abs{t_i}h_{R^i};x) \in (0,1) \quad
		\hbox{for}\ x\in C_i(t_i,\epsilon), \\
	&\chi_\epsilon(R^i,\abs{t_i}h_{R^i};x) = \chi(R^i,\abs{t_i}h_{R^i};x) \quad
		\hbox{for}\ x\not\in C_i(t_i,\epsilon),\\
	&\meas(C_i(t_i,\epsilon)) \to 0 \quad \hbox{as}\ \epsilon\to 0, 
	 {\hbox{ uniformly in $t \in \Sigma$}}.
	\end{align*}
Thus for $r\in {[1,\infty)}$ {and $s \in (0,\delta]$}
%	\begin{align}
%	\left\|\frac{1}{F(s)}F(s\chi_\epsilon(R^i,\abs{t_i}h_{R^i};\cdot))
%			-\chi(R^i,\abs{t_i}h_{R^i};\cdot)\right\|_r^r \nonumber 
%	&\leq \int_{C_i(t_i,\epsilon)} 
%			\pabs{\frac{1}{F(s)}F(s\chi_\epsilon(R^i,\abs{t_i}h_{R^i};x)) 
%					}^r\,dx 
%			\nonumber \\
%	&= \left(\max_{h\in [0,1]} \frac{\abs{F(hs)}}{\abs{F(s)}} 
%	\right)^r \meas(C_i(t_i,\epsilon))
%			\nonumber \\
%	&\to 0 \quad \hbox{as}\ \epsilon\to 0 \ \hbox{uniformly in}\ t\in\Sigma. \label{eq:4.14}
%	\end{align}
	\begin{eqnarray*}
	\lefteqn{ \left\|\frac{1}{F(s)}F(s\chi_\epsilon(R^i,\abs{t_i}h_{R^i};\cdot))
			-\chi(R^i,\abs{t_i}h_{R^i};\cdot)\right\|_r^r } \nonumber \\
	&\leq& \int_{C_i(t_i,\epsilon)} 
			\pabs{\frac{1}{F(s)}F(s\chi_\epsilon(R^i,\abs{t_i}h_{R^i};x)) 
					}^r\,dx 
			\nonumber \\
	&=& \left(\max_{h\in [0,1]} \frac{\abs{F(hs)}}{\abs{F(s)}} 
	\right)^r \meas(C_i(t_i,\epsilon))
			\nonumber \\
	&\to 0& \quad \hbox{as}\ \epsilon\to 0 \ \hbox{uniformly in}\ t\in\Sigma. \label{eq:4.14}
	\end{eqnarray*}
Here we use the fact that $\max_{h\in [0,1]} \frac{F(sh)}{F(s)} \leq C$, which follows from the local almost-monotonicity assumption in \hyperref[(CF4)]{\textnormal{(CF4)}}.
We note that \eqref{eq:4.14} implies, exploiting again \hyperref[(CF4)]{\textnormal{(CF4)}} % the local monotonicity
	\begin{align}
	&\pabs{ \frac{1}{(F(s))^2}\RRint I_\alpha(x-y)
		F(s\chi_\epsilon(R^i,\abs{t_i}h_{R^i};x))
		F(s\chi_\epsilon(R^j,\abs{t_j}h_{R^j};y)) %\, dxdy 
-a_{ij}(t)} \nonumber \\
	&\to 0 	\quad \hbox{as}\ \epsilon\to 0. \label{eq:4.15}
	\end{align}
By \eqref{eq:4.13} and \eqref{eq:4.15},
	\begin{equation}\label{eq:4.16}
		B_{ij}(s,t) \geq -(F(s))^2(a_{ij}(t)+o(1)) \quad \hbox{as}\ \epsilon\to 0.
	\end{equation}
Thus, it follows from \eqref{eq:4.10}--\eqref{eq:4.12} and \eqref{eq:4.16} that
\begin{align*}
	\mc{D}(s\gamma_{R,\epsilon}(t)) 
	&\geq (F(s))^2\left(\sum_{i=1}^na_{ii}(t)-\sum_{i\not=j}a_{ij}+o(1)\right) \\
	&\geq \frac 12 (F(s))^2 (A + o(1))
%A \quad \hbox{for}\ \epsilon>0 \ \hbox{small}.
\end{align*}
 This concludes the proof.
\QED

%%%%%%%%%%%%%%%%%%%%%%%%%%%%%%%%%%

\section{Asymptotic analysis of %symmetric 
mountain pass values}
%\label{section:4..2}
\label{sec_choquard_asympt_MP}

We end this Section with some key estimates on the asymptotic behaviour of $a_n(\lambda)$ as $\lambda\to\pm\infty$.

\begin{Proposition}\label{S:4.6}
Assume \hyperref[(F1)]{\textnormal{(F1)}}--\hyperref[(F4)]{\textnormal{(F4)}} and let $n \in \N^*$.
\begin{itemize}
\item[(i)] If \hyperref[(CF3)]{\textnormal{(CF3)}} holds, then $\lim_{\lambda\to +\infty} \frac{a_n(\lambda)}{e^{\lambda}}=+\infty$.
\item[(ii)] If \hyperref[(CF4)]{\textnormal{(CF4)}} holds, then $\lim_{\lambda \to -\infty} \frac{a_n(\lambda)}{e^{\lambda}}=0$.
\end{itemize}
\end{Proposition}

\claim Proof of (i) of Proposition \ref{S:4.6}.
Recall $q={\frac{N+\alpha}{N}}$, $p={\frac{N+\alpha+2}{N}}$, and write $\mu=e^{\lambda}$ (and consequently adapt the notations) for the sake of simplicity. 

Since $a_n(\mu)\geq a_1(\mu)$ for each $n \in \N^*$, it is sufficient to show the claim for $n=1$. By \hyperref[(CF3)]{\textnormal{(CF3)}}, for any $\delta>0$ there exists
$C_\delta>0$ such that
$$ \abs{F(s)} \leq \delta \abs s^p +C_\delta \abs s^q \quad \hbox{for all}\ s\in\R.$$
%Basta $F$. Serve sottocritico!
For $v\in H^1_r(\R^N)$, setting $u_s:=s^{N/2}v(s \cdot)$, we have
\begin{align}\label{eq:4.7}
	 \mc{D}(u_s) &= s^{-N-\alpha} \mc{D}(s^{N/2}v) \nonumber \\
	&\leq s^{-N-\alpha} 
	\intRN \left(I_\alpha*(\delta s^{{\frac{N}{2}}p}\abs v^p + C_\delta s^{{\frac{N}{2}}q}\abs v^q)\right)
	(\delta s^{{\frac{N}{2}}p}\abs v^p + C_\delta s^{{\frac{N}{2}}q}\abs v^q) \, dx\nonumber \\
	&= s^2 \intRN \Big(I_\alpha*(\delta \abs v^p + C_\delta s^{-1}\abs v^q)\Big) 
	(\delta \abs v^p + C_\delta s^{-1}\abs v^q) \, dx\nonumber\\
	&=: s^2 D_{\delta,C_\delta s^{-1}} (v),
\end{align}
where we write for $\delta>0$ and $A\geq 0$,
	\begin{align*}
	D_{\delta,A}(v) &:= \intRN \Big(I_\alpha*(\delta \abs v^p + A\abs v^q)\Big)
	(\delta \abs v^p + A \abs v^q)\, dx, \\ 
	\mc{J}_{\delta,A}(v) &:= \half\norm{\nabla v}_2^2 +\half\norm v_2^2 -\half D_{\delta,A}(v).
	\end{align*}
We also denote by $b(\delta,A)$ the MP value of $\mc{J}_{\delta,A}$. Taking into account the continuity and monotonicity property of $b(\delta,A)$ with respect of each variable $\delta$ and $A$ and observing that $\mc{J}_{\delta,A}$ satisfies the (PS) condition, we have
$$ b(\delta,A) \to b(\delta,0) \quad \hbox{as}\ A\to 0^+,$$
$$b(\delta,0) \to +\infty \quad \hbox{as}\ \delta\to 0^+.$$
Thus, from \eqref{eq:4.7} we have that
$$ \mc{J}(\mu,u_s) \geq s^2\left( \half\norm{\nabla v}_2^2 + \half \mu s^{-2} \norm v_2^2 
-\half D_{\delta,C_\delta s^{-1}}(v)\right).$$
Setting $s:=\sqrt\mu$, we obtain
$$\mc{J}(\mu, u_{\sqrt\mu}) \geq \mu \mc{J}_{\delta, C_\delta \mu^{-1/2}}(v)$$
and thus $ {\frac{a_1(\mu)}{\mu}} \geq b(\delta, C_\delta \mu^{-1/2})$, which implies
$$ \liminf_{\mu\to\infty} {\frac{a_1(\mu)}{\mu}} \geq \lim_{A\to 0} b(\delta,A) = b(\delta,0).$$
Since $\delta>0$ is arbitrary, we gain
$$ \lim_{\mu\to+\infty} {\frac{a_1(\mu)}{\mu}}=+\infty.
\QED
$$

\bigskip

\noindent
We deal now with the proof of (ii) of Proposition \ref{S:4.6}. We highlight that, when $F$ is even, the proof can be simplified (see \cite{CT1}, Proposition \ref{out}, and Proposition \ref{prop_F_posit_intorn}).

The proof will be based on the key Lemma \ref{S:4.7}. 
We start noticing that, by \hyperref[(CF4)]{\textnormal{(CF4)}} and Remark \ref{R:2.2}, for some $\delta_0>0$ 
	$$	F(s)>0 \quad \hbox{for}\ s\in (0,\delta_0],	$$
which implies
\begin{itemize}
\item[(i)] when $F$ is even, $F(s)>0$ for all $s\in [-\delta_0,\delta_0]\setminus\{ 0\}$;
\item[(ii)] when $F$ is odd, $F(s)<0$ for all $s\in [-\delta_0,0)$.
\end{itemize}
By \hyperref[(CF4)]{\textnormal{(CF4)}}, we also note that there exists $L_s>0$ with $L_s\to +\infty$ as
$s\to 0^+$ such that
	\begin{equation}\label{eq:4.8}	
		F(\sigma) \geq %\leq
 L_s \sigma^p \quad \hbox{for all}\ \sigma\in [0,s].
	\end{equation}

\claim Proof of (ii) of Proposition \ref{S:4.6}.
Let $\gamma_{R,\epsilon}$ defined in \eqref{eq_def_gamma_epsR}. For $s_0\in (0,\delta_0]$ and $\mu>0$, we consider the map
	$$	st \in D_n \mapsto ss_0 \gamma_{R,\epsilon}(t)(\cdot/\mu^{-\frac 12}) \in H_r^1(\R^N) .	$$
We have by Lemma \ref{S:4.7} (since $\eps>0$ is here fixed small, we write $A$ instead of $A+o_{\eps}(1)$)
	\begin{align*}
	&\mu^{-1} \mc{J}(\mu,ss_0\gamma_{R,\epsilon}(t)(\cdot/\mu^{-\frac 12})) \\
	&=\frac 12 \mu^{-\frac N2}(ss_0)^2 \norm{\nabla\gamma_{R,\epsilon}(t)}_2^2
		+\frac 12\mu^{-\frac N2 
		}(ss_0)^2 \norm{\gamma_{R,\epsilon}(t)}_2^2
		-\frac 12\mu^{-\frac N2 p}\mc{D}(ss_0\gamma_{R,\epsilon}(t)) \\
	&\leq \frac 12\mu^{-\frac N2}(ss_0)^2 \norm{\gamma_{R,\epsilon}(t)}_{H^1}^2
		-\frac 14 \mu^{-\frac N2 p} (F(ss_0))^2 A.
	\end{align*}
Thus for $\mu$ small
	$$	\mc{J}(\mu,s_0\gamma_{R,\epsilon}(t)(\cdot/\mu^{-\frac 12})) <0 
		\quad \hbox{for}\ t\in\Sigma,	$$
which implies that $st \mapsto s_0\gamma_{R,\epsilon}(t)(\cdot/\mu^{-\frac 12})$ is a path belonging to $\Gamma_n(\mu)$.
Moreover by \eqref{eq:4.8}
	\begin{align*}
	\mu^{-1}a_n(\mu) 
	&\leq \max_{s\in [0,1], \, t\in\Sigma} \mu^{-1}\mc{J}(\mu,ss_0\gamma_{R,\epsilon}(t)(\cdot/\mu^{-\frac 12}))
	 \\
	&\leq \max_{s\in [0,1], \, t\in\Sigma} 
		\frac 12\mu^{-\frac N2}(ss_0)^2 \norm{\gamma_{R,\epsilon}(t)}_{H^1}^2
		-\frac 14 \mu^{-\frac N2 p} (F(ss_0))^2 A \\
	&\leq \max_{s\in [0,1], \, t\in\Sigma} 
		\frac 12\mu^{-\frac N2}(ss_0)^2 \norm{\gamma_{R,\epsilon}(t)}_{H^1}^2
		-\frac 14 L_{s_0} (\mu^{-\frac N2} (ss_0)^2)^p A \\
	&\leq C_{s_0},
	\end{align*}
where
	$$	C_{s_0}:=\sup_{\tau\geq 0, \, t\in\Sigma} 
			\left( \frac 12{\tau} \norm{\gamma_{R,\epsilon}(t)}_{H^1}^2
			-\frac 14 L_{s_0} A {\tau}^{p}\right) \in \R.	$$
Thus we have
	$$	\limsup_{\mu\to 0^+} \mu^{-1} a_n(\mu) \leq C_{s_0}.	$$
Since $C_{s_0}\to 0$ as $s_0\to 0$, we have (ii) of Proposition \ref{S:4.6}.
\QED

%%%%%%%%%%%%%%%%%%%%%%%%%%%%%%%%%%

\section{The Pohozaev mountain}
\label{section:4.3}

In this Section we start studying the Lagrangian formulation, applying the previous asymptotic estimates to a Pohozaev geometry.
%\tor{ aggiusta}
We consider the functional 
$\mc{I}^m : \R \times H^1_r(\R^N) \to \R$ defined by
\begin{equation}\label{eq:2.2}
\mc{I}^m(\lambda, u):=\half \norm{\nabla u}_2^2 -\half \mc{D}(u)+ \frac{e^{\lambda}}{2} \bigl( \|u\|_2^2 -m \bigr), \quad (\lambda, u) \in \R\times H^1_r(\R^N).
\end{equation}
It is immediate that, for any $(\lambda, u) \in \R\times H^1_r(\R^N)$,
$$ \mc{I}^m(\lambda, u) = \mc{J}(\lambda, u) - \frac{e^{\lambda}}{2} m. $$
If \hyperref[(F1)]{\textnormal{(F1)}}-\hyperref[(F2)]{\textnormal{(F2)}} hold, by \cite[Theorems 2 and 3]{MS2} we have that each solution $u$ of \eqref{eq:2.4} belongs to $W^{2,2}_{loc}(\R^N)$ and it satisfies the Pohozaev identity
\begin{equation}\label{eq:2.5}
\frac{N-2}{2} \|\nabla u\|_2^2 + \frac{N }{ 2} e^{\lambda} \|u\|_2^2 - \frac{N + \alpha}{2} \, \mc{D}(u) =0
\end{equation}
or equivalently
$$\frac{1}{2^*_{\alpha}} \|\nabla u\|_2^2 + \frac{e^{\lambda}}{2^{\#}_{\alpha}} \|u\|_2^2 - \mc{D}(u)=0$$
where $2^*_{\alpha}=\frac{N+\alpha}{N-2}$ and $2^{\#}_{\alpha} = \frac{N+\alpha}{N}$ are the upper and the lower critical exponents; again we see that essentially the identity means $\frac{d}{d\theta} \mc{J}(\lambda, u(\cdot/e^{\theta}))_{|\theta=0}=0$.
Inspired by this relation, we also introduce the Pohozaev functional $\mc{P}:\R \times H^1_r(\R^N) \to \R$ by setting 
\begin{equation}\label{eq:2.6}
\mc{P}(\lambda , u):= \frac{N-2}{2}\|\nabla u\|_2^2 -\frac{N+ \alpha}{2} \mc{D}(u) + \frac{N}{2} e^{\lambda} \|u\|_2^2 , \quad (\lambda, u) \in \R\times H^1_r(\R^N).
\end{equation}
We consider the action of $\G:=\Z_2$ on $\R^n$, $n \in \N^*$, and on $\R\times H^1_r(\R^N)$, given by
$$(\pm1,\xi) \in \G\times \R^n\mapsto \pm \xi \in \R^n,$$
$$ (\pm1,\lambda, u) \in \G\times \big(\R\times H^1_r(\R^N)\big) \mapsto (\lambda, \pm u)\in \R\times H^1_r(\R^N).$$
We notice that, under the assumption \hyperref[(F5)]{\textnormal{(F5)}}, $\mc{I}^m$, $\mc{J}$ and $\mc{P}$ are invariant under this action, i.e. they are even in $u$:
$$\mc{I}^m(\lambda, -u)= \mc{I}^m(\lambda, u), \quad \mc{J}(\lambda, -u)=\mc{J}(\lambda, u), \quad \mc{P}(\lambda, -u)= \mc{P}(\lambda, u).$$
In addition, we observe by the Principle of Symmetric Criticality of Palais \cite{Pal0} that every critical point of $\mc{I}^m$ %(resp. $\mc{J}$)
 restricted to $\R \times H^1_r(\R^N)$ is actually a critical point of $\mc{I}^m$ %(resp. $\mc{J}$) 
on the whole $\R \times H^1(\R^N)$.
% This observation justifies our restriction onto the radial setting.
Finally, we denote by $P_2: \R \times H^1_r(\R^N)\to H^1_r(\R^N)$ the projection on the second component.

\smallskip

Moreover we consider the \emph{Pohozaev %level 
set}
$$ \Omega :=\big\{(\lambda,u) \in \R \times H^1_r(\R^N) \mid \mc{P}(\lambda,u)>0\big\} \cup\big\{(\lambda,0) \mid \lambda \in \R \big\};$$
%We notice that, 
under the assumption \hyperref[(F5)]{\textnormal{(F5)}}, $\Omega$ is symmetric with respect to the axis $\{(\lambda,0) \mid \lambda \in \R \}$, that is,
$$(\lambda, u) \in \Omega \implies (\lambda, -u) \in \Omega.$$
We start showing the following property, due to the fact that $\mc{D}(u)=o(\norm u_{H^1}^2)$ as $u \to 0$.

\begin{Lemma}
We have
	\begin{equation}\label{eq:4.17}	
\{ (\lambda,0) \mid \lambda\in \R\} \subset int(\Omega).
	\end{equation}
\end{Lemma}

\claim Proof.
%Since $\mc{D}(u)=o(\norm u_{H\delta |s|^q + ^1}^2)$ as $u \to 0$, the conclusion follows from the definition of $\mc{P}(\lambda,u)$. 
%We observe that by $\textnormal{(CF3)}$ for $\delta>0$ fixed, there exists $C_\delta>0$ such that
%$$ |F(s)| \lesssim \delta |u|^q + C_{\delta} |s|^p $$
%where we recall $q=\frac{N+\alpha}{N}$ and $p=\frac{N+ \alpha+2}{N}$.}
%%Basta su $F$. Serve sottocritico.
% Thus 
%$$ \|F(u_n)\|_{\frac{2N}{N+ \alpha}} \leq \delta \| |u_n|^q\|_{\frac{2N}{N+ \alpha}} + C_{\delta} \| |u_n|^p \|_{\frac{2N}{N+ \alpha}} 
%= \delta \| u_n\|_2^q + C_{\delta}\| u_n \|_{\frac{2Np}{N+ \alpha}}^p . $$
%Therefore %\tor{by $\textnormal{(CF2)}$}, 
%by Proposition \ref{prop_HLS} and Young's inequality we have %Qui non serve $L^2$-sottocritico per $f$, basta per $F$ ($f$ potrebbe essere $L^2$-critico anche)
%\begin{eqnarray*}
%\lefteqn{ \int_{\R^N} (I_\alpha*\abs{F(u_n)})\abs{F(u_n)}\,dx 
%\leq 	 C\norm{F(u_n)}_{\frac{2N}{N+\alpha}}^2 
%\leq C\left( \delta \| u_n\|_2^q + C_{\delta}\| u_n \|_{\frac{2Np}{N+ \alpha}}^p\right)^2}
%	 \\ % \norm{F(u_n)}_{\frac{2N}{N+\alpha}} } \\
%	%\cdot \left(\norm{u_n}_{\frac{2Np}{N+\alpha}}^{p} 
%	%+ \norm{u_n}_2^{q}\right) } \\
%	&&\leq CC'\delta \norm{u_n}_{\frac{2Np}{N+\alpha}}^{2p}
%	+ CC'(\delta+C_\delta)\left({\frac{\delta}{2}}
%	\norm{u_n}_{\frac{2Np}{N+\alpha}}^{2p}
%	+{\frac{1}{2\delta}}\norm{u_n}_2^{2q}\right)
%	+ CC'C_\delta \norm{u_n}_2^{2q} \\
%	&&\leq C''\delta\norm{u_n}_{\frac{2Np}{N+\alpha}}^{2p}
%	+C''_\delta\norm{u_n}_2^{2q} 
%\end{eqnarray*}
By
$$ |F(s)| \lesssim |s|^q + |s|^p $$
where $q=\frac{N+\alpha}{N}$ and $p=\frac{N+ \alpha+2}{N}<2^*$.
 Thus 
$$ \|F(u)\|_{\frac{2N}{N+ \alpha}} \lesssim \| |u|^q\|_{\frac{2N}{N+ \alpha}} + \| |u|^p \|_{\frac{2N}{N+ \alpha}} 
= \| u\|_2^q + \| u \|_{\frac{2Np}{N+ \alpha}}^p . $$
Therefore by Proposition \ref{prop_HLS} and Young's inequality we have 
\begin{align*}
 \int_{\R^N} (I_\alpha*\abs{F(u)})\abs{F(u)}\,dx 
&\lesssim 	 \norm{F(u)}_{\frac{2N}{N+\alpha}}^2 
\lesssim \left( \norm{u}_2^q +\norm{u}_{\frac{2Np}{N+ \alpha}}^p\right)^2
	 \\ 
	&\lesssim \norm{u}_2^{2q} +\norm{u}_{\frac{2Np}{N+ \alpha}}^{2p} \leq \norm{u}_{H^1}^{2q} + \norm{u}_{H^1}^{2p}
\end{align*}
thus
$$ \mc{P}(\lambda,u) \lesssim \norm{u}_{H^1}^2 - \norm{u}_{H^1}^{2q} -\norm{u}_{H^1}^{2p} > 0$$
for $\norm{u}_{H^1}$ small, $u\neq 0$.
\QED

\bigskip

By \eqref{eq:4.17} we detect the \emph{Pohozaev mountain}
$$
\partial \Omega =\big\{(\lambda,u) \in \R \times H^1_r(\R^N) \mid \mc{P}(\lambda,u)=0, \ u \nequiv 0 \big\}.
$$
We observe that $\partial \Omega \neq \emptyset$, for instance by \cite[Theorems 1 and 3]{MS2}.

\begin{Proposition}\label{S:4.9}
Assume \hyperref[(F1)]{\textnormal{(F1)}}--\hyperref[(F4)]{\textnormal{(F4)}} and \hyperref[(F5)]{\textnormal{(F5)}}. We have the following properties.
\begin{itemize}
\item[(i)] $\mc{J}(\lambda,u)\geq 0$ for all $(\lambda,u)\in \Omega$.
\smallskip
\item[(ii)] $\mc{J}(\lambda,u)\geq a_1(\lambda)>0$ for all $(\lambda,u)\in \partial\Omega$. 
\item[(iii)] 
Assume \hyperref[(CF3)]{\textnormal{(CF3)}}. For any $m>0$, we set
$$E^m:= \inf_{(\lambda,u)\in\partial\Omega}\mc{I}^m(\lambda,u), \quad \text{ and } \quad B^m:=\inf_{\lambda \in \R} \left(a_1(\lambda) -\frac{e^\lambda}{2}m \right).$$
Then $E^m \geq B^m >-\infty$. 
In particular $B^m \in \R$ and
$$ \mc{I}^m(\lambda,u) \geq B^m \quad \hbox{for every}\ (\lambda,u)\in\partial\Omega. $$
\end{itemize}
\end{Proposition}

\claim Proof.
We notice that for all $(\lambda,u)\in \Omega$ 
$$ {\mathcal J}(\lambda,u) \geq {\mathcal J}(\lambda,u) - \frac{{\mathcal P}(\lambda,u)}{N + \alpha} 
= \frac{\alpha +2}{2(N + \alpha)} \norm{\nabla u}_2^2+ \frac{\alpha}{2(N + \alpha)} e^{\lambda}\norm{u}_2^2 \geq 0 $$
and thus $(i)$ follows. 
Point $(ii)$ follows from the fact that for each $\lambda$ the mountain pass level $a_1(\lambda)$ coincides with the ground state energy level (see \cite[Section 4.2]{MS1} and Section \ref{sec_proper_pohozaev_gs} for details);
see also Remark \ref{rem_altern_proof_>0}. %\ref{...} for an approach 
%for a proof which does not require an existence theorem for the unconstrained problem. %\tor{scrivi qui quei dettagli}).
Focus on $(iii)$: the fact that $E^m \geq B^m$ is a direct consequence of $(ii)$, while the fact that $B^m >-\infty$ comes from Proposition \ref{S:4.6} (i).
\QED

\medskip

\begin{Remark}\label{rem_altern_proof_>0}
In order to show that $a_1(\lambda)>0$, without exploiting the existence result for the unconstrained problem, we argue as follows (see also \cite{CGT2}). 
Let $\gamma \in \Gamma_1(\lambda)$; by definition of $\Gamma_1(\lambda)$ and by Proposition \ref{S:4.9} (i) there exists $t^*$ such that $\gamma(t^*) \in \partial \Omega$ and $\gamma(t^*) \neq 0$, thus $\mc{P}(\lambda, \gamma(t^*))=0$. This means that
$$\mc{J}(\lambda, \gamma(t^*)) =\frac{\alpha +2}{2(N + \alpha)} \norm{\nabla \gamma(t^*)}_2^2+ \frac{\alpha \mu}{2(N + \alpha)} \norm{\gamma(t^*)}_2^2 \simeq \norm{\gamma(t^*)}_{H^1}^2 $$
thus
$$a(\lambda) \gtrsim \inf_{u \in (\partial \Omega)_{\lambda}} \norm{u}_{H^1}^2.$$
Since, by \eqref{eq:4.17}, $ (\partial \Omega)_{\lambda}$ is far from the line $(\lambda,0)$, we obtain that the right-hand side is strictly positive, which is the claim.
\end{Remark}

 From now on we assume \hyperref[(CF3)]{\textnormal{(CF3)}} to give sense to the quantity $B^m$. In view of Proposition \ref{S:4.9} (iii), 
we set for $m>0$ and $n\in\N^*$
\begin{align*}
\Gamma_n^m:=\big\{\Theta \in C(D_n, \R\times H^1_r(\R^N)) \mid \; & \text{$\Theta$ is $\G$-equivariant,} \ \mc{I}^m(\Theta(0)) \leq B^m-1 ,\\
& \Theta|_{\partial D_n}\notin \Omega, \ \mc{I}^m(\Theta|_{\partial D_n})\leq B^m-1 \big\}
\end{align*}
and
$$b_n^m := \inf_{\Theta \in \Gamma_n^m} \sup_{\xi \in D_n} \mc{I}(\Theta(\xi));$$
we point out that asking $\Theta=(\Theta_1, \Theta_2) \in \Gamma^m_n$ to be $\G$-equivariant means that $\Theta_1$ is even and $\Theta_2$ is odd, and in particular $\Theta_2(0)=0$ which implies $\Theta(0) \in \Omega$.

\begin{Proposition}\label{S:4.10}
Assume \hyperref[(F1)]{\textnormal{(F1)}}-\hyperref[(F2)]{\textnormal{(F2)}}-\hyperref[(CF3)]{\textnormal{(CF3)}}-\hyperref[(F4)]{\textnormal{(F4)}}-\hyperref[(F5)]{\textnormal{(F5)}}. We have the following properties.
\begin{itemize}
\item[(i)] For any $m>0$ and $n \in \N^*$, we have $\Gamma_n^m\neq \emptyset$ and
\begin{equation}\label{eq:4.18}
b_n^m\leq a_n(\lambda) - e^{\lambda} \frac{m}{2},
\end{equation}
for each $\lambda \in \R$. Moreover, $b_n^m$ increases with respect to $n$.
\item[(ii)] For any $k\in\N^*$ there exists $m_k\geq 0$, namely given by
\begin{equation}\label{eq:4.19}
m_k:= 2\inf_{\lambda\in\R}\frac{a_k(\lambda)}{e^{\lambda}},
\end{equation}
such that for $m>m_k$
	$$	b_n^m < 0 \quad \text{for}\ n=1,2,\dots, k.	$$
Moreover, $m_k$ is increasing with respect to $k$.
\item[(iii)] If \hyperref[(CF4)]{\textnormal{(CF4)}} holds, then $m_k=0$ for each $k \in \N^*$. That is, for each $m>0$ we have
	$$	b_n^m < 0 \quad \text{for all}\ n\in\N^*.	$$
\end{itemize}
\end{Proposition}

\claim Proof. 
For given $\lambda\in\R$ and $\zeta\in \Gamma_n(\lambda)$, we will find a $\psi \in \Gamma_n^m$ such that
\begin{equation} \label{eq:4.20}
\max_{\xi \in D_n} \mc{J}(\psi(\xi)) \leq \max_{\xi \in D_n} \mc{J}(\lambda, \zeta(\xi)),
\end{equation}
so that we have
$$b_n^m \leq \max_{\xi \in D_n} \mc{I}^m(\psi(\xi)) \leq \max_{\xi \in D_n} \mc{J}(\lambda, \zeta(\xi))-{\frac{e^\lambda}{ 2}}m$$
and, passing to the infimum over $\Gamma_n(\lambda)$, we gain \eqref{eq:4.18}. 
\\
To find $\psi\in\Gamma_n^m$ with \eqref{eq:4.20}, observe that, by definition of $\Gamma_n(\lambda)$ and compactness of $\zeta(\partial D_n)$, there exists $C>0$ such that $\mc{D}(\zeta(\xi))\geq C >0$ for $\xi \in \partial D_n$.
Thus, we have $\mc{I}^m(\lambda, \zeta(\xi)(\cdot/L)) \to - \infty$ and $\mc{P}(\lambda, \zeta(\xi)(\cdot/L))\to -\infty$ as $L \to +\infty$, uniformly for $\xi \in \partial D_n$. 
Thus, for $L\gg 1$ we obtain, for every $\xi \in \partial D_n$,
\begin{equation} \label{eq:4.21}
\mc{I}^m(\lambda, \zeta(\xi)(\cdot/L)) \leq B^m-1 \quad \hbox{ and } \quad 
\mc{P}(\lambda, \zeta(\xi)(\cdot/L))<0.
\end{equation}
We also note that $\mc{I}^m(\lambda + L,0)=-{\frac{e^{\lambda + L}}{ 2}}m \to -\infty$
as $L\to+\infty$. 
Thus, for $L\gg 1$, we find that the path $\psi: D_n \to \R\times H^1_r(\R^N)$
$$\psi(\xi) := \parag{
&(\lambda+L(1-2|\xi|), \,0) & \quad \hbox{ if $|\xi| \in [0,1/2]$,} \\ 
&\left(\lambda, \, \zeta\left({\frac{\xi}{\abs{\xi}}}(2\abs\xi-1)\right)(\cdot/L)\right) & \quad \hbox{ if $|\xi| \in (1/2,1]$}
}
$$	
satisfies $\psi(0)=(\lambda+L, 0)\in \R \times \{0\}$, $\mc{I}^m(\psi(0))\leq B^m-1$ and $\mc{I}^m(\psi(\xi))\leq B^m-1$ for $\xi\in\partial D_n$. 
Thus, by \eqref{eq:4.21}, we obtain $\psi \in \Gamma^m_n$ and \eqref{eq:4.20}
holds.
\\
The monotonicity of $b_n^m$ with respect to $n$ is a consequence of the definition. 
Point $(ii)$ follows from \eqref{eq:4.18} and $(iii)$ follows from Proposition \ref{S:4.6} (ii). 
\QED

\bigskip

As a corollary to Proposition \ref{S:4.10}, we have the following result.

\begin{Corollary}
For any $m>0$, we have 
$$B^m= E^m = b^m_1,$$
i.e. the first minimax value $b_1^m$ equals the Pohozaev minimum $E^m$ on the product space.
\end{Corollary}

\medskip

\claim Proof.
Since any path in $\Gamma^m_n$ passes through $\partial \Omega$, we have $b^m_n \geq E^m \geq B^m$ for each $n$. 
On the other hand, passing to the infimum \eqref{eq:4.18} we obtain $b_1^m \leq B^m$ and thus the claim. 
\QED

%\bigskip
%
%By Propositions \ref{S:3.2} and \ref{S:3.6}, $\mc{I}^m$ satisfies the $(PSP)_b$ condition for $b<0$ and the deformation lemma holds. 
%Let $m_k\geq 0$ be a number given in Proposition \ref{S:4.10}. For $m>m_k$ we can see that $b_n^m<0$ for $n=1,2,\cdots, k$ are critical values of $\mc{I}^m$.
%If $b_n^m$ are different, we can see the multiplicity of solutions. To deal with the case $b_n^m=b_{n'}^m$ for some $n\not=n'$, we need another family of minimax methods, which we consider in the following Section.

%%%%%%%%%%%%%%%%%%%%%%%%%%%%%%%%%%%%%%%%%%%%%

%\setcounter{equation}{0} %FOR Arxiv
\section{The Palais-Smale-Pohozaev condition} 
%\label{section:3}
\label{sec_PSP_choquard}

For every $b \in \R$ we set
$$ K_b^m := \big\{ (\lambda,u)\in \R\times H^1_r(\R^N) \mid \mc{I}^m(\lambda, u)=b,\,
 \partial_\lambda\mc{I}^m(\lambda, u)=0,\, \partial_u\mc{I}^m(\lambda, u)=0 \big\}.$$
As already observed, under \hyperref[(F1)]{\textnormal{(F1)}}-\hyperref[(F2)]{\textnormal{(F2)}} we have that $\mc{P}(\lambda, u)=0$ for each $(\lambda, u) \in K_b^m$. 
We notice also that, assuming \hyperref[(F5)]{\textnormal{(F5)}}, $K_b^m$ is invariant under the $\G$-action, that is
$$(\lambda, u) \in K_b^m \implies (\lambda, -u) \in K_b^m.$$

Under our assumptions on $F$, it seems difficult to verify the standard Palais-Smale condition for the functional $\mc{I}^m$.
Therefore we cannot recognize that $K_b^m$ is compact.

Inspired by \cite{HT0,IT0,CT2}, we introduce the Palais-Smale-Pohozaev condition, a weaker compactness condition that takes into account the scaling properties of $\mc{I}^m$ through the Pohozaev functional $\mc{P}$. 
Through this tool we will show that $K_b^m$ is compact when $b<0$.

\begin{Definition}\label{S:3.1}
	For $b \in \R$, we say that $(\lambda_n, u_n)_n \subset \R \times H^1_r(\R^N)$ is a \emph{Palais-Smale-Pohozaev sequence} for $\mc{I}^m$ at level $b$ (shortly a $(PSP)_b$ sequence) if
		\begin{equation}\label{eq:3.1}
	\mc{I}^m (\lambda_n, u_n) \to b,
	\end{equation}
	\begin{equation}\label{eq:3.2}	
\partial_{\lambda} 	\mc{I}^m(\lambda_n, u_n) \to 0, 
\end{equation}
	\begin{equation}\label{eq:3.3}
\norm{\partial_u 	\mc{I}^m(\lambda_n, u_n)}_{(H^1_r(\R^N))^*} \to 0,
\end{equation}
	\begin{equation}\label{eq:3.4}
	\mc{P}(\lambda_n, u_n) \to 0.
	\end{equation}
	We say that $\mc{I}^m$ satisfies the \emph{Palais-Smale-Pohozaev condition} at level $b$ (shortly the $(PSP)_b$ condition) if every $(PSP)_b$ sequence has a strongly convergent subsequence in $\R \times H^1_r(\R^N)$.
\end{Definition}

We show now the following result.

\begin{Proposition}\label{S:3.2}
Assume \hyperref[(F1)]{\textnormal{(F1)}}-\hyperref[(CF2)]{\textnormal{(CF2)}}-\hyperref[(CF3)]{\textnormal{(CF3)}} and let $b <0$. Then $\mc{I}^m$ satisfies the $(PSP)_b$ condition.
\end{Proposition}

\claim Proof.
Let $b <0$ and let $(\lambda_n, u_n) \subset \R \times H^1_r(\R^N)$ be a $(PSP)_b$ sequence, i.e. satisfying \eqref{eq:3.1}--\eqref{eq:3.4}. 
First we note that by \eqref{eq:3.2} we obtain
\begin{equation}\label{eq:3.5}
	e^{\lambda_n}\big(\|u_n\|_2^2-m\big)\to 0.
\end{equation}

\noindent
\textbf{Step 1:} \emph{$\lambda_n$ is bounded from below and $\|u_n \|_2^2 \to m$ as $n \to + \infty$.}
%\smallskip
%
%\noindent
\\
We have by \eqref{eq:3.4}, \eqref{eq:3.1} and \eqref{eq:3.5}
\begin{align*}
o(1) &= \mc{P}(\lambda_n, u_n) \\ 
&=
- \frac{\alpha +2}{2} \|\nabla u_n\|_2^2 +(N+ \alpha) \Bigl(\mc{I}^m(\lambda_n, u_n) - \frac{e^{\lambda_n}}{2} \bigl(\|u_n\|_2^2 -m \bigr) \Bigr) +\frac{N}{2} e^{\lambda_n} \|u_n\|_2^2 \\ 
&=
- \frac{\alpha +2}{2} \|\nabla u_n\|_2^2 +(N+ \alpha) (b + o(1)) +\frac{N}{2} e^{\lambda_n}m+o(1).
\end{align*}
Here we used \eqref{eq:3.5}. 
From the above identity, we derive boundedness of $\lambda_n$ from below, since $b <0$. This result joined to \eqref{eq:3.5} finally gives $\|u_n\|_2^2\to m$.

\smallskip

\noindent
\textbf{Step 2:} \emph{$\lambda_n$ and $\|\nabla u_n \|_2^2$ are bounded.} 
%
%\smallskip
%
%\noindent
\\
Since, by \eqref{eq:3.3}, $\varepsilon_n := \| \partial_u \mc{I}^m(\lambda_n, u_n)\|_{(H_r^1(\R^N))^*}\to 0$, we have
\begin{equation}\label{eq:3.6}
\|\nabla u_n\|_2^2 - \int_{\R^N} (I_\alpha \ast F(u_n)) f(u_n) u_n dx 
 + e^{\lambda_n} \|u_n\|_2^2 \leq \varepsilon_n\|u_n\|_{H^1}.
\end{equation}
We observe that by \hyperref[(CF3)]{\textnormal{(CF3)}} for $\delta>0$ fixed, there exists $C_\delta>0$ such that
$$ |F(s)| \leq \delta |s|^p + C_\delta |s|^q $$
where we recall $p=\frac{N+ \alpha+2}{N}$ and $q=\frac{N+\alpha}{N}$.
%Basta su $F$. Serve sottocritico.
 Thus 
$$ \|F(u_n)\|_{\frac{2N}{N+ \alpha}} \leq \delta \| |u_n|^p \|_{\frac{2N}{N+ \alpha}} 
+ C_\delta \| |u_n|^q\|_{\frac{2N}{N+ \alpha}} 
= \delta \| u_n \|_{\frac{2Np}{N+ \alpha}}^p + C_\delta \| u_n\|_2^q. $$
Therefore by \hyperref[(CF2)]{\textnormal{(CF2)}}, Proposition \ref{prop_HLS} and Young's inequality we have %Qui non serve $L^2$-sottocritico per $f$, basta per $F$ ($f$ potrebbe essere $L^2$-critico anche)
\begin{eqnarray*}
\lefteqn{ \int_{\R^N} \big(I_\alpha*\abs{F(u_n)}\big)\abs{f(u_n)u_n}\,dx }
\\	&&\leq
	 C\norm{F(u_n)}_{\frac{2N}{N+\alpha}}\norm{f(u_n)u_n}_{\frac{2N}{N+\alpha}} \\
	&&\leq C\left(\delta\norm{u_n}_{\frac{2Np}{N+\alpha}}^{p}
	+ C_\delta\norm{u_n}_2^q\right)
	\cdot C'\left(\norm{u_n}_{\frac{2Np}{N+\alpha}}^{p} 
	+ \norm{u_n}_2^{q}\right) \\
	&&\leq CC'\delta \norm{u_n}_{\frac{2Np}{N+\alpha}}^{2p}
	+ CC'(\delta+C_\delta)\left({\frac{\delta}{2}}
	\norm{u_n}_{\frac{2Np}{N+\alpha}}^{2p}
	+{\frac{1}{2\delta}}\norm{u_n}_2^{2q}\right)
	+ CC'C_\delta \norm{u_n}_2^{2q} \\
	&&\leq C''\delta\norm{u_n}_{\frac{2Np}{N+\alpha}}^{2p}
	+C''_\delta\norm{u_n}_2^{2q} 
\end{eqnarray*}
and thus, by the Gagliardo-Nirenberg inequality and \eqref{eq:3.6},
\begin{align*}
\|\nabla u_n \|_2^2 + e^{\lambda_n} \|u_n\|_2^2 	
&\leq \int_{\R^N} (I_\alpha \ast |F(u_n)|) |f(u_n) u_n| dx + \varepsilon_n\|u_n\|_{H^1} \\	 
&\leq 
 C''' \delta \|\nabla u_n\|^2_2 \|u_n\|_2^{2(p-1)} 
	+ C_\delta'' \|u_n \|_2^{\frac{2(N + \alpha)}{N}} + \varepsilon_n\|u_n\|_{H^1}.
\end{align*}
Since by Step 1 $\|u_n\|_2^2=m+o(1)$, we obtain
\begin{eqnarray*} 
\lefteqn{(1- C''' \delta(m+o(1))^{p-1})\|\nabla u_n\|_2^2 + e^{\lambda_n}(m+o(1)) }\\
 &&\leq C''_\delta (m+o(1))^{\frac{N+\alpha}{N}} + \varepsilon_n(\|\nabla u_n\|_2^2+m+o(1))^{1/2}.
\end{eqnarray*}
For $\delta$ small enough, we have the boundedness of $e^{\lambda_n}$ and $\|\nabla u_n\|_2$.
Hence $\lambda_n$ can not go to $+\infty$ and thus by Step 1 we infer that $\lambda_n$ is bounded.

\smallskip

\noindent
\textbf{Step 3:} \emph{$\lambda_n$ and $ u_n $ strongly converge.}
%
%\smallskip
%
%\noindent
\\
By Steps 1-2, the sequence $(\lambda_n, u_n)_n$ is bounded in $\R \times H^1_r(\R^N)$ and thus after extracting a subsequence, denoted in the same way, we may assume that $\lambda_n \to \lambda_0$ and $u_n \rightharpoonup u_0$ weakly in $H^1_r(\R^N)$ for some $(\lambda_0, u_0) \in \R \times H^1_r(\R^N)$.
Taking into account the assumptions \hyperref[(F1)]{\textnormal{(F1)}}--\hyperref[(F3)]{\textnormal{(F3)}} and the compact embedding of $H^1_r(\R^N)$ in $L^r(\R^N)$ for $r \in (2, 2^*)$, we have
by Proposition \ref{prop_converg_generiche_nonloc}
$$ \int_{\R^N} (I_\alpha \ast F(u_n)) f(u_n) u_0 \, dx \to \int_{\R^N} (I_\alpha \ast F(u_0)) f(u_0) u_0 \, dx $$
and 
$$ \int_{\R^N} (I_\alpha \ast F(u_n)) f(u_n) u_n\, dx \to \int_{\R^N} (I_\alpha \ast F(u_0)) f(u_0) u_0\, dx. $$
%Serve su $f$ (CF2), che implica sottocriticità Sobolev.
By \eqref{eq:3.3} we derive that $\langle \partial_u \mc{I}^m(\lambda_n, u_n), u_n \rangle \to 0$ and $\langle \partial_u \mc{I}^m(\lambda_n, u_n), u_0 \rangle \to 0$, and hence $(\nabla u_n, \nabla (u_n-u_0))_2 + (u_n, e^{\lambda_n} u_n - e^{\lambda_0} u_0)_2 \to 0$. Combining this with $u_n \wto u_0$ and $\lambda_n \to \lambda_0$ we get
$$\|\nabla u_n\|_2^2 + e^{\lambda_n} \|u_n\|^2_2 \to 
\|\nabla u_0\|_2^2 + e^{\lambda_0} \|u_0\|^2_2 $$ 
which implies $u_n \to u_0$ strongly in $H^1_r(\R^N)$.
\QED

\bigskip

As a straightforward consequence we obtain the following result.
\begin{Corollary}\label{S:3.3}
	Assume \hyperref[(F1)]{\textnormal{(F1)}}-\hyperref[(CF2)]{\textnormal{(CF2)}}-\hyperref[(CF3)]{\textnormal{(CF3)}} and let $b <0$. Then $K_b^m \cap (\R \times \{0\}) = \emptyset$ and $K_b^m$ is compact.
\end{Corollary}

\begin{Remark}
	We emphasize that the $(PSP)_b$ condition does not hold at level $b=0$. Indeed we can consider a $(PSP)_0$ unbounded sequence $(\lambda_n, 0)$ with $\lambda_n \to - \infty$.
\end{Remark}

%%%%%%%%%%%%%%%%%%%%%%%%%%%%%%%%%%

\section{Genus-shaped critical points}
%\subsection{\tb{An augmented functional}}
\label{section:3.2}

In this Section we essentially follow the lines of Sections \ref{sec_frac_augmented}--\ref{sec_frac_multiple_sol}. We give just an outline, avoiding details and proofs.

%\tor{SFOLTISCI (o sposta qualcosa)}

\subsection{Augmented functional} %A functional in the augmented space}} %A scaling augmented functional}}

We start by achieving a deformation lemma. In order to do this we define
$$M := \R \times \R \times H^1_r(\R^N)$$
and introduce the augmented functional $\mc{H}^m: M\to \R$
%\begin{equation}\label{eq:3.7}
%\mc{H}^m(\theta, \lambda, u):=\mc{I}^m(\lambda, u(e^{-\theta}\cdot)), \quad (\theta, \lambda, u) \in M.
%\end{equation}
%By the scaling properties of $\mc{I}^m$ we can recognize that
%$$ \mc{H}^m(\theta, \lambda, u)= 
%\frac{e^{(N-2)\theta}}{2}
%\norm{\nabla u}_2^2 - \frac{e^{(N+\alpha) \theta}}{2}
%\mc{D}(u) + 
%\frac{e^\lambda}{2} \bigl( e^{N \theta}\|u\|_2^2 -m \bigr) $$
\begin{align}\label{eq:3.7}
\mc{H}^m(\theta, \lambda, u):=& \,\mc{I}^m(\lambda, u(e^{-\theta}\cdot)) \\
= & \,
\frac{e^{(N-2)\theta}}{2}
\norm{\nabla u}_2^2 - \frac{e^{(N+\alpha) \theta}}{2}
\mc{D}(u) + 
\frac{e^\lambda}{2} \bigl( e^{N \theta}\|u\|_2^2 -m \bigr) \notag
\end{align}
for all $(\theta, \lambda,u ) \in M$, and thus
$\partial_\theta \mc{H}^m(\theta, \lambda, u) = 
		\mc{P}(\lambda, u(\cdot /e^\theta)).		$ 
We point out that, considered the action of $\G$ on $M$
$$\G \times M \to M ; \, (\pm1,\theta, \lambda, u) \mapsto (\theta, \lambda, \pm u) $$
and assumed \hyperref[(F5)]{\textnormal{(F5)}}, it results that $\mc{H}^m$ is $\G$-invariant.
%, i.e. it is even in $u$:
%$$\mc{H}^m(\theta, \lambda, -u)= \mc{H}^m(\theta, \lambda, u).$$
Introducing a metric on $M$ by 
$$ {\|(\alpha, \nu, h)\|}_{(\theta,\lambda, u)}^2:=\left|\left(\alpha, \nu,\norm{h(e^{-\theta} \cdot)}_{H^1}\right)\right|^2 $$
for any $(\alpha, \nu, h) \in T_{(\theta,\lambda,u)} M \equiv \R \times \R \times H^1_r(\R^N)$,
we regard $M$ as a Hilbert manifold.
We also denote the dual norm on 
$T^*_{(\theta,\lambda,u)}M$ by $\|\cdot \|_{(\theta,\lambda, u), *}$, 
and observe that both $\norm{\cdot}_{(\theta, \lambda, u)}$ and $\norm{\cdot}_{(\theta, \lambda, u),*}$ actually depend only on $\theta$. 

Denote now
	$	D:=(\partial_\theta,\partial_\lambda,\partial_u)	$
the gradient with respect to all the variables; a direct computation shows that
%\begin{eqnarray*}
%\lefteqn{D\mc{H}^m(\theta, \lambda, u)(\alpha,\nu,h) }\\
%&&=\mc{P}(\lambda, u(e^{-\theta} \cdot))\alpha + \partial_{\lambda} \mc{I}^m(\lambda, u(e^{-\theta}\cdot))\nu +
% \langle \partial_u \mc{I}^m(\lambda, u(e^{-\theta}\cdot)),\, h(e^{-\theta} \cdot)\rangle
% \end{eqnarray*}
for any $(\theta, \lambda, u)\in M$ %and $(\alpha, \nu, h) \in T_{(\theta, \lambda, u)}M$, and thus we obtain
\begin{eqnarray*}
\lefteqn{ \|{D\mc{H}^m(\theta, \lambda, u)\|}_{(\theta, \lambda, u),*}^2}\\
&& =\abs{\mathcal{P}(\lambda, u(e^{-\theta} \cdot))}^2 + \abs{\partial_{\lambda}\mc{I}^m(\lambda, u(e^{-\theta}\cdot))}^2 + \norm{\partial_u \mc{I}^m(\lambda, u(e^{-\theta}\cdot))}_{(H^1_r(\R^N))^*}^2 .
\end{eqnarray*}
We furthermore define 
%$$\tilde{K}_b^m := \big\{ (\theta, \lambda, u) \in M \mid \mc{H}^m(\theta, \lambda, u)=b,\, D \mc{H}^m(\theta, \lambda,u)=0\big\}$$
%the set of critical points of $\mc{H}^m$ at level $b$, and we deduce
%\begin{equation}\label{eq:3.8}
%\tilde{K}_b^m = \big\{(\theta, \lambda, u(e^{\theta} \cdot)) \mid (\lambda, u)\in K_b^m, \; \theta \in \R\big\}.
%\end{equation}
the set of critical points of $\mc{H}^m$ at level $b$ by
\begin{align}
\tilde{K}_b^m :=& \big\{ (\theta, \lambda, u) \in M \mid \mc{H}^m(\theta, \lambda, u)=b,\, D \mc{H}^m(\theta, \lambda,u)=0\big\} \notag \\
=&
\big\{(\theta, \lambda, u(e^{\theta} \cdot)) \mid (\lambda, u)\in K_b^m, \; \theta \in \R\big\}. \notag %\label{eq:3.8}
\end{align}
Finally we introduce the % standard 
distance between two points as %the infimum of the lengths of the curves connecting the two points, namely
\begin{eqnarray*}
\lefteqn{\dista_M \big((\theta_0, \lambda_0, h_0), (\theta_1, \lambda_1, h_1)\big) :=} \\
&&\inf \left\{ \int_0^1 \|\dot \gamma(t)\|_{\gamma(t)} \, dt \mid \gamma \in C^1([0,1],M), \; \gamma(0)= (\theta_0, \lambda_0, h_0), \gamma(1)= (\theta_1, \lambda_1, h_1)\right\}.
\end{eqnarray*}

As a consequence of Proposition \ref{S:3.2} we obtain the following.

\begin{Proposition}\label{S:3.5}
	Assume \hyperref[(F1)]{\textnormal{(F1)}}-\hyperref[(CF2)]{\textnormal{(CF2)}}-\hyperref[(CF3)]{\textnormal{(CF3)}} and let $b <0$. Then $\mc{H}^m$ satisfies the following Palais-Smale-type condition $(\widetilde{PSP})_b$: 
%That is, 
for each sequence $(\theta_n, \lambda_n, u_n)_n\subset M$ such that
	$$\mc{H}^m(\theta_n, \lambda_n, u_n) \to b,$$
	$$\norm{D \mc{H}^m(\theta_n, \lambda_n, u_n)}_{(\theta_n, \lambda_n, u_n),*} \to 0$$
	as $n \to +\infty$, we have, up to a subsequence,
	$$\dista_M((\theta_n, \lambda_n, u_n), \tilde{K}_b^m)\to 0.$$
\end{Proposition}

%%%%%%%%%%%%%%%%%%%%%%%%%%%%%%%%%%

\subsection{Deformation theory}
\label{section:3.3}

We write, for $b \in \R$
$$ [\mc{I}^m\leq b] := \big\{ (\lambda,u)\in\R\times H_r^1(\R^N) \mid \mc{I}^m(\lambda,u) \leq b\big\},$$
$$[\mc{H}^m \leq b]_M := \big\{ (\theta,\lambda,u)\in M \mid \mc{H}^m(\theta,\lambda,u) \leq b\big\}. $$ 
We state the following result.

\begin{Proposition}\label{S:3.6}
Assume \hyperref[(F1)]{\textnormal{(F1)}}-\hyperref[(CF2)]{\textnormal{(CF2)}}-\hyperref[(CF3)]{\textnormal{(CF3)}}. Let $b<0$, and let $\mc{O}$ be a neighborhood of $K_b^m$ with respect to the standard distance of $\R \times H^1_r(\R^N)$. Let $\bar{\varepsilon}>0$, then there exist $\varepsilon \in (0,\bar{\varepsilon})$ and $\eta: [0,1]\times (\R\times H^1_r(\R^N))\to \R\times H^1_r(\R^N)$ continuous such that
\begin{enumerate}
\item $\eta(0, \cdot,\cdot)=id_{\R\times H^1_r(\R^N)}$;
\item $\eta$ fixes $[\mc{I}^m\leq b-\bar{\varepsilon}]$, that is, $\eta(t, \cdot,\cdot)=id_{[\mc{I}^m\leq b-\bar{\varepsilon}]}$ for all $t \in [0,1]$;
\item $\mc{I}^m$ is non-increasing along $\eta$, and in particular $\mc{I}^m(\eta(t,\cdot, \cdot))\leq \mc{I}^m(\cdot, \cdot)$ for all $t \in [0,1]$;
\item if $K_b^m= \emptyset$, then $\eta(1, [\mc{I}^m\leq b+\varepsilon])\subseteq [\mc{I}^m\leq b-\varepsilon]$;
\item if $K_b^m\not=\emptyset$, then
$$\eta(1,[\mc{I}^m\leq b+\varepsilon]\setminus \mc{O}) \subseteq [\mc{I}^m\leq b-\varepsilon]$$
and
$$\eta(1, [\mc{I}^m\leq b+\varepsilon]) \subseteq [\mc{I}^m\leq b-\varepsilon]\cup \mc{O};$$
\item if \hyperref[(F5)]{\textnormal{(F5)}} holds, then $\eta(t, \cdot,\cdot)$ is $\G$-equivariant, i.e. for $\eta=(\eta_1, \eta_2)$ we have $\eta_1$ even and $\eta_2$ odd in $u$.
\end{enumerate}
\end{Proposition}

To prove this, we work first on the functional $\mc{H}$, for which we obtained a $(\widetilde{PSP})$ condition, which implies that 
for any $b<0$ there exists $\epsilon$, $\delta$, $\nu>0$ such that
	$$	\| D\mc{H}^m(\theta,\lambda,u)\|_{(\theta,\lambda,u),*} \geq \nu
	$$
for $(\theta,\lambda,u)\in M$ satisfying $\mc{H}^m(\theta,\lambda,u)\in [b-\epsilon,b+\epsilon]$ and $\dista_M((\theta,\lambda,u), \widetilde K_b^m)
\geq \delta$.

\begin{Proposition}\label{S:3.7}
Assume \hyperref[(F1)]{\textnormal{(F1)}}-\hyperref[(CF2)]{\textnormal{(CF2)}}-\hyperref[(CF3)]{\textnormal{(CF3)}}. Let $b<0$, and let $\tilde{\mc{O}}$ be a neighborhood of $\tilde{K}_b^m$ with respect to $\dist_M$. 
Let $\bar{\varepsilon}>0$,
then there exist $\varepsilon \in (0,\bar{\varepsilon})$ and $\tilde{\eta}: [0,1]\times M\to M$ continuous such that
\begin{enumerate}
\item $\tilde{\eta}(0, \cdot,\cdot,\cdot)=id_M$;
\item $\tilde{\eta}$ fixes $[\mc{H}^m \leq b-\bar{\varepsilon}]_M$, that is $\tilde{\eta}(t, \cdot,\cdot,\cdot)=id_{[\mc{H}^m \leq b-\bar{\varepsilon}]_M}$ for all $t\in [0,1]$;
\item $\mc{H}^m$ is non-increasing along $\tilde{\eta}$, and in particular $\mc{H}^m(\tilde{\eta}(t,\cdot,\cdot, \cdot))\leq \mc{H}^m(\cdot, \cdot, \cdot)$ for all $t\in [0,1]$; 
\item if $\tilde{K}_b^m= \emptyset$, then $\tilde{\eta}(1, [\mc{H}^m \leq b+\varepsilon]_M)\subseteq [\mc{H}^m \leq b-\varepsilon]_M$;
\item if $\tilde{K}_b^m \neq \emptyset$, then
$$\tilde{\eta}(1,[\mc{H}^m \leq b+\varepsilon]_M\setminus \tilde{\mc{O}}) \subseteq [\mc{H}^m \leq b-\varepsilon]_M$$
and
$$\tilde{\eta}(1,[\mc{H}^m \leq b+\varepsilon]_M) \subseteq [\mc{H}^m \leq b-\varepsilon]_M\cup \tilde{\mc{O}};$$
\item if \hyperref[(F5)]{\textnormal{(F5)}} holds, then $\tilde{\eta}(t, \cdot,\cdot)$ is $\G$-equivariant, i.e. for $\tilde{\eta}=(\tilde{\eta}_0,\tilde{\eta}_1,\tilde{\eta}_2)$ we have $\tilde{\eta}_0$, $\tilde{\eta}_1$ even and $\tilde{\eta}_2$ odd in $u$.
\end{enumerate}
\end{Proposition}
To get Proposition \ref{S:3.6} from Proposition \ref{S:3.7} we introduce
$$\pi: (\theta, \lambda, u)\in M \mapsto (\lambda, u(e^{-\theta}\cdot)) \in \R\times H^1_r(\R^N), $$
$$\iota: (\lambda, u) \in \R\times H^1_r(\R^N) \mapsto (0,\lambda, u) \in M,$$
which are a kind of rescaling projection and immersion satisfying
%. Observe that
$$\pi \circ \iota = id_{\R\times H^1_r(\R^N)}, \quad \pi(\tilde{K}_b^m) = K_b^m, $$%\hbox{(while $\iota \circ \pi \neq id_M$),}$$
$$\mc{H}^m\circ \iota=\mc{I}^m, \quad \mc{I}^m\circ \pi=\mc{H}^m.$$
%$$\pi(\tilde{K}_b^m) = K_b^m.$$
%By means of these operators we are able to prove that neighborhoods of $\tilde{K}_b^m$ are brought to neighborhoods of $K_b^m$.
For a deformation $\tilde{\eta}$ obtained in Proposition \ref{S:3.7} we thus define 
\begin{equation}\label{eq:3.9}
\eta(t,\lambda, u):=\pi(\tilde{\eta}(t,\iota(\lambda, u))), \quad (t, \lambda, u) \in [0,1] \times (\R \times H^1_r(\R^N)).
\end{equation}

%%%%%%%%%%%%%%%%%%%%%%%%%%%%%%%%%%

\subsection{Multiple critical points}
%\label{section:4.4}

For each $n\in \N^*$, define
\begin{align*}
\Lambda_n^m:=\{A=\Theta(\overline{D_{n+l}\setminus Y}) \mid \; & l\in \N %^*
, \; \Theta\in \Gamma_{n+l}^m, \\
& Y\subseteq D_{n+l}\setminus \{0\} \; \hbox{ is closed, symmetric in $0$}\\
& \qquad \qquad \qquad \quad \; \hbox{ and $\genus(Y)\leq l$} \}
\end{align*}
and
$$c_n^m:= \inf_{A\in \Lambda_n^m} \sup_{A} \mc{I}^m.$$

We notice that $\left\{ \Theta(D_n)\right\}_{\Theta \in \Gamma^m_n} \subset \Lambda^m_n$. 
In the following lemma, we observe that $\Lambda_n^m$ and $c_n^m$ inherit the properties of $\Gamma_n^m$ and $b_n^m$, also given by
$$ A \cap \partial \Omega\not=\emptyset \quad \hbox{for all}\ A\in \Lambda_1^m,	$$
together with the %and they enjoy an 
extra property (v). 

\begin{Proposition} \label{S:4.12}
Assume \hyperref[(F1)]{\textnormal{(F1)}}-\hyperref[(F2)]{\textnormal{(F2)}}-\hyperref[(CF3)]{\textnormal{(CF3)}}-\hyperref[(F4)]{\textnormal{(F4)}}. Let $n \in \N^*$ and $m>0$. Then
\begin{itemize}
\item[(i)] $\Lambda_n^m \neq \emptyset$.
\item[(ii)] $\Lambda_{n+1}^m\subseteq \Lambda_n^m$, and thus $c_n^m\leq c_{n+1}^m$.
\item[(iii)] $c_n^m\leq b_n^m$.
\item[(iv)] $B^m = E^m \leq c_1^m$.
\item[(v)] Let $A\in \Lambda_n^m$ and $Z\subset \R \times H^1_r(\R^N)$ be $\G$-invariant, closed, and such that $0 \notin \overline{P_2(Z)}$ and $\genus(\overline{P_2(Z)})\leq i<n$. Then $\overline{A\setminus Z} \in \Lambda_{n-i}^m$.
\end{itemize}
\end{Proposition}

%%%%%%%%%%%%%%%%%%%%%%%%%%%%%%%%%%

Fix $n\in \N^*$ and let $\Lambda_n^m$ and $c_n^m$ satisfying the properties of Proposition \ref{S:4.12}. We build now multiple solutions.

\begin{Proposition}\label{S:4.13}
Assume \hyperref[(F1)]{\textnormal{(F1)}}-\hyperref[(CF2)]{\textnormal{(CF2)}}-\hyperref[(CF3)]{\textnormal{(CF3)}}-\hyperref[(F4)]{\textnormal{(F4)}}-\hyperref[(F5)]{\textnormal{(F5)}}. Fix $k\in \N^*$ and assume $m>m_k$ (see \eqref{eq:4.19}). Then
$$c_1^m \leq c_2^m \leq \dots \leq c_k^m<0$$
are critical values of $\mc{I}^m$. Moreover
\begin{itemize}
\item[(i)]
 if, for some $q\in \N^*$,
$$c_n^m < c_{n+1}^m < \dots < c_{n+q}^m<0$$
then we have $q+1$ different nonzero critical values, and thus $q+1$ different pairs of nontrivial solutions of \eqref{eq:1.7};
\item[(ii)]
if instead, for some $q\in \N^*$,
\begin{equation}\label{eq:4.22}
	c_n^m = c_{n+1}^m = \dots = c_{n+q}^m =:b <0
\end{equation}
then 
\begin{equation}\label{eq:4.23}
\genus(P_2(K_b^m))\geq q+1
\end{equation}
and thus $\# P_2(K_b^m)=+\infty$, which means that we have infinite different solutions of \eqref{eq:1.7}. 
\end{itemize}
Summing up, we have at least $k$ different pairs of nontrivial solutions of \eqref{eq:1.7}.
\end{Proposition}

%%%%%%%%%%%%%%%%%%%%%%%%%%%%%%%%%%%%%%%%%%%%%

%\setcounter{equation}{0} %FOR Arxiv
\section{The unconstrained problem}
%\label{section:5}
\label{sec_unconstrained_choquard}

%COMMENT NOW
%\tr{Il risultato di esistenza di ground state si può estendere al caso frazionario sottocritico (non strettamente)? Non riesco a trovare una referenza (Byeon è sottocritico stretto, Chang chiedere $C^1$ per usare Pohozaev, Ambrosio chiede anche lui regolarità, Ikoma chiede a massa nulla (?), ...)}

In this Section we show how to exploit some of the developed tools also % sketch how 
to obtain infinitely many radial solutions for the \emph{unconstrained problem} \eqref{eq_Choquard_genericaF}, and give a sketch of the proof of Theorem \ref{S:1.3}.
Here we assume \hyperref[(F1)]{\textnormal{(F1)}}--\hyperref[(F5)]{\textnormal{(F5)}}.
We fix $\lambda \in \R$ and write $\mu=e^\lambda$; omitting $\lambda$, we denote $\mc{J}(\cdot):=\mc{J}(\lambda, \cdot): H^1_r(\R^N) \to \R$, i.e.
\begin{equation}
\mc{J}(u):=\half \norm{\nabla u}_2^2 -\half \mc{D}(u) + \frac{\mu}{2} \|u\|_2^2, \quad u \in H^1_r(\R^N).
\end{equation}
Similarly we write $\mc{P}(\cdot):=\mc{P}(\lambda, \cdot)$. 
For every $b \in \R$ we set
$$ K_b := \{ u\in H^1_r(\R^N) \mid \mc{J}(u)=b,\,
 \mc{J}'(u)=0 \}.$$
We have the following result.
\begin{Proposition}\label{S:5.1}
Assume \hyperref[(F1)]{\textnormal{(F1)}}--\hyperref[(F3)]{\textnormal{(F3)}} and let $b \in \R$. Then $\mc{J}$ satisfies the Palais-Smale-Pohozaev condition at level $b$ (shortly $(PSP)_b$), that is every sequence $(u_n)_n \subset H^1_r(\R^N)$ satisfying
\begin{equation}\label{eq:5.26}
	\mc{J} (u_n) \to b,
	\end{equation}
	\begin{equation}\label{eq:5.27}
\norm{\mc{J}'(u_n)}_{(H^1_r(\R^N))^*} \to 0,
\end{equation}
	\begin{equation}\label{eq:5.28}
	\mc{P}(u_n) \to 0,
	\end{equation}
	admits a strongly convergent subsequence in $H^1_r(\R^N)$.
In particular, $K_b$ %(\lambda)$ 
is compact in $H_r^1(\R^N)$. 
\end{Proposition}

\claim Proof.
First observe that, by \eqref{eq:5.26} and \eqref{eq:5.28} we obtain
\begin{equation}\label{eq:5.29}
 \frac{\alpha+2}{2}\norm{\nabla u_n}_2^2 + \frac{\alpha}{2} \mu\norm{u_n}_2^2 = (N+\alpha)b + o(1).
\end{equation}
We observe that $b\geq 0$ and the boundedness of $u_n$ in $H_r^1(\R^N)$. Thus by \hyperref[(F2)]{\textnormal{(F2)}}-\hyperref[(F3)]{\textnormal{(F3)}}, $\mc{D}'(u_n)$ has a strongly convergent subsequence in $(H_r^1(\R^N))^*$ and by \eqref{eq:5.27}, $u_n$ has a strongly convergent subsequence in $H_r^1(\R^N)$. Here we make use of Proposition \ref{prop_converg_generiche_nonloc}. 
%\tr{Serve Sobolev sottocritico per $f$, non basta per $F$!}
%
%COMMENT NOW
%\\ \tr{Più dettagli? per evidenziare che vale anche per $f$ sottoscritica non strettamente, e così possiamo dare il risultato anche per il frazionario, estendendo (di poco) anche Byeon-Kwon-Seok e Chang-Wang}
\QED

\bigskip

Set $ [\mc{J}\leq b] := \{ u\in H_r^1(\R^N) \mid \mc{J}_{\lambda}(u) \leq b\}$. 
Following the arguments of Section \ref{section:3.2} and \ref{section:3.3}, we prove the following deformation result by means of an augmented functional.

\begin{Proposition}\label{S:5.2}
Assume \hyperref[(F1)]{\textnormal{(F1)}}--\hyperref[(F3)]{\textnormal{(F3)}}. Let $b\in \R$ and let $\mc{O}$ be a neighborhood of $K_b(\lambda)$. Let $\bar{\varepsilon}>0$, then there exist $\varepsilon \in (0,\bar{\varepsilon})$ and $\eta: [0,1]\times H^1_r(\R^N)\to H^1_r(\R^N)$ continuous such that
\begin{enumerate}
\item $\eta(0, \cdot)=id_{H^1_r(\R^N)}$;
\item $\eta$ fixes $[\mc{J}\leq b-\bar{\varepsilon}]$, that is, $\eta(t, u)=u$ for all $t \in [0,1]$ and $\mc{J}(u)\leq b-\bar{\eps}$;
\item $\mc{J}$ is non-increasing along $\eta$, and in particular $\mc{J}(\eta(t,\cdot))\leq \mc{J}(\cdot)$ for all $t \in [0,1]$;
\item if $K_b= \emptyset$, then $\eta(1, [\mc{J}\leq b+\varepsilon])\subseteq [\mc{J}\leq b-\varepsilon]$;
\item if $K_b\not=\emptyset$, then
$$\eta(1,[\mc{J}\leq b+\varepsilon]\setminus \mc{O}) \subseteq [\mc{J}\leq b-\varepsilon]$$
and
$$\eta(1, [\mc{J}\leq b+\varepsilon]) \subseteq [\mc{J}\leq b-\varepsilon]\cup \mc{O};$$
\item if \hyperref[(F5)]{\textnormal{(F5)}} holds, then $\eta(t, \cdot)$ is $\G$-equivariant, i.e. it is odd.
\end{enumerate}
\end{Proposition}

As in Section \ref{section:4..1}, for any $n \in \N^*$ we define $\Gamma_n:=\Gamma_n(\lambda)$.
We note that $\Gamma_n\neq\emptyset$ is shown in Proposition \ref{S:4.1}.
Now our Theorem \ref{S:1.3} can be obtained through the arguments given in \cite{Rab0}.
Here we just give the definition of another minimax classes $\Lambda_n^m$, which ensures the multiplicity of solutions.
We set for $n\in\N^*$
%\begin{align*}
%\Lambda_n:=\{A=\Theta(\overline{D_{n+l}\setminus Y}) \mid \; & l\in \N^*, \; \Theta\in \Gamma_{n+l}(\lambda), \\
%& Y\subseteq D_{n+l}\setminus \{0\} \; \textnormal{ is closed,}\\
%& \textnormal{symmetric and $\genus(Y)\leq l$} \}
%\end{align*}
%and
%$$c_n:= \inf_{A\in \Lambda_n(\lambda)} \sup_{A} \mc{J}.$$
%Then we have $\{ \gamma(D_n)|\, \gamma\in\Gamma_n\}\subset \Lambda_n$ and we can also see that
%	$$	0 < c_1 \leq c_2\leq \cdots \leq c_n\leq \cdots.	$$
\begin{align*}
\Lambda_n:=\{A=\Theta(\overline{D_{n+l}\setminus Y}) \mid \; & l\in \N^*, \; \Theta\in \Gamma_{n+l}, \\
& Y\subseteq D_{n+l}\setminus \{0\} \; \hbox{ is closed, symmetric in $0$}\\
& \qquad \qquad \qquad \quad \; \hbox{ and $\genus(Y)\leq l$} \}
\end{align*}
and
$$c_n:= \inf_{A\in \Lambda_n} \sup_{A} \mc{J}.$$
Then we have $\{ \gamma(D_n)|\, \gamma\in\Gamma_n\}\subset \Lambda_n$ and we can also see that
	$$	0 < c_1 \leq c_2\leq \cdots \leq c_n\leq \cdots .	$$
Thus we have the following result.
\begin{Proposition} \label{S:5.3}
Assume \hyperref[(F1)]{\textnormal{(F1)}}--\hyperref[(F5)]{\textnormal{(F5)}}. Let $n \in \N^*$ and $m>0$. Then
\begin{itemize}
\item[(i)] $\Lambda_n \neq \emptyset$ and $c_n\leq c_{n+1}$.
\item[(ii)] Let $A\in \Lambda_n$ and $Z\subset H^1_r(\R^N)$ be $\G$-invariant, closed, 
and such that $0 \notin \overline{Z}$ and $\genus(\overline{Z})\leq i<n$. Then $\overline{A\setminus Z} \in \Lambda_{n-i}$.
\item[(iii)] $c_n$ is a critical value of $\mc{J}$. Moreover
	$$	c_n\to+\infty \quad \hbox{as}\ n\to+\infty.	$$
In particular, $\mc{J}$ has an unbounded sequence of critical values.
\end{itemize}
\end{Proposition}

\claim Proof.
Using Proposition \ref{S:5.2}, the proof can be given along the lines in \cite{Rab0}. See also \cite{CT2}.
\QED

\bigskip

\claim Proof of Theorem \ref{S:1.3}.
Theorem \ref{S:1.3} follows from Proposition \ref{S:5.3}.
\QED

%%%%%%%%%%%%%%%%%%%%%%%%%%%%%%%%%%%%%%%%%%%%%%%%%%%%%%%

%\phantomsection
%\addcontentsline{toc}{section}{Some open problems}
%
%\section*{Some open problems}
%
%\begin{itemize}
%
%\item Can we extend the existence results for the Choquard unconstrained problem (both existence and multiplicity) %\eqref{..} even 
%to general logarithmic-type nonlinearities in the spirit of \cite{Med2}? When the setting is local, in \cite{Med2} (see also \cite{GalSch} for recent results in the fractional setting) they make assumptions only on the positive part of the nonlinearity, letting the negative part be almost free. % (and including, in this way, also the "infinite mass" case). 
%What are the corresponding assumptions in the case of the Riesz potential?
%
%\item Can we apply our construction of multidimensional paths to other interesting kernels? For example positive kernels related to non-positive functionals, or sign-changing kernels, or kernels with no homogeneity properties.
%
%\end{itemize}

%%%%%%%%%%%%%%%%%%%%%%%%%%%%%%%%%%%%%%%%%%%%%%%%%%%%%%%
%%%%%%%%%%%%%%%%%%%%%%%%%%%%%%%%%%%%%%%%%%%%%%%%%%%%%%%
%%%%%%%%%%%%%%%%%%%%%%%%%%%%%%%%%%%%%%%%%%%%%%%%%%%%%%%

\chapter{Doubly nonlocal equations: qualitative and quantitative results} % for doubly nonlocal equations}
\label{chap_doubly}

This Chapter is dedicated to the study of the following \emph{fractional Choquard equation}
$$%\begin{equation}\label{eq_introduction}
(-\Delta)^s u + \mu u = \big(I_{\alpha}*F(u)\big)F'(u) \quad \hbox{in $\R^N$}
$$%\end{equation}
where $N\geq 2$, $s \in (0,1)$, $\alpha \in (0,N)$, $\mu>0$ and $F\in C^1(\R)$ is a general nonlinearity, in the spirit of Berestycki and Lions assumptions.
After having achieved existence of positive solutions and ground states, we will focus on the study of some qualitative properties of these solutions: boundedness, regularity, $L^1$-summability, positivity, radial symmetry and asymptotic decay. We will stress how the interplay between a fractional framework and a nonlocal nonlinearity, generally nonhomogeneous, obstructs the application of classical techniques. Some results generalize the ones presented in \cite{DSS1} and extend \cite{MS2,BKS}; in particular, some new results are stated also for the limiting case $s=1$. In addition, we will see that the interaction of the two nonlocalities arises a new critical threshold. % phenomena.

\medskip

This Chapter is mainly based on the papers:
\cite{CGT2} (see also \cite{CGT1}) for Section \ref{sec_doubl_exist}, 
\cite{CGT3} for Sections \ref{sec_doubl_exist}, \ref{sec_proper_pohozaev_gs}, \ref{sec_doubly_boundedness}--\ref{sec_L1_sum}, \ref{sec_doub_superl}, 
\cite{CG1} for Sections \ref{sec_proper_pohozaev_gs}, \ref{sec_regularity_fp}-\ref{sec_doub_C1C2}, \ref{sec_doubl_shape},
%and 
\cite{Gal0} for Sections \ref{sec_sublin_case}--\ref{sec_proof_main_2},
and \cite{CGT5} for Section \ref{sec_Pohozaev}.

%%%%%%%%%%%%%%%%%%%%%%%%%%%%%%%%%%%%%%%%%%%%%%%%%%%%%%%
%%%%%%%%%%%%%%%%%%%%%%%%%%%%%%%%%%%%%%%%%%%%%%%%%%%%%%%

\section{An example of double nonlocality: collapse of boson stars}
\label{sec_boson_stars}

In Sections \ref{sec_introd_frac_lap} and \ref{sec_intro_choquard} we highlighted the importance in physics of the fractional Laplacian and of the Hartree-type terms. Combinations of the two arise as well in different frameworks: 
for example equations of the type
\begin{equation}\label{eq_introduction}
(-\Delta)^s u + \mu u = \big(I_{\alpha}*F(u)\big)f(u) \quad \hbox{in $\R^N$}
\end{equation}
where $N\geq 2$, $s \in (0,1)$, $\alpha \in (0,N)$, $\mu>0$ and $f=F'\in C(\R)$, 
can be found in quantum chemistry \cite{ArMe,DOS,HMT0} (see also \cite{CFHMT} for some orbital stability results): here \eqref{eq_introduction} appears in the study of the mean field limit of weakly interacting molecules and in the physics of multi-particle systems. In particular the equation applies to the study of graphene \cite{LuMoMu}, where the nonlocal nonlinearity describes the short time interactions between particles. 
Doubly nonlocal equations appear also in the dynamics of populations \cite{CDV}, where small or large values of $s$ better model specific environments.
%one of the most relevant applications arises in relativistic physics, when the nonlinearity describes the short time interactions between particles. 
%

One of the main applications anyway arises in the study of exotic stars: minimization properties related to \eqref{eq_introduction} play indeed a fundamental role in the mathematical description of the dynamics of pseudorelativistic boson stars %\cite{ElSc0} 
and their gravitational collapse \cite{ElSc0,HL0,Len0,Len1,LeLe,LiYa0,FJL1, FJL2, FroL0,FroL1}, as well as the evolution of attractive fermionic systems, such as white dwarf stars \cite{HLLS}. 
%
%the minimization related to the problem (\ref{problem}) plays a fundamental role in the mathematical description of the gravitational collapse of boson stars \cite{LY0,FJL1, FJL2}. 
In fact, the study of the ground states to \eqref{eq_introduction} gives information on the size of the critical initial conditions for the solutions of the corresponding pseudorelativistic equation \cite{Len0}, where a critical value is given by the Chandrasekhar limiting mass. 
In particular, when $s=\frac{1}{2}$, $N=3$, $\alpha=2$ and $f(u)=u$, we obtain
\begin{equation}\label{eq_boson_massless}
\sqrt{-\Delta}u + \mu u = \left(\frac{1}{4 \pi |x|}*u^2\right) u \quad \hbox{in $\R^3$}
\end{equation}
related to the so called \emph{massless boson stars equation} \cite{FraL,LeLe,HL0}, where the pseudorelativistic operator $\sqrt{-\Delta + m}$ collapses to the square root of the Laplacian. 
Here $f(t)=|t|^{r-2}t$ with $r=2$ is $L^2$-critical: in this Chapter, when dealing with the mass-constrained problem, we essentially address the subcritical case $r\in (\frac{5}{3}, 2)$, but we believe that this result, together with the developed minimax tools, can be a first step towards the study of the $L^2$-mass critical (and supercritical) case, since for these problems the minimization approach is generally not well posed. Moreover, the high generality assumed on the function $f$ could be useful in the study of different physical problems.

\medskip

Mathematically, concerning the fractional Schr\"odinger equation with Hartree nonlinearity, we mention the papers \cite{DSS1,DSS2} where D'Avenia, Siciliano and Squassina considered the case of pure power nonlinearities and obtained existence and qualitative properties of the solutions.
 We mention also \cite{CFHMT,Geo0} for some orbital stability results, \cite{CHHO} for a Strichartz estimates approach, and \cite{CAB} for the unidimensional case. 
Other results can be found in \cite{SGY, BBMP, Luo0} for superlinear nonlinearities, in \cite{HRZ} for some local perturbation, in \cite{HR0, MuSr0} for critical equations and in \cite{YZ0} for concentration phenomena with strictly noncritical and monotone sources.

The existence of $L^2$-normalized solutions was investigated when $F(t)=|t|^p$ in \cite{Wu0} (see also \cite{GZ0,GS1} for $L^2$-supercritical Cauchy problems by scattering), while in \cite{ChLiu} it has been addressed the non-autonomous unconstrained case.
In \cite{DFQ}, symmetry and monotonicity of positive solutions are shown for the fractional Hartree equation for $\mu=0$ and a critical power nonlinearity, by means of the direct method of moving planes.
Regularity results for a class of doubly nonlocal equations on bounded domains are obtained in \cite{GDS}.

Some theoretical aspects related to the study of doubly nonlocal equations, both in the operator and in the source, remain open for general nonlinearities $F$, in particular when $F$ is not a power function or $F$ is odd.
%In the present Chapter we are interested to derive some qualitative properties of the solutions to \eqref{eq_introduction}, also in these special cases.

\medskip

In the present Chapter we are interested to derive some qualitative properties of the solutions to \eqref{eq_introduction}, also in these special cases.
In particular, after having stated existence of free and normalized solutions, we will focus our attention on the study of regularity of solutions (boundedness, $L^1$-summability, H\"older continuity, differentiability), moving then to positivity and symmetry of ground states, to tackle at the end the asymptotic behaviour at infinity. The precise statements will be presented throughout the Chapter: these results generalize some of the ones in \cite{DSS1} from the case of power functions to general nonlinearities of Berestycki-Lions type; moreover, we extend some results of \cite{MS2} to the fractional framework, and some results of \cite{BKS} to the Choquard framework.

The achieving of these results requires some technical effort in order to deal with the two nonlocalities and their interaction, as well as the nonhomogeneity and the nonregularity of the function $f$. In particular, we highlight some of the difficulties that arise in this general framework.

%\tr{Scrivi di più? SCRIVI i vari Teoremi} %CCOMMENT NOW

\smallskip

%\tb{
%We point out some difficulties which arise in this general framework, especially 
In the proof of the positivity, for instance (as well as in the proof of the existence), the presence of the fractional power of the Laplacian does not allow to use the fact that every solution satisfies the Pohozaev identity to conclude that, if $|u|$ is a solution, then it satisfies the Pohozaev identity; moreover, the conservation of the norm of the gradient does not hold anymore, i.e. $\norm{\nabla|u|}_2 = \norm{\nabla u}_2$ is not generally true in the fractional framework, and an inequality is needed. In addition, when dealing with $f$ even %-- contrary to the case $f$ odd -- 
other information about $u$ are lost through inequalities; this is not the case when dealing with %$s=1$ and 
$f$ odd \cite{MS2}.
Furthermore, the presence of the Choquard term, which scales differently from the $L^2$-norm term, does not allow to implement the classical minimization argument of \cite{CGM, BL1} (see \eqref{eq_min_probl}), which is useful to deal with the absolute value of $u$. Similarly, the nonhomogeneity of the nonlinearity $f$ obstructs the minimization approach of \cite{DSS1,MS0}.
Thus, a new approach is needed, and it relies on a fiber map %, introduced in Section \ref{sec_qualit}, 
which sends solutions to the Pohozaev set (see Proposition \ref{prop_sign_unconstr}). 
When dealing with $f$ even, this technique allows to treat also the case $s=1$, generalizing \cite{MS2}.
%}

Regarding the $L^1$-summability, the possibility of including a critical behaviour in zero (that is, $F(t) \sim t^{2^{\#}_{\alpha}}$) is not relevant when dealing with pure power functions \cite{DSS1,MS0}, since no solution exists in this case: this growth is instead relevant for general $f$ (for example, suitable sum of powers). Contrary to the case of noncritical nonlinearities, when $f$ is critical it is not possible to implement a simple bootstrap argument to achieve that every solution is in $L^1$: a new method is thus needed, and it is based on a suitable combination of bootstrap argument and fixed point theorems (see Proposition \ref{prop_u_L1}). The study of this case is new even for $s=1$, improving \cite{MS0,MS2}.

When studying the asymptotic behaviour of solutions, especially when $f$ has a sublinear growth, the interaction of the two nonlocalities is quite strong, and new phenomena arise: indeed, contrary to the local case $s=1$ \cite{MS0}, here the effect of the fractional Laplacian and of the Choquard term give rise to a new threshold depicting the qualitative profile of ground states at infinity (see Theorem \ref{thm_main}). From a technical point of view, new difficulties arise related to the explicit computation of the fractional Laplacian, and to the computation of concave powers, requiring a more delicate analysis and the implementation of new inequalities (see Sections \ref{sec_hypergeo} and \ref{sec_chain_rule}). This result is new even for power functions, improving \cite{DSS1}.

\smallskip

We refer to the following Sections for the detailed statements of the results.

\bigskip

The Chapter is organized as follows. 
In the remaining part of the Section we will briefly give a physical interpretation of equation \eqref{eq_introduction} in the framework of gravitational collapses.
In Section \ref{sec_doubl_exist} we will deal with existence of %positive 
solutions, both for the unconstrained problem and the constrained one, by highlighting some approach different from the ones developed in Chapters \ref{chap_fract_normal} and \ref{chap_choq_multi}; some properties related to the energy minimum levels and existence of positive solutions will be then investigated in Section \ref{sec_proper_pohozaev_gs}.
Section \ref{sec_doubl_regul} will be devoted to the study of regularity of positive solutions, including boundedness and H\"older regularity; moreover we will gain $L^1$-summability of solutions through a combination of bootstrap and fixed point maps arguments.
Then in Section \ref{sec_doubl_shape} we will exploit these results in order to gain positivity and radial symmetry of Pohozaev minima, by the implementation of maximum principles on some fiber maps.
Afterwards, we will investigate in Section \ref{sec_asymptotic} the asymptotic decay of ground states, focusing especially on the case of $f$ sublinear, which raises some new phenomenon. 
Finally in Section \ref{sec_Pohozaev} we furnish a proof of the Pohozaev identity in the doubly nonlocal framework, by assuming the solutions merely $C^{1}$.

%%%%%%%%%%%%%%%%%%
\subsubsection{Physical derivation}

Here we want to show how equation \eqref{eq_introduction}, in the particular case $N=3$, $s=\frac{1}{2}$, $\alpha=2$ and $F(u)=\frac{1}{2}u^2$, can be derived from a significant physical framework regarding boson stars. Aim of this Section is just to give an idea of the process, without any aim of accuracy or rigors. We refer to \cite{Chd0, ElSc0, FTY, FJL1, HS0, Len1, LeLe, LT0, LiYa0, MiSc0, Ngu1} (see also \cite{FraL, LiYa1,Ngu2, Thir1} and \cite{FJL2, HLLS, Len0}) for complete expositions on the topic.

The goal is to show how the equation \eqref{eq_boson_massless}
%$$\sqrt{-\Delta}u + \mu u = \big(I_2*u^2\big)u \quad \hbox{in $\R^3$}$$
is strictly connected with the self-gravitational collapse of boson stars. Actually a similar derivation holds also for neutron stars and white dwarfs, with some little complications.

Let us consider thus a group of $n$ \emph{bosons} (i.e. particles with entire spin, described by symmetric functions, and which do not respond to the Pauli exclusion principle). We assume these bosons to form a \emph{boson star}, i.e. $n\gg 0$ and we assume most of them (i.e. up to $o(N)$ particles) to be close one to each other and moving at a fast speed: these particles are at a same coherent state $\psi$ (and for this reason called \emph{condensate}) and create a trap for the remaining particles.

Due to the high speed of the particles, we cannot ignore the special relativistic effect; on the other hand, since the masses are not too big, we can ignore the effect of general relativity (this is not the case, instead, of neutron stars). %, see ).
Thus we consider the total relativistic energy
$$E^2 = (pc)^2 + (mc^2)^2$$
summation of the kinetic energy and the energy at rest; here $p$ is the momentum, $m$ the mass of the boson particle at rest, $c$ the speed of light. Thus we obtain, in momentum representation,
$$E = \sqrt{|\xi|^2 c^2 + m^2 c^4};$$
passing through a quantization $p \mapsto -i \hbar \nabla$ to the coordinate representation (and setting $\hbar :=1$) we obtain the \emph{pseudorelativistic} operator
$$E= \sqrt{- c^2\Delta+m^2 c^4}.$$
We observe that, letting $c\to +\infty$ (that is, the velocities are far from the one of the light) we obtain the operator $-\frac{1}{2m} \Delta$, that is the nonrelativistic operator (i.e. the classical Laplacian). %see 
We set instead, from now on, $c:=1$ for the sake of simplicity. Thus
$$E=\sqrt{-\Delta + m^2}$$
which is formally defined through the Fourier symbol $\mc{F}^{-1}\big((|\xi|^2+m^2)^{1/2}\widehat{u}\big)$ (see also Remark \ref{rem_Bessel_kernels}).

We consider now the interaction between the particles: this interaction can be treated classically as a Newton two-bodies interaction, and thus given by the quantity
$$- k \frac{1}{|x_i-x_j|}$$
where $k$ is a coupling constant (proportional to $G$, the gravitational constant).

Thus we come up with the Hamiltonian of the system
$$\mc{H}_n:= \sum_{i=1}^n \sqrt{-\Delta_i + m^2} - \frac{k}{2} \sum_{i \neq j}^n \frac{1}{|x_i-x_j|}.$$
When dealing with dwarf stars, additional pieces given by the interaction (due to the Pauli exclusion principle) appear; anyway, for $n\gg0$, one can ignore these pieces: this is called the \emph{Hartree approximation of Hartree-Fock theory}. %see 

Now we are interested in what happens when $n\gg 0$, that is, when the particles act like a single body, in what is called the \emph{mean field limit}:
$$n\to +\infty, \quad k\to 0, \quad nk =const;$$
formally this is given by assuming that the state of motion $\psi_n(t)$ can be factorized at each $t$ -- fact that is not generally true, even if one starts from a factorized state at $t=0$. We highlight that the powers appearing in the relation $nk=const$ are typical of the boson star framework, and are indeed different in other frameworks (for instance, in the case of white dwarfs, we have $n^{3/2} k = const$).

By letting $n\to +\infty$ one can formally prove that, in some precise sense, the motion of the (single body) boson star converges to the motion of the following (time-dependent) PDE
\begin{equation}\label{eq_boson_dynam}
i u_t = \sqrt{-\Delta + m^2} u - %k 
\left( \int_{\R^3} \frac{u^2(y)}{|x-y|} dy\right) u \quad \hbox{ in $(0, +\infty) \times \R^3$}
\end{equation}
that is (up to constant)
$$i u_t = \sqrt{-\Delta + m^2} u - %k 
(I_2*u^2) u \quad \hbox{ in $(0, +\infty) \times \R^3$.}$$
We focus now on this equation, and on the corresponding energy functional
$$E(u):=\frac{1}{2} \int_{\R^3} |(-\Delta+m)^{1/4}u|^2 - \frac{ %k
1}{2} \int_{\R^3} \big(I_2*u^2)u^2.$$
By exploiting the Hardy-Littlewood-Sobolev %inequality 
and the fractional Gagliardo-Nirenberg inequalities we obtain
\begin{equation}\label{eq_boson_best_const}
\int_{\R^3} \big(I_2*u^2)u^2 \leq \overline{C} \norm{(-\Delta)^{1/4}u}_2^2\norm{u}_2^2,
\end{equation}
which combined with the trivial inequality $|\xi|^2 + m^2 \geq |\xi|^2$ gives
$$E(u) \geq \frac{1}{2} \norm{(-\Delta)^{1/4}u}_2^2 \big( 1- %k
\overline{C} \norm{u}_2^2\big)$$
for some $\overline{C}>0$. Set
$$\int_{\R^3} u^2 =: M$$
the total mass of the boson star (interpreting $u^2(x)$ as the density in $x\in \R^3$) we obtain
$$E(u) \geq \frac{1}{2} \norm{(-\Delta)^{1/4}u}_2^2 \big( 1- %k
\overline{C} M\big);$$
from this we see that $E(u)$ could be or be not bounded from below on the sphere $\{u \in H^{1/2}(\R^3)\mid \norm{u}_2^2=M\}$ depending on the size of $M$: this is actually a phenomenon related to the \emph{$L^2$-critical growth} $2=\frac{3+2+2\frac{1}{2}}{3} = 2^m_{2,\frac{1}{2}}$ in $N=3$.
More precisely, one can prove that there exists a constant $M_*$, related to the best constant of the inequality \eqref{eq_boson_best_const}, such that
$$\inf_{\norm{u}_2^2 = M} E(u) \; \parag{&\geq 0& \quad \hbox{if $M < M_*$}, \\ &=-\infty& \quad \hbox{if $M>M_*$}.}$$
As a further consequence, one might study the dynamical properties of $u(x,t)=e^{it} u(x)$, solution of \eqref{eq_boson_dynam}, showing that
$$ u=u(x,t) \; \parag{&\hbox{exists for each $t>0$} & \quad \hbox{if $M < M_*$}, \\ &\hbox{explodes in finite time}& \quad \hbox{if $M>M_*$}.}$$
This is why $M_*$, called \emph{Chandrasekhar mass}, is related to the \emph{self-gravitational collapse} of boson stars (i.e., the collapse due to their own gravity). One could show that $M_*$ is related to a number of particles of the size of $\sim 10^{38}$, that is, the number of particles that can be approximately found in a mountain.

As already highlighted, $M_*$ is related to the best constant of the inequality \eqref{eq_boson_best_const}. And one can show that the optimizers $Q$ of this inequality satisfy the following equation 
$$\sqrt{-\Delta} Q + \mu Q = \big(I_2*Q^2\big)Q \quad \hbox{in $\R^3$}$$
for some $\mu>0$ (actually $M_*$ equals the $L^2$ norm squared of $Q$). And this is the equation we study.

%%%%%%%%%%%%%%%%%%%%%%%%%%%%%%%%%%%%%%%%%%%%%%%%%%%%%%%
%%%%%%%%%%%%%%%%%%%%%%%%%%%%%%%%%%%%%%%%%%%%%%%%%%%%%%%

\section{Different approaches for the existence problem}
\label{sec_doubl_exist}

In this Section we briefly sketch how to get existence of \emph{free} and \emph{constrained} solutions. The techniques are based on the ideas of Chapters \ref{chap_fract_normal} and \ref{chap_choq_multi}.
Anyway we present here a different approach to handle the boundary of $\R_+$, instead of considering the change of variable $\mu=e^{\lambda}$. With this aim, we will give a proof of some details, referring to Chapters \ref{chap_fract_normal}-\ref{chap_choq_multi} for all the other proofs.

This first Section is based on the paper \cite{CGT2} and \cite{CGT3} (see also \cite{CGT1}). For multiplicity results we refer to \cite{CGT5}.

%%%%%%%%%%%%%%%%%%%%%%%%%%%%%%%%%%%%%%%%%%%%%%%%%%%%%%%

%\subsection{Different approaches for the existence problem}

\medskip

The first goal we address is to study the unconstrained problem of \eqref{eq_introduction} when $f$ satisfies the following set of assumptions of Berestycki-Lions type \cite{BL1}:
%\begin{itemize}
%\item[(f1)] $f \in C(\R, \R)$;
%\item[(f2)] we have
% $$i) \; \limsup_{t \to 0} \frac{|tf(t)|}{|t|^{\frac{N+ \alpha}{N}}} <+\infty, \quad
% ii) \; \limsup_{ |t| \to + \infty} \frac{|t f(t)|}{|t|^{\frac{N+ \alpha}{N-2s}}} <+\infty;$$
%\item[(f3)] $F(t)= \int_0^t f(\tau) d\tau$ satisfies
% $$i) \; \lim_{t \to 0} \frac{F(t)}{|t|^{\frac{N+ \alpha}{N}}} =0, \quad
% ii) \; \lim_{ |t| \to + \infty} \frac{F(t)}{|t|^{\frac{N+ \alpha}{N-2s}}} =0;$$
%%\item[(f3)] $F(t)= \int_0^t f(\tau) d\tau$ satisfies
%% $$i) \; \lim_{t \to 0} \frac{t\tr{f}(t)}{|t|^{\frac{N+ \alpha}{N}}} =0, \quad
%% ii) \; \lim_{ |t| \to + \infty} \frac{t\tr{f}(t)}{|t|^{\frac{N+ \alpha}{N-2s}}} =0;$$
%\item[(f4)] there exists $t_0 \in \R$, $t_0 \neq 0$ such that $F(t_0) \neq 0$.
%\end{itemize}
\begin{itemize}
\item[(F1)] \label{(F1s)}
$f \in C(\R, \R)$;
\item[(F2)] \label{(F2s)}
we have
 $$i) \; \limsup_{t \to 0} \frac{|tf(t)|}{|t|^{2^{\#}_{\alpha}}} <+\infty, \quad
 ii) \; \limsup_{ |t| \to + \infty} \frac{|t f(t)|}{|t|^{2^*_{\alpha,s}}} <+\infty;$$
\item[(F3)] \label{(F3s)}
$F(t)= \int_0^t f(\tau) d\tau$ satisfies
 $$i) \; \lim_{t \to 0} \frac{F(t)}{|t|^{2^{\#}_{\alpha}}} =0, \quad
 ii) \; \lim_{ |t| \to + \infty} \frac{F(t)}{|t|^{2^*_{\alpha,s}}} =0;$$
%\item[(f3)] $F(t)= \int_0^t f(\tau) d\tau$ satisfies
% $$i) \; \lim_{t \to 0} \frac{t\tr{f}(t)}{|t|^{\frac{N+ \alpha}{N}}} =0, \quad
% ii) \; \lim_{ |t| \to + \infty} \frac{t\tr{f}(t)}{|t|^{\frac{N+ \alpha}{N-2s}}} =0;$$
\item[(F4)] \label{(F4s)}
there exists $t_0 \in \R$, $t_0 \neq 0$ such that $F(t_0) \neq 0$.
\end{itemize}
We observe again that \hyperref[(F3s)]{\textnormal{(F3)}} implies that we are in a \emph{noncritical} setting: indeed the exponents $2^{\#}_{\alpha}=\frac{N+\alpha}{N}$ and $2^*_{\alpha,s}=\frac{N+\alpha}{N-2s}$ 
have been addressed in \cite{MS0} as critical for Choquard-type equations when $s=1$, and then generalized to $s\in (0,1)$ in \cite{DSS1}; we will assume the noncriticality in order to obtain the existence of a solution, while all the qualitative results in the following Sections will be given in a \emph{possibly critical} setting.

This unconstrained case was studied by \cite{DSS1} for a power nonlinearity and by \cite{Bha0} in the case of combined local and nonlocal power-type nonlinearities; see also \cite{SGY, Luo0, GTC}. 
 
%We deal first with existence of a ground state for \eqref{eq_introduction}, obtaining the following result.

\medskip

We obtain the following result.

\begin{Theorem}\label{th_INT_exist_unconstrained}
Assume \hyperref[(F1s)]{\textnormal{(F1)}}--\hyperref[(F4s)]{\textnormal{(F4)}}. 
Then there exists a radially symmetric weak solution $u$ of \eqref{eq_introduction}, which satisfies the Pohozaev identity:
\begin{equation}\label{eq_INT_Pohozaev}
\frac{N-2s}{2}\int_{\R^N} |(-\Delta)^{s/2}u|^2 + \frac{N}{2} \mu \int_{\R^N} u^2 - \frac{N + \alpha}{2} \int_{\R^N} (I_{\alpha}*F(u)) F(u) =0
\end{equation}
or equivalently
$$\frac{1}{2^*_{\alpha,s}} \int_{\R^N} |(-\Delta)^{s/2}u|^2 + \frac{ \mu}{2^{\#}_{\alpha}} \int_{\R^N} u^2 - \int_{\R^N} (I_{\alpha}*F(u)) F(u) =0. $$
This solution is of Mountain Pass type. 
%and minimizes the energy among all the solutions satisfying \eqref{eq_INT_Pohozaev}.
\end{Theorem}

%We refer to Section \ref{sec_unconstrained} for the precise meaning of \emph{weak solution}, of \emph{Mountain Pass type} and \emph{energy}, according to a variational formulation of the problem.

We point out some difficulties which arise in this framework. Indeed, the presence of the fractional %power of the 
Laplacian does not allow to use the fact that every solution satisfies the Pohozaev identity to conclude that a Mountain Pass solution is actually a (Pohozaev) ground state, as in \cite{JT0} (see Remark \ref{rem_Pohozaev_c_p}). On the other hand, the presence of the Choquard term, which scales differently from the $L^2$-norm term, does not allow to implement the classical minimization argument by \cite{CGM, BL1}. Finally, the nonhomogeneity of the nonlinearity $f$ obstructs the minimization approach of \cite{MS2,DSS1}. Thus, we need a new approach to get existence of solutions, and this can be done in the spirit of Chapters \ref{chap_fract_normal}-\ref{chap_choq_multi}.
We omit the details, refering to \cite{CGT3}.

\bigskip

%COMMENT NOW
%\tor{-------------------}

%In this work, we consider the problem \eqref{problem_choquard}, which presents some nonlocal characteristics in the source, as well as in the fractional diffusion. 

The next goal is to study the constrained problem, 
%\tb{[RIPETIZIONE]
%In this paper, 
i.e. we study the existence of solutions $(\mu, u) \in (0,+\infty)\times H_r^s(\R^N)$ to the nonlocal problem
	\begin{equation} \label{problem_choquard}
\left \{
\begin{aligned}
(- \Delta)^s u & + \mu u =(I_\alpha*F(u))f(u) \quad \hbox{in}\ \R^N, \\ 
&\int_{\R^N} u^2 \, dx = m, 
\end{aligned}
\right.
\end{equation}
where %$s \in (0,1)$, $N\geq 2$, $\alpha\in (0,N)$, $F\in C^1(\R,\R)$ with $f=F'$, $m>0$, and 
$\mu>0$ is a Lagrange multiplier, part of the unknowns. 
%}
% of \eqref{problem_choquard}.

In particular we assume \hyperref[(F1s)]{\textnormal{(F1)}}, \hyperref[(F4s)]{\textnormal{(F4)}} together with the stronger assumptions
%\begin{itemize}
%%\item[(f1)] $f \in C(\R, \R)$;
%\item[(cf2)] there exists $C >0$ such that for every $t \in \R$, 
%$$|t f(t)| \leq C \big(|t|^{\frac{N + \alpha}{N}} + |t|^{\frac{N + \alpha+2s }{N}}\big);$$
%\item[(cf3)] $F(t)= \int_0^t f(\tau) \, d\tau$ satisfies
% $$\lim_{t \to 0} \frac{F(t)}{|t|^{\frac{N+ \alpha}{N}}} =0, \quad
% \lim_{ t \to + \infty} \frac{F(t)}{|t|^{\frac{N+ \alpha+2s}{N}}} =0;$$
%%\item[(cf3)] $F(t)= \int_0^t f(\tau) \, d\tau$ satisfies
%% $$\lim_{t \to 0} \frac{t\tr{f}(t)}{|t|^{\frac{N+ \alpha}{N}}} =0, \quad
%% \lim_{ t \to + \infty} \frac{F(t)}{|t|^{\frac{N+ \alpha+2s}{N}}} =0;$$
%%\item[(f4)] there exists $t_0 \in \R$, $t_0 \neq 0$ such that $F(t_0) \neq 0$.
%\end{itemize}
\begin{itemize}
%\item[(f1)] $f \in C(\R, \R)$;
\item[(CF2)] \label{(CF2s)}
%there exists $C >0$ such that for every $t \in \R$, 
%$$|t f(t)| \leq C \big(|t|^{\frac{N + \alpha}{N}} + |t|^{\frac{N + \alpha+2s }{N}}\big);$$
$$i) \; \limsup_{t \to 0} \frac{|tf(t)|}{|t|^{2^{\#}_{\alpha}}} <+\infty, \quad
 ii) \; \limsup_{ |t| \to + \infty} \frac{|t f(t)|}{|t|^{2^m_{\alpha,s}}} <+\infty;$$
\item[(CF3)] \label{(CF3s)}
%$F(t)= \int_0^t f(\tau) \, d\tau$ satisfies
 $$i) \; \lim_{t \to 0} \frac{F(t)}{|t|^{2^{\#}_{\alpha}}} =0, \quad
 ii) \; \lim_{ |t| \to + \infty} \frac{F(t)}{|t|^{2^m_{\alpha,s}}} =0;$$
%\item[(cf3)] $F(t)= \int_0^t f(\tau) \, d\tau$ satisfies
% $$\lim_{t \to 0} \frac{t\tr{f}(t)}{|t|^{\frac{N+ \alpha}{N}}} =0, \quad
% \lim_{ t \to + \infty} \frac{F(t)}{|t|^{\frac{N+ \alpha+2s}{N}}} =0;$$
%\item[(f4)] there exists $t_0 \in \R$, $t_0 \neq 0$ such that $F(t_0) \neq 0$.
\end{itemize}
we remark that the exponent $2^m_{\alpha,s}=\frac{N+ \alpha+2s}{N}$ appears as an $L^2$-critical exponent for the fractional Choquard equations and the conditions \hyperref[(F1s)]{\textnormal{(F1)}}-\hyperref[(CF2s)]{\textnormal{(CF2)}}-\hyperref[(CF3s)]{\textnormal{(CF3)}}-\hyperref[(F4s)]{\textnormal{(F4)}} correspond to $L^2$-subcritical growth.

For this general class of nonlinearities of the Berestycki--Lions type \cite{BL1,MS2} %, satisfying (f1)--(f4), 
we introduce a Lagrangian formulation:
namely, set $\R_+=(0,+\infty)$, a radially symmetric solution $(\mu,u)	\in \R_+ \times H^s_r(\R^N)$ of \eqref{problem_choquard} corresponds to a critical point 
of the functional $\mc{I}^m : \R_+ \times H^s_r(\R^N)\to \R$ defined by 
$$%\begin{equation*}%\label{eq:1.6}
\mc{I}^m(\mu, u):=\half \int_{\R^N} |(-\Delta)^{s/2} u|^2\, dx -\half \int_{\R^N} (I_\alpha*F(u))F(u)\, dx+ 
\frac{\mu}{2} \bigl( \|u\|_2^2 -m \bigr).
$$%\end{equation*}

Using a %new 
variant of the Palais--Smale condition \cite{HT0,IT0},
which takes into account the Pohozaev identity, we will prove a deformation theorem 
which enables us to detect minimax structures in the product space $\R_+ \times H^s_r(\R^N)$ by means of a \emph{Pohozaev mountain}.
Our deformation arguments show that solutions without Pohozaev identity are suitably deformable, and thus they \emph{do not influence the topology} of the sublevels of the functional. This information could be relevant in a fractional framework since it is not known if the Pohozaev identity holds for general continuous $f$ and general values of $s\in (0,1)$.

%\smallskip

\medskip

We state our main results.
\begin{Theorem}\label{S:1.1_doubly}
Assume \hyperref[(F1s)]{\textnormal{(F1)}}-\hyperref[(CF2s)]{\textnormal{(CF2)}}-\hyperref[(CF3s)]{\textnormal{(CF3)}}-\hyperref[(F4s)]{\textnormal{(F4)}}. Then there exists $m_0\geq 0$ such that, for any $m>m_0$, the problem \eqref{problem_choquard} has a radially symmetric solution, which satisfies the Pohozaev identity \eqref{eq_INT_Pohozaev}.
\end{Theorem}

%\newpage%€

\begin{Theorem}\label{S:1.12_doubly}
Assume \hyperref[(F1s)]{\textnormal{(F1)}}-\hyperref[(CF2s)]{\textnormal{(CF2)}}-\hyperref[(CF3s)]{\textnormal{(CF3)}}, %-\hyperref[(F4s)]{\textnormal{(F4)}}, 
together with an $L^2$-subcritical growth at zero, i.e.,
\begin{itemize}
\item[\textnormal{(CF4')}] \label{(CF4')}
$ \lim_{t \to 0} \frac{F(t)}{|t|^{2^m_{\alpha,s}}} = + \infty$. %\tr{CF4'}
\end{itemize}

 Then, for any $m >0$, the problem \eqref{problem_choquard} has a radially symmetric solution, which satisfies the Pohozaev identity \eqref{eq_INT_Pohozaev}. 
\end{Theorem}

We naively notice that \hyperref[(CF4')]{\textnormal{(CF4')}} automatically implies \hyperref[(F4s)]{\textnormal{(F4)}}. 
We remark that, as in the local unconstrained case \cite{JT0}, the Mountain Pass solutions obtained in the above theorems are ground state solutions, that is, they have the least energy among all solutions; see Section \ref{sec_ground_state} for details. This fact gives a strong indication on the stability properties of the found solution \cite{FZ0,CFHMT}. 

\medskip

Here we find solutions satisfying automatically the Pohozaev identity: in Section \ref{sec_Pohozaev} we will prove that a general $C^1$ solution actually satisfies such relation.
%we refer to the paper in preparation \cite{CGT5} for some proof of the identity in the case of $C^1$ solutions, together with multiplicity results. 

%%%%%%%%%%%%%%%%%%%%%%%%

 \subsection{Dealing with the boundary} % Palais-Smale-Pohozaev condition}
\label{section_PSP}

In what follows, we will often denote %by $q$, the lower Hardy--Littlewood--Sobolev critical exponent and, by $p$, the $L^2$-critical exponent appearing in (f2)--(f3), i.e.,
$$q=2^{\#}_{\alpha}= \frac{N+\alpha}{N}, \quad p= 2^m_{\alpha,s}= \frac{N+ \alpha+2s}{N}.$$
Consider the functional
$$\mc{J}_{\mu}(u):= \half \int_{\R^N} |(-\Delta)^{s/2} u|^2\,dx -\half {\mathcal D}(u) + \frac{\mu}{2} \|u\|_2^2$$
with $\mc{D}(u)=\int_{\R^N} (I_\alpha*F(u))F(u)$. 
We notice that, by the Principle of Symmetric Criticality of Palais, the critical points of $\mc{J}_{\mu}$ are \emph{weak solutions} of \eqref{eq_introduction}.
Moreover, inspired by the Pohozaev identity \eqref{eq_INT_Pohozaev}, 
%$$ %\begin{equation*}\label{eq_Poh_ident} %eq_INT_Pohozaev
%\frac{N-2s}{2}\norm{(-\Delta)^{s/2}u}_{2}^2 + \frac{N}{2} \mu \|u\|_2^2 = \frac{N + \alpha}{2} \, {\mathcal D}(u)
%$$ %\end{equation*}
we define also the Pohozaev functional $\mc{P}_{\mu}: H^s_r(\R^N) \to \R$ by 
$$\mc{P}_{\mu}(u):=\frac{N-2s}{2}\norm{(-\Delta)^{s/2}u}_{2}^2 -\frac{N+ \alpha}{2}{\mathcal D}(u)+ \frac{N}{2} \mu \|u\|_2^2.$$

Here we highlight how to deal with the boundary of $\R_+$ without implementing the change of variable $\mu=e^{\lambda}$. More details can be found in \cite{CGT2}.

\smallskip

As a matter of fact, we notice that $\R_+\times H_r^s(\R^N)$ with the standard metric induced by $\R\times H_r^s(\R^N)$ is not complete, and thus it is not suitable for a deformation argument. 
Since $(\R_+,\frac{1}{x^2}dx^2 )$ is instead complete, it is natural to introduce a related metric on $\R_+\times H_r^s(\R^N)$. That is, we regard 
$$R:=\R_+\times H_r^s(\R^N)$$
as a Riemannian manifold with the metric
	$$	\big((\nu_1, w_1), (\nu_2,w_2)\big)_{T_{(\mu,u)}R} :=
	\frac{1}{\mu^2} \nu_1\nu_2 +(w_1,w_2)_{H_r^s}
	$$
for $(\nu_1,w_1)$, $(\nu_2,w_2)\in T_{(\mu,u)}R$, $(\mu,u)\in R$; 
it is standard to see that $(R, (\cdot,\cdot)_{TR})$ is a complete Riemannian manifold.
We regard thus $\mc{I}^m$ as a functional defined on $R$, and obtain
	$$	\norm{\big(\partial_\mu\mc{I}^m(\mu,u), 
		\partial_u\mc{I}^m(\mu,u)\big)}_{(T_{(\mu,u)}R)^*}^2
	= \mu^2 \abs{\partial_\mu\mc{I}^m(\mu,u)}^2
	+ \norm{\partial_u\mc{I}^m(\mu,u)}_{(H_r^s)^*}^2.
	$$

\begin{Definition}\label{defPSP}
For $b \in \R$, we say that $(\mu_j, u_j)_j \subset R= \R_+ \times H^s_r(\R^N)$ is a \emph{Palais-Smale-Pohozaev sequence} at level $b$ (in short, the $(PSP)_b$ sequence) if, as $j \to +\infty$,
\begin{eqnarray*}
	&&\mc{I}^m(\mu_j, u_j) \to b, \\
	&&	\norm{ \big(\partial_\mu\mc{I}^m(\mu,u), 
		\partial_u\mc{I}^m(\mu,u)\big)}_{(T_{(\mu,u)}M)^*}
		\to 0, \label{4.P} \\
	&&{\mathcal P}(\mu_j, u_j) \to 0, %\label{ps4} 
\end{eqnarray*}
or equivalently
\begin{eqnarray}
	&&\mc{I}^m(\mu_j, u_j) \to b, \label{ps1} \\
	&&\mu_j\cdot\partial_\mu \mc{I}^m(\mu_j, u_j) \to 0, \label{ps2}\\
	&&\partial_u \mc{I}^m(\mu_j, u_j) \to 0 \quad \hbox{strongly in $(H^s_r(\R^N))^*$}, \label{ps3} \\
	&&{\mathcal P}(\mu_j, u_j) \to 0. \label{ps4} 
\end{eqnarray}
We say that $\mc{I}^m$ satisfies the \emph{$(PSP)_b$ condition} if, 
for any $(PSP)_b$ sequence $(\mu_j, u_j)_j \subset \R_+ \times H^s_r(\R^N)$, it happens that $(\mu_j, u_j)_j$ has a strongly convergent subsequence in $\R_+ \times H^s_r(\R^N)$.
\end{Definition}

\begin{Remark}
Clearly, setting
	$$	\widetilde{\mc{I}}^m(\lambda,u) :=\mc{I}^m(e^\lambda,u), \quad \widetilde{\mc{I}}: 
		\R\times H_r^s(\R^N)\to \R,
	$$
we can observe that $\widetilde{\mc{I}}^m$ satisfies the $(PSP)_b$ in the sense of Definitions \ref{PSPcondition2} and \ref{S:3.1} if and only if $\mc{I}^m$ satisfies the $(PSP)_b$ condition in the sense of Definition \ref{defPSP} with $\mu_j:=e^{\lambda_j}$.
\end{Remark}

%$(\lambda_j,u_j)\subset \R\times H_r^s(\R^N)$ satisfies
%	\begin{eqnarray*}
%	&\partial_\lambda\widetilde\mc{I}^m(\lambda_j,u_j)\to 0, \\
%	&\norm{\partial_u\widetilde\mc{I}^m(\lambda_j,u_j)}_{(H_r^s)^*}
%		\to 0
%	\end{eqnarray*}
%if and only if $(\mu_j,u_j):=(e^{\lambda_j},u_j)$ satisfies \eqref{4.P}.

%We remark that this compactness condition takes the scaling properties of $\mc{I}^m$ into consideration through the Pohozaev functional ${\mathcal P}$. 
% We now show the following crucial result.

\smallskip

For the sake of completeness, we give here some details on the proof of the $(PSP)_b$ condition at strictly negative levels.
We emphasize again indeed that the $(PSP)_{b}$ condition does not hold at level $b=0$: it is sufficient to consider an infinitesimal sequence $(\mu_j, 0)$ 
with $\mu_j \to 0$.

\begin{Theorem}\label{PSP2}
Assume \hyperref[(F1s)]{\textnormal{(F1)}}-\hyperref[(CF2s)]{\textnormal{(CF2)}}-\hyperref[(CF3s)]{\textnormal{(CF3)}}. Let $b <0$. Then $\mc{I}^m$ satisfies the $(PSP)_b$ condition on $\R_+\times H_r^s(\R^N)$.
\end{Theorem}

\claim Proof. 
Let $b <0$ and $(\mu_j, u_j)_j \subset \R \times H^s_r(\R^N)$ be a sequence satisfying \eqref{ps1}--\eqref{ps4}. 
First we note that, by \eqref{ps2}, we have
\begin{equation}\label{ps5}
	\mu_j\big(\|u_j\|_2^2-m\big)\to 0.
\end{equation}

\smallskip

\noindent
\textbf {Step 1:} \emph{$\liminf_{j\to\infty} \mu_j>0$ and $\norm{u_j}_2^2\to m$.} 
%Forse basta questo step, togliere gli altri? COMMENT NOW
%mdpi: is the italic and noindent necessary in the whole text? %I confirm.
%\smallskip
%
%
\\
By \eqref{ps4} and \eqref{ps1}, we have
\begin{align*}
	o(1) &= {\mathcal P}(\mu_j, u_j) = \frac{N-2s}{2} \|
	(-\Delta)^{s/2} u_j \|_2^2 +\\ 
& \quad+ (N+ \alpha) \Bigl(\mc{I}^m(\mu_j, u_j) -\frac{1}{2}
\|(-\Delta)^{s/2} u_j\|_2^2 
- \frac{\mu_j}{2} \bigl(\|u_j\|_2^2 -m \bigr) \Bigr) +
\frac{N}{2} \mu_j \|u_j\|_2^2 \\ 
&=
- \frac{\alpha +2s}{2} 
\|(-\Delta)^{s/2} u_j\|_2^2 +(N+ \alpha) (b + o(1)) +\frac{N}{2} \mu_j m +o(1);
\end{align*}
here we have used \eqref{ps5}. Since $b<0$, we have $\liminf_{j\to\infty} \mu_j>0$.
Thus \eqref{ps5} implies $\norm{u_j}_2^2\to m$.

\smallskip

\noindent
\textbf{Step 2:} \emph{$\|(-\Delta)^{s/2} u_j \|_2^2$ and $\mu_j$ are bounded.} 
%
%\smallskip
%
\\
Since $\epsilon_j\equiv \| \partial_u \mc{I}^m(\mu_j, u_j)\|_{(H_r^s(\R^N))^*}\to 0$,
we have
\begin{equation}\label{limit}
\|(-\Delta)^{s/2} u_j\|_2^2 - \int_{\R^N} (I_\alpha \ast F(u_j)) f(u_j) u_j \, dx 
 + \mu_j \|u_j\|_2^2 \leq \epsilon_j\|u_j\|_{H^s_r}.
\end{equation}

Note that $\frac{2Np}{N+ \alpha} \in (2, 2_s^*)$. Moreover, we observe that, by \hyperref[(CF3s)]{\textnormal{(CF3)}}, for $\delta>0$ fixed,
%, when fixed, %mdpi: please check intended meaning has been retained. %I confirm
there exists $C_\delta>0$ such that
$$
|F(t)| \leq \delta |t|^p + C_\delta |t|^{\frac{N+ \alpha}{N}}, \quad t \in \R,
$$
where $p=\frac{N+ \alpha+2s}{N}$, and thus 
$$
\|F(u_j)\|_{\frac{2N}{N+ \alpha}} \leq \delta \| |u_j|^p \|_{\frac{2N}{N+ \alpha}} 
+ C_\delta \| |u_j|^\frac{N+\alpha}{N}\|_{\frac{2N}{N+ \alpha}} 
= \delta \| u_j \|_{\frac{2Np}{N+ \alpha}}^p + C_\delta \| u_j\|_2^{\frac{N+ \alpha}{N}}.
$$
Therefore, by \hyperref[(CF2s)]{\textnormal{(CF2)}} we have
%\begin{eqnarray*}
%\lefteqn{\int_{\R^N} (I_\alpha*\abs{F(u_j)})\abs{f(u_j)u_j}\,dx }
%\\	&\leq&
%	C\norm{F(u_j)}_{\frac{2N}{N+\alpha}}\norm{f(u_j)u_j}_{\frac{2N}{N+\alpha}} \\
%	&\leq& C \left(\delta\norm{u_j}_{\frac{2Np}{N+\alpha}}^{p}
%	+ C_\delta\norm{u_j}_2^\frac{N+\alpha}{N} \right)
%	\cdot C' \left(\norm{u_j}_{\frac{2Np}{N+\alpha}}^{p}
%	+ \norm{u_j}_2^\frac{N+\alpha}{N} \right) \\
%	&=& CC'\delta \norm{u_j}_{\frac{2Np}{N+\alpha}}^{2p}
%	+ CC'(\delta+C_\delta)\norm{u_j}_{\frac{2Np}{N+\alpha}}^p
%	\norm{u_j}_2^\frac{N+\alpha}{N}
%	+ CC'C_\delta \norm{u_j}_2^{\frac{2(N+\alpha)}{N}} \\
%	&=& CC'\delta \norm{u_j}_{\frac{2Np}{N+\alpha}}^{2p}
%	+ CC'(\delta+C_\delta) \left(\frac{\delta}{2}
%	\norm{u_j}_{\frac{2Np}{N+\alpha}}^{2p}
%	+\frac{1}{2\delta}\norm{u_j}_2^{\frac{2(N+\alpha)}{N}} \right)
%	+ CC'C_\delta \norm{u_j}_2^{\frac{2(N+\alpha)}{N}} \\
%	&\leq& C''\delta\norm{u_j}_{\frac{2Np}{N+\alpha}}^{2p}
%	+C''_\delta\norm{u_j}_2^{\frac{2(N+\alpha)}{N}} 
%\end{eqnarray*}
\begin{eqnarray*}
\lefteqn{\int_{\R^N} (I_\alpha*\abs{F(u_j)})\abs{f(u_j)u_j}\,dx }
\\	&\leq&
	C\norm{F(u_j)}_{\frac{2N}{N+\alpha}}\norm{f(u_j)u_j}_{\frac{2N}{N+\alpha}} \\
	&\leq& C' \left(\delta\norm{u_j}_{\frac{2Np}{N+\alpha}}^{p}
	+ C_\delta\norm{u_j}_2^\frac{N+\alpha}{N} \right)
	\cdot \left(\norm{u_j}_{\frac{2Np}{N+\alpha}}^{p}
	+ \norm{u_j}_2^\frac{N+\alpha}{N} \right) \\
	&=& C'\delta \norm{u_j}_{\frac{2Np}{N+\alpha}}^{2p}
	+ C'(\delta+C_\delta)\norm{u_j}_{\frac{2Np}{N+\alpha}}^p
	\norm{u_j}_2^\frac{N+\alpha}{N}
	+ C'C_\delta \norm{u_j}_2^{\frac{2(N+\alpha)}{N}} \\
	&=& C'\delta \norm{u_j}_{\frac{2Np}{N+\alpha}}^{2p}
	+ C'(\delta+C_\delta) \left(\frac{\delta}{2}
	\norm{u_j}_{\frac{2Np}{N+\alpha}}^{2p}
	+\frac{1}{2\delta}\norm{u_j}_2^{\frac{2(N+\alpha)}{N}} \right)
	+ C'C_\delta \norm{u_j}_2^{\frac{2(N+\alpha)}{N}} \\
	&\leq& C''\delta\norm{u_j}_{\frac{2Np}{N+\alpha}}^{2p}
	+C''_\delta\norm{u_j}_2^{\frac{2(N+\alpha)}{N}} 
\end{eqnarray*}
and thus, by the fractional Gagliardo--Nirenberg inequality \eqref{GN}, with $r = \frac{2Np}{N+\alpha}$ and $\beta =\frac{1}{p}$, we derive 
\begin{eqnarray*}
\lefteqn{\|(-\Delta)^{s/2} u_j \|_2^2 + \mu_j \|u_j\|_2^2 
\leq \int_{\R^N} (I_\alpha \ast |F(u_j)|) |f(u_j) u_j| \, dx + \epsilon_j\|u_j\|_{H^s_r} }\\	 
&\leq &
 C'' \delta \|
 (-\Delta)^{s/2} u_j\|^2_2 \|u_j\|_2^{2(p-1)}
	+ C_\delta'' \|u_j \|_2^{\frac{2(N + \alpha)}{N}} + \epsilon_j\|u_j\|_{H^s_r}.
\end{eqnarray*}
Since $\|u_j\|_2^2=m+o(1)$, we get
\begin{eqnarray*} 
\lefteqn{\big(1- C'' \delta(c+o(1))^{p-1}\big)
	\|(-\Delta)^{s/2} u_j\|_2^2 + \mu_j\big(m+o(1)\big)} \\
 &\leq& C''_\delta \big(m+o(1)\big)^\frac{N+\alpha}{N} + \epsilon_j\big(\|
 (-\Delta)^{s/2} u_j\|_2^2+m+o(1)\big)^{1/2}.
\end{eqnarray*}
For a small enough $\delta$, we have the boundedness of $\|(-\Delta)^{s/2} u_j\|_2$ and $\mu_j$.

\smallskip

\noindent
\textbf{Step 3:} \emph{Convergence in $\R_+ \times H^s_r(\R^N)$.} 
%
%\smallskip
%
\\
By Steps 1-2, the sequence $(\mu_j, u_j)_j$ is bounded in $\R_+ \times H^s_r(\R^N)$ and thus, after extracting a subsequence denoted in the same way, we may assume that $\mu_j \to \mu_0>0$ and $u_j \wto u_0$ weakly in $H^s_r(\R^N)$ for some $(\mu_0, u_0) \in \R_+\times H^s_r(\R^N)$. 
 
\smallskip

\noindent
\textbf{Step 4:} \emph{Conclusion.} 
%\smallskip
\\
Taking into account the assumptions \hyperref[(F1s)]{\textnormal{(F1)}}--\hyperref[(F4s)]{\textnormal{(F4)}}, we obtain by Proposition \ref{prop_converg_generiche_nonloc}
$$
\int_{\R^N} (I_\alpha \ast F(u_j)) f(u_j) u_0 \, dx
\to 
\int_{\R^N} (I_\alpha \ast F(u_0)) f(u_0) u_0 \, dx
$$
and 
$$
\int_{\R^N} (I_\alpha \ast F(u_j)) f(u_j) u_j\, dx
\to 
\int_{\R^N} (I_\alpha \ast F(u_0)) f(u_0) u_0\, dx.
$$
Thus, we derive that 
$
\langle \partial_u \mc{I}^m(\mu_j, u_j), u_j \rangle \to 0$ 
and 
$
\langle \partial_u \mc{I}^m(\mu_j, u_j), u_0 \rangle \to 0$,
and hence
$$\|(- \Delta)^{s/2} u_j \|_2^2 + \mu_0 \|u_j\|^2_2 \to 
\|(- \Delta)^{s/2} u_0 \|^2_2
 + \mu_0 \|u_0\|^2_2 
$$ 
which implies $u_j \to u_0$ strongly in $H^s_r(\R^N)$.
\QED

%\medskip

%\begin{Remark}
%We emphasize that the $(PSP)_{b}$ condition does not hold at level $b=0$; it is sufficient to consider an infinitesimal sequence $(\mu_j, 0)$ 
%with $\mu_j \to 0$.
%\end{Remark}

%COMMENT NOW
\bigskip
%\tr{Come diventa questa PSP per l'aumentato? Non lo facciamo nell'altro articolo. Forse è da scrivere qui}

Now we define a metric on the Hilbert manifold
$$M := \R \times R = \R \times \R_+ \times H^s_r(\R^N)$$
by setting 
%\begin{align*}
%{\|(\alpha, \nu, h)\|}_{(\theta,\mu, u)}^2 :=& \pabs{\left(\alpha, \nu,\norm{h(e^{-\theta} \cdot)}_{H^s_r(\R^N)}\right)}^2 \\
%=& \, \alpha^2 + \frac{1}{\mu^2}\nu^2 + e^{N\theta} \norm{h}_2^2+ e^{(N-2s)\theta} \norm{(-\Delta)^{s/2} h}_2^2 
%\end{align*}
$${\|(\alpha, \nu, h)\|}_{(\theta,\mu, u)}^2 := \, \alpha^2 + \frac{1}{\mu^2}\nu^2 + e^{N\theta} \norm{h}_2^2+ e^{(N-2s)\theta} \norm{(-\Delta)^{s/2} h}_2^2 $$
for any $(\alpha, \nu, h) \in T_{(\theta,\mu,u)} M \equiv \R \times \R_+ \times H^s_r(\R^N)$.
We also denote the dual norm on 
$T^*_{(\theta,\mu,u)}M$ by $\|\cdot \|_{(\theta,\mu, u), *}$. 
We notice that 
${\|(\cdot, \cdot, \cdot)\|}_{(\theta,\mu, u)}^2$ depends both on $\theta$ and $\mu$ (but not on $u$). % and we can write ${\|(\cdot, \cdot, \cdot )\|}_{(\theta,\mu, \cdot)}^2$.
%Moreover for any $(\alpha, \nu, h) \in T_{(\theta,\mu,u)}M$ and $\beta \in \R$ we have
%\begin{equation}\label{eq_shift_norma}
%{\|(\alpha, \nu, h(e^\beta x))\|}_{(\theta + \beta,\cdot, \cdot)}^2=
%{\|(\alpha, \nu, h)\|}_{(\theta,\cdot, \cdot)}^2.
%\end{equation}
Furthermore we define the standard distance between two points $\dist_M$ as the infimum of length of curves connecting the two points.
%, namely
%$$ \dist_M\big((\theta_0, \lambda_0, h_0), (\theta_1, \lambda_1, h_1)\big) := 
%\inf_{\gamma \in \mathcal{G}} \int_0^1 \|\dot \gamma(t)\|_{\gamma(t)} dt $$
%where 
%$$\mathcal{G} :=\left\{\gamma \in C^1([0,1],M) \, \middle | \, \gamma(0)= (\theta_0, \lambda_0, h_0), \gamma(1)= (\theta_1, \lambda_1, h_1) \right\}.$$
%
%Observe that, if $\sigma$ is a path connecting $(\alpha_0, \nu_0, h_0)$ and $(\alpha_1, \nu_1, h_1)$, then by \eqref{eq_shift_norma} $\tilde{\sigma}(t):=(\sigma_1(t)+\beta, \sigma_2(t), (\sigma_3(t))(e^{\beta}\cdot))$ is a path connecting $(\alpha_0 +\beta, \nu_0, h_0(e^{\beta}\cdot))$ and $(\alpha_1+\beta, \nu_1, h_1(e^{\beta}\cdot))$ with same length, and hence
%\begin{equation}\label{eq_shift_distance}
%\dist_M\big((\alpha_0, \nu_0, h_0), (\alpha_1, \nu_1, h_1)\big) = \dist_M\big((\alpha_0 +\beta, \nu_0, h_0(e^{\beta}\cdot)), (\alpha_1+\beta, \nu_1, h_1(e^{\beta}\cdot))\big).
%\end{equation}

On $M$ we consider the augmented functional
$$\mc{H}^m(\theta, \mu, u) := \mc{I}^m(\mu, u(e^{-\theta}\cdot)) ;$$
denoted $D:=(\partial_\theta,\partial_\mu,\partial_u)$, 
%; a direct computation shows that
%$$D\mc{H}^m(\theta, \mu, u)(\alpha,\nu,h) 
%= \mc{P} (\mu, u(e^{-\theta} \cdot))\alpha +\partial_{\mu} \mc{I}^m(\mu, u(e^{-\theta}\cdot))\nu
%+\partial_u \mc{I}^m(\lambda, u(e^{-\theta}\cdot))h(e^{-\theta} \cdot)$$
%and 
 we obtain
\begin{eqnarray*}
\lefteqn{\|{D\mc{H}^m(\theta, \mu, u)\|}_{(\theta, \mu, u),*}^2 }\\
%	&=& \left|\left(\mathcal{P}(\lambda, u(e^{-\theta} \cdot)), \partial_{\lambda}\mc{I}^m(\lambda, u(e^{-\theta}\cdot)), \norm{\partial_u \mc{I}^m(\lambda, u(e^{-\theta}\cdot))}_{(H^s_r(\R^N))^*}\right)\right|^2 \\
	&=& \abs{\mathcal{P}(\mu, u(e^{-\theta} \cdot))}^2 + \mu^2 \abs{\partial_{\mu}\mc{I}^m(\mu, u(e^{-\theta}\cdot))}^2 
		+ \norm{\partial_u \mc{I}^m(\mu, u(e^{-\theta}\cdot))}_{(H^s_r)^*}^2 .
\end{eqnarray*}
Finally, defined
$$\tilde{K}_b :=\big \{ (\theta, \lambda, u) \in M \mid \mc{H}^m(\theta, \lambda, u)=b,\, D \mc{H}^m(\theta, \lambda,u)=0\big\}$$
the set of critical points at level $b$ of $\mc{H}^m$, we deduce the following.
%\begin{equation}\label{eq_confronto_K}
%\tilde{K}_b = \big\{(\theta, \lambda, u(e^{\theta} \cdot)) \mid (\lambda, u)\in K^{PSP}_b, \; \theta \in \R\big\}.
%\end{equation}
\begin{Proposition}%\label{PSPtilde}
%	Assume \textnormal{(g1)--(g3)}. 
	Let $b \in \R$, $b <0$. Then the functional $\mc{H}^m$ satisfies the following Palais-Smale type condition $(\widetilde{PSP})_b$. 
	That is, for each sequence $(\theta_j, \mu_j, u_j)_j$ such that
	$$\mc{H}^m(\theta_j, \mu_j, u_j) \to b,$$
	$$\norm{D \mc{H}^m(\theta_j, \mu_j, u_j)}_{(\theta_j, \mu_j, u_j),*} \to 0,$$
	we have, up to a subsequence,
	$$\dist_M((\theta_j, \mu_j, u_j), \tilde{K}_b)\to 0.$$
\end{Proposition}
%We note that $(\widetilde{PSP})_b$ condition is different from the standard Palais-Smale condition and it ensures the compactness of $(\theta_n,\lambda_n,u_n)$ after a suitable scaling. 
%By \eqref{eq_confronto_K} we also highlight that, if $\tilde K_b\not=\emptyset$, then $\tilde K_b$ is not compact. % (see \eqref{eq_confronto_K}).
%
%\medskip
%
%\claim Proof.	
%Let $(\theta_n, \lambda_n, u_n)$ as in $(\widetilde{PSP})_b$. Then set $\tilde{u}_n := u_n(e^{-\theta_n} \cdot)$ we have
%$$ \mc{P}(\lambda_n, \tilde{u}_n)\to 0,$$
%$$\partial_{\lambda}\mc{I}^m(\lambda_n, \tilde{u}_n) \to 0,$$
%$$ \norm{\partial_u \mc{I}^m(\lambda_n,\tilde{u}_n)}_{(H^s_r(\R^N))^*}\to 0,$$
%and thus by Proposition \ref{PSP} the sequence $(\lambda_n, \tilde{u}_n)$ is convergent (up to subsequences) to a $(\lambda, \tilde{u})\in K^{PSP}_b$. 
%Observe that, for each $n$, set $v_n : =\tilde{u}(e^{\theta_n}\cdot)$, we have $(\theta_n, \lambda, v_n)\in \tilde{K}_b$. 
%Therefore by \eqref{eq_shift_distance}
%\begin{align*}
%	\dist_M((\theta_n, \lambda_n, u_n), \tilde{K}_b) &\leq \dist_M((\theta_n, \lambda_n, u_n), (\theta_n, \lambda, v_n)) \\
%	&= \dist_M((0,\lambda_n, \tilde{u}_n), (0, \lambda, \tilde{u})) \\
%	&\leq \sqrt{|\lambda_n -\lambda|^2 + \norm{\tilde{u}_n-\tilde{u}}_{H^s_r(\R^N)}^2} \to 0,
%\end{align*}
%which reaches the claim.
%\QED
Through the use of the augmented functional we can obtain again a deformation result. We write here the statement for the unconstrained case (similarly to Proposition \ref{S:5.2}), since it will be used afterwards. %, it takes the following form. 
Set
%$$\mc{J}_{\mu}(u):= \half \int_{\R^N} |(-\Delta)^{s/2} u|^2\,dx -\half {\mathcal D}(u) + \frac{\mu}{2} \|u\|_2^2$$
%and
$$K_b^{PSP}:=\big\{ u \in H^s_r(\R^N) \mid \mc{J}_{\mu}(u)=0, \; \mc{J}'_{\mu}(u)=0, \; \mc{P}_{\mu}(u)=0\big\}.$$
\begin{Lemma}\label{lem_doubl_deform}
For any $b\in \R$, $\bar{\eps}>0$ and any $U$ open neighborhood of $K_b^{PSP}$, there exist an $\epsilon \in (0, \bar \epsilon)$ and a continuous map $\eta:[0,1] \times H^s_r(\R^N) \to H^s_r(\R^N)$ such that 
\begin{itemize}
	\item[$(1^o)$] 
	$
	 \eta(0,u)=u \quad \forall u \in H^s_r(\R^N)
	$;
	\item[$(2^o)$] 
		$
	\eta(t, u)=u \quad \forall (t,u) \in [0,1] \times [\mc{J}_{\mu}\leq b- \bar\epsilon]
	$;
	\item[$(3^o)$] 
	$
	\mc{J}_{\mu}( \eta(t, u)) \leq \mc{J}_{\mu}(u) \quad 
	\ \forall (t,u) \in [0,1] \times H^s_r(\R^N)
	$;
	\item[$(4^o)$] 
	$ \eta (1, [\mc{J}_{\mu} \leq b+ \epsilon] \setminus U) \subset [\mc{J}_{\mu} \leq b- \epsilon]
	$;
				\item[$(5^o)$] 
	$ \eta (1, [\mc{J}_{\mu} \leq b+ \epsilon]) \subset [\mc{J}_{\mu} \leq b - \epsilon] \cup U
	$;
	\item[$(6^o)$] 
	if $K_b^{PSP}= \emptyset$, then
	$ \eta (1, [\mc{J}_{\mu} \leq b+ \epsilon]) \subset [\mc{J}_{\mu} \leq b- \epsilon].
	$
\end{itemize}
\end{Lemma}
The remaining part of the proof follows the lines of the previous Chapters, so that we obtain the existence of a (normalized) Mountain Pass solution: this proves Theorems %\ref{th_INT_exist_unconstrained}, 
\ref{S:1.1_doubly} and \ref{S:1.12_doubly}.

\subsection{Existence of $L^2$-ground states}
\label{sec_ground_state}

In this Section we show (with an approach different from Section \ref{sec_L2_min_prima}) how to obtain the existence of an $L^2$ ground state, by assuming that this energy level is negative and by exploiting Ekeland variational principle together with our Palais-Smale-Pohozaev condition; then we relate this solution to our Mountain Pass solution of Theorem \ref{S:1.1_doubly}. %See also Section \ref{sec_L2_min_prima} for a different approach.

More precisely, for any $m>0$, we introduce the functional ${\mathcal L} : \mc{S}_m\to \R$ defined by 
	\begin{equation}\label{func}
		{\mathcal L}(u):=\half \int_{\R^N} |(-\Delta)^{s/2} u|^2\, dx -\half \mc{D}(u)
	\end{equation}
	on the sphere
	$$
	\mc{S}_m := \big\{u \in H^s_r(\R^N) \mid \|u\|^2_2= m \big\}
	$$
and we consider the $L^2$ ground state level
$$\kappa_m :=\inf_{u \in \mc{S}_m} {\mathcal L}(u).$$
We have the following result.
\begin{Proposition} \label{minimizing}
Under the assumption of Theorem \ref{S:1.1_doubly}, % \ref{S:1.1_doubly}, 
we have, for any $m>m_0$,
\begin{itemize}
\item[\textnormal{(i)}] $-\infty<\kappa_m<0$ and $\kappa_m$ is attained;
\item[\textnormal{(ii)}] $\kappa_m=b_m$, where $b_m$ is defined in \eqref{crti}.
\end{itemize}
Moreover, in the assumptions of Theorem \ref{S:1.12_doubly}, %\ref{S:1.12_doubly}, 
$m_0=0$.
\end{Proposition}

\claim Proof.
We split in some steps. 

\noindent
\textbf{Step 1:} \emph{$\kappa_m>-\infty$.} 
%
%\smallskip
%
%\noindent
\\
By arguing as in Step 2 of Theorem \ref{PSP} we obtain
%
%By (g1)-(g2), for any $\delta>0$ there exists $C_\delta>0$ such that
%	$$	\mc{L}(u) \geq \half\norm{(-\Delta)^{s/2}u}_2^2
%		-\frac\delta{p+1}\norm u_{p+1}^{p+1}
%		-C_\delta \norm u_2^2.
%	$$
%By the fractional Gagliardo-Nirenberg inequality \eqref{...} we have, for $u\in \mc{S}_m$
%	\begin{align*}
%	\mc{L}(u) &\geq \half\norm{(-\Delta)^{s/2}u}_2^2
%		-\frac{C\delta}{p+1}\norm{(-\Delta)^{s/2}u}_2^2
%			\norm u_2^{p-1}
%		-C_\delta \norm u_2^2 \\
%		&=\left(\half-\frac{C\delta}{p+1}m^{\frac{p-1}2}\right)
%			\norm{(-\Delta)^{s/2}u}_2^2 -C_\delta m.
%	\end{align*}
$$ \mc{L}(u) \geq \left(\half-\delta C m^{2(p-1)}\right)
			\norm{(-\Delta)^{s/2}u}_2^2 -C_\delta m^{2 \frac{N+\alpha}{N}}.
$$
Choosing $\delta>0$ small so that $\half-\delta C m^{2(p-1)}>0$, we have $\kappa_m\geq -C_\delta m^{2 \frac{N+\alpha}{N}}>-\infty$.

\smallskip

\noindent
 \textbf{Step 2:} \emph{For $m>m_0$, $\kappa_m<0$.} 
\\
Since the solution $u_*\in\mc{S}_m$ obtained in Theorem \ref{S:1.1_doubly} satisfies, for $m>m_0$, 
	$$	0>b_m=\mc{L}(u_*) \geq \kappa_m,
	$$
we have the claim. %$\kappa_m<0$ for $m>m_0$.

\smallskip

\noindent
 \textbf{Step 3:} \emph{For $m>m_0$, $\kappa_m$ is attained.}
%
%\smallskip
%
%\noindent
\\
To show the existence of a minimizer of $\mc{L}$ on $\mc{S}_m$, we use a linear action $ \Phi:\, \R\to L(H_r^s(\R^N))$ defined by
	$$	\Phi_\theta v := e^{\frac N2 \theta}v(e^\theta \cdot).
	$$
We note that $\mc{S}_m$ is invariant under 
$\Phi_\theta$, that is, $\Phi_\theta(\mc{S}_m) = \mc{S}_m$. 
Let 
$$N:=\R\times \mc{S}_m$$
and on the tangent bundle $TN=\R\times T\mc{S}_m =\coprod_{(\theta,u)\in N} (\R\times T_u\mc{S}_m)$
we introduce a $C^2$-metric
	$$	\norm{(\kappa,v)}_{(\theta,u)} 
	:=\left( \kappa^2+\norm{\Phi_\theta v}_{H^s(\R^N)}^2
		\right)^{1/2}
	$$
for all $(\theta, u) \in N$ and $(\kappa, v) \in TN$. We also introduce $\tilde{\mc{L}}:\, N \to \R$ by
\begin{align*}
\tilde{\mc{L}}(\theta,u)
		&:= \mc{L}(\Phi_\theta u)\\
		&= \half e^{2s\theta}\norm{(-\Delta)^{s/2}u}_2^2
-\half e^{-(N+\alpha)\theta} \mc{D}\big(e^{\frac{N}{2}\theta}u_0\big). % \int_{\R^N} (I_\alpha*F(e^{\frac{N}{2}\theta}u_0))F(e^{\frac{N}{2}\theta}u_0)\,dx.
		%-e^{-N\theta}\int_{\R^N} G(e^{\frac N2\theta} u).
\end{align*}
We note that 
	$$	\inf_{(\theta,u)\in N}\tilde{\mc{L}}(\theta,u)=\kappa_m.
	$$
Since $\kappa_m \in \R$ by Step 1, applying Ekeland's principle, there exists a sequence $(\theta_j,u_j)_{j=1}^\infty\subset N$ such that 
	\begin{eqnarray*}
	&\tilde{\mc{L}}(\theta_j,u_j)\to \kappa_m, \\
	&\norm{D\tilde{\mc{L}}(\theta_j,u_j)}_{T^*_{(\theta_j,u_j)}N}\to 0. 
	\end{eqnarray*}
That is, noticing that $T_u\mc{S}_m \equiv \big\{ v\in H_r^s(\R^N) \mid \int_{\R^N} uv =0\big \}$, we have
%	\begin{eqnarray*}
%	&\partial_\theta\tilde{\mc{L}}(\theta_j,u_j)\to 0,\\
%	&\norm{\partial_u \tilde{\mc{L}}
%		(\theta_j,u_j)}_{T^*_{u_j}\mc{S}_m}
%	=\sup_{\substack{ 
%	v\in T_{u_j}\mc{S}_m \\
%			\norm{\Phi_{\theta_j}v}_{H^s(\R^N)} \leq 1 } }
%	\abs{\partial_u\tilde{\mc{L}}(\theta_j,u_j)v}
%	\to 0.	
%	\end{eqnarray*}
$$ \partial_\theta\tilde{\mc{L}}(\theta_j,u_j)\to 0,$$
$$\norm{\partial_u \tilde{\mc{L}}
		(\theta_j,u_j)}_{T^*_{u_j}\mc{S}_m}
	=\sup_{\substack{ 
	v\in T_{u_j}\mc{S}_m \\
			\norm{\Phi_{\theta_j}v}_{H^s(\R^N)} \leq 1 } }
	\abs{\partial_u\tilde{\mc{L}}(\theta_j,u_j)v}
	\to 0.	$$
Setting $\tilde u_j :=\Phi_{\theta_j}u_j$, we thus have
	\begin{eqnarray}
	&\norm{\tilde u_j}_2^2=m, \label{b.0}\\
	&\mc{L}(\tilde u_j)
	=\half\norm{(-\Delta)^{s/2}\tilde u_j}_2^2 -\frac{1}{2}\mc{D}(\tilde u_j)
%		-\int_{\R^N} G(\tilde u_j)
 \to \kappa_m, \label{b.1}\\
	&s\norm{(-\Delta)^{s/2}\tilde u_j}_2^2
+ \frac{N+\alpha}{2}\mc{D}(\tilde u_j) % (I_\alpha*F(\tilde u_j))F(\tilde u_j)\, dx
		-\frac{N}{2} \int_{\R^N} (I_\alpha*F(\tilde u_j))f(\tilde u_j)\tilde u_j
%		+N\int_{\R^N} G(\tilde u_j)
%		-\frac N2 \int_{\R^N} g(\tilde u_j)\tilde u_j
	\to 0	\label{b.2}
 \end{eqnarray}
and for a suitable $\mu_j\in\R$
 \begin{equation}\label{b.3}
 \mc{L}'(\tilde u_j) \tilde{v}
 +\mu_j\int_{\R^N} \tilde u_j \tilde{v}= o(1)\norm{\tilde{v}}_{H^s(\R^N)}
 \quad\hbox{for all}\ \tilde{v}\in H_r^s(\R^N).
 \end{equation}
By using \eqref{b.1} and arguing as in Step 1 we see that $\tilde{u}_j$ is bounded in $H^s_r(\R^N)$. Thus, choosing $\tilde{v}= \tilde{u}_j$ in \eqref{b.3}, we have
	$$	\norm{(-\Delta)^{s/2}\tilde u_j}_2^2
- \int_{\R^N} (I_\alpha*F(\tilde u_j))f(\tilde u_j)\tilde u_j
	%-\int_{\R^N} g(\tilde u_j)\tilde u_j 
+ \mu_j m
	=o(1),
	$$
which, joined to \eqref{b.2}, gives a Pohozaev identity in the limit
%	$$	\left(1-\frac{2s}N\right)\norm{(-\Delta)^{s/2}\tilde u_j}_2^2
%		-2\int_{\R^N} G(\tilde u_j) +\mu_j m =o(1),
%	$$
\begin{equation}\label{eq_Ekeland_Pohozaev}
	\frac{N-2s}{2} \norm{(-\Delta)^{s/2}\tilde u_j}_2^2
		-\frac{N+\alpha}{2} \mc{D}(\tilde u_j)+ \frac{N}{2} \mu_j m =o(1).
\end{equation}
From this relation and \eqref{b.1} we have
%	$$	-\frac{2s}N \norm{(-\Delta)^{s/2}\tilde u_j}_2^2
%		+ 2\mc{L}(\tilde u_j)
%		+\mu_j m =o(1).
%	$$
%Thus, by \eqref{b.1},
%	$$	\mu_j \geq -\frac{2\kappa_m}m + o(1)
%	$$
$$\mu_j = \frac{2}{Nm} \left( \frac{\alpha+2s}{2} \norm{(-\Delta)^{s/2}\tilde u_j}_2^2 - (N+\alpha) \kappa_m\right) + o(1)$$
which implies, by Step 1, that $\mu_j>0$ for $j$ large.

%Hence we can write $\mu_j=e^{\tilde\lambda_j}$ for some $\tilde\lambda_j\in\R$. 

Relations \eqref{b.0}, \eqref{b.1}, \eqref{b.3} and \eqref{eq_Ekeland_Pohozaev} imply that $(\tilde \mu_j,\tilde u_j)$ is a $(PSP)_{\kappa_m}$ sequence.
% satisfies \eqref{prima}--\eqref{terza4.14} with $b=\kappa_m<0$.
Thanks to the Palais-Smale-Pohozaev condition given in Proposition \ref{PSP2}, $(\tilde \mu_j,\tilde u_j)$ has a strongly convergent subsequence to some $(\tilde\mu_*,\tilde u_*) \in N$, which shows the existence of a minimizer $\tilde{u}_*$. 
Thus (i) is proved.

\smallskip

\noindent
\textbf{Step 4:} \emph{For $m>m_0$, $\kappa_m=b_m$.} 
%
%\smallskip
%
%\noindent
\\
In Step 2 we showed $b_m\geq \kappa_m$. On the other hand by the argument in Step 3, for the minimizer $\tilde u_*$ of $\mc{L}$ on $\mc{S}_m$, there exists $\tilde \mu_*\in\R$ such that
	\begin{eqnarray*}
	&\mc{I}^m(\tilde\mu_*,\tilde u_*)=\kappa_m,
	\quad 
	\partial_u\mc{I}^m(\tilde\mu_*,\tilde u_*)=0,
		\\
	&\partial_{\mu}\mc{I}^m(\tilde\mu_*,\tilde u_*)=0,
	\quad 
	\mc{P}(\tilde\mu_*,\tilde u_*)=0.
	\end{eqnarray*}
Set $\xi_*(t):= \tilde{u}_*(\cdot/t)$ we have, by the Pohozaev identity, $\mc{I}^m(\tilde{\mu}_*, \xi_*(t))\to -\infty$ as $t\to +\infty$; thus, up to a rescaling, we obtain $\xi_*\in\Gamma^m$ and 
	$$	\max_{t\in [0,1]}\mc{I}^m(\xi_*(t)) %\leq 
=		 \mc{I}^m(\tilde{\mu}_*,\tilde{u}_*)
		= \kappa_m,
	$$
%by arguing as in Proposition \ref{tom}
% \tr{VEDI meglio} %COMMENT NOW
%%the argument in previous Sections, 
%there exists $\xi_*\in\Gamma^m$ such that
%	$$	\max_{t\in [0,1]}\mc{I}^m(\xi_*(t)) %\leq 
%=
%		 \mc{I}^m(\tilde \lambda_*,\tilde u_*)
%		= \kappa_m,
%	$$
which implies $b_m \leq \kappa_m$ and the proof %of Proposition \ref{minimizing} 
is completed. 
%
%Set $\zeta_\lambda(t):=u_\lambda(\cdot/t)$ for $t>0$ and $\zeta_\lambda(0):=0$ and note that, since $u_{\lambda}$ satisfies the Pohozaev identity, we have $\mc{I}^m (\lambda,\zeta_\lambda(t)) \to-\infty$ and $\mc{P}(\lambda,\zeta_\lambda(t)) \to-\infty$ as $t\to+\infty$. 
%We can find $\gamma_\lambda:=\zeta_\lambda(L\cdot)$ for $L\gg 1$ satisfying
\QED

%%%%%%%%%%%%%%%%%%%%%%%%%%%%%%%%%%%%%%%%%%%%%%%%%%%%%%%

\section{Preliminary properties of Pohozaev energy levels}
\label{sec_proper_pohozaev_gs}

As highlighted, the goal of this Chapter is to study qualitative properties of solutions and, in particular, of Pohozaev minima. In this Section, thus, we start by observing that the solution found in Theorem \ref{th_INT_exist_unconstrained} is actually a Pohozaev minimum. % ground state.
Since, afterwards, we will be interested in studying symmetric properties of general ground states, in this Section we highlight the dependence of some sets and energy levels from the subspace of radially symmetric functions. Moreover, we show existence of positive solutions.

%\smallskip

\subsubsection{Energy levels in radially symmetric spaces}
%We define the functional $\mc{J}_{\mu}: H^s(\R^N) \to \R$ by setting
%$$\mc{J}_{\mu}(u):= \half \int_{\R^N} |(-\Delta)^{s/2} u|^2\,dx -\half {\mathcal D}(u) + \frac{\mu}{2} \|u\|_2^2.$$
%where
%$$\mc{D}(u) := \int_{\R^N} (I_\alpha*F(u))F(u)\, dx.$$
We introduce the set of paths
$$\Gamma_r(\mu) := \big\{ \gamma \in C\big([0,1], H^s_r(\R^N)\big) \mid \gamma(0)=0, \, \mc{J}_{\mu}(\gamma(1))<0\big\}$$
and the \emph{Mountain Pass} (MP for short) value 
\begin{equation}\label{eq_MP_value}
 a_r(\mu) := \inf_{\gamma\in\Gamma_\mu}\max_{t\in [0,1]} \mc{J}_{\mu}(\gamma(t)).
 \end{equation}
Then we introduce
$$p_r(\mu):= \inf \big\{ \mc{J}_{\mu}(u) \mid u \in H^s_r(\R^N) \setminus \{0\}, \; \mc{P}_{\mu}(u)=0\big\}$$
the \emph{least energy} of $\mc{J}_{\mu}$ on the Pohozaev set of radially symmetric functions.

%We prove now that the found solution is actually a ground state over the Pohozaev set.

\begin{Proposition}\label{prop_MP=PM}
The Mountain Pass level and the Pohozaev minimum level coincide, that is
$$a_r(\mu)=p_r(\mu)>0.$$
In particular, the solution found in Theorem \ref{th_INT_exist_unconstrained} is a Pohozaev minimum.
\end{Proposition}

\claim Proof.
Let $u \in H^s_r(\R^N)\setminus\{0\}$ such that $\mc{P}_{\mu}(u)=0$; observe that $\mc{D}(u)>0$. We define $\bar{\gamma}(t):=u(\cdot/t)$ for $t \neq 0$ and $\bar{\gamma}(0):=0$ so that $t\in (0,+\infty) \mapsto \mc{J}_{\mu}(\bar{\gamma}(t))$ is negative for large values of $t$, and it attains the maximum in $t=1$. After a suitable rescaling we have $\bar{\gamma} \in \Gamma_r(\mu)$ and thus
\begin{equation}\label{eq_dis_P-MP}
\mc{J}_{\mu}(u) = \max_{t \in [0,1]}\mc{J}_{\mu}(\bar{\gamma}(t)) \geq a_r(\mu).
\end{equation}
Passing to the infimum in \eqref{eq_dis_P-MP} we have $p_r(\mu) \geq a_r(\mu)$. Let now $\gamma \in \Gamma_r(\mu)$. By definition we have $\mc{J}_{\mu}(\gamma(1))<0$, thus by
$$ \mc{P}_{\mu}(v) = N \mc{J}_{\mu}(v) - s\norm{(-\Delta)^{s/2} v}_2^2 - \frac{\alpha}{2} \mc{D}(v), \quad v \in H^s_r(\R^N),$$
we obtain $\mc{P}_{\mu}(\gamma(1))<0$. 
In addition, since $\mc{D}(u)=o(\norm u_{H^s}^2)$ as $u \to 0$ and $\gamma(t)\to 0$ as $t \to 0$ in $H^s_r(\R^N)$, we have
$$\mc{P}_{\mu}(\gamma(t))>0 \quad \hbox{for small $t>0$}.$$
Thus there exists a $t^*$ such that $\mc{P}_{\mu}(\gamma(t^*))=0$, and hence
$$p_r(\mu) \leq \mc{J}_{\mu}(\gamma(t^*)) \leq \max_{t \in [0,1]} \mc{J}_{\mu}(\gamma(t));$$
passing to the infimum we come up with $p_r(\mu) \leq a_r(\mu)$, and hence the claim.
\QED

\bigskip

We pass to investigate more in details Pohozaev minima, showing that it is a general fact that they are solutions of equation \eqref{eq_introduction}.

\begin{Proposition}\label{prop_poho_min_sol}
Every Pohozaev minimum is a solution of \eqref{eq_introduction}, i.e.
$$\mc{J}_{\mu}(u)=p_r(\mu) \; \hbox{ and } \; \mc{P}_{\mu}(u)=0$$
imply
$$\mc{J}'_{\mu}(u)=0.$$
As a consequence
$$p_r(\mu)= \inf \big\{ \mc{J}_{\mu}(u) \mid u \in H^s_r(\R^N) \setminus \{0\}, \; \mc{P}_{\mu}(u)=0, \; \mc{J}'_{\mu}(u)=0\big\}.$$
\end{Proposition}

\claim Proof.
Let $u$ be such that $\mc{J}_{\mu}(u)=p_r(\mu)$ and $\mc{P}_{\mu}(u)=0$. 
In particular, considered $\gamma(t):=u(\cdot/t)$, we have that $\mc{J}_{\mu}(\gamma(t))$ is negative for large values of $t$ and its maximum value is $p(\mu)$ attained only in $t=1$. 

Assume by contradiction that $u$ is not critical. Let $I:=[1-\delta, 1+\delta]$ be such that $\gamma(I) \cap K_{p(\mu)}=\emptyset$, and set $\bar{\eps} := p(\mu) - \max_{t \notin I} \mc{J}_{\mu}(\gamma(t)) >0$.
Let now $U$ be a neighborhood of $K_{p(\mu)}$ verifying $\gamma(I) \cap U = \emptyset$: by the Deformation Lemma \ref{lem_doubl_deform}
there exists an $\eta:[0,1]\times H^s_r(\R^N) \to H^s_r(\R^N)$ at level $p_r(\mu)\in\R$ with properties $(1^o)$-$(6^o)$.
 Define then $\tilde{\gamma}(t):= \eta(1, \gamma(t))$ a deformed path. 
 
 For $t \notin I$ we have $\mc{J}_{\mu}(\gamma(t))< p_r(\mu) - \bar{\eps}$, and thus by $(2^o)$ we gain
\begin{equation}\label{eq_dim_P_sol_1}
\mc{J}_{\mu}(\tilde{\gamma}(t)) = \mc{J}_{\mu}(\gamma(t)) < p_r(\mu) - \bar{\eps}, \quad \hbox{ for $t \notin I$}.
\end{equation}
Let now $t \in I$: we have $\gamma(t) \notin U$ and $\mc{J}_{\mu}(\gamma(t))\leq p_r(\mu) \leq p_r(\mu)+ \eps$, thus by $(4^o)$ we obtain
\begin{equation}\label{eq_dim_P_sol_2}
\mc{J}_{\mu}(\tilde{\gamma}(t)) \leq p_r(\mu) - \eps.
\end{equation}
 Joining \eqref{eq_dim_P_sol_1} and \eqref{eq_dim_P_sol_2} we have
$$\max_{t \geq 0} \mc{J}_{\mu}(\tilde{\gamma}(t)) < p_r(\mu)=a_r(\mu)$$
which is an absurd, since after a suitable rescaling it results that $\tilde{\gamma} \in \Gamma_r(\mu)$, thanks to $(3^o)$. 
\QED

\begin{Remark}\label{rem_Pohozaev_c_p}
We point out that it is not known, even in the case of local nonlinearities \cite{BKS}, if 
$$p_r(\mu) \stackrel{?} = \inf \big\{ \mc{J}_{\mu}(u) \mid u \in H^s_r(\R^N) \setminus \{0\}, \; \mc{J}_{\mu}'(u)=0\big\}.$$
On the other hand, by assuming that every solution of \eqref{eq_introduction} satisfies the Pohozaev identity (see e.g. \cite[Proposition 2]{SGY} and \cite[Eq (6.1)]{DSS1} and Section \ref{sec_Pohozaev}), %and also the paper in preparation \cite{CGT5}), 
the claim holds true. 
We point out that the equality may hold even if it is not true that every solution satisfies the Pohozaev identity. 
The fact that Deformation Lemma \ref{lem_doubl_deform} % we proved (see Theorem \ref{defarg}) 
allows to deform the functional near critical points not satisfying the Pohozaev identity might be useful in the investigation of these facts.
\end{Remark}

\subsubsection{Energy levels in the whole space}

We pass studying general Pohozaev minima on the whole space $H^s(\R^N)$. 
We start defining the least energy of $\mc{J}_{\mu}$ on the Pohozaev set, and call every minimizer a \emph{Pohozaev minimum} (or \emph{ground state})
\begin{equation}\label{eq_poh_minim}
p(\mu):= \inf \big\{ \mc{J}_{\mu}(u) \mid u \in H^s(\R^N) \setminus \{0\}, \; \mc{P}_{\mu}(u)=0\big\}.
\end{equation}

We start by showing that Proposition \ref{prop_poho_min_sol} holds also in a nonradial setting, providing here the proof.
To do this, we get advantage of the minimax paths and level of $\mc{J}_{\mu}$. Set
$$a(\mu):= \inf_{\gamma \in \Gamma_{\mu}} \sup_{t \in [0,1]} \mc{J}_{\mu}(\gamma(t))$$
where
$$\Gamma(\mu):=\big\{ \gamma \in C\big([0,1], H^s(\R^N)\big) \mid \gamma(0)=0, \mid \mc{J}_{\mu}(\gamma(1))<0\big\}.$$
Notice that, with the same proof of Proposition \ref{prop_MP=PM} we obtain
\begin{equation}\label{eq_MP=PM}
a(\mu)=p(\mu)>0.
\end{equation}

%\begin{Lemma}\label{lem_MP=PM}
%The Mountain Pass level and the Pohozaev minimum level of $\mc{J}_{\mu}$ coincide, that is
%$$b(\mu)=p(\mu)>0.$$
%\end{Lemma}
%
%\begin{proof}
%Let $u \in H^s(\R^N)\setminus\{0\}$ such that $\mc{P}_{\mu}(u)=0$; observe that $\mc{D}(u)>0$.
%We define $\gamma(t):=u(\cdot/t)$ for $t \neq 0$ and $\gamma(0):=0$ so that $t\in (0,+\infty) \mapsto \mc{J}_{\mu}(\gamma(t))$ is negative for large values of $t$, and it attains the maximum in $t=1$.
%After a suitable rescaling we have $\gamma \in \Gamma_{\mu}$ and thus
%\begin{equation}\label{eq_dis_P-MP}
%\mc{J}_{\mu}(u) = \max_{t \in [0,1]}\mc{J}_{\mu}(\gamma(t)) \geq b(\mu).
%\end{equation}
%Passing to the infimum in equation \eqref{eq_dis_P-MP} we have $p(\mu) \geq b(\mu)$. Let now $\gamma \in \Gamma_{\mu}$. By definition we have $\mc{J}_{\mu}(\gamma(1))<0$, thus by
%$$ \mc{P}_{\mu}(v) = N \mc{J}_{\mu}(v) - s\norm{(-\Delta)^{s/2} v}_2^2 - \frac{\alpha}{2} \mc{D}(v), \quad v \in H^s(\R^N),$$
%we obtain $\mc{P}_{\mu}(\gamma(1))<0$.
%In addition, since $\mc{D}(u)=o(\norm u_{H^s}^2)$ as $u \to 0$ and $\gamma(t)\to 0$ as $t \to 0$ in $H^s(\R^N)$, we have
%$$\mc{P}_{\mu}(\gamma(t))>0 \quad \hbox{for small $t>0$}.$$
%Thus there exists a $t^*$ such that $\mc{P}_{\mu}(\gamma(t^*))=0$, and hence
%$$p(\mu) \leq \mc{J}_{\mu}(\gamma(t^*)) \leq \max_{t \in [0,1]} \mc{J}_{\mu}(\gamma(t));$$
%passing to the infimum we come up with $p(\mu) \leq b(\mu)$, and hence the claim.
%\end{proof}

\begin{Proposition}\label{prop_poho_min_Rsol}
Assume \hyperref[(F1s)]{\textnormal{(F1)}}-\hyperref[(F2s)]{\textnormal{(F2)}}. Then every Pohozaev minimum of $\mc{J}_{\mu}$ is a solution of \eqref{eq_introduction}, i.e.
$$\mc{J}_{\mu}(u)=p(\mu) \; \hbox{ and } \; \mc{P}_{\mu}(u)=0$$
imply
$$\mc{J}'_{\mu}(u)=0.$$
As a consequence
$$p(\mu)= \inf \big\{ \mc{J}_{\mu}(u) \mid u \in H^s(\R^N) \setminus \{0\}, \; \mc{P}_{\mu}(u)=0, \; \mc{J}'_{\mu}(u)=0\big\}.$$
\end{Proposition}

\claim Proof.
Assume by contradiction that $\mc{J}_{\mu}'(u)\neq 0$. Thus $\mc{J}_{\mu}'$ remains far from zero in a neighborhood of $u$, that is there exist $\delta>0$ and $\lambda>0$ such that
$$v \in B_{3 \delta}(u) \implies \norm{\mc{J}'_{\mu}(v)}_* \geq \lambda.$$
Consider the path $\gamma(t):=u(\cdot/t)$; it is straightforward to show that $t \in \R_+ \mapsto \mc{J}_{\mu}(\gamma(t))$ is negative for $t\gg0$ and it has a unique strict maximum, equal to $\mc{J}_{\mu}(u)=p(\mu)>0$, attained in $t=1$.
Let now $I:=[1-\omega, 1+\omega]$, $\omega$ small, be such that
$$S:=\gamma(I)\subset B_{\delta}(u);$$
we can also assume that
$$\max_{t\notin I}\mc{J}_{\mu}(\gamma(t))\in \big(0, p(\mu)\big).$$
Introduce moreover
$$0<\eps< \min\left\{ \frac{p(\mu)-\max_{t\notin I}\mc{J}_{\mu}(\gamma(t))}{2}, \, \frac{\lambda \delta}{8} \right \}.$$
By writing $S_{2\delta}:= \{ v \in H^s(\R^N) \mid d(v, S)\leq 2 \delta\}$, we see that
$$v \in S_{2 \delta} \implies \norm{\mc{J}'_{\mu}(v)}_* \geq\frac{8 \eps}{\delta}$$
and in particular
$$v \in \mc{J}_{\mu}^{-1}\big([p(\mu)-2\eps, p(\mu)+2\eps]\big) \cap S_{2 \delta} \implies \norm{\mc{J}'_{\mu}(v)}_* \geq\frac{8 \eps}{\delta},$$
where we observe that $p(\mu)-2\eps>0$.
We are thus in the assumptions of \cite[Lemma 2.3]{Wil0}, and we have the existence of a \emph{local} continuous deformation $\eta : [0,1]\times H^s(\R^N) \to H^s(\R^N)$ such that (we write $\mc{J}_{\mu}^b:= \mc{J}_{\mu}^{-1}\big((-\infty, b]\big)$)
\begin{itemize}
\item[(a)] $\eta(0,v)=v$,
\item[(b)] $\eta(t,v)=v$ if $v \notin \mc{J}_{\mu}^{-1}\big([p(\mu)-2\eps, p(\mu)+2\eps]\big) \cap S_{2 \delta}$,
\item[(c)] $\mc{J}_{\mu}(\eta(\cdot, v))$ is non increasing for each $v\in H^s(\R^N)$,
\item[(d)] $\eta(1,\mc{J}_{\mu}^{p(\mu)+\eps}\cap S)\subset \mc{J}_{\mu}^{p(\mu)-\eps}$.
\end{itemize}
We thus define a deformed path
$$\tilde{\gamma}(t):=\eta(1, \gamma(t)).$$
Consider first $t\notin I$. By (c) and the definition of $\eps$, we have
$$\mc{J}_{\mu}(\tilde{\gamma}(t))\leq \mc{J}_{\mu}(\gamma(t)) < p(\mu) - 2 \eps<p(\mu).$$
Assume instead $t\in I$. Then $\gamma(t) \in \gamma(I)=S$ and $\mc{J}_{\mu}(\gamma(t)) \leq \mc{J}_{\mu}(\gamma(1))= p(\mu) \leq p(\mu)+\eps$, thus by (d) we have
$$\mc{J}_{\mu}(\tilde{\gamma}(t)) \leq p(\mu)- \eps < p(\mu).$$
Joining together the two inequalities we obtain
\begin{equation}\label{eq_dim_less_pmu}
\max_{t >0
} \mc{J}_{\mu}(\tilde{\gamma}(t)) < p(\mu).
\end{equation}
On the other hand, we have $\tilde{\gamma}(0)=\eta(1,\gamma(0))=\eta(1,0)=0$ since $0 \notin \mc{J}_{\mu}^{-1}\big([p(\mu)-2\eps, p(\mu)+2\eps]\big)$, and $\mc{J}_{\mu}(\tilde{\gamma}(t)) \leq \mc{J}_{\mu}(\eta(1, \gamma(t))\leq \mc{J}_{\mu}(\gamma(t))<0$ for $t \gg 0$.
Up to a rescaling, we can assume
$\tilde{\gamma} \in \Gamma(\mu)$ and hence, by \eqref{eq_MP=PM} %Lemma \ref{lem_MP=PM}
\begin{equation*}\label{eq_dim_great_pmu}
\max_{t \in [0,1]} \mc{J}_{\mu}(\tilde{\gamma}(t))\geq a(\mu) = p(\mu),
\end{equation*}
which is in contradiction with \eqref{eq_dim_less_pmu}. The proof is thus concluded.
\QED

\bigskip

\begin{Remark}\label{rem_Pohozaev_c_p_2}
As in Remark \ref{rem_Pohozaev_c_p}, we point out that it is not known, even in the case of local nonlinearities, 
if 
$$p(\mu) \stackrel{?} = \inf \big\{ \mc{J}_{\mu}(u) \mid u \in H^s(\R^N) \setminus \{0\}, \; \mc{J}_{\mu}'(u)=0\big\},$$
unless some additional assumptions on $s$ or $f$ are assumed.
\end{Remark}

In Corollary \ref{cor_uguagl_p} we will state some relation between $p_r(\mu)$ and $p(\mu)$.

\medskip

Most of the qualitative properties that we will investigate, will be stated in the case of positive solutions. Thus it is important to highlight the existence of a solution of constant sign.
% We show now that, under the same assumptions of Theorem \ref{th_exist_unconstrained}, we can find a solution with constant sign.

\begin{Proposition}\label{prop_ex_posiv_solut}
Assume \hyperref[(F1s)]{\textnormal{(F1)}}--\hyperref[(F4s)]{\textnormal{(F4)}} and that $F\nequiv 0$ on $(0,+\infty)$ (i.e., $t_0$ in assumption \hyperref[(F4s)]{\textnormal{(F4)}} can be chosen positive). Then there exists a positive radially symmetric solution of \eqref{eq_introduction}, which is minimum over all the positive functions on the Pohozaev set.
\end{Proposition}

\claim Proof.
Let us define
$$\tilde{f}:= \chi_{(0, +\infty)} f.$$
We have that $\tilde{f}$ still satisfies \hyperref[(F1s)]{\textnormal{(F1)}}--\hyperref[(F4s)]{\textnormal{(F4)}}. Thus, by Theorem \ref{th_INT_exist_unconstrained} there exists a solution $u$ of
$$(- \Delta)^s u + \mu u = (I_\alpha*\tilde{F}(u))\tilde{f}(u) \quad \hbox{in $\mathbb{R}^N$}$$
where $\tilde{F}(t):=\int_0^t \tilde{f}(\tau) d \tau$, $\tilde{F} = \chi_{(0,+\infty)} F$. 
We show now that $u$ is positive. 
%We start observing the following: by \eqref{eq_semin_gagl} we have
%\begin{eqnarray*}
%\norm{(-\Delta)^{s/2} |u|}_2^2 &=& C_{N,s} \int_{\R^{2N}} \frac{\big(|u(x)|-|u(y)|\big)^2}{|x-y|^{N+2s}} \, dx \, dy \\
%&=& C_{N,s} \int_{\R^{2N}} \frac{|u|^2(x) + |u|^2(y) - 2|u|(x)|u|(y)}{|x-y|^{N+2s}} \, dx \, dy \\
%&\leq& C_{N,s} \int_{\R^{2N}} \frac{u^2(x) + u^2(y) - 2u(x)u(y)}{|x-y|^{N+2s}} \, dx \, dy \\
%&=& C_{N,s} \int_{\R^{2N}} \frac{\big(u(x)-u(y)\big)^2}{|x-y|^{N+2s}} \, dx \, dy = \norm{(-\Delta)^{s/2}u}_2^2,
%\end{eqnarray*}
%thus
%$$\norm{(-\Delta)^{s/2} |u|}_2 \leq \norm{(-\Delta)^{s/2} u}_2.$$
%In particular, written $u=u_+ - u_-$, by the previous argument we have 
Recall by Lemma \ref{lem_dis_modul_1} that $u_- = \frac{|u|-u}{2} \in H^s_r(\R^N)$. Thus, chosen $u_-$ as test function, we obtain
$$\int_{\R^N} (-\Delta)^{s/2} u \, (-\Delta)^{s/2}u_- \, dx + \mu \int_{\R^N} u \, u_- \, dx = \int_{\R^N} (I_{\alpha} * \tilde{F}(u))\tilde{f}(u) u_- \, dx.$$
By definition of $\tilde{f}$ and \eqref{eq_semin_gagl} we have
\begin{equation}\label{eq_splitting_posit}
C_{N,s}\int_{\R^N \times \R^N} \frac{(u(x)-u(y))(u_-(x)-u_-(y))}{|x-y|^{N+2s}} \, dx \, dy - \mu \int_{\R^N} u_-^2 \, dx =0.
\end{equation}
Splitting the domain, we gain
\begin{eqnarray*}
\lefteqn{\int_{\R^N \times \R^N} \frac{(u(x)-u(y))(u_-(x)-u_-(y))}{|x-y|^{N+2s}} \, dx \, dy =}\\
&&- \int_{\{u(x)\geq 0\} \times \{u(y)<0\}} \frac{(u_+(x)+u_-(y))(u_-(y))}{|x-y|^{N+2s}} \, dx \, dy -\\
&& -\int_{\{u(x)<0\} \times \{u(y)\geq 0\}} \frac{(u_-(x)+u_+(y))(u_-(x))}{|x-y|^{N+2s}} \, dx \, dy -\\
&&- \int_{\{u(x)<0\} \times \{u(y)<0\}} \frac{(u_-(x)-u_-(y))^2}{|x-y|^{N+2s}} \, dx \, dy.
\end{eqnarray*}
Since the left-hand side of \eqref{eq_splitting_posit} is sum of nonpositive pieces, we have $u_- \equiv 0$, that is $u\geq 0$. Hence $\tilde{f}(u)=f(u)$ and $\tilde{F}(u)=F(u)$, which imply that $u$ is a positive solution of \eqref{eq_introduction}.
\QED

%%%%%%%%%%%%%%%%%%%%%%%%%%%%%%%%%%%%%%%%%%%%%%%%%%%%%%%
%%%%%%%%%%%%%%%%%%%%%%%%%%%%%%%%%%%%%%%%%%%%%%%%%%%%%%%

\section{Regularity}
\label{sec_doubl_regul}

In this Section we investigate regularity of solutions, focusing in particular on boundedness, H\"older regularity and $L^1$-summability.

The discussed results generalize some of the ones in \cite{DSS1} to the case of general, not homogeneous, nonlinearities; in particular, we do not even assume $f$ to satisfy Ambrosetti-Rabinowitz type conditions nor monotonicity conditions. 
Moreover, we improve the results in \cite{SGY, Luo0} since we do not assume $f$ to be superlinear, and we have no restriction on the parameter $\alpha$.

Some of these results extend the ones in \cite{MS2, BKS} to the fractional, Choquard framework.

%%%%%%%%%%%%%%%%%%%%%%%%%%%%%%%%%%
\subsection{Boundedness by splitting}
\label{sec_doubly_boundedness}

%In this Section we prove some regularity results for \eqref{eq_introduction}. 
%We split the proof of Theorem \ref{th_INT_regular} in different steps.

Here we prove that solutions of \eqref{eq_introduction} are bounded. 
%We split the proof of Theorem \ref{th_INT_regular} in different steps.
In particular, when dealing with sign-changing solutions, we will consider also the following stronger assumption: 
\begin{itemize} %\label{eq_cond_superl}
%\item[(f5)] $\limsup_{t \to 0} \frac{|tf(t)|}{|t|^{2}} <+\infty$,
\item[(F6)] \label{(F6)}
%$\limsup_{t \to 0} \frac{|tf(t)|}{|t|^{2}} <+\infty$, 
$\limsup_{t \to 0} \frac{|f(t)|}{|t|} <+\infty$, 
\end{itemize}
which says that $f$ is linear or superlinear in the origin. 
Observe that 
$$\hbox{\hyperref[(F6)]{\textnormal{(F6)}} $\implies$ \hyperref[(F2s)]{\textnormal{(F2,i)}} and \hyperref[(F3s)]{\textnormal{(F3,i)}}.}$$

\begin{Theorem}\label{th_INT_regular}
Assume \hyperref[(F1s)]{\textnormal{(F1)}}-\hyperref[(F2s)]{\textnormal{(F2)}}. Let $u\in H^s(\R^N)$ be a weak positive solution of \eqref{eq_introduction}.
Then $u \in %L^1(\R^N) \cap 
L^{\infty}(\R^N)$. 
The same conclusion holds for generally (possibly sign-changing) %signed 
solutions by assuming also \hyperref[(F6)]{\textnormal{(F6)}}.
\end{Theorem}

We start from the following lemma, that can be found in \cite[Lemma 3.3]{MS2}.
\begin{Lemma}[\cite{MS2}]\label{lem_primoMS}
Let $N \geq 2$ and $\alpha \in (0, N)$. Let $\lambda \in [0,2]$ and $q,r, h,k \in [1, +\infty)$ be such that
$$1+ \frac{\alpha}{N} - \frac{1}{h} - \frac{1}{k} = \frac{\lambda}{q} + \frac{2-\lambda}{r}.$$
Let $\theta \in (0,2)$ satisfying
$$\min\{q,r\} \left(\frac{\alpha}{N} - \frac{1}{h}\right) < \theta < \max\{q,r\} \left( 1- \frac{1}{h}\right),$$
$$\min\{q,r\} \left(\frac{\alpha}{N} - \frac{1}{k}\right) < 2-\theta < \max\{q,r\} \left( 1- \frac{1}{k}\right).$$
Let $H \in L^h(\R^N)$, $K \in L^k(\R^N)$ and $u \in L^q(\R^N) \cap L^r(\R^N)$. Then
$$\int_{\R^N} \left( I_{\alpha} * \big(H|u|^{\theta}\big) \right) K|u|^{2-\theta} \, dx \leq C \norm{H}_h \norm{K}_k \norm{u}_q^{\lambda} \norm{u}_r^{2-\lambda}$$
for some $C>0$ (depending on $\theta$). 
\end{Lemma}

By a proper use of Lemma \ref{lem_primoMS} we obtain now an estimate on the Choquard term depending on $H^s$-norm of the function.

\begin{Lemma}\label{lem_secondMS}
Let $N\geq 2$, $s \in (0,1)$ and $\alpha \in (0,N)$. 
Let moreover $\theta \in (\frac{\alpha}{N}, 2- \frac{\alpha}{N})$ and $H, K \in L^{\frac{2N}{\alpha}}(\R^N) + L^{\frac{2N}{\alpha+2s}}(\R^N) $. Then for every $\eps>0$ there exists $C_{\eps, \theta}>0$ such that
$$\int_{\R^N} \left(I_{\alpha}*\big(H|u|^{\theta}\big)\right)K|u|^{2-\theta} \, dx \leq \eps^2 \norm{(-\Delta)^{s/2} u}_2^2 + C_{\eps, \theta}\norm{u}_2^2$$
for every $u \in H^s(\R^N)$.
\end{Lemma}

\claim Proof.
Observe that $2-\theta \in (\frac{\alpha}{N}, 2- \frac{\alpha}{N})$ as well. We write
$$H= H^* + H_* \in L^{\frac{2N}{\alpha}}(\R^N) + L^{\frac{2N}{\alpha+2s}}(\R^N),$$
$$K= K^* + K_* \in L^{\frac{2N}{\alpha}}(\R^N) + L^{\frac{2N}{\alpha+2s}}(\R^N).$$
We split $\int_{\R^N} \left(I_{\alpha}*\big(H|u|^{\theta}\big)\right)K|u|^{2-\theta} \, dx $ in four pieces and choose
$$q=r=2, \quad h=k=\frac{2N}{\alpha}, \quad \lambda=2,$$
$$q = 2, \; r=\frac{2N}{N-2s}, \quad h=\frac{2N}{\alpha}, \; k=\frac{2N}{\alpha+2s}, \quad \lambda=1,$$
$$q = 2, \; r=\frac{2N}{N-2s}, \quad h=\frac{2N}{\alpha+2s}, \; k=\frac{2N}{\alpha}, \quad \lambda=1,$$
$$q=r=\frac{2N}{N-2s}, \quad h=k=\frac{2N}{\alpha+2s}, \quad \lambda=0,$$
in Lemma \ref{lem_primoMS}, to obtain
%\begin{align*}
%\int_{\R^N} \left(I_{\alpha}*\big(H|u|^{\theta}\big)\right)K|u|^{2-\theta} \, dx \lesssim& \norm{H^*}_{\frac{2N}{\alpha}} \norm{K^*}_{\frac{2N}{\alpha}} \norm{u}_2^2 + \norm{H^*}_{\frac{2N}{\alpha}}\norm{K_*}_{\frac{2N}{\alpha+2s}} \norm{u}_2 \norm{u}_{\frac{2N}{N-2s}}+ \\
%&+ \norm{H_*}_{\frac{2N}{\alpha+2s}} \norm{K^*}_{\frac{2N}{\alpha}}\norm{u}_2 \norm{u}_{\frac{2N}{N-2s}} + \norm{H_*}_{\frac{2N}{\alpha+2s}}\norm{K_*}_{\frac{2N}{\alpha+2s}} \norm{u}_{\frac{2N}{N-2s}}^2.
%\end{align*}
\begin{eqnarray*}
\lefteqn{\int_{\R^N} \left(I_{\alpha}*\big(H|u|^{\theta}\big)\right)K|u|^{2-\theta} \, dx \lesssim } \\
&& \norm{H^*}_{\frac{2N}{\alpha}} \norm{K^*}_{\frac{2N}{\alpha}} \norm{u}_2^2 + \norm{H^*}_{\frac{2N}{\alpha}}\norm{K_*}_{\frac{2N}{\alpha+2s}} \norm{u}_2 \norm{u}_{\frac{2N}{N-2s}}+ \\
&&+ \norm{H_*}_{\frac{2N}{\alpha+2s}} \norm{K^*}_{\frac{2N}{\alpha}}\norm{u}_2 \norm{u}_{\frac{2N}{N-2s}} + \norm{H_*}_{\frac{2N}{\alpha+2s}}\norm{K_*}_{\frac{2N}{\alpha+2s}} \norm{u}_{\frac{2N}{N-2s}}^2.
\end{eqnarray*}
Recalled that $\frac{2N}{N-2s}=2^*_s$ and the Sobolev embedding \eqref{eq_embd_homog}, we obtain
%\begin{align}
%\int_{\R^N} \left(I_{\alpha}*\big(H|u|^{\theta}\big)\right)K|u|^{2-\theta} \, dx \lesssim & \left(\norm{H^*}_{\frac{2N}{\alpha}}\norm{K^*}_{\frac{2N}{\alpha}}\right) \norm{u}_2^2 + \left(\norm{H_*}_{\frac{2N}{\alpha+2s}} \norm{K_*}_{\frac{2N}{\alpha+2s}}\right) \norm{(-\Delta)^{s/2} u}_2^2 + \notag \\
%& + \left(\norm{H^*}_{\frac{2N}{\alpha}}\norm{K_*}_{\frac{2N}{\alpha+2s}} + \norm{H_*}_{\frac{2N}{\alpha+2s}} \norm{K^*}_{\frac{2N}{\alpha}}\right) \norm{u}_2 \norm{(-\Delta)^{s/2} u}_2.\label{eq_dim_mixed_terms_HK}
%\end{align}
\begin{eqnarray}
\lefteqn{\int_{\R^N} \left(I_{\alpha}*\big(H|u|^{\theta}\big)\right)K|u|^{2-\theta} \, dx \lesssim} \notag \\
&& \left(\norm{H^*}_{\frac{2N}{\alpha}}\norm{K^*}_{\frac{2N}{\alpha}}\right) \norm{u}_2^2 + \left(\norm{H_*}_{\frac{2N}{\alpha+2s}} \norm{K_*}_{\frac{2N}{\alpha+2s}}\right) \norm{(-\Delta)^{s/2} u}_2^2 + \notag \\
&& + \left(\norm{H^*}_{\frac{2N}{\alpha}}\norm{K_*}_{\frac{2N}{\alpha+2s}} + \norm{H_*}_{\frac{2N}{\alpha+2s}} \norm{K^*}_{\frac{2N}{\alpha}}\right) \norm{u}_2 \norm{(-\Delta)^{s/2} u}_2.\label{eq_dim_mixed_terms_HK}
\end{eqnarray}
%where $\lesssim$ denotes an inequality up to a constant.
We want to show now that, since $ \frac{2N}{\alpha}>\frac{2N}{\alpha+2s}$, we can choose the decomposition of $H$ and $K$ such that the $L^{\frac{2N}{\alpha+2s}}$-pieces are arbitrary small (see \cite[Lemma 2.1]{BK0}). Indeed, let
$$H = H_1 + H_2 \in L^{\frac{2N}{\alpha}}(\R^N) + L^{\frac{2N}{\alpha+2s}}(\R^N)$$
be a first decomposition. Let $M>0$ to be fixed, and write
$$H = \left( H_1 + H_2 \chi_{\{|H_2| \leq M\}}\right) + H_2 \chi_{\{|H_2| >M\}}.$$
Since $ H_2 \chi_{\{|H_2| \leq M\}} \in L^{\frac{2N}{\alpha+2s}}(\R^N) \cap L^{\infty}(\R^N)$ and $\frac{2N}{\alpha} \in (\frac{2N}{\alpha+2s}, \infty)$, we have $H_2 \chi_{\{|H_2| \leq M\}} \in L^{\frac{2N}{\alpha}}(\R^N)$, and thus
$$H^*:= H_1 + H_2 \chi_{\{|H_2| \leq M\}} \in L^{\frac{2N}{\alpha}}(\R^N), \quad H_* : =H_2 \chi_{\{|H_2| >M\}} \in L^{\frac{2N}{\alpha+2s}}(\R^N).$$
On the other hand
$$\norm{H_*}_{\frac{2N}{\alpha+2s}} = \left(\int_{|H_2|>M} |H_2|^{\frac{2N}{\alpha+2s}} \, dx\right)^{\frac{\alpha+2s}{2N}}$$
which can be made arbitrary small for $M\gg 0$. In particular we choose the decomposition so that
$$\left(\norm{H_*}_{\frac{2N}{\alpha+2s}} \norm{K_*}_{\frac{2N}{\alpha+2s}}\right) \lesssim \eps^2$$
and thus
$$C'(\eps):\approx \left(\norm{H^*}_{\frac{2N}{\alpha}}\norm{K^*}_{\frac{2N}{\alpha}}\right) .$$
In the last term of \eqref{eq_dim_mixed_terms_HK} we use the generalized Young's inequality $ab \leq \frac{\delta}{2} a^2 + \frac{1}{2\delta} b^2$, with
$$\delta :=\eps^2 \left(\norm{H^*}_{\frac{2N}{\alpha}}\norm{K_*}_{\frac{2N}{\alpha+2s}} + \norm{H_*}_{\frac{2N}{\alpha+2s}} \norm{K^*}_{\frac{2N}{\alpha}}\right)^{-1}$$
so that
%$$ \left(\norm{H^*}_{\frac{2N}{\alpha}}\norm{K_*}_{\frac{2N}{\alpha+2s}} + \norm{H_*}_{\frac{2N}{\alpha+2s}} \norm{K^*}_{\frac{2N}{\alpha}}\right) \norm{u}_2 \norm{(-\Delta)^{s/2} u}_2 \leq \tfrac{1}{2}\eps^2 \norm{u}_2^2 
%+ C''(\eps) \norm{(-\Delta)^{s/2}u}_2^2.$$
\begin{eqnarray*}
\lefteqn{ \left(\norm{H^*}_{\frac{2N}{\alpha}}\norm{K_*}_{\frac{2N}{\alpha+2s}} + \norm{H_*}_{\frac{2N}{\alpha+2s}} \norm{K^*}_{\frac{2N}{\alpha}}\right) \norm{u}_2 \norm{(-\Delta)^{s/2} u}_2} \\
&\qquad \qquad \qquad \leq& \tfrac{1}{2}\eps^2 \norm{u}_2^2 
+ C''(\eps) \norm{(-\Delta)^{s/2}u}_2^2. 
\end{eqnarray*}
Merging the pieces, we have the claim.
\QED

\bigskip

The following technical result can be found in \cite[Lemma 3.5]{GDS}. 

\begin{Lemma}[\cite{GDS}]\label{lem_truncation}
Let $a, b \in \R$, $r \geq 2$ and $k\geq 0$. Set $T_k: \R \to [-k,k]$ the truncation in $k$, that is
$$T_k(t):= \parag{ -k && \quad \hbox{ if $t\leq -k$}, \\ t&& \quad \hbox{ if $t \in (-k,k)$}, \\ k&& \quad \hbox{ if $t \geq k$},}$$
and write $a_k:= T_k(a)$, $b_k:= T_k(b)$. Then
$$\frac{4(r-1)}{r^2} \left(|a_k|^{r/2} - |b_k|^{r/2}\right)^2 \leq (a-b)\left( a_k |a_k|^{r-2} - b_k |b_k|^{r-2}\right).$$
\end{Lemma}

Notice that the (optimal) Sobolev embedding tells us that $H^s(\R^N) \hookrightarrow L^{2^*_s}(\R^N)$. In what follows we show that $u$ belongs to some $L^r(\R^N)$ with $r> 2^*_s$; %=\frac{2N}{N-2s}$; 
we highlight that we make no use of the Caffarelli-Silvestre $s$-harmonic extension method, and work directly in the fractional framework.

\begin{Proposition}\label{prop_uinLr}
Let $H, K \in L^{\frac{2N}{\alpha}}(\R^N) + L^{\frac{2N}{\alpha+2s}}(\R^N) $. Assume that $u \in H^s(\R^N)$ solves
$$(-\Delta)^s u + \mu u = (I_{\alpha}*(Hu))K, \quad \hbox{ in $\R^N$}$$
in the weak sense. Then
$$u \in L^r(\R^N) \quad \hbox{for all $r \in \left[2, \frac{N}{\alpha} \frac{2N}{N-2s}\right)$}.$$
Moreover, for each of these $r$, we have
$$\norm{u}_r \leq C_r \norm{u}_2$$
with $C_r>0$ not depending on $u$.
\end{Proposition}

\claim Proof.
By Lemma \ref{lem_secondMS} there exists $\lambda> \mu$ (that we can assume large) such that
\begin{equation}\label{eq_H_lambda}
\int_{\R^N} \left(I_{\alpha}*\big(H|u|\big)\right)K|u| \, dx \leq \frac{1}{2} \norm{(-\Delta)^{s/2} u}_2^2 + \frac{\lambda}{2}\norm{u}_2^2.
\end{equation}
Let us set
$$H_n:= H \chi_{\{|H|\leq n\}}, \quad K_n:= K \chi_{\{|K|\leq n\}}, \quad \hbox{ for $n \in \N$}$$
and observe that
$$H_n, \; K_n \in L^{\frac{2N}{\alpha}}(\R^N),$$
$$H_n \to H, \quad K_n \to K \quad \hbox{ almost everywhere, as $n\to +\infty$}$$
and
\begin{equation}\label{eq_H_nH}
|H_n| \leq |H|, \quad |K_n| \leq |K| \quad \hbox{ for every $n \in \N$}.
\end{equation}
We thus define the bilinear form
$$a_n(\varphi, \psi):= \int_{\R^N} (-\Delta)^{s/2} \varphi \, (-\Delta)^{s/2} \psi \, dx + \lambda \int_{\R^N} \varphi \psi \, dx - \int_{\R^N} \left(I_{\alpha}*\big(H_n \varphi\big)\right) K_n \psi \, dx$$
for every $\varphi, \psi \in H^s(\R^N)$. Since, by \eqref{eq_H_nH} and \eqref{eq_H_lambda}, we have
\begin{equation}\label{eq_coercive}
a_n(\varphi, \varphi) \geq \frac{1}{2} \norm{(-\Delta)^{s/2} \varphi}_2^2 + \frac{\lambda}{2} \norm{\varphi}_2^2 \geq \frac{1}{2} \norm{\varphi}_{H^s(\R^N)}^2
\end{equation}
for each $\varphi \in H^s(\R^N)$, we obtain that $a_n$ is coercive. Set
$$f:= (\lambda-\mu) u \in H^s(\R^N)$$
we obtain by Lax-Milgram theorem that, for each $n \in \N$, there exists a unique $u_n \in H^s(\R^N)$ solution of
$$a_n(u_n, \varphi)= (f, \varphi)_2, \quad \varphi \in H^s(\R^N),$$
that is
\begin{equation}\label{eq_troncat}
(-\Delta)^s u_n + \lambda u_n - \big(I_{\alpha}*(H_n u_n)\big)K_n =(\lambda-\mu) u, \quad \hbox{ in $\R^N$}
\end{equation}
in the weak sense; moreover the theorem tells us that
$$\norm{u_n}_{H^s} \leq \frac{\norm{f}_2}{1/2}= 2(\lambda-\mu) \norm{u}_2$$
(since $1/2$ appears as coercivity coefficient in \eqref{eq_coercive}), and thus $u_n$ is bounded. Hence $u_n \wto \bar{u}$ in $H^s(\R^N)$ up to a subsequence for some $\bar{u}$. This means in particular that $u_n \to \bar{u}$ almost everywhere pointwise.
 
 Thus we can pass to the limit in
$$\int_{\R^N} (-\Delta)^{s/2} u_n \, (-\Delta)^{s/2} \varphi + \lambda \int_{\R^N} u_n \varphi - \int_{\R^N} \left(I_{\alpha}*\big(H_n u_n\big)\right) K_n \varphi = (\lambda-\mu) \int_{\R^N} u \varphi ;$$
we need to check only the Choquard term. 
We first see by the continuous embedding that $u_n \wto \bar{u}$ in $L^q(\R^N)$, for $q \in [2, 2^*_s]$. Split again $H=H^*+ H_*$, $K=K^* + K_*$ and work separately in the four combinations; we assume to work generally with $\tilde{H} \in \{H^*, H_*\}$, $\tilde{H}\in L^{\beta}(\R^N)$ and $\tilde{K} \in \{K^*, K_*\}$, $\tilde{K}\in L^{\gamma}(\R^N)$, where $\beta, \gamma \in \{ \frac{2N}{\alpha}, \frac{2N}{\alpha+2s}\}$. 
Then one can easily prove that $\tilde{H}_n u_n \wto \tilde{H} \bar{u}$ in $L^r(\R^N)$ with $\frac{1}{r} = \frac{1}{\beta} + \frac{1}{q}$.
By the continuity and linearity of the Riesz potential we have $I_{\alpha} * (H_n u_n) \wto I_{\alpha} * (H \bar{u})$ in $L^h(\R^N)$, where $\frac{1}{h}= \frac{1}{r} - \frac{\alpha}{n}$. 
As before, we obtain $\left(I_{\alpha}*\big(H_n u_n\big)\right) K_n \wto \left(I_{\alpha}*\big(H \bar{u}\big)\right) K$ in $L^k(\R^N)$, where $\frac{1}{k} = \frac{1}{\gamma} + \frac{1}{h}$.
Simple computations show that if $\beta=\gamma=\frac{2N}{\alpha}$ and $q=2$, then $k'=2$; if $\beta=\frac{2N}{\alpha}$, $\gamma= \frac{2N}{\alpha+2s}$ (or viceversa) and $q=2$, then $k'=2^*_s$; if $\beta=\gamma=\frac{2N}{\alpha+2s}$ and $q=2^*_s$, then $k'=2^*_s$. 
Therefore $H^s(\R^N) \subset L^{k'}(\R^N)$ and we can pass to the limit in all the four pieces, obtaining
$$\int_{\R^N} \left(I_{\alpha}*\big(H_n u_n\big)\right) K_n \varphi \, dx \to \int_{\R^N} \left(I_{\alpha}*\big(H \bar{u}\big)\right) K \varphi \, dx.$$
Therefore, $\bar{u}$ satisfies
$$(-\Delta)^s \bar{u} + \lambda \bar{u} - \big(I_{\alpha}*(H \bar{u})\big)K =(\lambda-\mu) u, \quad \hbox{ in $\R^N$}$$
as well as $u$. But we can see this problem, similarly as before, with a Lax-Milgram formulation and obtain the uniqueness of the solution. Thus $\bar{u}=u$ and hence, as $n\to +\infty$, 
$$u_n \wto u \quad \hbox{ in $H^s(\R^N)$}$$
and almost everywhere pointwise. 
Let now $k\geq 0$ and write
$$u_{n,k}:= T_k(u_n)\in L^2(\R^N) \cap L^{\infty}(\R^N)$$
where $T_k$ is the truncation introduced in Lemma \ref{lem_truncation}. Let $r\geq 2$. We have $|u_{n,k}|^{r/2} \in H^s(\R^N)$, by exploiting \eqref{eq_semin_gagl} and the fact that $h(t):=(T_k(t))^{r/2}$ is a Lipschitz function with $h(0)=0$.
By \eqref{eq_semin_gagl} and by Lemma \ref{lem_truncation} we have
\begin{eqnarray*}
\lefteqn{\frac{4(r-1)}{r^2} \int_{\R^N} |(-\Delta)^{s/2} (|u_{n,k}|^{r/2}) |^2 %\, dx 
= C_{N,s} \int_{\R^{2N}} \frac{ \frac{4(r-1)}{r^2}\left(|u_{n,k}(x)|^{r/2} - |u_{n,k}(y)|^{r/2}\right)^2}{|x-y|^{N+2s}} %\, dx \, dy 
} \qquad \quad \\
&&\leq C_{N,s} \int_{\R^{2N}} \frac{\big(u_n(x)-u_n(y)\big)\left(u_{n,k}(x)|u_{n,k}(x)|^{r-2} - u_{n,k}(y) |u_{n,k}(y)|^{r-2}\right)}{|x-y|^{N+2s}} %\, dx \, dy .
\end{eqnarray*}
Set
$$\varphi:= u_{n,k}|u_{n,k}|^{r-2}$$
it results that $\varphi \in H^s(\R^N)$, since again $h(t):=T_k(t) |T_k(t)|^{r-2}$ is a Lipschitz function with $h(0)=0$.
Thus we can choose it as a test function in \eqref{eq_troncat} and obtain, by \eqref{eq_sobolev_polaraz}, %polarizing the identity \eqref{eq_semin_gagl},
\begin{eqnarray*}
\lefteqn{
 \frac{4(r-1)}{r^2} \int_{\R^N} |(-\Delta)^{s/2} (|u_{n,k}|^{r/2}) |^2 %\, dx
 \leq C_{N,s} \int_{\R^{2N}} \frac{\big(u_n(x)-u_n(y)\big)\left(\varphi(x) - \varphi(y)\right)}{|x-y|^{N+2s}} %\, dx \, dy
} \qquad \quad \\
 &&= - \lambda \int_{\R^N} u_n \varphi % \, dx
 + \int_{\R^N} \left(I_{\alpha}*(H_n u_n)\right) K_n \varphi% \, dx 
+ (\lambda-\mu) \int_{\R^N} u \varphi %\, dx
 \end{eqnarray*}
 and since $u_n \varphi \geq |u_{n,k}|^r$ we gain
 \begin{eqnarray}
 \lefteqn{
 \frac{4(r-1)}{r^2} \int_{\R^N} |(-\Delta)^{s/2} (|u_{n,k}|^{r/2}) |^2 %\, dx 
\leq } \notag\\
 &&\leq - \lambda \int_{\R^N} |u_{n,k}|^r %\, dx
+\int_{\R^N} \big(I_{\alpha}*(H_n u_n)\big) K_n \varphi %\, dx 
+ (\lambda-\mu) \int_{\R^N} u \varphi % \, dx
. \label{eq_dim_4(r-1)}
 \end{eqnarray}
 Focus on the Choquard term on the right-hand side. We have, by using \eqref{eq_H_nH},
% \begin{eqnarray}
%\lefteqn{\int_{\R^N} \big(I_{\alpha}*(H_n u_n)\big) K_n \varphi \leq} \notag\\
%&\leq& \int_{\R^N} \big(I_{\alpha}*(|H_n| |u_n|\chi_{\{|u_n|\leq k\}})\big) |K_n| |u_{n,k}|^{r-1} + \int_{\R^N} \big(I_{\alpha}*(|H_n| |u_n| \chi_{\{|u_n|>k\}})\big) |K_n| |u_{n,k}|^{r-1} \notag\\
%&\leq & \int_{\R^N} \big(I_{\alpha}*(|H_n| |u_{n,k}|)\big) |K_n| |u_{n,k}|^{r-1} +\int_{\R^N} \big(I_{\alpha}*(|H_n||u_{n}|\chi_{\{|u_n|>k\}})\big) |K_n| |u_n|^{r-1} \notag\\
%&%\stackrel{\eqref{eq_H_nH}}
%\leq & \int_{\R^N} \big(I_{\alpha}*(|H| |u_{n,k}|)\big) |K| |u_{n,k}|^{r-1} + \int_{\R^N} \big(I_{\alpha}*(|H_n||u_{n}|\chi_{\{|u_n|>k\}})\big) |K_n| |u_n|^{r-1} \notag\\
%& =:& (I)+(II). \label{eq_dim_(I-II)}
% \end{eqnarray}
%
 \begin{eqnarray}
\lefteqn{\int_{\R^N} \big(I_{\alpha}*(H_n u_n)\big) K_n \varphi \leq} \notag\\
&\leq& \int_{\R^N} \big(I_{\alpha}*(|H_n| |u_n|\chi_{\{|u_n|\leq k\}})\big) |K_n| |u_{n,k}|^{r-1} + \notag\\
&&+ \int_{\R^N} \big(I_{\alpha}*(|H_n| |u_n| \chi_{\{|u_n|>k\}})\big) |K_n| |u_{n,k}|^{r-1} \notag\\
&\leq & \int_{\R^N} \big(I_{\alpha}*(|H_n| |u_{n,k}|)\big) |K_n| |u_{n,k}|^{r-1} +\int_{\R^N} \big(I_{\alpha}*(|H_n||u_{n}|\chi_{\{|u_n|>k\}})\big) |K_n| |u_n|^{r-1} \notag\\
&%\stackrel{\eqref{eq_H_nH}}
\leq & \int_{\R^N} \big(I_{\alpha}*(|H| |u_{n,k}|)\big) |K| |u_{n,k}|^{r-1} + \int_{\R^N} \big(I_{\alpha}*(|H_n||u_{n}|\chi_{\{|u_n|>k\}})\big) |K_n| |u_n|^{r-1} \notag\\
& =:& (I)+(II). \label{eq_dim_(I-II)}
 \end{eqnarray}
 
Focus on $(I)$. Consider $r \in [2, \frac{2N}{\alpha})$, so that $\theta:= \frac{2}{r} \in (\frac{\alpha}{N}, 2-\frac{\alpha}{N})$. Choose moreover $v:= |u_{n,k}|^{r/2} \in H^s(\R^N)$ and $\eps^2:= \frac{2(r-1)}{r^2}>0$.
 Thus, observed that if a function belongs to a sum of Lebesgue spaces then its absolute value does the same (see Remark \ref{rem_somma_lebesg}), by Lemma \ref{lem_secondMS} we obtain
 \begin{equation}\label{eq_dim_(I)}
 (I) \leq \frac{2(r-1)}{r^2} \norm{(-\Delta)^{s/2}(|u_{n,k}|^{r/2})}_2^2 + C(r) \norm{|u_{n,k}|^{r/2}}_2^2.
 \end{equation}
 Focus on $(II)$. Assuming $r< \min\{\frac{2N}{\alpha}, \frac{2N}{N-2s}\}$, we have $u_n \in L^r(\R^N)$ 
 and $H_n \in L^{\frac{2N}{\alpha}}(\R^N)$, thus
 $$|H_n| |u_n| \in L^{a}(\R^N), \quad \hbox{with $\frac{1}{a} = \frac{\alpha}{2N} + \frac{1}{r}$}$$
 for the H\"older inequality. Similarly
 $$|K_n| |u_n|^{r-1} \in L^{b}(\R^N), \quad \hbox{with $\frac{1}{b} = \frac{\alpha}{2N} + 1-\frac{1}{r}$}.$$
 Thus, since $\frac{1}{a} + \frac{1}{b}= \frac{N+\alpha}{N}$, we have by the Hardy-Littlewood-Sobolev inequality (see Proposition \ref{prop_HLS}) that
\begin{eqnarray*}
\lefteqn{
\int_{\R^N} \big(I_{\alpha}*(|H_n||u_{n}|\chi_{\{|u_n|>k\}})\big) |K_n| |u_n|^{r-1} \, dx} \\
&& \leq C\left( \int_{\{|u_n|>k\}} \abs{|H_n| |u_n|}^a \, dx\right)^{1/a} \left( \int_{\R^N}\abs{|K_n||u_n|^{r-1}}^{b}\, dx\right)^{1/b} .
\end{eqnarray*}
 With respect to $k$, the second factor on the right-hand side is bounded, while the first factor goes to zero thanks to the dominated convergence theorem, thus
\begin{equation}\label{eq_dim_(II)}
(II) = o_k(1), \quad \hbox{ as $k\to +\infty$}.
\end{equation}
 Joining \eqref{eq_dim_4(r-1)}, \eqref{eq_dim_(I-II)}, \eqref{eq_dim_(I)}, \eqref{eq_dim_(II)} we obtain
\begin{eqnarray*}
\lefteqn{ \frac{2(r-1)}{r^2} \int_{\R^N} |(-\Delta)^{s/2} (|u_{n,k}|^{r/2}) |^2\, dx \leq} \\
&& \leq - \lambda \int_{\R^N} |u_{n,k}|^r
 \, dx + C(r) \int_{\R^N} |u_{n,k}|^r
 \, dx + (\lambda-\mu) \int_{\R^N} u \varphi \, dx + o_k(1).
\end{eqnarray*}
 That is, by Sobolev inequality \eqref{eq_embd_homog} 
$$C'(r)\left( \int_{\R^N} |u_{n,k}|^{\frac{r}{2} 2^*_s} % \, dx 
\right)^{2/2^*_s} \leq (C(r)-\lambda) \int_{\R^N} |u_{n,k}|^r %\, dx 
+ (\lambda-\mu)
 \int_{\R^N} |u| \, |u_{n,k}|^{r-1} %\, dx 
+ o_k(1).$$
Letting $k\to +\infty$ by the monotone convergence theorem (since $u_{n,k}$ are monotone with respect to $k$ 
and $u_{n,k} \to u_n$ pointwise)
we have
\begin{equation}\label{eq_moser}
C'(r)\left( \int_{\R^N} |u_n|^{\frac{r}{2} 2^*_s} %\, dx
 \right)^{2/2^*_s} \leq (C(r)-\lambda) \int_{\R^N} |u_{n}|^r %\, dx 
+ (\lambda-\mu) \int_{\R^N} |u| \, |u_{n}|^{r-1} %\, dx
\end{equation}
and thus $u_n \in L^{\frac{r}{2}2^*_s}(\R^N)$. Notice that $\frac{r}{2} \in \big[1, \min\{\frac{N}{\alpha}, \frac{N}{N-2s}\}\big)$. If $N-2s < \alpha$ we are done. Otherwise, set $r_1:=r$, we can now repeat the argument with
$$r_2 \in \left( \frac{2N}{N-2s}, \min\left\{ \frac{2N}{\alpha}, 2\left(\frac{N}{N-2s}\right)^2\right\}\right).$$
Again, if $\frac{2N}{\alpha} < 2\left(\frac{N}{N-2s}\right)^2$ we are done, otherwise we repeat the argument. Inductively, we have
$$\left(\frac{N}{N-2s}\right)^m \to +\infty, \quad \hbox{as $m\to +\infty$}$$
thus $\frac{2N}{\alpha}<2\left(\frac{N}{N-2s}\right)^m$ after a finite number of steps. For such $r=r_m$, consider again \eqref{eq_moser}: 
by the almost everywhere convergence of $u_n$ to $u$ and Fatou's lemma
\begin{eqnarray*}
\lefteqn{C''(r) \left(\int_{\R^N} |u|^{\frac{r}{2} 2^*_s}\right)^{2/2^*_s} \, dx \leq \liminf_{n} C''(r) \left( \int_{\R^N} |u_n|^{\frac{r}{2} 2^*_s} \, dx \right)^{2/2^*_s}} \\ 
&&\leq \liminf_n \left( (C(r)-\lambda) \int_{\R^N} |u_{n}|^r \, dx + (\lambda-\mu) \int_{\R^N} |u| \, |u_{n}|^{r-1} \, dx\right) \\
&&\leq (C(r)-\lambda ) \limsup_n\int_{\R^N} |u_{n}|^r \, dx + (\lambda-\mu) \limsup_n \int_{\R^N} |u| \, |u_{n}|^{r-1} \, dx.
\end{eqnarray*}
Being $u_n$ equibounded in $H^s(\R^N)$ and thus in $L^{2^*_s}(\R^N)$, by the iteration argument we have that it is equibounded also in $L^r(\R^N)$; in particular, the bound is given by $\norm{u}_2$ times a constant $C(r)$. Thus the right-hand side is a finite quantity, and we gain $u \in L^{\frac{r}{2}2^*_s}(\R^N)$, which is the claim.
\QED

\bigskip

The following lemma states that $I_{\alpha}*g \in L^{\infty}(\R^N)$ whenever $g$ lies in $L^q(\R^N)$ with $q$ in a neighborhood of $\frac{N}{\alpha}$; in particular, it extends Proposition \ref{prop_HLS} (see also Remark \ref{rem_conv_welldef}). %prop_HLS} to the case $h=\infty$ and $r \approx \frac{N}{\alpha}$).

In addition, it shows the decay at infinity of the Riesz potential, which will be useful in Section \ref{sec_asymptotic}.

\begin{Proposition}\label{prop_conv_C0}
Assume that \hyperref[(F1s)]{\textnormal{(F1)}}-\hyperref[(F2s)]{\textnormal{(F2)}} hold. Let $u\in H^s(\R^N)$ be a solution of \eqref{eq_introduction}. Then 
$u\in L^q(\R^N)$ for $q \in \big[2,\frac{N}{\alpha} \frac{2N}{N-2s}\big)$, and
$$I_{\alpha} * F(u) \in C_0(\R^N),$$
that is, continuous and zero at infinity. 
In particular, 
$$I_{\alpha} * F(u) \in L^{\infty}(\R^N)$$
and
$$\big(I_{\alpha} * F(u)\big)(x) \to 0 \quad \hbox{as $|x| \to +\infty$}.$$
\end{Proposition}

\claim Proof.
We first check to be in the assumptions of Proposition \ref{prop_uinLr}. Indeed, by \hyperref[(F1s)]{\textnormal{(F1)}}-\hyperref[(F2s)]{\textnormal{(F2)}} and the fact that $u\in H^s(\R^N)\subset L^2(\R^N) \cap L^{2^*_s}(\R^N)$ we obtain that
$$H:= \frac{F(u)}{u}, \quad K:=f(u)$$
lie in $L^{\frac{2N}{\alpha}}(\R^N) + L^{\frac{2N}{\alpha+2s}}(\R^N)$, since bounded by functions in this sum space (see Remark \ref{rem_somma_lebesg}). 
Now by Proposition \ref{prop_uinLr} we have $u \in L^q(\R^N)$ for $q \in [2, \frac{N}{\alpha} \frac{2N}{N-2s})$; 
the claim follows by Remark \ref{rem_conv_welldef}.
%To gain the information on the convolution, we want to use Young's Theorem, which states that if $g, h$ belong to two Lebesgue spaces with conjugate (finite) indexes, then $g*h \in C_0(\R^N)$. 
%We first split
%$$I_{\alpha}*F(u) = (I_{\alpha}\chi_{B_1})*F(u) + (I_{\alpha}\chi_{B_1^c})*F(u)$$
%where
%$$ I_{\alpha}\chi_{B_1} \in L^{r_1}(\R^N), \quad \hbox{ for $r_1 \in [1, \frac{N}{N-\alpha})$},$$
%$$ I_{\alpha}\chi_{B_1^c} \in L^{r_2}(\R^N), \quad \hbox{ for $r_2 \in (\frac{N}{N-\alpha}, \infty]$}.$$
%We need to show that $F(u) \in L^{q_1}(\R^N)\cap L^{q_2}(\R^N)$ for some $q_i$ satisfying
%$$\frac{1}{q_i} + \frac{1}{r_i} = 1, \quad i=1,2$$
%that is
%$$\frac{q_1}{q_1-1} \in \left[1, \frac{N}{N-\alpha}\right), \quad \frac{q_2}{q_2-1}\in \left(\frac{N}{N-\alpha}, \infty\right]$$
%or equivalently $q_2 < \frac{N}{\alpha} < q_1$.
%Recall that
%$$|F(u)| \leq C\left(|u|^{\frac{N+\alpha}{N}} + |u|^{\frac{N+\alpha}{N-2s}}\right).$$
%Note that $u \in L^q(\R^N)$ for $q \in [2, \frac{N}{\alpha} \frac{2N}{N-2s})$ implies
%$$|u|^{\frac{N+\alpha}{N}} , |u|^{\frac{N+\alpha}{N-2s}} \in L^{q_1}(\R^N) \cap L^{q_2}(\R^N)$$
%for some $q_2<\frac{N}{\alpha}<q_1$. Thus we have the claim.
\QED

\bigskip

Once obtained the boundedness of the Choquard term, we can finally gain the boundedness of the solution.

\begin{Proposition}\label{prop_u_bounded}
Assume that \hyperref[(F1s)]{\textnormal{(F1)}}-\hyperref[(F2s)]{\textnormal{(F2)}} hold. Let $u\in H^s(\R^N)$ be a positive solution of \eqref{eq_introduction}. Then $u\in L^{\infty}(\R^N)$.
\end{Proposition}

\claim Proof.
By Lemma \ref{prop_conv_C0} we obtain
$$a:= I_{\alpha}*F(u) \in L^{\infty}(\R^N).$$
Thus $u$ satisfies the following nonautonomous problem, with a local nonlinearity
$$(-\Delta)^{s/2} u + \mu u = a(x) f(u), \quad \hbox{ in $\R^N$}$$
with $a$ bounded. In particular
$$(-\Delta)^{s/2} u = g(x,u):=- \mu u + a(x) f(u), \quad \hbox{ in $\R^N$}$$
where
$$|g(x,t)| \leq \mu |t| + C \norm{a}_{\infty} \left(|t|^{\frac{\alpha}{N}} + |t|^{\frac{\alpha+2s}{N-2s}}\right).$$
Set $\gamma:= \max\{1, \frac{\alpha+2s}{N-2s}\} \in [1, 2^*_s)$, we thus have
$$|g(x,t)| \leq C(1 + |t|^{\gamma}).$$
Hence we are in the assumptions of \cite[Proposition 5.1.1]{DMV} and we can conclude. 
\QED

\bigskip

\claim Proof of Theorem \ref{th_INT_regular}.
The first part of the claim comes from Proposition \ref{prop_u_bounded}. In the case of sign-changing %signed 
solutions, we may apply Proposition \ref{prop_prop_u_Linf} with
$$g(x,t) := \big(I_{\alpha}*F(u)\big)(x) f(t) - \mu u,$$
whenever $u$ is a fixed solution and \hyperref[(F6)]{\textnormal{(F6)}} holds (together with \hyperref[(F1s)]{\textnormal{(F1)}}--\hyperref[(F2s)]{\textnormal{(F2)}}), thanks to Proposition \ref{prop_conv_C0}.
\QED

%%%%%%%%%%%%%%%%%%%%%%%%%%%%%%%
\subsection{H\"older regularity: strong solutions}

Gained the boundedness, % of the solutions, 
we obtain now that solutions are H\"older continuous and satisfy the equation in the strong sense.
% also some additional regularity, which 
This extra regularity will be also implemented in some bootstrap argument for the $L^1$-summability, see Section \ref{sec_L1_sum}.

\begin{Proposition}\label{prop_u_holder}
Assume that \hyperref[(F1s)]{\textnormal{(F1)}}-\hyperref[(F2s)]{\textnormal{(F2)}} hold. Let $u\in H^s(\R^N)\cap L^{\infty}(\R^N)$ be a weak solution of \eqref{eq_introduction}. Then $u \in H^{2s}(\R^N) \cap C^{0,\gamma}(\R^N)$ for any $\gamma \in (0, \min\{1,2s\})$, and $u$ is a strong solution, i.e. $u$ satisfies \eqref{eq_introduction} almost everywhere.

In addition, if $s\in (\tfrac{1}{2}, 1)$, then $u \in C^{1,\gamma}(\R^N)$ for any $\gamma \in (0, 2s-1)$.
\end{Proposition}

\claim Proof.
By Proposition \ref{prop_u_bounded}, Proposition \ref{prop_conv_C0} and \hyperref[(F2s)]{\textnormal{(F2)}} we have that $u\in L^{\infty}(\R^N)$ satisfies
$$(-\Delta)^s u = g \in L^{\infty}(\R^N)$$
where $g(x):= (I_{\alpha}*F(u))(x) f(u(x))- \mu u(x)$. 
We prove first that $u\in H^{2s}(\R^N)$. Indeed, we already know that $f(u)$, $F(u)$ and $I_{\alpha}*F(u)$ belong to $L^{\infty}(\R^N)$. By Remark \ref{rem_buona_posit_f}, we obtain
$$f(u) \in L^{\frac{2N}{\alpha+2s}}(\R^N)\cap L^{\infty}(\R^N), \quad F(u) \in L^{\frac{2N}{N+\alpha}}(\R^N) \cap L^{\infty}(\R^N),$$
$$I_{\alpha}*F(u) \in L^{\frac{2N}{N-2s}}(\R^N) \cap L^{\infty}(\R^N), \quad (I_{\alpha}*F(u))f(u) \in L^2(\R^N)\cap L^{\infty}(\R^N).$$
In particular,
$$g=(I_{\alpha}*F(u))f(u)-\mu u \in L^2(\R^N).$$
Since $u$ is a weak solution, we have, fixed $\varphi \in H^s(\R^N)$,
\begin{equation}\label{eq_dim_weak_form}
\int_{\R^N} (-\Delta)^{s/2} u \, (-\Delta)^{s/2} \varphi \, dx = \int_{\R^N} g \, \varphi \, dx.
\end{equation}
Since $g \in L^2(\R^N)$, we can apply Plancharel theorem and obtain
\begin{equation}\label{eq_dim_plancharel}
\int_{\R^N} |\xi|^{2s} \widehat{u} \, \widehat{\varphi} \, d\xi = \int_{\R^N} \widehat{g} \, \widehat{\varphi} \, d \xi.
\end{equation}
Since $H^s(\R^N) = \mc{F}(H^s(\R^N))$ and $\varphi$ is arbitrary, we gain
$$|\xi|^{2s} \widehat{u} = \widehat{g} \in L^2(\R^N).$$
By definition, we obtain $u \in H^{2s}(\R^N)$, which concludes the proof.
Observe moreover that $\mc{F}^{-1}\big((1+|\xi|^{2s})\widehat{u}\big) = u +g \in L^2(\R^N) \cap L^{\infty}(\R^N)$, thus by definition $u \in H^{2s}(\R^N) \cap W^{2s, \infty}(\R^N)$. By the embedding \eqref{eq_immers_Holder} (see also Proposition \ref{prop_reg_fraz}) we obtain $u \in C^{0,\gamma}(\R^N)$ if $2s<1$ and $\gamma \in (0, 2s)$, while $u \in C^{1,\gamma}(\R^N)$ if $2s>1$ and $\gamma \in (0,2s-1)$. 

It remains to show that $u$ is an almost everywhere pointwise solution. Thanks to the fact that $u\in H^{2s}(\R^N)$, we use again \eqref{eq_dim_plancharel}, where we can apply Plancharel theorem (that is, we are integrating by parts \eqref{eq_dim_weak_form}) and thus
$$\int_{\R^N}(-\Delta)^s u \, \varphi \, dx = \int_{\R^N} g \, \varphi \, dx.$$
Since $\varphi \in H^s(\R^N)$ is arbitrary, we obtain
$$(-\Delta)^s u = g \quad \hbox{ almost everywhere.}$$
This concludes the proof.
\QED

\bigskip

We observe, by the proof, that if $s\in (\tfrac{1}{2}, 1)$, then %$u \in C^{1,\gamma}(\R^N)$ for any $\gamma \in (0, 2s-1)$, and 
$u$ is a classical solution, with $(-\Delta)^s u \in C(\R^N)$ and equation \eqref{eq_introduction} satisfied pointwise.
We will further investigate these aspects in Section \ref{sec_regularity_fp}.

%%%%%%%%%%%%%%%%%%%%%%%%%%%%%%%
\subsection{$L^1$-summability: fixed point maps}
\label{sec_L1_sum}

We deal now with the summability of $u$ in Lebesgue spaces $L^r(\R^N)$ for $r<2$.
We observe that the information $u\in L^1(\R^N)\cap L^2(\R^N)$ is new even in the power-type setting: indeed in \cite{DSS1} the authors, in order to ensure existence of solutions, assume the nonlinearity to be not 
%lower 
%\tr{(?)} %COMMENT NOW
critical, while here we can include the possibility of criticality.
Moreover, this result is new even for $s=1$, improving \cite{MS2}. The $L^1$-summability will be then used also to gain the asymptotic behaviour of the solutions %when $f$ is sublinear, 
in Section \ref{sec_asymptotic}. %\ref{sec_sublin_case}.

\begin{Remark}\label{rem_L1}
We start noticing that, if a solution $u$ belongs to some $L^q(\R^N)$ with $q<2$, then $u\in L^1(\R^N)$. 
Assume thus $q \in (1,2)$ and let $u \in L^q(\R^N) \cap L^{\infty}(\R^N)$, then we have
$$f(u) \in L^{\frac{qN}{\alpha}}(\R^N)\cap L^{\infty}(\R^N), \quad F(u) \in L^{\frac{qN}{N+\alpha}}(\R^N) \cap L^{\infty}(\R^N),$$
$$I_{\alpha}* F(u) \in L^{\frac{qN}{N+ \alpha(1-q)}}(\R^N) \cap L^{\infty}(\R^N), \quad (I_{\alpha}*F(u)) f(u) \in L^{\frac{qN}{N + \alpha(2-q)}}(\R^N) \cap L^{\infty}(\R^N).$$
Thanks to Proposition \ref{prop_u_holder}, $u$ satisfies \eqref{eq_introduction} almost everywhere, thus we have
$$\mc{F}^{-1}\big((|\xi|^{2s} + \mu)\, \widehat{u}\big) = (-\Delta)^s u + \mu u = (I_{\alpha}*F(u))f(u) \in L^{\frac{qN}{N + \alpha(2-q)}}(\R^N)$$
%which equivalently means that the Bessel operator verifies
%$$\mc{F}^{-1}\big((|\xi|^{2} + 1)^s \,\widehat{u}\big) \in L^{\frac{qN}{N + \alpha(2-q)}}(\R^N).$$
%T
hence by the properties of the Bessel operator \eqref{eq_Bessel_guadagna} we obtain that $u$ itself lies in the same Lebesgue space, that is
$$u \in L^{\frac{qN}{N + \alpha(2-q)}}(\R^N).$$
If $\frac{qN}{N + \alpha(2-q)}<1$, we mean that $(I_{\alpha}*F(u)) f(u) \in L^1(\R^N) \cap L^{\infty}(\R^N)$, and thus $u \in L^1(\R^N) \cap L^{\infty}(\R^N)$. We convey this when we deal with exponents less than $1$.

If $q <2$, then
$$\frac{qN}{N + \alpha(2-q)}<q$$
and we can implement a bootstrap argument to gain $u \in L^1(\R^N)$. More precisely
$$\parag{ &q_0\in [1, 2)& \\ &q_{n+1} = \frac{q_n N}{N + \alpha(2-q_n)}& }$$
where $q_n \to 0$ (but we stop at $1$). 
%Thus, in order to implement the argument, we need to show that $u\in L^q(\R^N)$ for some $q<2$.
\end{Remark}

We show now that $u\in L^1(\R^N)$. 
It is easy to see that, if the problem is (strictly) not lower-critical, i.e., \hyperref[(F2s)]{\textnormal{(F2)}} holds together with
	$$	\lim_{t\to 0}\frac{F(t)}{\abs{t}^{\beta}}=0
	$$
for some $\beta\in (\frac{N+\alpha}{N}, \frac{N+\alpha}{N-2s})$, then $u\in L^1(\R^N)$. 
Indeed $u\in H^s(\R^N)\cap L^\infty(\R^N)\subset
L^2(\R^N)\cap L^\infty(\R^N)$ and
	$$	(I_\alpha*F(u))f(u) \in L^q(\R^N),
	$$
where $\frac{1}{q}=\frac{\beta}{2}-\frac{\alpha}{2N}$; noticed that $q<2$, we can implement the bootstrap argument of Remark \ref{rem_L1}.

We will show that the same conclusion can be reached by assuming only \hyperref[(F2s)]{\textnormal{(F2)}}. 

\begin{Proposition}\label{prop_u_L1}
Assume that \hyperref[(F1s)]{\textnormal{(F1)}}-\hyperref[(F2s)]{\textnormal{(F2)}} hold. Let $u\in H^s(\R^N)\cap L^{\infty}(\R^N)$ be a weak solution of \eqref{eq_introduction}. 
Then $u \in L^1(\R^N)$.
\end{Proposition}

\claim Proof. % of Proposition \ref{prop_u_L1}.
For a given solution $u\in H^s(\R^N)\cap L^\infty(\R^N)$ we
set again
 $$ H:=\frac{F(u)}{u}, \quad K:=f(u).
 $$
Since $u\in L^2(\R^N)\cap L^\infty(\R^N)$, by \hyperref[(F2s)]{\textnormal{(F2)}} we have
$H$, $K\in L^{\frac{2N}{\alpha}}(\R^N)$. For $n\in\N$, we set
 $$ H_n:=H\chi_{\{ |x|\geq n\} }.
 $$
Then we have
 \begin{equation}\label{kt.a}
 \norm{H_n}_{\frac{2N}{\alpha}}\to 0 \quad
 \hbox{as}\ n\to\infty.
 \end{equation}
Since $\hbox{supp} (H-H_n)\subset \big\{ |x|\leq n\big\}$ is a
bounded set, we have for any $\beta \in [1,{\frac{2N}{\alpha}}]$ 
 \begin{equation}\label{kt.b}
 H-H_n \in L^\beta(\R^N) \quad \hbox{for all}\ n\in\N.
 \end{equation}
We write our equation \eqref{eq_introduction} as
 $$ (-\Delta)^s u+\mu u =(I_\alpha*H_nu)K +R_n \quad \hbox{in $\mathbb{R}^N$},
 $$
where we introduced the function $R_n$ by
 $$ R_n:= (I_\alpha*(H-H_n)u)K.
 $$
Now we consider the following linear equation:
 \begin{equation}\label{kt.c}
 (-\Delta)^s v+\mu v =(I_\alpha*H_nv)K +R_n \quad \hbox{in $\mathbb{R}^N$}.
 \end{equation}
We have the following facts:
\begin{itemize}
\item[(i)] The given solution $u$ solves \eqref{kt.c}.
\item[(ii)] By the property \eqref{kt.b} with $\beta \in (\frac{2N}{N+\alpha},{\frac{2N}{\alpha}})$, there exists $q_1
\in (1,2)$, namely $\frac{1}{q_1}=\frac{1}{\beta}+\half-\frac{\alpha}{2N}$, 
such that $R_n\in L^{q_1}(\R^N)\cap L^2(\R^N)$.
\item[(iii)] By the property \eqref{kt.a}, for any $r\in (\frac{2N}{2N-\alpha}, 2] \subset (1,2]$
 $$ v\in L^r(\R^N) \mapsto A_n(v):=(I_\alpha*H_n v)K\in L^r(\R^N)
 $$
is well defined and verifies
 \begin{equation}\label{kt.d}
 \norm{A_n(v)}_r \leq C_{r,n}\norm v_r.
 \end{equation}
Here $C_{r,n}$ satisfies $C_{r,n}\to 0$ as $n\to\infty$.
\end{itemize}
We show only (iii). Since $v\in L^r(\R^N)$, by Hardy-Littlewood-Sobolev inequality and H\"older inequality we obtain
 $$ \norm{A_n(v)}_r \leq C_r\norm{H_n}_{\frac{2N}{\alpha}}
 \norm K_{\frac{2N}{\alpha}} \norm v_r,
 $$
where $C_r>0$ is independent of $n$, $v$. Thus by \eqref{kt.a} we have
$C_{r,n}:=C_r\norm{H_n}_{\frac{2N}{\alpha}}\norm K_{\frac{2N}{\alpha}}\to 0$ 
as $n\to\infty$.

Now we show $u\in L^{q_1}(\R^N)$, where $q_1\in (1,2)$ is
given in (ii).
Since $((-\Delta)^s+\mu)^{-1}:\, L^r(\R^N)\to L^r(\R^N)$ is
a bounded linear operator for $r\in(1,2]$ (see \eqref{eq_operat_contin_K}), \eqref{kt.c} can be rewritten as
 $$ v=T_n(v), $$
where
 $$ T_n(v):=((-\Delta)^s+\mu)^{-1}\big(A_n(v)+R_n\big). $$
By choosing $\beta \in (2, \frac{2N}{\alpha})$ we have $q_1 \in (\frac{2N}{2N-\alpha},2)\subset(1,2)$, thus we observe that for $n$ large, $T_n$ is a contraction in $L^2(\R^N)$ and in $L^{q_1}(\R^N)$. We fix such an $n$.

Since $T_n$ is a contraction in $L^2(\R^N)$, we can see that
$u\in H^s(\R^N)$ is a unique fixed point of $T_n$.
In particular, we have
 $$ u = \lim_{k\to \infty} T_n^k (0) \quad \hbox{in} \ L^2(\R^N).
 $$
On the other hand, since $T_n$ is a contraction in $L^{q_1}(\R^N)$,
$(T_n^k (0))_{k=1}^\infty$ also converges in $L^{q_1}(\R^N)$. Thus the
limit $u$ belongs to $L^{q_1}(\R^N)$.

Since $q_1<2$ we can use the bootstrap argument of Remark \ref{rem_L1} to get $u\in L^1(\R^N)$, and reach the claim.
\QED

\bigskip

With similar arguments we obtain also the following result for $s=1$.
\begin{Proposition}%\label{prop_u_L1}
Let $s=1$ and assume $N\geq 3$ and \hyperref[(F1)]{\textnormal{(F1)}}-\hyperref[(F2)]{\textnormal{(F2)}}. %with $s=1$. 
Let $u\in H^1(\R^N)\cap L^{\infty}(\R^N)$ be a weak solution of \eqref{eq_Choquard_genericaF}
Then $u \in L^1(\R^N)$.
\end{Proposition}

%%%%%%%%%%%%%%%%%%%%%%%%%%%%%%%
\subsection{$C^{\gamma}$-regularity: classical solutions}
\label{sec_regularity_fp}

We continue the analysis of the regularity started in Proposition \ref{prop_u_holder} and we infer the following %main %regularity 
result. 
This extra regularity will be exploited in the discussion of the positivity of solutions, see Section \ref{sec_positivity}; some more results about the regularity of the solutions will be stated in Section \ref{sec_doub_C1C2}.

Consider the condition
\begin{itemize}
\item[\textnormal{(F7)}] \label{(F7)}
%$s \in (\tfrac{1}{2},1)$ or $s\in (0,\tfrac{1}{2}]$ with 
$f\in C^{0,\sigma}_{loc}(\R)$ for some $\sigma \in (0,1]$. %)$.
\end{itemize}

\begin{Proposition}
\label{prop_regolar} %\label{th_INT_regular_N}
%Assume \textnormal{(f5)} in addition to \textnormal{(f1)-(f2)}. 
Assume \hyperref[(F1s)]{\textnormal{(F1)}}-\hyperref[(F2s)]{\textnormal{(F2)}}. 
%Moreover, assume
%\begin{itemize}
%\item[\textnormal{(F7)}] \label{(F7)}
%%$s \in (\tfrac{1}{2},1)$ or $s\in (0,\tfrac{1}{2}]$ with 
%$f\in C^{0,\sigma}_{loc}(\R)$ for some $\sigma \in (0,1]$. %)$.
%\end{itemize}
Let $u\in H^s(\R^N)\cap L^{\infty}(\R^N)$ be a weak solution of \eqref{eq_introduction}. 
If $s\in (\tfrac{1}{2}, 1)$, then $u \in C^{1,\gamma}(\R^N)$ for any $\gamma \in (0, 2s-1)$ and $u$ is a classical solution.

Assume now instead $s \in (0,1)$ and in addition \hyperref[(F7)]{\textnormal{(F7)}}. 
Then $u$ is a classical solution, that is a pointwise solution lying in
\begin{itemize}
\item $C^{0, \gamma}(\R^N) \cap H^{2s}(\R^N)$ for some $\gamma>2s$, if $2s < 1$,
\item $C^{1, \gamma-1}(\R^N) \cap H^{2s}(\R^N)$ for some $\gamma>2s$, if $2s \geq 1$.
\end{itemize}
More specifically, %if $f\in C_{loc}^{0,\sigma}(\R)$ for some $\sigma \in (0,1)$, 
set $\omega:=\min\{\sigma, 2s \sigma, \alpha\}$, we have
%\begin{itemize}
%\item if $\omega + 2s \leq 1$, then $u \in C^{0, \gamma}(\R^N)$ for every $\gamma \in (0,\omega+2s)$,
%\item if $\omega + 2s > 1$, then $u \in C^{1,\gamma}(\R^N)$ for every $\gamma \in (0,\omega+2s-1)$.
%\end{itemize}
\begin{itemize}
\item if $\omega + 2s \in (0,1]$, then $u \in C^{0, \gamma}(\R^N)$ for each $\gamma \in (0, \omega+2s]\cap (0,1)$, %$\gamma \leq \omega+2s$, $\gamma <1$,
\item if $\omega + 2s \in (1,2]$, then $u \in C^{1,\gamma-1}(\R^N)$ for each $\gamma \in (0,\omega+2s] \cap (0,2)$, %$\gamma \leq \omega+2s$, $\gamma <2$,
\item if $\omega + 2s \in (2,3)$, then $u \in C^{2,\omega+2s-2}(\R^N)$. % for each $\gamma \leq \omega+2s$,
\end{itemize}
\end{Proposition}

\claim Proof.
%We notice that the first regularity part of the claim is actually already known for $s\in (\frac{1}{2}, 1)$ thanks to Theorem \ref{teo_u_holder}.
%
Start noticing that by Proposition \ref{prop_conv_C0} we have 
%by the assumption, $u \in L^2(\R^N) \cap L^{\infty}(\R^N)$, and thus 
$I_{\alpha}*F(u) \in C_0(\R^N)$; in particular $I_{\alpha}*F(u)$ is pointwise finite. %finite almost everywhere.
Moreover, by \hyperref[(F2)]{\textnormal{(F2)}} we have $F(u) \in L^{\frac{2N}{N+\alpha}}(\R^N)\cap L^{\infty}(\R^N)$. If we choose
$$\parag{q \in [\tfrac{2N}{N+\alpha}, \infty) \quad \hbox{ if $\alpha \in (0,1]$}, \\ q \in [\tfrac{2N}{N+\alpha}, \tfrac{N}{\alpha-1}) \quad \hbox{ if $\alpha \in (1,N)$}}$$
we obtain
$$F(u) \in L^q(\R^N),
\quad \tfrac{N}{q} < \alpha < 1+ \tfrac{N}{q}$$
and thus we can apply Proposition \ref{prop_regol_riesz} to conclude that
$$I_{\alpha}* F(u) \in C^{0, \alpha - \frac{N}{q}}(\R^N).$$
In particular, by suitable choices of $q$, we gain
$$I_{\alpha}* F(u) \in C^{0, \omega}(\R^N) \quad \hbox{ for every $\omega \in (0,\min\{1,\alpha\})$}.$$
Notice that up to now we did not use the regularity on $f$.
Assume \hyperref[(F7)]{\textnormal{(F7)}} now.
By Proposition \ref{prop_u_holder} we have that $u$ is bounded and $u \in C^{0, \gamma}(\R^N)$ for every $\gamma \in (0, \min\{1,2s\})$.
By composition, we obtain
$$f(u) \in C^{0,\theta}(\R^N), \quad \hbox{for $\theta \in (0, \min\{\sigma, 2s\sigma\})$}.$$
Chosen
$$\omega\equiv\theta \in (0, \min\{\sigma, 2s \sigma, \alpha\})$$
then, since both $f(u)$ and $I_{\alpha}*F(u)$ are bounded and H\"older continuous, we have
$$(I_{\alpha}*F(u))f(u) \in C^{0,\omega}(\R^N).$$
At this point we can use Proposition \ref{prop_reg_fraz} to gain
%\begin{itemize}
%\item if $\otmega + 2s \leq 1$, then $u \in C^{0, \omega +2s}(\R^N)$,
%\item if $\omega + 2s > 1$, then $u \in C^{1,\omega + 2s-1}(\R^N)$,
%\end{itemize}
\begin{itemize}
\item if $\omega + 2s \in (0,1]$, then $u \in C^{0, \gamma}(\R^N)$ for each $\gamma \leq \omega+2s$, $\gamma <1$,
\item if $\omega + 2s \in (1,2]$, then $u \in C^{1,\gamma-1}(\R^N)$ for each $\gamma \leq \omega+2s$, $\gamma <2$,
\item if $\omega + 2s \in (2,3)$, then $u \in C^{2,\omega+2s-2}(\R^N)$, % for each $\gamma \leq \omega+2s$,
\end{itemize}
and thus the regularity claim. 
Finally, again by Proposition \ref{prop_u_holder} %Theorem \ref{teo_u_holder} 
$u$ satisfies \eqref{eq_introduction} almost everywhere; moreover, by the achieved regularity and Proposition \ref{prop_well_posed}, %prop_buon_def_lap}, 
we have that all the appearing functions in \eqref{eq_introduction} are continuous; thus the equation must be satisfied everywhere pointwise.
This concludes the proof.
\QED

%\bigskip
%
%%\smallskip
%
%\claim Proof of Theorem \ref{th_INT_regular_N}.
%Taking together Proposition \ref{prop_u_holder} %Theorem \ref{th_INT_regular_V} 
%and Proposition \ref{prop_regolar}, we get the claim.
%\QED
%
%\bigskip

\subsection{$C^1$ and $C^2$ regularity}
\label{sec_doub_C1C2}

We prove now that, under some more restrictive conditions on $s$, $\alpha$ and $\sigma$, where $f \in C^{0,\sigma}_{loc}(\R^N)$, we can prove that $u \in C^1(\R^N)$.
We notice that partial results for $s \in [\frac{1}{4}, 1)$ are already contained in Proposition \ref{prop_u_holder} %Theorem \ref{teo_u_holder} 
and Proposition \ref{prop_regolar}.
This $C^1$-regularity will be implemented then in the study of the Pohozaev identity in Section \ref{sec_Pohozaev}.

\begin{Proposition}
\label{prop_u_C1}
Assume \hyperref[(F1s)]{\textnormal{(F1)}}-\hyperref[(F2s)]{\textnormal{(F2)}}.
%Assume \textnormal{(f1)-(f2)} and $f\in C^{0,\sigma}_{loc}(\R)$ for some $\sigma \in (0,1)$. 
Let $u\in H^s(\R^N) \cap L^{\infty}(\R^N)$ be a weak solution of \eqref{eq_introduction}. Then
\begin{itemize}
\item[i)] if $s\in (\tfrac{1}{2}, 1)$, then $u \in C^{1,\gamma}(\R^N)$ for any $\gamma \in (0, 2s-1)$.
\end{itemize}
Assume now \hyperref[(F7)]{\textnormal{(F7)}} in addition. Then
\begin{itemize}
%\item[ii)] if $s \in [\frac{1}{2}, 1)$, then $u \in C^{1,\gamma}(\R^N)$ for any $\gamma \in (0, \omega+2s-1)$, where $\omega:=\min\{\sigma, \alpha\}$;
\item[ii)] if $s \in [\frac{1}{2}, 1)$ and $\omega:=\min\{\sigma, \alpha\}\leq 2-2s$, then $u \in C^{1,\gamma}(\R^N)$ for any $\gamma \in (0, \omega+2s-1)$; if instead $\omega>2-2s$, then $u \in C^{2,\omega+2s-2}(\R^N)$;
\item[iii)] if $s \in [\frac{1}{4}, \frac{1}{2})$, $\alpha > 1-2s$ and $\sigma > \frac{1-2s}{2s}$, then $u \in C^{1,\gamma}(\R^N)$ for every $\gamma \in (0,\omega+2s-1)$, where $\omega:=\min\{2s \sigma, \alpha\}$;
\item[iv)] if $\alpha<2$ %\textcolor{blue}{$\alpha <1$} 
and $\sigma > 1-2s$,
then $u\in C^{1, \gamma}(\R^N)$ for every $\gamma \in (0,1)$.
\end{itemize}
%Moreover, if $s\in (\tfrac{1}{2}, 1)$, $\alpha >2(1-s)$ and $\sigma > 2(1-s)$, then $u \in C^{2,\omega+2s-2}(\R^N)$ where $\omega:=\min\{\sigma, \alpha\}$.
\end{Proposition}

\claim Proof.
We need to check only the fourth case.
We aim to prove
$$(I_{\alpha}*F(u))f(u) \in C^{0,\omega}_{loc}(\R^N), \quad \hbox{for some $\omega + 2s >1$}$$
in order to apply Proposition \ref{prop_reg_fraz}. 
We want to show thus that $I_{\alpha}*F(u)$ is H\"older continuous; more precisely, we will show that it belongs to $C^{0,\omega}(\R^N)$ for some $\omega$ that increases according to $\gamma$, where $u\in C^{0,\gamma}(\R^N)$, so that we can employ a bootstrap argument.

Thanks to Proposition \ref{prop_regolar}, set
\begin{align*} % Package amsmath Error: \begin{split} won't work here.
\theta_0:=& \min\{\sigma +2s, 2s\sigma+2s, \alpha +2s, 1\} \\
=& \min\{2s\sigma+2s, \alpha +2s, 1\}
\end{align*}
we have
$$u \in C^{0, \gamma}(\R^N)\cap L^{\infty}(\R^N), \quad \hbox{ for $\gamma \in (0, \theta_0)$}.$$
By composition we obtain
$$f(u) \in C^{0,\gamma}(\R^N), \quad \hbox{for $\gamma \in (0, \sigma \theta_0)$},$$
and
$$F(u) \in C^{0,\gamma}(\R^N), \quad \hbox{for $\gamma \in (0, \theta_0)$},$$
which implies, by Proposition \ref{prop_regol_riesz} (possible since $\alpha<2$) and Remark \ref{rem_DuPl} (recall that $I_{\alpha}*F(u)\in L^{\infty}(\R^N)$), that
$$I_{\alpha}*F(u) \in C^{0, \alpha + \gamma}(\R^N), \quad \hbox{for $\gamma \in (0,\min\{\theta_0, 1-\alpha\})$}$$
that is
$$I_{\alpha}*F(u) \in C^{0, \gamma}(\R^N), \quad \hbox{for $\gamma \in (0,\min\{\theta_0+\alpha, 1\})$}.$$
Since $\omega_0:=\sigma \theta_0 < \min\{\theta_0+\alpha, 1\}$, we have
$$(I_{\alpha}*F(u))f(u) \in C^{0,\gamma}(\R^N) \quad \hbox{for $\gamma \in (0, \omega_0)$}.$$
We implement now the bootstrap argument.
By Proposition \ref{prop_reg_fraz} we gain
\begin{itemize}
\item if $\theta_1:=\omega_0 + 2s > 1$, then $u \in C^{1,\gamma}(\R^N)$ for $\gamma \in (0, \theta_1-1)$;
\item if $\theta_1=\omega_0 + 2s \leq 1$, then $u \in C^{0, \gamma}(\R^N)$ for $\gamma \in (0, \theta_1)$.
\end{itemize}
In the first case, we stop. Otherwise,
$$(I_{\alpha}*F(u))f(u) \in C^{0, \omega_1}(\R^N), \quad \omega_1:= \sigma \theta_1,$$
and
\begin{itemize}
\item if $\theta_2:=\omega_1 + 2s > 1$, then $u \in C^{1,\gamma}(\R^N)$ for $\gamma \in (0, \theta_2-1)$;
\item if $\theta_2=\omega_1 + 2s \leq 1$, then $u \in C^{0, \gamma}(\R^N)$ for $\gamma \in (0, \theta_2)$.
\end{itemize}
We proceed inductively by setting
$$\parag{&\omega_i:= \sigma \theta_i,&\\&\theta_i:=\omega_{i-1} +2s,& }$$
that is
$$\theta_i = \sigma \theta_{i-1} + 2s.$$
We need to show that $\theta_i>1$ at some point. We observe that
$$\theta_i> \theta_{i-1} \iff \theta_{i-1} < \frac{2s}{1-\sigma}.$$
If for some $i$ we have
$$\theta_i \geq \frac{2s}{1-\sigma}>1$$
then we stop. Otherwise, $\theta_i $ is increasing, and thus its limit $\theta_i \to l$ satisfies
$$l = \sigma l + 2s$$
which means that $l=+\infty$ or $l= \frac{2s}{1-\sigma}>1$. This concludes the proof.
\QED

\bigskip

\begin{Remark}\label{rem_regul_refined_1}
We notice that, by assuming $F\geq 0$, we can use Proposition \ref{prop_regol_riesz} to implement a bootstrap argument (namely $\omega_0:= \min\{\sigma, \frac{N}{N+\alpha}\} \theta_0$ with the notations of the above proof) to get additional regularity for a generic $\alpha \in (0,N)$. We leave the details to the interested reader. Similar arguments can be developed by assuming
$$|F(t)-F(s)| \lesssim |t-s|^{\theta} |f(t)-f(s)|, \quad \hbox{for $t,s \in \R$}$$
for some $\theta \in (0,1]$.
\end{Remark}

\begin{Remark}\label{rem_regul_refined_2}
Let us consider $u>0$ radially symmetric decreasing and assume $f \in C^1((0,+\infty))$ with
$$|f'(t)| \lesssim |t|^{-\frac{N-\alpha}{N}} + |t|^{\frac{2\alpha}{N-\alpha}}, \quad \hbox{for $t>0$}$$
and
$$|f(t)-f(s)|\lesssim |t-s|^{\theta} |f'(t)-f'(s)|, \quad \hbox{for $t,s>0$}$$
for some $\theta \in (0,1]$, then we can refine the regularity argument of Remark \ref{rem_regul_refined_1} by exploiting some asymptotic estimates.
% thanks to the the found asymptotic decay. 
Indeed, better regularity on $(I_{\alpha}*F(u)) f(u)$ can be deduced as follows: for $x,y \in \R^N$ we have
%\begin{align*}
%\pabs{\big(I_{\alpha}*F(u)\big)(x)f(u(x)) - \big(I_{\alpha}*F(u)\big)(y)f(u(y))} \leq & \, \pabs{\big(I_{\alpha}*F(u)\big)(x) - \big(I_{\alpha}*F(u)\big)(y)} \norm{f(u)}_{\infty} + \\
%&+ |u(x)-u(y)|^{\theta} \pabs{\big(I_{\alpha}*F(u)\big)(y)}\pabs{f'(u(x))-f'(u(y))} \\
%\end{align*}
\begin{eqnarray*}
\lefteqn{\pabs{\big(I_{\alpha}*F(u)\big)(x)f(u(x)) - \big(I_{\alpha}*F(u)\big)(y)f(u(y))}} \\
 &\leq  \, \pabs{\big(I_{\alpha}*F(u)\big)(x) - \big(I_{\alpha}*F(u)\big)(y)} \norm{f(u)}_{\infty} + |u(x)-u(y)|^{\theta} \pabs{\big(I_{\alpha}*F(u)\big)(y)}\pabs{f'(u(x))-f'(u(y))} 
\end{eqnarray*}
where, by Corollary \ref{corol_stima_Riesz} and Remark \ref{rem_dec_N}
$$\pabs{\big(I_{\alpha}*F(u)\big)(y)} |u(x)|^{-\frac{N-\alpha}{N}} \leq \frac{1+|x|^{N-\alpha}}{1+|y|^{N-\alpha}} \leq C $$
whenever $|x-y|\leq 1$. Thus the H\"older regularity exponent of $(I_{\alpha}*F(u)) f(u)$ directly depends on the ones of $I_{\alpha}*F(u)$, $u$ and on $\theta$. We leave the details to the interested reader. 
\end{Remark}

Finally, we exploit the $L^1$-summability in order to further investigate the $C^2$-regularity of the solution $u$. %have the following result which states that $u\in C^2(\R^N)$.
We notice that some results are already contained in Proposition \ref{prop_regolar}, whenever $s \in (\frac{1}{2}, 1)$, with some restriction on the regularity of $f$ and on $\alpha$: for instance, if $f \in C^{0,1}_{loc}(\R)$ we need $\alpha + 2s > 2$ (e.g., $\alpha \geq 1$). We prove now that, if $f\in C^1(\R)$, then no restriction on $\alpha$ is needed.
Notice that $f \in C^1(\R)$ implies $f$ non sublinear in zero, that is \hyperref[(F6)]{\textnormal{(F6)}}. %\eqref{eq_cond_superl}.

\begin{Proposition}\label{prop_u_C2}
Assume \hyperref[(F1s)]{\textnormal{(F1)}}-\hyperref[(F2s)]{\textnormal{(F2)}} and $s \in (\frac{1}{2}, 1)$.
Let $u\in H^s(\R^N) \cap L^{\infty}(\R^N)$ be a weak solution of \eqref{eq_introduction}. Then we have
\begin{itemize}
\item if \hyperref[(F7)]{\textnormal{(F7)}} holds with $\omega:=\min\{\sigma,\alpha\}>2-2s$, then $u \in C^{2,\omega+2s-2}(\R^N)$,
\item if $f \in C^1(\R)$, then $u \in C^{2,\gamma-2}(\R^N)$ for every $\gamma< 2s+1$.
% \tr{$u \in C^{2,\gamma-2}_{loc}(\R^N)$ for every $\gamma\leq 2s+1$.}
\end{itemize}
%and 
%Let $u$ be a positive weak solution of \eqref{eq_introduction}.
%Then $u \in C^{2,\gamma}_{loc}(\R^N)$ for $\gamma \in (0, 2s-1)$.
\end{Proposition}

\medskip

\claim Proof.
We need to prove only the second point.
First we show that $I_{\alpha}*F(u)$ is in $C^1(\R^N)$. Indeed, considered $\eta\in C^{\infty}_c(\R^N)$ a smooth mollification of $\chi_{B_1}$, we have
$$\big(I_{\alpha}\eta\big)*F(u) \in C^1(\R^N)$$
since $I_{\alpha}\eta \in L^1(\R^N)$ has compact support and $u \in C^1(\R^N)$ (by Proposition \ref{prop_u_holder}), % Theorem \ref{teo_u_holder}, 
while
$$ \big(I_{\alpha}(1-\eta)\big)*F(u) \in C^1(\R^N)$$
since $I_{\alpha}(1-\eta)$ has support far from the origin and thus belongs to $ C^1_b(\R^N)$, while $F(u) \in L^1(\R^N)$ by Proposition \ref{prop_u_L1} (since $u \in L^1(\R^N) \cap L^2(\R^N) \supset L^{\frac{N+\alpha}{N}}(\R^N)$, see also Remark \ref{rem_L1}). %, thanks to Theorem \ref{...}.

In particular
$$(-\Delta)^s u = - \mu u + (I_{\alpha}*F(u))f(u) \in C^1(\R^N).$$
Since $2s>1$, we gain $u \in H^{2s}(\R^N) \hookrightarrow H^1(\R^N)$, and in particular $\partial_j u \in L^2(\R^N)\cap C^{0, \gamma}(\R^N)\subset L^{\infty}(\R^N)$ for each $j =1 \dots N$.
Moreover we have
$$\partial_j \big(I_{\alpha}*F(u)\big) = \big(I_{\alpha}\eta\big)*\big( \partial_j F(u)\big) + \big(\partial_j (I_{\alpha}(1-\eta))\big)*F(u).$$
We want to show that the derivative can be moved to $F(u)$ in the second term. Indeed, set $h:=I_{\alpha}(1-\eta)$ for brevity, and let $\phi_n$ be a cut-off function with $\phi_n\equiv 1 $ in $B_n$ and support in $B_{n+1}$; thus
%$$
%\int_{\R^N} \partial_j h(x-y)F(u(y))\phi_n(y) = \int_{\R^N} h(x-y) \partial_j F(u(y))\phi_n(y) 
%+ \int_{\R^N} h(x-y) F(u(y)) \partial_j\phi_n(y);
%$$
\begin{align*}
\int_{\R^N} \partial_j h(x-y)F(u(y))\phi_n(y) =& \int_{\R^N} h(x-y) \partial_j F(u(y))\phi_n(y) + \\
&+ \int_{\R^N} h(x-y) F(u(y)) \partial_j\phi_n(y);
\end{align*}
being $\phi_n \to 1$, $\partial_j \phi_n \to 0$ as $n\to +\infty$, and $|h|, |\partial_j h| \leq C$ together with $F(u) \in L^1(\R^N)$ and $\partial_j F(u) = f(u) \partial_j u \in L^1(\R^N)$ (notice that $f(u) \in L^2(\R^N)$ since $u \in L^2(\R^N) \cap L^{\infty}(\R^N)$), by dominated convergence theorem we reach the claim.
Thus we obtain
$$\partial_j\big((I_{\alpha}*F(u))f(u)\big)= \big(I_{\alpha}*(f(u)\partial_j u)\big)f(u) + \big(I_{\alpha}*F(u)\big) f'(u)\partial_j u.$$
Since $u\in L^{\infty}(\R^N)$ and $f'$ is continuous, we have $f'(u)$ is bounded. Thus the right hand side belongs to $L^2(\R^N) \cap L^{\infty}(\R^N)$.

If we prove that
\begin{equation}\label{eq_scambio_der}
\partial_j ((-\Delta)^s u) = (-\Delta)^s(\partial_j u)
\end{equation}
then we have
$$(-\Delta)^s (\partial_j u) =- \mu \partial_j u + \partial_j \big( (I_{\alpha}*F(u)) f( u)\big) \in L^{\infty}(\R^N);$$
by Proposition \ref{prop_reg_fraz} and again $2s>1$, we obtain that $\partial_j u \in C^{1,\gamma}(\R^N)$ for any $\gamma \in (0, 2s-1)$, which is the claim.

We deal thus with \eqref{eq_scambio_der}. Since $\partial_j ((-\Delta)^s u) \in L^2(\R^N)$, we can evaluate the Fourier transform $\mc{F}\big(\partial_j ((-\Delta)^s u) \big)$, and since $(-\Delta)^s u \in C^1(\R^N)$ we have
$$\mc{F}\big(\partial_j ((-\Delta)^s u) \big) = i \xi_j \mc{F}\big( (-\Delta)^s u\big) = i \xi_j \big(|\xi|^{2s} \mc{F}(u)\big).$$
Since $u\in L^2(\R^N) \cap C^1(\R^N)$ we obtain
$$\mc{F}\big(\partial_j ((-\Delta)^s u) \big) = |\xi|^{2s} \big(i \xi_j \mc{F}(u)\big) = |\xi|^{2s} \mc{F}(\partial_j u);$$
taking back the Fourier transform, we obtain \eqref{eq_scambio_der}.
This concludes the proof.
\QED

%%%%%%%%%%%%%%%%%%%%%%%%%%%%%%%%%%%%%%%%%%%%%%%%%%%%%%%
%%%%%%%%%%%%%%%%%%%%%%%%%%%%%%%%%%%%%%%%%%%%%%%%%%%%%%%
\section{Shape of %Qualitative properties of 
ground states} %Positivity and symmetry}
\label{sec_doubl_shape}

In this Section we exploit the regularity of the solutions gained in Proposition \ref{prop_regolar} %Theorem \ref{th_INT_regular_N} 
to deduce the following theorem concerning the sign and the symmetry of the ground state solutions.

\begin{Theorem}\label{th_INT_positiv}
Assume $N\geq 2$ and \hyperref[(F7)]{\textnormal{(F7)}} in addition to \hyperref[(F1s)]{\textnormal{(F1)}}-\hyperref[(F2s)]{\textnormal{(F2)}}. Assume moreover
\begin{itemize}
\item[\textnormal{(F8)}] \label{(F8)}
\begin{itemize}
\item[(i)] $f$ is odd or even,
\item[(ii)] $f$ has constant sign on $(0, +\infty)$.
\end{itemize}
\end{itemize}
Then every Pohozaev minimum of \eqref{eq_introduction} has strict constant sign (strictly positive or negative), is radially symmetric and decreasing.
\end{Theorem}
%We refer to Section \ref{sec_prelim} (see \eqref{eq_poh_minim}) for the exact definition of Pohozaev minimum.
This last result is obtained also for \emph{constrained} problem with fixed mass, see Remark \ref{rem_constrain} for details.

\begin{Remark}
We observe that the qualitative results in Theorem \ref{th_INT_positiv} holds also for least energy solutions, when the Pohozaev identity holds for every solution, see Section \ref{sec_Pohozaev} (see also %. This is the case for instance of \cite{SGY} dealing with functions $f\in C^1(\R)$ (see also 
\cite[Eq (6.1)]{DSS1} and \cite{SGY}). % the paper in preparation \cite{CGT5}).
\end{Remark}

This theorem extends the result in Theorem \ref{thm_MVS_cita1} to the fractional case; 
in particular, \cite{MS2} deals with the case $F$ even. %: this restriction is not only a matter of symmetry, but also of sign.
Here we address also the study of the case $F$ odd: as already highlighted in Chapter \ref{chap_choq_multi}, this case is generally less studied in literature, even if (in the nonlocal framework) this assumption makes the functional symmetric as well as the odd case.
Mathematically, $F$ odd reveals to be more challenging, since the interactions in the nonlocal term among positive and negative contributions is stronger and more difficult to manage. 

Specifically, we highlight that this result is new even in the limiting local case $s=1$, $N\geq 3$, when $F$ is odd, extending some results in \cite{MS2}.
Notice that in this framework the regularity results hold for $f$ merely continuous, and moreover every Pohozaev minimum is a least energy solution, since every solution satisfies the Pohozaev identity \eqref{eq:2.5}. 

\begin{Theorem}\label{th_INT_LOC_positiv}
Let $s=1$ and assume $N\geq 3$ and \hyperref[(F1)]{\textnormal{(F1)}}-\hyperref[(F2)]{\textnormal{(F2)}}. %with $s=1$. 
Assume moreover \hyperref[(F8)]{\textnormal{(F8)}}.
%\begin{itemize}
%\item $f$ is odd or even,
%\item $f$ has constant sign on $(0, +\infty)$.
%\end{itemize}
Then every least energy solution of \eqref{eq_Choquard_genericaF} has strict constant sign (strictly positive or negative), is radially symmetric and decreasing.
\end{Theorem}

%%%%%%%%%%%%%%%%%%%%%%%%%%%%%%%%%%%%%%%%%%%%%%%%%%%%%%%
\subsection{Positivity through fibers}
\label{sec_positivity}

We want to show now that every Pohozaev ground state has constant sign.
This result requires some additional symmetric condition on $f$.

We start by providing some trivial but useful inequalities, consequence of Lemma \ref{lem_dis_modul_1}.

\begin{Lemma}\label{lem_dis_modul}
Let $u \in H^s(\R^N)$. Then
$$\norm{(-\Delta)^{s/2} |u|}_2 \leq \norm{(-\Delta)^{s/2} u}_2.$$
Assume moreover that
\begin{itemize}
\item $f$ is odd, or
\item $f$ is even, and $F$ has constant sign on % $f$ has constant sign on 
$(0, +\infty)$,
\end{itemize}
% and that $f$ is even or odd, 
then
$$\mc{D}(|u|)\geq \mc{D}(u);$$
if $f$ is odd, equality holds. As a consequence
$$\mc{J}_{\mu}(|u|)\leq \mc{J}_{\mu}(u), \quad \mc{P}_{\mu}(|u|)\leq\mc{P}_{\mu}(u).$$
\end{Lemma}

%\begin{proof}
%By \eqref{eq_semin_gagl} we have
%\begin{align*} % Package amsmath Error: \begin{split} won't work here.
%\norm{(-\Delta)^{s/2} |u|}_2^2 =& C_{N,s} \int_{\R^{2N}} \frac{\big(|u(x)|-|u(y)|\big)^2}{|x-y|^{N+2s}} \, dx dy \\
%=& C_{N,s} \int_{\R^{2N}} \frac{|u|^2(x) + |u|^2(y) - 2|u|(x)|u|(y)}{|x-y|^{N+2s}} \, dx dy \\
%\leq& C_{N,s} \int_{\R^{2N}} \frac{u^2(x) + u^2(y) - 2u(x)u(y)}{|x-y|^{N+2s}} \, dx dy \\
%=& C_{N,s} \int_{\R^{2N}} \frac{\big(u(x)-u(y)\big)^2}{|x-y|^{N+2s}} \, dx dy = \norm{(-\Delta)^{s/2}u}_2^2,
%\end{align*}
%thus the first claim.
%Focus on the second claim: if $f$ is even, then
%\begin{align*} % Package amsmath Error: \begin{split} won't work here.
%\mc{D}(|u|) =& \int_{\R^{2N}} I_{\alpha}(x-y) F(|u(x)|) F(|u(y)|) \, dx dy \\
%=& \int_{\{u>0\} \times \{u>0\}} I_{\alpha}(x-y) F(u(x)) F(u(y)) \, dx dy - \\
%& -\int_{\{u<0\} \times \{u>0\}} I_{\alpha}(x-y) F(u(x)) F(u(y)) \, dx dy - \\
%&- \int_{\{u>0\} \times \{u<0\}} I_{\alpha}(x-y) F(u(x)) F(u(y)) \, dx dy + \\
%& + \int_{\{u<0\} \times \{u<0\}} I_{\alpha}(x-y) F(u(x)) F(u(y)) \, dx dy \\
%\geq& \int_{\{u>0\} \times \{u>0\}} I_{\alpha}(x-y) F(u(x)) F(u(y)) \, dx dy +\\
%&+ \int_{\{u<0\} \times \{u>0\}} I_{\alpha}(x-y) F(u(x)) F(u(y)) \, dx dy + \\
%&+ \int_{\{u>0\} \times \{u<0\}} I_{\alpha}(x-y) F(u(x)) F(u(y)) \, dx dy +\\
%& + \int_{\{u<0\} \times \{u<0\}} I_{\alpha}(x-y) F(u(x)) F(u(y)) \, dx dy \\
%=& \int_{\R^{2N}} I_{\alpha}(x-y) F(u(x)) F(u(y)) \, dx dy = \mc{D}(u),
%\end{align*}
%which concludes the proof, observing that equality holds if $f$ is odd.
%\end{proof}

To prove the positivity of Pohozaev ground states, we need to get information about the absolute value of the function.
This analysis is simplified when dealing with local operators $s=1$ (since $\norm{\nabla|u|}_2 = \norm{\nabla u}_2$ and the Pohozaev identity holds for every solution, see \cite{MS2}), or when dealing with local nonlinearities (since the source scales in the argument in the same way as $|u|^2$ and an equivalent minimization approach can be exploited, see \cite{BL1}), or when dealing with homogeneous nonlinearities (since another minimization approach holds, see \cite{MS0,DSS1}).
In order to implement a different approach, we start observing the following fact.

\smallskip

For every $u \in H^s(\R^N)$, $u\nequiv 0$, we define the fiber $g_u:(0,+\infty)\to \R$ as follows
$$g_u(t):= \mc{J}_{\mu}(u(\cdot/t)) = \frac{t^{N-2s}}{2} \norm{(-\Delta)^{s/2}u}_2^2 + \mu \frac{t^N}{2} \norm{u}_2^2 - \frac{t^{N+\alpha}}{2} \mc{D}(u) , \quad t \in (0,+\infty) .$$
By a straightforward computation we notice that
$$g_u'(1) = \mc{P}_{\mu}(u).$$
Since $N+\alpha > N > N-2s$ it is immediate showing that there exists a single critical point for $g_u$, that we call $\lambda(u)$, which is a global maximum. That is
$$g'_u(\lambda(u))=0, \quad g_u(\lambda(u))\geq g_u(t) \; \hbox{ for each $t \in (0,+\infty)$%for all $t>0$
}.$$
Noticed that $\lambda(u)>0$, we set
$$v:= u(\cdot/\lambda(u));$$
by the fact that $g_u(\lambda(u) t) = g_v(t)$, we obtain $g'_v(1)=0$, that is
$$\mc{P}_{\mu}(v)=0.$$
In other words, the scaling through $\lambda(u)$ brings $u$ to the Pohozaev manifold; moreover, the energy is maximized in $\lambda(u)$ all over the scaling.

\begin{Proposition}\label{prop_sign_unconstr}
Assume 
\begin{itemize}
\item $f$ is odd, or
\item $f$ is even, and $F$ has constant sign on % $f$ has constant sign on 
$(0, +\infty)$,
\end{itemize}
%$f$ odd or even 
in addition to \hyperref[(F1s)]{\textnormal{(F1)}}-\hyperref[(F2s)]{\textnormal{(F2)}}. Assume moreover \hyperref[(F7)]{\textnormal{(F7)}}.
Let $u$ be a Pohozaev minimum of \eqref{eq_introduction}.
Then $u$ has strict constant sign (strictly positive or negative).
\end{Proposition}

\claim Proof.
Since $u$ satisfies $\mc{P}_{\mu}(u)=0$, we obtain $\lambda(u)=1$ and thus
\begin{equation}\label{eq_magg_gu}
g_u(t)\leq g_u(1) \quad \hbox{ for each $t \in (0,+\infty)$}.
\end{equation}
Consider $|u|$ and $\lambda(|u|)$. Define
$$v:=|u|(\cdot/\lambda(|u|))$$
which satisfies $\mc{P}_{\mu}(v)=0$. Since $u$ is a Pohozaev minimum we obtain
$$\mc{J}_{\mu}(u) \leq \mc{J}_{\mu}(v).$$
We then use Lemma \ref{lem_dis_modul} to gain
\begin{align*} % Package amsmath Error: \begin{split} won't work here.
\mc{J}_{\mu}(u) &\leq \mc{J}_{\mu}(v)\\
&= \frac{(\lambda(|u|))^{N-2s}}{2} \norm{(-\Delta)^{s/2}|u|}_2^2 + \mu \frac{(\lambda(|u|))^N}{2} \norm{|u|}_2^2 - \frac{(\lambda(|u|))^{N+\alpha}}{2} \mc{D}(|u|)\\
&\leq \frac{(\lambda(|u|))^{N-2s}}{2} \norm{(-\Delta)^{s/2}u}_2^2 + \mu \frac{(\lambda(|u|))^N}{2} \norm{u}_2^2 - \frac{(\lambda(|u|))^{N+\alpha}}{2} \mc{D}(u)\\
&= g_u(\lambda(|u|)).
\end{align*}
We finally use \eqref{eq_magg_gu} with $t=\lambda(|u|)$ and obtain
$$\mc{J}_{\mu}(u)\leq \mc{J}_{\mu}(v) \leq g_u(\lambda(|u|))\leq g_u(1) = \mc{J}_{\mu}(u).$$
Thus
$$\mc{J}_{\mu}(v)=\mc{J}_{\mu}(u)=p(\mu)$$
which, together with $\mc{P}_{\mu}(v)=0$, implies that $v$ is also a Pohozaev minimum of \eqref{eq_introduction}.
By Proposition \ref{prop_poho_min_Rsol} we obtain that $v$ is a weak solution of \eqref{eq_introduction}, positive by definition.
Thus by Proposition \ref{prop_u_bounded} %Theorem \ref{th_INT_regular_V} 
we have $v\in H^s(\R^N) \cap L^{\infty}(\R^N)$; this implies, by Proposition \ref{prop_regolar}, that $v$ is a classical solution, and in particular well defined pointwise.
Thus, if by contradiction there exists an $x_0 \in \R^N$ such that $v(x_0)=0$, then computing
$$(-\Delta)^{s} v(x_0) + \mu v(x_0)= (I_{\alpha}*F(v))(x_0)f(v(x_0))$$
we obtain, by definition of fractional Laplacian %\eqref{eq_def_lap_fraz} 
and $f(0)=0$,$$
-\int_{\R^N} \frac{v(y)}{|x_0-y|^{N+2s}} \, dy= 0$$
and hence $v\equiv 0$, which is absurd. 
Thus $|u|\neq 0$. Being $v\in L^{\infty}(\R^N)$, we obtain $u\in L^{\infty}(\R^N)$, and hence $u$ continuous by Proposition \ref{prop_u_holder}. % Theorem \ref{teo_u_holder}. 
As a consequence, $u$ does not change sign.
%Since $v$ is continuous, we have that $v$ does not change sign, and thus $u$ has constant sign. 
This concludes the proof.
%Questo argomento non vale per s=1! Serve, come chiedono Moroz e Van, che il membro di destra sia positivo (in ogni punto), da cui un principio del massimo...
\QED

\bigskip

\begin{Remark}
We point out that, without assuming \hyperref[(F7)]{\textnormal{(F7)}}, we can achieve
$$ p(\mu) = \inf \big\{ \mc{J}_{\mu}(u) \mid u \in H^s(\R^N) \setminus \{0\}, \; \mc{P}_{\mu}(u)=0, \; \hbox{$u$ positive}\big\}$$
and the same for $p_r(\mu)$. 
%Indeed, it is sufficient to argue as in the first part of Proposition \ref{prop_sign_unconstr}, but working with a minimizing sequence. 
Indeed, let $(u_n)_n \subset H^s(\R^N) \setminus \{0\}$, $\mc{P}_{\mu}(u_n) =0$, $\mc{J}_{\mu}(u_n) \to p(\mu)$ be a minimizing sequence. 
Set $v_n:=|u_n|(\cdot/\lambda(|u_n|))$ we have $\mc{P}_{\mu}(u_n) =0$ and, arguing as in the first part of Proposition \ref{prop_sign_unconstr}, we obtain
$$ \lim_{n \to +\infty} \mc{J}_{\mu}(u_n) = p(\mu) \leq \mc{J}_{\mu}(v_n) \leq g_{u_n}(\lambda(|u_n|)) \leq \mc{J}_{\mu}(u_n);$$
thus $\mc{J}_{\mu}(v_n) \to p(\mu)$, which means that $v_n$ is a positive minimizing sequence.
%\\ \tr{AGGIUNGERE DETTAGLI NELLA TESI?} % CCOMMENT NOW
\end{Remark}

%%%%%%%%%%%%%%%%%%%%%%%%%%%%%%%%%%%%%%%%%%%%%%%%%%%%%%%
\subsection{Radial symmetry}

The solution found in Theorem \ref{th_INT_exist_unconstrained} is radially symmetric by construction. We show now that, under some condition on $f$, every Pohozaev ground state is actually radially symmetric.
To this aim, we will exploit the polarization introduced in Section \ref{sec_absol_pol}.
We remark that other techniques could be investigated (with different assumptions on $f$, see e.g. \cite{LoMa08,WY0}), but this goes beyond the scope of this thesis.

\begin{Proposition}\label{prop_radial_sym}
Assume that $f$ has constant sign on $(0, +\infty)$ in addition to \hyperref[(F1s)]{\textnormal{(F1)}}-\hyperref[(F2s)]{\textnormal{(F2)}}.
Let $u$ be a positive Pohozaev minimum of \eqref{eq_introduction}.
Then $u$ is radially symmetric and decreasing with respect to some point.
\end{Proposition}

\claim Proof.
Let $u^H$ be the polarization of $u$ with respect to a closed half-space $H\subset \R^N$. By Proposition \ref{prop_prel_polar_lapl} we have
$$\norm{(-\Delta)^{s/2} u^H}_2 \leq \norm{(-\Delta)^{s/2} u}_2.$$
Assume moreover that $f \geq 0$ on $(0,+\infty)$ (if we substitute $f$ with $-f$ the Hartree-type terms are conserved). Observed that $F$ is nondecreasing on $(0,+\infty)$, we have by \eqref{eq_prel_pol_F_incr}
$$F(v^H)=(F(v))^H \quad \hbox{whenever $v \geq 0$}.$$
Thanks to these facts, we can argue as in \cite[Section 5.3]{MS2} to reach that $\mc{J}_{\mu}(u^H)=\mc{J}_{\mu}(u)$, which implies $\mc{D}(u^H) \geq \mc{D}(u)$; on the other hand, the inverse inequality is always true, and hence $\mc{D}(u^H) = \mc{D}(u)$; again by the argument in \cite{MS2} we have the claim.
%\\ \tr{INSERIRE DETTAGLI NELLA TESI?} % CCOMMENT NOW. Forse meglio di no, poiché da questi si può dedurre un metodo alternativo per la positività. Inoltre la radialità non è un risultato importante.
\QED

\medskip

\begin{Corollary}\label{cor_uguagl_p}
In the assumptions of Theorem \ref{th_INT_positiv}, every Pohozaev minimum of \eqref{eq_introduction} has constant sign, is radially symmetric and decreasing. Moreover, assuming also \hyperref[(F3s)]{\textnormal{(F3)}}-\hyperref[(F4s)]{\textnormal{(F4)}}, we have
\begin{equation*}\label{eq_ugua_Poh}
p_r(\mu) = p(\mu) = \inf \big\{ \mc{J}_{\mu}(u) \mid u \in H^s_r(\R^N) \setminus \{0\}, \; \mc{P}_{\mu}(u)=0, \; \hbox{$u$ positive}\big\}.
\end{equation*}
\end{Corollary}

\medskip

\claim Proof Theorems \ref{th_INT_positiv} and \ref{th_INT_LOC_positiv}.
The claim of Theorem \ref{th_INT_positiv} is contained in Corollary \ref{cor_uguagl_p}. The proof of Theorem \ref{th_INT_LOC_positiv} can be obtained arguing in the same way, obtaining regularity of solutions by standard results (see e.g. \cite{MS2}).
\QED

\medskip

%\begin{Remark}
%We point out that, without assuming \textnormal{(f5)}, we can achieve
%$$ p(\mu) = \inf \big\{ \mc{J}_{\mu}(u) \mid u \in H^s(\R^N) \setminus \{0\}, \; \mc{P}_{\mu}(u)=0, \; \hbox{$u$ positive}\big\}$$
%and the same for $p_r(\mu)$. Indeed, it is sufficient to argue as in the first part of Proposition \ref{prop_sign_unconstr}, but working with a minimizing sequence.
%\end{Remark}

\begin{Remark}\label{rem_constrain}
In Section \ref{sec_ground_state} we found a Mountain Pass solution $(\bar{\mu}, \bar{u})$ for the $L^2$-mass prescribed problem 
%in advanced, that is $\norm{u}_2^2=m>0$, and free (unknown) $\mu>0$; in this case 
by assuming $L^2$-subcriticality of the nonlinearity. This solution is a ground state of
$$ \mc{L}(u)=\half \int_{\R^N} |(-\Delta)^{s/2} u|^2\, dx -\half \int_{\R^N} (I_\alpha*F(u))F(u)\, dx, $$
restricted to the set
$$\mc{S}_m= \{ u \in H^s_{r}(\R^N) \mid \norm{u}_2^2=m\};$$
moreover, this solution $(\bar{\mu}, \bar{u})$ is a minimum over the Pohozaev set in the product space, that is
$$\mc{L}(\mu, u)= \inf_{\mc{P}_{\nu}(v)=0} \mc{L}(\nu, v).$$
This property easily implies that $\bar{u}$ is a ground state (in the unconstrained case) of $\mc{J}_{\bar{\mu}}$ over the Pohozaev set, that is 
%\tr{Questa proprietà si può mettere in evidenza parte} %CCOMMENT NOW
$$\mc{J}_{\bar{\mu}}(\bar{u})=p_r(\bar{\mu}).$$ 
Thus, the positivity result in Proposition \ref{prop_sign_unconstr} applies to $\bar{u}$.

Actually, all the positivity and symmetry results gained in this Section hold also for this constrained mass problem, up to simple adaptations.
Indeed, in this case the proof of the positivity is even easier, since
$$u \in \big\{\norm{v}_2^2=m \big\} \implies |u| \in \big\{\norm{v}_2^2=m\big\},$$
which means that if $u$ is a ground state, then $|u|$ is a ground state as well. 
We highlight again that this simplified approach can be not implemented in the unconstrained case.
In addition, under these symmetric and regularity assumptions \hyperref[(F7)]{\textnormal{(F7)}}-\hyperref[(F8)]{\textnormal{(F8)}}, %on $f$, 
also this $L^2$-minimum is actually an $L^2$-minimum % least energy solution 
all over the whole $H^s(\R^N)$ (and not only restricting the functional on $H^s_r(\R^N)$).
%\\ \tr{SCRIVI MEGLIO} %CCOMMENT NOW
\end{Remark}

\section{Asymptotic decay}
\label{sec_asymptotic}

In this Section we exploit the $L^1$-summability of the solutions to study the asymptotic behaviour of solutions for $|x|\to +\infty$. Recall that $2^{\#}_{\alpha}=\frac{N+\alpha}{N}$ and $2^*_{\alpha,s} = \frac{N+\alpha}{N-2s}$. % (and $2^*_{\alpha} = \frac{N+\alpha}{N-2}$).

When $s=1$ and $f(u)=|u|^{r-2} u$, that is
\begin{equation}\label{eq_MS}
-\Delta u + \mu u = \big(I_{\alpha}*|u|^r\big)|u|^{r-2} u \quad \hbox{in $\R^N$}
\end{equation}
Cingolani, Clapp and Secchi in \cite[Proposition A.2]{CCS1} obtained an exponential decay of positive solutions whenever $r\geq 2$, which means that the effect of the classical Laplacian prevails. 
 Afterwards, Moroz and Van Schaftingen in \cite{MS0} (see also \cite{MS1,MS3} and \cite{CLO0, ClSa}) extended the previous analysis in the case of ground state solutions to all the possible values of $r$ in the range $ [2^{\#}_{\alpha}, 2^*_{\alpha,1}]$, % $[\frac{N+\alpha}{N}, \frac{N+\alpha}{N-2}]$, 
in particular by finding a polynomial decay when $f$ is sublinear (i.e., the Choquard term effect prevails). 
They prove the following result \cite[Theorem 4]{MS0}.

\begin{Theorem}[\cite{MS0}]\label{thm_MVS}
Let $s=1$ and let $u \in H^1(\R^N)$ be a nonnegative ground state of \eqref{eq_MS}, and $r \in [2^{\#}_{\alpha}, 2^*_{\alpha,1}]$. %[\frac{N+\alpha}{N}, \frac{N+\alpha}{N-2}]$. 
Assume $\mu=1$. 
Then
\begin{itemize}
\item if $r>2$, then
$$\lim_{|x|\to +\infty} u(x) |x|^{\frac{N-1}{2}} e^{|x|} \in (0, +\infty);$$
\item if $r=2$, then 
$$\lim_{|x|\to +\infty} u(x) |x|^{\frac{N-1}{2}} e^{\int_{\nu}^{|x|} \sqrt{1- \frac{\nu^{N-\alpha}}{t^{N-\alpha}}} dt } \in (0, +\infty)$$
for some explicit $\nu=\nu(u)$;
\item if $r<2$, then
$$ \lim_{|x|\to +\infty} u(x) |x|^{\frac{N-\alpha}{2-r}} = C(N, \alpha, r, u) \in (0, +\infty)$$
where 
\begin{equation}\label{eq_costant_stima}
C(N,\alpha, r, u) := \big( C_{N,\alpha} \norm{u}_r^r\big)^{\frac{1}{2-r}}
\end{equation}
with $C_{N,\alpha}:=\frac{\Gamma(\frac{N-\alpha}{2})}{2^{\alpha} \pi^{N/2} \Gamma(\frac{\alpha}{2})}$.
\end{itemize}
\end{Theorem}

Notice that, when $\mu\neq 1$, $\mu$ influences both the limiting constants and the speed of the exponential decays. 
We refer also to \cite[Section 8.2]{DG0} for some results on convolution equations with non-variational structure. 

\smallskip

The case of the fractional Choquard equation $s \in (0,1)$ with homogeneous $f$, that is
\begin{equation}\label{eq_fr_ch_hom}
(-\Delta)^s u + \mu u = \big(I_{\alpha}*|u|^r\big)|u|^{r-2} u \quad \hbox{in $\R^N$},
\end{equation}
 has been studied by D'Avenia, Siciliano and Squassina in \cite{DSS1} (see also \cite{DSS2} and \cite{MZng0, ZY0} for other related results). 
In this paper the authors gain existence of ground states, multiplicity and qualitative properties of solutions. 
In particular they obtain asymptotic decay of solutions whenever the source is linear or superlinear, that is when $r\geq 2$ (see also \cite{BBMP} for the $p$-fractional Laplacian counterpart): in this case the rate is polynomial, as one can expect dealing with the fractional Laplacian; more specifically, it does not depend on $\alpha$, and they prove the following theorem. % (see also \cite[(1.12) and Theorem 1.4]{Gre0} for other results on asymptotic decays).

\begin{Theorem}[\cite{DSS1}]\label{thm_DSS1}
Let $u\in H^s(\R^N)$ be a solution of \eqref{eq_fr_ch_hom}, and assume $r \in [2, 2^*_{\alpha,s}]$. %\frac{N+\alpha}{N-2s}]$. 
Then
$$0 < \liminf_{|x|\to +\infty} |u(x)| |x|^{N+2s} \leq \limsup_{|x|\to +\infty} |u(x)| |x|^{N+2s} <+\infty.$$
\end{Theorem}

In this Section we study the asymptotic profile of solutions of equation \eqref{eq_introduction}, starting by 
%first we focus on 
the case $f$ linear or superlinear. 
In the remeaning part we will develop the more tricky case of $f$ sublinear: the found decay is of polynomial type, with a rate possibly slower than $\sim\frac{1}{|x|^{N+2s}}$; %
%T
the result is new even for homogeneous functions $f(u)=|u|^{r-2}u$, $r\in [\frac{N+\alpha}{N},2)$, and, 
differently from the local case $s=1$ in \cite{MS0}, new phenomena arise connected to a new $s$-sublinear threshold that we detect on $r$.

\smallskip

This Section is mainly based on papers \cite{CGT3} and \cite{Gal1}.

\medskip

We show first some conditions which imply the decay at infinity of the solutions. 
%\tr{forse da spostare prima? Vedi meglio} %COMMENT NOW

\begin{Lemma}\label{lem_u_to0}
Assume that \hyperref[(F1s)]{\textnormal{(F1)}}-\hyperref[(F2s)]{\textnormal{(F2)}} hold. Let $u \in H^s(\R^N)$ be a weak solution of \eqref{eq_introduction}. Assume 
%$$u \in L^{\frac{N}{2s}}(\R^N)\cap L^{\infty}(\R^N)$$
%and
$$(I_{\alpha}*F(u))f(u)\in L^2(\R^N) \cap %L^{\frac{N}{2s}}(\R^N)\cap 
L^{\infty}(\R^N).$$
 Then we have
$$%\begin{equation*}%\label{eq_conv_unif_0_veps2}
u(x) \to 0 \quad \hbox{as $|x|\to +\infty$}.
$$%\end{equation*}
\end{Lemma}

\claim Proof.
Being $u$ solution of
%$$(-\Delta)^s u + u = (1-\mu) u + \big(I_{\alpha}*F(u)\big) f(u) =: \chi \quad \hbox{in $\mathbb{R}^N$},$$
$$(-\Delta)^s u + \mu u = \big(I_{\alpha}*F(u)\big) f(u) =: \chi \quad \hbox{in $\mathbb{R}^N$},$$
where $\chi \in L^2(\R^N) \cap L^{\frac{N}{2s}}(\R^N) \cap L^{\infty}(\R^N)$, 
we have the representation formula (being $\chi \in L^2(\R^N)$)
$$ u= \mc{K} * \chi $$
where $\mc{K}:=\mc{K}_{2s, \mu}$ is the Bessel kernel; we recall that $\mc{K}$ is positive, it satisfies $\mc{K}(x) \leq \frac{C}{|x|^{N+2s}}$ for $|x| \geq 1$ and $\mc{K} \in L^q(\R^N)$ for $q \in [1, 1 + \tfrac{2s}{N-2s})$ (see Lemma \ref{lem_Bessel_kernel}). 
Let us fix $\eta>0$; 
we have, for $x \in \R^N$,
\begin{align*}
u(x) =& \int_{\R^N} \mc{K}(x-y) \chi(y) dy \\
=& \int_{|x-y|\geq 1/\eta} \mc{K}(x-y) \chi(y)dy +\int_{|x-y|< 1/\eta} \mc{K}(x-y) \chi(y)dy.
\end{align*}
As regards the first piece
$$ \int_{|x-y|\geq 1/\eta} \mc{K}(x-y) \chi(y)dy \leq \norm{\chi}_{\infty} \int_{|x-y|\geq 1/\eta} \frac{C}{|x-y|^{N+2s}} dy \leq C \eta^{2s}$$
while for the second piece, 
 fixed a whatever $q \in (1, \min\{2,\frac{N}{N-2s}%1 + \tfrac{2s}{N-2s}
\})$ and its conjugate exponent $q' \in ( \max\{2,\frac{N}{2s}\}, + \infty)$
%fixed a whatever $q \in (1, 1 + \tfrac{2s}{N-2s})$ and its conjugate exponent $q'> \frac{N}{2s}$, 
we have by H\"older inequality
$$
\int_{|x-y|< 1/\eta} \mc{K}(x-y) \chi(y)dy \leq \norm{\mc{K}}_q \norm{\chi}_{L^{q'}(B_{1/\eta}(x))}\\
$$
where the second factor can be made small for $|x| \gg 0$. 
Joining the pieces, we conclude the proof. %have \eqref{eq_conv_unif_0_veps2}. 
\QED

\bigskip

Notice that $u \in L^2(\R^N)\cap L^{\infty}(\R^N)$ implies the assumptions of Lemma \ref{lem_u_to0}. % (see Remark \ref{rem_L1}).

%Notice that $u \in L^{\frac{N+2\alpha}{2s+\alpha}}(\R^N)\cap L^{\infty}(\R^N)$ implies the assumptions of Lemma \ref{lem_u_to0} (see Remark \ref{rem_L1}).
%If $N \geq 4s$, then $u \in L^2(\R^N) \cap L^{\infty}(\R^N) \subset L^{\frac{N+2\alpha}{2s+\alpha}}(\R^N)$. %\cap L^{\infty}(\R^N)$. 
%Differently, if $N \in (2s, 4s)$, we need to exploit the extra assumption $u \in L^1(\R^N)$. 

%%%%%%%%%%%%%%%%%%
\subsection{The (super)linear case}
\label{sec_doub_superl}

%In particular, 
By assuming the condition in zero \hyperref[(F6)]{\textnormal{(F6)}} for the function $f$, we obtain % assumed in leads to 
the following polynomial decay, as stated in paper \cite{CGT3}. % of the solutions.
\begin{Theorem}\label{th_INT_decay}
Assume \hyperref[(F1s)]{\textnormal{(F1)}}-\hyperref[(F2s)]{\textnormal{(F2)}} and \hyperref[(F6)]{\textnormal{(F6)}}. Let $u\in H^s(\R^N)$ be a positive weak solution of \eqref{eq_introduction}. Then there exists $C', C''>0$ such that
$$\frac{C'}{1+|x|^{N+2s}} \leq u(x)\leq \frac{ C''}{1+|x|^{N+2s}},\quad \textit{ for $x \in \R^N$}.$$
\end{Theorem}

\bigskip

%We observe that the assumptions of the Lemma are fulfilled by assuming that $u$ is bounded thanks to Proposition \ref{prop_u_L1}. 
We are now ready to prove the polynomial decay of the solutions.

\medskip

\claim Proof of Theorem \ref{th_INT_decay}.
%Observe that, by \hyperref[(F6)]{\textnormal{(F6)}} and Lemma \ref{lem_u_to0}, we have 
%\begin{equation}\label{eq_dim_V_lim}
%\frac{f(u)}{u} \in L^{\infty}(\R^N).
%\end{equation}
%Moreover, by Theorem \ref{th_INT_regular} and Proposition \ref{prop_u_L1} we have $u \in L^1(\R^N) \cap L^{\infty}(\R^N) \subset L^{\frac{N}{2s}}(\R^N) \cap L^{\frac{N+2\alpha}{2s+\alpha}}(\R^N)\cap L^{\infty}(\R^N)$, which in particular implies $(I_{\alpha}*F(u))f(u)\in L^{\frac{N}{2s}}(\R^N)\cap L^{\infty}(\R^N)$ (see Remark \ref{rem_L1}). 
%Thus we can apply Proposition \ref{prop_conv_C0} and obtain 
Observe that, by \hyperref[(F6)]{\textnormal{(F6)}} and Theorem \ref{th_INT_regular} % Lemma \ref{lem_u_to0}, we have 
\begin{equation}\label{eq_dim_V_lim}
\frac{f(u)}{u} \in L^{\infty}(\R^N).
\end{equation}
Moreover, by 
%Theorem \ref{th_INT_regular} and Proposition \ref{prop_u_L1} we have $u \in L^1(\R^N) \cap L^{\infty}(\R^N) \subset L^{\frac{N}{2s}}(\R^N) \cap L^{\frac{N+2\alpha}{2s+\alpha}}(\R^N)\cap L^{\infty}(\R^N)$, which in particular implies $(I_{\alpha}*F(u))f(u)\in L^{\frac{N}{2s}}(\R^N)\cap L^{\infty}(\R^N)$ (see Remark \ref{rem_L1}). 
%Thus we can apply 
Proposition \ref{prop_conv_C0} we obtain 
\begin{equation}\label{eq_conv_inf_block}
(I_{\alpha}*F(u))(x) \frac{f(u(x))}{u(x)} \to 0 \quad \hbox{ as $|x| \to + \infty$}. 
\end{equation}
As a consequence, by \eqref{eq_conv_inf_block} and the positivity of $u$, we have for some $R'\gg 0$
$$(-\Delta)^s u + \tfrac{1}{2} \mu u = (I_{\alpha}*F(u))f(u) - \tfrac{1}{2} \mu u = \left( (I_{\alpha}*F(u))\tfrac{f(u)}{u} - \tfrac{1}{2} \mu \right) u \leq 0 \quad \hbox{in $\mathbb{R}^N\setminus B_{R'}$}.$$
Similarly
$$(-\Delta)^s u + \tfrac{3}{2} \mu u = (I_{\alpha}*F(u))f(u) + \tfrac{1}{2} \mu u = \left( (I_{\alpha}*F(u))\tfrac{f(u)}{u} + \tfrac{1}{2} \mu \right) u \geq 0 \quad \hbox{in $\mathbb{R}^N\setminus B_{R'}$}.$$
Notice that we always intend differential inequalities in the weak sense.
%, that is tested with functions in $H^s(\R^N)$ with supports contained in the reference domain (e.g., $\R^N \setminus B_{R'}$). 

In addition, by Lemma \ref{lem_esist_sol_part} we have that there exist two positive functions $\underline{W}'$, $\overline{W}'$ and three positive constants $R''$, $C'$ and $C''$ depending only on $\mu$, such that
$$ \parag{
& (-\Delta)^s \underline{W}' + \frac{3}{2}\mu \, \underline{W}' = 0 \quad \hbox{in $\mathbb{R}^N \setminus B_{R''}$},& \\ 
&\frac{C'}{|x|^{N+2s}}< \underline{W}' (x), \quad \hbox{ for $|x|>2R''$}.&}$$
and
$$ \parag{& (-\Delta)^s \overline{W}' + \frac{1}{2}\mu \, \overline{W}' = 0 \quad \hbox{in $\mathbb{R}^N \setminus B_{R''}$},& \\ 
& \overline{W}'(x) < \frac{C''}{|x|^{N+2s}}, \quad \hbox{ for $|x|>2R''$}.&}$$
Set $R:=\max\{ R', 2R''\}$. Let $\underline{C}_1$ and $\overline{C}_1$ be some lower and upper bounds for $u$ on $B_R$, $\underline{C}_2:=\min_{B_R} \overline{W}'$ and $\overline{C}_2:= \max_{B_R} \underline{W}'$, all strictly positive. Define
$$\underline{W}:= \underline{C}_1 \overline{C}_2 ^{-1} \underline{W}', \quad \overline{W}:= \overline{C}_1 \underline{C}_2^{-1} \overline{W}'$$
so that
$$\underline{W}(x)\leq u(x) \leq \overline{W}(x), \quad \hbox{ for $|x|\leq R$}.$$
Thanks to the comparison principle in Lemma \ref{lem_comp_prin}, and redefining $C'$ and $C''$, we obtain
$$ \frac{C'}{|x|^{N+2s}} <\underline{W}(x) \leq u(x) \leq \overline{W}(x) < \frac{C''}{|x|^{N+2s}}, \quad \hbox{ for $|x|>R$}.$$
By the boundedness of $u$, we obtain the claim.
\QED

\bigskip

We see that, for non sublinear $f$ (that is, \hyperref[(F6)]{\textnormal{(F6)}}), the decay is essentially given by the fractional operator. It is important to remark that, contrary to the limiting local case $s=1$ (Theorem \ref{thm_MVS}), the Choquard term in case of linear $f$ seems not to affect the decay of the solution.

\begin{Remark}
We observe that the conclusion of the proof of Theorem \ref{th_INT_decay} can be substituted %(and extended to sign changing solutions) 
by exploiting the results in \cite{FLS} through a Kato's inequality (see also \cite[Theorem 3.2]{Amb-1}). Indeed write $V:= -(I_{\alpha}*F(u)) \frac{f(u)}{u}$, which is bounded and zero at infinity as observed in \eqref{eq_dim_V_lim}--\eqref{eq_conv_inf_block}, and gain
$$(-\Delta)^s u + V(x) u = - \mu u \quad \hbox{in $\mathbb{R}^N$}. $$
Up to dividing for $\norm{u}_2$, we may assume $\norm{u}_2=1$. Thus we are in the assumptions of \cite[Lemma C.2]{FLS} and obtain, for constant sign or sign-changing solutions of \eqref{eq_introduction},
$$|u(x)| \leq \frac{C_1}{(1 + |x|^2)^{\frac{N+2s}{2}}}$$
together with
$$|u(x)| = \frac{C_2}{|x|^{N+2s}} + o\left( \frac{1}{|x|^{N+2s}}\right) \quad \hbox{ as $|x| \to +\infty$}$$
for some $C_1, C_2>0$.
\end{Remark}

%\tr{L'ARTICOLO SOTTOMESSO LO INSERISCO MEGLIO ALLA FINE (aspetto eventuali risposte dai referee)} %COMMENT NOW

\subsection{The sublinear case: fractional Laplacian versus Riesz potential} %introduction} %JDE: avoiding a detailed literature survey or a summary of the results.
\label{sec_sublin_case}

We focus now on the case $f$ sublinear: 
we aim to study the fractional Choquard case $s\in (0,1)$, $\alpha \in (0,N)$, in presence of general, sublinear nonlinearities. 
We point out that the arguments in \cite{MS0} cannot be directly adapted to the fractional framework: for instance, we see that the explicit computation of the fractional Laplacian of some comparison function is not possible, and the choice of the comparison functions itself is hindered by some growth condition typical of the nonlocal framework; moreover, it is not obvious that all the weak solutions are pointwise solutions, and neither one can deduce that the concave power of a pointwise solution is indeed a solution (of a different equation) itself.

\smallskip

We start by presenting the case of homogeneous powers $f$, which has an interest on its own.
Since in the superlinear case the rate of convergence is of the type $\sim \frac{1}{|x|^{N+2s}}$, in the sublinear case we generally expect a slower decay. 
Actually this is what we find, as the following theorem states.

\begin{Theorem}\label{thm_homog0_ws}
Let $u\in H^s(\R^N)$, strictly positive, radially symmetric and decreasing, be a weak solution of \eqref{eq_fr_ch_hom}.
 Let $r \in [2^{\#}_{\alpha}, %\frac{N+\alpha}{N}, 
2)$ and set 
$$\beta:= \min\left\{ \frac{N-\alpha}{2-r}, N+2s\right\} \geq N.$$
Then
$$0 < \liminf_{|x| \to +\infty} u(x)|x|^{\beta} \leq \limsup_{|x|\to +\infty} u(x) |x|^{\beta} < + \infty .$$
%$$\liminf_{|x| \to +\infty} u(x)|x|^{\beta} \in (0, +\infty).$$
%If in addition $\mu> r-1$, then 
%$$\limsup_{|x|\to +\infty} u(x) |x|^{\beta} \in (0, +\infty).$$
\end{Theorem}

We refer to Remark \ref{rem_more_gen} and Corollary \ref{corol_homogen2} for some comments and generalizations on the assumptions. 
%Thanks to the qualitative analysis of the previous Sections, 
This result in particular applies to ground states solutions. % (see Definition \ref{def_pohozaev_min2}). 

\begin{Corollary}\label{corol_homog0_gs}
Let $u$ be a positive ground state of \eqref{eq_fr_ch_hom}. 
Then the conclusions of Theorem \ref{thm_homog0_ws} hold.
\end{Corollary}

We highlight that the found decay of the ground states might give information, when $r<2$, also on the twice Gateaux differentiability of the corresponding functional and on the nondegenaracy of the ground state solution itself, see \cite{MS0} (see also \cite[Section 3.3.5]{MS3}). 
Moreover %We further highlight that 
this information on the decay may be exploited to study fractional Choquard equations with potentials $V=V(x)$ approaching, as $|x|\to +\infty$, some $V_{\infty}>0$ from above or oscillating, in the spirit of \cite{MPS1, MPS2}. It might be further used, for example, in the semiclassical analysis of concentration phenomena, see e.g. Chapter \ref{chap_concentr}.

\smallskip

Joining the results in Theorem \ref{thm_DSS1} and Theorem \ref{thm_homog0_ws} we obtain the following picture of the asymptotic decay of fractional Choquard equations.
\begin{Corollary}
Let $u$ be a positive ground state of \eqref{eq_fr_ch_hom}, with $r \in [2^{\#}_{\alpha}, 2^*_{\alpha,s}]$ %[ \frac{N+\alpha}{N}, \frac{N+\alpha}{N-2s}]$ 
and $\mu>r-1$. 
\begin{itemize}
\item If $r \in [2^{\#}_{\alpha}, %\frac{N+\alpha}{N}, 
\frac{N+\alpha+4s}{N+2s}]$, then 
$$0 < \liminf_{|x|\to +\infty} u(x) |x|^{\frac{N-\alpha}{2-r}} \leq \limsup_{|x|\to +\infty} u(x) |x|^{\frac{N-\alpha}{2-r}} <+\infty;$$
in particular, $\frac{N-\alpha}{2-r}=N$ in the lower critical case $r=2^{\#}_{\alpha}$. % \frac{N+\alpha}{N}$.
\item If $r \in [\frac{N+\alpha+4s}{N+2s}, 2^*_{\alpha,s}]$, %\frac{N+\alpha}{N-2s}]$, 
then
$$0 < \liminf_{|x|\to +\infty} u(x) |x|^{N+2s} \leq \limsup_{|x|\to +\infty} u(x) |x|^{N+2s} <+\infty.$$
\end{itemize}
\end{Corollary}
By the previous Corollary we see that the exponent 
$$r^*_{\alpha,s}:= \frac{N+\alpha+4s}{N+2s},$$
$r^*_{\alpha,s} \in (2^{\#}_{\alpha}, %\frac{N+\alpha}{N}, 
2)$, separates the cases where the fractional Laplacian influences more the rate of convergence (which does not depend on $\alpha$), from the cases where the asymptotic behaviour is dictated by the Choquard term (which does not depend on $s$). 
This phenomenon seems to highlight a difference between the fractional and the local case, where the separating exponent is $r=2$ (see Theorem \ref{thm_MVS}): indeed, when $r \in (r^*_{\alpha,1}, 2)$, the arbitrary big (as $r \to 2$) polynomial behaviour $\sim \frac{1}{|x|^{\frac{N-\alpha}{2-r}}}$ keeps being slower than the exponential decay induced by the classical Laplacian; this is not the case when compared with the polynomial decay induced by the fractional Laplacian, and this is why this new phenomenon appears in this range. 
Thus $r^*_{\alpha,s}$ can be seen as a kind of \emph{$s$-subquadratic} threshold for the growth of $F$; 
set instead
$$p^*_{\alpha,s}:=r^*_{\alpha,s}-1 = \frac{\alpha+2s}{N+2s},$$
it can be seen as a \emph{$s$-sublinear} threshold for the growth of $f$. 
Notice that
$$ r^*_{\alpha,s} \stackrel{s \to 0} \to 2^{\#}_{\alpha}, %\frac{N+\alpha}{N}, 
\quad r^*_{\alpha,s} \stackrel{\alpha \to N} \to 2,$$
while
$$ r^*_{\alpha,s} \stackrel{s \to 1} \to \frac{N+\alpha+4}{N+2} \in (2^{\#}_{\alpha}, %\frac{N+\alpha}{N}, 
2), \quad r^*_{\alpha,s} \stackrel{\alpha \to 0} \to \frac{N+4s}{N+2s} \in (1,2).$$
%$$ p^*_{\alpha,s} \stackrel{s \to 0} \to \frac{\alpha}{N}, \quad p^*_{\alpha,s} \stackrel{\alpha \to N} \to 1,$$
%while
%$$ p^*_{\alpha,s} \stackrel{s \to 1} \to \frac{\alpha+2}{N+2} \in \Big(\frac{\alpha}{N}, 1\Big), \quad p^*_{\alpha,s} \stackrel{\alpha \to 0} \to \frac{2s}{N+2s} \in (0,1).$$
We refer also to the recent paper \cite[Theorem 1.4]{Gre0} where asymptotic decay results are studied in a different framework (still involving the fractional Laplacian and the Riesz potential); here a threshold different from the classical case $s=1$ is detected as well.

\smallskip

When %$\mu=1$ and 
$r \in [2^{\#}_{\alpha}, %\frac{N+\alpha}{N}, 
r^*_{\alpha,s})$ we are also able to find a \emph{sharp decay} for $u$.
\begin{Corollary}\label{corol_sharp_decay}
Let $u\in H^s(\R^N)$, strictly positive, radially symmetric and decreasing, be a weak solution of \eqref{eq_fr_ch_hom}; in particular, $u$ may be a ground state.
If $r \in [2^{\#}_{\alpha}, %\frac{N+\alpha}{N}, 
r^*_{\alpha,s})$, %and $\mu=1$, 
we have
$$ \lim_{|x|\to +\infty} u(x) |x|^{\frac{N-\alpha}{2-r}} = \left( \frac{C_{N,\alpha}\norm{u}_r^r}{\mu} \right)^{\frac{1}{2-r}}; %C(N, \alpha, r, u)
$$
%where $C(N,\alpha,r,u)$ is defined in \eqref{eq_costant_stima}.
notice that, if $\mu=1$, the constant is coherent with \eqref{eq_costant_stima}.
\end{Corollary}
%We notice that, if $\mu=1$, the constant in \eqref{eq_corol_sharp_decay} is coherent with \eqref{eq_costant_stima}. 
We finally highlight that, for $s \in (0,1]$, the rate of convergence of the solutions for $r \leq r^*_{\alpha,s}$ is $\sim \frac{1}{|x|^{\frac{N-\alpha}{2-r}}}$: for bigger values of $r$, the rate stabilizes to $\sim \frac{1}{|x|^{N+2s}}$ when $s<1$, while it keeps getting faster when $s=1$ (up to the threshold $r=2$, where it gets constantly exponential). 
It might be interesting to investigate other possible phenomena on fractional Choquard equations when $r$ is above and below this exponent $r^*_{\alpha,s}$, or also possible phenomena in $(r^*_{\alpha,1},2)$ for the local Choquard equation.
%and in $(1, r^*_{s,0}+p^*_{s,0}=\frac{N+6s}{N+2s})$ for the fractional local equation.

\begin{Remark}
We notice that, fixed a positive solution $u$, by setting
$$\rho:= I_{\alpha}* u^r$$
equation \eqref{eq_fr_ch_hom} can be rewritten as
$$(-\Delta)^s u + \mu u = \rho(x) u^{r-1}.$$
When $\mu=0$ and $\rho(x) \leq \frac{1}{|x|^{\gamma}}$ with $\gamma > N$, this fractional sublinear equation ($r\in (0,2)$) has been studied in \cite{PT0} (see also \cite[Theorem 4.4]{GM0} where they extend the result to $\gamma>2s$): here the authors find an estimate from above of the asymptotic decay of the solutions, which is strictly slower than $\sim \frac{1}{|x|^N}$. 
Notice that, in our case, $\rho= I_{\alpha}*u^r$ decays at most as $\sim \frac{1}{|x|^{N-\alpha}}$ (see \eqref{eq_decay_Riesz}), 
and we discuss the strict positive mass case $\mu>0$.
%strongly rely on the strictly positivity of $\mu$. 
See also \cite{DSS1, Le0} for more results on the zero mass case.
\end{Remark}

\medskip

%eq_introduction = eq_main_intr

We pass now to more general nonlinearities, and study \eqref{eq_introduction}. 
%For the whole paper we assume the following conditions on $f$ in order to give sense to appearing integrals: 
%\begin{itemize}
%\item[\textnormal{(f1)}] \label{(f1)}
% $f= F' \in C(\R, \R)$;
%\item[\textnormal{(f2)}] \label{(f2)}
%$f$ satisfies
% $$i) \; \limsup_{t \to 0} \frac{|tf(t)|}{|t|^{\frac{N+ \alpha}{N}}} <+\infty, \quad
% ii) \; \limsup_{ |t| \to + \infty} \frac{|t f(t)|}{|t|^{\frac{N+ \alpha}{N-2s}}} <+\infty,$$
%or equivalently there exists $C >0$ such that for every $t \in \R$, 
%$$|t f(t)| \leq C \big(|t|^{\frac{N + \alpha}{N}} + |t|^{\frac{N+ \alpha}{N-2s}}\big).$$
%\end{itemize}
%In particular, \hyperref[(f2)]{\textnormal{(f2)}} implies
We will assume \hyperref[(F1s)]{\textnormal{(F1)}}-\hyperref[(F2s)]{\textnormal{(F2)}}, which in particular imply
 \begin{equation}\label{eq_condiz_F}
i) \; \limsup_{t \to 0} \frac{|F(t)|}{|t|^{2^{\#}_{\alpha}}} <+\infty, \quad
 ii) \; \limsup_{ |t| \to + \infty} \frac{|F(t)|}{|t|^{2^*_{\alpha,s}}} <+\infty,
\end{equation}
or equivalently that there exists $C >0$ such that for every $t \in \R$, 
$$|F(t)| \leq C \big(|t|^{2^{\#}_{\alpha}} + |t|^{2^*_{\alpha,s}}\big).$$

In addition we consider $f$ \emph{sublinear in the origin}, given by the following assumptions:

\begin{itemize}
\item[\textnormal{(F9)}] \label{(F9)}
 there exists $r \in [2^{\#}_{\alpha}, 2)$ such that
$$\limsup_{t \to 0^+} \frac{|f(t)|}{t^{r-1}} \in [0, +\infty),$$
i.e., for some $\bar{C}>0$ and $\delta \in (0,1)$ we have 
\begin{equation}\label{eq_cond_2sublin_f}
|f(t)| \leq \bar{C} t^{r-1} \quad \hbox{for $t\in (0,\delta)$};
\end{equation}
\item[\textnormal{(F10)}] \label{(F10)} 
there exists $r \in [2^{\#}_{\alpha}, 2)$ such that
$$\liminf_{t \to 0^+} \frac{f(t)}{t^{r-1}} \in (0, +\infty),$$
i.e., for some $\underline{C}>0$ and $\delta \in (0,1)$ we have 
\begin{equation}\label{eq_cond_4sublin_f}
f(t) \geq \underline{C} t^{r-1} \quad \hbox{for $t \in (0,\delta)$}.
\end{equation}
\end{itemize}

\noindent
A sufficient condition for \hyperref[(F9)]{\textnormal{(F9)}} 
%\hyperref[(f4)]{(f4)} \nameref{(f3)}{(f3)} \href{(f3)}{(f3)} \ref{(f3)} \pageref{(f3)} \textnormal{(f3)} 
is clearly given by
\begin{equation}\label{eq_strongf3}
\limsup_{t \to 0^+} \frac{f(t)}{t^{r-1}}=0 \quad \hbox{for some $r \in [2^{\#}_{\alpha}, 2)$},
\end{equation}
which means that $\bar{C}$ can be taken arbitrary small in \eqref{eq_cond_2sublin_f}; in particular it includes logarithmic nonlinearities $f(t)=t \log(t^2)$, %(see \cite{GalSch}), 
where $r$ can be chosen arbitrary close to $2$. 
A sufficient condition for \hyperref[(F10)]{\textnormal{(F10)}} is instead given (for example) by a local Ambrosetti-Rabinowitz condition of the type
$$f(t) t \geq r F(t) >0 \quad \hbox{for $t \in (0,\delta)$}.$$
The restriction in \hyperref[(F9)]{\textnormal{(F9)}} and \hyperref[(F10)]{\textnormal{(F10)}} to right neighborhoods of zero is due to the fact we deal with positive solutions.

\smallskip

We eventually come up with the following generalization of Theorem \ref{thm_homog0_ws}.
\begin{Theorem}\label{thm_main}
Assume \hyperref[(F1s)]{\textnormal{(F1)}}-\hyperref[(F2s)]{\textnormal{(F2)}}, and let $u\in H^s(\R^N)$, strictly positive, radially symmetric and decreasing, be a weak solution of \eqref{eq_introduction}. Let $r \in [2^{\#}_{\alpha}, %1+ \frac{\alpha}{N}, 
2)$ and set
$$\beta:= \min\left\{ \frac{N-\alpha}{2-r}, N+2s\right\} \geq N.$$
\begin{itemize}
\item[\textnormal{(i)}] Assume \hyperref[(F9)]{\textnormal{(F9)}}. %and $\mu> (r-1) \bar{C}^{\frac{1}{r-1}}$. 
Then 
$$\limsup_{|x|\to +\infty} u(x) |x|^{\beta} \in (0, +\infty).$$
\item[\textnormal{(ii)}] Assume \hyperref[(F10)]{\textnormal{(F10)}}, $f$ locally H\"older continuous and $\int_{\R^N}F(u)>0$ (e.g. $F\geq 0$ on $(0,+\infty)$). 
Then
$$\liminf_{|x| \to +\infty} u(x)|x|^{\beta} \in (0, +\infty).$$
\end{itemize}
If both conditions in (i) and (ii) hold, together with $\overline{C}=\underline{C}$ (i.e., $f$ is a power near the origin) and $r \in [\frac{N+\alpha}{N}, \frac{N+\alpha+4s}{N+2s})$, then %when $\beta < N+2s$ 
we have the sharp decay
\begin{equation}\label{eq_sharp_decay_F} 
%\lim_{|x|\to +\infty} u(x) |x|^{\frac{N-\alpha}{2-r}} = \left( \frac{\underline{C} C_{N,\alpha} \int_{\R^N} F(u)}{\mu } \right)^{\frac{1}{2-r}} 
\lim_{|x|\to +\infty} u(x) |x|^{\frac{N-\alpha}{2-r}} = \left( \frac{C_{N,\alpha} \big(\lim_{t\to 0^+} \frac{f(t)}{t^{r-1}}\big) \int_{\R^N} F(u)}{\mu } \right)^{\frac{1}{2-r}} 
\end{equation}
where $C_{N,\alpha}>0$ is given in \eqref{eq_costant_stima}.
\end{Theorem}
\begin{Remark}\label{rem_more_gen}
We highlight that the conclusions of Theorem \ref{thm_main} (as well as of Theorem \ref{thm_homog0_ws}) hold in more general cases. Indeed:
\begin{itemize}
%\item as regards (i), if $r \in (2^{\#}_{\alpha}, %1+ \frac{\alpha}{N},
%2)$ we can drop the condition on $\mu$ by possibly considering a slower decay: see Corollary \ref{corol_slower_deca2};
\item The case
%$$\lim_{t\to 0} \frac{|f(t)|}{|t|}=+\infty$$
%in a non-strict sense (for example $f(t) \sim t \log(t^2)$ near zero, is included, and as we expect the decay is of order $\sim \frac{1}{|x|^{N+2s}}$ . See Corollary \ref{corol_nonstrict}.
%}
$$\lim_{t\to 0^+} \frac{f(t)}{t}=+\infty$$
in a non-strict sense (i.e. $\lim_{t\to 0} \frac{|f(t)|}{|t|^{r-1}}=0$ for each $r \in [1+\frac{\alpha}{N}, 2)$, for example $f(t) \sim - t \log(t^2)$) is included, and as we expect the decay is of order $\sim \frac{1}{|x|^{N+2s}}$. %It is sufficient to apply the argument of Remark \ref{rem_dec_2sN} (since $f(t) \geq \underline{C} t$ for $t$ small and positive), and the results in Proposition \ref{prop_estim_above} (after having chosen a whatever $r \in [\frac{N+\alpha+4s}{N+2s}, 2)$). 
See Corollary \ref{corol_nonstrict}.
\item The conclusions hold also without assuming radial symmetry and monotonicity of $u$, but by assuming a priori that
$$\limsup_{|x|\to +\infty} |u(x)| |x|^{\omega} < +\infty$$
for some $\omega> \frac{N^2}{N+\alpha}$: see Remark \ref{rem_cond_asym_dec}. When $u\in L^q(\R^N)$, $q < \frac{N+\alpha}{N}$, is radially symmetric and decreasing, this is the case with $\omega=\frac{N}{q}$ (see Remark \ref{rem_dec_N}); in particular, if $q=1$, we have $\omega=N$. Notice that $u$ is automatically radially symmetric and decreasing when $u\in C^{1,1}_{loc}(\R^N)$, $f(u)=|u|^{r-2}u$ and $\omega>\frac{\alpha}{r-1}$ thanks to \cite[Theorem 1]{Le1} (see also \cite[Theorem 1.3]{WY0}).
\item In light of the previous remark, we highlight that the estimate from above actually holds true also for nonnegative solutions $u \geq 0$; see Proposition \ref{prop_estim_above}; moreover, it can be further extended to $|u|$ in the case of changing sign solutions, by applying a Kato's inequality \cite[Theorem 3.2]{Amb-1}.
%\item the estimate from above actually holds true also for nonnegative solutions $u \geq 0$; see Proposion \ref{prop_estim_above};
\item The conclusions hold also for solutions $u\in L^1(\R^N) \cap C(\R^N)$ in the viscosity sense, without assuming $f$ H\"older continuous (which is needed in (ii) only to pass from weak to viscosity solutions): see Section \ref{sec_est_bel_vis}.
\item When \hyperref[(F10)]{\textnormal{(F10)}} holds, we actually have $F(t) \geq \underline{C} \frac{t^r}{r}$ for $t \in (0,\delta)$; thus, 
% $\int_{\R^N\setminus} F(u) \geq \frac{\overline{C}}{r} \norm{u}_{L^r(\R^N \setminus B_R)}^r>0$ for $R$ large; thus, 
being also $u\in L^{\infty}(\R^N)$, the condition $\int_{\R^N} F(u)>0$ means that $F$ is not \emph{too negative} in $[\delta, \norm{u}_{\infty}]$. 
We highlight that the energy term $\int_{\R^N} \big(I_{\alpha}*F(u)\big) F(u)$ is always positive (see e.g. Proposition \ref{prop_positivity_Riesz}).
\item We find some estimates on the asymptotic constants, which are coherent, when $r \in [2^{\#}_{\alpha}, %\frac{N+\alpha}{N}, 
r^*_{\alpha,s})$, with the one found in Theorem \ref{thm_MVS} and Theorem \ref{thm_main}: %(when $F(u)= \frac{1}{r}|u|^r$, $\bar{C}=\underline{C}=1$ and $\mu=1$): 
see Propositions \ref{prop_estim_above} and \ref{prop_below_2}, and Corollary \ref{corol_sharp_decay}. 
We notice that \eqref{eq_fr_ch_hom} is obtained by \eqref{eq_introduction} formally choosing $f(t)=\sqrt{r}|t|^{r-2}t$. In our proofs -- up to well posedness and regularity -- we do not use that $F$ is the primitive of $f$: in particular, we do not apply \hyperref[(F9)]{\textnormal{(F9)}} and \hyperref[(F10)]{\textnormal{(F10)}} to $F$. Thus we can arbitrary move constants from $f$ to $F$ in our arguments to adjust -- for example -- the value of $\underline{C}$, and this allows to gain the result for every $\mu>0$ (see also Corollary \ref{corol_slower_deca2}).
\end{itemize}
\end{Remark}

\smallskip

Our results apply in particular to Pohozaev minima of the equation, % (see Definition \ref{def_pohozaev_min2}), 
whenever some symmetric assumption is assumed on $f$, that is \hyperref[(F7)]{\textnormal{(F7)}}-\hyperref[(F8)]{\textnormal{(F8)}}.
We notice that, since every Pohozaev minimum has constant sign, it is not restrictive to assume a priori the sign of $u$.

\begin{Corollary}\label{corol_main}
Assume \hyperref[(F1s)]{\textnormal{(F1)}}-\hyperref[(F2s)]{\textnormal{(F2)}} and \hyperref[(F7)]{\textnormal{(F7)}}-\hyperref[(F8)]{\textnormal{(F8)}}. 
%\hyperref[(f5)]{\textnormal{(f5)}}. 
Let $u$ be a (positive) Pohozaev minimum of \eqref{eq_introduction}. 
Then the conclusions of Theorem \ref{thm_main} hold.
\end{Corollary}

\medskip

We finally want to highlight that our results may be adapted to the local case $s=1$, extending Theorem \ref{thm_MVS} to general nonlinearities, studied in \cite{MS2}. 
We leave the details to the reader, observing that in this case the rate of decaying is simply given by $\beta= \frac{N-\alpha}{2-r}$, since, as already observed, the solutions of the homogeneous linear (associated) equation decay exponentially.

\begin{Theorem}\label{thm_main_loc}
Let $s=1$ and $N\geq 3$, and assume \hyperref[(F1)]{\textnormal{(F1)}}-\hyperref[(F2)]{\textnormal{(F2)}}.
% (where the upper critical exponent is substituted by $2^*_{\alpha,1}=\frac{N+\alpha}{N-2}$). 
Let $u\in H^1(\R^N)$, strictly positive, radially symmetric and decreasing, be a solution of \eqref{eq_Choquard_genericaF};
%$$-\Delta u + \mu u = \big(I_{\alpha}*F(u)\big)f(u) \quad \hbox{on $\R^N$};$$
in particular, $u$ may be a ground state. Let $r \in [2^{\#}_{\alpha}, %1+ \frac{\alpha}{N}, 
2)$.
\begin{itemize}
\item[\textnormal{(i)}] Assume \hyperref[(F9)]{\textnormal{(F9)}} and $\mu> (r-1) \bar{C}^{\frac{1}{r-1}}$. Then 
$$\limsup_{|x|\to +\infty} u(x) |x|^{ \frac{N-\alpha}{2-r}} \in (0, +\infty).$$
\item[\textnormal{(ii)}] Assume \hyperref[(F10)]{\textnormal{(F10)}} and $\int_{\R^N}F(u)>0$ (e.g. $F\geq 0$ on $(0,+\infty)$). 
Then
$$\liminf_{|x| \to +\infty} u(x)|x|^{\frac{N-\alpha}{2-r}} \in (0, +\infty).$$
\end{itemize}
%\tr{aggiusta label sezioni}
If both conditions (i) and (ii) hold, together with $\overline{C}=\underline{C}$, then \eqref{eq_sharp_decay_F} holds.
\end{Theorem}

\medskip

In both the estimates from above and below in Theorem \ref{thm_main} we rely on some comparison principle and the use of some auxiliary function whose fractional Laplacian is related to the Gauss hypergeometric function.
For the estimate from above we succeed in working with the weak formulation of the problem; %, by adding a restriction on $\mu$ in order to control the interaction among the linear and the nonlinear term.
%; we notice that in \cite{MS0} this restriction is automatically fulfilled since $\bar{C}=\mu=1$ and $r<2$. 
on the other hand, in order to deal with the estimate from below, we find the necessity of working with $u^{2-r}$, where $2-r \in (0,1)$: this concave power of the solution may fail to lie in $H^s(\R^N)$, and thus we cannot treat the problem with its weak formulation. 
The pointwise formulation seems to arise some problems as well, since the fractional Laplacian of $u^{2-r}$ needs some restrictive assumption on $\alpha,s,N$ and $r$ in order to be well defined. 
This is why we work with a viscosity formulation of the problem: in this case, to pass from weak to viscosity solutions, we ask only a bit of H\"older regularity on $f$. 
We remark that the estimate from above may be treated with the viscosity formulation as well.

\medskip

The remaining part of the Chapter is organized as follows. %\tr{aggiusta label sezioni}
%In Section \ref{sec_prelim} we introduce definitions and notations, collecting some existence and comparison results in Appendix \ref{app_ex_comp}. 
In Section \ref{sec_frac_aux} we recall the suitable auxiliary function introduced in Section \ref{sec_hypergeo}, and % (see Appendix \ref{app_frac_dec} for some related asymptotic property) and 
establish some asymptotic behaviour on suitable comparison functions; other preliminary estimates are studied in Section \ref{sec_prel_est}.
Then in Section \ref{sec_above} we deal with the estimate from above, by working with the weak formulation, while in Section \ref{sec_est_bel_vis} we study the asymptotic behaviour from below, by exploiting a viscosity formulation. 
Finally in Section \ref{sec_proof_main_2} we conclude the proofs of Theorem \ref{thm_main} and its corollaries. % of the main results.

\subsection{Fractional auxiliary functions}
\label{sec_frac_aux}

In order to implement some comparison argument, in Section \ref{sec_hypergeo} we introduced the function
$$h_{\beta}(x):= \frac{1}{(1+|x|^2)^{\frac{\beta}{2}}}, $$
which behaves, at infinity, like $\sim\frac{1}{|x|^{\beta}}$, $\beta >0$, but lies in $H^s(\R^N)$, avoiding the pole in the origin when $\beta \geq N$.
This function 
%, we 
%%it would be useful to work with the fractional Laplacian of the function $\frac{1}{|x|^{\beta}}$, $\beta >0$. 
%%On the other hand, $(-\Delta)^s \frac{1}{|\,\cdot\,|^{\beta}}$ is well defined only for $\beta \in (0,N)$, because of the presence of a pole in the origin. 
%%Thus we 
%%make the following choice, by considering, 
%consider, for any $\beta>0$, 
%
%which 
verifies 
$$%\begin{equation*}%\label{eq_frac_hyper}
(-\Delta)^s h_{\beta} (x) = C_{\beta, N, s} \, {}_2F_1\left(\frac{N}{2}+s, \frac{\beta}{2} + s, \frac{N}{2}; -|x|^2\right)
$$%\end{equation*}
where $C_{\beta, N, s}:= 2^{2s} \frac{\Gamma\big(\frac{N}{2}+s\big) \Gamma\big(\frac{\beta}{2} + s\big)}{\Gamma\big(\frac{N}{2}\big) \Gamma\big(\frac{\beta}{2}\big)}>0 $ 
and ${}_2F_1$ 
denotes the Gauss hypergeometric function.
Notice that we will be interested in $\beta \in (0, N+2s].$
In Section \ref{sec_hypergeo} % \ref{app_frac_dec} 
we collected some results on Gauss hypergeometric functions and their asymptotic behaviour at infinity. We use now this auxiliary function to study some comparison at infinity. % function.

\begin{Lemma}[Comparison for weak equation]\label{lem_confr_2xbeta}
Let $u\in C(\R^N)$ be a weak %, strictly positive 
solution of
$$(-\Delta)^s u + \lambda u = \gamma h_{\beta} \quad \hbox{in $\R^N\setminus B_{\rho}(0)$}$$
for some $\lambda, \gamma>0$, $\rho>0$ and 
$$\beta \in \left(\frac{N}{2}, N+2s\right].$$
Then 
%$$0<\liminf_{|x|\to+\infty} u(x) |x|^{\beta} \leq \limsup_{|x|\to +\infty} u(x) |x|^{\beta} <\infty.$$
$$ \limsup_{|x|\to +\infty} u(x) |x|^{\beta} <\infty.$$
%More precisely, 
Moreover, if $\beta \in (\frac{N}{2}, N+2s)$, we have
$$\lim_{|x|\to +\infty} u(x) |x|^{\beta} = \frac{\gamma}{\lambda}.$$
\end{Lemma}

\claim Proof.
We start noticing that, since $\beta > \frac{N}{2}$, then the equation is well posed from a weak point of view.
By Lemma \ref{lem_esist_sol_part} there exists a continuous function $w \in H^{2s}(\R^N)$, such that
$$(-\Delta)^s w + \lambda w =0 \quad \hbox{in $\R^N\setminus B_{\rho}(0)$}$$
in the weak sense and pointwise, and moreover, for some $C_1'',C_2''>0$,
$$\frac{C_1''}{|x|^{N+2s}} <w(x) \leq \frac{C_2''}{|x|^{N+2s}}, \quad \hbox{for every $|x|>\rho$}.$$
Let thus define, for some $\tau, \sigma \in \R$ and $\theta \in [\beta, N+2s]$ %$\theta > \beta $ 
to be chosen,
$$v_{\tau,\sigma}(x):= \frac{\gamma}{\lambda}h_{\beta}(x)+ \sigma h_{\theta}(x) + \tau w(x)$$
for every $x \in \R^N$. 
We have, for $|x|>\rho$,
\begin{align*}
(-\Delta)^s v_{\tau,\sigma}(x) + \lambda v_{\tau,\sigma}(x) &= \gamma h_{\beta}(x) + \left( \frac{\gamma}{\lambda} (-\Delta)^s h_{\beta}(x) + \sigma (-\Delta)^s h_{\theta}(x) + \lambda \sigma h_{\theta}(x)\right) \\
&=: \gamma h_{\beta}(x) + g_{\sigma, \theta}(x).
\end{align*}
By Lemma \ref{lem_calcol_2potenz} we obtain
\begin{itemize}
\item if $\beta \in (\frac{N}{2}, N)\setminus \{ N-2s\}$,
$$g_{\sigma, \theta}(x) \sim \frac{\gamma}{\lambda} C'_{\beta, N, s} h_{\beta +2s}(x) + \sigma C'_{\theta, N, s} h_{\theta+2s}(x) + \lambda \sigma h_{\theta}(x) \quad \hbox{as $|x|\to +\infty$};$$
in this case we assume $\theta \in (\beta, \min\{N, \beta+2s\})\setminus \{N-2s\}$;

\item if $\beta =N$,
$$g_{\sigma, \theta}(x) \sim \frac{\gamma}{\lambda} C'_{N, N, s} \log(x) h_{N+2s}(x) + \sigma C'_{\theta, N, s} h_{N+2s}(x) + \lambda \sigma h_{\theta}(x) \quad \hbox{as $|x|\to +\infty$};$$
in this case we assume $\theta \in (N, N+2s)$;

\item otherwise 
$$g_{\sigma, \theta}(x) \sim \frac{\gamma}{\lambda} C'_{\beta, N, s} h_{N+2s}(x) + \sigma C'_{\theta, N, s} h_{N+2s}(x) + \lambda \sigma h_{\theta}(x) \quad \hbox{as $|x|\to +\infty$},$$
and in this case 
\begin{itemize}
\item if $\beta=N-2s$ (possible only if $N>4s$), we choose $\theta \in (N, N+2s)$;
\item if $\beta \in (N,N+2s)$, we choose $\theta \in (\beta, N+2s)$;
\item if $\beta =N+2s$, we simply assume $\theta = %> 
N+2s$.
\end{itemize}
\end{itemize}
Assume first $\beta < N+2s$. 
By the abovementioned choices of $\theta>\beta$ we obtain
$$g_{\sigma, \theta}(x) \sim \lambda \sigma h_{\theta}(x) \quad \hbox{as $|x|\to +\infty$}. $$
In particular, fixed $\eps>0$, for some $R =R_{\eps} (\gamma, \lambda, \beta, \theta, \sigma) \gg 0$ (we may assume $R>\rho$) we obtain
$$(1-\eps) \lambda \sigma h_{\theta}(x) \leq g_{\sigma, \theta}(x) \leq (1+\eps) \lambda \sigma h_{\theta}(x) \quad \hbox{for $|x|\geq R$}$$
if $\sigma >0$, and
$$(1+\eps) \lambda \sigma h_{\theta}(x) \leq g_{\sigma, \theta}(x) \leq (1-\eps) \lambda \sigma h_{\theta}(x) \quad \hbox{for $|x|\geq R$}$$
if $\sigma <0 $. Notice that $R$ does not depend on $\tau$. Thus
$$(-\Delta)^s v_{\tau,\overline{\sigma}}(x) + \lambda v_{\tau,\overline{\sigma}}(x) \geq \gamma h_{\beta}(x) + (1-\eps) \lambda \overline{\sigma} h_{\theta}(x) \geq \gamma h_{\beta}(x) \quad \hbox{in $\R^N\setminus B_R(0)$}$$
by choosing a whatever $\overline{\sigma}>0$, and
$$(-\Delta)^s v_{\tau,\underline{\sigma}}(x) + \lambda v_{\tau,\underline{\sigma}}(x) \leq \gamma h_{\beta}(x) + (1-\eps) \lambda \underline{\sigma} h_{\theta}(x) \leq \gamma h_{\beta}(x) \quad \hbox{in $\R^N\setminus B_R(0)$}$$
by choosing a whatever $\underline{\sigma}<0$. Summing up
\begin{equation}\label{eq_confr_sistem}
\parag{ 
& (-\Delta)^s v_{\tau,\overline{\sigma}}(x) + \lambda v_{\tau,\overline{\sigma}}(x) \geq \gamma h_{\beta}(x) & \quad \hbox{in $\R^N\setminus B_R(0)$}, \\ 
& (-\Delta)^s v_{\tau,\underline{\sigma}}(x) + \lambda v_{\tau,\underline{\sigma}}(x) \leq \gamma h_{\beta}(x) & \quad \hbox{in $\R^N\setminus B_R(0)$}. }
\end{equation}
We choose now $\overline{\tau}>0$ such that
\begin{equation*}\label{eq_confr_tau_sop}
v_{\overline{\tau},\overline{\sigma}} - u \geq 0 \quad \hbox{on $B_R(0)$}.
\end{equation*}
Indeed, we impose
$$\frac{\gamma}{\lambda}h_{\beta}(x)+ \overline{\sigma} h_{\theta}(x) + \tau w(x) \geq u(x) \quad \hbox{on $B_R(0)$}$$
that is
$$\tau w(x) \geq u(x) - \frac{\gamma}{\lambda}h_{\beta}(x) - \overline{\sigma} h_{\theta}(x) \quad \hbox{on $B_R(0)$}$$
which is satisfied if we impose (recall that $\overline{\sigma}>0$)
$$ \tau \min_{B_R} w \geq \max_{B_R} u - \frac{\gamma}{\lambda} h_{\beta}(R) \geq u(x) - \frac{\gamma}{\lambda}h_{\beta}(x)- \overline{\sigma} h_{\theta}(x) \quad \hbox{on $B_R(0)$}$$
that is
$$ \overline{\tau} \geq \frac{ \max_{B_R} u - \frac{\gamma}{\lambda} h_{\beta}(R)}{\min_{B_R} w }.$$
Similarly, we choose $\underline{\tau}\in \R%>0
$ such that
\begin{equation*}\label{eq_confr_tau_sot}
v_{\underline{\tau},\underline{\sigma}} - u \leq 0 \quad \hbox{on $B_R(0)$},
\end{equation*}
given by
$$ \underline{\tau} \leq %\geq 
\frac{ \min_{B_R} u - \frac{\gamma}{\lambda} h_{\beta}(R) }{\max_{B_R} w }.$$
We notice that both the minimum and the maximum of $w$ in the ball are finite and strictly positive, since $w>0$ is continuous. Thus, summing up
\begin{equation}\label{eq_confr_tau}
\parag{ & v_{\overline{\tau},\overline{\sigma}} - u \geq 0& \quad \hbox{on $B_R(0)$}, \\ &v_{\underline{\tau},\underline{\sigma}} - u \leq 0& \quad \hbox{on $B_R(0)$}.}
\end{equation}
By joining \eqref{eq_confr_sistem} with the assumption on $u$, we obtain
\begin{equation}\label{eq_confr_u_sistem}
\parag{ & (-\Delta)^s (v_{\overline{\tau},\overline{\sigma}} -u)(x)+ \lambda (v_{\overline{\tau},\overline{\sigma}}-u)(x) \geq 0 & \quad \hbox{in $\R^N\setminus B_R(0)$}, \\ 
& (-\Delta)^s (v_{\underline{\tau},\underline{\sigma}}-u)(x) + \lambda (v_{\underline{\tau},\underline{\sigma}}-u)(x) \leq 0 & \quad \hbox{in $\R^N\setminus B_R(0)$}. }
\end{equation}
By the weak version of the Comparison Principle (Lemma \ref{lem_comp_prin}) we obtain
$$
\parag{ & v_{\overline{\tau},\overline{\sigma}} - u \geq 0& \quad \hbox{on $\R^N$}, \\ &v_{\underline{\tau},\underline{\sigma}} - u \leq 0& \quad \hbox{on $\R^N$}.}
$$
that is
$$ \frac{\gamma}{\lambda}h_{\beta}(x)+ \underline{\sigma} h_{\theta}(x) + \underline{\tau} w(x) \leq u(x) \leq \frac{\gamma}{\lambda}h_{\beta}(x)+ \overline{\sigma} h_{\theta}(x) + \overline{\tau} w(x) $$
and hence, by the assumption on $w$,
$$ \frac{\gamma}{\lambda}h_{\beta}(x)+ \underline{\sigma} h_{\theta}(x) + \underline{\tau} \frac{C_1''}{|x|^{N+2s}} \leq u(x) \leq \frac{\gamma}{\lambda}h_{\beta}(x)+ \overline{\sigma} h_{\theta}(x) + \overline{\tau}\frac{C_2''}{|x|^{N+2s}}$$
for each $x \in \R^N$, $x\neq 0$. Thus
%$$ \frac{\gamma}{\lambda}\frac{|x|^{\beta}}{(1+|x|^2)^{\frac{\beta}{2}}}+ \underline{\sigma} \frac{|x|^{\beta}}{(1+|x|^2)^{\frac{\theta}{2}}} + \underline{\tau} \frac{C_1''}{|x|^{N+2s-\beta}} 
%\leq u(x)|x|^{\beta} 
%\leq \frac{\gamma}{\lambda} \frac{|x|^{\beta}}{(1+|x|^2)^{\frac{\beta}{2}}}+ \overline{\sigma} \frac{|x|^{\beta}}{(1+|x|^2)^{\frac{\theta}{2}}} + \overline{\tau}\frac{C_2''}{|x|^{N+2s-\beta}},$$
\begin{eqnarray*}
\lefteqn{\frac{\gamma}{\lambda}\frac{|x|^{\beta}}{(1+|x|^2)^{\frac{\beta}{2}}}+ \underline{\sigma} \frac{|x|^{\beta}}{(1+|x|^2)^{\frac{\theta}{2}}} + \underline{\tau} \frac{C_1''}{|x|^{N+2s-\beta}} 
\leq } \\
&& \leq u(x)|x|^{\beta} 
\leq \frac{\gamma}{\lambda} \frac{|x|^{\beta}}{(1+|x|^2)^{\frac{\beta}{2}}}+ \overline{\sigma} \frac{|x|^{\beta}}{(1+|x|^2)^{\frac{\theta}{2}}} + \overline{\tau}\frac{C_2''}{|x|^{N+2s-\beta}},
\end{eqnarray*}
which gives the claim passing to the limit $|x|\to +\infty$, since $\theta>\beta$ and $\beta < N+2s$. 

\smallskip

Assume now $\beta=N+2s$, and choose $\theta = %>
\beta=N+2s$. Now we have
$$g_{\sigma, \theta}(x) \sim \overline{C}_{\sigma} h_{N+2s}(x) \quad \hbox{as $|x| \to +\infty$}$$
where
$$\overline{C}_{\sigma}:=\frac{\gamma}{\lambda} C'_{N+2s,N,s} + \sigma C'_{N+2s,N,s} %C'_{\theta, N, s}
+ \lambda \sigma;$$
%$$\overline{C}_{\sigma}:=\left(\frac{\gamma}{\lambda} + \sigma \right)C'_{N+2s,N,s} %C'_{\theta, N, s}
%+ \lambda \sigma;$$
recall that $C'_{N+2s,N,s}, C'_{\theta,N,s}<0$.
We can choose proper $\overline{\sigma} \in \R$ such that $\overline{C}_{\overline{\sigma}}<0$
%$\overline{\sigma} <- \frac{\frac{\gamma}{\lambda}C'_{N+2s,N,s}}{C_{\theta,N,s}}<0 $ and $\underline{\sigma} \in (-\frac{\gamma}{\lambda} \frac{C'_{N+2s,N,s}}{C_{\theta,N,s}}, 0)$ such that \eqref{eq_confr_sistem} still hold.
and thus the first equation in \eqref{eq_confr_sistem} still hold.
Since the sign of $\overline{\sigma}$ may be %is 
now different, we choose
$$ \overline{\tau} \geq \frac{ \max_{B_R} u - \frac{\gamma}{\lambda} h_{\beta}(R)- \min\{\overline{\sigma},0\}%\overline{\sigma}
}{\min_{B_R} w }.$$
We come up then with the same proof, obtaining
$$ \limsup_{|x|\to +\infty} u(x)|x|^{\beta} \leq \frac{\gamma}{\lambda}+ \overline{\sigma} + \overline{\tau}C_2''.$$
%$$ \frac{\gamma}{\lambda}+ + \underline{\tau} C_1'' \leq \liminf_{|x|\to +\infty} u(x)|x|^{\beta} \leq \limsup_{|x|\to +\infty} u(x)|x|^{\beta} \leq \frac{\gamma}{\lambda}+ \overline{\tau}C_2''.$$
Notice that the appearing constants depend on $u, \gamma,\lambda, \rho, \beta, N, s$.
\QED

\bigskip

\begin{Lemma}[Comparison for pointwise equation]\label{lem_confr_3xbeta}
Let $u \in C(\R^N)$ be a pointwise %, strictly positive 
solution of
$$(-\Delta)^s u + \lambda u = \gamma h_{\beta} \quad \hbox{in $\R^N\setminus B_{\rho}(0)$}$$
for some $\lambda, \gamma>0$, $\rho>0$ and 
$$\beta \in (0, N+2s].$$
Then
$$0<\liminf_{|x|\to+\infty} u(x) |x|^{\beta} \leq \limsup_{|x|\to +\infty} u(x) |x|^{\beta} <\infty.$$
More precisely, if $\beta \in (0, N+2s)$, we have
$$\lim_{|x|\to +\infty} u(x) |x|^{\beta} = \frac{\gamma}{\lambda}.$$
\end{Lemma}

\claim Proof.
The proof goes as the previous Lemma, with the difference that at the end we apply the pointwise version of the Comparison Principle (Lemma \ref{lem_comp_3prin}).
\QED

\subsection{A preliminary range}
\label{sec_prel_est}

We start with some observations.

\begin{Remark}\label{rem_dec_N}
Let $u\in L^q(\R^N)$, for some $q \in [1,+\infty)$, be continuous and such that $|u|$ is radially symmetric and decreasing. Then, for every $x \in \R^N$,
\begin{align*}
|u(|x|)|^q |x|^N &= N |u(|x|)| \int_0^{|x|} t^{N-1} dt = N \int_0^{|x|} |u(|x|)|^q t^{N-1} dt \\
&\leq N \int_0^{|x|} |u(t)|^q t^{N-1} dt = \frac{N}{\omega_{N-1}} \int_{B_{|x|}(0)} |u(y)|^q dy \leq \frac{N}{\omega_{N-1}} \norm{u}_{L^q(\R^N)}^q
\end{align*}
where $\omega_{N-1}$ denotes the area of the $N-1$ dimensional sphere. Thus 
$$|u(x)| \leq \frac{C^2_u}{|x|^{\frac{N}{q}}}, \quad x \neq 0$$
where $C^2_u := C_N \norm{u}_q^q >0$. In particular, if $u \in L^1(\R^N)$, we have
$$|u(x)| \leq \frac{C^2_u}{|x|^N}, \quad x \neq 0.$$
\end{Remark}

%By Remarks \ref{rem_dec_2sN} and \ref{rem_dec_N}, we obtain that every strictly positive, continuous, radially symmetric and decreasing solution of \eqref{eq_main_prel} in $L^1(\R^N)$ satisfies
%\begin{equation}\label{eq_stima_grezza}
%\frac{C^1_u}{|x|^{N+2s}} \leq u(x) \leq \frac{C^2_u}{|x|^N} \quad \hbox{for $|x|\geq R \gg 0$},
%\end{equation}
% whenever $f$ satisfies \hyperref[(F1s)]{\textnormal{(F1)}}-\hyperref[(F2s)]{\textnormal{(F2)}} and \hyperref[(F10)]{\textnormal{(F10)}}.
%Thus the goal is to improve this asymptotic decay in the case of sublinear nonlinearities.

We keep with some preliminary results.
%; in particular they show that the decay of $I_{\alpha}*F(u)$ ``stabilizes'' when $u$ satisfies \eqref{eq_stima_grezza}, and thus no good ``bootstrap'' argument can be implemented.

\begin{Lemma}\label{lem_decay_Ialpha}
Let $u\in L^1 (\R^N)$ continuous be such that $|u|$ is radially symmetric and decreasing. Let $f$ satisfy \hyperref[(F1s)]{\textnormal{(F1)}}-\hyperref[(F2s)]{\textnormal{(F2,i)}}, and let $\theta \in (N, N+\alpha ]$. Then there exists $C=C(N,\alpha)>0$ such that
$$\pabs{\big(I_{\alpha}*F(u)\big)(x) - I_{\alpha}(x) \int_{\R^N} F(u)} \leq C \norm{F(u)}_{\infty,\theta} I_{\alpha}(x) \left( \frac{1}{1+|x|} + \frac{1}{1+ |x|^{\theta-N}}\right)$$
for each $x \in \R^N$, $x \neq 0$.
\end{Lemma}

\claim Proof.
First notice that $u\in L^{\infty}(\R^N)$, $F(u) \in L^{\infty}(\R^N)$, and that $I_{\alpha}*F(u)$ and $\int_{\R^N} F(u) $ are finite and well defined. 
By Remark \ref{rem_dec_N} we have
$$|u(x)| \leq \frac{C_u^2}{|x|^N} \to 0.$$
Thus $\big|F(u(x))\big| |x|^{\theta}$ is bounded on a ball $B_R$ (since $F(u)$ is bounded), and it is bounded on the complement of this ball since
$$\big|F(u(x))\big| |x|^{\theta} = \frac{\big|F(u(x))\big|}{|u(x)|^{\frac{N+\alpha}{N}}} |u(x)|^{\frac{N+\alpha}{N}} |x|^{\theta} \leq \frac{\big|F(u(x))\big|}{|u(x)|^{\frac{N+\alpha}{N}}} \frac{C}{|x|^{N+\alpha -\theta}}$$
by considering the growth condition \hyperref[(F2)]{\textnormal{(F2,i)}} of $F$ in zero (when $R\gg 0$, not depending on $\theta$) and the restriction on $\theta$. Thus
$$\sup_{x \in \R^N} \big|F(u(x))\big| |x|^{\theta} < +\infty$$
and Lemma \ref{lem_stima_1_Riesz} applies with $g(x):= F(u(x))$, which concludes the proof. We further 
notice that
$$\norm{F(u)}_{\infty,\theta}\leq \norm{F(u)}_{\infty}(1+ R^{\theta}) + \left(\limsup_{t \to 0} \frac{|F(t)|}{|t|^{\frac{N+\alpha}{N}}} \right)\frac{1+R^{\theta}}{R^{N+\alpha}}$$
for any $\theta \in (N, N+\alpha]$ and any $R\gg 0$ (not depending on $\theta$, but depending on $u$).
\QED

\begin{Remark}\label{rem_cond_asym_dec}
 In what follows, for the sake of exposition we will restrict our analysis to the space of radially symmetric and decreasing functions in $L^1(\R^N)$, but we highlight that this assumption is needed only to get the a priori asymptotic decay of Remark \ref{rem_dec_N}. 
By the above proof, actually we see that we may ask only
$$|u(x)| \leq \frac{C}{|x|^{\omega}}$$
for some $\omega$ such that
$$\omega > \frac{N^2}{N+\alpha}.$$
In particular $\omega=N$, obtained in Remark \ref{rem_dec_N}, fits this condition. 
Alternatively, one may assume this a priori asymptotic decay on $u$ (and adapt the restrictions on $\theta$ by $\theta \in (N, \frac{N+\alpha}{N} \omega]$).
\end{Remark}

\begin{Corollary}\label{corol_stima_Riesz}
Let $u\in L^1 (\R^N)$ continuous be such that $|u|$ is radially symmetric and decreasing. 
Let $f$ satisfy \hyperref[(F1s)]{\textnormal{(F1)}}-\hyperref[(F2s)]{\textnormal{(F2,i)}}, and let $\theta \in (N, N+\alpha]$. 
Then for any $\eps>0$, there exists $R_{\eps}=R_{\eps}(N, \alpha,\theta)\gg 0 $ such that
$$\Big|\big(I_{\alpha}*F(u)\big)(x) \Big| \leq I_{\alpha}(x) \left( \pabs{\int_{\R^N} F(u)} + \eps \norm{F(u)}_{\infty,\theta}\right)$$
and
$$\big(I_{\alpha}*F(u)\big)(x) \geq I_{\alpha}(x) \left(\int_{\R^N} F(u)-\eps \norm{F(u)}_{\infty, \theta}\right) $$
for each $|x|\geq R_{\eps}$.
\end{Corollary}

\begin{Remark}\label{rem_dec_2sN}
In Section \ref{sec_doub_superl} it was showed that the solutions decay as fast as $\sim \frac{1}{|x|^{N+2s}}$ when the nonlinearity is linear or superlinear. 
In the sublinear case, we expect a slower decay. 
Indeed, assume \hyperref[(F1s)]{\textnormal{(F1)}}-\hyperref[(F2s)]{\textnormal{(F2)}} and \hyperref[(F10)]{\textnormal{(F10)}}, and let $u$ be a strictly positive solution of \eqref{eq_main_prel}.
%Since by Theorem \ref{th_INT_regular} and Proposition \ref{prop_u_L1} we have $u \in L^1(\R^N) \cap L^{\infty}(\R^N) \subset %L^{\frac{N}{2s}}(\R^N) \cap 
%L^{\frac{N+2\alpha}{2s+\alpha}}(\R^N)\cap L^{\infty}(\R^N)$, which in particular implies $(I_{\alpha}*F(u))f(u)\in L^2(\R^N) \cap L^{\frac{N}{2s}}(\R^N)\cap L^{\infty}(\R^N)$, we can apply
By 
 Lemma \ref{lem_u_to0} we %and 
 obtain $u(x)\to 0$ as $|x|\to +\infty$\footnote{If $u$ is assumed radially symmetric and $s \in (\frac{1}{2},1)$, this is actually a consequence of the Strauss radial lemma. If $u$ is radially symmetric and decreasing, it is a consequence of the monotonicity and (a whatever) summability.}. 
 Thus there exists $R\gg 0$ such that 
% $|u(x)|\leq \delta <1$ for $|x|\geq R$ and thus
$0 \leq u(x)\leq \delta <1$ for $|x|\geq R$ and thus
$$f(u(x))\geq \underline{C} u^{r-1}(x) \quad \hbox{for $|x|\geq R$}$$
together with
%$$|u(x)|^{r-1} \geq |u(x)| \quad \hbox{ for $|x|\geq R$}.$$
$$u^{r-1}(x) \geq u(x) \quad \hbox{ for $|x|\geq R$}.$$
%Hence %being $F(u(x))= \int_0^{u(x)} f(\tau) \geq \underline{C} \int_0^{u(x)} \tau^{r-1} \geq 0$,
If we assume $(I_{\alpha}*F(u))(x) \geq 0$ for $|x|\geq R$, we gain
$$(-\Delta)^s u + \mu u \geq (I_{\alpha}*F(u))u \quad \hbox{on $\R^N\setminus B_R(0)$}$$
which implies
$$(-\Delta)^s u + \tfrac{3}{2}\mu u \geq \left(I_{\alpha}*F(u) + \tfrac{1}{2}\mu\right) u \quad \hbox{on $\R^N\setminus B_R(0)$}.$$
By Proposition \ref{prop_conv_C0} we have that $(I_{\alpha}*F(u))(x)\to 0$ as $|x|\to +\infty$,%
\footnote{Notice that the claim is obtained by assuming that $u$ is a weak solution or, alternatively, that $u$ is in the right Lebesgue spaces; in particular it is true if $u\in L^1(\R^N) \cap L^{\infty}(\R^N)$ a priori.}, 
 thus for some $R'\geq R \gg0$ we get
$$(-\Delta)^s u + \tfrac{3}{2}\mu u \geq 0 \quad \hbox{on $\R^N\setminus B_{R'}(0)$}.$$
At this point (being $u$ strictly positive) we conclude as in the proof of Theorem \ref{th_INT_decay} and obtain
$$u(x) \geq \frac{C^1_u}{|x|^{N+2s}} \quad \hbox{for $|x| \geq R$}$$ 
for some constant $C^1_u = C_{N,\alpha, R, \mu} \min_{B_R} u>0$ and some sufficiently large $R\gg0$. 
\end{Remark}

By Remarks \ref{rem_dec_2sN} and \ref{rem_dec_N}, we obtain that every strictly positive, continuous, radially symmetric and decreasing solution of \eqref{eq_main_prel} in $L^1(\R^N)$ satisfies
\begin{equation}\label{eq_stima_grezza}
\frac{C^1_u}{|x|^{N+2s}} \leq u(x) \leq \frac{C^2_u}{|x|^N} \quad \hbox{for $|x|\geq R \gg 0$},
\end{equation}
 whenever $f$ satisfies \hyperref[(F1s)]{\textnormal{(F1)}}-\hyperref[(F2s)]{\textnormal{(F2)}} and \hyperref[(F10)]{\textnormal{(F10)}}, together with $\int_{\R^N} F(u)>0$: indeed in this case, by Lemma \ref{lem_decay_Ialpha}, we have $\big(I_{\alpha}*F(u)\big)(x) \sim I_{\alpha}(x) \int_{\R^N} F(u) >0$ for $|x|$ large.
Thus the goal is to improve the asymptotic decay \eqref{eq_stima_grezza} in the case of sublinear nonlinearities.

\begin{Remark}\label{rem_dim_induzione}
By Lemma \ref{lem_confr_2xbeta}, Corollary \ref{corol_stima_Riesz}, and a bootstrap argument we can give some first qualitative proofs of the main result. Indeed, by
$$(-\Delta)^s u + \mu u = 
\big(I_{\alpha}*F(u)\big) f(u) \lesssim
 I_{\alpha}(x) u^{r-1} \lesssim 
 \frac{1}{|x|^{N-\alpha}} \frac{1}{|x|^{\gamma_0(r-1)}} = 
 \frac{1}{|x|^{\gamma_1}}$$
where $\gamma_0:=N$ and $\gamma_1:=\gamma_0(r-1)+ N-\alpha$, and a comparison argument, we obtain $u(x) \lesssim \frac{1}{|x|^{\gamma_1}}$. By induction, set
$$\gamma_{i+1} := \gamma_i(r-1)+ N-\alpha$$
we obtain $u(x) \lesssim \frac{1}{|x|^{\gamma_i}}$ and $\gamma_i \nearrow \frac{N-\alpha}{2-r}$ (but the argument works only for $\gamma_i \leq N+2s$). Similarly, 
$$(-\Delta)^s u + \mu u = 
\big(I_{\alpha}*F(u)\big) f(u) \gtrsim
 I_{\alpha}(x) u^{r-1} \gtrsim
 \frac{1}{|x|^{N-\alpha}} \frac{1}{|x|^{\gamma_0(r-1)}} = 
 \frac{1}{|x|^{\gamma_1}}$$
 where now $\gamma_0:=N+2s$, which implies $u(x) \gtrsim \frac{1}{|x|^{\gamma_i}}$ and $\gamma_i \searrow \frac{N-\alpha}{2-r}$ if $r \leq r^*_{s,\alpha}$ (while the case $r \geq r^*_{s,\alpha}$ implying $\gamma_i \nearrow \frac{N-\alpha}{2-r}$ cannot be set in motion).

In order to pass to the limit we have to take care %The roughness of the results, together with the difficulty in taking care 
of the bounding constants (or, equivalently, of the radii related to the complements of the balls where the inequalities hold), which is not an easy task; see anyway Corollary \ref{corol_slower_deca2}. This suggests the implementation of more direct approaches, as done in following Sections.
\end{Remark}

\subsection{Estimate from above}
\label{sec_above}

First, we deal with the estimate from above. In this case we succeed in arguing in the weak sense with no additional assumption on $f$.

\begin{Proposition}\label{prop_estim_above}
Assume \hyperref[(F1s)]{\textnormal{(F1)}} and \hyperref[(F9)]{\textnormal{(F9)}}. 
Let $u\in H^s(\R^N)\cap L^1(\R^N)$, continuous, nonnegative, %strictly positive, 
radially symmetric and decreasing, be a weak solution of \eqref{eq_main_prel}. 
Assume moreover
$$\mu> (r-1) \bar{C}^{\frac{1}{r-1}}.$$
Then, set 
$\beta:=\min\left\{\frac{N-\alpha}{2-r}, N+2s\right\}$, 
we have, for some $C_u \geq 0$, 
$$\limsup_{|x|\to +\infty} u(x) |x|^{\beta} \leq C_u;$$
if $\beta< N+2s$, the constant $C_u$ depends on $u$ in the following way:
$$C_u= \frac{(2-r) \left(C_{N,\alpha} \pabs{ \int_{\R^N} F(u)}\right)^{\frac{1}{2-r}} }{\mu - (r-1)\bar{C}^{\frac{1}{r-1}} }
$$
where $C_{N,\alpha} >0$ is given in \eqref{eq_def_Riesz}.
\end{Proposition}

\begin{Remark}
We observe that
$$r= 1+ \frac{\alpha}{N} \implies \beta=N,$$
$$r \in \left( 1+ \frac{\alpha}{N}, 1 + \frac{\alpha+2s}{N+2s}\right) \implies \beta = \frac{N-\alpha}{2-r} \in (N, N+2s),$$
$$r \in \left[1+\frac{\alpha+2s}{N+2s}, 2\right) \implies \beta = N+2s;$$
actually, as already observed, the asymptotic decay with $\beta=N+2s$ applies for general $r \in [1+\frac{\alpha+2s}{N+2s}, 1+\frac{\alpha+2s}{N-2s}]$, including linear and superlinear cases, thanks to the results in Section \ref{sec_doub_superl}. 
We notice that, when $r> 1+ \frac{\alpha}{N}$, we are actually improving \eqref{eq_stima_grezza}.
\end{Remark}

\claim Proof.
We start noticing that, by the Young product inequality, we obtain 
$$(I_{\alpha}*F(u))f(u) \leq \frac{1}{a} \big|I_{\alpha}*F(u)\big|^a + \frac{1}{b} |f(u)|^b$$
when $a,b>0$, $\frac{1}{a} + \frac{1}{b}=1$. In particular we choose $b=\frac{1}{r-1}$ and thus $a= \frac{1}{2-r}>0$ (possible thanks to the sublinearity restriction on $r$); with this choice, by \eqref{eq_cond_2sublin_f} and the fact that $u(x)\to 0$ as $|x|\to +\infty$, we obtain
$$(I_{\alpha}*F(u))f(u) \leq (2-r) \big|I_{\alpha}*F(u)\big|^{\frac{1}{2-r}} + (r-1)\bar{C}^{\frac{1}{r-1}} u$$
for $|x|\geq R$, where $R=R(u)\gg 0$ is sufficiently large. 
By Corollary \ref{corol_stima_Riesz}, for a whatever fixed $\theta \in (N, N+\alpha]$ and any $\eps>0$
we obtain
\begin{align*}
(I_{\alpha}*F(u))f(u) &\leq (2-r) \left(I_{\alpha}(x) \left( \pabs{ \int_{\R^N} F(u)} + \eps \norm{F(u)}_{\infty, \theta}\right)\right)^{\frac{1}{2-r}} + (r-1)\bar{C}^{\frac{1}{r-1}} u \\
&= (2-r) \left( \pabs{ \int_{\R^N} F(u)}+ \eps \norm{F(u)}_{\infty, \theta}\right)^{\frac{1}{2-r}} \frac{C_{N,\alpha}^{\frac{1}{2-r}}}{|x|^{\frac{N-\alpha}{2-r}}} + (r-1)\bar{C}^{\frac{1}{r-1}} u
\end{align*}
for every $|x| \geq R_{\eps}=R_{\eps}(u,N, \alpha, \theta)$, 
thus
$$(-\Delta)^s u + \mu u \leq (2-r) C_{N,\alpha}^{\frac{1}{2-r}} \left( \pabs{ \int_{\R^N} F(u)}+ \eps \norm{F(u)}_{\infty,\theta}\right)^{\frac{1}{2-r}} \frac{1}{|x|^{\frac{N-\alpha}{2-r}}} + (r-1)\bar{C}^{\frac{1}{r-1}} u .$$
Notice that $F(u) \nequiv 0$ (otherwise, by the equation, $u\equiv 0$ and the claim is trivial), thus we set
%Set
$$\gamma_{u, \eps}:= (2-r) C_{N,\alpha}^{\frac{1}{2-r}} \left( \pabs{ \int_{\R^N} F(u)} + \eps \norm{F(u)}_{\infty, \theta}\right)^{\frac{1}{2-r}} %\geq 
>0$$
and
$$\lambda:=\mu - (r-1)\bar{C}^{\frac{1}{r-1}} >0$$
we obtain
$$(-\Delta)^s u + \lambda u \leq \frac{\gamma_{u,\eps}}{|x|^{\beta}} \quad \hbox{ in $\R^N \setminus B_{R_{\eps}}(0)$};$$
notice that we use the fact that $\frac{1}{|x|^{\frac{N-\alpha}{2-r}}} \leq \frac{1}{|x|^{\beta}}$ for $|x|$ large.
%; we can assume $\gamma_{u,\eps}>0$ (otherwise the claim is trivial). 
For each $\delta>1$ we have $\frac{1}{|x|^{\beta}} \leq \delta h_{\beta}(x)$ when $|x|> R_{\delta} :=(\delta^{\frac{2}{\beta}}-1)^{-\frac{1}{2}}$; we may choose $R_{\delta,\eps}>\max\{R_{\delta}, R_{\eps}\}$. Thus
\begin{equation}\label{eq_u_2gammau}
(-\Delta)^s u + \lambda u \leq \delta \gamma_{u,\eps} h_{\beta}(x) \quad \hbox{ in $\R^N\setminus B_{R_{\delta, \eps}}(0)$}.
\end{equation}
We have $ h_{\beta} \in L^2(\R^N)$, since $2\frac{N-\alpha}{2-r}>N$.
By Lemma \ref{lem_lax_milgr}, being $u \in H^s(\R^N)$, %\cap C(\R^N)$, 
there exists $v\in H^s(\R^N)$ such that
$$\parag{ &(- \Delta)^s v + \lambda v = \delta \gamma_{u,\eps} h_{\beta}(x)& \quad \hbox{ in $\R^N \setminus B_{R_{\delta,\eps}}(0)$}, \\ &v = u& \quad \hbox{on $B_{R_{\delta,\eps}}(0)$}.}$$
Joining the first equation with \eqref{eq_u_2gammau} we obtain
$$(-\Delta)^s (u-v) + \lambda(u-v) \leq 0 \quad \hbox{ in $\R^N \setminus B_{R_{\delta,\eps}}(0)$}$$
and thus, by the weak version of the Comparison Principle (Lemma \ref{lem_comp_prin}) we gain
\begin{equation}\label{eq_stima_uv_2comp}
u \leq v \quad \hbox{on $\R^N$}.
\end{equation}
By Lemma \ref{lem_confr_2xbeta}, if $\beta<N+2s$, we can estimate $v$ by
$$\limsup_{|x|\to +\infty} v(x)|x|^{\beta} \leq \frac{ \delta \gamma_{u,\eps}}{\lambda}.$$
This relation, combined with \eqref{eq_stima_uv_2comp}, gives
$$\limsup_{|x|\to +\infty} u(x) |x|^{\beta} \leq \frac{ \delta \gamma_{u,\eps}}{\lambda}$$
for each $\delta>1$. In particular, as $\delta \to 1^+$ and $\eps \to 0^+$, we obtain the claim. 
If $\beta=N+2s$, we argue similarly (without moving $\delta$ and $\eps$).
\QED

\bigskip

Notice that, if we assume \eqref{eq_strongf3}, then one can choose every $\bar{C}>0$, and thus in particular every $\mu>0$ is allowed. 
%In \cite{MS0} the condition on $\mu$ is not evident, since $\mu=1$ is fixed and being $f(t)=|t|^{r-2}t$ they have $\bar{C}=1$ and thus $\mu>(r-1)\bar{C}^{\frac{1}{r-1}}$ is automatically fulfilled for every $r<2$.

\begin{Corollary}
Assume \hyperref[(F1s)]{\textnormal{(F1)}} and the condition \eqref{eq_strongf3}. 
Let $u\in H^s(\R^N) \cap L^1(\R^N)$, continuous, nonnegative, % strictly positive, 
radially symmetric and decreasing, be a solution of \eqref{eq_main_prel}. 
Then, set $\beta:=\min\left\{\frac{N-\alpha}{2-r}, N+2s\right\}$, 
we have
$$\limsup_{|x|\to +\infty} u(x) |x|^{\beta} \leq C_u;$$
if $\beta< N+2s$ 
the constant $C_u>0$ depends on $u$ in the following way:
$$C_u := \frac{(2-r) \left(C_{N,\alpha} \pabs{ \int_{\R^N} F(u)}\right)^{\frac{1}{2-r}} }{\mu }.$$
\end{Corollary}

We observe that the previous estimate from above is still valid by considering viscosity solutions $u \in L^1(\R^N) \cap C(\R^N)$, see Section \ref{sec_est_bel_vis}. We leave the details to the reader.

\subsection{Estimate from below} %: viscosity solutions} 
\label{sec_est_bel_vis}

Next, we deal with the estimate from below. 
%We need first some preliminary results, in order to 
Here we need to deal with the fractional Laplacian of the concave power of a function: since it might happen that $u^{\theta} \notin H^s(\R^N)$ when $u \in H^s(\R^N)$ and $\theta \in (0,1)$, the weak formulation seems not to be appropriate. 
Similarly, $(-\Delta)^s u^{\theta}$ might be not well defined pointwise, even if $u$ is regular enough. 
Notice that knowing a priori that $u$ is continuous, radially symmetric and decreasing seems of no use. 
The idea is thus to treat the problem via viscosity formulation, and we do it by exploiting the concave chain rule obtained in Section \ref{sec_chain_rule}.

\begin{Proposition}\label{prop_below_2}
Assume \hyperref[(F1s)]{\textnormal{(F1)}}-\hyperref[(F2s)]{\textnormal{(F2,i)}} and the sublinear condition \hyperref[(F10)]{\textnormal{(F10)}}. 
Let $u\in L^1(\R^N) \cap C(\R^N)$, strictly positive, radially symmetric and decreasing, be a viscosity solution of \eqref{eq_main_prel}. Assume $\int_{\R^N}F(u)>0$.
Then, 
$$\liminf_{|x|\to +\infty} u(x) |x|^{\frac{N-\alpha}{2-r}} \geq C'_u$$
where 
$$C'_u := \left( \frac{\underline{C} C_{N,\alpha} \int_{\R^N} F(u) dx}{\mu } \right)^{\frac{1}{2-r}} $$
and $C_{N,\alpha}>0$ is given in \eqref{eq_def_Riesz}. 
Moreover, set $\beta:=\min\left\{\frac{N-\alpha}{2-r}, N+2s\right\}$, 
we have, for some $C''_u>0$, 
$$\liminf_{|x|\to +\infty} u(x) |x|^{\beta} \geq C''_u;$$
if $\frac{N-\alpha}{2-r}\leq N+2s$ (i.e. $\beta=\frac{N-\alpha}{2-r}$), we have $C''_u:=C'_u$, otherwise we have $C''_u:=C^1_u$ (see Remark \ref{rem_dec_2sN}).
\end{Proposition}

\claim Proof.
First notice that, by the assumptions, $u\in L^1(\R^N) \cap L^{\infty}(\R^N)$ and thus, by Remark \ref{rem_conv_welldef}, $I_{\alpha}*F(u)$ is pointwise well defined.

By Corollary \ref{corol_concav_u2}, since $2-r \in (0,1-\frac{\alpha}{N}]\subset (0,1)$ we have 
$$(-\Delta)^s u^{2-r} \geq \frac{2-r}{u^{r-1}} \Big(-\mu u + \big(I_{\alpha}*F(u))f(u)\Big) $$
on $\R^N$, in the viscosity sense.
Thus
$$ (-\Delta)^s u^{2-r} + \mu (2-r) u^{2-r} \geq (2-r) \frac{\big(I_{\alpha}*F(u))f(u)}{u^{r-1}} .$$
For a fixed $\theta \in (N, N+\alpha]$ and any $\eps >0$ small, by Corollary \ref{corol_stima_Riesz} and \eqref{eq_cond_4sublin_f} (since $u(x) \to 0$ as $|x|\to +\infty$, being $u$ decreasing and in $L^1(\R^N)$) we obtain -- we use here that $\int_{\R^N} F(u)>0$ --
$$\big(I_{\alpha}*F(u)\big)f(u) \geq \underline{C} \left(\int_{\R^N} F(u)-\eps \norm{F(u)}_{\infty, \theta}\right)
I_{\alpha} u^{r-1} \quad \hbox{in $\R^N \setminus B_{R_{\eps}}(0)$}$$
for some $R_{\eps} \gg 0$, thus
$$ (-\Delta)^s u^{2-r} + \mu (2-r) u^{2-r} \geq %\left(
 (2-r) \underline{C}
\left(\int_{\R^N} F(u)-\eps \norm{F(u)}_{\infty,\theta}\right)
%\right)
 I_{\alpha} \quad \hbox{in $\R^N \setminus B_{R_{\eps}}(0)$};$$
that is, exploiting $\frac{1}{|x|^{N-\alpha}} \geq \frac{1}{(1+|x|^2)^{\frac{N-\alpha}{2}}}$, we get
$$(-\Delta)^s u^{2-r} + \lambda' u^{2-r} \geq \gamma'_{u,\eps} h_{N-\alpha} \quad \hbox{ in $\R^N \setminus B_{R_{\eps}}(0)$}$$
in the viscosity sense, where
$$\gamma'_{u,\eps} := (2-r) \underline{C} C_{N,\alpha} 
\left(\int_{\R^N} F(u)-\eps \norm{F(u)}_{\infty,\theta}\right) > 0$$
and
$$\lambda':=\mu (2-r).$$

We observe that $u^{2-r} \in L^{\infty}(B_R(0)) \cap C(B_{R_{\eps}}(0))$, while $h_{N-\alpha} \in L^{\infty}(\R^N) \cap C^{\sigma}_{loc}(\R^N)$
 (for any $\sigma$), thus by Lemma \ref{lem_lax_3milgr}, there exists $v\in C^{\omega}_{loc}(\R^N)$, for some $\omega>2s$ such that
$$\parag{ 
&(- \Delta)^s v + \lambda' v = \gamma'_{u,\eps} h_{N-\alpha} 
& \quad \hbox{ in $\R^N \setminus B_{R_{\eps}}(0)$}, \\ &v = u^{2-r} 
& \quad \hbox{on $B_{R_{\eps}}(0)$},}$$
pointwise. 
Thus
$$(-\Delta)^s (u^{2-r}-v) + \lambda' (u^{2-r}-v) \geq 0 \quad \hbox{ in $\R^N \setminus B_{R_{\eps}}(0)$}$$
in the viscosity sense, with
$$u^{2-r} - v \geq 0 \quad \hbox{on $B_{R_{\eps}}(0)$}.$$
Observe that, by Lemma \ref{lem_confr_3xbeta}, we have $v(x) \to 0$ as $|x|\to +\infty$.
Since $(u^{r-2}-v)(x) \to 0$ as $|x|\to +\infty$, by the viscosity version of the Comparison Principle (Lemma \ref{lem_comp_3prin}) we obtain
$$u^{2-r} \geq v \quad \hbox{on $\R^N$}.$$
By Lemma \ref{lem_confr_3xbeta} we gain
$$\liminf_{|x|\to +\infty} v(x) |x|^{N-\alpha} \geq \frac{\gamma'_{u,\eps}}{\lambda'} .$$
Combining the previous inequalities and 
sending $\eps \to 0^+$, we have the first claim. We conclude by adapting Remark \ref{rem_dec_2sN} to the viscosity case (notice that $u\in L^1(\R^N) \cap L^{\infty}(\R^N)$).
\QED

\bigskip

The above estimate applies, in particular, to pointwise solutions.

\begin{Corollary}
Assume \hyperref[(F1s)]{\textnormal{(F1)}}-\hyperref[(F2s)]{\textnormal{(F2,i)}} and the sublinear condition \hyperref[(F10)]{\textnormal{(F10)}}. 
Let $u\in L^1(\R^N) \cap C^{\gamma}_{loc}(\R^N)$ for some $\gamma>2s$, strictly positive, radially symmetric and decreasing, be a pointwise 
solution of \eqref{eq_main_prel}, such that $\int_{\R^N}F(u)>0$. 
Then the conclusions of Proposition \ref{prop_below_2} hold.
\end{Corollary}

By the results of the previous Sections (see Proposition \ref{prop_regolar}), we gain sufficient conditions in order to state that a weak solution is a pointwise solution. 

\begin{Corollary}\label{cor_est_bel_wea}
Assume \hyperref[(F1s)]{\textnormal{(F1)}}-\hyperref[(F2s)]{\textnormal{(F2,i)}} together with \hyperref[(F7)]{\textnormal{(F7)}}, and the sublinear condition \hyperref[(F10)]{\textnormal{(F10)}}. 
Let $u\in H^s(\R^N) \cap L^1(\R^N) \cap C(\R^N)$, strictly positive, radially symmetric and decreasing, be a weak solution
of \eqref{eq_main_prel}, such that $\int_{\R^N}F(u)>0$. %Assume moreover that $f \in C^{0,\sigma}_{loc}(\R)$ for some $\sigma \in (0,1)$. 
Then $u$ is a classical solution and 
 the conclusions of Proposition \ref{prop_below_2} hold.
\end{Corollary}

Notice that, by the sublinearity in zero, the H\"older exponent $\sigma$ can lie only in $(0, r-1]$. 
We conjecture anyway that the conclusion of Corollary \ref{cor_est_bel_wea} holds in more general cases, by assuming merely $f$ continuous.

\subsection{An $s$-sublinear threshold} %Some corollaries} %Proof of Theorem \ref{thm_main}}
\label{sec_proof_main_2}

We can sum up some of the results of the previous Sections in what follows.
\begin{Corollary}
Assume \hyperref[(F1s)]{\textnormal{(F1)}}-\hyperref[(F2s)]{\textnormal{(F2,i)}} and the sublinear conditions \hyperref[(F9)]{\textnormal{(F9)}}-\hyperref[(F10)]{\textnormal{(F10)}}, in particular
$$0< \liminf_{t \to 0} \frac{f(t)}{|t|^{r-1}} \leq \limsup_{t \to 0} \frac{f(t)}{|t|^{r-1}} \leq \bar{C} < +\infty.$$ 
Let $u\in H^s(\R^N) \cap L^1(\R^N) \cap C(\R^N)$, 
strictly positive, radially symmetric and decreasing, be a weak solution of \eqref{eq_main_prel}. 
Finally assume \hyperref[(F7)]{\textnormal{(F7)}}, i.e. 
$$f\in C^{0,\sigma}(\R) \; \hbox{ for some $\sigma \in (0,r-1]$},$$
and $\int_{\R^N}F(u)>0$. 
Then, if
$$\mu> (r-1) \bar{C}^{\frac{1}{r-1}}$$
we have
$$0 < \liminf_{|x| \to +\infty} u(x)|x|^{\beta} \leq \limsup_{|x|\to +\infty} u(x) |x|^{\beta} < +\infty$$
where 
$\beta:= \min\left\{ \frac{N-\alpha}{2-r}, N+2s\right\}.$
\end{Corollary}

We notice that, by assuming
$$\limsup_{t \to 0} \frac{f(t)}{|t|^{r-1}} \in (0, +\infty)$$
we obtain that
$$ \limsup_{t \to 0} \frac{f(t)}{|t|^{r-\eps-1}} =0 $$
for any $\eps>0$.
Thus we may directly extend the estimate from above to a whatever $\mu>0$ by paying the cost of a slower decay at infinity; this was essentially contained already in Remark \ref{rem_dim_induzione}. Notice that we still need $r-\eps \geq 2^{\#}_{\alpha}$.

\begin{Corollary}\label{corol_slower_deca2}
Assume \hyperref[(F1s)]{\textnormal{(F1)}}-\hyperref[(F2s)]{\textnormal{(F2,i)}} and the sublinear conditions \hyperref[(F9)]{\textnormal{(F9)}}-\hyperref[(F10)]{\textnormal{(F10)}}, in particular
$$0< \liminf_{t \to 0} \frac{f(t)}{|t|^{r-1}} \leq \limsup_{t \to 0} \frac{f(t)}{|t|^{r-1}} < +\infty$$ 
with $r \in (2^{\#}_{\alpha},2)$. 
Let $u\in H^s(\R^N) \cap L^1(\R^N)\cap C(\R^N)$, 
strictly positive, radially symmetric and decreasing, be a weak solution of \eqref{eq_main_prel}. 
Finally assume \hyperref[(F7)]{\textnormal{(F7)}}, i.e. 
$$f\in C^{0,\sigma}(\R) \; \hbox{ for some $\sigma \in (0,r-1]$},$$
and $\int_{\R^N}F(u)>0$. 
Then, if $\mu>0$ is arbitrary and $\eps>0$ is small, we have
$$0 < \liminf_{|x| \to +\infty} u(x)|x|^{\beta_0} \leq \limsup_{|x|\to +\infty} u(x) |x|^{\beta_{\eps}} < +\infty$$
where
$$\beta_{\eps}:= \min\left\{ \frac{N-\alpha}{2-r+\eps}, N+2s\right\}.$$
\end{Corollary}

We can now conclude the proof of the main theorem.

\medskip

\claim Proof of Theorem \ref{thm_main}.
First, we show how to remove the restriction on $\mu$ in Proposition \ref{prop_estim_above}. Indeed, for any $\kappa>0$ we can write $\big(I_{\alpha}*F(u)\big) f(u) \equiv \big(I_{\alpha}*F_{\kappa}(u)\big) f_{\eta}(u)$, where $f_{\kappa}:=\frac{1}{\kappa} f$ and $F_{\kappa}:=\kappa F$. We can thus rewrite \hyperref[(f3)]{\textnormal{(f3)}} as
$$|f_{\kappa} (t)|\leq \frac{1}{\kappa} \overline{C} t^{r-1} \quad \hbox{for $t \in (0,\delta)$}.$$
Since in Proposition \ref{prop_estim_above} we did not use that $F$ is the primitive of $f$ (in particular, we did not apply \hyperref[(f3)]{\textnormal{(f3)}} to $F$), fixed a whatever $\mu>0$ we can choose $\kappa$ such that
$$\mu > (r-1) \left(\frac{\overline{C}}{\kappa}\right)^{\frac{1}{r-1}} >0,$$
that is a large $\kappa$ given by $\kappa > \left(\frac{r-1}{\mu}\right)^{r-1} \overline{C}$, and obtain 
$$\limsup_{|x|\to +\infty} u(x) |x|^{\beta} \leq C_{u,\kappa}$$
where, if $\beta < N+2s$, 
\begin{equation*}\label{eq_constant_k}
C_{u,\kappa}:= \frac{(2-r) \left( C_{N,\alpha} \kappa \pabs{ \int_{\R^N} F(u)}\right)^{\frac{1}{2-r}} }{\mu -(r-1)\left(\frac{\bar{C}}{\kappa}\right)^{\frac{1}{r-1}} }.
\end{equation*}
We notice, as we expect, that as $\mu \to 0$ then $\kappa \to +\infty$ and $C_{u,\kappa}\to +\infty$, while $C'_u$ defined in Proposition \ref{prop_below_2} is invariant under $\kappa$-transformations.

We show now the sharp decay. Indeed, we search for a $\kappa$ such that $C_{u,\kappa}=C'_u$. By a straightforward analysis of $g(\kappa):=C_{u,\kappa}-C'_u$, $\kappa> \left(\frac{r-1}{\mu}\right)^{r-1} \overline{C}$, we find a (unique, explicit) zero $\kappa^*$ (which actually is a point of minimum) if only if $\overline{C}=\underline{C}$, i.e. if $f$ is exactly a power near the origin.

\smallskip

By the results of the previous Sections (Theorem \ref{th_INT_regular}, Proposition \ref{prop_u_holder}, Proposition \ref{prop_u_L1}), we have that every positive solution is bounded, and every bounded solution is in $H^{2s}(\R^N)\cap C(\R^N) \cap L^1(\R^N)$. 
%Thus the main theorem is a subcase of the previous results.
By the previous results we conclude the proof.
\QED

\bigskip

The conditions on $f$ in the previous results imply that $f$ is sublinear, but in a strict sense. We see that the results actually generalize to sublinear functions in a non strict sense.
\begin{Corollary}\label{corol_nonstrict}
Assume \hyperref[(F1s)]{\textnormal{(F1)}}-\hyperref[(F2s)]{\textnormal{(F2,i)}}. Assume moreover that $f$ is sublinear in a non-strict sense, i.e.
%$$\lim_{t\to 0} \frac{|f(t)|}{|t|}=+\infty$$
$$\lim_{t\to 0^+} \frac{f(t)}{t}=+\infty$$
but
$$\lim_{t\to 0} \frac{f(t)}{|t|^{r-1}}=0 \quad \hbox{for each $r \in (2^{\#}_{\alpha}, 2)$}.$$
Let $u\in H^s(\R^N) \cap L^1(\R^N)\cap C(\R^N)$, 
strictly positive, radially symmetric and decreasing, be a weak solution of \eqref{eq_main_prel}. 
Finally assume
$$f\in C^{0,\sigma}(\R) \; \hbox{ for some $\sigma \in (0,r-1]$}.$$
Then, if $\mu> 0$, we have
$$0 < \liminf_{|x| \to +\infty} u(x)|x|^{N+2s} \leq \limsup_{|x|\to +\infty} u(x) |x|^{N+2s} < +\infty.$$
\end{Corollary}
\claim Proof.
The estimate from below comes from the argument in Remark \ref{rem_dec_2sN} (since $f(t) \geq \underline{C} t$ for $t$ small and positive). The estimate from above comes from Proposition \ref{prop_estim_above}, after having chosing a whatever $r \in [r^*_{\alpha,s}, %\frac,ha+4s}{N+2s}, 
2)$.
\QED

\bigskip

\claim Proof of Corollary \ref{corol_main}.
By the results in the previous Sections (Theorem \ref{th_INT_positiv}), we have that every Pohozev minimum has constant sign -- e.g., it is strictly positive -- (if $f$ is odd or even, and H\"older continuous), and it is radially symmetric and decreasing (if in addition $f$ has constant sign on $(0,+\infty)$). 
Thus we conclude by the previous results.
\QED

\bigskip

All the previous theorems particularly apply to homogeneous nonlinearities $f(u)=|u|^{r-2} u$; notice that in this case we have $f \in C^{r-1}_{loc}(\R^N)$.

\begin{Corollary}\label{corol_homogen2}
Let $u\in H^s(\R^N)$, strictly positive, radially symmetric and decreasing, be a solution of
$$(-\Delta)^s u + \mu u = (I_{\alpha}*|u|^r)|u|^{r-2} u \quad \hbox{in $\R^N$}$$
with $r \in [2^{\#}_{\alpha}, 2)$.
Set, for every $\eps\geq 0$,
$$\beta_{\eps}:= \min\left\{ \frac{N-\alpha}{2-r+\eps}, N+2s\right\}.$$
We have
\begin{itemize}
\item if $\mu>r-1$ then
$$0 < \liminf_{|x| \to +\infty} u(x)|x|^{\beta_0} \leq \limsup_{|x|\to +\infty} u(x) |x|^{\beta_0} < +\infty;$$
\item if $r \in (2^{\#}_{\alpha},2)$ and $\mu \in (0, r-1]$ then, for any $\eps>0$ small,
$$0 < \liminf_{|x| \to +\infty} u(x)|x|^{\beta_0} \leq \limsup_{|x|\to +\infty} u(x) |x|^{\beta_{\eps}} < +\infty.$$
\end{itemize}
\end{Corollary}

\medskip

\claim Proof of Theorem \ref{thm_homog0_ws}, Corollary \ref{corol_homog0_gs} and Corollary \ref{corol_sharp_decay}.
Theorem \ref{thm_homog0_ws} is a direct consequence of the above result. 
By \cite[Theorems 3.2 and 4.2]{DSS1} we have that every ground state satisfies all the assumptions of the previous results; thus we have the claims of Corollary \ref{corol_homog0_gs} and Corollary \ref{corol_sharp_decay}.
% When $\mu=1$ and $r \in ( %2^{\#}_{\alpha}, 
%2^{\#}_{\alpha}, 
%r^*_{\alpha,s})$, the asymptotic coefficients appearing in Propositions \ref{prop_estim_above} and \ref{prop_below_2} are both equal to $\big( C_{N,\alpha} \norm{u}_r^r\big)^{\frac{1}{2-r}}$; this concludes the proof of Corollary \ref{corol_sharp_decay}.
\QED

%\medskip
%
%\begin{Remark}
%After the publication of this thesis and \cite{Gal1}, we were notified of the recent paper \cite{DeYa0} where the authors obtain the same asymptotic estimate as in Theorem \ref{thm_homog0_ws}, in the noncritical setting $r \in (2^{\#}_{\alpha}, 2)$: differently from here, where we employ a direct approach, in \cite{DeYa0} the authors employ an iterative approach, essentially in the spirit of Remark \ref{rem_dim_induzione}. Moreover, we achieve here the asymptotic constant (see Corollary \ref{corol_sharp_decay}).
%%\label{rem_regul_refined}
%\end{Remark}

%%%%%%%%%%%%%%

\section{The Pohozaev identity} % \tr{[taken from old notes, to be checked and completed]}} %Proof of Proposition \ref{prop_pohozaev} (Pohozaev identity)} 
\label{sec_Pohozaev}

%\small
%\emph{Disclaimer:} this Section was not present in the version of the Thesis filed at Università degli Studi di Bari (28/02/2023, even if the result was there announced), but it has been added before the dissertation of the Thesis itself (21/03/2023).
%
%\normalsize 
%
%\medskip

%In \cite[equation (6.1)]{DSS1}, in presence of power nonlinearities, it is proved that every weak solution $u$ is $C^2$ and thus satisfies the Pohozaev identity \eqref{eq_INT_Pohozaev}.
In \cite[equation (6.1)]{DSS1}, in presence of power nonlinearities, it is proved that every weak solution $u$ is $C^2$ and thus satisfies the Pohozaev identity \eqref{eq_INT_Pohozaev},
and this relation is extended to general superlinear nonlinearities $f \in C^1(\R)$ %which have a lower growth condition higher than $2$, 
in \cite[Proposition 2]{SGY}. %; in all these cases, the solutions are proved to be $u \in C^2(\R^N)$.
Here we want to further extend the identity to more general nonlinearities and to more general solutions $u\in C^1$, without employing the Caffarelli-Silvestre $s$-harmonic extension. 

This Section is mainly based on the paper \cite{CGT5}.

\smallskip 

First, we collect the results of the previous Sections to highlight the conditions that ensure the right regularity of the solutions.
\begin{Corollary}
\label{corol_reg_pohz}
Assume that \hyperref[(F1s)]{\textnormal{(F1)}}-\hyperref[(F2s)]{\textnormal{(F2)}} hold. 
Let $u\in H^s(\R^N)\cap L^{\infty}(\R^N)$
be a weak solution of \eqref{eq_introduction}. Assume in addition one of the following
\begin{itemize}
\item $s \in (\frac{1}{2}, 1)$,
\item $s=\frac{1}{2}$ and \hyperref[(F7)]{\textnormal{(F7)}},
\item $s \in [\frac{1}{4}, \frac{1}{2})$, $\alpha \in ( 1-2s, N)$ and \hyperref[(F7)]{\textnormal{(F7)}} with $\sigma \in (\frac{1-2s}{2s}, 1]$,
\item $s \in (0, \frac{1}{2})$, $\alpha \in (0,2)$ and \hyperref[(F6)]{\textnormal{(F7)}} with $\sigma \in (1-2s, 1]$.
\end{itemize}
Then $u \in C^{1,\gamma}(\R^N)$ for some $\gamma \in (0,1)$. If $s \in (\frac{1}{2}, 1)$ and \hyperref[(F7)]{\textnormal{(F7)}} holds too, then we can choose $\gamma \in (2s-1,1)$.
\end{Corollary}

Thus we want to prove the following result.

\begin{Theorem}
\label{thm_pohoz_frac_choq}
Let $u\in H^s(\R^N)\cap L^{\infty}(\R^N)$
be a weak solution of \eqref{eq_introduction}, and assume \hyperref[(F1s)]{\textnormal{(F1)}}-\hyperref[(F2s)]{\textnormal{(F2)}}.
Assume moreover \hyperref[(F7)]{\textnormal{(F7)}} and one of the following:
\begin{itemize}
\item $s \in [\frac{1}{2}, 1)$, % and \hyperref[(F6)]{\textnormal{(F6)}},
\item $s \in [\frac{1}{4}, \frac{1}{2})$, $\alpha \in ( 1-2s, N)$ and %\hyperref[(F6)]{\textnormal{(F6)}} with 
$\sigma \in (\frac{1-2s}{2s}, 1]$,
\item $s \in (0, \frac{1}{2})$, $\alpha \in (0,2)$ and %\hyperref[(F6)]{\textnormal{(F6)}} with 
$\sigma \in (1-2s, 1]$.
\end{itemize}
Then $u \in C^{1,\gamma}(\R^N)$ for some $\gamma \in (\max\{0,2s-1\},1) $, and $u$ satisfies the Pohozaev identity \eqref{eq_INT_Pohozaev}, 
%\begin{equation}\label{eq:2.5}
%\frac{N-2s}{2} \|(-\Delta)^{s/2} u\|_2^2 + \frac{N }{ 2} e^{\lambda} \|u\|_2^2 - \frac{N + \alpha}{2} \, \mc{D}(u) =0
%\end{equation}
or equivalently
$$\frac{1}{2^*_{\alpha,s}} \|(-\Delta)^{s/2} u\|_2^2 + \frac{\mu}{2^{\#}_{\alpha}} \|u\|_2^2 - \mc{D}(u)=0.$$
%where $\mc{D}(u) = \int_{\R^N} (I_\alpha*F(u))F(u)$. 
The result in particular applies to positive weak solutions $u \in H^s(\R^N)$ of \eqref{eq_introduction}.
\end{Theorem}

%The Theorem will be a consequence of the following statements.
%\smallskip

We start by the following integration by parts rule, inspired by \cite[Lemma 4.2]{DFW1}, obtained 
under a pointwise well posedness of the fractional Laplacian and the existence of a weak gradient. %, we obtain the. %, DFW2}). 

\begin{Proposition}
\label{prop_int_kern_div}
Let $s\in (0,1)$. 
%Let $F \in Lip_{loc}(\R^N)$ with $F(u) \in \R^N$ and $F'=f$.
Let $u \in \dot{H}^s(\R^N) \cap C^{\gamma}_{loc}(\R^N) \cap Lip_{loc}(\R^N)$ for some $\gamma >2s$, and assume \eqref{eq_spazio_gener_s}.
%Assume moreover that
%$$(-\Delta)^s u = f(u) \quad \hbox{ in $\R^N$}$$
%in the pointwise sense. %, with $F(u) \in L^1_{loc}(\R^N)$. 
Let moreover $X \in C^1_c %Lip_c
(\R^N,\R^N)$ %(\tor{also $C^1_c(\R^N)$ or more is okay, if needed}) 
be a vector field, and define, for $x, y\in \R^N$, $x \neq y$, 
%$$\mc{K}_X(x,y):= \dive\big(X(x)) + X(y)\big) - (N+2s) \frac{(X(x)-X(y))\cdot(x-y)}{|x-y|^2}$$
%$$\mc{K}_X(x,y):= \frac{C_{N,s}}{4} \left( (N+2s) \frac{(X(x)-X(y))\cdot(x-y)}{|x-y|^2} -\dive\big(X(x)) + X(y)\big)\right)$$
%$$\mc{K}_{X}^s(x,y):= \frac{\dive\big(X(x) + X(y)\big)}{2}- \frac{N+2s}{2}\frac{(X(x)-X(y))\cdot(x-y)}{|x-y|^2}$$
$$\mc{K}_{X}^s(x,y):= \frac{\big(\dive (X)\big)(x) + \big(\dive (X)\big)(y)}{2}- \frac{N+2s}{2}\frac{(X(x)-X(y))\cdot(x-y)}{|x-y|^2}$$
%$$\mc{K}_{X}^s(x,y):= C'_{N,s} \left( \frac{\dive\big(X(x) + X(y)\big)}{N+2s}- \frac{(X(x)-X(y))\cdot(x-y)}{|x-y|^2}\right)$$
the \emph{fractional divergence kernel} related to $X$. %here $C'_{N,s}:= \frac{N+2s}{4}$. %C_{N,s}$. 
Then it holds
$$%\frac{C_{N,s}}{4}
\frac{C_{N,s}}{2}\int_{\R^N} \int_{\R^N} \frac{|u(x)-u(y)|^2}{|x-y|^{N+2s}} \mc{K}_{X}^s(x,y) dx dy = -
 \int_{\R^N} (-\Delta)^s u \, (\nabla u \cdot X ) \, dx;$$
noticed that the left-hand side is the weighted Gagliardo seminorm with weight $\mc{K}_{X}^s$, set 
$$\mc{G}_u^s(x,y):= \frac{C_{N,s}}{2} \frac{|u(x)-u(y)|^2}{|x-y|^{N+2s}}$$
we can write
$$%[u]_{\R^N, \mc{K}_{X}^s}^2 =% - \frac{4}{C_{N,s}}
(\mc{G}_u^s, \mc{K}_X^s)_{L^2(\R^{2N})} =
- \big((-\Delta)^{s} u \, \nabla u,% \nabla u \cdot 
X\big)_{L^2(\R^N)}.$$
\end{Proposition}

\claim Proof. % of Proposition \ref{prop_int_kern_div}.
For the proof, we follow the lines of \cite{Dje0}. % (see also \cite{DeNDj}).
%For the reader's convenience, we write here the details. 
%
We start noticing that, being $u \in \dot{H}^s(\R^N)$, by the assumptions we have
$$\mc{G}_u^s \in L^1(\R^{2N}), \quad \mc{K}^s_X \in L^{\infty}(\R^{2N})$$
so that the product is summable. By dominated convergence theorem, the symmetry of the kernel, and the Fubini theorem, we obtain
\begin{eqnarray*}
%\frac{1}{C'_{N,s}} 
\lefteqn{
\frac{2}{C_{N,s}}(\mc{G}_u^s, \mc{K}_X^s)_{L^2(\R^{2N})} } \\
%[u]_{\R^N, \mc{K}^s_X}^2 
&=& %\frac{1}{C'_{N,s}}
 \lim_{\eps \to 0} \iint_{|x-y|>\eps} \frac{|u(x)-u(y)|^2}{|x-y|^{N+2s}} \mc{K}_{X}^s(x,y) dx dy \\
&=& %2 \lim_{\eps \to 0} \iint_{|x-y|>\eps}\frac{|u(x)-u(y)|^2}{|x-y|^{N+2s}} \left( \frac{\dive(X)(x)}{N+2s} - \frac{(x-y)\cdot X(x)}{|x-y|^2} \right)dxdy \\
 \lim_{\eps \to 0} \iint_{|x-y|>\eps}\frac{|u(x)-u(y)|^2}{|x-y|^{N+2s}} \left( \dive(X)(x) - (N+2s) \frac{(x-y)\cdot X(x)}{|x-y|^2} \right)dxdy \\
&=& %2 \lim_{\eps \to 0} \int_{\R^N} \left( \int_{\R^N \setminus B_{\eps}(y)}\frac{|u(x)-u(y)|^2}{|x-y|^{N+2s}} \left( \frac{\dive(X)(x)}{N+2s} - \frac{(x-y)\cdot X(x)}{|x-y|^2} \right) dx \right) dy.
 \lim_{\eps \to 0} \int_{\R^N} \left( \int_{\R^N \setminus B_{\eps}(y)}\frac{|u(x)-u(y)|^2}{|x-y|^{N+2s}} \left( \dive(X)(x)-\right. \right. \\ 
&& \left. \left. - (N+2s)\frac{(x-y)\cdot X(x)}{|x-y|^2} \right) dx \right) dy.
\end{eqnarray*}
Exploiting that, for $x\neq y$, $\nabla_x \frac{1}{|x-y|^{N+2s}} =- (N+2s) \frac{x-y}{|x-y|^{N+2s+2}}$, and the divergence theorem (possible because $\frac{X}{|\cdot-y|^{N+2s}} \in C^1_c % Lip_c
(\overline{\R^N \setminus B_{\eps}(y)})$ and $u \in Lip_{loc}(\R^N) \subset W^{1,\infty}(\supp(X))$, see \cite[Theorem 4.6]{EG0})
%see \href{https://encyclopediaofmath.org/wiki/Integration_by_parts}{Link} or [Evans-Gariepy, Theorem 4.6])
\begin{eqnarray*}
%\frac{1}{2C'_{N,s}}
\lefteqn{
\tfrac{2}{C_{N,s}} (\mc{G}_u^s, \mc{K}_X^s)_{L^2(\R^{2N})} } \\ % [u]_{\R^N, \mc{K}^s_X}^2 
%=& \lim_{\eps \to 0} \int_{\R^N} \left(\int_{\R^N \setminus B_{\eps}(y)}|u(x)-u(y)|^2 \left( \frac{\dive(X)(x)}{(N+2s)|x-y|^{N+2s}} - \frac{(x-y)\cdot X(x)}{|x-y|^{N+2s+2}} \right) dx \right) dy \\
&=& \lim_{\eps \to 0} \int_{\R^N} \left(\int_{\R^N \setminus B_{\eps}(y)}|u(x)-u(y)|^2 \left( \frac{\dive(X)(x)}{|x-y|^{N+2s}} - \right. \right. \\
&& \left. \left. -(N+2s)\frac{(x-y)\cdot X(x)}{|x-y|^{N+2s+2}} \right) dx \right) dy \\
&=&%\frac{1}{N+2s}
\lim_{\eps \to 0} \int_{\R^N} \left(\int_{\R^N \setminus B_{\eps}(y)}|u(x)-u(y)|^2 \dive_x \left( \frac{X}{|x-y|^{N+2s}}\right)(x) dx \right) dy \\
%=&\frac{1}{N+2s}\lim_{\eps \to 0} \int_{\R^N} \left(\int_{|x-y|>\eps}\nabla_y \big(|u(x)-u(y)|^2\big) \frac{X}{|x-y|^{N+2s}}(y) dy \right) dx + \\
%& +\frac{1}{N+2s}\lim_{\eps \to 0} \int_{\R^N} \left(\int_{|x-y|=\eps} |u(x)-u(y)|^2 \frac{X}{|x-y|^{N+2s}}(y) d\sigma(y) \right) dx
&=&- 2 %\frac{2}{N+2s}
\lim_{\eps \to 0} \int_{\R^N} \left(\int_{\R^N \setminus B_{\eps}(y)}(u(x)-u(y)) \nabla u(x) \cdot \frac{X(x)}{|x-y|^{N+2s}} dx \right) dy + \\
&& +%\frac{1}{N+2s}
\lim_{\eps \to 0} \int_{\R^N} \left(\int_{\partial B_{\eps}(y)} |u(x)-u(y)|^2 \frac{X(x)}{|x-y|^{N+2s}} \cdot \frac{x-y}{|x-y|} d\sigma(x) \right) dy \\
&=&- 2 %\frac{2}{N+2s}
\lim_{\eps \to 0} \int_{\R^N} \left(\int_{\R^N \setminus B_{\eps}(y)} \frac{u(x)-u(y)}{|x-y|^{N+2s}} \nabla u(x) \cdot X(x) dx \right) dy + \\
&& + %\frac{1}{N+2s}
\lim_{\eps \to 0} \frac{1}{\eps^{N+2s+1}} \int_{\R^N} \left(\int_{\partial B_{\eps}(y)} |u(x)-u(y)|^2 X(x) \cdot (x-y) d\sigma(x) \right) dy.
\\
&=:& \, -2 \lim_{\eps \to 0} I_{\eps} + \lim_{\eps \to 0} E_{\eps};
\end{eqnarray*}
here we split the limits since we will prove the existence of both. 

For the first integral, we notice that $x \mapsto \int_{\R^N \setminus B_{\eps}(x)} \frac{|u(x)-u(y)|}{|x-y|^{N+2s}} |\nabla u(x) \cdot X(x)| dy \leq C_{\eps} \norm{u}_{\infty} |\nabla u (x)\cdot X(x)| \in L^1(\R^N)$
%$\nabla u \cdot X \in L^{\infty}(\R^{2N})$ with compact support \tr{and \tr{... }
%$$\frac{u(x)-u(y)}{|x-y|^{N+2s}} \leq \frac{1}{|x-y|^{N+2s-1}} \in L^1(|x-y|>\eps), \quad \nabla u \cdot X \in L^{\infty}(\R^{2N}),$$
so that we can apply Fubini theorem, % \tr{(?)},} 
then we perform a symmetrization substitution and apply again Fubini theorem, %(possible because $x \mapsto \int_{\R^N \setminus B_{\eps}(x)} \frac{ |2u(x)-u(x+y)-u(x-y)|}{|y|^{N+2s}} |\nabla u(x) \cdot X(x)| dy \in L^1(\R^N)$) 
and finally dominated convergence theorem (since $y \mapsto \frac{ 2u(x)-u(x+y)-u(x-y)}{|y|^{N+2s}} \in L^1(\R^N)$ by Proposition \ref{prop_diff_represent}), obtaining 
\begin{eqnarray*}
%-\frac{N+2s}{2}
\lefteqn{
C_{N,s} \lim_{\eps \to 0} I_{\eps}} \\
 &=& C_{N,s} \lim_{\eps \to 0} \iint_{|x-y|>\eps} \frac{u(x)-u(y)}{|x-y|^{N+2s}} \nabla u(x) \cdot X(x) dx dy \\
%&\tr{\; \stackrel{?}= \, } \lim_{\eps \to 0} \int_{\R^N} \nabla u(x) \cdot X(x) \left(\int_{\R^N \setminus B_{\eps}(x)} \frac{u(x)-u(y)}{|x-y|^{N+2s}} dy \right) dx \\
%\int_{\R^N} \left(\int_{\partial B_{\eps}(y)} |u(x)-u(y)|^2 \frac{X(x)}{|x-y|^{N+2s}} \cdot \frac{x-y}{|x-y|} d\sigma(x) \right) dy
%&= \lim_{\eps \to 0} \int_{\R^N} \left(\int_{\R^N \setminus B_{\eps}(0)} \frac{2u(x)-u(x+y)-u(x-y)}{|x|^{N+2s}} \nabla u(x) \cdot X(x) dx \right) dy \\
&=& \frac{C_{N,s}}{2}\lim_{\eps \to 0} \iint_{|y|>\eps} \frac{2u(x)-u(x+y)-u(x-y)}{|y|^{N+2s}} \nabla u(x) \cdot X(x) dx dy \\
&=& \lim_{\eps \to 0} \int_{\R^N} \nabla u(x) \cdot X(x) \left(\frac{C_{N,s}}{2}\int_{\R^N \setminus B_{\eps}(0)} \frac{ 2u(x)-u(x+y)-u(x-y)}{|y|^{N+2s}} dy \right) dx \\
%&\tr{\stackrel{?}=} \frac{1}{C_{N,s}} \int_{\R^N} \nabla u(x) \cdot X(x) \left( \lim_{\eps \to 0} C_{N,s}\int_{\R^N \setminus B_{\eps}(x)} \frac{u(x)-u(y)}{|x-y|^{N+2s}} dy \right) dx \\
&=& %\frac{1}{C_{N,s}} 
\int_{\R^N} \nabla u \cdot X (-\Delta)^s u.
\end{eqnarray*}

For the second integral, notice that 
%if $\supp(X)\subset B_R(0)$, then the set $\{(x,y) \in \R^{2N} \mid x \in B_R(0), |x-y|=\eps\}$
 the set $\{(x,y) \in \R^{2N} \mid x \in \supp(X), \, |x-y|=\eps\}$
%$$A_{R,\eps}:=\{(x,y) \mid x \in B_R(0), |x-y|=\eps\}$$
is bounded.
%and its measure is of order $\sim \eps^N$ \tr{MMM forse no, forse $\sim \eps^{N-1}$, quindi dobbiamo sfruttare anche la Lipschitzianeità di $X$}); 
Thus the integrand (being bounded) is summable, which allows us to implement the Fubini theorem and obtain, by exploiting also a symmetrization argument, % of the argument,
%\begin{align*}
%(N+2s)E_{\eps} =& \lim_{\eps \to 0}\frac{1}{\eps^{N+2s+1}} \iint_{A_{R,\eps}} |u(x)-u(y)|^2 X(x) \cdot (x-y) d\sigma(x) \times dy \\
%\lesssim & \lim_{\eps \to 0}\frac{C}{\eps^{N+2s+1}} \iint_{A_{R,\eps}} |x-y|^{2\gamma+1} d\sigma(x) \times dy \\
%\lesssim & \lim_{\eps \to 0} \eps^{2(\gamma-s)} =0.
%\end{align*}
\begin{align*}
(N+2s)E_{\eps} &= %\lim_{\eps \to 0}
\frac{1}{\eps^{N+2s+1}} \iint_{|x-y|=\eps} |u(x)-u(y)|^2 X(x) \cdot (x-y) d\sigma(x) \times dy \\
&= % \lim_{\eps \to 0}
\frac{1}{2\eps^{N+2s+1}} \iint_{|x-y|=\eps}|u(x)-u(y)|^2 (X(x)-X(y)) \cdot (x-y) d\sigma(x) \times dy .
\end{align*}
If $\supp(X)\subset B_R(0)$, then out of the set
$$A_{R,\eps}:=\{(x,y) \in B_R(0)\times B_R(0) \mid |x-y|=\eps\}$$
the integrand is null. 
Thus, being $u \in Lip_{loc}(\R^N)$ (actually it is sufficient $u \in C^{0,\theta}(\R^N)$ for some $\theta>s$) % (actually $\theta=1$)
and $X \in Lip(\R^N, \R^N)$, we get
\begin{align*}
E_{\eps} &\lesssim %\lim_{\eps \to 0}
\frac{1}{\eps^{N+2s+1}} \iint_{A_{R,\eps}} |x-y|^4 d\sigma(x) \times dy \\
&= % \lim_{\eps \to 0}
\eps^{-N-2s +3} m_{2N-1}(A_{R,\eps}).
\end{align*}
Observed that $m_{2N-1}(A_{R,\eps}) \lesssim m_N(B_R) m_{N-1}(\partial B_{\eps}) \sim \eps^{N-1}$, we obtain $E_{\eps} \lesssim \eps^{-2s+2} \to 0$. %(being $\theta>s$) that the error $E_{\eps} \to 0$.
Joining the pieces, we reach the claim.
\QED

\medskip

\begin{Corollary}
In the assumptions of Proposition \ref{prop_int_kern_div}, let $G \in C^1%Lip_{loc}
(\R^N)$ with $G(u) \in L^1(\R^N)$ and 
$$(-\Delta)^s u = g(u) \quad \hbox{ in $\R^N$}$$
in the pointwise sense, where $G'=g$. Then 
\begin{align*} %\frac{C_{N,s}}{4}
\frac{C_{N,s}}{2} \int_{\R^N} \int_{\R^N} \frac{|u(x)-u(y)|^2}{|x-y|^{N+2s}} \mc{K}_{X}^s(x,y) dx dy &= -
 \int_{\R^N} \nabla G(u) \cdot X \, dx \\
 &= \int_{\R^N} G(u) \, \dive(X) dx, 
 \end{align*}
 i.e.
 $$%[u]_{\R^N, \mc{K}_{X}^s}^2 =% - \frac{4}{C_{N,s}}
(\mc{G}_u^s, \mc{K}_X^s)_{L^2(\R^{2N})} =
(G(u), \dive(X))_{L^2(\R^N)}.$$
\end{Corollary}

\smallskip

We deal now with the Riesz kernel right-hand side of the equation.

\begin{Proposition}
\label{prop_int_riesz_div}
%assume also $u \in H^{1}(\R^N)$ (e.g., $u \in H^{2s}(\R^N)$)
Let $\alpha \in (0,N)$ and $H \in Lip_{loc}(\R^N)\cap L^{\infty}(\R^N)$ be such that 
$$ %I_{\alpha}*|H|<\infty, \quad 
(I_{\alpha}*|H|)|H| \in L^1(\R^N), \quad (I_{\alpha}*|H|)|\nabla H| \in L^1_{loc}(\R^N).$$
%$H \in W^{1,1}(\R^N)$ %and $F'=f$ and assume \eqref{eq_Choquard_genericaF} 
%$$(-\Delta)^s u + = f(u) \quad \hbox{ in $\R^N$}$$
%is satisfied in the pointwise sense. 
Let moreover $X \in C^1_c %Lip_c
(\R^N,\R^N)$ be a vector field %(\tor{also $C^1_c(\R^N)$ or more is okay, if needed}) 
and set, for $x, y\in \R^N$, $x \neq y$, %(with an abuse of notation)
%$$\mc{K}_{X}^{\alpha}(x,y):= C'_{N,\alpha} \left( \frac{\dive\big(X(x) + X(y)\big)}{N-\alpha}- \frac{(X(x)-X(y))\cdot(x-y)}{|x-y|^2}\right)$$
%$$\mc{K}_{X}^{\alpha}(x,y):= \frac{\dive\big(X(x) + X(y)\big)}{2}- \frac{N-\alpha}{2} \frac{(X(x)-X(y))\cdot(x-y)}{|x-y|^2}.$$
$$\mc{K}_{X}^{-\frac{\alpha}{2}}(x,y):= \frac{\big(\dive (X)\big)(x) + \big(\dive (X)\big)(y)}{2}- \frac{N-\alpha}{2} \frac{(X(x)-X(y))\cdot(x-y)}{|x-y|^2}.$$
%where $C'_{N,\alpha}:= \frac{N-\alpha}{2}$. % C_{N,\alpha}$. %\frac{N-\alpha}{4} C_{N,\alpha}$.
Then 
\begin{align*} %\frac{C_{N,s}}{4}
\int_{\R^N} \int_{\R^N} I_{\alpha}(x-y) H(x) H(y) \mc{K}_{X}^{-\frac{\alpha}{2}}(x,y) dx dy &= -
 \int_{\R^N}\big (I_{\alpha}*H\big) \nabla H \cdot X \, dx, %\\
% &= \int_{\R^N} G(u) \, \dive(X)
 \end{align*}
 i.e. set
$$\mc{R}_H^{\alpha}(x,y):= I_{\alpha}(x-y) H(x) H(y)$$
we have
$$ (\mc{R}_{H}^{\alpha}, \mc{K}_{X}^{-\frac{\alpha}{2}})_{L^2(\R^{2N})} = %[u]_{\R^N, \mc{K}_{X}^s}^2 =% - \frac{4}{C_{N,s}}
%- (\nabla F(u), X)_2 =
- \big((I_{\alpha}*H) \nabla H, %\nabla F(u) \cdot 
X\big)_{L^2(\R^N)}.
% \int_{\R^N}\big (I_{\alpha}*F(u)\big) \nabla F(u) \cdot X \, dx.
$$
 % (G(u), \dive(X))_{L^2(\R^N)}.$$
\end{Proposition}

\claim Proof. 
We proceed as in the proof of Proposition \ref{prop_int_kern_div}.
% (\tor{%if correct, 
%we can cut something here, by writing down only the justifications for the formal passages ($... \in L^1$ etc). Or viceversa, we can keep this proof and cut a bit the one about the fractional Laplacian.}).
We start noticing that %, being $u \in ... $, %\dot{H}^s(\R^N)$,
$$\mc{R}_{H}^{\alpha} \in L^1(\R^{2N}), \quad \mc{K}_{X}^{-\frac{\alpha}{2}} \in L^{\infty}(\R^{2N})$$
by the assumptions, so that the product is summable. Thus
%By dominated convergence theorem, the symmetry of the kernel, and the Fubini theorem, we have
\begin{eqnarray*}
%\frac{1}{C'_{N,\alpha}}
\lefteqn{ 
 (\mc{R}_{H}^{\alpha}, \mc{K}_{X}^{-\frac{\alpha}{2}})_{L^2(\R^{2N})} } \\ %[u]_{\R^N, \mc{K}^s_X}^2 
%&=%\frac{1}{C'_{N,\alpha}}
% \lim_{\eps \to 0} \iint_{|x-y|>\eps} I_ {\alpha}(x-y) H(x) H(y) \mc{K}_{X}^{-\frac{\alpha}{2}}(x,y) dx dy \\
%&= %2 \lim_{\eps \to 0} \iint_{|x-y|>\eps} I_ {\alpha}(x-y) H(x) H(y) \left( \frac{\dive(X)(x)}{N-\alpha} - \frac{(x-y)\cdot X(x)}{|x-y|^2} \right)dxdy \\
% \lim_{\eps \to 0} \iint_{|x-y|>\eps} I_ {\alpha}(x-y) H(x) H(y) \left( \dive(X)(x) - (N-\alpha)\frac{(x-y)\cdot X(x)}{|x-y|^2} \right)dxdy \\
&=&%2 \lim_{\eps \to 0} \int_{\R^N} \left( \int_{\R^N \setminus B_{\eps}(y)} I_ {\alpha}(x-y) H(x) H(y)\left( \frac{\dive(X)(x)}{N-\alpha} - \frac{(x-y)\cdot X(x)}{|x-y|^2} \right) dx \right) dy.
 \lim_{\eps \to 0} \int_{\R^N} \left( \int_{\R^N \setminus B_{\eps}(y)} I_ {\alpha}(x-y) H(x) H(y)\left( \dive(X)(x) - \right. \right. \\
&& \left. \left. - (N-\alpha)\frac{(x-y)\cdot X(x)}{|x-y|^2} \right) dx \right) dy.
\end{eqnarray*}
Since $H \in Lip_{loc}(\R^N) \subset W^{1,\infty}(\supp(X))$, we have
\begin{eqnarray*}
%\frac{1}{2C'_{N,\alpha}} 
\lefteqn{\frac{1}{C_{N,\alpha}} (\mc{R}_{H}^{\alpha}, \mc{K}_{X}^{-\frac{\alpha}{2}})_{L^2(\R^{2N})} } \\% \\
&=& - %\frac{1}{N-\alpha}
\lim_{\eps \to 0} \int_{\R^N} \left(\int_{\R^N \setminus B_{\eps}(y)}\frac{1}{|x-y|^{N-\alpha}} H(y) \nabla H(x) \cdot X(x) dx \right) dy + \\
&& + %\frac{1}{N-\alpha}
\lim_{\eps \to 0} \frac{1}{\eps^{N-\alpha+1}} \int_{\R^N} \left(\int_{\partial B_{\eps}(y)} H(x) H(y) X(x) \cdot (x-y) d\sigma(x) \right) dy.
\\
&=:& \, %\frac{1}{N-\alpha} \left(-I_{\eps} +E_{\eps}\right);
-\lim_{\eps \to 0} I_{\eps} +\lim_{\eps \to 0} E_{\eps}. %;
\end{eqnarray*}
%here we split the limits since we will prove the existence of both. 

For $I_{\eps}$ %the first integral, 
we notice that $(I_{\alpha}*|H|)|\nabla H| |X| \in L^1(\R^N)$, thus 
$ (x,y) \mapsto I_{\alpha}(x-y) H(y) \nabla H(x) \cdot X(x) \in L^1(\R^{2N})$
and we can apply (twice) Fubini theorem; moreover $I_{\alpha}(x-\cdot) H \in L^1(\R^N)$, and we can apply dominated convergence theorem.
Hence we obtain
\begin{align*}
C_{N,\alpha} \lim_{\eps \to 0} I_{\eps} &= \lim_{\eps \to 0} \iint_{|x-y|>\eps} I_{\alpha}(x-y) H(y) \nabla H(x) \cdot X(x) dx dy \\
&= \lim_{\eps \to 0} \int_{\R^N} \nabla H(x)\cdot X(x) \left(\int_{\R^N \setminus B_{\eps}(x)} I_{\alpha}(x-y) H(y) dy \right) dx \\
&= \int_{\R^N} \nabla H(x) \cdot X(x) \left( \lim_{\eps \to 0}\int_{\R^N \setminus B_{\eps}(x)} I_{\alpha}(x-y) H(y) dy \right) dx \\
&= \int_{\R^N} \nabla H(x) \cdot X \, (I_{\alpha}*H).
\end{align*}

We can write $E_{\eps}$ instead as
\begin{align*}
E_{\eps} 
%&= %\lim_{\eps \to 0}
%\frac{1}{\eps^{N-\alpha+1}} \iint_{|x-y|=\eps} H(x) H(y) X(x) \cdot (x-y) d\sigma(x) \times dy \\
&= % \lim_{\eps \to 0}
\frac{1}{2\eps^{N-\alpha+1}} \iint_{|x-y|=\eps} H(x) H(y) (X(x)-X(y)) \cdot (x-y) d\sigma(x) \times dy .
\end{align*}
If $\supp(X)\subset B_R(0)$, set $A_{R,\eps}:=\{(x,y) \in B_R(0)\times B_R(0) \mid |x-y|=\eps\}$
and observed that $H \in L^{\infty}(\R^N)$, % and $X \in Lip(\R^N, \R^N)$ 
we obtain %$u \in C^{0,\theta}(\R^N)$ for some $\theta$ (actually $\theta=1$)
$$ E_{\eps} \lesssim \frac{1}{\eps^{N-\alpha+1}} \iint_{A_{R,\eps}} |x-y|^2 d\sigma(x) \times dy = \eps^{-N+\alpha+1} m_{2N-1}(A_{R,\eps}) \lesssim \eps^{\alpha} \to 0,$$
being $\alpha>0$. This concludes the proof.
%Since 
%Observed that $m_{2N-1}(A_{R,\eps}) \lesssim m_N(B_R) m_{N-1}(\partial B_{\eps}) \sim \eps^{N-1}$, we gain (being $\alpha>0$) that the error $E_{\eps} \to 0$.
%Joining the pieces, we get the claim.
\QED

\medskip

\begin{Theorem}
\label{thm_pohoz_X}
Let $s\in (0,1)$ and $\alpha \in (0,N)$ and assume that \hyperref[(F1s)]{\textnormal{(F1)}}-\hyperref[(F2s)]{\textnormal{(F2)}} hold. 
%Let $F \in Lip_{loc}(\R^N)$ with $F(u) \in \R^N$ and $F'=f$.
Let $u \in H^s(\R^N) \cap C^{\gamma}_{loc}(\R^N) \cap Lip_{loc}(\R^N)$ for some $\gamma >2s$, be 
%Assume that $u$ is 
a pointwise solution of \eqref{eq_introduction}. %, and that 
Then $u$ satisfies the Pohozaev identity \eqref{eq_INT_Pohozaev}.
\end{Theorem}

\claim Proof.
We apply Proposition \ref{prop_int_kern_div} and Proposition \ref{prop_int_riesz_div} with $H=F(u)$; notice that the assumptions on $u$ and $F$ imply the needed conditions on $H$ % integrability conditions 
(in particular we highlight that $f(u) \in L^{\frac{2N}{N+\alpha}}_{loc}(\R^N)$). Thus, for a generic $X \in C^1_c(\R^N, \R^N)$ we obtain
\begin{align*}
(\mc{G}_u^s, \mc{K}_X^s)_{L^2(\R^{2N})} & =- \big((-\Delta)^{s} u \, \nabla u,X\big)_{L^2(\R^N)} \\
&= \mu \big(u \, \nabla u,X\big)_{L^2(\R^N)} - \big((I_{\alpha}*F(u)) f(u) \, \nabla u,X\big)_{L^2(\R^N)} \\
&= \frac{\mu}{2} \big(\nabla(u^2),X\big)_{L^2(\R^N)} - \big((I_{\alpha}*F(u)) \nabla F(u),X\big)_{L^2(\R^N)} \\
&=-\frac{\mu}{2} \big(u^2,\dive(X)\big)_{L^2(\R^N)} + (\mc{R}_{F(u)}^{\alpha}, \mc{K}_{X}^{-\frac{\alpha}{2}})_{L^2(\R^{2N})}.
\end{align*}
In particular, we apply the result to 
$$X_n(x):= \varphi_n(x) \, x,$$
where $\varphi_n$ is a cut-off function with $\varphi_n \equiv 1$ in $B_n(0)$, $\supp(\varphi_n) \subset B_{n+1}(0)$, $\norm{\varphi_n}_{\infty}=1$ and $|x| |\nabla \varphi_n(x)|\leq C$ for each $x \in \R^N$ and $n \in \N$; for instance, defined such $\varphi_1$, we can set $\varphi_n:=\varphi_1(\cdot/n)$ and obtain
$$|x| |\nabla \varphi_n(x)|= |x/n| |\nabla \varphi_1 (x/n)| \leq \norm{|x| |\nabla \varphi_1(x)|}_{\infty}.$$
In particular, $x\mapsto x \varphi_n(x)$ is equi-Lipschitz. Noticed that $\dive(X_n) = N \varphi_n + \nabla \varphi_n \cdot x $
we gain
\begin{eqnarray*}
\lefteqn{ (\mc{G}_u^s, \mc{K}_{X_n}^s)_{L^2(\R^{2N})} } \\
%=& \, C'_{N,s}\int_{\R^N} \int_{\R^N} \frac{|u(x)-u(y)|^2}{|x-y|^{N+2s}} \left( \frac{\dive\big(X(x) + X(y)\big)}{N+2s}- \frac{(X(x)-X(y))\cdot(x-y)}{|x-y|^2}\right) dx dy \\
&=& %\, C'_{N,s}\int_{\R^N} \int_{\R^N} \frac{|u(x)-u(y)|^2}{|x-y|^{N+2s}} \left( \frac{N(\varphi_n(x)+\varphi_n(y))}{N+2s} - \frac{(\varphi_n(x)\, x - \varphi_n(y) y)\cdot(x-y)}{|x-y|^2}\right) dx dy +\\
%\frac{C_{N,s}}{2}\int_{\R^N} \int_{\R^N} \frac{|u(x)-u(y)|^2}{|x-y|^{N+2s}} \left( \frac{N(\varphi_n(x)+\varphi_n(y))}{2} - \frac{N+2s}{2}\frac{(\varphi_n(x)\, x - \varphi_n(y) y)\cdot(x-y)}{|x-y|^2}\right) dx dy +\\
\frac{C_{N,s}}{2}\int_{\R^N} \int_{\R^N} \frac{|u(x)-u(y)|^2}{|x-y|^{N+2s}}\left( \frac{N(\varphi_n(x)+\varphi_n(y))}{2} \right)dx dy -\\
&&-\frac{C_{N,s}}{2}\int_{\R^N} \int_{\R^N} \frac{|u(x)-u(y)|^2}{|x-y|^{N+2s}} \left( \frac{N+2s}{2}\frac{(\varphi_n(x)\, x - \varphi_n(y) y)\cdot(x-y)}{|x-y|^2} \right) dx dy +\\
%&+C'_{N,s}\int_{\R^N} \int_{\R^N} \frac{|u(x)-u(y)|^2}{|x-y|^{N+2s}} \left(\frac{\nabla \varphi_n(x) \cdot x + \nabla \varphi_n(y) \cdot y}{N+2s} \right) dx dy \\
&&+\frac{C_{N,s}}{2}\int_{\R^N} \int_{\R^N} \frac{|u(x)-u(y)|^2}{|x-y|^{N+2s}} \left(\frac{\nabla \varphi_n(x) \cdot x + \nabla \varphi_n(y) \cdot y}{2} \right) dx dy \\
&\to& \, %C'_{N,s}
\frac{C_{N,s}}{2} \int_{\R^N} \int_{\R^N} \frac{|u(x)-u(y)|^2}{|x-y|^{N+2s}} \left( N- \frac{N+2s}{2}\right) % \left( \frac{2N}{N+2s} -1\right) dx dy
= \frac{N-2s}{2} %C_{N,s} 
[u]_{\R^N}^2
\end{eqnarray*}
where we used $\varphi_n \to 1$, $\nabla \varphi_n \to 0$ and dominated convergence theorem.
Similarly
\begin{align*} 
(\mc{R}_{F(u)}^{\alpha}, \mc{K}_{X_n}^{\alpha})_{L^2(\R^{2N})} &\to %\, C'_{N,\alpha}
\int_{\R^N} \int_{\R^N} I_{\alpha}(x-y) F(u(x)) F(u(y)) \left( N - \frac{N-\alpha}{2} \right) dx dy \\ %\left( \frac{2N}{N-\alpha} -1\right) dx dy \\
&= \frac{N+\alpha}{2} \int_{\R^N} \big(I_{\alpha}*F(u)\big) F(u).
\end{align*}
and
$$\frac{\mu}{2} \big(u^2,\dive(X_n)\big)_{L^2(\R^N)} \to \mu \frac{N}{2} \norm{u}_2^2.$$
Joining the pieces, we have the claim. %\tr{CHECK the constants $\frac{1}{4}$}.
\QED

\bigskip

\medskip

\claim Proof of Theorem \ref{thm_pohoz_frac_choq}.
The theorem is a consequence of Corollary \ref{corol_reg_pohz} and Theorem \ref{thm_pohoz_X}.
\QED

\begin{Remark}
We comment the name of $\mc{K}_{X}^s$. Indeed, up to a multiplicative constant, we have, for any $\beta \in (0,1)$ and $X \in Lip_c(\R^N, \R^N)$, by \cite[equations (2.9c) and (2.11)]{CoSt} (see also \cite{{Silh0}})
\begin{align*}
\int_{\R^N} \int_{\R^N} \frac{\mc{K}_{X}^s(x,y)}{|x-y|^{N+\beta-1}} =& \int_{\R^N} \left( %\frac{2}{N+2s}
 \int_{\R^N} \frac{\dive(X(y))}{|x-y|^{N+\beta-1}}dy\right)dx- \\
&- \frac{N+2s}{2} \int_{\R^N} \left( \int_{\R^N} \frac{(X(x)-X(y))\cdot(x-y)}{|x-y|^{N+\beta+1}}dy \right)dx \\
=& (N+\beta-1) %\frac{2(N+\beta-1)}{N+2s} 
\int_{\R^N} \dive^{\beta}(X)(x) %dx 
- \frac{N+2s}{2} \int_{\R^N} \dive^{\beta}(X)(x) %dx
\\
=& (N+2\beta-2-2s) % \frac{N+2\beta-2-2s}{N+2s}
 \int_{\R^N} \dive^{\beta}(X) ; 
\end{align*}
in particular
$$\int_{\R^N} \int_{\R^N} \frac{\mc{K}_{X}^s(x,y)}{|x-y|^{N+s-1}} = (N-2) %\frac{N-2}{N+2s} 
\int_{\R^N} \dive^{s}(X). $$
We refer also to \cite[Chapter 3]{Dje0} where $\mc{K}_{X}^s$ is seen as the derivative of a suitable family of \emph{deformations}.

\end{Remark}

%%%%%%%%%%%%%%%%%%%%%%%%%%%%%%%%%%%%%%%%%%%%%%%%%%%%%%%

%\phantomsection
%\addcontentsline{toc}{section}{Some open problems}
%
%\section*{Some open problems}
%
%\begin{itemize}
%
%\item Can we extend some regularity results %about the doubly nonlocal problem \eqref{..} even 
%to sign-changing solutions? In this case, when $f$ is sublinear, we expect some regularity only out of the set $\{u=0\}$.
%
%\item Are there other interesting phenomena related to the $s$-sublinear threshold $p^*_{\alpha,s}$?
%% appearing in the study of the asymptotic behaviour of the ground state of the doubly nonlocal equation \eqref{..}?
%
%\item Can we extend some asymptotic decay results to doubly nonlocal equations in non-variational form? In particular, where a term of the type $(I_{\alpha}*u^p)u^q$ appears, with $p \neq q+1$. In this case it could be interesting to see if the threshold related to the change of decay is related only to $q$, or if $p$ has also a role as well.
%
%\item Can we implement a more direct proof of the Pohozaev identity, by exploiting the fractional divergence and gradient operators?
%
%\end{itemize}

%%%%%%%%%%%%%%%%%%%%%%%%%%%%%%%%%%%%%%%%%%%%%%%%%%%%%%%
%%%%%%%%%%%%%%%%%%%%%%%%%%%%%%%%%%%%%%%%%%%%%%%%%%%%%%%
%%%%%%%%%%%%%%%%%%%%%%%%%%%%%%%%%%%%%%%%%%%%%%%%%%%%%%%

\chapter{Concentration phenomena: the effect of the fractional operator}
\label{chap_concentr}

In this Chapter we investigate how the nonlocalities interact with concentration phenomena. 
%
%Firstly 
We consider the fractional, semiclassical nonlinear Schr\"odinger equation
$$ \eps^{2s}(- \Delta)^s v+ V(x) v= f(v), \quad x \in \R^N$$
where $s \in (0,1)$, $N \geq 2$, $V \in C(\R^N,\R)$ is a positive potential and $f$ is a nonlinearity satisfying Berestycki-Lions type conditions. 
For $\eps>0$ small, we prove the existence of at least $\cupl(K)+1$ positive solutions, where $K$ is a set of local minima in a bounded potential well and $\cupl(K)$ denotes the cup-length of $K$. 
Due to the generality of $f$, we cannot implement a Lyapunov-Schimdt reduction, nor we can bound our functional on a Nehari manifold: thus, by means of variational methods, our approach is to analyze the topological difference between two levels of an indefinite functional in a neighborhood of expected solutions.
Since the nonlocality comes in the decomposition of the space directly, we introduce also a new \emph{fractional center of mass}, via a suitable seminorm.
Some other delicate aspects arise strictly related to the presence of the nonlocal operator: in particular, $L^{\infty}$-boundedness, regularity and polynomial decay have to be specifically investigated.
% . By using regularity results based on fractional De Giorgi classes, 
We show then that the found solutions decay polynomially and concentrate around some point of $K$ as $\eps \to 0$.

\medskip

The main discussion (Section \ref{sec_conc_intro}--\ref{sec_cup-length}) will be focused on the case $f$ Sobolev-subcritical. This is based mainly on paper \cite{CG0}.
Afterwards, in Section \ref{sec_conc_critical}, we will see how to treat the case $f$ critical; the argument will be based mainly on paper \cite{Gal0}.

%\medskip
%
%The study of concentration on points for the Choquard equation %($s=1$, see equation \eqref{eq_concent_choquard}) 
%%$$ -\eps^{2} \Delta v+ V(x) v= \frac{1}{\eps^{\alpha}}\big(I_{\alpha}*F(v)\big) f(v), \quad x \in \R^N$$
%has been covered, in the full generality, by Cingolani and Tanaka in \cite{CT0} (see also \cite{MS-23, BCS, CCS2}) by means of some suitable tail minimization function. 
%Anyway, in the last Section about open problems, % \ref{sec_conc_cerchi_sfere} 
%we will give a sketch of an ongoing work, in collaboration with Denis Bonheure and Silvia Cingolani \cite{BoCiGa0, BG0}, about a Schr\"odinger-Newton system, where the nonlocality of the nonlinearity plays an heavy role in the concentration on spheres and disks.
	
%%%%%%%%%%%%%%%%%%%%%%%%%%%%%%%%%%%%%%%%%%%%%%%%%

\section{From classical to quantum: semiclassical states 
%\tb{spostare qualcosa all'inizio}%COMMENT NOW
}
%\label{section:1}
\label{sec_conc_intro}

In Section \ref{sec_introd_frac_lap} we highlighted the physical relevance of the fractional Laplacian operator. In particular we mentioned the study of standing waves of
%Special solutions of the equation
 the fractional nonlinear Schr\"odinger (fNLS for short) equation 
	\begin{equation}\label{eq_fNLS_2}
	i \hbar \partial_t \psi = \hbar^{2s} (-\Delta)^s \psi + V(x) \psi- f(\psi), \quad (t,x) \in (0,+\infty) \times \R^N
	\end{equation}
% \eqref{fNLS} 
%are given by the standing waves, 
i.e. factorized functions $\psi (t,x) = e^{\frac{i \mu t}{\hbar}} v(x)$, $\mu\in \R$. Instead of considering the fixed case $\hbar=1$, we focus now on the study of 
%For 
small $\hbar >0$: in this case % these 
standing waves are usually called %referred to as 
\emph{semiclassical} states and the transition from quantum physics to classical physics is somehow described letting $\hbar \to 0$.

Roughly speaking, when $s = 1$ the centers of mass $q_{\eps}=q_{\eps}(t)$ of the soliton solutions in \eqref{eq_fNLS_2}, under suitable assumptions and initial conditions, converge as $\hbar\to 0$ to the solution of the Newton's equation of motion
\begin{equation}\label{eq_newton}
\ddot{q}(t)=-\nabla V(q(t)), \quad t \in (0, +\infty);
\end{equation}
for $s\in (0,1)$ a suitable power-type modification of equation \eqref{eq_newton} is needed. Here, considering small $\hbar$ roughly means that the size of the support of the soliton in \eqref{eq_fNLS_2} is considerably smaller than the size of the potential $V$; 
for details we refer to \cite{BrJe,FGJS,JFG,BGM2}, and to \cite{SeSq} for the fractional case (see also \cite{BDGS} for the Choquard case).

Similar problems %This kind of equations 
arise also in the study of superconductivity in Ginzburg-Landau vortices, see \cite{BBH} and references therein; here the point of concentration is indeed a point where a vortex is formed.

\medskip

%Without loss of generality, shifting $E$ to
%$0$ and denoting $\hbar \equiv \varepsilon$, the search of semiclassical states leads to investigate the following nonlocal equation
%\begin{equation} \label{eq_concent_genericaf}
%%\label{eq:1.1}
%\varepsilon^{2s} (-\Delta)^s v + V(x)v = f(v), \quad x \in \R^N 
%\end{equation}
%when $V$ is positive and $\varepsilon > 0$ is small.

%[RIPETIZIONE]
Without loss of generality, shifting $\mu$ to $0$ and denoting $\hbar \equiv \eps$, the search for semiclassical states leads to the investigation of the following nonlocal equation
\begin{equation} \label{eq_concent_genericaf}
%\label{eq_principale}
\eps^{2s} (-\Delta)^s v + V(x) v = f(v), \quad x \in \R^N
\end{equation}
where $V$ is positive and $\eps > 0$ is small. 
Setting $u:= v(\eps \cdot)$, we observe that \eqref{eq_concent_genericaf} can be rewritten as
\begin{equation}\label{eq_cambio_var}
(-\Delta)^s u + V(\eps x) u= f(u), \quad x \in \R^N,
\end{equation}
thus the equation
\begin{equation}\label{eq_limite_iniz} % \label{eq_intr_limit_eq}
(-\Delta)^s U + a U = f(U), \quad x \in \R^N
\end{equation}
becomes a formal limiting equation, as $\eps \to 0$, of \eqref{eq_cambio_var}, for some $a>0$. Indeed, if $x_0 \in \R^N$ and $r>0$, 
$$\sup_{x \in B(x_0, \eps r)} \pabs{V(\eps x)- V(x_0)} \to 0 \quad \hbox{as $\eps \to 0$}.$$
Solutions of \eqref{eq_concent_genericaf} usually exhibit concentration behaviour as $\eps\to 0$: by \emph{concentrating solutions} we mean a family $v_{\eps}$ of solutions of \eqref{eq_concent_genericaf} which converges, up to rescaling, to a ground state of \eqref{eq_limite_iniz} and whose maximum points converge to some point $x_0 \in \R^N$ given by the topology of $V$ (see Theorem \ref{teo_main} for a precise statement). 
This point $x_0$ reveals, generally, to be a critical point of $V$ -- i.e. an equilibrium of \eqref{eq_newton} -- as shown in \cite{Wan0, FMV}.

\medskip

In the limiting case $s=1$ the semiclassical analysis of NLS equations has been largely investigated, starting from the seminal paper \cite{FlWe0}: 
by means of a finite Lyapunov-Schmidt dimensional reduction argument, Floer and Weinstein proved the existence of positive spike solutions to the homogeneous 3D cubic NLS equation, concentrating at each nondegenerate critical point of the potential $V$ (see also \cite{Oh2}); here the nondegeneracy of the ground states of the limiting problem \eqref{eq_limite_iniz} is crucial. 
%here the authors implement a Lyapunov-Schmidt dimensional reduction argument to gain existence of solutions for homogeneous sources, relying on the nondegeneracy of the ground states of the limiting problem \eqref{eq_limite_iniz}. 
Successively, refined variational techniques were implemented to study singularly perturbed elliptic problems in entire space: several existence results of positive spike solutions to the NLS equation in a semiclassical regime are derived under different assumptions on the potential and the nonlinear terms. 
 We confine to mention \cite{Rab1, Wan0, DF0, ABC,DF2, BoVS0, BJ0, BT1, BT2, DPR} and references therein.

Starting from the work \cite{DaWe0}, in \cite{CiLa1,CiLa2,AmMaSe,CJT,DLY0, Jin0, Che0} topological invariants were used to derive multiplicity results in singularly perturbed frameworks, in the spirit of well known results of Bahri, Coron \cite{BaCo} and Benci, Cerami \cite{BeCe1} for semilinear elliptic problems with Dirichlet boundary condition. Precisely, in \cite{CiLa1} it has been proved that the number of positive solutions of the stationary NLS equation is influenced by the topological richness of the set of global minima of $V$. 
Some years later, using a perturbative approach, Ambrosetti, Malchiodi and Secchi \cite{AmMaSe} obtained a multiplicity result for the NLS equation with power nonlinearity, assuming that the set of critical points of $V$ is nondegenerate in the sense of Bott.
More recently, in \cite{CJT} Cingolani, Jeanjean and Tanaka improved the result in \cite{CiLa1}, relating the number of semiclassical standing waves solutions to the \emph{cup-length} of $K$, where $K$ is a set of local (possibly degenerate) minima of the potential, under almost optimal assumptions on the nonlinearity (see also the recent paper \cite{CT0} in the context of nonlinear Choquard equations).

\smallskip

When $s \in (0,1)$, the search of semiclassical standing waves for the fNLS equation has been firstly considered by D\'avila, Del Pino and Wei in \cite{DPW} under the assumptions $f(t)= |t|^{p-2}t$, with $2 < p < 2^*_s$, where $2^*_s: = \frac{2N}{N-2s}$ is the Sobolev critical exponent, and $V \in C^{1,\alpha}(\R^N)$ is bounded. 
Using a Lyapunov-Schmidt reduction inspired by \cite{FlWe0, Oh2}, they showed the existence of a positive spike solution whose maximum point concentrates at some nondegenerate critical point of $V$: this approach relies on the nondegeneracy property of the linearization at the positive ground state shown by Frank, Lenzmann and Silvestre \cite{FLS}.
Successively, inspired by \cite{DF0, BJ0}, variational techniques were employed to derive existence of spike solutions concentrating at local minima of $V$, 
see \cite{Seo0, AlMi1, Amb0} and references therein (see also \cite{Sec0} where global assumptions on $V$ are considered). 

A first multiplicity result for the (fNLS) equation is obtained 
in \cite{FiSi0}, inspired by \cite{CiLa1}. Precisely, letting $K$ be the set of global minima of $V$, Figueiredo and Siciliano proved that the number of positive solutions of \eqref{eq_concent_genericaf}, when $f$ satisfies monotonicity and Ambrosetti-Rabinowitz condition, is at least given by the Ljusternik-Schnirelmann category of $K$: here the search of solutions of \eqref{eq_concent_genericaf} can be reduced to the study of the (global) level sets of the Nehari manifold, where the energy functional is restricted, and to deformation arguments valid on Hilbert manifolds without boundary. 
See also \cite{AA0} where the Ambrosetti-Rabinowitz condition is dropped. 
In \cite{Che0}, moreover, Chen implemented a Lyapunov-Schmidt reduction for nondegenerate critical points of $V$ and power-type functions $f$ in order to get multiplicity results related to the cup-length, extending the results of \cite{AmMaSe}.

\medskip

In this first part of the Chapter we are interested to prove multiplicity of positive solutions for the fNLS equation \eqref{eq_concent_genericaf} when $\eps$ is small, without monotonicity and Ambrosetti-Rabinowitz conditions on $f$, nor nondegeneracy and global conditions on $V$, concentrating at a local minimum of $V$.

On the potential $V$ we assume
\begin{itemize}
	\item[(V1)] \label{(V1)}
$\, V\in C(\R^N)\cap L^{\infty}(\R^N)$, $\underline{V}:=\inf_{\R^N} V>0$ (see also Remark \ref{rem_ipotesi_V});
	\item[(V2)] \label{(V2)}
$\,$there exists a bounded domain $\Omega\subset \R^N$ such that
	$$m_0:= \inf_{\Omega} V < \inf_{\partial \Omega} V;$$
\end{itemize}
\vspace{-\topsep}
by the strict inequality and the continuity of $V$, we can assume that $\partial \Omega$ is regular.
We define $K$ as the set of local minima
\begin{equation}\label{eq_def_K}
K:=\{ x \in \Omega \mid V(x)=m_0\}.
\end{equation}

On $f$ we assume the following \emph{subritical} assumptions 
\begin{itemize}
	\item[(f1)] \label{(f1)}
Berestycki-Lions type assumptions with respect to $m_0$, that is
	\begin{itemize}
		\item[(f1.1)] $\, f\in C(\R, \R)$;
		\item[(f1.2)] $\, \lim_{t \to 0} \frac{f(t)}{t}=0$;
		\item[(f1.3)] $\, \lim_{t \to +\infty} \frac{f(t)}{|t|^p}=0$ for some $p \in (1, 2^*_s-1)$, where we recall $2^*_s=\frac{2N}{N-2s}$; 
		\item[(f1.4)] $\, F(t_0)> \frac{1}{2} m_0 t_0^2$ for some $t_0>0$, where $F(t):= \int_0^t f(s) ds$;
	\end{itemize}
	\item[(f2)] \label{(f2)}
$f(t)=0$ for $t\leq 0$.
\end{itemize}

On $f$ we further assume 
\begin{itemize}
	\item[(f3)] \label{(f3)}
$f\in C^{0,\gamma}_{loc}(\R)$ for some $\gamma \in (1-2s, 1)$ if $s\in (0,1/2]$.
\end{itemize}

\begin{Remark}
We remark that \hyperref[(f3)]{\textnormal{(f3)}} is needed only to get a Pohozaev identity (see Proposition \ref{Pohozaev-prop}). 
See also Remark \ref{rem_remove_Pohozaev} below.
\end{Remark}

%\smallskip

%Setting $u:= v(\eps \cdot)$, \eqref{eq_concent_genericaf} can be rewritten as
%\begin{equation}\label{eq_cambio_var}
%(-\Delta)^s u + V(\eps x) u= f(u), \quad \textnormal{ $x \in \R^N$},
%\end{equation}
%thus the equation
%\begin{equation}\label{eq_intr_limit_eq}
%(-\Delta)^s U + a U = f(U), \quad x \in \R^N,
%\end{equation}
%for some $a>0$, becomes a formal limiting equation, as $\eps \to 0$, for \eqref{eq_cambio_var}.

It is standard that weak solutions to \eqref{eq_cambio_var} correspond to critical points of the $C^1$-energy functional
$$I_{\eps}(u):=\frac{1}{2} \int_{\R^N} \abs{(-\Delta)^{s/2} u}^2 dx+ \frac{1}{2} \int_{\R^N} V(\eps x) u^2 dx- \int_{\R^N}F(u) dx, \quad u\in H^s(\R^N). $$
%where $H^s(\R^N)$ is the fractional Sobolev space. 
We remark that, because of the general assumptions on $f$, we can not take advantage of the boundedness of the functional from above and below, nor of Nehari type constraint. 
Therefore in the present paper we combine reduction methods and \emph{penalization arguments} in a nonlocal setting: 
in particular, as in \cite{CJT, CT0}, the analysis of the topological changes between two level sets of the indefinite energy functional $I_{\eps}$ in a small \emph{neighborhood $\mc{X}_{\eps, \delta}$
of expected solutions} is essential in our approach. 
With the aid of $\eps$-independent pseudo-differential estimates, we detect 
such a neighborhood, which will be positively invariant
under a pseudo-gradient flow, and we develop our deformation argument in the context of nonlocal operators. 
To this aim we introduce two maps $\Phi_\varepsilon$ and $\Psi_\varepsilon$ between topological pairs: we emphasize that to define such maps, a center of mass $\Upsilon$ and a functional $P_a$ which is inspired by the Pohozaev identity are crucial.

With respect to the local case, several difficulties arise linked to special features of the nonlocal nature of the problem: among them we have 
the polynomial decay of the least energy solutions of the limiting problems, the weak regularizing effect of the fractional Laplacian, the lack of general comparison arguments, the differences between the supports of a function and of its Fourier transform, and the lack of the standard Leibniz formula (see e.g. \cite{Sec0, BWZ}). 
Moreover we highlight that, for fractional equations, the nonlocal part strongly influences the decomposition of the space and this makes quite delicate to use truncating test functions and perform the localization of the centers of mass. 

In the present Chapter we need to implement new ideas to overcome the above obstructions; in particular we introduce a new \emph{fractional local center of mass} by means of a suitable seminorm, stronger than the usual Gagliardo seminorm in a bounded set.

\vspace{2mm}
Our main result is the following theorem.
\begin{Theorem}
\label{teo_concen_esist}
%\label{claim:1.1}
	Suppose $N\geq 2$ and that \hyperref[(V1)]{\textnormal{(V1)}}-\hyperref[(V2)]{\textnormal{(V2)}}, \hyperref[(f1)]{\textnormal{(f1)}}--\hyperref[(f3)]{\textnormal{(f3)}} hold. 		
	Let $K$ be defined by \eqref{eq_def_K}.
	Then, for sufficiently small $\eps >0$, equation \eqref{eq_concent_genericaf} has at least
	$\cupl(K)+1$ positive solutions, which belong to $C^{0, \sigma}(\R^N) \cap L^{\infty}(\R^N)$ for some $\sigma \in (0,1)$.
\end{Theorem}

Here $\cupl(K)$ denotes the \emph{cup-length} of $K$ defined by the Alexander-Spanier cohomology with coefficients in some field $\F$ (see Appendix \ref{chap_app_alg_top}). %Section \ref{sec_cup-length}). 
Notice that the cup-length of a set $K$ is strictly related to the \emph{category} of $K$, see Lemma \ref{lem_coll_cat_cupl} and Remark \ref{examples}.

\begin{Remark}\label{rem_ipotesi_V}
Observe that, arguing as in \cite{BJ0} and \cite{BT1}, we could omit the assumption that $V$ is bounded from above in Theorem \ref{teo_concen_esist}. For the sake of simplicity, we assume here the boundedness of $V$.
\end{Remark}

The regularity statement in Theorem \ref{teo_concen_esist} relies on some recent regularity results based on fractional De Giorgi classes and tail functions (see Section \ref{sec_regol_degiorgi}); notice that the fact that the noncriticality is \emph{strict} in \hyperref[(f1)]{\textnormal{(f1.3)}} (that is $p<2^*_s-1$) is here needed.
Through these results we are able to prove also the concentration of the solutions. % theorem. 

%\begin{Theorem}
%%\label{mean-concen}
%\label{teo_concen_conc}
%	Let $(\eps_n)_{n }$ with $\eps_n \to 0^+$ as $n \to +\infty$. For sufficiently large $n\in \N$, in the assumptions of Theorem \ref{teo_concen_esist}, let $v_{\eps_n}$ be one of the $\cupl(K)+1$ solutions of equation \eqref{eq_concent_genericaf}.
%	Then, up to a subsequence, $(v_{\eps_n})_{n\in \N}$ \emph{concentrates} in $K$ as $n \to +\infty$. Namely, for each $n\in \N$ there exists a maximum point $x_{\eps_n}\in \R^N$ of $v_{\eps_n}$ such that 
%$$\lim_{n\to +\infty}d(x_{\eps_n}, K) =0;$$
%\tor{more precisely, there exists a point $x_0 \in K$ sucht that $x_{\eps_n} \to x_0$.}
%	In addition, $ v_{\eps_n}(\eps_n \cdot+x_{\eps_n})$ converges in $H^s(\R^N)$ and uniformly on compact sets to a least energy solution of
%	\begin{equation}\label{eq_least_energy_m0}
%	(- \Delta)^s U + m_0 U = f(U), \quad U >0, \quad U \in H^{s}(\R^N),
%	\end{equation}
%	and, for some positive $C', C''$ independent on $n \in \N$, we have the uniform polynomial decay
%	$$\frac{C'}{1+|\frac{x-x_{\eps_n}}{\eps_n}|^{N+2s}}\leq v_{\eps_n}(x) \leq \frac{C''}{1+|\frac{x-x_{\eps_n}}{\eps_n}|^{N+2s}},\quad \textit{ for $x \in \R^N$}.$$
%\end{Theorem}

\begin{Theorem}
%\label{mean-concen}
\label{teo_concen_conc}
	 %For sufficiently large $n\in \N$, 
In the assumptions of Theorem \ref{teo_concen_esist}, 
let $v_{\eps}$ be one of the $\cupl(K)+1$ family of solutions of equation \eqref{eq_concent_genericaf}. Then, $(v_{\eps})_{\eps>0}$ \emph{concentrates} in $K$ as $\eps \to 0$, i.e. there exist a maximum points $x_{\eps}\in \R^N$ of $v_{\eps}$ such that 
$$\lim_{\eps \to 0}d(x_{\eps}, K) =0;$$
moreover, for some positive $C', C''$, we have the uniform polynomial decay
	$$\frac{C'}{1+|\frac{x-x_{\eps}}{\eps}|^{N+2s}}\leq v_{\eps}(x) \leq \frac{C''}{1+|\frac{x-x_{\eps}}{\eps}|^{N+2s}},\quad \textit{ for $x \in \R^N$}.$$
	% Namely, for each $n\in \N$ 
%\tor{more precisely, .}
In addition, let $(\eps_n)_{n }$ with $\eps_n \to 0^+$ as $n \to +\infty$. Then, up to a subsequence,
there exists a point $x_0 \in K$ sucht that $x_{\eps_n} \to x_0$ as $n \to +\infty$, and $ v_{\eps_n}(\eps_n \cdot+x_{\eps_n})$ converges in $H^s(\R^N)$ and uniformly on compact sets to a least energy solution of
	\begin{equation}\label{eq_least_energy_m0}
	(- \Delta)^s U + m_0 U = f(U), \quad U >0, \quad U \in H^{s}(\R^N).
	\end{equation}
\end{Theorem}

\smallskip

This first part of the Chapter is organized as follows. 
%\tr{aggiusta label sezioni} %COMMENT NOW
In Section \ref{sec_fract_fram} we recall the mixed Gagliardo seminorm introduced in Section \ref{sec_prel_sobolev}, 
% some basic notions on the nonlocal framework. 
while in Section \ref{sez_limit_eq} we show the uniform polynomial decay of the solutions of \eqref{eq_limite_iniz}, and we introduce a new fractional center of mass $\Upsilon$, by means of a suitable seminorm. 
% we postpone to Appendix \ref{sez_poly_dec} the proof of the uniform polynomial decay of the solutions of \eqref{eq_limite_iniz}, since it shares some argument with the proof of Theorem \ref{teo_concen_conc}. 
Section \ref{sez_singur_perturb} is the main core of the Chapter, where we introduce a penalized functional and prove a deformation lemma on a neighborhood of expected solutions; moreover, we build suitable maps $\Phi_{\eps}$, $\Psi_{\eps}$ essential in the proof of the multiplicity of solutions. Then in Section \ref{sec_cup-length} we prove Theorem \ref{teo_concen_esist} by the use of the deformation lemma and the built maps applied to the theory of relative category and relative cup-length. %, of which we briefly recall the definitions. 
Finally we prove Theorem \ref{teo_concen_conc} by using regularity results based on fractional De Giorgi classes.

\medskip

Afterwards, in Section \ref{sec_conc_critical} we move to the study of the critical case.

%\tr{Introduci corsivo per i funzionali? Cambia $f$ con $g$?} %COMMENT NOW

\subsection{A tail-controlling mixed norm}
\label{sec_fract_fram}

In this Chapter %, together with the notation of Section \ref{...}, 
we will make use of the following norm
$$\norm{u}_{H_{\eps}^s(\R^N)}^2:= \norm{(-\Delta)^{s/2} u}_2^2 + \int_{\R^N} V(\eps x) u^2 dx$$
which is an equivalent norm on $H^s(\R^N)$ (once $\eps$ is fixed), thanks to the positivity and the boundedness of $V$; % and \eqref{eq_equiv_norm_laplac}; 
the space $H^s_{\eps}(\R^N)$ is defined straightforwardly.
Moreover we will make use of the \emph{mixed} Gagliardo seminorm (introduced in Section \ref{sec_mixed_norm})
$$[u]_{A_1, A_2}^2= \int_{A_1} \int_{A_2} \frac{|u(x)-u(y)|^2}{|x-y|^{N+2s}} dx \, dy, \quad [u]_A:=[u]_{A,A}$$
for any $A_1, A_2,A \subset \R^N$ and $u \in H^s(\R^N)$.
%by using that $\varphi_u(x,y):= \frac{|u(x)-u(y)|}{|x-y|^{N/2+s}}$ satisfies $\varphi_{u+v} \leq \varphi_u + \varphi_v$ and $[u]_{A_1, A_2} = \norm{\varphi_u}_{L^2(A_1 \times A_2)}$, we have that $[u]_{A_1, A_2}$ is actually a seminorm. 
For any $u \in H^s(\R^N)$ and $A \subset \R^N$ it will be useful to work also with the following norms:
\begin{equation}\label{eq_def_Htilde}
\norm{u}_{A}^2:= \norm{u}_{L^2(A)}^2 + [u]_{A,\R^N}^2
\end{equation}
and
\begin{equation}\label{eq_def_triplanorm}
\tnorm{u}_{A}:=\norm{u}_{A} + \norm{u}_{L^{p+1}(A)},
\end{equation}
where $p$ is introduced in assumption \hyperref[(f1)]{\textnormal{(f1.3)}}. We highlight that $\norm{u}_{\R^N}=\norm{u}_{H^s(\R^N)}$, but generally $\norm{u}_A \geq \norm{u}_{H^s(A)}$ for $A \neq \R^N$. 
By $H^s(A) \hookrightarrow L^{p+1}(A)$ the norms $\norm{\cdot}_A$ and $\tnorm{\cdot}_A$ are equivalent: on the other hand, the constant such that $\tnorm{u}_A \leq C_A \norm{u}_A$ depends on $A$, thus not useful for expanding sets % $\eps$-dependent sets 
$A=A(\eps)$. This is why we will make direct use of $\tnorm{\cdot}_A$. 

%Moreover, we denote by 
%$$[u]_{C^{0,\sigma}(A)}:= \sup_{\substack{x, y \in A \\ x\neq y}} \frac{\abs{u(x)-u(y)}}{\abs{x-y}^{\sigma}}$$
%the usual seminorm in H\"older spaces for $\sigma \in (0,1]$. We allow $\sigma>1$ by simply writing $C^{\sigma}(\R^N)$.

%[TOGLIERE] %CCOMMENT NOW
Before ending this Section, we highlight that, by the assumptions on $f$,
%the assumptions on the function $f$ imply the following standard fact, which will be extensively used throughout the paper: 
for each $q\geq p$ and $\beta>0$ there exists a $C_{\beta}>0$ such that
\begin{equation}\label{eq_prop_f}
\abs{f(t)}\leq \beta |t| + C_{\beta} |t|^q \quad \hbox{ and } \quad \abs{F(t)}\leq C\left(\beta |t|^2 + C_{\beta} |t|^{q+1}\right).
\end{equation}

%\medskip
%
%\textbf{Further notation}. We highlight that, all throughout the paper, we will assume $N\geq 2$, and the constants $C, C'$ appearing in inequalities may change from a passage to another. 
%To avoid cumbersome notations, we will not stress the dependence of such constants, which will be based only on the fixed quantities. 
%Moreover, by $\approx$ we will mean \emph{approximately equal to}.
%
%We will write $B_R(x_0)$ for the ball, centered in $x_0 \in \R^N$, with radius $R>0$; in particular, $B_R:=B_R(0)$. Moreover 
%$$A_{\delta}:= \{x \in X \mid d(x,A)\leq \delta\}$$
%for any $A\subset (X,d)$ metric space. In addition, we will sometimes write $\complement(A):= \R^N \setminus A$ for $A \subset \R^N$ to avoid cumbersome notation. 
%
%Finally, for every $A\subset B \subset \R^N$, we will write
%$$A\prec \phi \prec B$$
%to indicate a Urysohn-type regular function $\phi \in C^{\infty}_c(\R^N)$ such that
%$$\phi_{|A}=1 \quad \textnormal{ and } \quad \phi_{|\R^N \setminus B}=0.$$

\section{Limiting equation}
\label{sez_limit_eq}

In this Section we further investigate the autonomous equation studied in Section \ref{sec_frac_unconstr}.

\subsection{A single equation} \label{subsec_single_eq}

Consider 
\begin{equation}\label{eq_limit_eq}
(-\Delta)^s U + a U = f(U), \quad x \in \R^N
\end{equation}
with $a>0$. Weak solutions of \eqref{eq_limit_eq} are known to be characterized as critical points of the $C^1$-functional $L_a: H^s(\R^N) \to \R$
$$L_a(U):= \frac{1}{2} \norm{(-\Delta)^{s/2}U}_2^2 + \frac{a}{2} \norm{U}_2^2 - \int_{\R^N} F(U) dx, \quad U \in H^s(\R^N).$$
Set moreover the Pohozaev functional $P_a: H^s(\R^N) \setminus \{0\} \to \R_+$
$$P_a(U):= \left( 2^*_s % \frac{2N}{N-2s} 
\; \frac{\int_{\R^N} F(U) dx - \frac{a}{2} \norm{U}_2^2}{\norm{(-\Delta)^{s/2} U}_2^2}\right)_+^{\frac{1}{2s}}, \quad U\in H^s(\R^N), \; U \neq 0.$$
We further set
$$C_{po,a}:=\inf \big\{ L_a(U) \mid U \in H^s(\R^N) \setminus \{0\}, \; P_a(U)=1\big\}$$
the Pohozaev minimum energy, 
and
$$E_a:=\inf \big\{L_a(U) \mid U \in H^s(\R^N) \setminus \{0\}, \; L'_a(U)=0 \big\}$$
the least energy for $L_a$.

We recall the following result by Section \ref{sec_frac_unconstr}, where we further highlight the regularity, the positivity and the decay at infinity (see \cite[Theorems 1.1--1.3]{BKS} and \cite[Theorem 1.5]{FQT}).

\begin{Theorem}%[Existence of a Pohozaev minimum] %[Regularity and polynomial decay]
\label{thm_esist_Poh_min}
\label{thm_reg_dec}
Assume \hyperref[(f1)]{\textnormal{(f1)}} with respect to $a>0$ and \hyperref[(f2)]{\textnormal{(f2)}}. 
\begin{itemize}
\item There exists a positive minimizer for $C_{po,a}>0$, which is a weak solution of \eqref{eq_limit_eq}.
%; this minimizer is generally not unique, even up to translation.
\item Every weak solution $U \in H^s(\R^N)$ of \eqref{eq_limit_eq} is actually a strong solution, i.e. $U$ satisfies \eqref{eq_limit_eq} almost everywhere. Moreover
$U\in H^{2s}(\R^N)\cap C^{\sigma}(\R^N)$
for every $\sigma \in (0,2s)$.
\item Every weak solution $U \in H^s(\R^N)$ of \eqref{eq_limit_eq} is strictly positive 
%$$U>0 \quad \textnormal{and} \quad \lim_{|x|\to +\infty} U(x)=0.$$
and decays polynomially at infinity, that is there exist positive constants $C_a', C_a''$ such that
\begin{equation}\label{eq_bound_funct}
\frac{ C_a'}{1+|x|^{N+2s}} \leq U(x)\leq \frac{ C_a''}{1+|x|^{N+2s}},\quad \textit{ for $x \in \R^N$}.
\end{equation}
Observe that the bounding functions in \eqref{eq_bound_funct} belong to $L^q(\R^N)$ for any $q\in [1,+\infty]$. 
\item If \hyperref[(f3)]{\textnormal{(f3)}} holds, then the Pohozaev identity
$$P_a(U)=1$$
holds for each nontrivial solution $U$ of \eqref{eq_limit_eq}. As a consequence
\begin{equation}\label{eq_ug_poh_gr.st}
E_a = C_{po,a}.
\end{equation}
\end{itemize}
\end{Theorem}

We observe that to reach the Pohozaev identity we need the solutions to be regular enough, fact that is given by \hyperref[(f3)]{\textnormal{(f3)}}. 
The functional $P_a$ will be of key importance for estimating $L_a$ from below, see Lemma \ref{lemma_stima_g0} and Lemma \ref{lemma_stima_g}. 

We highlight the polynomial decay of solutions of \eqref{eq_limit_eq}. This decay is much less slower than the one, exponential, of the local case $s=1$. An alternative proof, which underlines some uniformity in $a$, can be found in Proposition \ref{prop_polyn_dec}.

\begin{Remark} 
\label{rem_decay_of_frac}
Since the equation is satisfied almost everywhere, we have also, by \eqref{eq_prop_f}
\begin{align*}
\abs{(-\Delta)^s U(x)} &\leq \abs{f(U(x))} + a|U(x)| \leq \beta |U(x)| + C_{\beta} |U(x)|^p + a |U(x)| \\
&\leq C\left (\frac{1}{1+|x|^{(N+2s)p}} + \frac{1}{1+|x|^{N+2s}} \right) \leq C\frac{1}{1+|x|^{N+2s}}
\end{align*}
for almost every $x\in \R^N$. %\tb{(A noi però serve $(-\Delta)^{s/2}U(x)$ - NON SERVE) }
\end{Remark}

\smallskip

We end this Section by a technical lemma, which allows to link the level $L_{a}(u)$ of a whatever function $u$, having a functional $P_a(u) \approx 1$, with the ground state $E_{a}$; in particular, it provides a useful lower bound for $L_a$.

\begin{Lemma}\label{lemma_stima_g0}
Let $u\in H^s(\R^N)$ and define
$$g(t):=\frac{1}{2s} \left(N t^{N-2s} - (N-2s) t^N\right), \quad t \in \R.$$
\begin{itemize}
\item[(a)] If $P_a(u)\in \Big(0, \left(\frac{N}{N-2s}\right)^{\frac{1}{2s}}\Big)$, which we highlight is a neighborhood of $1$, then
$$L_a(u)\geq g(P_a(u))E_a.$$
\item[(b)] If $u=U_a\left(\frac{\cdot-q}{t}\right)$ for some $q \in \R^N$ and $t \in \R$, with $U_a$ being a ground state of \eqref{eq_limit_eq}, then the above inequality is indeed an equality, that is
$$ L_a\left ( U_a\left(\frac{\cdot-q}{t}\right) \right) = g(t)E_a. $$
\end{itemize}%
We highlight that the function $g$ verifies
$$g(t)\leq 1 \quad \hbox{ and } \quad g(t)=1 \iff t=1.$$
\end{Lemma}

\claim Proof.
Let $\sigma:=P_a(u)$ and set $v:=u(\sigma \cdot)$. Then $P_a(v)=1$. A straightforward computation shows, by using $P_a(v)=1$ and $g(\sigma)>0$, that
$$L_a(u) = g(\sigma) L_a(v) \geq g(\sigma) C_{po,a} = g(\sigma) E_a.$$
We see that, if $u=U_a\left(\frac{\cdot-q}{t}\right)$, then by $P_a(U_a)=1$ we have $\sigma=P_a\left( U_a\left(\frac{\cdot-q}{t}\right) \right)= t$ and thus $v=U_a(\cdot-q)$, which by translation invariance is again a ground state of \eqref{eq_limit_eq}; thus $L_a(v)=C_{po,a}$. This concludes the proof.
\QED

\medskip

\subsection{A family of equations: minimal radius map}
\label{subsec_family}

In this Section we study equation \eqref{eq_limit_eq} for variable values of $a>0$. Introduce the notation
$$\Omega[a,b]:=V^{-1}([m_0+a,m_0+b])\cap \Omega$$
and similarly $\Omega(a,b)$ and mixed-brackets combinations. 
We choose now a small $\nu_0>0$ such that the minimum $m_0$ is not heavily perturbed, namely
\begin{itemize}
\item Berestycki-Lions type assumptions \hyperref[(f1)]{\textnormal{(f1)}} hold with respect to $a \in [m_0, m_0+\nu_0]$ -- i.e., in particular, $F(t_0)> \frac{1}{2} (m_0+\nu_0) t_0^2$;
\item assumption \hyperref[(V2)]{\textnormal{(V2)}} holds with respect to $m_0+\nu_0$, i.e. $m_0 + \nu_0 < \inf_{\partial \Omega} V$;
\item $\overline{\Omega[0,\nu_0]}\subset K_d \subset \Omega$ for a sufficiently small $d>0$ subsequently fixed, see Lemma \ref{lem_cupl_jK};
\item other conditions subsequently stated, see e.g. \eqref{eq_stima_l0} and Lemma \ref{lemma_stima_g}.
\end{itemize}
We observe that, by construction, for $a \in [m_0, m_0+\nu_0]$ the considerations of Section \ref{subsec_single_eq} apply. Moreover, by scaling arguments on $C_{po,a}$, we notice that 
$$a\in [m_0, m_0+\nu_0] \mapsto E_a \in (0, +\infty)$$
is strictly increasing and that, up to choosing a smaller $\nu_0$, we have%
\footnote{ % on $\R$.
%%Let $\mu<\nu$, and $v$ be a $\nu$-Pohozaev minimum (i.e. $J_{\nu}(v)=a(\nu)$). Let rescale $v$ such that it belongs to the $\mu$-Pohozaev, i.e. $u:=v(\lambda \cdot)$: straightforwards computation shows the explicit value of $\lambda$ and that (since $\mu<\nu$) $\lambda <1$. Thus $a(\mu) \leq J_{\mu}(u) < J_{\nu}(v) = a(\nu)$. The strict inequalities comes from the Pohozaev identities, the scaling and from $\lambda<1$.}
%%Let $\nu_0$ to be fixed, and $v$ be a $(m_0+\nu_0)$-Pohozaev minimum (i.e. $L_{m_0+\nu_0}(v)=C_{po,m_0+\nu_0}$). Let rescale $v$ in such a way it belongs to the $m_0$-Pohozaev set, i.e. $u:=v(\theta \cdot)$ for some explicit $\theta$: computation shows $\theta <1$. Thus $E_{m_0+\nu_0}=C_{p0,m_0} \leq L_{m_0}(u) = \theta L_{m_0+\nu_0,v) =\theta C_{po,m_0+\nu_0} = \theta E_{m_'}$. 
Let $\nu_0$ to be fixed, and $v$ be a $m_0$-Pohozaev minimum (i.e. $L_{m_0}(v)=C_{po,m_0}$). 
Let rescale $v$ in such a way it belongs to the $(m_0+\nu_0)$-Pohozaev set, i.e. $u:=v(\cdot/\theta)$ for some explicit $\theta$: 
computation shows $\theta >1$ and $\theta \to 1$ as $\nu_0 \to 0$. 
Thus %$E_{m_0+\nu_0}=
$C_{po,m_0+\nu_0} \leq L_{m_0+\nu_0}(u) = \theta^{N-2s} L_{m_0}(v) =\theta^{N-2s} C_{po,m_0}$. % = \theta^{-N+2s} E_{m_0}$. 
By choosing $\nu_0$ small we have $\theta^{N-2s} < 2$. 
%%The strict inequalities comes from the Pohozaev identities, the scaling and from $\theta<1$.
} 
$E_{m_0+\nu_0} < 2E_{m_0}$ and thus we can find an $l_0=l_0(\nu_0)\in \R$ such that
\begin{equation}\label{eq_stima_l0}
E_{m_0+\nu_0} < l_0 < 2E_{m_0}.
\end{equation}
As a final step in the proof of the main Theorem, we will make $\nu_0$ and $l_0$ moving such that $\nu_0^n \to 0$ and $l_0^n\to E_{m_0}$. 
We now define the set of \emph{almost} ground states of \eqref{eq_limit_eq}
$$S_a:= \left \{ U \in H^s(\R^N)\setminus\{0\} \mid L'_a(U)=0,\; L_a(U)\leq l_0, \; U(0)=\max_{\R^N} U\right \}\neq \emptyset.$$
We observe that we set the last condition in order to fix solutions in a point and prevent them to escape to infinity; the idea is to gain thickness and compactness (see \cite{BJ0, CSS, CJT}): notice indeed that, in the case of a proper ground state $U\in S_a$, then $U$ is radially symmetric (see also \eqref{eq_S_radial}). 
We further define
$$\widehat{S}:= \bigcup_{a\in [m_0, m_0+\nu_0]} S_a.$$
The following properties of the set $\widehat{S}$ will be of key importance in the whole paper.

\begin{Lemma} \label{lemma_stima_unif} 
The following properties hold.
\begin{itemize}
\item[(a)] There exist positive constants $C', C''$ such that, for each $U\in \widehat{S}$ we have
\begin{equation}\label{eq_stima_uniformeU}
\frac{C'}{1+|x|^{N+2s}} \leq U(x)\leq \frac{C''}{1+|x|^{N+2s}},\quad \textit{ for $x \in \R^N$}.
\end{equation}
\item[(b)] $\widehat{S}$ is compact. Since it does not contain the zero function, we have
$$r^*:= \min_{U \in \widehat{S}} \norm{U}_{H^s(\R^N)} >0;$$
the maximum is attained as well.
\item[(c)] We have
$$\lim_{R\to +\infty} \norm{U}_{\R^N \setminus B_R} =0, \quad \textit{uniformly for $U\in \widehat{S}$},$$
where the norm $\norm{\cdot}_{\R^N\setminus B_R}$ is defined in \eqref{eq_def_Htilde}. 
Moreover, if $(U_n)_n \subset \widehat{S}$, and $(\theta_n)_n \subset \R^N$ is bounded, % included in a compact set, 
then
$$\lim_{n \to +\infty} \norm{U_n(\cdot + \theta_n)}_{\R^N \setminus B_n}=0.$$
\end{itemize}
\end{Lemma}

\claim Proof.
We divide the proof in some steps.

\noindent
\textbf{Step 1.} We see that $\widehat{S}$ is bounded. Indeed, by the Pohozaev identity, we have
$$\norm{(-\Delta)^{s/2}U}_2^2 = \frac{N}{s} L_a(U) \leq \frac{N}{s} l_0.$$
By \eqref{eq_embd_homog} we have that also $\norm{U}_{2^*_s}$ is uniformly bounded. 
Since $L'_a(U)U=0$, we have by \eqref{eq_prop_f}
$$\norm{(-\Delta)^{s/2}U}_2^2 + a\norm{U}_2^2 = \int_{\R^N} f(U)U dx\leq \beta \norm{U}_2^2 + C_{\beta} \norm{U}_{2^*_s}^{2^*_s}$$
which implies, by choosing $\beta < a$, that also $\norm{U}_2^2$ is bounded. 

\noindent
\textbf{Step 2.} There exist uniform $C>0$ and $\sigma \in (0,1)$ such that
\begin{equation}\label{eq_dim_stime_uniformi}
\norm{U}_{\infty} \leq C, \quad [U]_{C^{0,\sigma}_{loc}(\R^N)} \leq C
\end{equation}
for any $U \in \widehat{S}$. We postpone %to the Appendix 
the proof of \eqref{eq_dim_stime_uniformi}, as well as the proof of the uniform pointwise estimate \eqref{eq_stima_uniformeU} (where we use \hyperref[(f3)]{\textnormal{(f3)}}), since they will carry some arguments used in the proof of Theorem \ref{teo_concen_conc} in Section \ref{sec_concentration}; see Proposition \ref{prop_polyn_dec}.

\medskip

We show now $(b)$, which is a refinement of the fact that $S_a$ itself is compact.

\noindent
\textbf{Step 3.} We observe first that $\widehat{S}$ is closed. Indeed, if $U_k \in S_{a_k} \subset \widehat{S}$ converges strongly to $U$, then up to a subsequence we have $a_k \to a\in [m_0,m_0+\nu_0]$ and, by the strong convergence, we have that the condition
$$E_{m_0} \leq L_a(U)\leq l_0 $$
holds, which in particular implies that $U\nequiv 0$.
Moreover, exploiting the weak convergence $U_k \wto U$, and the almost everywhere convergence (together with the estimate on $f$, the uniform estimate \eqref{eq_stima_uniformeU} and the dominated convergence theorem), we obtain %by classical arguments 
that for each $v\in H^s(\R^N)$
%\tor{Non posso usare Proposition \ref{prop_converg_generiche_loc} perché non sono sulle radiali. Come ragiono? Devo usare l'uniforme limitatezza in $L^{\infty}$ (o la stima uniforme polinomiale)? VEDI. Inoltre, va bene per sottocritico non stretto?}
\begin{eqnarray*}
\lefteqn{0=L'_{a_k}(U_k)v = \int_{\R^N} (-\Delta)^{s/2} U_k (-\Delta)^{s/2} v \, dx + a_k \int_{\R^N} U_k v \, dx - \int_{\R^N} f(U_k) v \, dx} \\
&\to& \int_{\R^N} (-\Delta)^{s/2} U(-\Delta)^{s/2} v \, dx+ a \int_{\R^N} U v \, dx- \int_{\R^N} f(U) v \, dx= L'_a(U)v,
\end{eqnarray*}
that is, $L'_a(U)=0$. As regards the maximum in zero, we need a pointwise convergence. In order to get it, we exploit the fact that, by \eqref{eq_dim_stime_uniformi}, $U_k$ are uniformly bounded in $L^{\infty}(\R^N)$ and in $C^{0, \sigma}_{loc}(\R^N)$ and we apply Ascoli-Arzel\`{a} theorem to get local uniform convergence. 
This shows that $U\in S_a \subset \widehat{S}$. 

\smallskip

\noindent
\textbf{Step 4.} Let now $U_k \in S_{a_k} \subset \widehat{S}$. 
By the boundedness, up to a subsequence we have $U_k \wto U \in H^s(\R^N)$ and $a_k \to a \in [m_0,m_0+\nu_0]$. We need to show that $\norm{U_k}_{H^s(\R^N)}\to \norm{U}_{H^s(\R^N)}$; the closedness of $\widehat{S}$ will conclude the proof.

As observed in Step 3, we have by the weak convergence that $L'_{a_k}(U_k) U_k=0=L'_a(U)U$; hence, if $R>0$ is some radius to be fixed, we gain
\begin{eqnarray*}
\lefteqn{\left |\left(\norm{(-\Delta)^{s/2}U_k}_2^2 + a_k \norm{U_k}_2^2 \right) - \left(\norm{(-\Delta)^{s/2}U}_2^2 + a \norm{U}_2^2 \right) \right|}\\
 &\leq& \Bigg| \int_{|x| \leq R} f(U_k)U_k dx - \int_{|x|\leq R} f(U) U dx \Bigg | + \\
 &&+ \int_{|x| > R} |f(U_k) U_k| dx + \int_{|x|>R} |f(U)U)| dx =: (I) + (II).
\end{eqnarray*}
Fix now a small $\eta>0$. As regards $(II)$, we have by \eqref{eq_prop_f} %(recall that the definition of $p$ is given in condition (f1.3)) %\tor{va bene anche per $p+1=2^*_s$}
$$\int_{|x|>R} |f(U_k)U_k| dx\leq \beta \int_{|x|>R} |U_k|^2 dx+ C_{\beta} \int_{|x|>R} |U_k|^{p+1} dx< \eta$$
for sufficiently (uniformly in $k$) large $R>0$ thanks to \eqref{eq_stima_uniformeU}; up to taking a larger $R$, it holds also for $U$.

Fixed this $R>0$, focusing on $(I)$, by %classical arguments (see 
Proposition \ref{prop_converg_generiche_loc} %Va bene per il sottocritico non stretto!
we have 
$$ \Bigg| \int_{|x| \leq R} f(U_k)U_k dx - \int_{|x|\leq R} f(U) U dx \Bigg |<\eta$$
for sufficiently large $k =k(R)$. 
Merging together, we obtain
$$\norm{(-\Delta)^{s/2}U_k}_2^2 + a_k \norm{U_k}_2^2 \to \norm{(-\Delta)^{s/2}U}_2^2 + a \norm{U}_2^2$$
which with elementary passages leads to the claim. 

\medskip

\noindent
\textbf{Step 5.}
Finally, we prove $(c)$. By contradiction, there exists an $\eta>0$ such that, for each $n \in \N$ there exists a $U_n \in \widehat{S}$ which satisfies
$$\norm{U_n}_{\R^N \setminus B_n} > \eta.$$
By the compactness, we have, up to a subsequence, $U_n \to U \in \widehat{S}$ as $n \to +\infty$. 
Thus (notice that $\int_{\R^N} \frac{|U(x)-U(\cdot)|^2}{|x-\cdot|^{N+2s}} dx \in L^1(\R^N)$ and absolute integrability of the integral applies)
$$\eta < \norm{U_n}_{\R^N \setminus B_n} \leq \norm{U_n -U}_{H^s(\R^N)} + \norm{U}_{\R^N \setminus B_n} \to 0$$
which is an absurd. 

For the second part, we argue similarly. Indeed, up to a subsequence, $U_n \to U$ in $H^s(\R^N)$ and $\theta_n \to \theta$ in $\R^N$, thus
\begin{eqnarray*}
\lefteqn{\norm{U_n(\cdot-\theta_n)}_{\R^N\setminus B_n}}\\
&\leq& \norm{U_n- U}_{H^s(\R^N)} +\norm{\tau_{\theta_n}U - \tau_{\theta}U}_{H^s(\R^N)} + \norm{U(\cdot-\theta)}_{\R^N\setminus B_n}
\to0,
\end{eqnarray*}
where $\tau_{\theta}$ is the translation. This concludes the proof.
\QED
%
%\bigskip
%

\begin{Remark}
The compactness of the set $\widehat{S}$ of (almost) ground states is somehow expected by thinking at the power case $f(u)=|u|^{p-2}u$: in this case, indeed, the ground state is unique (and nongenerate) \cite{Kwo0, FLS}. On the other hand, in the general case (for examples for suitable sums of powers), uniqueness seems not to be the case \cite{DPG, WeWu}.
%[Wei-Wu '21])
%[D\'avila-del Pino-Guerra '12]
\end{Remark}

Gained compactness, we turn back considering the set of all the solutions (with no restrictions in zero), that is
$$\widehat{S}':= \bigcup_{p\in \R^N} \tau_p(\widehat{S});$$
we observe that $\widehat{S}'$ is bounded. Moreover we define an open $r$-neighborhood of $\widehat{S}'$, reminiscent of the perturbation approach in \cite{FlWe0, ABC, DPW},
$$S(r):= \{ u \in H^s(\R^N) \mid d(u,\widehat{S}')< r\}$$ 
that is
%\begin{eqnarray*}
%\lefteqn{S(r)=}\\
%&&\Big\{ u=U(\cdot-p)+\varphi \in H^s(\R^N) \mid U \in \widehat{S}, \; p \in \R^N, \; \varphi \in H^s(\R^N), \; \norm{\varphi}_{H^s(\R^N)}< r \Big\}.
%\end{eqnarray*}
$$
S(r)=
\Big\{ u=U(\cdot-p)+\varphi \in H^s(\R^N) \mid U \in \widehat{S}, \; p \in \R^N, \; \varphi \in H^s(\R^N), \; \norm{\varphi}_{H^s(\R^N)}< r \Big\}.
$$
In order to re-gain some compactness, we aim to detect and somehow bound the point of translation and the size of the error. To this last goal, we
%We further 
define a \emph{minimal radius} map $\widehat{\rho}: H^s(\R^N) \to \R_+$ by
$$\widehat{\rho}(u):= \inf \left \{ \norm{u-U(\cdot-y)}_{H^s(\R^N)} \mid U \in \widehat{S}, \; y \in \R^N \right\}.$$
We observe
\begin{equation}\label{eq_stima_widerho}
u \in S(r) \implies \widehat{\rho}(u) < r,
\end{equation}
and in addition
$$\widehat{\rho}(u)=\inf\{t\in \R_+ \mid u \in S(t)\},$$
where the infimum on the right-hand side is not attained. Finally, $\widehat{\rho}\in Lip(H^s(\R^N), \R)$ with Lipschitz constant equal to $1$, that is, for every $u,v \in H^s(\R^N)$,
\begin{equation}\label{eq_lipsc_rho}
\abs{\widehat{\rho}(u)-\widehat{\rho}(v)} \leq \norm{u-v}_{H^s(\R^N)}.
\end{equation}
The detection of the point of translation will be instead more tricky, and will be investigated in Section \ref{sec_conc_center}.

\smallskip

We end this Section with two technical lemmas. The first one is a direct consequence of Lemma \ref{lemma_stima_g0}, and allows to link the level $L_{m_0}(u)$ of a whatever function $u\in S(r)$ with the ground state $E_{m_0}$, once $r$ is sufficiently small; this further gives a lower bound for the functional $L_{m_0}$.

\begin{Lemma}\label{lemma_stima_g}
Up to taking a smaller $\nu_0=\nu_0(l_0)>0$,
there exists a sufficiently small $r'=r'(\nu_0, r^*)>0$ such that, for every $u \in S(r')$, we have
$$L_{m_0}(u)\geq g(P_{m_0}(u))E_{m_0}.$$
\end{Lemma}

\claim Proof.
By Lemma \ref{lemma_stima_g0} (a), we know that the inequality holds if $P_{m_0}(u)$ is in a neighborhood of the value $1$. Observe that $P_{a}(U)=1$ if $U\in S_{a}$: by continuity and compactness, $P_{m_0}(U) \approx 1$ if $U\in S_a$ and $a \approx m_0$. In particular, by choosing a small value of $\nu_0$, $P_{m_0}(U) \approx 1$ for $U \in \widehat{S}$. Indeed
$$
P_{m_0}(U)=1+\frac{1}{2s} \frac{N}{N-2s} \; (a-m_0)\frac{ \norm{U}_2^2}{\norm{(-\Delta)^{s/2} U}_2^2} + o(1);
$$
the addendum on the right-hand side can be bounded by the maximum and the minimum over $\widehat{S}$ (notice that $(-\Delta)^{s/2}U$ cannot be zero) and thus we can find a uniform small $\nu_0$. 
Again by continuity we have that $P_{m_0}(u) \approx 1$ for $u\in S(r')$, $r'$ sufficiently small. Indeed
$$
P_{m_0}(u)= 1-\frac{1}{2s} \frac{N}{N-2s}\frac{1}{\norm{(-\Delta)^{s/2}U}_2^2}\big(L_{m_0}(U(\cdot-p)+\varphi)-L_{m_0}(U(\cdot-p)) \big)+ o(1).
$$
This concludes the proof.
\QED

\bigskip

We notice that the condition $E_{m_0+\nu_0}< l_0$ keeps holding by decreasing $\nu_0$, so no ambiguity in the $l_0$-depending choice of $\nu_0$ in Lemma \ref{lemma_stima_g} arises. 
We focus now on the second lemma. 

\begin{Lemma}\label{lem_tec2}
There exist $\nu_1 \in (0, \nu_0)$ and $\delta_0=\delta_0(\nu_1)>0$ such that
$$L_{m_0+\nu_1}(U) \geq E_{m_0}+\delta_0, \quad \textit{uniformly for $U \in \widehat{S}$}.$$
\end{Lemma}

\claim Proof.
Observe first that, since $\widehat{S}$ is compact also in $L^2(\R^N)$, we have also finite and strictly positive minimum $\underline{M}$ and maximum $\overline{M}$ with respect to $\norm{\cdot}_2$. Consider $\nu_1 \in (0, \nu_0)$ such that 
$$E_{m_0+\nu_1}-E_{m_0} > \frac{1}{2}(\nu_0-\nu_1)\overline{M};$$
we notice that such $\nu_1$ exists since, as $\nu_1\to \nu_0^+$, the left hand side positively increases while the right hand side goes to zero. Let now $a\in [m_0, m_0+\nu_0]$; we consider two cases. 
If $a\in [m_0,m_0+\nu_1]$ we argue as follow: for $U\in S_a$ we have
\begin{align*}
L_{m_0+\nu_1}(U) &= L_a(U) + \tfrac{1}{2}(m_0 + \nu_1 - a)\norm{U}_2^2 \\
&\geq E_a + \tfrac{1}{2}(m_0 + \nu_1 - a)\underline{M}=:(I)+(II);
\end{align*}
now, the quantity $(I)$ is minimum when $a=m_0$, while $(II)$ is minimum when $a=m_0+\nu_1$; if both could apply at the same time, we would have as a minimum the quantity $E_{m_0}$. Since it is not possible, we obtain
$$\inf_{U\in \bigcup_{a\in [m_0, m_0+\nu_1]} S_a} L_{m_0+\nu_1}(U) > E_{m_0}.$$
If $a\in (m_0+\nu_1, m_0+\nu_0]$ instead, we have
$$ L_{m_0+\nu_1}(U) \geq E_a - \tfrac{1}{2}(a-(m_0 + \nu_1))\overline{M} \geq E_{m_0+\nu_1} - \tfrac{1}{2}(\nu_0-\nu_1)\overline{M} $$
and thus, by the property on $\nu_1$,
$$\inf_{U\in \bigcup_{a\in [m_0+\nu_1, m_0+\nu_0]} S_a} L_{m_0+\nu_1}(U) \geq E_{m_0+\nu_1} - \tfrac{1}{2}(\nu_0-\nu_1)\overline{M}>E_{m_0}.$$
This concludes the proof, by taking as $\delta_0>0$ the smallest of the two differences.
\QED

\subsection{Fractional center of mass}
\label{sec_conc_center}

As in \cite{CJT}, inspired by \cite{BeCe0, DF2, BT1}, we want to define a \emph{barycentric map} $\Upsilon$ which, given a function $u=U(\cdot-p) + \varphi \in S(r)$, gives an estimate on the maximum point $p$ of $U(\cdot-p)$; since $\varphi$ is small and $U$ decays (polynomially) at infinity, $p$ is, in some ways, the \emph{center of mass} of $u$. The idea will be to bound $\Upsilon(u)$ in order to re-gain compactness.

Since the nonlocality comes into the very definition of the ambient space, we need the use of the norm \eqref{eq_def_Htilde}, which we notice being stronger than the one induced by the Gagliardo seminorm.

\begin{Lemma} \label{lem_def_bar}
Let $r^*$ be as in Lemma \ref{lemma_stima_unif}. Then there exist a sufficiently large $R_0>0$, a sufficiently small radius $r_0\in (0, r^*)$ and a continuous map
$$\Upsilon: S(r_0) \to \R^N$$
such that, for each $u=U(\cdot-p)+\varphi \in S(r_0)$ we have
$$|\Upsilon(u)-p|\leq 2R_0.$$
Moreover, $\Upsilon$ is continuous and
$-\Upsilon$ is \emph{shift-equivariant}, that is, $\Upsilon(u(\cdot+\xi))=\Upsilon(u)-\xi$ for every $u \in S(r_0)$ and $\xi\in \R^N$.
\end{Lemma}

\claim Proof.
Recalled that $r^*=\min_{U\in \widehat{S}} \norm{U}_{H^s(\R^N)}>0$, we have by Lemma \ref{lemma_stima_unif}
\begin{equation}\label{eq_dim_cent_mass}
\norm{U}_{\R^N\setminus B_{R_0}} <\tfrac{1}{8}r^*, \quad \hbox{ uniformly for $U \in \widehat{S}$}
\end{equation}
for $R_0 \gg 0$. Thus
$$r^* \leq \norm{U}_{H^s(\R^N)} \leq \norm{U}_{B_{R_0}} + \norm{U}_{\R^N\setminus B_{R_0}} < \norm{U}_{B_{R_0}} + \tfrac{1}{8}r^*$$
which implies
$$\norm{U}_{B_{R_0}} > \tfrac{7}{8} r^* \quad \hbox{and} \quad \norm{U}_{\R^N\setminus B_{R_0}}<\tfrac{1}{8} r^* $$
for each $U \in \widehat{S}$. Consider now a cutoff function $\psi \in C^{\infty}_c(\R_+)$ such that
$$[0,\tfrac{1}{4} r^*] \prec \psi \prec [\tfrac{1}{2} r^*,+\infty);$$
let $r_0\in (0, \frac{1}{8} r^*)$ and define, for each $u\in S(r_0)$ and $q\in \R^N$, a \emph{density} function
{\tiny }$$d(q,u):= \psi \left(\inf_{\tilde{U}\in \widehat{S}} \norm{u-\tilde{U}(\cdot-q)}_{B_{R_0}(q)}\right).$$
Notice that $d(\cdot, u)$ is an integrable function: indeed $[u(\cdot + \xi)]_{B_{R_0}(q)} = [u]_{B_{R_0}(q+\xi)}$ for $\xi \in \R^N$, and $q\mapsto \norm{u-\tilde{U}(\cdot-q)}_{B_{R_0}(q)}=\norm{\tau_q u-\tilde{U}}_{B_{R_0}}$ is continuous by
$$\abs{\norm{\tau_q u-\tilde{U}}_{B_{R_0}}-\norm{\tau_pu -\tilde{U}}_{B_{R_0}}} \leq \norm{\tau_q u - \tau_p u}_{B_{R_0}} \leq \norm{\tau_q u - \tau_p u}_{H^s(\R^N)} \to 0$$
as $p\to q$, and the infimum over continuous functions is upper semicontinuous. 

If we show that $d(\cdot,u)\geq 0$ is not identically zero and it has compact support, then it will be well defined the quantity
$$\Upsilon(u):=\frac{\int_{\R^N} q \, d(q,u) dq}{\int_{\R^N} d(q,u) dq}.$$
We show first that $d(\cdot,u)$ has compact support. Indeed if $u=U(\cdot-p)+\varphi$ and $\tilde{U}\in \widehat{S}$ is arbitrary, then
\begin{eqnarray*}
\lefteqn{\norm{u-\tilde{U}(\cdot-q)}_{B_{R_0}(q)}} \\
&\geq & \norm{\tilde{U}(\cdot-q)}_{B_{R_0}(q)} - \norm{U(\cdot-p)}_{B_{R_0}(q)} - \norm{\varphi}_{B_{R_0}(q)} \\
&\geq & \norm{\tilde{U}}_{B_{R_0}} - \norm{U}_{B_{R_0}(q-p)} -\norm{\varphi}_{H^s(\R^N)} \geq \tfrac{6}{8} r^* - \norm{U}_{B_{R_0}(q-p)};
\end{eqnarray*}
take now $q\notin B_{2R_0}(p)$: if $x \in B_{R_0}(q-p)$, by the fact that $|x-(q-p)|<R_0$ and $|q-p|\geq 2R_0$, we obtain that $\abs{x} \geq R_0$, that is, $x \in \R^N \setminus B_{R_0}$. Therefore by \eqref{eq_dim_cent_mass}
\begin{equation}\label{eq_dim_58} 
\norm{u-\tilde{U}(\cdot-q)}_{B_{R_0}(q)} \geq \tfrac{6}{8} r^* - \norm{U}_{\R^N \setminus B_{R_0}} \geq \tfrac{5}{8} r^* > \tfrac{1}{2} r^* 
\end{equation}
thus $\inf_{U\in \widehat{S}} \norm{u-\tilde{U}(\cdot-q)}_{B_{R_0}(q)}\geq \frac{1}{2} r^*$ and hence $d(q,u)=0$ for $q\notin B_{2R_0}(p)$; this means that $\supp(d(\cdot,u))\subset B_{2R_0}(p)$.

We show next that $d(\cdot, u)$ is equal to $1$ on a ball. Indeed if $u=U(\cdot-p)+\varphi$
\begin{eqnarray*}
\lefteqn{\inf_{\tilde{U}\in \widehat{S}} \norm{u-\tilde{U}(\cdot-q)}_{B_{R_0}(q)} \leq \norm{u-U(\cdot-q)}_{B_{R_0}(q)}}\\
& \leq& \norm{U(\cdot-p)-U(\cdot-q)}_{B_{R_0}(q)} + \norm{\varphi}_{B_{R_0}(q)} \leq \norm{\tau_{p-q}U-U}_{B_{R_0}} + \tfrac{1}{8}r^*.
\end{eqnarray*}
We can make the first term as small as we want by taking $|p-q|$ small, that is 
$$\inf_{U\in \widehat{S}} \norm{u-\tilde{U}(\cdot-q)}_{B_{R_0}(q)} \leq \tfrac{1}{4}r^*$$
 for $q\in B_r(p)$, $r$ small, which implies $d(q,u)=1$.

By the fact that $B_r(p)\subset \supp(d(\cdot, u))\subset B_{2R_0}(p)$ we have the well posedness of $\Upsilon(u)$ and
$$\Upsilon(u) = \frac{\int_{B_{2R_0}(p)} q \, d(q,u) dq}{\int_{B_{2R_0}(p)} d(q,u) dq}.$$
The main property comes straightforward, as well as the shift equivariance. We show now the continuity. 
Indeed, assume $\norm{u-v}_{H^s(\R^N)}\leq \frac{1}{8}r^*$. Then, by \eqref{eq_dim_58},
$$\norm{v-\tilde{U}(\cdot -q)}_{B_{R_0}(q)} \geq \norm{u-\tilde{U}(\cdot -q)}_{B_{R_0}(q)} - \norm{v-u}_{B_{R_0}(q)} \geq \tfrac{1}{2} r^*$$
and again we can conclude that $\supp(d(\cdot, v))\subset B_{2R_0}(p)$ for each $\norm{u-v}_{H^s(\R^N)}\leq \tfrac{1}{8} r^*$, where $p$ depends only on $u$. Moreover, observe that $\int_{B_{2R_0}(p)} d(q,u) dq \geq \int_{B_{r}(p)}1 \,dq \geq |B_{r}|=:C_1$ not depending on $u$ and $p$ (and similarly $C_2:= |B_{2R_0}|$), and that $d(q, \cdot)$ is Lipschitz 
(since $\psi$ and the norm are so, and the infimum over a family of Lipschitz functions is still Lipschitz). 
Thus we have
\begin{align*}
\abs{\Upsilon(u)-\Upsilon(v)} &\leq \frac{\int_{B_{2R_0}(p)} \abs{q} \, \abs{d(q,u)-d(q,v)} dq}{\int_{B_{2R_0}(p)} d(q,u) dq} +\\
& \quad+\int_{B_{2R_0}(p)} \abs{q} \, d(q,v) dq\frac{\int_{B_{2R_0}(p)} \abs{d(q,v)-d(q,u)} dq}{\int_{B_{2R_0}(p)} d(q,u) dq \int_{B_{2R_0}(p)} d(q,v) dq}\\
&\leq \int_{B_{2R_0}(p)} \abs{q} dq \frac{1}{C_1}\left(1+ \frac{C_2}{C_1}\right)\norm{u-v}_{H^s(\R^N)}
=: C_p \norm{u-v}_{H^s(\R^N)}.
\end{align*}
Since $C_p$ can be bounded above by a constant of the type $C(1+\Upsilon(u))$, we have
$$\norm{u-v}_{H^s(\R^N)}\leq r_0 \implies \abs{\Upsilon(u)-\Upsilon(v)} \leq C(1+\Upsilon(u)) \norm{u-v}_{H^s(\R^N)};$$
in particular this implies the continuity.
\QED

\section{Singularly perturbed equation}
\label{sez_singur_perturb}

We come back now to our equation
\begin{equation}\label{eq_princ_epsx}
(-\Delta)^s u + V(\eps x) u= f(u), \quad x \in \R^N.
\end{equation}
It is known that the solutions of \eqref{eq_princ_epsx} can be characterized as critical points of the functional $I_{\eps}: H^s(\R^N) \to \R$
$$I_{\eps}(u):=\frac{1}{2} \norm{(-\Delta)^{s/2} u}_2^2 + \frac{1}{2} \int_{\R^N} V(\eps x) u^2 dx - \int_{\R^N}F(u)dx, \quad u \in H^s(\R^N) $$
where $I_{\eps}\in C^1(H^s(\R^N),\R)$, since $\norm{\cdot}_{H_{\eps}^s(\R^N)}$ is a norm.

We start with a technical result. Let $\nu_1$ be as in Lemma \ref{lem_tec2}; we want to show that the claim of the lemma 
continues holding, for $\eps$ small, if we replace $L_{m_0+\nu_1}$ with $I_{\eps}$, and $\widehat{S}$ with $S(r_0')\cap \{\eps \Upsilon(u) \in \Omega[\nu_1,\nu_0]\}$, $r_0'$ small.

\begin{Lemma}\label{lemma_stimabasso_I}
Let $\nu_1$ and $\delta_0$ be as in Lemma \ref{lem_tec2}. 
Then there exist $\delta_1 \in (0, \delta_0)$ and $r_0'=r_0'(\delta_1) \in (0,r_0)$ sufficiently small, such that for every $\eps$ small we have 
$$I_{\eps}(u)\geq E_{m_0} + \delta_1$$
for each $u \in \{u \in S(r_0')\mid \eps \Upsilon(u) \in \Omega[\nu_1,+\infty)\} \supset \{u \in S(r_0')\mid \eps \Upsilon(u) \in \Omega[\nu_1,\nu_0]\}$.
\end{Lemma}

\claim Proof.
First we improve Lemma \ref{lem_tec2} for $L_a$ in the direction of the nonautonomous equation. Indeed, by the assumption, we have $V(\eps \Upsilon(u))\geq m_0+\nu_1$, that is
$$L_{V(\eps \Upsilon(u))}(U) \geq L_{m_0+\nu_1}(U) \geq E_{m_0} + \delta_0$$
for any $U\in \widehat{S}$. Moreover, if $u=\tilde{U}(\cdot - p) + \tilde{\varphi}\in S(r_0)$ then, by Lemma \ref{lem_def_bar}, $\eps p \in \Omega_{2\eps R_0}\subset \Omega_{2R_0}$ which is compact. 
 By uniform continuity of $V$ and boundedness from above of $\widehat{S}$, we have
\begin{equation}\label{eq_dim_nonaut_to_aut}
L_{V(\eps p)}(U) \geq E_{m_0} + \delta_0/2
\end{equation}
for all $U\in \widehat{S}$ and $\eps$ small enough. 

Let now $r_0'$ to be fixed and $u=U(\cdot-p)+\varphi\in S(r_1)$. Then we have
$$I_{\eps}(u) = I_{\eps}(U(\cdot-p)+\varphi) = I_{\eps}(U(\cdot-p)) + I'_{\eps}(v)\varphi$$
for some $v\in H^s(\R^N)$ in the segment $[U(\cdot-p), u]$. Notice that $v$ lies in a ball of radius $\max \widehat{S}+r_0'$ and $I'_{\eps}$ sends bounded sets in bounded sets (uniformly on $\eps$); thus there exists a constant $C$, not depending on $U$, $p$ and $\varphi$, such that
\begin{equation}\label{eq_dim_stimaI}
I_{\eps}(u) \geq I_{\eps}(U(\cdot-p)) - C \norm{\varphi}_{H^s(\R^N)} \geq I_{\eps}(U(\cdot-p)) - \delta_1/2
\end{equation}
for $\norm{\varphi}_{H^s(\R^N)}< r_0'$ sufficiently small. 
Recalled that $\eps p \in \Omega_{2 R_0}$ we have, by the uniform continuity of $V$ and the uniform estimate \eqref{eq_stima_uniformeU}, for sufficiently small $\eps$, 
$$I_{\eps}(U(\cdot-p)) \geq L_{V(\eps p)}(U) - \delta_1/2$$
and the claim comes from \eqref{eq_dim_stimaI} and \eqref{eq_dim_nonaut_to_aut}, since, for $\delta_1< \delta_0/4$, 
$$I_{\eps}(u) \geq E_{m_0}+\delta_0/2- \delta_1 \geq E_{m_0} + \delta_1. 
\QED
$$

\bigskip

Before introducing the penalized functional, we state another technical lemma, which gives a (trivial, but useful) lower bound for $I'_{\eps}(v)v$ for small values of $v\in H^s(\R^N)$.

\begin{Lemma}\label{lemma_stima_I}
There exists $r_1>0$ sufficiently small %(see also the Remark below) 
and a constant $C>0$ such that
\begin{equation}\label{eq_stima_I}
I'_{\eps}(v)v \geq C\norm{v}_{H^s(\R^N)}^2
\end{equation}
for every $\eps>0$ and $v\in H^s(\R^N)$ with $\norm{v}_{H^s(\R^N)}\leq r_1$.
\end{Lemma}

\claim Proof.
We have, by \eqref{eq_prop_f} with $\beta < \frac{1}{2}\underline{V}$,
\begin{align}
I'_{\eps}(v)v &\geq \norm{(-\Delta)^{s/2}v}_2^2 + \int_{\R^N} \underline{V} v^2 dx - \beta \norm{v}_2^2 - C_{\beta} \norm{v}_{p+1}^{p+1} \notag \\
&\geq \norm{(-\Delta)^{s/2}v}_2^2 + \frac{1}{2}\underline{V} \norm{v}_2^2 - C_{\beta} \norm{v}_{p+1}^{p+1} \label{eq_dim_da_migliorare}\\
&\geq C\norm{v}_{H^s(\R^N)}^2 - C_{\beta} \norm{v}_{H^s(\R^N)}^{p+1}
\geq C'\norm{v}_{H^s(\R^N)}^2 \notag
\end{align}
where the last inequality holds for $\norm{v}_{H^s(\R^N)}$ small, since $p+1>2$.
\QED

\medskip

\begin{Remark} 
For a later use, we observe that one can improve \eqref{eq_dim_da_migliorare} by
\begin{equation}\label{eq_stima_2p}
\norm{(-\Delta)^{s/2}v}_2^2 + \frac{1}{2}\underline{V} \norm{v}_2^2 - 2^p C_{\beta} \norm{v}_{p+1}^{p+1} \geq C\norm{v}_{H^s(\R^N)}^2
\end{equation}
up to choosing a smaller $r_1$.
\end{Remark}

\subsection{A mass-concentrating penalization} %penalized functional}

We want to study now a penalized functional (see \cite{ByWa0,BJ0, CJT}), that is $I_{\eps}$ plus a term which forces solutions to stay in $\Omega$.
%; we highlight that if we had information on $V$ at infinity (for example, trapping potentials) we would not have need of this penalization. 

Since $V>m_0$ on $\partial \Omega$, we can find an annulus around $\partial \Omega$ where this relation keeps holding, that is
$$V(x)>m_0, \quad \hbox{ for $x\in \overline{\Omega_{2h_0}\setminus \Omega}$}$$
for $h_0$ sufficiently small. We then define the \emph{mass-concentrating} penalization functional $Q_{\eps}: H^s(\R^N) \to \R$
$$Q_{\eps}(u):=\left(\frac{1}{\eps^{\alpha}} \norm{u}^2_{L^2(\R^N\setminus (\Omega_{2h_0}/\eps))}-1\right)_+^{\frac{p+1}{2}}, \quad u \in H^s(\R^N)$$
where $\alpha \in (0, \min\{1/2, s\})$. 

We observe that, for every $u,v \in H^s(\R^N)$,
$$Q'_{\eps}(u)v = \frac{(p+1)}{\eps^{\alpha}} \left(\frac{1}{\eps^{\alpha}} \norm{u}^2_{L^2(\R^N\setminus (\Omega_{2h_0}/\eps))}-1\right)_+^{\frac{p-1}{2}} \int_{\R^N\setminus (\Omega_{2h_0}/\eps)} u v \, dx $$
and it is straightforward to prove the following estimate
\begin{equation} \label{eq_stima_Q}
Q'_{\eps}(u)u \geq (p+1) Q_{\eps}(u).
\end{equation}
We thus set
$$J_{\eps}:=I_{\eps} + Q_{\eps}$$
the \emph{penalized} functional. It results that $Q_{\eps}$ and $J_{\eps}$ are in $C^1(H^s(\R^N), \R)$.

We want to find critical points of $J_{\eps}$ and show, afterwards, that these critical points, under suitable assumptions, are critical points of $I_{\eps}$ too, since $Q_{\eps}$ will be identically zero. 
Let $\eps=1$: observed that $Q_1(u)$ vanishes if $u$ have much mass inside $\Omega$, we see that $J_1(u)=I_1(u)$ holds when the mass of $u$ \emph{concentrates} in $\Omega$; this is why we say that $Q_{\eps}$ \emph{forces} $u$ to stay in $\Omega$. Similarly, as $\eps \to 0$, much less mass must be found outside $\Omega/\eps$.

We start by two technical lemmas. The first one gives a sufficient condition to pass from weak to strong convergent sequences in a Hilbert space, similarly to the convergence of the norms.

\begin{Lemma}\label{lemma_conv_forte} 
Fix $\eps>0$ and let $(u_j)_j\subset H^s(\R^N)$ be such that
\begin{equation}\label{eq_ip_conv_J'}
\norm{J'_{\eps}(u_j)}_{(H^s(\R^N))^*} \to 0 \quad \textit{ as $j \to +\infty$}.
\end{equation}
Assume moreover that $u_j \wto u_0$ in $H^s(\R^N)$ as $j \to +\infty$, and that
\begin{equation}\label{eq_doppio_limite}
\lim_{R,j \to +\infty} \norm{u_j}_{L^q(\R^N \setminus B_R)}=0
\end{equation}
for $q=2$ and $q=p+1$. 
Then $u_j \to u_0$ in $H^s(\R^N)$ as $j \to +\infty$.
\end{Lemma}

\claim Proof. 
We have by the weak lower semicontinuity of the norm
\begin{equation}\label{eq_semicont_inf}
\liminf_{j\to +\infty} \norm{u_j}_{H_{\eps}^s(\R^N)} \geq \norm{u_0}_{H_{\eps}^s(\R^N)}.
\end{equation}
Moreover
%\begin{eqnarray*}
%\lefteqn{\norm{u_j}_{H_{\eps}^s(\R^N)}^2 = \int_{\R^N} f(u_j) u_j dx + I'_{\eps}(u_j)u_j } \\
%&=& \Bigg( \int_{\R^N} f(u_j) u_j dx -\int_{\R^N} f(u_0)u_0 dx \Bigg) +\Big( I'_{\eps}(u_j)u_j -I'_{\eps}(u_0)u_0\Big) +\\
%&&+ I'_{\eps}(u_0)u_0 + \int_{\R^N} f(u_0)u_0 dx =: (I) + (II) + \norm{u_0}_{H_{\eps}^s(\R^N)}^2;
%\end{eqnarray*}
\begin{align*}
\norm{u_j}_{H_{\eps}^s(\R^N)}^2 =& \int_{\R^N} f(u_j) u_j dx + I'_{\eps}(u_j)u_j \\
=& \, \Bigg( \int_{\R^N} f(u_j) u_j dx -\int_{\R^N} f(u_0)u_0 dx \Bigg) +\Big( I'_{\eps}(u_j)u_j -I'_{\eps}(u_0)u_0\Big) +\\
&+ I'_{\eps}(u_0)u_0 + \int_{\R^N} f(u_0)u_0 dx =: (I) + (II) + \norm{u_0}_{H_{\eps}^s(\R^N)}^2;
\end{align*}
if we prove that
$$\limsup_{j \to +\infty} \big( (I) + (II) \big) \leq 0$$
we are done, because together with \eqref{eq_semicont_inf} we obtain
$$\norm{u_j}_{H_{\eps}^s(\R^N)} \to \norm{u_0}_{H_{\eps}^s(\R^N)} \quad \hbox{ as $j \to +\infty$},$$
which implies the claim, since $H^s_{\eps}(\R^N)$ is a Hilbert space.

Focus on $(I)$; we have
\begin{eqnarray*}
\lefteqn{\int_{\R^N} (f(u_j)u_j - f(u_0)u_0) dx = \int_{B_R} (f(u_j)u_j - f(u_0)u_0) dx +}\\
&& + \int_{\R^N\setminus B_R} (f(u_j)u_j - f(u_0)u_0) dx =:(I_1) + (I_2).
\end{eqnarray*}
The piece $(I_2)$ can be made small for $j$ and $R$ sufficiently large, by exploiting the estimates on $f$, assumption \eqref{eq_doppio_limite} and the absolute continuity of the Lebesgue integral for $u_0$. %\tor{VEDI dettagli - vanno bene per sottocritico non stretto?}. 
For such large $R$ and $j$, up to taking a larger $j$, we can make the piece $(I_1)$ small by Proposition \ref{prop_converg_generiche_loc}. % classical arguments. %Va bene per sottocritico non stretto!

Focus now on $(II)$; we first observe that by exploiting H\"older inequalities and again classical arguments we have $I'_{\eps}(u_j)u_0 \to I'_{\eps}(u_0)u_0$.
Thus we have, by \eqref{eq_ip_conv_J'},
 \begin{eqnarray*}
\lefteqn{\limsup_{j \to +\infty} \left(I'_{\eps}(u_j) u_j - I'_{\eps}(u_0)u_0 \right)
= -\liminf_{j \to +\infty} \left(Q'_{\eps}(u_j) u_j - Q'_{\eps}(u_j)u_0 \right)} \\
&=& -\liminf_{j \to +\infty} \Bigg( \left(\frac{1}{\eps^{\alpha}} \norm{u_j}^2_{L^2(\R^N\setminus (\Omega_{2h_0}/\eps))}-1\right)_+^{\frac{p-1}{2}} \cdot \\
&& \cdot \Bigg(\int_{\R^N\setminus (\Omega_{2h_0}/\eps)}u_j^2 dx -\int_{\R^N\setminus (\Omega_{2h_0}/\eps)} u_j u_0 dx \Bigg) \Bigg) \leq 0
 \end{eqnarray*}
where the last inequality is due to the following fact: observe first that $u_j \wto u_0$ in $H^s_{\eps}(\R^N)\hookrightarrow L^2(\R^N)$ thus (by restriction) $u_j \wto u_0$ in $L^2(\R^n \setminus (\Omega_{2h_0}/\eps))$; by definition of weak convergence and by the lower semicontinuity of the norm, we have
$$\liminf_{j \to +\infty} \int_{\R^N\setminus (\Omega_{2h_0}/\eps)}u_j^2 dx \geq \int_{\R^N\setminus (\Omega_{2h_0}/\eps)}u_0^2 dx = \lim_{j \to +\infty} \int_{\R^N\setminus (\Omega_{2h_0}/\eps)}u_j u_0 dx ,$$
that is
$$\liminf_{j \to +\infty} \Bigg(\int_{\R^N\setminus (\Omega_{2h_0}/\eps)}u_j^2dx -\int_{\R^N\setminus (\Omega_{2h_0}/\eps)} u_j u_0 dx \Bigg) \geq 0.$$
Noticed that $a_n \geq 0$ and $\liminf_n b_n \geq 0$ imply $\liminf_n (a_n b_n) \geq 0$, we conclude.
\QED

\bigskip

The second Lemma is a lower bound for $J_{\eps}$ with respect to the functional $L_{m_0}$.
We highlight that in what follows we understand that the case $m_0=\underline{V}$, i.e. $m_0$ \emph{global} minimum, gives rise to a not-perturbed result.

\begin{Lemma}\label{lemma_stima_J}
Set $C_{min}:= \frac{1}{2}(m_0-\underline{V})\geq 0$ we have, for $\eps$ small and $u\in H^s(\R^N)$,
$$J_{\eps}(u) \geq L_{m_0}(u) - C_{min} \eps^{\alpha}.$$
\end{Lemma}

\claim Proof.
We have, recalling that $m_0$ is the infimum of $V$ over $\Omega_{2h_0}$ and $\underline{V}$ is the infimum over $\R^N$,
\begin{align*}
J_{\eps}(u) &= L_{m_0}(u) + \frac{1}{2} \int_{\R^N} (V(\eps x) - m_0) u^2 dx + Q_{\eps}(u)\\
&\geq L_{m_0}(u) + \frac{1}{2} \int_{\R^N\setminus(\Omega_{2h_0} / \eps)} (V(\eps x) - m_0) u^2 dx + Q_{\eps}(u)\\
&\geq L_{m_0}(u) - C_{min} \norm{u}_{L^2(\R^N\setminus (\Omega_{2h_0}/\eps))}^2 + Q_{\eps}(u).
\end{align*}
If $\norm{u}_{L^2(\R^N\setminus (\Omega_{2h_0}/\eps))}^2 \leq 2\eps^{\alpha}$ we have the claim by the positivity of $Q_{\eps}(u)$. 
If instead $\norm{u}_{L^2(\R^N\setminus (\Omega_{2h_0}/\eps))}^2 \geq 2\eps^{\alpha}$, then
$$ Q_{\eps}(u) \geq \left(\frac{1}{2\eps^{\alpha}} \norm{u}^2_{L^2(\R^N\setminus (\Omega_{2h_0}/\eps))}\right)^{\frac{p+1}{2}}
\geq \frac{1}{2\eps^{\alpha}} \norm{u}^2_{L^2(\R^N\setminus (\Omega_{2h_0}/\eps))}. $$
Thus
$$ J_{\eps}(u) \geq L_{m_0} + \left (\frac{1}{2\eps^{\alpha}} - C_{min} \right) \norm{u}_{L^2(\R^N\setminus (\Omega_{2h_0}/\eps))}^2 
\geq L_{m_0} $$
for $\eps$ small. This concludes the proof.
\QED

\subsection{Critical points and truncated Palais-Smale condition}

In order to get critical points of $J_{\eps}$ we want to implement a deformation argument. As usual, we need a uniform estimate from below of $\norm{J'_{\eps}(u)}_{(H^s(\R^N))^*}$, and this is the next goal. 

First, by the strict monotonicity of $E_a$, let us fix $l_0'=l_0'(\nu_1)>0$ such that
$$E_{m_0}< l_0' < E_{m_0 + \nu_1};$$
as well as $\nu_0$ and $l_0$, even $l_0'$ will be let vary as $(l_0')^n \to E_{m_0}$ in the proof of the existence.

\begin{Lemma}\label{lemma_stima_basso} 
Let $r_0$ and $r_1$ be as in Lemma \ref{lem_def_bar} and Lemma \ref{lemma_stima_I}. There exists $r_2'\in (0,\min\{r_0, r_1\})$ sufficiently small with the following property: let $0<\rho_1<\rho_0\leq r_2'$ and $(u_{\eps})_{\eps}\subset S(r_2')$ be such that
\begin{eqnarray}
&\quad \norm{J'_{\eps}(u_{\eps})}_{(H^s(\R^N))^*}\to 0 \quad \textit{ as $\eps \to 0$},& \label{eq_PS2}\\ 
&J_{\eps}(u_{\eps})\leq l_0' < E_{m_0+\nu_1}, \quad \textit{ for any $\eps >0$},& \label{eq_PS1} 
\end{eqnarray}
with the additional assumption
$$(\widehat{\rho}(u_{\eps}))_{\eps} \subset [0,\rho_0], \quad (\eps\Upsilon(u_{\eps}))_{\eps} \subset \Omega[0,\nu_0].$$
Then, for $\eps$ small
$$\widehat{\rho}(u_{\eps})\in [0,\rho_1], \quad \eps\Upsilon(u_{\eps})\in \Omega[0,\nu_1].$$
\end{Lemma}

We notice, by \eqref{eq_PS2} and \eqref{eq_PS1}, that $(u_{\eps})_{\eps}$ resembles a particular (truncated) Palais-Smale sequence. As an immediate consequence of the Lemma, set the sublevel
$$J_{\eps}^c:= \{ u \in H^s(\R^N) \mid J_{\eps}(u) \leq c\}$$
 we have the following theorem.

\begin{Theorem}\label{theorem_stima_J'}
There exists $r_2'\in (0,\min\{r_0, r_1\})$ sufficiently small with the following property: if $0<\rho_1<\rho_0\leq r_2'$, then there exists a $\delta_2=\delta_2(\rho_0,\rho_1)>0$ such that, for $\eps$ small
$$\norm{J'_{\eps}(u)}_{(H^s(\R^N))^*}\geq \delta_2$$
for any 
\begin{align*}
u &\in \left\{u\in S(r_2')\cap J_{\eps}^{l_0'} \mid (\widehat{\rho}(u), \eps \Upsilon(u))\in \big((0,\rho_0]\times \Omega[0,\nu_0]\big)\setminus\big([0,\rho_1]\times \Omega[0,\nu_1]\big)\right\} \\
&\supset \left\{u\in S(r_2')\cap J_{\eps}^{l_0'} \mid \rho_1<\widehat{\rho}(u)\leq \rho_0, \; \eps \Upsilon(u)\in \Omega(\nu_1,\nu_0]\right\}.
\end{align*}
\end{Theorem}

\begin{Remark}\label{remark_triplanorma}
Arguing as in the last part of the proof of Lemma \ref{lemma_stima_unif}, noticed that $\widehat{S}$ is compact not only in $H^s(\R^N)$ but also in $L^q(\R^N)$ for $q \in [2, 2^*_s]$, if $(U_n)_n\subset \widehat{S}$ and $(\theta_n)_n \subset \R^N$ is included in a compact, we have
$$\lim_{n\to +\infty} \tnorm{U_n(\cdot + \theta_n)}_{\R^N \setminus B_n}=0,$$
where the norm $\tnorm{\cdot}_{\R^N\setminus B_n}$ is defined in \eqref{eq_def_triplanorm}.
\end{Remark}

\claim Proof of Lemma \ref{lemma_stima_basso}. 
We use the notation, for $h>0$,
$$\Omega^{\eps}_{h}:=(\Omega_{\eps h})/\eps = (\Omega/\eps)_h$$
and notice that if $h<h'$ then $\Omega/\eps \subset\Omega^{\eps}_{h} \subset \Omega^{\eps}_{h'}$. 
Let $r_2'<\min\{r_0, r_1\}$ to be fixed.

\medskip

\noindent 
\textbf{Step 1.} \emph{An estimate for $u_{\eps}$.}
\\
We have, for $u_{\eps}= U_{\eps}(\cdot-p_{\eps}) + \varphi_{\eps} \in S(r_2')$, 
\begin{align*}
\tnorm{u_{\eps}}_{\R^N\setminus (\Omega/\eps)} &\leq 
\tnorm{U_{\eps}(\cdot-p_{\eps})}_{\R^N\setminus (\Omega/\eps)} + \tnorm{\varphi_{\eps}}_{\R^N\setminus (\Omega/\eps)}\\
&\leq \tnorm{U_{\eps}(\cdot-p_{\eps}+\Upsilon(u_{\eps}))}_{\R^N\setminus (\Omega/\eps-\Upsilon(u_{\eps}))} + C
\norm{\varphi_{\eps}}_{H^s(\R^N)}\\
&\leq \tnorm{U_{\eps}(\cdot-p_{\eps}+\Upsilon(u_{\eps}))}_{\R^N\setminus (\Omega/\eps-\Upsilon(u_{\eps}))} + C
 r_2'.
\end{align*}
By the fact that $\eps \Upsilon(u_{\eps})\in \Omega[\nu_1,\nu_0]\subset \Omega$, we have that $0 \in \Omega/\eps - \Upsilon(u_{\eps})$ and thus $\Omega/\eps - \Upsilon(u_{\eps})$ expands in $\R^N$ as $\eps \to 0$. Moreover by Lemma \ref{lem_def_bar} we have $\theta_{\eps}:= \Upsilon(u_{\eps})- p_{\eps} \in B_{2R_0}$ compact. By Remark \ref{remark_triplanorma}, for $\varepsilon$ small we have
\begin{equation} \label{eq_stima_princ_ueps}
\tnorm{u_{\eps}}_{\R^N\setminus (\Omega/\eps)}\leq (1+C)r_2' =C' r_2'.
\end{equation}
Let 
$$n_{\eps}:= \left[ \frac{\sqrt{1+4h_0/\eps}+1}{2}\right]\in \N$$
which by definition satisfies $\eps n_{\eps}(n_{\eps}+1)\leq h_0$ and $n_{\eps}\to+ \infty$ as $\eps \to 0$. We have
$$\sum_{i=1}^{n_{\eps}} \norm{u_{\eps}}_{L^2(\Omega_{n_{\eps}(i+1)}^{\eps}\setminus \Omega_{n_{\eps} i}^{\eps})}^2 
\leq 
 \norm{u_{\eps}}_{L^2(\Omega_{n_{\eps}(n_{\eps}+1)}^{\eps}\setminus \Omega_{n_{\eps}}^{\eps})}^2 \leq (C' r_2')^2$$
and similarly 
$$\sum_{i=1}^{n_{\eps}} [u_{\eps}]_{\Omega_{n_{\eps}(i+1)}^{\eps}\setminus \Omega_{n_{\eps} i}^{\eps}, \R^N}^2 \leq (C'
 r_2')^2, \qquad \sum_{i=1}^{n_{\eps}} \norm{u_{\eps}}_{L^{p+1}(\Omega_{n_{\eps}(i+1)}^{\eps}\setminus \Omega_{n_{\eps} i}^{\eps})}^{p+1} \leq (C'
 r_2')^{p+1}$$
thus, for some $C=C(r_2')$,
$$\sum_{i=1}^{n_{\eps}} \left( 
\norm{u_{\eps}}_{\Omega_{n_{\eps}(i+1)}^{\eps}\setminus \Omega_{n_{\eps} i}^{\eps}}^2 + \norm{u_{\eps}}_{L^{p+1}(\Omega_{n_{\eps}(i+1)}^{\eps}\setminus \Omega_{n_{\eps} i}^{\eps})}^{p+1} \right) \leq C. $$
This implies that there exists $i_{\eps} \in \{1, \dots, n_{\eps}\}$ such that
\begin{equation}\label{eq_stima_anello}
\norm{u_{\eps}}_{A^{\eps}}^2 + \norm{u_{\eps}}_{L^{p+1}(A^{\eps})}^{p+1} \leq \frac{C}{n_{\eps}}\to 0 \quad \hbox{ as $\eps \to 0$},
\end{equation}
where 
$$A^{\eps}:= \Omega_{n_{\eps}(i_{\eps}+1)}^{\eps}\setminus \Omega_{n_{\eps} i_{\eps}}^{\eps}$$
and $C$ depends on $r_2'$ (we will omit this dependence).

\smallskip

\noindent \textbf{Step 2.} \textit{Split the sequence.}
\\
Consider cutoff functions $\varphi_{\eps} \in C^{\infty}_c(\R^N)$ 
$$\Omega^{\eps}_{n_{\eps} i_{\eps}} \prec \varphi_{\eps} \prec \Omega^{\eps}_{n_{\eps}(i_{\eps}+1)}$$
such that $\norm{\nabla \varphi_{\eps}}_{\infty} \leq \frac{C}{n_{\eps}}= o(1)$ as $\eps \to 0$ (which is possible because the distance between $\Omega^{\eps}_{n_{\eps} i_{\eps}}$ and $\Omega^{\eps}_{n_{\eps}(i_{\eps}+1)}$ is $n_{\eps}\to +\infty$).

Define 
$$u_{\eps}^{(1)}:=\varphi_{\eps} u_{\eps}, \quad u_{\eps}^{(2)}:=(1-\varphi_{\eps}) u_{\eps}\quad \hbox{ and } \; u_{\eps}=u_{\eps}^{(1)}+u_{\eps}^{(2)};$$
notice that both $\supp \big(u_{\eps}^{(1)} u_{\eps}^{(2)}\big)$ and $\supp(F(u_{\eps}) -F(u_{\eps}^{(1)}) - F(u_{\eps}^{(2)}))$ are contained in $A^{\eps}$, that is where we gained the estimate of the norm. Moreover, since $$
\supp(u_{\eps}^{(1)}) \subset \Omega_{n_{\eps}(i_{\eps}+1)}^{\eps} \subset \Omega_{\eps n_{\eps}(n_{\eps}+1)}/\eps \subset \Omega_{2h_0}/\eps $$
we have, by definition of $Q_{\eps}$, that $Q_{\eps}(u_{\eps}^{(1)})=0$, $Q_{\eps}(u_{\eps}) = Q_{\eps}(u_{\eps}^{(2)})$ and
\begin{equation}\label{eq_relaz_der_Qeps}
Q'_{\eps}(u_{\eps}^{(1)})=0, \quad Q'_{\eps}(u_{\eps}) = Q'_{\eps}(u_{\eps}^{(2)}).
\end{equation}

\smallskip

\noindent 
\textbf{Step 3.} \emph{Relations of the functionals.}
\\
We show that 
$$\abs{I_{\eps}(u_{\eps}) - I_{\eps}(u_{\eps}^{(1)}) - I_{\eps}(u_{\eps}^{(2)})}\to 0 \quad \hbox{ as $\eps \to 0$}$$ 
from which
\begin{equation}\label{eq_relaz_funzionali}
J_{\eps}(u_{\eps}) = I_{\eps}(u_{\eps}^{(1)}) + I_{\eps}(u_{\eps}^{(2)}) + Q_{\eps}(u_{\eps}^{(2)}) + o(1).
\end{equation}
Indeed 
%\begin{eqnarray*}
%\lefteqn{\abs{I_{\eps}(u_{\eps}) - I_{\eps}(u_{\eps}^{(1)}) - I_{\eps}(u_{\eps}^{(2)})}
% \leq \Bigg|\int_{ \R^N} (-\Delta)^{s/2}u_{\eps}^{(1)} (-\Delta)^{s/2}u_{\eps}^{(2)}dx \Bigg| +}\\
% &&+ \int_{ A^{\eps}} \abs{V(\eps x)u_{\eps}^{(1)} u_{\eps}^{(2)}}dx +\int_{ A^{\eps}} \abs{F(u_{\eps}) -F(u_{\eps}^{(1)}) - F(u_{\eps}^{(2)}) }dx \\
%&=:& (I)+(II)+(III).
%\end{eqnarray*}
\begin{align*}
\abs{I_{\eps}(u_{\eps}) - I_{\eps}(u_{\eps}^{(1)}) - I_{\eps}(u_{\eps}^{(2)})}
 \leq& \; \Bigg|\int_{ \R^N} (-\Delta)^{s/2}u_{\eps}^{(1)} (-\Delta)^{s/2}u_{\eps}^{(2)} \Bigg| + \int_{ A^{\eps}} \abs{V(\eps x)u_{\eps}^{(1)} u_{\eps}^{(2)}} \\
&+\int_{ A^{\eps}} \abs{F(u_{\eps}) -F(u_{\eps}^{(1)}) - F(u_{\eps}^{(2)}) } \\
=:& \, (I)+(II)+(III).
\end{align*}

The second piece can be easily estimated by the boundedness of $\varphi_{\eps}$ and $V$, and the information on the $L^2$-norm given by \eqref{eq_stima_anello}, i.e $(II) \leq \frac{C}{n_{\eps}}$. 
Similarly, as regards $(III)$, we estimate each single piece separately, in the same way: use \eqref{eq_prop_f} and the information on the $L^2$-norm and the $L^{p+1}$-norm given by \eqref{eq_stima_anello}, obtaining $(III)\leq \frac{C}{n_{\eps}}$. 

%\vskip-0.5cm
Focus instead on (I). Recall that $(u_{\eps})_{\eps}\subset S(r_2')$, and thus $\norm{u_{\eps}}_2$ is bounded. We have 
\setlength{\abovedisplayshortskip}{0pt}
%\vskip-2cm
\begin{align*}
(I) &\leq C \int_{\R^{2N}} \frac{\abs{u_{\eps}^{(1)}(x)-u_{\eps}^{(1)}(y)} \abs{u_{\eps}^{(2)}(x)-u_{\eps}^{(2)}(y)}}{|x-y|^{N+2s}} dx \, dy \\
&\leq 2C \int_{\Omega^{\eps}_{n_{\eps} i_{\eps}}\times \complement(\Omega^{\eps}_{n_{\eps}(i_{\eps}+1)})} \frac{\abs{u_{\eps}^{(1)}(x)-u_{\eps}^{(1)}(y)} \abs{u_{\eps}^{(2)}(x)-u_{\eps}^{(2)}(y)}}{|x-y|^{N+2s}} dx \, dy+ \\
& \quad+ 2C \int_{A^{\eps}\times \R^N} \frac{\abs{u_{\eps}^{(1)}(x)-u_{\eps}^{(1)}(y)} \abs{u_{\eps}^{(2)}(x)-u_{\eps}^{(2)}(y)}}{|x-y|^{N+2s}} dx \, dy \\
&=: 2C \big((I_1) + (I_2)\big)
\end{align*}
%\vskip-2cm
since on $\Omega^{\eps}_{n_{\eps} i_{\eps}} \times \Omega^{\eps}_{n_{\eps} i_{\eps}}$ and $\complement(\Omega^{\eps}_{n_{\eps}(i_{\eps}+1)}) \times \complement(\Omega^{\eps}_{n_{\eps}(i_{\eps}+1)})$ the integrand is null. Focusing on $(I_1)$
%\vskip-5cm
\begin{align*}
(I_1) &=\int_{(\Omega^{\eps}_{n_{\eps} i_{\eps}}\times \complement(\Omega^{\eps}_{n_{\eps}(i_{\eps}+1)}))\cap \{|x-y|>n_{\eps}\}} \frac{\abs{u_{\eps}(x)u_{\eps}(y)}}{|x-y|^{N+2s}} dx \, dy\\
&\leq \frac{1}{2}\int_{(\Omega^{\eps}_{n_{\eps} i_{\eps}}\times \complement(\Omega^{\eps}_{n_{\eps}(i_{\eps}+1)}))\cap \{|x-y|>n_{\eps}\}} \frac{u^2_{\eps}(x) + u^2_{\eps}(y)}{|x-y|^{N+2s}}dx \, dy \\
&= \frac{1}{2} \int_{\Omega^{\eps}_{n_{\eps} i_{\eps}}}u^2_{\eps}(x) \int_{\complement(\Omega^{\eps}_{n_{\eps}(i_{\eps}+1)})\cap \{|x-y|>n_{\eps}\}} \frac{1}{|x-y|^{N+2s}} dy \, dx+\\
& \quad+ \frac{1}{2} \int_{\complement(\Omega^{\eps}_{n_{\eps}(i_{\eps}+1)})}u^2_{\eps}(y) \int_{\Omega^{\eps}_{n_{\eps} i_{\eps}}\cap \{|x-y|>n_{\eps}\}} \frac{1}{|x-y|^{N+2s}} dx \, dy \\
&\leq C\norm{u_{\eps}}_2^2 \int_{|x-y|>n_{\eps}} \frac{1}{|x-y|^{N+2s}} dx \, dy 
\leq \frac{C}{{n_{\eps}}^{2s}} \to 0 \quad \hbox{ as $\eps \to 0$}.
\end{align*}
Focusing on $(I_2)$ we have

\begin{eqnarray*}
\lefteqn{(I_2) }\\
&\leq& \int_{A^{\eps}\times \R^N} \frac{1}{|x-y|^{N+2s}} \Big( \abs{(\varphi_{\eps}(x) -\varphi_{\eps}(y))u_{\eps}(x)} + \abs{ \varphi_{\eps}(y)(u_{\eps}(x)-u_{\eps}(y))} \Big) \cdot \\
&& \cdot \Big( \abs{(\varphi_{\eps}(y)-\varphi_{\eps}(x))u_{\eps}(x)} + \abs{(1-\varphi_{\eps}(y))(u_{\eps}(x) -u_{\eps}(y))} \Big) dx \, dy \\
&\leq& \int_{A^{\eps}\times \R^N} \frac{\abs{\varphi_{\eps}(x) -\varphi_{\eps}(y)}^2 \abs{u_{\eps}(x)}^2}{|x-y|^{N+2s}} dx \, dy
 + \int_{A^{\eps}\times \R^N} \frac{\abs{u_{\eps}(x) -u_{\eps}(y)}^2}{|x-y|^{N+2s}}dx \, dy+\\
&&+ 2\int_{A^{\eps}\times \R^N} \frac{\abs{\varphi_{\eps}(x) -\varphi_{\eps}(y)} \abs{u_{\eps}(x)} \abs{u_{\eps}(x) -u_{\eps}(y)}}{|x-y|^{N+2s}}dx \, dy\\
&\leq& \int_{A^{\eps}\times \R^N} \frac{\abs{\varphi_{\eps}(x) -\varphi_{\eps}(y)}^2 \abs{u_{\eps}(x)}^2}{|x-y|^{N+2s}} %dx \, dy
+\int_{A^{\eps}\times \R^N} \frac{\abs{u_{\eps}(x) -u_{\eps}(y)}^2}{|x-y|^{N+2s}} %dx \, dy
+\\
&& +2 \Bigg( \int_{A^{\eps}\times \R^N} \frac{\abs{\varphi_{\eps}(x) -\varphi_{\eps}(y)}^2 \abs{u_{\eps}(x)}^2}{|x-y|^{N+2s}}% dx \, dy 
\Bigg)^{\frac{1}{2}} \Bigg( \int_{A^{\eps}\times \R^N} \frac{\abs{u_{\eps}(x) -u_{\eps}(y)}^2}{|x-y|^{N+2s}} %dx \, dy 
\Bigg)^{\frac{1}{2}} \\
&=:& A^2 + B^2 + 2AB
\end{eqnarray*}
and we see that both $A$ and $B$ go to zero: $B=[u_{\eps}]_{A^{\eps}, \R^N} \leq \frac{C}{n_{\eps}^{1/2}}$ by \eqref{eq_stima_anello}, while for $A$ we exploit that $\norm{\nabla \varphi_{\eps}}_{\infty}\to 0$. 
Indeed, let $\alpha_{\eps}:=\frac{1}{\norm{\nabla \varphi_{\eps}}_{\infty}}$; we have
\begin{align*}
A^2 &= \int_{A^{\eps}} \abs{u_{\eps}(x)}^2 \int_{|x-y|\leq \alpha_{\eps}} \frac{\abs{\varphi_{\eps}(x) -\varphi_{\eps}(y)}^2}{|x-y|^{N+2s}} dy \, dx+\\
 & \quad + \int_{A^{\eps}} \abs{u_{\eps}(x)}^2 \int_{|x-y|> \alpha_{\eps}} \frac{\abs{\varphi_{\eps}(x) -\varphi_{\eps}(y)}^2}{|x-y|^{N+2s}} dy \, dx\\
&\leq C \norm{u_{\eps}}_{L^2(A^{\eps})}^2 \Bigg(\norm{\nabla \varphi_{\eps}}_{\infty}^2 \int_{|z|\leq \alpha_{\eps}} \frac{1}{|z|^{N+2s-2}} dz + 4 \norm{\varphi_{\eps}}_{\infty}^2 \int_{|z|> \alpha_{\eps}} \frac{1}{|z|^{N+2s}}dz \Bigg)\\
&\leq \frac{C}{n_{\eps}} \alpha_{\eps}^{-2s} \left(\norm{\nabla \varphi_{\eps}}_{\infty}^2 \alpha_{\eps}^2 + 1\right)
= \frac{C}{n_{\eps}} \norm{\nabla \varphi_{\eps}}_{\infty}^{2s}\leq \frac{C}{n_{\eps}^{2s+1}}\to 0 \quad \hbox{ as $\eps \to 0$}.
\end{align*}
Thus $(I_2) \leq \frac{C}{n_{\eps}^{2s+1}} + \frac{C}{n_{\eps}^{s+1}} + \frac{C}{n_{\eps}}\leq \frac{C'}{n_{\eps}}\to 0$, 
which reaches the claim. 

\smallskip

\noindent 
\textbf{Step 4.} \emph{Relations of the derivatives.}
\\
We have
\begin{equation}\label{eq_conv_I'}
\norm{I'_{\eps}(u_{\eps}) - I'_{\eps}(u_{\eps}^{(1)}) - I'_{\eps}(u_{\eps}^{(2)})}_{(H^s(\R^N))^*}\to 0 \quad \hbox{ as $\eps \to 0$},
\end{equation}
from which, joined to \eqref{eq_relaz_der_Qeps},
\begin{equation}\label{eq_relaz_derivate}
J'_{\eps}(u_{\eps} ) = I'_{\eps}(u_{\eps}^{(1)}) + I'_{\eps}(u_{\eps}^{(2)}) + Q'_{\eps}(u_{\eps}^{(2)}) + o(1).
\end{equation}
Indeed by H\"older inequality, for any $v\in H^s(\R^N)$,
$$ \abs{I'_{\eps}(u_{\eps})v - I'_{\eps}(u_{\eps}^{(1)})v - I'_{\eps}(u_{\eps}^{(2)})v} \leq \int_{A^{\eps}} \abs{f(u_{\eps})-f(u_{\eps}^{(1)}) - f(u_{\eps}^{(2)})} |v| dx $$
and again we argue in the same way as in the third piece of Step 3, observing that, by \eqref{eq_prop_f}, $|f(u_{\eps})| |v| \leq \beta |u_{\eps}| |v| + C_{\beta} |u_{\eps}|^{p} |v|$ thus
\begin{align*} 
\int_{A^{\eps}} |f(u_{\eps})| |v| dx &\leq \beta \norm{u_{\eps}}_{L^2(A^{\eps})} \norm{v}_2 + C_{\beta}\norm{u_{\eps}}_{L^{p+1}(A^{\eps})}^p \norm{v}_{p+1} \\
&\leq C\left( \beta \norm{u_{\eps}}_{L^2(A^{\eps})} + C_{\beta}\norm{u_{\eps}}_{L^{p+1}(A^{\eps})}^p \right) \norm{v}_{H^s}
\end{align*}
and hence the claim. In particular, $\abs{\big(I'_{\eps}(u_{\eps}) - I'_{\eps}(u_{\eps}^{(1)}) - I'_{\eps}(u_{\eps}^{(2)})\big)u_{\eps}^{(2)}} \leq \frac{C}{n_{\eps}}$.
We see also that 
\begin{equation}\label{eq_o_picc}
 I'_{\eps}(u_{\eps}^{(1)})u_{\eps}^{(2)} = o(1).
\end{equation}
Indeed
\begin{eqnarray*}
\lefteqn{\abs{I'_{\eps}(u_{\eps}^{(1)})u_{\eps}^{(2)}} \leq \Bigg|\int_{ \R^N} (-\Delta)^{s/2}u_{\eps}^{(1)} (-\Delta)^{s/2}u_{\eps}^{(2)} dx\Bigg|+ } \\
&&+ \int_{ A^{\eps}} \abs{V(\eps x)u_{\eps}^{(1)} u_{\eps}^{(2)}dx } + \int_{A^{\eps}} \abs{f(u_{\eps}^{(1)})u_{\eps}^{(2)}dx }
=: (I) + (II) + (III)
\end{eqnarray*}
where for $(I)$ and $(II)$ we argue as in Step 3 obtaining $(I)+ (II) \leq \frac{C}{n_{\eps}^{2s}} + \frac{C}{n_{\eps}}$, while for $(III)$ we argue as in \eqref{eq_relaz_derivate} obtaining $(III) \leq \frac{C}{n_{\eps}} $.

\smallskip

\noindent 
\textbf{Step 5.} \emph{Convergence of $u_{\eps}^{(2)}$.}
\\
Observing that the support of $u_{\eps}^{(2)}$ is outside $\Omega/\eps$, we have with arguments similar to Step 3 that, by \eqref{eq_stima_princ_ueps},
\begin{equation}\label{eq_dim_stima_r2}
\norm{u_{\eps}^{(2)}}_{H^s(\R^N)}\leq r_1.
\end{equation}
Indeed, focusing only on the nonlocal part, we have (recall that $\supp\big(u_{\eps}^{(2)}\big)\subset \complement(\Omega/\eps)$)
\begin{align*}
\int_{\R^{2N}} \frac{\abs{u_{\eps}^{(2)}(x)-u_{\eps}^{(2)}(y)}^2}{|x-y|^{N+2s}}
&\leq 2 \int_{\complement(\Omega/\eps)\times \R^N} \frac{\abs{u_{\eps}^{(2)}(x)-u_{\eps}^{(2)}(y)}^2}{|x-y|^{N+2s}} dx \, dy \\
&\leq 4 \int_{\complement(\Omega/\eps)\times \R^{N}} \frac{\abs{\varphi_{\eps}(x)-\varphi_{\eps}(y)}^2\abs{u_{\eps}(x)}^2}{|x-y|^{N+2s}}dx \, dy+\\
& \quad + 4 \int_{\complement(\Omega/\eps)\times \R^N} \frac{\abs{u_{\eps}(x)-u_{\eps}(y)}^2}{|x-y|^{N+2s}}dx \, dy
\end{align*}
and we use again the final argument in Step 3 and \eqref{eq_stima_princ_ueps} to gain, for $\eps $ small,
$$\norm{u_{\eps}^{(2)}}_{H^s(\R^N)}^2 \leq(C+ o(1)) \norm{u_{\eps}}_{L^2(\complement(\Omega/\eps))}^2 + C [u_{\eps}]_{\complement(\Omega/\eps), \R^N}^2 \leq (Cr_2')^2$$
where $C$ does not depend on $r_2'$. We choose thus $r_2'$ such that \eqref{eq_dim_stima_r2} holds.

This allows us to use Lemma \ref{lemma_stima_I}. 
By joining \eqref{eq_PS2}, \eqref{eq_relaz_derivate}, \eqref{eq_o_picc}, \eqref{eq_stima_I}, \eqref{eq_stima_Q} we obtain 
\begin{align*}
o(1) &= J'_{\eps}(u_{\eps} )u_{\eps}^{(2)} = I'_{\eps}(u_{\eps}^{(1)}) u_{\eps}^{(2)} + I'_{\eps}(u_{\eps}^{(2)})u_{\eps}^{(2)} + Q'_{\eps}(u_{\eps}^{(2)})u_{\eps}^{(2)} + o(1) \\
&\geq C \norm{u_{\eps}^{(2)}}_{H^s(\R^N)}^2 + (p+1)Q_{\eps}(u_{\eps}^{(2)}) +o(1)
\end{align*}
or more precisely (we highlight this for a later use), for some $C=C(r_2')$,
\begin{equation}
 o(1) = J'_{\eps}(u_{\eps} )u_{\eps}^{(2)} \geq C \norm{u_{\eps}^{(2)}}_{H^s(\R^N)}^2 + (p+1)Q_{\eps}(u_{\eps}^{(2)}) -
 \left( \frac{C}{n_{\eps}}+\frac{C}{n_{\eps}^{2s}}\right), \label{eq_dim_stima_mista}
\end{equation}
which implies (since $Q_{\eps}$ is positive) that
\begin{equation}\label{eq_conv_u2}
\norm{u_{\eps}^{(2)}}_{H^s(\R^N)}\to 0 \quad \hbox{ as $\eps\to 0$}
\end{equation}
and
\begin{equation}\label{eq_conv_Q2}
Q_{\eps}(u_{\eps}^{(2)})\to 0 \quad \hbox{ as $\eps\to 0$}.
\end{equation}
As a further consequence, \eqref{eq_conv_u2} and the boundedness of $V$ imply 
\begin{equation}\label{eq_conv_I2}
I_{\eps}(u_{\eps}^{(2)})\to 0 \quad \hbox{and} \quad I'_{\eps}(u_{\eps}^{(2)})\to 0 \quad \hbox{ as $\eps \to 0$}.
\end{equation}

\smallskip

\noindent 
\textbf{Step 6.} \emph{Convergence of $I_{\eps}'(u_{\eps}^{(1)})$.}
\\
In particular we obtain from \eqref{eq_conv_I2}, together with \eqref{eq_relaz_funzionali} and \eqref{eq_conv_Q2}, that
\begin{equation}\label{eq_confr_JI}
J_{\eps}(u_{\eps}) = I_{\eps}(u_{\eps}^{(1)}) + o(1).
\end{equation}
We want now to show that
\begin{equation}\label{eq_conv_I1}
I'_{\eps}(u_{\eps}^{(1)})\to 0 \quad \hbox{in $(H^s(\R^N))^*$} \quad \hbox{ as $\eps \to 0$}.
\end{equation}
Start observing that \eqref{eq_conv_I2} together with \eqref{eq_conv_I'} give
\begin{equation}\label{eq_info_I1}
I'_{\eps}(u_{\eps}) = I'_{\eps}(u_{\eps}^{(1)}) + o(1);
\end{equation}
let now $v \in H^s(\R^N)$ and evaluate $I'_{\eps}(u_{\eps}^{(1)})v$. We want to exploit \eqref{eq_info_I1} together again with the assumption \eqref{eq_PS2}. In order to do this we need to pass from $u_{\eps}^{(2)}$ to $u_{\eps}$, but getting rid of $Q'_{\eps}(u_{\eps})$ on which we have no information. Thus we introduce a cutoff function $\tilde{\varphi}\in C^{\infty}_c(\R^N)$ such that
$$\Omega_{\frac{3}{2}h_0} \prec \tilde{\varphi} \prec \Omega_{2h_0}$$
and hence
$$\supp(u_{\eps}^{(1)}) \subset \Omega_{h_0}/\eps\subset \Omega_{\frac{3}{2}h_0}/\eps\subset \{\tilde{\varphi}(\eps \cdot)\equiv 1\}\subset \supp(\tilde{\varphi}(\eps \cdot)) \subset \Omega_{2h_0}/\eps.$$
Thus we have
%\begin{eqnarray*}
%\lefteqn{ I'_{\eps}(u_{\eps}^{(1)})v \stackrel{(*)}= I'_{\eps}(u_{\eps}^{(1)})(\tilde{\varphi}(\eps \cdot) v) +(1+\norm{v}_2) o(1)} \\
%&=& I'_{\eps}(u_{\eps})(\tilde{\varphi}(\eps \cdot) v) - \left(I'_{\eps}(u_{\eps})(\tilde{\varphi}(\eps \cdot) v) - I'_{\eps}(u_{\eps}^{(1)})(\tilde{\varphi}(\eps \cdot) v) \right) +(1+\norm{v}_2) o(1) \\
%&=& J'_{\eps}(u_{\eps})(\tilde{\varphi}(\eps \cdot) v) - \left(I'_{\eps}(u_{\eps}) - I'_{\eps}(u_{\eps}^{(1)}) \right)(\tilde{\varphi}(\eps \cdot) v) +(1+\norm{v}_2) o(1).
%\end{eqnarray*}
\begin{align*}
 I'_{\eps}(u_{\eps}^{(1)})v &\stackrel{(*)}= I'_{\eps}(u_{\eps}^{(1)})(\tilde{\varphi}(\eps \cdot) v) +(1+\norm{v}_2) o(1) \\
&= I'_{\eps}(u_{\eps})(\tilde{\varphi}(\eps \cdot) v) - \left(I'_{\eps}(u_{\eps})(\tilde{\varphi}(\eps \cdot) v) - I'_{\eps}(u_{\eps}^{(1)})(\tilde{\varphi}(\eps \cdot) v) \right) +(1+\norm{v}_2) o(1) \\
&= J'_{\eps}(u_{\eps})(\tilde{\varphi}(\eps \cdot) v) - \left(I'_{\eps}(u_{\eps}) - I'_{\eps}(u_{\eps}^{(1)}) \right)(\tilde{\varphi}(\eps \cdot) v) +(1+\norm{v}_2) o(1).
\end{align*}
Indeed, we justify $(*)$ as done in Step 3: notice that $u_{\eps}^{(1)}\equiv 0$ outside $\Omega_{h_0}/\eps$ and $1-\tilde{\varphi}(\eps \cdot)\equiv 0$ in $\Omega_{\frac{3}{2}h_0}/\eps$, so in the annulus $(\Omega_{\frac{3}{2} h_0}/\eps) \setminus (\Omega_{h_0}/\eps)$ both $u_{\eps}^{(1)}$ and $1-\tilde{\varphi}(\eps \cdot)$ are zero; notice also that $\complement(\Omega_{\frac{3}{2} h_0}/\eps)$ and $\Omega_{h_0}/\eps$ get far one from the other as $\eps \to 0$. Thus we have
%\begin{eqnarray*}
%\lefteqn{\abs{I'_{\eps}(u_{\eps}^{(1)})((1-\tilde{\varphi}(\eps \cdot))v)} =2 \int_{(\Omega_{h_0}/\eps)\times \complement(\Omega_{\frac{3}{2}h_0}/\eps)} \frac{\abs{u_{\eps}^{(1)}(x)} \abs{(1-\tilde{\varphi}(\eps y))v(y)}}{\abs{x-y}^{N+2s}} dx \, dy} \\
%&&\leq 2\int_{(\Omega_{h_0}/\eps)\times \complement(\Omega_{\frac{3}{2}h_0}/\eps)} \frac{\abs{u_{\eps}(x)} \abs{v(y)}}{\abs{x-y}^{N+2s}}dx \, dy 
%\leq (1+\norm{v}_2) o(1) 
%\end{eqnarray*}
\begin{align*}
\abs{I'_{\eps}(u_{\eps}^{(1)})((1-\tilde{\varphi}(\eps \cdot))v)} &= 2 \int_{(\Omega_{h_0}/\eps)\times \complement(\Omega_{\frac{3}{2}h_0}/\eps)} \frac{\abs{u_{\eps}^{(1)}(x)} \abs{(1-\tilde{\varphi}(\eps y))v(y)}}{\abs{x-y}^{N+2s}} dx \, dy \\
&\leq 2\int_{(\Omega_{h_0}/\eps)\times \complement(\Omega_{\frac{3}{2}h_0}/\eps)} \frac{\abs{u_{\eps}(x)} \abs{v(y)}}{\abs{x-y}^{N+2s}}dx \, dy 
\leq (1+\norm{v}_2) o(1) 
\end{align*}
where in the last passage we argue as for $(I_1)$ in Step 3. Notice that $o(1)$ does not depend on $v$. 
Thus we obtain
\begin{eqnarray*}
\lefteqn{\abs{I'_{\eps}(u_{\eps}^{(1)})v}} \\
 &\leq& \left( \norm{J'_{\eps}(u_{\eps})}_{(H^s(\R^N))^*} +\norm{I'_{\eps}(u_{\eps}) - I'_{\eps}(u_{\eps}^{(1)})}_{(H^s(\R^N))^*}\right) \norm{\tilde{\varphi}(\eps \cdot)v}_{H^s(\R^N)} \\
 &&+ (1+\norm{v}_{H^s(\R^N)}) o(1) \\
 &\leq& \left( \norm{J'_{\eps}(u_{\eps})}_{(H^s(\R^N))^*} + \norm{I'_{\eps}(u_{\eps}) - I'_{\eps}(u_{\eps}^{(1)})}_{(H^s(\R^N))^*}\right) \norm{v}_{H^s(\R^N)} \big(C +o(1)\big)+\\
 &&+(1+\norm{v}_{H^s(\R^N)}) o(1)
\end{eqnarray*}
where in the last inequality we argue as in Step 5 (again $o(1)$ does not depend on $v$). 
Concluding, we have, by choosing $\norm{v}_{H^s(\R^N)}=1$, that
\begin{eqnarray*}
\lefteqn{\norm{I'_{\eps}(u_{\eps}^{(1)})}_{(H^s(\R^N))^*}}\\
&& \leq C \left( \norm{J'_{\eps}(u_{\eps})}_{(H^s(\R^N))^*} + \norm{I'_{\eps}(u_{\eps}) - I'_{\eps}(u_{\eps}^{(1)})}_{(H^s(\R^N))^*}\right) \big(1 +o(1)\big) + o(1) \to 0
\end{eqnarray*}
by using \eqref{eq_PS2} and \eqref{eq_info_I1}.

\smallskip

\noindent 
\textbf{Step 7.} \emph{Weak convergence of $u_{\eps}^{(1)}$.}
\\
Set $q_{\eps}:=\Upsilon(u_{\eps})$ to avoid cumbersome notation. Since $(\eps q_{\eps})_{\eps} \subset \Omega[0, \nu_0]\subset \Omega$ bounded in $\R^N$, we have that up to a subsequence
$$\eps q_{\eps} \to p_0 \in \overline{\Omega[0, \nu_0]} \subset K_d
\subset \Omega.$$
Moreover, by estimating the norm of $u_{\eps}^{(1)}$ with the norm of $u_{\eps}$ (as done before, in Step 5, for $u_{\eps}^{(2)}$), where $u_{\eps}$ belongs to $S(r_2')$ bounded in $H^s(\R^N)$, we have that also $u_{\eps}^{(1)}$ is a bounded sequence, and thus is so $u_{\eps}^{(1)}(\cdot + q_{\eps})$, which implies, up to a subsequence
$$u_{\eps}^{(1)}(\cdot + q_{\eps}) \wto \tilde{U} \; \; \hbox{in $H^s(\R^N)$} \quad \hbox{ as $\eps \to 0$}.$$
For each $v\in %\tr{C^{\infty}_c(\R^N)}$ %
H^s(\R^N)$ 
we apply this weak convergence to the following equalities, derived from \eqref{eq_conv_I1},
\begin{align*}
o(1) &= I'_{\eps}(u_{\eps}^{(1)}) v(\cdot-q_{\eps}) = \int_{\R^N} (-\Delta)^{s/2} u_{\eps}^{(1)}(y+q_{\eps}) (-\Delta)^{s/2}v(y) dy + \\
&\quad+\int_{\R^N} V(\eps y + \eps q_{\eps}) u_{\eps}^{(1)}(y+q_{\eps}) v(y) dy - \int_{\R^N} f(u_{\eps}^{(1)}(y+q_{\eps})) v(y) dy\\
&=: (I) + (II) + (III).
\end{align*}
For $(I)$ and $(III)$ we obtain by the weak convergence and by %classical arguments \tor{non posso usare direttamente 
Proposition \ref{prop_converg_generiche_loc} %, né decadimenti polinomiali uniformi, \tr{cosa uso?}} %Non serve sottocritico strettamente!
$$(I)\to \int_{\R^N} (-\Delta)^{s/2} \tilde{U} \,(-\Delta)^{s/2}v \, dy \quad \hbox{ and } \quad (III)\to - \int_{\R^N} f(\tilde{U}) v \, dy \quad \hbox{ as $\eps \to 0$}.$$
For $(II)$ instead we have
%\begin{eqnarray*}
%\lefteqn{\abs{(II) -\int_{\R^N} V(p_0) \tilde{U}v \, dy}}\\ 
%&\leq& \Bigg(\int_{\R^N} (V(\eps y + \eps q_{\eps}) -V(p_0))^2 v^2(y) dy\Bigg)^{1/2} \norm{u_{\eps}^{(1)}(\cdot+q_{\eps})}_2+ \\
%&&+V(p_0)\Bigg|\int_{\R^N} u_{\eps}^{(1)}(y+q_{\eps}) v(y) dy - \int_{\R^N} \tilde{U}v \, dy \Bigg|
%\end{eqnarray*}
\begin{align*}
\pabs{(II) -\int_{\R^N} V(p_0) \tilde{U}v % \, dy
}
\leq& \, \Bigg(\int_{\R^N} (V(\eps y + \eps q_{\eps}) -V(p_0))^2 v^2(y) %dy
\Bigg)^{1/2} \norm{u_{\eps}^{(1)}(\cdot+q_{\eps})}_2+ \\
&+V(p_0)\Bigg|\int_{\R^N} u_{\eps}^{(1)}(y+q_{\eps}) v(y) % dy
 - \int_{\R^N} \tilde{U}v % \, dy 
\Bigg|
\end{align*}
where the first term goes to zero (thanks to the boundedness of $u_{\eps}^{(1)}$) by the dominated convergence theorem, while the second thanks to the weak convergence. Thus we finally obtain
$$L'_{V(p_0)}(\tilde{U})v = \int_{\R^N} (-\Delta)^{s/2} \tilde{U} \, (-\Delta)^{s/2}v \, dy+\int_{\R^N} V(p_0) \tilde{U}v \, dy- \int_{\R^N} f(\tilde{U}) v \, dy=0$$
for each %\tr{$v \in C^{\infty}_c(\R^N)$, and thus by density for each} 
$v\in H^s(\R^N)$, that is
\begin{equation}\label{eq_critic_Vp0}
L'_{V(p_0)}(\tilde{U})=0.
\end{equation}

\smallskip

\noindent 
\textbf{Step 8.} \emph{Strong convergence of $u_{\eps}^{(1)}$.}
\\
We want to show the strong convergence of $u_{\eps}^{(1)}(\cdot+q_{\eps})$, that is
\begin{equation}\label{eq_dim_strong_conv_tildeU}
u_{\eps}^{(1)}(\cdot + q_{\eps}) \to \tilde{U} \quad \hbox{ in $H^s(\R^N)$} \quad \hbox{ as $\eps \to 0$}.
\end{equation}
Set $\tilde{w}_{\eps}:= u_{\eps}^{(1)}(\cdot+q_{\eps}) -\tilde{U} \wto 0$, again by \eqref{eq_conv_I1} we have
\begin{align*}
o(1) &= I'_{\eps}(u_{\eps}^{(1)})\tilde{w}_{\eps}(\cdot-q_{\eps}) \\
&= L'_{V(p_0)}(\tilde{U})\tilde{w}_{\eps} +\left(\norm{(-\Delta)^{s/2} \tilde{w}_{\eps}}_2^2 + \int_{\R^N} V(\eps y + \eps q_{\eps}) \tilde{w}_{\eps}^2 dy\right)+\\
&\quad +\int_{\R^N} (V(\eps y + \eps q_{\eps})-V(p_0))\tilde{U}\tilde{w}_{\eps}dy + \int_{\R^N} (f(\tilde{U})-f(\tilde{U}+\tilde{w}_{\eps}))\tilde{w}_{\eps} dy\\
&=: \Bigg(\norm{(-\Delta)^{s/2} \tilde{w}_{\eps}}_2^2 + \int_{\R^N} V(\eps y + \eps q_{\eps}) \tilde{w}_{\eps}^2 dy\Bigg) + (I) \\
&\geq \norm{(-\Delta)^{s/2} \tilde{w}_{\eps}}_2^2 +\underline{V} \norm{\tilde{w}_{\eps}}_2^2 + (I)
\end{align*}
where we have used \eqref{eq_critic_Vp0}. We obtain by the boundedness of $V$ and \eqref{eq_prop_f} %Proposition \ref{prop_converg_generiche_loc}
\begin{align*}
(I) \geq& - 2 \norm{V}_{\infty} \int_{\R^N}|\tilde{U}| |\tilde{w}_{\eps}|dy-\int_{\R^N} \left(2\beta |\tilde{U|} + C_{\beta}(2^p+1)|\tilde{U}|^p\right)|\tilde{w}_{\eps}| dy -\\
&- \int_{\R^N} \left(\beta |\tilde{w}_{\eps}|^2 + 2^p C_{\beta} |\tilde{w}_{\eps}|^{p+1}\right) dy
= o(1) - \beta \norm{\tilde{w}_{\eps}}_2^2 - 2^p C_{\beta} \norm{\tilde{w}_{\eps}}_{p+1}^{p+1};
\end{align*}
in the last passage we have used that $\tilde{w}_{\eps}\wto 0$ in $H^s(\R^N)$, thus by Remark \ref{rem_conv_deb_assol} %classical properties
$|\tilde{w}_{\eps}|\wto 0$ in $H^s(\R^N)$ and hence in $L^2(\R^N)$ and in $L^{p+1}(\R^N)$ (observing that $\tilde{U}^p \in L^{1+\frac{1}{p}}(\R^N)$).
Merging together all the things we have, by \eqref{eq_stima_2p} and choosing $\beta < \frac{1}{2} \underline{V}$, %va bene anche per $p+1=2^*_s$
$$o(1) \geq \norm{(-\Delta)^{s/2} \tilde{w}_{\eps}}_2^2 +(\underline{V}-\beta) \norm{\tilde{w}_{\eps}}_2^2 - 2^p C_{\beta} \norm{\tilde{w}_{\eps}}_{p+1}^{p+1}\geq C \norm{\tilde{w}_{\eps}}_{H^s(\R^N)}^2$$
and thus $\tilde{w}_{\eps} \to 0$ strongly in $H^s(\R^N)$, that is the claim.

\smallskip

\noindent 
\textbf{Step 9.} \emph{Localization.}
\\
Observe first that $\tilde{U} \nequiv 0$. Indeed, if not, by \eqref{eq_conv_u2}, \eqref{eq_dim_strong_conv_tildeU} and translation invariance of the norm we would have 
\begin{align*}
r^*&\leq \liminf_{\eps \to 0} \norm{U_{\eps}}_{H^s(\R^N)} \\
&\leq \liminf_{\eps \to 0} \norm{U_{\eps}(\cdot-p_{\eps})+\varphi_{\eps}}_{H^s(\R^N)} + \liminf_{\eps \to 0} \norm{\varphi_{\eps}}_{H^s(\R^N)} \\
&\leq \lim_{\eps \to 0} \Big( \norm{u_{\eps}^{(1)}}_{H^s(\R^N)} + \norm{u_{\eps}^{(2)}}_{H^s(\R^N)}\Big)+ r_2' \leq r_2' <r^*,
\end{align*}
impossible. 
By \eqref{eq_dim_strong_conv_tildeU} we obtain also
$$I_{\eps}(u_{\eps}^{(1)}) \to L_{V(p_0)}(\tilde{U}) \quad \hbox{ as $\eps \to0$}.$$
Thus we find, by using also \eqref{eq_confr_JI} and \eqref{eq_PS1},
$$L_{V(p_0)}(\tilde{U}) = I_{\eps}(u_{\eps}^{(1)}) + o(1) = J_{\eps}(u_{\eps}) + o(1) \leq l_0'+ o(1)$$
and hence, letting $\eps \to 0$,
\begin{equation}\label{eq_stima_Lvp0}
L_{V(p_0)}(\tilde{U}) \leq l_0'< E_{m_0+\nu_1}.
\end{equation}
Moreover by \eqref{eq_critic_Vp0} and $\tilde{U}\nequiv 0$, we have
$$E_{V(p_0)}\leq L_{V(p_0)}(\tilde{U});$$
joining together the two previous inequalities we find $E_{V(p_0)} < E_{m_0+\nu_1}$ which implies, by the monotonicity of $E_a$, that
$$V(p_0) < m_0+\nu_1.$$
Joining this information to the fact that $p_0 \in \Omega$ (and in particular $V(p_0) \geq m_0$) we have $p_0 \in \Omega[0, \nu_1)$, that is
\begin{equation}\label{eq_local_qeps}
\eps \Upsilon(u_{\eps}) \to p_0 \in \Omega[0,\nu_1) \quad \hbox{ as $\eps \to 0$}.
\end{equation}
Exploiting again \eqref{eq_stima_Lvp0} (observe that $E_{m_0+\nu_1}< E_{m_0+\nu_0} < l_0$) together with $L'_{V(p_0)}(\tilde{U})=0$, and $\tilde{U}\neq 0$, we have that $\tilde{U}$ belongs to $S_{V(p_0)}$ up to translations, that is
$$U:= \tilde{U}(\cdot - y_0) \in S_{V(p_0)}\subset \widehat{S}$$
for some suitable $y_0 \in \R^N$. So, set
$$p_{\eps}:=q_{\eps} + y_0$$
we have
\begin{equation}\label{eq_conv_u1}
\norm{u_{\eps}^{(1)} - U(\cdot - p_{\eps})}_{H^s(\R^N)} \to 0 \quad \hbox{ as $\eps \to 0$}.
\end{equation}
For a later use observe also that
\begin{equation}\label{eq_conv_p_eps}
\eps p_{\eps} \to p_0 \in \Omega[0, \nu_1) \quad \hbox{ as $\eps \to 0$}.
\end{equation}

\smallskip

\noindent 
\textbf{Step 10.} \emph{Conclusions.}
\\
By \eqref{eq_local_qeps} we have that
$$\eps \Upsilon(u_{\eps})\in \Omega[0,\nu_1) $$
definitely for $\eps$ small. This is the first part of the claim. 
Moreover, by \eqref{eq_conv_u1} and \eqref{eq_conv_u2} we gain
\begin{equation}\label{eq_conv_uu}
\norm{u_{\eps} - U(\cdot - p_{\eps})}_{H^s(\R^N)} \to 0 \quad \hbox{ as $\eps \to0$}
\end{equation}
and thus, since $\widehat{\rho}(u_{\eps})\leq \norm{u_{\eps} - U(\cdot - p_{\eps})}_{H^s(\R^N)}$ by definition, also $\widehat{\rho}(u_{\eps})\to 0$ and hence
$$\widehat{\rho}(u_{\eps})\in [0, \rho_1]$$
definitely for $\eps$ small. This concludes the proof.
\QED

\bigskip

In the next proposition we see that solutions of $J'_{\eps}(u)=0$ are, under suitable assumptions, also solutions of $I'_{\eps}(u)=0$.

\begin{Corollary}\label{corol_passaggio_critici}
Let $(u_{\eps})_{\eps}$ be a sequence of critical points of $J_{\eps}$, that is $J_{\eps}'(u_{\eps})=0$,
satisfying 
$$u_{\eps}\in S(r_2'), \quad J_{\eps}(u_{\eps})\leq l_0' \quad \textit{and} \quad \eps \Upsilon(u_{\eps})\in \Omega[0,\nu_0] $$
for any $\eps>0$. 
Then, for $\eps$ sufficiently small, we have
$$Q_{\eps}(u_{\eps})=0, \quad \textit{and} \quad Q'_{\eps}(u_{\eps})=0.$$
In particular $I'_{\eps}(u_{\eps})=0$, 
which means that $u_{\eps}$ is a solution of \eqref{eq_princ_epsx}.
\end{Corollary}

\claim Proof.
By the proof of Lemma \ref{lemma_stima_basso}, we notice, since $1-\varphi_{i_\eps}^{\eps}\equiv 1$ outside $\Omega_{i_{\eps}+1}^{\eps}$, and thus outside $\Omega_{h_0}/{\eps}$, that
\begin{equation}\label{eq_dim_uguagl_norme}
\norm{u_{\eps}}_{L^2(\R^N \setminus (\Omega_{h_0}/\eps))} = \norm{u_{\eps}^{(2)}}_{L^2(\R^N \setminus (\Omega_{h_0}/\eps))} \leq \norm{u_{\eps}^{(2)}}_{H^s(\R^N)}
\end{equation}
and hence $\norm{u_{\eps}}_{L^2(\R^N \setminus (\Omega_{h_0}/\eps))}\to 0$ by \eqref{eq_conv_u2}. 

Through a careful analysis of the Steps 3--5 of the proof, that is by \eqref{eq_dim_stima_mista} and \eqref{eq_dim_uguagl_norme}, we see, more precisely, that
$$\norm{u_{\eps}}^2_{L^2(\R^N \setminus (\Omega_{h_0}/\eps))}\leq \frac{C}{n_{\eps}}+ \frac{C}{n_{\eps}^{2s}} +o(1)$$
where $C=C(r_2')$ and $o(1)$ depends on the rate of convergence of $J'_{\eps}(u_{\eps})$. Thus, since we assume $J'_{\eps}(u_{\eps})\equiv 0$, we gain uniformity, i.e., called $\alpha^*:=\min\{1, 2s\}$, we obtain
\begin{equation}\label{eq_conv_u_a0_unif}
\norm{u_{\eps}}^2_{L^2(\R^N \setminus (\Omega_{h_0}/\eps))} \leq \frac{C}{n_{\eps}^{\alpha^*}} \sim \eps^{\alpha^*/2}
\end{equation}
As a consequence
$$\frac{1}{\eps^{\alpha}}\norm{u_{\eps}}^2_{L^2(\R^N \setminus (\Omega_{h_0}/\eps))} \to 0 \quad \hbox{ as $\eps \to 0$}$$
for $\alpha \in (0, \alpha^*/2)$, and hence $Q_{\eps}(u_{\eps})\equiv Q'_{\eps}(u_{\eps})\equiv 0$ for $\eps$ sufficiently small.
\QED

\bigskip

We want to show now a (truncated) Palais-Smale-like condition.

\begin{Proposition}\label{prop_Palais_Smale}
There exists $r_2'' \in (0, \min\{r_0, r_1\})$ sufficiently small with the following property: let $\eps>0$ fixed and let $(u_j)_{j}\subset S(r_2'')$ be such that 
\begin{equation}\label{eq_PS3}
\norm{J'_{\eps}(u_j)}_{(H^s(\R^N))^*}\to 0 \quad \textit{ as $j\to +\infty$}
\end{equation}
with the additional assumption
$$(\eps\Upsilon(u_j))_j \subset \Omega[0,\nu_0].$$
Then $(u_j)_j$ admits a strongly convergent subsequence in $H^s(\R^N)$.
\end{Proposition}

\medskip

\claim Proof.
Let $r_2''$ to be fixed. Since $S(r_2'')$ is bounded, up to a subsequence we can assume $u_j \wto u_0$ in $H^s(\R^N)$. We want to show that
$$\lim_{R \to +\infty} \lim_{j \to +\infty} \norm{u_j}_{L^q(\R^N \setminus B_R)}=0$$
for $q=2$ and $q=p+1$ and conclude by Lemma \ref{lemma_conv_forte}. 

Arguing similarly to Step 1 of the proof of Lemma \ref{lemma_stima_basso}, i.e. exploiting Remark \ref{remark_triplanorma}, we obtain for $L \gg 0$, uniformly in $j \in \N$,
$$\tnorm{u_j}_{\R^N \setminus B_L} \leq C r_2'';$$
indeed we work with the set $B_L - \Upsilon(u_j)$ which expands to $\R^N$ as $L, j\to +\infty$, since $\Upsilon(u_j) \in \Omega/\eps$, a fixed bounded set. 
Moreover, for any $n \in \N$, we have
$$\sum_{i=1}^n \norm{u_j}^2_{L^2(B_{L+ni}\setminus B_{L+n(i-1)})} \leq (C r_2'')^2$$
and similarly for the Gagliardo seminorm and the $(p+1)$-norm, thus for some $i_{j,n}\in \{1, \dots, n\}$
$$ \norm{u_j}^2_{ A^{j,n}} + \norm{u_j}^{p+1}_{L^{p+1}(A^{j,n})} \leq \frac{C}{n}$$
where $A^{j,n}:=B_{L+ni_{j,n}}\setminus B_{L+n(i_{j,n}-1)}$. Again similarly to Step 2 of the proof of Lemma \ref{lemma_stima_basso}, we introduce $\psi_{j,n}$ such that
$$B_{L + n(i_{j,n} -1)} \prec \psi_{j,n} \prec B_{L+ni_{j,n}}$$ 
and $\norm{\nabla \psi_{j,n}}_{\infty}=o(1)$ as $n\to +\infty$;
moreover we set
$$\tilde{u}_{j,n} :=(1- \psi_{j,n}) u_j.$$
Observe that $\chi_{B_L} \leq \psi_{j,n}$, thus $\supp(\tilde{u}_{j,n}) \subset \complement(B_L)$. 
Arguing as in Step 5 and 3 of the proof of Lemma \ref{lemma_stima_basso} we obtain
\begin{eqnarray*}
\lefteqn{ 
\int_{\R^{2N}} \frac{\abs{\tilde{u}_{j,n}(x)-\tilde{u}_{j,n}(y)}^2}{|x-y|^{N+2s}} \leq 4 \int_{\complement(B_L) \times \R^{N}} \frac{\abs{\psi_{j,n}(x)-\psi_{j,n}(y)}^2\abs{u_j(x)}^2}{|x-y|^{N+2s}} dx \, dy +}\\
 &&+4 \int_{\complement(B_L)\times \R^N} \frac{\abs{u_j(x)-u_j(y)}^2}{|x-y|^{N+2s}} dx \, dy \leq o(1) \norm{u_j}_{L^2(\complement(B_L))}^2 + C [u_j]_{\complement(B_L), \R^N}^2 
\end{eqnarray*}
thus $\norm{\tilde{u}_{j,n}}_{H^s(\R^N)} \leq C r_2''$ and hence, choosing $r_2''$ sufficiently small, we have 
$$\norm{\tilde{u}_{j,n}}_{H^s(\R^N)} \leq r_1;$$
by Lemma \ref{lemma_stima_I}, for $q \in \{2, p+1\}$, we obtain
$$\norm{u_j}_{L^q(\R^N \setminus B_{L+ni_{j,n}})}^2 = \norm{\tilde{u}_{j,n}}_{L^q(\R^N \setminus B_{L+ni_{j,n}})}^2 \leq C \norm{\tilde{u}_{j,n}}_{H^s(\R^N)}^2 \leq C I'_{\eps}(\tilde{u}_{j,n})\tilde{u}_{j,n}.$$
Thus the claim comes if we show that 
$$I'_{\eps}(\tilde{u}_{j,n})\tilde{u}_{j,n} \to 0 \quad \hbox{ as $j,n\to +\infty$}.$$
Indeed we have
\begin{eqnarray*}
\lefteqn{I'_{\eps}(\tilde{u}_{j,n})\tilde{u}_{j,n} = I'_{\eps}(u_j) \tilde{u}_{j,n} - \int_{\R^{2N}} (-\Delta)^{s/2}(\psi_{j,n} u_j) (-\Delta)^{s/2}((1-\psi_{j,n})u_j) dx-}\\
 && -\int_{A^{j,n}} V(\eps x) \psi_{j,n} (1-\psi_{j,n}) u_j^2dx - \int_{A^{j,n}} (f((1-\psi_{j,n})u_j)- f(u_j)) (1-\psi_{j,n})u_j dx\\
&&=: J'_{\eps}(u_j)\tilde{u}_{j,n} - Q'_{\eps}(u_j)\tilde{u}_{j,n} + (I)
\leq o(1) + (I)
\end{eqnarray*}
where we have used that $J'_{\eps}(u_j) \to 0$ (as $j\to +\infty$, uniformly in $n\in \N$), the boundedness of $\norm{\tilde{u}_{j,n}}_{H^s(\R^N)}$ and the positivity of $ Q'_{\eps}(u_j)\tilde{u}_{j,n}$. The term $(I)$ can be estimated in the same way as done in Steps 3-4 of the proof of Lemma \ref{lemma_stima_basso} (fixed $j\in \N$, and $n\to +\infty$), and hence we reach the claim.
%Sto usando Proposition \ref{prop_converg_generiche_loc}? Credo di no
\QED

\subsection{Deformation lemma on a neighborhood of expected solutions}
\label{sec_neighbo_expect}

We want to define now a neighborhood of expected solutions (see \cite{CJT}), which will be invariant under a suitable deformation flow. 
Consider $r_3:= \min\{ r', r_0', r_2', r_2''\}$ (see Lemma \ref{lemma_stima_g}, Lemma \ref{lemma_stimabasso_I}, Theorem \ref{theorem_stima_J'} and Proposition \ref{prop_Palais_Smale}), and let us define
$$R(\delta, u):= \delta - \frac{\delta_2}{2} (\widehat{\rho}(u) - \rho_1)_+ \leq \delta$$
and
$$\mathcal{X}_{\eps, \delta}:= \Big\{ u \in S(\rho_0) \mid \eps \Upsilon(u) \in \Omega[0,\nu_0), \; J_{\eps}(u) < E_{m_0} + R(\delta, u) \Big\}$$
where
$$0< \rho_1 < \rho_0 < r_3,$$
$\eps$ is sufficiently small and 
\begin{equation}\label{eq_propr_delta}
\delta \in \left(0, \min \left\{ \tfrac{\delta_2}{4}(\rho_0- \rho_1), \, \delta_1, \, l_0'-E_{m_0}\right\}\right);
\end{equation}
here $\delta_1$ and $\delta_2$ are the ones that appear in Lemma \ref{lemma_stimabasso_I} and Theorem \ref{theorem_stima_J'}. Notice that the \emph{height} of the sublevel in $\mc{X}_{\eps, \delta}$ depends on $u$ itself; this will be used to gain a deformation which preserves $\mc{X}_{\eps, \delta}$. 

\noindent We begin by pointing out some geometrical features of the neighborhood $\mc{X}_{\eps, \delta}$. 
\begin{itemize}
\item $\mc{X}_{\eps, \delta}$ is open. Indeed, $S(\rho)$ and $\{J_{\eps}(u) < E_{m_0} + R(\delta, u)\}$ are open, and $\Omega[0,\nu_0)= \Omega(-\gamma, \nu_0)$ for a whatever $\gamma>0$ (since $V$ cannot go under $m_0$ in $\Omega$) and thus open. Moreover it is nonempty (see e.g. Section \ref{sec_def_phi_eps}).
\item If $v\in\mathcal{X}_{\eps,\delta} \subset S(\rho_0)$, then by \eqref{eq_stima_widerho} we have 
$\widehat{\rho}(v)<\rho_0.$

\item 
If $v\in \mathcal{X}_{\eps, \delta} \subset \{ \eps \Upsilon(v) \in \Omega[0, \nu_0]\}$, then
\begin{equation}\label{eq_appart_insieme_piccolo}
\eps \Upsilon(v) \in \Omega [0, \nu_1).
\end{equation}
Indeed, if not, i.e. $\eps \Upsilon(v) \in \Omega [\nu_1, \nu_0]$, then by Lemma \ref{lemma_stimabasso_I} we have
$$J_{\eps}(v) \geq I_{\eps}(v) \geq E_{m_0} + \delta_1 > E_{m_0} + \delta \geq E_{m_0} + R(\delta, v)$$
which is an absurd.

\item If $R(\delta, v) \geq -\delta$ then
\begin{equation}\label{eq_cons_conten}
v \in S(\rho_0).
\end{equation}
Indeed $\frac{\delta_2}{4}(\widehat{\rho}(u)-\rho_1)_+ \leq \delta$ implies, by the restriction on $\delta$,
$$(\widehat{\rho}(u)-\rho_1)_+ < \rho_0-\rho_1.$$
If $\widehat{\rho}(u)<\rho_1$ then clearly $u \in S(\rho_1)\subset S(\rho_0)$. If instead $\widehat{\rho}(u)\geq \rho_1$, then 
$\widehat{\rho}(u) < \rho_0$, 
which again implies $u \in S(\rho_0).$
\end{itemize} 
We further define the set of critical points of $J_{\eps}$ lying in the neighborhood of expected solutions
$$K_c:= \left\{ u \in \mathcal{X}_{\eps, \delta} \mid J'_{\eps}(u)=0, \; J_{\eps}(u)=c\right\},$$
the sublevel
 $$\mc{X}_{\eps, \delta}^c:= \mc{X}_{\eps, \delta}\cap J_{\eps}^c$$
 and the strip level
 $$ (\mc{X}_{\eps, \delta})^c_d:= \{ u \in \mc{X}_{\eps, \delta} \mid d \leq J_{\eps}(u) \leq c\}, $$
for every $c,d \in \R$. 
We present now a deformation lemma with respect to $K_c$, for $c$ sufficiently close to $E_{m_0}$. 

\begin{Lemma}\label{lem_def_lem}
Let $c \in (E_{m_0}-\delta, E_{m_0}+\delta)$.
Then there exists a deformation at level $c$, which leaves the set $\mathcal{X}_{\eps, \delta}$ invariant.
That is, for every $U$ neighborhood of $K_c$ ($U=\emptyset$ if $K_c=\emptyset$), there exist a small $\omega>0$ and a continuous deformation $\eta: [0,1]\times \mc{X}_{\eps, \delta} \to \mc{X}_{\eps,\delta}$ such that
\begin{itemize}
\item[(i)] \ $\eta(0, \cdot)=id$;
\item[(ii)] \, $J_{\eps}(\eta(\cdot, u))$ is non-increasing;
\item[(iii)] \, $\eta(t, u)=u$ for every $t \in [0,1]$, if $J_{\eps}(u) \notin (E_{m_0}-\delta, E_{m_0}+\delta)$;
\item[(iv)] \, $\eta(1,\mc{X}_{\eps, \delta}^{c+\omega}\setminus U)\subset \mc{X}_{\eps, \delta}^{c-\omega}$;
\item[(v)] \, $\eta(\cdot, u)$ is a semigroup.
\end{itemize}
\end{Lemma}

\claim Proof. 
Let $\mathcal{V}: \{ u \in H^s(\R^N) \mid J'_{\eps}(u)\neq 0\} \to H^s(\R^N)$ be a locally Lipschitz pseudo-gradient vector field associated to $J_{\eps}$, and let $\phi \in Lip_{loc}(H^s(\R^N), \R)$ be a cutoff function such that 
$\supp (\phi) \subset (\mc{X}_{\eps, \delta})^{E_{m_0}+\delta}_{E_{m_0}-\delta}$ 
and $\phi=1$ in a small neighborhood of $c$. 
We consider the Cauchy problem
\begin{equation}\label{eq_dim_CP}
\parag{&\dot{\eta} = - \phi(\eta) \frac{\mathcal{V}(\eta)}{\norm{\mathcal{V}(\eta)}_{H^s(\R^N)}},& \\ &\eta(0, u)=u.&}
\end{equation}
The proof keeps on classically, obtaining a deformation $\eta: [0,1]\times \mathcal{X}_{\eps, \delta} \to H^s(\R^N)$. 
We want to prove now that $\eta$ goes into $\mathcal{X}_{\eps, \delta}$.

Let $u\in \mathcal{X}_{\eps, \delta}$. We need to show that $\eta(t,u)\in \mathcal{X}_{\eps, \delta}$ for every $t>0$. 
Since $\mathcal{X}_{\eps, \delta}$ is open, $\eta(s,u)$ continues staying in $\mathcal{X}_{\eps, \delta}$ for $s$ small. Thus assume that 
$$\eta(s,u) \in \mathcal{X}_{\eps,\delta}, \quad \hbox{ for every $0\leq s<t_0$}$$
for some $t_0>0$, and we want to show that $\eta(t_0, u)\in\mathcal{X}_{\eps, \delta}$. 
Notice first that, by using (iii), (v), (ii) and the continuity of $\eta$ and $J_{\eps}$ 
we can assume that 
\begin{equation}\label{eq_dim_stima[)}
J_{\eps}(\eta(t_0, u)) \in [E_{m_0}-\delta, E_{m_0}+\delta).
\end{equation}

\noindent
\textbf{Step 1:} \emph{$\eps \Upsilon(\eta(t_0,u)) \in \Omega[0, \nu_0)$.} 
\\
By \eqref{eq_appart_insieme_piccolo} we have
$$\eps \Upsilon(\eta(s,u)) \in \Omega[0,\nu_1), \quad \hbox{ for every $0\leq s<t_0$},$$
and thus by continuity 
$\eps \Upsilon(\eta(t_0,u)) \in \overline{\Omega[0,\nu_1)} \subset \Omega[0,\nu_0].$

\smallskip

\noindent
\textbf{Step 2:} \emph{$J_{\eps}(\eta(t_0,u)) < E_{m_0} + R(\delta, \eta(t_0, u))$.} 
\\
If $\widehat{\rho}(\eta(t_0, u)) \leq \rho_1$ then $R(\delta, \eta(t_0, u))=\delta$ and we directly have the claim, recalled that $J_{\eps}(\eta(t_0, u))<E_{m_0}+\delta$ by \eqref{eq_dim_stima[)}. 
Assume instead $\widehat{\rho}(\eta(t_0, u)) > \rho_1$. By continuity, there exists $t_1 \in (0, t_0)$ such that we have
$$\widehat{\rho}(\eta(s, u)) > \rho_1, \quad \hbox{ for every $s\in [t_1, t_0]$.} $$
In particular 
$$\parag{&\eta(s, u) \in S(\rho_0) \subset S(r_3) \subset S(r_2'),& \\ 
&J_{\eps}(\eta(s, u)) < E_{m_0}+\delta < E_{m_0+\nu_1},& \\ 
&\widehat{\rho}(\eta(s, u)) \in (\rho_1, \rho_0],& \\ 
&\eps \Upsilon(\eta(s, u)) \in \Omega[0, \nu_0]&
}
$$
for $s\in [t_1, t_0]$. 
Then by Theorem \ref{theorem_stima_J'} we have
$$\norm{J'_{\eps}(\eta(s, u))}_{(H^s(\R^N))^*} \geq \delta_2, \quad \hbox{ for every $s\in [t_1, t_0]$.} $$
We can thus compute with standard argument, by using \eqref{eq_dim_CP}, the properties of the pseudo-gradient and \eqref{eq_lipsc_rho}, 
\begin{align*}
J_{\eps}(\eta(t_0, u)) &\leq J_{\eps}(\eta(t_1, u)) - \frac{\delta_2}{2}\big(\widehat{\rho}(\eta(t_0, u)) - \widehat{\rho}(\eta(t_1, u))\big)\\
&< E_{m_0} + \delta - \frac{\delta_2}{2}\big(\widehat{\rho}(\eta(t_1,u)) - \rho_1\big) - \frac{\delta_2}{2}\big(\widehat{\rho}(\eta(t_0, u)) - \widehat{\rho}(\eta(t_1, u))\big)\\
&= E_{m_0} + R(\delta, \eta(t_0, u)),
\end{align*}
that is the claim.

\smallskip

\noindent
\textbf{Step 3:} \emph{$\eta(t_0, u)\in S(\rho_0)$.} 
\\
By the previous point we have
$J_{\eps}(\eta(t_0,u)) \leq E_{m_0} + R(\delta, \eta(t_0, u)).$
Since \eqref{eq_dim_stima[)} implies $J_{\eps}(\eta(t_0,u))\geq E_{m_0}-\delta$, then by \eqref{eq_cons_conten} we have $\eta(t_0,u) \in S(\rho_0)$, 
and thus the claim. 
\QED

\subsection{Maps homotopic to the embedding}

We search now for two maps $\Phi_{\eps}$, $\Psi_{\eps}$ such that, for a sufficiently small $\sigma_0\in (0,1)$ and a sufficiently small $\hat{\delta}=\hat{\delta}(\sigma_0)\in (0, \delta)$ (see \eqref{eq_propr_delta}), defined
$$ I:=[1-\sigma_0, 1+\sigma_0],$$
we have, for small $\eps$,
$$I \times K \; \stackrel{\Phi_{\eps}} \to \; \mathcal{X}_{\eps, \delta}^{E_{m_0}+\hat{\delta}} \;\stackrel{\Psi_{\eps}} \to \; I \times K_d $$
with the additional condition
$$\partial I \times K \;\stackrel{\Phi_{\eps}} \to \;\mathcal{X}_{\eps, \delta}^{E_{m_0}-\hat{\delta}} \;\stackrel{\Psi_{\eps}} \to\; (I\setminus \{1\}) \times K_d;$$
then we will prove that $\Psi_{\eps} \circ \Phi_{\eps}$ is homotopic to the identity. While the first property is useful for category arguments to gain multiplicity of solutions, the second additional condition will be essential for developing \emph{relative} category (and cup-length) arguments and controlling the sublevels of the functional \emph{below} the expected critical level.

\subsubsection{Definition of $\Phi_{\eps}$}\label{subsubsec_Phi}
\label{sec_def_phi_eps}

Let us fix a ground state $U_0 \in S_{m_0}\subset \widehat{S}$, i.e. $L_{m_0}(U_0)=E_{m_0}$ (see Theorem \ref{thm_esist_Poh_min} and \eqref{eq_ug_poh_gr.st}). Define, for $p \in K$ and $t \in I$ ($\sigma_0$ to be fixed)
$$\Phi_{\eps}(t,p):= U_0\left(\tfrac{\cdot-p/\eps}{t}\right)\in H^s(\R^N).$$
We show now that, for $\eps$ small, $\Phi_{\eps}(t,p) \in \mathcal{X}_{\eps, \delta}^{E_{m_0}+\hat{\delta}} $.
\begin{itemize}
\item 
$\Phi_{\eps}(t,p) \in S(\rho_1) \subset S(\rho_0)$:

indeed, recalled that the dilation $t\in \R \mapsto U_0(\cdot/t)\in H^s(\R^N)$ is continuous, we have
$$\norm{U_0\left(\tfrac{\cdot-p/\eps}{t}\right)-U_0(\cdot-p/\eps)}_{H^s(\R^N)}= \norm{U_0\left( \cdot/t\right)-U_0}_{H^s(\R^N)}< \rho_1$$
for $t \in I$ and sufficiently small $\sigma_0=\sigma_0(U_0)$ (not depending on $\eps$). Thus, setting $\varphi_t:= U_0\left(\frac{\cdot-p/\eps}{t}\right)-U_0(\cdot-p/\eps)$ we have
$$\Phi_{\eps}(t,p) = U_0(\cdot-p/\eps) + \varphi_t$$
with $U_0 \in \widehat{S}$, $p/\eps \in \R^N$ and $\norm{\varphi_t}_{H^s(\R^N)}<\rho_1$, which is the claim.

\item $\eps \Upsilon(\Phi_{\eps}(t,p)) \in \Omega[0, \nu_0)$:

indeed, by the previous point and Lemma \ref{lem_def_bar}, we have
$$\abs{\Upsilon(\Phi_{\eps}(t,p)) - p/\eps} < 2R_0$$
hence $\abs{\eps \Upsilon(\Phi_{\eps}(t,p)) - p} < 2\eps R_0$, and since $p \in K$
$$d(\eps \Upsilon(\Phi_{\eps}(t,p)), K) < 2\eps R_0.$$
For sufficiently small $\eps$, we have $K_{2\eps R_0} \subset \Omega[0, \nu_0)$, and thus the claim. 
In particular, for a later use observe that
\begin{equation}\label{eq_posizione_Ups}
\eps\Upsilon(\Phi_{\eps}(t,p)) = p + o(1).
\end{equation}

\item $J_{\eps}(\Phi_{\eps}(t,p)) < E_{m_0}+R(\delta, \Phi_{\eps}(t,p))$:

indeed $\Phi_{\eps}(t,p) \in S(\rho_1)$, thus $\widehat{\rho}(\Phi_{\eps}(t,p)) < \rho_1$, which implies $R(\delta, \Phi_{\eps}(t,p) ) = \delta$ and the claim comes from the following point, since $\hat{\delta}<\delta$.

\item $J_{\eps}(\Phi_{\eps}(t,p)) < E_{m_0}+\hat{\delta}$:

indeed we have by Lemma \ref{lemma_stima_g0} (b)
%\begin{eqnarray}
%\lefteqn{J_{\eps}(\Phi_{\eps}(t,p))\notag}\\ 
%&=& L_{m_0}\left(\Phi_{\eps}(t,p)\right) + \frac{1}{2}\int_{\R^N} (V(\eps x)- m_0) \Phi_{\eps}^2(t,p)dx+ Q_{\eps}\left(\Phi_{\eps}(t,p)\right) \notag\\
%&=:& L_{m_0}\left(U_0\left(\tfrac{\cdot-p/\eps}{t}\right)\right) + (I) +(II) = g(t) E_{m_0} + o(1)
% \label{eq_dim_stima_g(t)}\\
%&\leq & E_{m_0} +o(1) 
%\notag
%\end{eqnarray}
\begin{align}
J_{\eps}(\Phi_{\eps}(t,p))\notag
&= L_{m_0}\left(\Phi_{\eps}(t,p)\right) + \frac{1}{2}\int_{\R^N} (V(\eps x)- m_0) \Phi_{\eps}^2(t,p)dx+ Q_{\eps}\left(\Phi_{\eps}(t,p)\right) \notag\\
&=: L_{m_0}\left(U_0\left(\tfrac{\cdot-p/\eps}{t}\right)\right) + (I) +(II) = g(t) E_{m_0} + o(1)
 \label{eq_dim_stima_g(t)}\\
&\leq E_{m_0} +o(1) 
\notag
\end{align}
where we used $g(t) \leq 1$. Indeed, as regards $(I)$ we have
$$(I)= \frac{1}{2}\int_{\R^N} (V(\eps x + p)- m_0) U_0^2(x/t) dx\to 0 \quad \hbox{ as $\eps \to 0$}$$
by exploiting that $p \in K$ and the dominated convergence theorem, together with the boundedness of $V$. 
Focusing on $(II)$ instead, we have 
$$(II) = \left( \frac{1}{\eps^{\alpha}} \norm{U_0(\cdot/t)}_{L^2(\R^N \setminus ((\Omega_{2h_0} - p)/\eps)}^2 -1 \right)_+^{\frac{p+1}{2}};$$
since $p \in K \subset \Omega \subset \Omega_{2h_0}$, we have $0 \in \Omega_{2h_0}-p$ and moreover $B_r \subset \Omega_{2h_0}-p$ for some ball $B_r$; notice that $B_r/\eps$ covers the whole $\R^N$ as $\eps \to 0$. Therefore, by the polynomial estimate we have
$$ \norm{U_0(\cdot/t)}_{L^2(\R^N \setminus ((\Omega_{2h_0} - p)/\eps)}^2 \leq C \norm{\tfrac{1}{1+|x|^{N+2s}}}_{L^2(\R^N \setminus (B_r/\eps))} ^2 \leq C \eps^{N+4s},$$
and hence $(II) \to 0$ as $\eps \to 0$, since $\alpha < N+4s$. 
Therefore, by choosing a sufficiently small $\eps$, we obtain
$$J_{\eps}(\Phi_{\eps}(t,p)) \leq E_{m_0} + \tfrac{1}{2} \hat{\delta} < E_{m_0} + \hat{\delta}.$$
\end{itemize}
\vspace{-\topsep}
Finally, we show the additional condition.

\begin{itemize}
\item $J_{\eps}(\Phi_{\eps}(1\pm \sigma_0,p)) < E_{m_0}-\hat{\delta}$:

indeed, looking at \eqref{eq_dim_stima_g(t)} we see that, for small $\eps$,
$$J_{\eps}(\Phi_{\eps}(1\pm \sigma_0 ,p)) < g(1\pm \sigma_0) E_{m_0} + \hat{\delta};$$
since $ g(1\pm \sigma_0)<1$, we can find a small $\hat{\delta}< \frac{1-g(1\pm \sigma_0)}{2} E_{m_0}$ (not depending on $\eps$) such that
\begin{equation}\label{eq_buona_post_Phieps}
J_{\eps}(\Phi_{\eps}(1\pm \sigma_0 ,p)) < g(1\pm \sigma_0) E_{m_0} + \hat{\delta} < E_{m_0} - \hat{\delta}
\end{equation}
and thus the claim.
\end{itemize}
\vspace{-\topsep}

\subsubsection{Definition of $\Psi_{\eps}$}

Define a truncation
$$T(t):= \parag{1-\sigma_0 & \quad \hbox{if $t \leq 1-\sigma_0$},\\ t & \quad \hbox{if $t \in (1-\sigma_0 ,1+\sigma_0)$}, \\ 1+\sigma_0 & \quad \hbox{if $t \geq 1+\sigma_0$}}$$
for $t \in \R$, and
$$\Psi_{\eps}(u):= \big( T(P_{m_0}(u)), \eps \Upsilon(u)\big)$$
for every $u \in \mathcal{X}_{\eps, \delta}^{E_{m_0}+\hat{\delta}}$. By the definition of $T$ and property \eqref{eq_appart_insieme_piccolo}, we have directly
$$\Psi_{\eps}(u) \in I \times \Omega[0, \nu_1] \subset I \times \Omega[0, \nu_0] \subset I \times K_d.$$
Assume now $u\in \mathcal{X}_{\eps, \delta}^{E_{m_0}-\hat{\delta}}$. We have, by using Lemma \ref{lemma_stima_J} and Lemma \ref{lemma_stima_g}, 
$$E_{m_0}-\hat{\delta} \geq J_{\eps}(u) \geq L_{m_0}(u) - C_{min} \eps^{\alpha} \geq g(P_{m_0}(u)) E_{m_0} - C_{min}\eps^{\alpha}$$
and hence
$$E_{m_0} \geq g(P_{m_0}(u)) E_{m_0} + \hat{\delta} - C_{min}\eps^{\alpha}> g(P_{m_0}(u)) E_{m_0}$$
where the last inequality holds for $\eps$ small, not depending on $u$. Thus
$$g(P_{m_0}(u))< 1$$
and this must imply, by the properties of $g$, that $P_{m_0}(u) \neq 1$, and in particular
$$T(P_{m_0}(u))\neq 1.$$
This reaches the goal.

\subsubsection{An homotopy to the identity}

Introduce the notation of \emph{topological pair} from the algebraic topology: we write, for $ B \subset A$ and $B' \subset A'$,
$$f: (A,B) \to (A', B')$$
whenever
$$f\in C(A,A')\quad \hbox{ and } \quad f(B) \subset B'.$$
Observed that $\Phi_{\eps}$ and $\Psi_{\eps}$ are continuous, we can rewrite the stated properties as
$$\Phi_{\eps} : \Big( I\times K, \, \partial I \times K \Big) \to \Big( \mathcal{X}_{\eps, \delta}^{E_{m_0}+\hat{\delta}}, \, \mathcal{X}_{\eps, \delta}^{E_{m_0}-\hat{\delta}}\Big),$$
$$\Psi_{\eps}: \Big( \mathcal{X}_{\eps, \delta}^{E_{m_0}+\hat{\delta}}, \, \mathcal{X}_{\eps, \delta}^{E_{m_0}-\hat{\delta}}\Big) \to \Big( I \times K_d, \, (I\setminus \{1\}) \times K_d\Big)$$
and
$$\Psi_{\eps}\circ \Phi_{\eps}: \Big( I \times K, \, \partial I \times K \Big) \to \Big( I \times K_d, \, (I \setminus \{1\}) \times K_d \Big),$$
where a straightforward computation shows
$$ (\Psi_{\eps}\circ \Phi_{\eps}) (t, p) = \left( t, \, \eps \Upsilon\left( U_0\left(\tfrac{\cdot-p/\eps}{t}\right)\right)\right),$$
thus actually $\Psi_{\eps}\circ \Phi_{\eps}: \Big( I \times K, \, \partial I \times K \Big) \to \Big( I \times K_d, \, \partial I \times K_d \Big)$.
Clearly, we notice that the inclusion map has the same property, that is set $j(t,p):=(t,p)$ we have
$$j: \Big( I \times K, \, \partial I \times K \Big) \to \Big( I \times K_d, \, \partial I \times K_d\Big) \subset \Big( I \times K_d, \, (I \setminus \{1\}) \times K_d\Big).$$
We want to show that these maps are homotopic, information useful in the theory of relative cup-length.

\begin{Proposition} \label{prop_esist_homot}
For sufficiently small $\eps$, the maps $\Psi_{\eps}\circ \Phi_{\eps}$ and $j$ are homotopic, that is there exists a continuous map 
$H: [0,1] \times I \times K \to I \times K_d$
such that
$$H(\theta, \cdot, \cdot): \Big( I \times K, \, \partial I \times K \Big) \to \Big( I \times K_d, \, \partial I \times K_d\Big) \subset \Big( I \times K_d, \, (I \setminus \{1\}) \times K_d\Big)$$
for each $\theta \in [0,1]$, with
$H(0, \cdot, \cdot) = \Psi_{\eps}\circ \Phi_{\eps}$ and $ H(1, \cdot, \cdot) = j.$
\end{Proposition}

\claim Proof.
Noticed that also $\Psi_{\eps}\circ \Phi_{\eps}$ fixes the first variable, it is sufficient to link the second variables through a segment, that is
$$H(\theta, t, p) := \left( t, \, (1-\theta) \eps \Upsilon\left( U_0\left(\tfrac{\cdot-p/\eps}{t}\right)\right) + \theta p\right),$$
with $\theta \in [0,1]$. We must check that $H$ is well defined, since $K_d$ is not a convex set, generally. Indeed we have, by \eqref{eq_posizione_Ups}
$$(1-\theta) \eps \Upsilon\left( U_0\left(\tfrac{\cdot-p/\eps}{t}\right)\right) + \theta p = (1-\theta) p + o(1) + \theta p = p + o(1).$$
Since $p \in K$, for sufficiently small $\eps$ we have that $p+o(1) \in K_d$, and thus the claim.
\QED

\bigskip

Before coming up to multiplicity results, we highlight that existence of a single solution could be obtained without any use of algebraic tools. Notice that we need only the map $\Phi_{\eps}$ and the first component of $\Psi_{\eps}$.
%: the idea is to wisely use the deformation lemma and the definition of $E_{m_0}$, together with the existence of Pohozaev minima, dilatations and the intermediate value theorem.
%We omit the details, gaining instead the existence as a byproduct of the multiplicity.

%\begin{Remark}\label{rem_sola_esist}
%Before coming up to multiplicity results, we highlight that existence of a single solution could be obtained without any use of algebraic tools: 
%the idea is to wisely use the deformation lemma and the definition of $E_{m_0}$, together with the existence of Pohozaev minima, dilatations and the intermediate value theorem.
%We omit the details, gaining instead the existence as a byproduct of the multiplicity.
%\tor{Inserisci}
%\end{Remark}
%

\medskip

 %Uso solo $\phi_{\eps}$, non $\psi_{\eps}$!
\claim Proof (existence).
Let $p\in K$. First observe that we can slightly change the map $\Phi_{\eps}$ such that 
\begin{equation}\label{eq_dim_miglior_sigma0}
\Phi_{\eps}(1\pm \sigma_0, p)\in \mc{X}_{\eps, \delta}^{E_{m_0}-\delta};
\end{equation}
indeed (see \eqref{eq_dim_stima_g(t)} and \eqref{eq_buona_post_Phieps}), it is sufficient to take a smaller $\eps>0$ and $\hat{\delta} < \delta < \frac{1-g(1\pm \sigma_0)}{2} E_{m_0} < E_{m_0}$, where we point out that $\sigma_0$ depends only on $U_0$ and $\rho_1$ (and thus not on $\delta$).

Let $c= E_{m_0}$; by contradiction, assume $K_c \neq 0$. Thus, by the Lemma \ref{lem_def_lem}, there exists a deformation $\eta$ related to the regular value $c$. 
By Lemma \ref{lemma_stima_J} we have, for each $\sigma \in I$,
$$L_{m_0}(\eta(1, \Phi_{\eps}(\sigma, p))) \leq J_{\eps}(\eta(1,\Phi_{\eps}(\sigma, p))) +C_{min} \eps^{\alpha} \leq 
%E_{m_0}-c+ C_{min}\eps^{\alpha}$$
E_{m_0}- \delta + C_{min}\eps^{\alpha}$$
where in the last inequality we have used that $\Phi_{\eps}(\sigma, p) \in \mc{X}_{\eps, \delta}^{E_{m_0}+\hat{\delta}}
%\tr{\subset \mc{X}_{\eps, \delta}^{E_{m_0}+c}}
\subset \mc{X}_{\eps, \delta}^{c+\delta}
$. Thus, for $\eps$ small, we have
$$L_{m_0}(\eta(1, \Phi_{\eps}(\sigma, p))) < E_{m_0} \quad \hbox{ for each $\sigma \in I$}.$$
To conclude, we need to find a $\tilde{\sigma} \in I$ such that
\begin{equation*}%\label{eq_dim_cond_poho}
P_{m_0}(\eta(1, \Phi_{\eps}(\tilde{\sigma}, p))) = 1
\end{equation*}
since this implies $L_{m_0}(\eta(1, \Phi_{\eps}(\tilde{\sigma}, p)))\geq C_{po, m_0} = E_{m_0}$ and thus an absurd.

Indeed, by \eqref{eq_dim_miglior_sigma0} we have %(notice that $\partial I= \{1\pm \sigma_0\}$)
$$P_{m_0}(\eta(1%t
, \Phi_{\eps}(1\pm \sigma_0,p)))= P_{m_0}(\Phi_{\eps}(1\pm \sigma_0, p))= 1\pm \sigma_0 %\neq 1.
$$
%By the homotopy invariance of the degree we have
%$$\deg(H(\cdot, 1), I, 1)=\deg(H(\cdot, 0), I, 1)= \deg(P_{m_0}(\Phi_{\eps}(\cdot, p)), I, 1)$$
%where again $P_{m_0}(\Phi_{\eps}(\sigma, p))= \sigma$ and thus
%$$\deg(P_{m_0}(\eta(1, \Phi_{\eps}(\cdot, p))), I, 1) = \deg(id, I, 1)=1\neq 0.$$
%This concludes the proof.
and the claim follows by the intermediate value theorem.
\QED

%%%%%%%%%%%%%%%%%%%%%%%%%%%%%%%%%%%%%%%%
%%%%%%%%%%%%%%%%%%%%%%%%%%%%%%%%%%%%%%%%
\section{Existence of multiple solutions}
\label{sec_cup-length}

We finally come up to the existence of multiple solutions. 
Here the algebraic notions of \emph{relative category} and \emph{relative cup-length} (built on the Alexander-Spanier cohomology with coefficients in some field $\mathbb{F}$) are of key importance. We refer to the Appendix \ref{chap_app_alg_top} for definitions, comments and properties of these algebraic tools.

\medskip

\claim Proof of Theorem \ref{teo_concen_esist}.
By construction of the neighborhood $\mc{X}_{\eps, \delta}$ and Corollary \ref{corol_passaggio_critici} (recall that $\rho_0 < r_3 \leq r_2'$ and that $J_{\eps}(u) < E_{m_0} + R(\hat{\delta}, u) \leq E_{m_0} + \hat{\delta} \leq l_0'$ for $u \in \mc{X}_{\eps, \delta}$), we have
$$
\Big\{ u \in (\mathcal{X}_{\eps, \delta})^{E_{m_0}+\hat{\delta}}_{E_{m_0}-\hat{\delta}} \mid J'_{\eps}(u)=0 \Big\} \subset \left\{ u \in H^s(\R^N) \mid I'_{\eps}(u)=0 \right\}.
$$
Thus we obtain
\begin{eqnarray*}
\lefteqn{ \# \{ u \hbox{ solutions of \eqref{eq_princ_epsx}}\} \geq \# \Big\{ u \in (\mathcal{X}_{\eps, \delta})^{E_{m_0}+\hat{\delta}}_{E_{m_0}-\hat{\delta}} \mid J'_{\eps}(u)=0 \Big\} } \\
&\stackrel{(i)}\geq& \cat\Big (\mathcal{X}_{\eps, \delta}^{E_{m_0}+\hat{\delta}},\, \mathcal{X}_{\eps, \delta}^{E_{m_0}-\hat{\delta}}\Big) 
\stackrel{(ii)}\geq \cupl\Big (\mathcal{X}_{\eps, \delta}^{E_{m_0}+\hat{\delta}},\, \mathcal{X}_{\eps, \delta}^{E_{m_0}-\hat{\delta}}\Big) +1 \\
&\stackrel{(iii)}\geq& \cupl(K) +1
\end{eqnarray*}
that is the claim, up to the proof of (i)--(iii). Indeed, (i) is obtained classically from the Deformation Lemma \ref{lem_def_lem} as in Section \ref{sec_mult_cat_cupl}. 
Inequality (ii) is given by the algebraic-topological Lemma \ref{lem_coll_cat_cupl}. 
Point (iii) is instead due to the existence of the homotopy gained in Proposition \ref{prop_esist_homot} and properties of the cup-length: indeed,
by \eqref{eq_propr_cupl} in Lemma \ref{lem_propr_cupl} (a), we have
$$\cupl \Big (\mathcal{X}_{\eps, \delta}^{E_{m_0}+\hat{\delta}},\, \mathcal{X}_{\eps, \delta}^{E_{m_0}-\hat{\delta}}\Big) \geq \cupl (\Psi_{\eps} \circ \Phi_{\eps});$$
moreover, since $\Psi_{\eps} \circ \Phi_{\eps}$ is homotopic to the immersion $j$ thanks to Proposition \ref{prop_esist_homot}, we have by Lemma \ref{lem_propr_cupl} (b)
$$ \cupl (\Psi_{\eps} \circ \Phi_{\eps}) = \cupl(j),$$
which leads to the conclusion thanks to Lemma \ref{lem_cupl_jK}. See Remark \ref{rem_regolar} for the proof of regularity.
\QED

\subsection{Concentration in the potential well} %$K$} 
\label{sec_concentration}

We prove now the polynomial decay and the concentration of the found solutions in $K$. 
To deal with uniform bound, we will make use of the fractional De Giorgi class recalled in Section \ref{sec_regol_degiorgi}.

\medskip

\claim Proof of Theorem \ref{teo_concen_conc}. 
For $\eps$ sufficiently small, let $u_{\eps}$ be one of the $\cupl(K)+1$ critical points of $J_{\eps}$ built in Theorem \ref{teo_concen_conc}, which by Corollary \ref{corol_passaggio_critici} is also a solution of \eqref{eq_princ_epsx}, positive by \hyperref[(f2)]{\textnormal{(f2)}}. In particular, since it satisfies the assumptions of Lemma \ref{lemma_stima_basso}, looking at the proof (see \eqref{eq_conv_uu} and \eqref{eq_conv_p_eps}) we obtain that
$$\norm{u_{\eps} - U(\cdot - p_{\eps})}_{H^s(\R^N)}\to 0$$
with $U \in S_{V(p_0)}$, $p_{\eps}\in \R^N$ and
$$\eps p_{\eps} \to p_0 \in \Omega[0, \nu_1).$$

\smallskip

\noindent
\textbf{Step 1.}
Notice that we have found these solutions by fixing $\nu_0$, $l_0$ and $l_0'$. Let them move, throughout three sequences $\nu_0^n \searrow 0$, $l_0^n \searrow E_{m_0}$, and $(l_0')^n \searrow E_{m_0}$, 
and find the corresponding (sufficiently small) $\eps_n >0$ such that $\cupl(K)+1$ solutions exist; let $u_{\eps_n}$ be one of those and $p_{\eps_n}$ as before. It is not reductive to assume $\eps_n \to 0$ as $n \to +\infty$; by a diagonalization-like argument we obtain
\begin{equation} \label{eq_dim_conv_u_x}
u_{\eps_n}(\cdot+p_{\eps_n})\to U \; \; \hbox{ in $H^s(\R^N)$}, \; \hbox{ for some $U$ least energy solution of \eqref{eq_least_energy_m0}},
\end{equation}
$$\eps_n p_{\eps_n} \to p_0 \in K,$$
as $n \to +\infty$.

\smallskip

\noindent
\textbf{Step 2.} %\label{pag_dim_concent}
From now on we write $\eps\equiv \eps_n$ to avoid cumbersome notation. 
By $I_{\eps}'(u_{\eps})=0$ we obtain
\begin{equation}\label{eq_dim_perconfr}
(-\Delta)^s u_{\eps} + V(\eps x) u_{\eps} = f(u_{\eps}) , \quad x \in \R^N,
\end{equation}
thus (recall that $u_{\eps}$ is positive), by choosing $\beta < \underline{V}$ in \eqref{eq_prop_f}, 
$$(-\Delta)^s u_{\eps} \leq -\underline{V} u_{\eps} + f(u_{\eps})\leq (\beta-\underline{V})u_{\eps} + C_{\beta} u_{\eps}^p \leq C_{\beta} u_{\eps}^p, \quad x \in \R^N.$$
Therefore by Theorem \ref{thm_degiorgi_class} we have, choosing $q=p+1 %p
$, $d_1=0$ and $d_2=C_{\beta}$,
$$u_{\eps} \in \DG_+^{s,2} \Big( B_{R_0}(x_0), 0, H, 0, 1-\frac{p+1}{2^*_s}, 2s, R_0\Big),$$
%\tr{qui è importante che sia strettamente sottocritico!} %COMMENT NOW
%
with $H=H(N, s, p, \beta)$ 
and $R_0$ depending on $N, s, p, C_{\beta}$ and a uniform upper bound of the $H^s$-norms of $u_{\eps }$.
 
We can thus use now \cite[Proposition 6.1]{Coz0}: observing that $d(x_0, \partial B_{R_0}(x_0))= R_0$, and that $\mu = 1-\frac{p+1}{2^*_s}$, we obtain, for any $\omega\in (0,1]$ and $R\in (0, \frac{R_0}{2})$, 
$$\sup_{B_R(x_0)} u_{\eps} \leq \frac{C}{(N-2s)^{\frac{1}{2\mu }}} \frac{1}{\omega^{\frac{1}{2\mu }}} \frac{1}{(2R)^{N/2}} \norm{u_{\eps}}_{L^2(B_{2R}(x_0))} + \omega \Tail(u_{\eps};x_0, R)$$
that is, rewriting the constant $C=C(N, s, p, \beta)$,
$$\sup_{B_R(x_0)} u_{\eps} \leq C\frac{1}{\omega^{\frac{1}{2\mu }}} \frac{1}{R^{N/2}} \norm{u_{\eps}}_{L^2(B_{2R}(x_0))} + \omega \Tail(u_{\eps};x_0, R).$$

\smallskip

\noindent
\textbf{Step 3.}
We have
\begin{align*}
\norm{u_{\eps}}_{L^{\infty}(\R^N)} &= \sup_{x_0 \in \R^N} \sup_{B_R(x_0)} u_{\eps} \\
&\leq \sup_{x_0 \in \R^N} \left(C\frac{1}{\omega^{\frac{1}{2\mu }}} \frac{1}{R^{N/2}} \norm{u_{\eps}}_{L^2(B_{2R}(x_0))} + \omega \Tail(u_{\eps};x_0, R)\right).
\end{align*}
Observe that, by definition of Tail function \eqref{eq_def_tail} and H\"older inequality,
\begin{align*}
\Tail(u_{\eps};x_0, R) &\leq (1-s) R^{2s} \norm{u_{\eps}}_{L^2(\R^N \setminus B_R(x_0))} \norm{\tfrac{1}{|x-x_0|^{N+2s %2N+4s
}}}_{L^2(\R^N \setminus B_R(x_0))} \\
&\leq \frac{C}{R^{N/2}}\norm{u_{\eps}}_{L^2(\R^N)}.
\end{align*}
Thus
\begin{align*}
\norm{u_{\eps}}_{L^{\infty}(\R^N)} &\leq \frac{C}{R^{N/2}}\sup_{x_0 \in \R^N} \left(\omega^{-\frac{1}{2\mu }} \norm{u_{\eps}}_{L^2(B_{2R}(x_0))} + \omega \norm{u_{\eps}}_{L^2(\R^N \setminus B_R(x_0))}\right)\\
&\leq \frac{C}{R^{N/2}}\left(\omega^{-\frac{1}{2\mu }} + \omega\right) \norm{u_{\eps}}_{L^2(\R^N)}
\end{align*}
which is uniformly bounded by the properties on $u_{\eps}$. Hence $u_{\eps}$ are uniformly bounded in $L^{\infty}(\R^N)$.

In addition, by the estimates on $V$, $f$ and $u_{\eps}$, we have 
$$g_{\eps}(x):= -V(\eps x) u_{\eps}(x) + f(u_{\eps}(x)) \in L^{\infty}(\R^N)$$
with bound uniform in $\eps$; since
$$(-\Delta)^s u_{\eps} = g_{\eps}(x), \quad x \in \R^N,$$
by \cite[Theorem 8.2]{Coz0} there exists $\sigma \in (0,1)$, not depending on $u_{\eps}$, and $C=C(N,s)$, such that, for each $R>1$ and $x_0 \in \R^N$, 
%\begin{eqnarray*}
%\lefteqn{[u_{\eps}]_{C^{0,\sigma}(B_R(x_0))}}\\
%&\leq& \frac{C}{R^{\sigma}} \left( \norm{u_{\eps}}_{L^{\infty}(B_{4R}(x_0))} + \Tail(u_{\eps}; x_0, 4R) + R^{2s} \norm{g_{\eps}}_{L^{\infty}(B_{8R}(x_0))}\right)
%\leq C'.\label{eq_stima_uniforme_C0alpha}
%\end{eqnarray*}
\begin{align}
[u_{\eps}]_{C^{0,\sigma}(B_R(x_0))}
&\leq \frac{C}{R^{\sigma}} \left( \norm{u_{\eps}}_{L^{\infty}(B_{4R}(x_0))} + \Tail(u_{\eps}; x_0, 4R) + R^{2s} \norm{g_{\eps}}_{L^{\infty}(B_{8R}(x_0))}\right) \notag \\
&\leq C'.\label{eq_stima_uniforme_C0alpha}
\end{align}
We highlight that, since the constant is uniform in $R>1$, we obtain $u_{\eps} \in C^{0, \sigma}(\R^N)$.

\smallskip

\noindent
\textbf{Step 4.}
By the local uniform estimate on $u_{\eps}$ we could gain $\norm{u_{\eps}}_{L^{\infty}(\R^N\setminus \left(\Omega_{h_0}/\eps)_{2R_0}\right)}\to 0$, but this lack of uniformity on the domain can be improved. 
Thus we exploit the tightness of $\tilde{u}_{\eps}$ to reach the claim, where
$$\tilde{u}_{\eps}:=u_{\eps}(\cdot + p_{\eps}).$$
Indeed, by Step 2, and \eqref{eq_dim_conv_u_x} we have
$$\parag{ &(-\Delta)^s \tilde{u}_{\eps} +V(\eps x+ \eps p_{\eps}) \tilde{u}_{\eps} = f(\tilde{u}_{\eps}), \quad x \in \R^N,& \\ &\norm{\tilde{u}_{\eps}}_{\infty}\leq C,& \\ &\tilde{u}_{\eps} \to U \quad \hbox{ in $H^s(\R^N)$ as $\eps \to 0$}, \quad \hbox{ $U$ least energy solution of \eqref{eq_least_energy_m0}}.&}$$
In particular, it is standard to show that $f(\tilde{u}_{\eps}) \to f(U)$ in $L^2(\R^N)$, $\norm{f(\tilde{u}_{\eps})}_{\infty}\leq C$ and $U, f(U) \in L^{\infty}(\R^N)$. By interpolation we thus obtain
$$\chi_{\eps}:=\tilde{u}_{\eps} + f(\tilde{u}_{\eps})\to \chi:=U+f(U) \quad \hbox{ in $L^{q'}(\R^N)$}$$
for every $q'\in [2, +\infty)$, and $\norm{\chi_{\eps}}_{\infty} \leq C.$ 
Proceeding as in the proof of Lemma \ref{lem_u_to0} we gain
\begin{equation}\label{eq_conv_unif_0_veps}
\tilde{u}_{\eps}(x) \to 0 \; \; \hbox{ as $|x|\to +\infty$}, \quad \hbox{uniformly in $\eps$}.
\end{equation}
For the reader's convenience, we give some details. 
%\tr{diminuisci dettagli} %COMMENT NOW
 Indeed, being $\tilde{u}_{\eps}$ solution of
$$(-\Delta)^s \tilde{u}_{\eps} + \tilde{u}_{\eps} = \chi_{\eps} - V(\eps x + \eps p_{\eps}) \tilde{u}_{\eps}, \quad x \in \R^N,$$
we have the representation formula
$$\tilde{u}_{\eps}= \mc{K}_{2s} * (\chi_{\eps} - V(\eps x + \eps p_{\eps}) \tilde{u}_{\eps})$$
where $\mc{K}_{2s}$ is the Bessel kernel. 
Let us fix $\eta>0$; since $V$, $\tilde{u}_{\eps}$ and $\mc{K}_{2s}$ are positive, we have, for $x \in \R^N$,
\begin{align*}
\tilde{u}_{\eps}(x) &= \int_{\R^N} \mc{K}_{2s}(x-y) \big(\chi_{\eps}(y) - V(\eps x + \eps p_{\eps})\tilde{u}_{\eps}(y)\big)dy \\
&\leq \int_{|x-y|\geq 1/\eta} \mc{K}_{2s}(x-y) \chi_{\eps}(y)dy +\int_{|x-y|< 1/\eta} \mc{K}_{2s}(x-y) \chi_{\eps}(y)dy.
\end{align*}
As regards the first piece
$$ \int_{|x-y|\geq 1/\eta} \mc{K}_{2s}(x-y) \chi_{\eps}(y)dy \leq \norm{\chi_{\eps}}_{\infty} \int_{|x-y|\geq 1/\eta} \frac{C}{|x-y|^{N+2s}} dy \leq C \eta^{2s}$$
while for the second piece, fixed a whatever $q \in (1, \min\{2,\frac{N}{N-2s}%1 + \tfrac{2s}{N-2s}
\})$ and its conjugate exponent $q' \in ( \max\{2,\frac{N}{2s}\}, + \infty)$, we have by H\"older inequality
\begin{align*}
\int_{|x-y|< 1/\eta} \mc{K}_{2s}(x-y) \chi_{\eps}(y)dy &\leq \norm{\mc{K}_{2s}}_q \norm{\chi_{\eps}}_{L^{q'}(B_{1/\eta}(x))}\\
& \leq \norm{\mc{K}_{2s}}_q\left( \norm{\chi_{\eps}-\chi}_{q'}+ \norm{\chi}_{L^{q'}(B_{1/\eta}(x))}\right)
\end{align*}
where the first norm can be made small for $\eps<\eps_0=\eps_0(\eta)$, while the second 
 for $|x| \gg 0$ (uniformly in $\eps$). On the other hand, for $\eps\geq \eps_0$ (and thus for a finite number of elements, since we recall we are working with $\eps \equiv \eps_n$ small) the quantity $ \norm{\chi_{\eps}}_{L^{q'}(B_{1/\eta}(x))}$ can be made small for $|x|\gg 0$, uniformly in $\eps$. Joining the pieces, we have \eqref{eq_conv_unif_0_veps}.

\smallskip

\noindent
\textbf{Step 5.}
Let now $y_{\eps}\in \R^N$ be a maximum point for $u_{\eps}$, which exists by the boundedness of $u_{\eps}$ and its continuity (see \eqref{eq_stima_uniforme_C0alpha}). 
Therefore $z_{\eps}:=y_{\eps}-p_{\eps}$ 
 is a maximum point for $\tilde{u}_{\eps}$. In particular
$$\tilde{u}_{\eps}(z_{\eps})=\max_{\R^N} \tilde{u}_{\eps}=\norm{\tilde{u}_{\eps}}_{\infty}\not \to 0 \quad \hbox{ as $\eps \to 0$}$$
since on the contrary we would have $\tilde{u}_{\eps}\to 0$ almost everywhere, which is in contradiction with the fact that $\tilde{u}_{\eps}\to U\nequiv 0$ almost everywhere (up to a subsequence). As a consequence, thanks to \eqref{eq_conv_unif_0_veps}, we have that $z_{\eps}$ is bounded (up to a subsequence). That is, again up to a subsequence,
$$z_{\eps} \to \overline{p}$$
for some $\overline{p}\in \R^N$. In particular
$$\eps y_{\eps} = \eps z_{\eps} + \eps p_{\eps} \to p_0 \in K$$
and, by the fact that
$$U(\cdot + z_{\eps}) \to U(\cdot + \overline{p}) =: \overline{U} \quad \hbox{ in $H^s(\R^N)$}$$
we have 
$u_{\eps}(\cdot + y_{\eps}) \to \overline{U}$ in $H^s(\R^N)$, $\overline{U}$ least energy solution of \eqref{eq_least_energy_m0}. 
We set
$$\overline{u}_{\eps}:= u_{\eps}(\cdot + y_{\eps}),$$
$$\overline{u}_{\eps}\to \overline{U} \; \; \hbox{ in $H^s(\R^N)$}, \quad \hbox{ $\overline{U}$ least energy solution of \eqref{eq_least_energy_m0}};$$
in addition, $\overline{u}_{\eps}$ is positive by \hyperref[(f2)]{\textnormal{(f2)}}, and in the same way we obtained \eqref{eq_conv_unif_0_veps} we obtain also
\begin{equation}\label{eq_conv_unif_0_veps2}
\overline{u}_{\eps}(x) \to 0 \; \; \hbox{ as $|x|\to +\infty$}, \quad \hbox{uniformly in $\eps$}.
\end{equation}
Moreover, by exploiting the uniform estimates in $L^{\infty}(\R^N)$ and $C^{0, \sigma}_{loc}(\R^N)$ we obtain by Ascoli-Arzell\`{a} theorem also that $\overline{u}_{\eps}\to \overline{U}>0$ in $L^{\infty}_{loc}(\R^N)$, with $\overline{U}$ continuous; this easily implies, for every $r>0$, that
\begin{equation}\label{eq_stima_basso_min}
\min_{B_r} \overline{u}_{\eps} \geq \frac{1}{2} \min_{B_r} \overline{U} >0
\end{equation}
for $\eps$ small, depending on $\overline{U}$ and $r$.

\smallskip

\noindent
\textbf{Step 6.}
By \eqref{eq_conv_unif_0_veps2} we have, for $R'$ large (uniform in $\eps$), that
$$\overline{u}_{\eps}(x)\leq \eta', \quad \hbox{ for $|x|>R'$}$$
for every $\eps>0$, where $\eta'>0$ is preliminary fixed. 
As a consequence, by \hyperref[(f1)]{\textnormal{(f1.2)}}, we gain
$$-\frac{1}{2} \overline{V} \overline{u}_{\eps}(x)\leq f(\overline{u}_{\eps}(x)) \leq \frac{1}{2}\underline{V} \overline{u}_{\eps}(x),\quad \hbox{ for $|x|>R'$},$$
where $\overline{V}:=\norm{V}_{\infty}$. We obtain by \eqref{eq_dim_perconfr}
$$(-\Delta)^s \overline{u}_{\eps} + \tfrac{1}{2}\underline{V} \overline{u}_{\eps} \leq f(\overline{u}_{\eps}) - \tfrac{1}{2} \underline{V}\overline{u}_{\eps} \leq 0, \quad x \in \R^N \setminus B_{R'},$$
$$(-\Delta)^s \overline{u}_{\eps} + \tfrac{3}{2} \overline{V} \overline{u}_{\eps} \geq f(\overline{u}_{\eps}) + \tfrac{1}{2} \overline{V}\overline{u}_{\eps}\geq 0, \quad x \in \R^N \setminus B_{R'}.$$
Notice that we always intend differential inequalities in the weak sense.
%, that is tested with functions in $H^s(\R^N)$ with supports contained (e.g.) in $\R^N \setminus B_{R'}$. 
%
In addition, by Lemma \ref{lem_esist_sol_part} we have that there exist two positive functions $\underline{W}'$, $\overline{W}'$ and three positive constants $R''$, $C'$ and $C''$ depending only on $V$, such that
$$ \parag{
& (-\Delta)^s \underline{W}' + \frac{3}{2}\overline{V} \, \underline{W}' = 0, \quad x \in \R^N \setminus B_{R''},& \\ 
&\frac{C'}{|x|^{N+2s}}< \underline{W}' (x), \quad \hbox{ for $|x|>2R''$},&}$$
and
$$ \parag{& (-\Delta)^s \overline{W}' + \frac{1}{2}\underline{V} \, \overline{W}' = 0, \quad x \in \R^N \setminus B_{R''},& \\ 
& \overline{W}'(x) < \frac{C''}{|x|^{N+2s}}, \quad \hbox{ for $|x|>2R''$},&}$$
Set $R:=\max\{ R', 2R''\}$. Let $\underline{C}_1$ and $\overline{C}_1$ be some uniform lower and upper bounds for $\overline{u}_{\eps}$ on $B_R$, $\underline{C}_2:=\min_{B_R} \overline{W}'$ and $\overline{C}_2:= \max_{B_R} \underline{W}'$, all strictly positive. Define
$$\underline{W}:= \underline{C}_1 \overline{C}_2 ^{-1} \underline{W}', \quad \overline{W}:= \overline{C}_1 \underline{C}_2^{-1} \overline{W}'$$
so that
$$\underline{W}\leq \overline{u}_{\eps} \leq \overline{W}, \quad \hbox{ for $|x|\leq R$}.$$
Through a Comparison Principle (see Lemma \ref{lem_comp_prin}), and redefining $C'$ and $C''$, we obtain
$$ \frac{C'}{|x|^{N+2s}} <\underline{W}(x) \leq\overline{u}_{\eps}(x) \leq \overline{W}(x) < \frac{C''}{|x|^{N+2s}}, \quad \hbox{ for $|x|>R$}.$$
By the uniform boundedness of $\overline{u}_{\eps}$ and \eqref{eq_stima_basso_min} we also obtain
$$\frac{C'}{1+|x|^{N+2s}}< \overline{u}_{\eps}(x) < \frac{C''}{1+|x|^{N+2s}}, \quad \hbox{ for $x \in \R^N $}.$$
Recalling the definition of $\overline{u}_{\eps}$, we have finally obtained a sequence of solutions such that
$$\parag{&u_{\eps_n}(y_{\eps_n})=\max_{\R^N} u_{\eps_n}, &\\& 
d(\eps_n y_{\eps_n}, K)\to 0, &\\& 
\frac{C'}{1+|x-y_{\eps_n}|^{N+2s}} \leq u_{\eps_n}(x) \leq \frac{C''}{1+|x-y_{\eps_n}|^{N+2s}}, \quad \hbox{ for $x \in \R^N$}, &\\& 
\norm{u_{\eps_n}(\cdot + y_{\eps_n}) - \overline{U}}_{H^s(\R^N)} \to 0, \quad \hbox{ for some $\overline{U}$ least energy solution of \eqref{eq_least_energy_m0}},}$$
where the limits are given by $n\to +\infty$. Furthermore, by the uniform estimates in $L^{\infty}(\R^N)$ and the local uniform estimates in $C^{0, \sigma}_{loc}(\R^N)$ of $u_{\eps_n}$, together with the locally-compact version of Ascoli-Arzel\`{a} theorem, we have that the last convergence is indeed uniform on compacts. 
Thus, recalled that $v_{\eps_n}= u_{\eps_n}(\cdot/\eps_n)$ are solutions of the original problem \eqref{eq_concent_genericaf}, defined $x_{\eps_n}:=\eps_n y_{\eps_n}$ we obtain, as $n \to +\infty$,
 $$\parag{&v_{\eps_n}(x_{\eps_n})=\max_{\R^N} v_{\eps_n}, &\\& 
 d(x_{\eps_n}, K)\to 0, &\\& 
 \frac{C'}{1+|\frac{x-x_{\eps_n}}{\eps_n}|^{N+2s}}\leq v_{\eps_n}(x) \leq \frac{C''}{1+|\frac{x-x_{\eps_n}}{\eps_n}|^{N+2s}}, \quad \hbox{ for $x \in \R^N$}, &\\& 
 \norm{v_{\eps_n}(\eps_n \cdot + x_{\eps_n}) - \overline{U}}_{X} \to 0, \quad X= H^s(\R^N) \hbox{ and } X=L^{\infty}_{loc}(\R^N), }$$
for some $\overline{U}$ least energy solution of \eqref{eq_least_energy_m0}. This concludes the proof.
\QED

\begin{Remark}\label{rem_regolar}
We observe that Steps 2 and 3 apply to a whatever family of equations $(u_{\eps})_{\eps>0}$, that is why the regularity statement in Theorem \ref{teo_concen_esist} holds true. 
Moreover, the uniform concentration in $K$ and the uniform polynomial decay are obtained by a contradiction argument.
\end{Remark}

%%%%%%%%%%%%%
\subsubsection{Proof of Lemma \ref{lemma_stima_unif}: polynomial decay of $\widehat{S}$}

By adapting some argument of the proof of Theorem \ref{teo_concen_conc} we can now complete the proof of Lemma \ref{lemma_stima_unif}.

\begin{Proposition}[Polynomial decay]
\label{prop_polyn_dec}
Assume \hyperref[(f1)]{\textnormal{(f1)}}--\hyperref[(f3)]{\textnormal{(f3)}}. Let $a>0$ and let $U$ be a weak solution of
$$(-\Delta)^s U + a U = f(U), \quad x \in \R^N.$$
Then there exist positive constants $C_a', C''_a$ such that
$$\frac{ C_a'}{1+|x|^{N+2s}} \leq U(x)\leq \frac{ C_a''}{1+|x|^{N+2s}},\quad \textit{ for $x \in \R^N$}.$$
These constants can be chosen uniform for $U \in \widehat{S}$.
\end{Proposition}

\claim Proof. 
The proof is similar to the one carried out in Theorem \ref{teo_concen_conc}.

Indeed, as in Step 2 and Step 3, we obtain the uniform boundedness in $L^{\infty}(\R^N)$. We point out that the values $C_{\delta}$, $H$, $C$ and $R_0$ depend on $a \in [m_0, m_0+\nu_0]$, since they depend on $\delta$ and we must have $\delta < a$; on the other hand, it is sufficient to take $\delta < m_0$ to gain uniformity. 
The same can be said on the uniform boundedness in $C^{0, \sigma}(\R^N)$ and for the constants $R'', C', C''$ related to the comparison functions $\underline{W}', \overline{W}'$, thanks to Lemma \ref{lem_esist_sol_part}. 
As we will show, this allows us to gain that
\begin{equation}\label{eq_conv_unif_0}
\lim_{|x|\to +\infty} U(x)=0 \quad \hbox{ uniformly for $U\in \widehat{S}$},
\end{equation}
which leads, as in Step 6 of the proof, to
$$\abs{f(U(x))} \leq \frac{1}{2}a U(x),\quad \hbox{ for $|x|>R'$}$$
where $R'$ does not depend on $a\in [m_0, m_0+\nu_0]$. In addition, compactness of $\widehat{S}$ and a simple contradiction argument lead to $\min_{B_R} U \geq C>0$ uniformly for $U\in \widehat{S}$. If we prove \eqref{eq_conv_unif_0}, we conclude as in Step 6. 

Let us prove \eqref{eq_conv_unif_0}. By contradiction, there exist $(x_k)_k\subset \R^N$, $|x_k|\to +\infty$, $(U_k)_k \subset \widehat{S}$ and $\theta>0$ such that $U_k(x_k)> \theta>0$. Define
$$V_k:= U_k(\cdot+x_k).$$
Since both are bounded sequences in $H^s(\R^N)$, we have $U_k \wto U$ and $V_k \wto V$ in $H^s(\R^N)$; moreover, by the uniform $L^{\infty}(\R^N)$ and $C^{0,\sigma}_{loc}(\R^N)$ estimates and Ascoli-Arzel\`{a} theorem, we have also that the convergences are pointwise. In particular by
$$U_k(0)\geq U_k(x_k) > \theta, \quad V_k(0)=U_k(x_k)>\theta$$
we obtain
$$U(0)\geq \theta >0, \quad V(0) \geq \theta>0.$$
As a consequence, $U$ and $V$ are not trivial. Let now $(a_k)_k \subset \R$ be such that $U_k \in S_{a_k}$; up to a subsequence we have $a_k \to a\in [m_0, m_0+\nu_0]$. Observed that also $V_k$ are solutions of $L'_{a_k}(V_k)=0$, we obtain, as in Step 3 of Lemma \ref{lemma_stima_unif} (see also Step 7 of the proof of Lemma \ref{lemma_stima_basso}), that $U$ and $V$ are (nontrivial) solutions of $L_a'(U)=0$. Hence
$$E_{m_0} \leq E_a \leq L_a(U), \quad E_{m_0} \leq E_a\leq L_a(V).$$
By the Pohozaev identity (applied to $U_k$) we have the following chain of inequalities, once fixed $R>0$ and $k\gg 0$ such that $|x_k|\geq 2R$,
\begin{align*}
l_0 &\geq \liminf_{k\to +\infty} L_{a_k}(U_k) = \frac{s}{N} \liminf_{k\to +\infty} \int_{\R^N}\abs{ (-\Delta)^{s/2}U_k}^2 dx \\
&\geq \frac{s}{N} \liminf_{k\to +\infty} \Bigg( \int_{B_R} \abs{ (-\Delta)^{s/2}U_k}^2 dx + \int_{B_R} \abs{ (-\Delta)^{s/2}V_k}^2 dy \Bigg) \\
&\geq \frac{s}{N} \Bigg( \int_{B_R} \abs{ (-\Delta)^{s/2}U}^2 dx + \int_{B_R} \abs{ (-\Delta)^{s/2}V}^2 dy\Bigg)
\end{align*}
where in the last passage we have used that $U_k \wto U$ in $H^s(\R^N)$, thus $(-\Delta)^{s/2} U_k \wto (-\Delta)^{s/2} U$ in $L^2(\R^N)$, hence (by restriction) in $L^2(B_R)$, and the weak lower semicontinuity of the norm. 
Thus, by choosing $R$ sufficiently large, we have, again by the Pohozaev identity (applied to $U$ and $V$, we use \hyperref[(f3)]{\textnormal{(f3)}})
\begin{align*}
l_0 &\geq \frac{s}{N} \Bigg( \int_{\R^N} \abs{ (-\Delta)^{s/2}U}^2 dx+ \int_{\R^N} \abs{ (-\Delta)^{s/2}V}^2 dy\Bigg) - \eta\\
&= L_a(U) + L_a(V) - \eta 
\geq 2E_{m_0}- \eta
\end{align*}
which leads to a contradiction if we choose $\eta \in (0, 2E_{m_0} - l_0)$, possible thanks to \eqref{eq_stima_l0}.
\QED

\medskip

\begin{Remark}\label{rem_remove_Pohozaev}
Actually, \hyperref[(f3)]{\textnormal{(f3)}} can be dropped, and we highlight here some modifications to the previous proofs.
\begin{itemize}
\item Define %for $a \in [m_0, b]$ (where $E_b:= l_0$), 
$S_a:=\{ U \in H^s_r(\R^N) \setminus \{0\} \mid L_a'(U)=0, \; L_a(U) \leq l_0, \; P_a(U)=1\}$.

We comment the proof of the compactness (Lemma \ref{lemma_stima_unif}).
The nonemptiness si given by the existence of a ground state with $L_a(U) = E_a \leq l_0$, which is automatically radially symmetric. 
The boundedness of $\widehat{S}$ is given by the extra condition on the Pohozaev; compactness is now enduced by the radial symmetry (and the fact that $\int_{\R^N} g(u_n) u_n \to \int_{\R^N} g(u) u$, see Proposition \ref{prop_converg_generiche_loc}), and the strong convergence implies that the Pohozaev identity is preserved. 

Finally, the proof of uniform asymptotic decay (Proposition \ref{prop_polyn_dec}) is modified in the following way: after having shown that, for $U_k \in \widehat{S}$ and $|x_k|\to +\infty$, $V_k=U_k(\cdot + x_k)$ satisfies $V_k \wto V \nequiv 0$ (in $H^s(\R^N)$, thus in $L^p(\R^N)$, $p \in (2, 2^*_s)$), while $U_k \wto U$ in $H^s_r(\R^N)$, by the compactness we have $U_k \to U$ in $L^p(\R^N)$, $p \in (2, 2^*_s)$. Thus, for every $\varphi \in L^{p'}(\R^N)$ we have $\varphi(\cdot-x_k) \wto 0 \in L^{p'}(\R^N)$ and hence
$$\int_{\R^N} V_k \varphi = \int_{\R^N} U_k \varphi(\cdot - x_k) \to 0$$
i.e. $V_k \wto 0 $ in $L^p(\R^N)$, thus $V \equiv 0$, impossible.
\item Lemma \ref{lemma_stima_g}, Lemma \ref{lemma_stima_g0} (and whenever the Pohozaev identity is used for $\widehat{S}$), can be proved thanks to the extra condition in $S_a$.
\end{itemize}
Her we kept the original definition of $S_a$ (i.e. $U \in H^s(\R^N)$ such that $\max U = U(0)$), since this approach can be adapted also to frameworks where radial symmetry is not a feature of the limiting problem.
\end{Remark}

%%%%%%%%%%%%%%%%%%%%%%%%%%%%%%%%%%%%%%%%%%%%%%%%%%%%%%%
%%%%%%%%%%%%%%%%%%%%%%%%%%%%%%%%%%%%%%%%%%%%%%%%%%%%%%%

\section{The critical case}
\label{sec_conc_critical}

Goal of this Section is to study equation \eqref{eq_concent_genericaf}, that is %the singularly perturbed nonlinear Schr\"odinger equation 
$$ \varepsilon^{2s}(- \Delta)^s v+ V(x) v= f(v), \quad x \in \mathbb{R}^N,$$
where %$s \in (0,1)$, $N \geq 2$, $V \in C(\mathbb{R}^N,\mathbb{R})$ is a positive potential and 
now $f$ is assumed critical and satisfying general Berestycki-Lions type conditions. 
When $\eps>0$ is small, we obtain again existence and multiplicity of semiclassical solutions, relating the number of solutions to the cup-length of the set of local minima of $V$; these solutions are proved to concentrate in the potential well, exhibiting a polynomial decay. 
In particular, we improve the result in \cite{HeZo}.
%Furthermore,
Finally, we prove the previous results also in the limiting local setting $s=1$ and $N\geq 3$, with an exponential decay of the solutions.

\bigskip

Here, thus, we assume \hyperref[(V1)]{\textnormal{(V1)}}-\hyperref[(V2)]{\textnormal{(V2)}} %, that we recall here
%Namely, we focus on possibly degenerate local minima of $V$, that is $V$ satisfies
%\begin{itemize}
%\item[(V1)] $ V\in C(\R^N, \R)\cap L^{\infty}(\R^N)$, $\underline{V}=\inf_{\R^N} V>0$,
%\item[(V2)] there exists a bounded domain $\Omega\subset \R^N$ such that
%$$m_0= \inf_{\Omega} V < \inf_{\partial \Omega} V,$$
%with set of local minima
%\begin{equation}\label{eq_def_K}
%K=\{ x \in \Omega \mid V(x)=m_0\},
%\end{equation}
%\end{itemize}
where we recall
$$m_0= \inf_{\Omega} V$$
with %set of local minima
\begin{equation}\label{eq_def_K_2}
K=\{ x \in \Omega \mid V(x)=m_0\},
\end{equation}
and \hyperref[(f1)]{\textnormal{(f1)}}-\hyperref[(f3)]{\textnormal{(f3)}}, where now \hyperref[(f1)]{\textnormal{(f1.3)}} is substituted with a critical (not pure) growth, i.e.%we assume general Berestycki-Lions assumptions on $f$, i.e.
%\begin{itemize}
%\item[(f1)] $f\in C(\R, \R)$, and 
%$f\in C^{0,\gamma}_{loc}(\R, \R)$ for some $\gamma \in (1-2s, 1)$ if $s\in (0,1/2]$, \<
%item[(f2)] $f(t)\equiv 0$ for $t \leq 0$,
%\item[(f3)] $\, \lim_{t \to 0} \frac{f(t)}{t}=0$,
%\item[(f4)] $\lim_{t \to +\infty} \frac{f(t)}{t^{2^*_s-1}}=a>0$, where $2^*_s=\frac{2N}{N-2s}$,
%\item[(f5)] for some $C>0$ and $\max\{2^*_s-2, 2\} <p< 2^*_s$, i.e. satisfying
%\begin{equation}\label{eq_cond_p}
%p \in \parag{&\Big(\frac{4s}{N-2s}, \frac{2N}{N-2s}\Big)& \quad N \in (2s, 4s), \\ &\Big(2, \frac{2N}{N-2s}\Big)& \quad N\geq 4s, }
%\end{equation}
%it results that
%$$f(t) \geq a t^{2^*_s-1} + C t^{p-1} \quad \hbox{for $t \geq 0$}.$$
%\end{itemize}
\begin{itemize}
	\item[(f1')] \label{(f1')}
Berestycki-Lions type assumptions with respect to $m_0>0$, that is
	\begin{itemize}
		\item[(f1.1)] $\, f\in C(\R, \R)$;
		\item[(f1.2)] $\, \lim_{t \to 0} \frac{f(t)}{t}=0$;
		\item[(f1.3')] $\, \lim_{t \to +\infty} \frac{f(t)}{t^{2^*_s-1}}=a>0$, where $2^*_s=\frac{2N}{N-2s}$, and moreover
for some $C>0$ and $\max\{2^*_s-2, 2\} <p< 2^*_s$, i.e. satisfying
\begin{equation}\label{eq_cond_p}
p \in \parag{&\Big(\frac{4s}{N-2s}, \frac{2N}{N-2s}\Big)& \quad N \in (2s, 4s), \\ &\Big(2, \frac{2N}{N-2s}\Big)& \quad N\geq 4s, }
\end{equation}
it results that
$$f(t) \geq a t^{2^*_s-1} + C t^{p-1} \quad \hbox{for $t \geq 0$};$$
		\item[(f1.4)] $\, F(t_0)> \frac{1}{2} m_0 t_0^2$ for some $t_0>0$. %, where $F(t):= \int_0^t f(s) ds$;
	\end{itemize}
%	\item[(f2)] $f(t)\equiv 0$ for $t\leq 0$.
%	\item[(f3)] $f\in C^{0,\gamma}_{loc}(\R)$ for some $\gamma \in (1-2s, 1)$ if $s\in (0,1/2]$.
\end{itemize}

%\tor{Caso sottocritico non strettamente? VEDI}

See also Remark \ref{rem_ipotesi_reg} for some weakening and comments on the assumptions \hyperref[(V1)]{\textnormal{(V1)}}, \hyperref[(f1')]{\textnormal{(f1.3')}} and \hyperref[(f3)]{\textnormal{(f3)}}. %(f1) and (f5). 
Notice that the stronger condition on $p$ in the first line of \eqref{eq_cond_p} is verified, whenever $N\geq 2$, only if $N=2$ and $s \in (\frac{1}{2}, 1]$, or $N=3$ and $s\in (\frac{3}{4}, 1]$. 
We point out that the condition $C>0$ in \hyperref[(f1')]{\textnormal{(f1.3')}} %(f4) 
is of key importance: indeed, for pure critical nonlinearities of the type 
$f(t)=|t|^{2^*_s-2}t,$
the limiting problem \eqref{eq_limite_iniz}, that is
$$(-\Delta)^s u + m_0 u = |u|^{2^*_s-2} u, \quad x \in \R^N$$
does not admit any variational solution \cite{DSS1}. % exhibits a quite different scenario. 

\smallskip

The existence of a solution in a critical, fractional setting, in the case of local minima \hyperref[(V1)]{\textnormal{(V1)}}-\hyperref[(V2)]{\textnormal{(V2)}} and general Berestycki-Lions assumptions \hyperref[(f1')]{\textnormal{(f1')}}-\hyperref[(f2)]{\textnormal{(f2)}}-\hyperref[(f3)]{\textnormal{(f3)}}, %(f2)--(f5),
 has been faced in \cite{JLZ} by assuming $V \in C^1(\R^N)$, and moreover in \cite{He2} by means of penalization methods.

Inspired by \cite{Rab1}, multiplicity of solutions of \eqref{eq_concent_genericaf} in the case of global minima of $V$ was studied in \cite{SZY} for power-type nonlinearities. Moreover, in \cite{LTZW} the authors consider functions of the type
\begin{equation}\label{eq_cond_g}
f(t)=g(t)+ |t|^{2^*_s-2}t,
\end{equation}
where $g$ is subcritical and satisfies a monotonicity condition which allows to implement the Nehari manifold tool, and they relate the number of solutions to the Lusternik-Schnirelmann category of the set of global minima.

Existence of multiple solutions for local minima of $V$ has been investigated, in the spirit of \cite{DF0}, by \cite{HeZo} with sources of the type \eqref{eq_cond_g}, where now $g$ satisfies also an Ambrosetti-Rabinowitz condition: this assumption enables to employ Mountain Pass and Palais-Smale arguments, combined with a penalization scheme. 
Again, the authors are able to find $\cat(K)$ solutions, where $K$ is the set of local minima of $V$ and $\cat(K)$ denotes its Lusternik-Schnirelmann category.

\smallskip

In the present Section we prove a multiplicity result for equation \eqref{eq_concent_genericaf} under almost optimal assumptions of $f$, showing the concentration of the solutions around local minima of $V$.

\smallskip

In particular, we prove the following result.
\begin{Theorem}\label{teo_main}
Assume $s \in (0,1)$, $N\geq 2$ and that \hyperref[(V1)]{\textnormal{(V1)}}-\hyperref[(V2)]{\textnormal{(V2)}}, \hyperref[(f1')]{\textnormal{(f1')}}-\hyperref[(f2)]{\textnormal{(f2)}}-\hyperref[(f3)]{\textnormal{(f3)}} %(f1)--(f5)} 
hold. 		
Let $K$ be defined by \eqref{eq_def_K_2}.
Then, for small $\eps>0$ equation \eqref{eq_concent_genericaf} has at least $\cupl(K)+1$ positive solutions, which belong to $C^{0, \sigma}(\R^N) \cap L^{\infty}(\R^N)$ for some $\sigma \in (0,1)$. 
Moreover, each of these sequences $v_{\eps}$ concentrates in $K$ as $\eps \to 0$: namely, there exist $x_{\eps} \in \R^N$ global maximum points of $v_{\eps}$, such that
$$\lim_{\eps\to 0}d(x_{\eps}, K) =0$$
and
$$\frac{C'}{1+|\frac{x-x_{\eps}}{\eps}|^{N+2s}}\leq v_{\eps}(x) \leq \frac{C''}{1+|\frac{x-x_{\eps}}{\eps}|^{N+2s}} \quad \hbox{ for $x \in \R^N$}$$
where $C', C''>0$ are uniform in $\eps>0$. 
Finally, for every sequence $\eps_n \to 0^+$ there exist a ground state solution $U$ of \eqref{eq_limite} and a point $x_0 \in K$ such that, up to a subsequence,
$$x_{\eps_n} \to x_0 \in K$$
and
$$v_{\eps_n}(\eps_n \cdot + x_{\eps_n}) \to U \quad \hbox{as $n\to +\infty$}$$
in $H^s(\R^N)$ and locally on compact sets.
\end{Theorem}

We highlight that Theorem \ref{teo_main} extends the existence results in \cite{He2, LTZW} to a multiplicity result, and it improves the multiplicity theorem in \cite{HeZo}, since we do not assume monotonicity nor Ambrosetti-Rabinowitz conditions on the nonlinearity. Moreover, no nondegeneracy and global conditions on $V$ are considered.

\begin{Remark}\label{rem_ipotesi_reg}
As observed in Remark \ref{rem_ipotesi_V}, assumption \hyperref[(V1)]{\textnormal{(V1)}} in Theorem \ref{teo_main} can be relaxed without assuming the boundedness of $V$ (see also \cite{BJ0, BT1}). 
Moreover, the condition $$p>\max\{2^*_s-2, 2\}$$ in \hyperref[(f1')]{\textnormal{(f1.3')}} %(f5)
can be relaxed in $p>2$ by paying the cost of considering a sufficiently large $C\gg 0$; see for instance \cite{SZ2, He2}. 
Finally, we remark that \hyperref[(f3)]{\textnormal{(f3)}}, %(f1.1)%(f1)
 instead of the mere continuity of $f$, is needed only to get a Pohozaev identity by means of the regularity of solutions (see Proposition \ref{Pohozaev-prop}). 
 See also Remark \ref{rem_caso_subcr_nonstr} for further comments.
\end{Remark}

\smallskip

The idea of the present Section is the following: first, we gain compactness and uniform $L^{\infty}$-bounds on the set of ground states of the critical limiting problem \eqref{eq_limite_iniz}; to this aim we employ a Moser's iteration argument adapted to the fractional framework, without the use of the $s$-harmonic extension, and appropriate for weak solutions (see Proposition \ref{prop_prop_u_Linf}). 
The criticality of the problem, as well as the absence of a chain rule, make the argument more delicate. 
The gained uniformity allows then the introduction of a suitable truncation on the nonlinearity $f$; the new truncated function reveals thus to be subcritical.

Therefore, we can apply to the truncated problem the approach of the previous Sections: we employ a penalization argument on a neighborhood of expected solutions, perturbation of the ground states of a limiting problem, and this neighborhood results to be invariant under the action of a deformation flow. 
Compactness is restored also by the use of the new fractional center of mass, which engages the new strong seminorm; % stronger than the usual Gagliardo one; 
the topological machinery between two level sets of the associated indefinite energy functional is then built also through the use of the Pohozaev functional. 
The number of solutions is thus related to the cup-length of $K$ and these solutions are proved to exhibit a polynomial decay and to converge to a ground state of the limiting equation. 
This last convergence allows finally to prove that these solutions solve the original critical problem \eqref{eq_concent_genericaf}.

We point out that the techniques employed in the previous Sections cannot be applied directly to the critical framework: indeed, the embedding of $H^s(\R^N)$ in $L^{2^*_s}(\R^N)$ is not compact, even if we reduce to radially symmetric functions or to bounded domains; in particular, the criticality obstructs the convergence of truncated Palais-Smale sequences related to the penalized functional, which is a key point in the proof. 
Moreover, the regularity results given by \cite{Coz0}, exploited in the concentration and in the decay of the solutions, do not apply; in particular, $L^{\infty}$-bounds and compactness of the set of ground states of the limiting problem have to be specifically investigated.

\medskip

We highlight that the conclusions of Theorem \ref{teo_main} hold also for $s=1$ and $N\geq 3$, as we state in Theorem \ref{teo_main_locale}. 
Regarding this local framework, Theorem \ref{teo_main_locale} is the critical counterpart of the result in \cite{CJT}: again, we point out that the arguments exploited in the subcritical setting of \cite{CJT} cannot be directly implemented in our framework, because of the lack of compactness. 
In the critical case, previous results were given by \cite{ADS, ZhZo2, Amb4}: in particular we extend here the existence result in \cite{ZCZ} to a multiplicity result, and we improve the multiplicity theorem in \cite{WLXF} in the sense that we do not need to work with global minima of $V$ nor we need monotonicity on $f$.
In this setting, the solutions decay exponentially and enjoy more regularity. 
Notice that in such a case \hyperref[(f3)]{\textnormal{(f3)}} is no more needed. % (f1) means $f$ merely continuous.

\medskip

This last part of the Chapter % Section 
is organized as follows. 
%\tr{aggiusta label sezioni} %COMMENT NOW
In Section \ref{sec_recalls} %we recall some notions on the fractional Sobolev space, and then 
we obtain compactness of the set of ground states and a crucial $L^{\infty}$-bound on the critical limiting problem. 
In Section \ref{sec_truncated} we use this uniform estimate to introduce a truncation which brings the problem back to the subcritical case, and we prove Theorem \ref{teo_main}. Finally, in Section \ref{sec_locale} we deal with the local case.

%%%%%%%%%%%%%%%%%%%%%%%%%%%%%%%%%%%%%%%%%%%%%%%%%

\subsection{Uniform $L^{\infty}$-bound}
\label{sec_recalls}

\medskip

\noindent Let us recall some crucial results on the limiting critical problem \eqref{eq_limite_iniz}, that is
\begin{equation}\label{eq_limite} 
(-\Delta)^s U + m_0 U = f(U), \quad x \in \R^N.
\end{equation}
We recall the energy %$C^1$-functional
 $\mc{L}: H^s(\R^N) \to \R$
$$\mc{L}(U):=\frac{1}{2} \int_{\R^N} |(-\Delta)^{s/2}U|^2 \, dx + \frac{ m_0}{2} \int_{\R^N} U^2 \, dx - \int_{\R^N} F(U) \, dx, \quad U \in H^s(\R^N),$$
the related least energy % \emph{least energy}
$$E_{m}:=\inf \big\{ \mc{L}(U) \mid U\in H^s(\R^N)\setminus\{0\}, \; \mc{L}'(U)=0\big\},$$
%Moreover we define 
and the Mountain Pass level %\emph{Mountain Pass level} 
$$C_{mp}:=\inf_{\gamma \in \Gamma} \sup_{t \in [0,1]} \mc{L}(\gamma(t))$$
with
$$\Gamma:=\big \{ \gamma \in C\big([0,1], \, H^s(\R^N)\big) \mid \gamma(0)=0, \; \mc{L}(\gamma(1))<0\big \}.$$
We introduce also the following \emph{minimization problem}
\begin{equation}\label{eq_min_probl}
C_{min}:=\inf \big\{ \mc{T}(U) \mid U \in H^s(\R^N), \; \mc{V}(U)=1\big\}
\end{equation}
where
$$\mc{T}(U):= \int_{\R^N} |(-\Delta)^{s/2}U|^2 \, dx, \quad \mc{V}(U):=\int_{\R^N} \left(F(U)-\frac{m_0}{2}U^2\right)\, dx.$$
Notice that $\mc{L} = \frac{1}{2} \mc{T} - \mc{V}$. The following collection of results states the equivalence of the previous problems and the existence of a solution.

%Uno studio su di una generalizzazione dell'ensemble di Wishart

\smallskip

\begin{Proposition}\label{prop_esist_ground}
Assume \hyperref[(f1')]{\textnormal{(f1')}}-\hyperref[(f2)]{\textnormal{(f2)}}-\hyperref[(f3)]{\textnormal{(f3)}}. Then
there exists a ground state solution for the problem \eqref{eq_limite}, that is a function $U$ which solves the equation and such that
$$\mc{L}(U)=E_{m}.$$
Moreover, every ground state is also a Mountain Pass solution and (up to scaling) also a solution for the minimization problem \eqref{eq_min_probl}, and viceversa; in addition the following relations hold
$$ E_{m}=C_{mp},$$
\begin{equation}\label{eq_rel_ener_cmin}
E_{m}= \frac{s}{N} (2^*_s)^{-\frac{N}{2^*_s s}} (C_{min})^{\frac{N}{2s}},
\end{equation}
and every ground state is positive. 
Finally, recalled that $\mc{S}$ is the best Sobolev constant for the embedding \eqref{eq_embd_homog}, we have that the following upper bound holds
\begin{equation}\label{eq_rel_best_const}
C_{min} < \left(\frac{2^*_s}{a}\right)^{\frac{2}{2^*_s}} \mc{S}
\end{equation}
where $a>0$ appears in assumption \hyperref[(f1')]{\textnormal{(f1.3')}}. %\textnormal{(f4)-(f5)}.
\end{Proposition}

\claim Proof.
The positivity is a straightforward consequence of assumption \hyperref[(f2)]{\textnormal{(f2)}}. Existence of a ground state solution can be achieved through the use of \eqref{eq_rel_best_const} and minimization of $C_{min}$ as classically made by \cite{BL1} (see also \cite[Lemma 1]{BJM0}). 
The equivalence with the Mountain Pass formulation is instead discussed as in \cite{JT0}. 
We refer to \cite[Proposition 2.4 and Remark 1.3]{JLZ} for the precise statement and to \cite[Section 4.1 and Remark 1.2]{ZMS}, \cite[Section 2]{LMS} for details.

Moreover, as observed in Remark \ref{rem_ipotesi_reg}, to get the existence of a ground state, the restriction on the range of $p$ in assumption \hyperref[(f1')]{\textnormal{(f1.3')}} % (f5) 
can be substituted, by arguing as in \cite[Lemma 3.3]{SZY}, with the request that $C$ is sufficiently large (see also \cite[Proposition 2.8]{He2} and references therein).

We refer also to \cite[Theorem 3.1.3, Theorem 3.1.5]{Amb5}.
\QED

\bigskip

Thanks to Proposition \ref{prop_esist_ground} we can define
$$\widehat{S}:= \big\{ U \in H^s(\R^N)\setminus \{0\} \mid \hbox{$U$ ground state solution of \eqref{eq_limite}, $ U(0)=\max_{\R^N} U$}\big\}.$$

We observe that, by the fractional version of the P\'olya-Szeg\H o inequality \cite{Par0}, every minimizer of $C_{min}$ (i.e. every ground states of \eqref{eq_limite}) is actually radially symmetric decreasing up to a translation (see also Remark \ref{rem_byeon_addit} %Theorem \ref{Byeonminimizers} 
and \cite[Proposition B.3]{BKS}). 
Thus, the request in $\widehat{S}$ for $U$ to have a maximum in zero is equivalent to the radial symmetry of $U$; that is
\begin{equation}\label{eq_S_radial}
\widehat{S}= \big\{ U \in H^s(\R^N)\setminus \{0\} \mid \hbox{$U$ radially symmetric ground state solution of \eqref{eq_limite}}\big\}.
\end{equation}

\begin{Proposition}\label{prop_compat_S}
Every $U\in \widehat{S}$ satisfies the Pohozaev identity, i.e.
\begin{equation}\label{eq_pohozaev}
 \int_{\R^N} |(-\Delta)^{s/2}U|^2 \, dx - 2^*_s \int_{\R^N} \left(F(U)-\frac{m_0}{2} U^2 \right)\, dx=0.
\end{equation}
Moreover, the set $\widehat{S}$ is compact.
\end{Proposition}

\claim Proof.
Once one observes that $U \in L^{\infty}(\R^N)$, which follows from Proposition \ref{prop_prop_u_Linf}, 
%ose proof is an easy adaptation of Proposition \ref{prop_upperbound} below (focusing on a single $U\in \widehat{S}$), 
the proof of \eqref{eq_pohozaev} is gained by means of regularity results and explicit computations on the $s$-harmonic extension problem; the arguments can be easily adapted from \cite[Proposition 1.1]{BKS} to the critical case.

\smallskip

Let us show the boundedness of $\widehat{S}$. For any $U \in \widehat{S}$, the embedding \eqref{eq_embd_homog} and the Pohozaev identity \eqref{eq_pohozaev} lead to
$$\norm{U}_{2^*_s} \leq \mc{S}^{-\frac{1}{2}} \norm{(-\Delta)^{s/2} U}_2 = \mc{S}^{-\frac{1}{2}} \frac{N}{s} \mc{L}(U) = \mc{S}^{-\frac{1}{2}} \frac{N}{s} E_{m};$$
moreover equation \eqref{eq_limite} and assumption \hyperref[(f1')]{\textnormal{(f1')}} %(f3)-(f4) imply
imply
$$\norm{(-\Delta)^{s/2}U}_2^2 + m_0 \norm{U}_2^2 = \int_{\R^N} f(U) U \, dx \leq \delta \norm{U}_2^2 + C_{\delta} \norm{U}_{2^*_s}^{2^*_s}$$
for $\delta < m_0$ and some $C_{\delta}>0$. The combination of the two bounds leads to the claim.

\smallskip

Let thus focus on compactness; we use some ideas from \cite{ZhZo1}. 
Let $U_n$ be a sequence in $\widehat{S}$; by \eqref{eq_S_radial} we assume $(U_n)_n \subset H^s_r(\R^N)$, where %, where $H^s_r(\R^N)$ denotes the subset of $H^s(\R^N)$ consisting in radially symmetric functions. 
%We recall that 
$H^s_r(\R^N) \hookrightarrow \hookrightarrow L^q(\R^N)$ for $q \in (2, 2^*_s)$. By the boundedness of $\widehat{S}$ we can assume $U_n \wto U$ in $H^s_r(\R^N)$. 
Set 
$$\sigma:= \left( \frac{1}{2^*_s} C_{min}\right)^{\frac{1}{2s}}$$
and
$$V_n := U_n(\sigma \cdot ), \quad V:= U(\sigma \cdot)$$
we have, by exploiting the Pohozaev identity, that $V_n$ are solutions of the minimization problem \eqref{eq_min_probl}, that is
$$\mc{T}(V_n)=C_{min}, \quad \mc{V}(V_n)=1.$$
Thus we have $V_n \wto V$ in $H^s_r(\R^N)$, and hence $V_n \to V$ in $L^q(\R^N)$, $q \in (2, 2^*_s)$, and $V_n \to V$ almost everywhere. By the lower semicontinuity of the norm we obtain
\begin{equation}\label{eq_dim_stima_T}
\mc{T}(V) \leq C_{min};
\end{equation}
hence, to conclude the proof, it is sufficient to show that $\mc{V}(V)=1$, since this implies also that $U=V(\sigma^{-1} \cdot)$ lies in $\widehat{S}$.

Set
$$W_n:=V_n-V$$
we have by the Brezis-Lieb Lemma (since $(-\Delta)^{s/2} V_n \wto (-\Delta)^{s/2} V$ in the Hilbert space $L^2(\R^N)$) 
\begin{align}
 \mc{T}(W_n) &= \mc{T}(V_n) - \mc{T}(V) + o(1) \nonumber \\
 &= C_{min}- \mc{T}(V) + o(1) \label{eq_dim_T}\\
 &\leq C_{min} + o(1). \label{eq_dim_2_T}
\end{align}
Moreover, rewrite $\mc{V}(W_n)$ as
\begin{equation}\label{eq_dim_riscr_V}
\mc{V}(W_n) = \int_{\R^N} \left(F(W_n) - \frac{a}{2^*_s} W_n^2\right) \, dx+ \frac{a}{2^*_s} \norm{W_n}_{2^*_s}^{2^*_s} - \frac{m_0}{2} \norm{W_n}_2^2.
\end{equation}
Again by the Brezis-Lieb Lemma (since $V_n \wto V$ in $L^q(\R^N)$, $q=2, 2^*_s$ and $V_n \to V$ almost everywhere) we have
\begin{equation}\label{eq_dim_norm}
\norm{W_n}_q^q = \norm{V_n}_q^q - \norm{V}_q^q + o(1), \quad q=2, 2^*_s.
\end{equation}
Set
$$g(t):= f(t)-a t^{2^*_s-1}$$
we have that $g$ is subcritical at infinity by \hyperref[(f1')]{\textnormal{(f1.3')}}, %(f4), 
and superlinear in zero by \hyperref[(f1')]{\textnormal{(f1.2)}}; %(f3); 
thus, set $G(t):= \int_0^t g(\tau)d \tau$, by Proposition \ref{prop_converg_generiche_loc} we have
\begin{equation}\label{eq_dim_G}
\int_{\R^N} G(W_n) \, dx= o(1), \quad \int_{\R^N} G(V_n) \, dx= \int_{\R^N} G(V) \, dx+ o(1).
\end{equation}
Therefore by \eqref{eq_dim_riscr_V}--\eqref{eq_dim_G} we obtain
\begin{align}
\mc{V}(W_n) &= \mc{V}(V_n) - \mc{V}(V) + o(1) \nonumber\\
&= 1- \mc{V}(V) + o(1). \label{eq_dim_V}
\end{align}
Finally, through a simple scaling argument, we observe that
\begin{equation}\label{eq_dim_stima_gen}
\mc{T}(u) \geq C_{min} (\mc{V}(u))^{\frac{2}{2^*_s}} \quad \hbox{for every $\mc{V}(u)\geq 0$}.
\end{equation}
We pass to prove that $\mc{V}(V)=1$ by contradiction.
\\ \textbf{Case $\mc{V}(V)>1$.} In this case, by \eqref{eq_dim_stima_gen} we have
$$\mc{T}(V) \geq C_{min} (\mc{V}(V))^{\frac{2}{2^*_s}} > C_{min}$$
which contradicts \eqref{eq_dim_stima_T}.
\\ \textbf{Case $\mc{V}(V)<0$.} Then, by \eqref{eq_dim_V} we have that $$\mc{V}(W_n) \geq 1- \frac{1}{2} \mc{V}(V) >1 \quad \hbox{for $n\gg 0$}.$$ Thus, by \eqref{eq_dim_stima_gen} we obtain
$$\mc{T}(W_n) \geq C_{min}(\mc{V}(W_n))^{\frac{2}{2^*_s}} \geq C_{min} \left( 1- \frac{1}{2} \mc{V}(V)\right)^{\frac{2}{2^*_s}}$$
which contradicts \eqref{eq_dim_2_T}.
\\ \textbf{Case $\mc{V}(V) \in (0,1)$.} Again by \eqref{eq_dim_V} we have that $$\mc{V}(W_n) \geq \frac{1}{2} \left(1- \mc{V}(V)\right) > 0 \quad \hbox{for $n\gg 0$.}$$ Thus by \eqref{eq_dim_T}, \eqref{eq_dim_stima_gen} and \eqref{eq_dim_V} we gain
\begin{align*}
C_{min} &= \lim_n \big( \mc{T}(W_n) + \mc{T}(V) \big) \geq C_{min} \lim_n \left( (\mc{V}(W_n))^{\frac{2}{2^*_s}} + (\mc{V}(V))^{\frac{2}{2^*_s}}\right)\\
&= C_{min} \left( (1-\mc{V}(V))^{\frac{2}{2^*_s}} + (\mc{V}(V))^{\frac{2}{2^*_s}}\right) \\ 
&> C_{min} \big( (1-\mc{V}(V)) + \mc{V}(V)\big) = C_{min}
\end{align*}
which is an absurd.
\\ \textbf{Case $\mc{V}(V)=0$.} By \eqref{eq_dim_V} we have 
\begin{equation}\label{eq_dim_2_V}
\mc{V}(W_n)=1+o(1),
\end{equation}
and thus by \eqref{eq_dim_stima_gen} $\mc{T}(W_n) \geq C_{min}(1+o(1))^{\frac{2}{2^*_s}}$. This, combined with \eqref{eq_dim_2_T}, gives
\begin{equation}\label{eq_dim_3_T}
\mc{T}(W_n) = C_{min} + o(1).
\end{equation}
Combining \eqref{eq_dim_2_V}, \eqref{eq_dim_riscr_V} and \eqref{eq_dim_G} we obtain
$$1+o(1) = \mc{V}(W_n) = \frac{a}{2^*_s} \norm{W_n}_{2^*_s}^{2^*_s} - \frac{m_0}{2} \norm{W_n}_2^2$$
that is
\begin{align}
\norm{W_n}_{2^*_s}^{2^*_s} &=\frac{2^*_s}{a} + \frac{ 2^*_s m_0 }{2 a} \norm{W_n}_2^2 + o(1) \nonumber\\
&\geq\frac{2^*_s}{a} + o(1). \label{eq_dim_stima_crit}
\end{align}
By \eqref{eq_dim_3_T}, the Sobolev embedding \eqref{eq_embd_homog} and \eqref{eq_dim_stima_crit} we gain
$$C_{min} + o(1) = \mc{T}(W_n) = \norm{(-\Delta)^{s/2}W_n}_2^2 \geq \mc{S} \norm{W_n}_{2^*_s}^2 \geq \mc{S} \left(\frac{2^*_s}{a} + o(1) \right)^{\frac{2}{2^*_s}}.$$
Letting $n\to +\infty$ we finally have
$$C_{min} \geq \left(\frac{2^*_s}{a}\right)^{\frac{2}{2^*_s}} \mc{S}$$
which is in contradiction with \eqref{eq_rel_best_const}. This concludes the proof.
\QED

\bigskip

As a key property to employ the truncation argument, and to detect a handy neighborhood of approximating solutions, we have the following result. % and Proposition \ref{prop_prop_u_Linf}). \tor{sfoltisci da una delle due parti}
\begin{Proposition}\label{prop_upperbound}
The following bound holds
$$\sup_{U \in \widehat{S}} \norm{U}_{\infty} < \infty.$$
\end{Proposition}

\claim Proof.
Assume by contradiction that there exists $(U_n)_n \subset \widehat{S}$ such that $\norm{U_n}_{\infty} \to +\infty$ as $n \to +\infty$. 
By the compactness of $\widehat{S}$ in Proposition \ref{prop_compat_S} we may assume that $U_n$ is positive and convergent in $H^s(\R^N)$; in particular $U_n$ converges in $L^{2^*_s}(\R^N)$ and is equibounded a.$\,$e. pointwise by a function in $L^{2^*_s}(\R^N)$. If we prove that
$$\sup_n \norm{U_n}_{\infty} < +\infty$$
we get a contradiction and conclude the proof. 
In order to do this, we argue as in the proof of Proposition \ref{prop_prop_u_Linf}, uniformly in $n$ for $U_n=U_n^+$; the idea is %We use 
a Moser's iteration argument in a critical, fractional framework, appropriate for weak solutions. We refer to \cite{Gal0} for details.
\QED

%%%%%%%%%%%%%%%%%%%%%%%%%%%%%%%%%%%%%%%%%%%%%%%%%

\subsection{The truncated problem}
\label{sec_truncated}

In virtue of Proposition \ref{prop_upperbound}, let
$$M:= \sup_{U \in \widehat{S}} \norm{U}_{\infty} +1.$$
We preliminary observe that we can find a $t_0 \in [0,M]$ such that
\begin{equation}\label{eq_t0}
F(t_0) > \frac{1}{2} m_0 t_0^2.
\end{equation}
Indeed 
fixed a whatever $U\in \widehat{S}$, by the Pohozaev identity \eqref{eq_pohozaev} we have (notice that $(-\Delta)^{s/2}U$ cannot identically vanish)
$$ \int_{\R^N} \left(F(U)-\frac{m_0}{2} U^2 \right)\, dx=\frac{1}{2^*_s} \norm{(-\Delta)^{s/2}U}_2^2 >0$$
and thus there exists an $x_0 \in \R^N$ such that
$$F(U(x_0))> \frac{m_0}{2} U(x_0)^2;$$
setting $t_0:=U(x_0)\in [0,M]$ we have the claim. 

We thus set
$$k:= \sup_{t \in [0, M]} f(t) \in (0, +\infty),$$
where we observe that the strict positivity is due to the fact that $F(t_0)>0$. 
Moreover we define the \emph{truncated} nonlinearity $f_k : \R \to \R$
$$f_k(t):= \min\{ f(t), k\}, \quad t \in \R.$$
We have the following properties on $f_k: \R \to \R$:
\begin{itemize}
\item $f_k(t) \leq f(t)$ for each $t \in\R$,
\item $f_k(t)=f(t)$ whenever $|t|\leq M$,
\item $f_k(U)=f(U)$ for every $U \in \widehat{S}$.
\end{itemize}	
Notice that the same relations hold also for $F$ and $$F_k(t):= \int_0^t f_k(\tau)d \tau.$$
We have that $f_k$ is subcritical, i.e. $f_k$ satisfies assumptions \hyperref[(f1)]{\textnormal{(f1)}}--\hyperref[(f3)]{\textnormal{(f3)}};
here $p \in (1, 2^*_s-1)$ is however fixed and $t_0\in [0,M]$ is the one appearing in \eqref{eq_t0}; notice that $t_0$ does not depend on $k$. 
%
% in particular
%%(f1.1)-(f1.2)-\hyperref[(f2)]{\textnormal{(f2)}}-\hyperref[(f3)]{\textnormal{(f3)}} %(f1)--(f3) 
%together with
%\begin{itemize}
%\item[(fk4)] $\, \lim_{t \to +\infty} \frac{f_k(t)}{t^q}=0$ for some %(actually any) 
%$q \in (1, 2^*_s-1)$, 
%\item[(fk5)] $\, F_k(t_0)> \frac{1}{2} m_0 t_0^2$ for some $t_0>0$; 
%\end{itemize}
%here $q \in (1, 2^*_s-1)$ is however fixed and $t_0\in [0,M]$ is the one appearing in \eqref{eq_t0}; notice that $t_0$ does not depend on $k$. 

Consider now the truncated problem
\begin{equation}\label{eq_principale_k}
\eps^{2s} (-\Delta)^s v+ V(x) v = f_k(v), \quad x \in \R^N
\end{equation}
and the corresponding limiting truncated problem
\begin{equation}\label{eq_limite_k}
(-\Delta)^s U + m_0 U = f_k(U), \quad x \in \R^N.
\end{equation}
Notice again that, since $f_k$ satisfies \hyperref[(f2)]{\textnormal{(f2)}}, all the ground states of \eqref{eq_limite_k} are positive. Thus define
$$\widehat{S}_k:= \big\{U \in H^s(\R^N)\setminus \{0\} \mid \hbox{$U$ ground state solution of \eqref{eq_limite_k}, $ U(0)=\max_{\R^N} U$}\big\}.$$
We have that the following key relation holds.
\begin{Proposition}\label{prop_uguagl_insiem}
It results that $\widehat{S}=\widehat{S}_k$. Moreover, the least energy levels coincide.
\end{Proposition}

\claim Proof.
Let us denote by $\mc{L}_k$, $\Gamma_k$, $\mc{V}_k$, $E_{m}^k=C_{mp}^k$, $C_{min}^k$ the quantities of the truncation problem analogous to the ones introduced in Section \ref{sec_recalls} for the critical problem.

First observe that, by $\mc{L}_k \geq \mc{L}$, we have $\Gamma_k \subset \Gamma$ and
\begin{equation}\label{eq_dis_mp}
C_{mp}^k \geq C_{mp};
\end{equation}
moreover for any $V \in \widehat{S}$ we have also $\mc{L}'_k(V)=0$, and hence
\begin{equation}\label{eq_dis_min}
\min_{V \in \widehat{S}} \mc{L}_k(V) \geq \min_{\mc{L}'_k(V)=0} \mc{L}_k(V) = E_{m}^k.
\end{equation}
Let now $U \in \widehat{S}$. We have by \eqref{eq_dis_mp} and \eqref{eq_dis_min}
$$C_{mp}^k \geq C_{mp} = \mc{L}(U) = E_{m}= \min_{V \in \widehat{S}} \mc{L}(V)= \min_{V \in \widehat{S}} \mc{L}_k(V) \geq E_{m}^k.$$
Therefore
$$\mc{L}_k(U)=\mc{L}(U) = C_{mp}^k = E_{m}^k$$
which, together with $\mc{L}'_k(U)=\mc{L}'(U)=0$, gives that $U \in \widehat{S}_k$. Hence $\widehat{S}\subset \widehat{S}_k$. 
As a further consequence we gain
\begin{equation}\label{eq_energ_ugual}
 E_{m}^k = E_{m}.
 \end{equation}
We show now that $\widehat{S}_k \subset \widehat{S}$. 
By \eqref{eq_energ_ugual}, \eqref{eq_rel_ener_cmin} and the analogous relation on the subcritical problem, we have
$$C_{min}^k=C_{min},$$
thus, 
by rescaling, it is sufficient to prove that every minimizer of $C_{min}^k$ is also a minimizer of $C_{min}$. Let thus $U$ be a minimizer for $C_{min}^k$, i.e. $\mc{T}(U)=C_{min}^k$ and $\mc{V}_k(U)=1$. 
Since $\mc{T}(U)=C_{min}$, it suffices to prove that $\mc{V}(U)=1$. By definition, we have 
$$\mc{V}(U)\geq \mc{V}_k(U)=1.$$
On the other hand, set $\theta:= (\mc{V}(U))^{\frac{1}{N}}$ we obtain, by scaling, that $\mc{V}(U(\theta \cdot))=1$ and thus
$$\mc{T}(U)=C_{min} \leq \mc{T}(U(\theta \cdot))= \theta^{-\frac{N+2s}{N}} \mc{T}(U)$$
from which we achieve
$$\mc{V}(U)\leq 1.$$
This concludes the proof.
\QED

\bigskip

%%%%%%%%%%%%%%%%%%%%%%%%%%%%%%%%%%%%%%%%%%%%%%%%%
%
%
%\subsection{Proof of Theorem \ref{teo_main}}
%\label{sec_proof}

We are now ready to prove Theorem \ref{teo_main}.

\medskip

\claim Proof of Theorem \ref{teo_main}.

\noindent
\textbf{Step 1.} We first look at the truncated problem \eqref{eq_principale_k}. 
Indeed, by Theorems \ref{teo_concen_esist} and \ref{teo_concen_conc} we obtain the existence of $\cupl(K)+1$ sequences of solutions of \eqref{eq_principale_k} satisfying the properties of Theorem \ref{teo_main} for $\eps>0$ small. 
For each of these sequences $v_{\eps}$ of solutions of \eqref{eq_principale_k}, called $x_{\eps}\in \R^N$ a global maximum point of $v_{\eps}$, 
we obtain
$$\lim_{\eps\to 0}d(x_{\eps}, K) =0$$
and
\begin{equation*}\label{eq_stima_polinom}
\frac{C'}{1+|\frac{x-x_{\eps}}{\eps}|^{N+2s}}\leq v_{\eps}(x) \leq \frac{C''}{1+|\frac{x-x_{\eps}}{\eps}|^{N+2s}} \quad \hbox{ for $x \in \R^N$}
\end{equation*}
where $C', C''>0$ are uniform in $\eps>0$. 

Moreover, for every sequence $\eps_n \to 0^+$ there exist $U \in \widehat{S}_k$ and an $x_0 \in \R^N$ such that, up to subsequences,
\begin{equation}\label{eq_conv_U}
v_{\eps_n}(\eps_n \cdot + x_{\eps_n}) \to U(\cdot + x_0), \quad \hbox{as $n\to +\infty$} 
\end{equation}
in $H^s(\R^N)$ and locally on compact sets.

\smallskip

\noindent
\textbf{Step 2.} Notice that by Proposition \ref{prop_uguagl_insiem} we have $U \in \widehat{S}$, thus $U(\cdot + x_0)$ is a ground state of \eqref{eq_limite}. 
We prove now that $v_{\eps}$ are solutions of the original equation, which is given by
\begin{equation}\label{eq_dim_stimaM}
\norm{v_{\eps}}_{\infty} < M \quad \hbox{ definitely for $\eps$ small}.
\end{equation}
Assume by contradiction that \eqref{eq_dim_stimaM} does not hold: thus there exists a sequence $\eps_n \to 0$ such that
$$\norm{v_{\eps_n}}_{\infty}\geq M \quad \hbox{for each $n \in \N$}.$$
By the previous Step, there exists an $U\in \widehat{S}_k$ and an $x_0 \in \R^N$ such that, up to subsequence, \eqref{eq_conv_U} holds. 
In particular, by the pointwise convergence we obtain
$$\norm{v_{\eps_n}}_{\infty}= v(x_{\eps_n})\to U(x_0) \leq \norm{U}_{\infty} < M$$
which implies
$$\norm{v_{\eps_n}}_{\infty} < M$$
definitely for $n\gg 0$, which is an absurd. Thus \eqref{eq_dim_stimaM} holds. As a consequence
$$f_k(v_{\eps})=f(v_{\eps})$$
and hence $v_{\eps}$ are solutions of the original problem \eqref{eq_concent_genericaf}, satisfying the desired properties.
\QED

\begin{Remark}
We point out that the found solutions are perturbations of ground states of the truncated limiting problem \eqref{eq_limite_k} which are, on the other hand, coinciding with the ground states of the critical limiting problem \eqref{eq_limite} thanks to Proposition \ref{prop_uguagl_insiem}. 
One may think to search directly the solutions as perturbation of functions in the compact set $\widehat{S}$, but actually the direct approach in a critical setting reveals several problems, such as the convergence of the Palais-Smale sequences. 
A different and direct approach is given in \cite{Amb4} by means of Concentration-Compactness techniques, but in the assumptions that $f$ satisfies monotonicity and Ambrosetti-Rabinowitz conditions.
\end{Remark}

\begin{Remark}\label{rem_caso_subcr_nonstr}
%\tr{DA CONTROLLARE} \\ %CCCOMMENT NOW
We see that actually the ideas of this Section adapts to study the case of $f$ negatively critical $a<0$, or subcritical $a=0$ (but not in the strict sense of \hyperref[(f1)]{\textnormal{(f1.3)}} treated in the previous Sections), % \ref{...}, 
that is
$$\lim_{t \to +\infty} \frac{f(t)}{|t|^{2^*_s-1}}=a \in (-\infty,0], %=0, %In realtà \limsup
$$
instead of \hyperref[(f1')]{\textnormal{(f1.3')}}, %(f4)}-\textnormal{(f5)}, 
filling the gap between the papers \cite{CG0} and \cite{Gal0}. This case covers functions of the type $f(t)= |t|^{p-2}t - |t|^{2^*_s-2}t$ and $f(t)= \frac{|t|^{2^*_s-2}t}{\log(t^2+2)}$.

In order to achieve this result, we sketch the steps:
\begin{itemize}
\item We substitute the existence result Proposition \ref{prop_esist_ground} with the one by \cite{ChWa0}, observing that $f \in C^1$ is needed only to get the Pohozaev identity, thus our assumptions \hyperref[(f3)]{\textnormal{(f3)}} %(f1)} 
is enough (see also Remark \ref{rem_esist_nonstrett_unc}). % and \cite{GalSch}).
\item The uniform $L^{\infty}$-bound of Proposition \ref{prop_upperbound} can be easily adapted.
\item The uniform $C^{0,\sigma}$-bound can be obtained as in the Step 3 of the proof of Theorem \ref{teo_concen_conc}.
\item The compactness result Proposition \ref{prop_compat_S} can be obtained as in the strict-subcritical case Lemma \ref{lemma_stima_unif}. When $a=0$, the proof follows verbatim, otherwise we adapt Step 4 in the following way.
\\Let $f=f^+-f^-$. If $a<0$ it means that $f^+(t)=0$ for $t\gg0$; in particular $f^+$ is subcritical. By knowing $U_k \wto U$ in $H^s(\R^N)$, $a_k \to a$ and $L'_{a_k}(U_k) U_k=0=L'_a(U)U$ we want to show that $\norm{(-\Delta)^{s/2}U_k}_2^2 + a_k \norm{U_k}_2^2 \to \norm{(-\Delta)^{s/2}U}_2^2 + a \norm{U}_2^2$. Observe that, by \eqref{eq_stima_uniformeU}, Proposition \ref{prop_converg_generiche_loc} and Fatou's Lemma, for any $\eta>0$ there exists $R\gg 0$ such that
$$ \pabs{ \int_{B_R^c} f(U_k)U_k }, \pabs{ \int_{B_R^c} f(U) U }<\eta \quad \hbox{for each $k \in \N$},$$
$$ \int_{B_R} f^+(U_k)U_k \to \int_{B_R} f^+(U) U,$$
$$ \liminf_{k \to +\infty} \int_{B_R} f^-(U_k)U_k \geq \int_{B_R} f^-(U) U.$$
Thus
%\begin{align*}
%\limsup_{k} \left(\norm{(-\Delta)^{s/2}U_k}_2^2 + a_k \norm{U_k}_2^2 \right) &\leq \limsup_{k} \int_{\R^N} f(U_k) U_k \\
%& \leq \limsup_k \int_{B_R^c} f(U_k)U_k dx + \limsup_k \int_{B_R} f^+ (U_k)U_k dx - \liminf_k \int_{B_R} f^- (U_k)U_k dx
%\end{align*}
\begin{eqnarray*}
\lefteqn{\limsup_{k} \left(\norm{(-\Delta)^{s/2}U_k}_2^2 + a_k \norm{U_k}_2^2 \right) \leq \limsup_{k} \int_{\R^N} f(U_k) U_k} \\
& \leq & \limsup_k \int_{B_R^c} f(U_k)U_k + \limsup_k \int_{B_R} f^+ (U_k)U_k - \liminf_k \int_{B_R} f^- (U_k)U_k \\
& \leq & \eta + \int_{B_R} f^+ (U)U - \int_{B_R} f^- (U)U %= \eta + \int_{B_R} f (U)U \\
 = \eta + \int_{\R^N} f (U)U - \int_{B_R^c} f(U) U \\
&\leq& 2\eta + \int_{\R^N} f (U)U = 2 \eta + \norm{(-\Delta)^{s/2}U}_2^2 + a \norm{U}_2^2.
\end{eqnarray*}
Letting $\eta \to 0$ we obtain
$$\limsup_{k} \left(\norm{(-\Delta)^{s/2}U_k}_2^2 + a_k \norm{U_k}_2^2 \right) \leq \norm{(-\Delta)^{s/2}U}_2^2 + a \norm{U}_2^2 $$
which, together with the semicontinuity of the norm, gives the claim.
\QED
\end{itemize}
The remaining part of the proof follows the lines of the critical case treated in this Section.
\end{Remark}

%%%%%%%%%%%%%%%%%%%%%%%%%%%%%%%%%%%%%%%%%%%%%%%%%

\subsection{The local case}
\label{sec_locale}

The arguments presented in Theorem \ref{teo_main} apply, with suitable modifications, also to local nonlinear Schr\"odinger equations. We give here some details. 
Condition \hyperref[(f1')]{\textnormal{(f1')}} %--\hyperref[(f3)]{\textnormal{(f3)}} %(f5) 
rewrites in the local case $s=1$ as
%\begin{itemize}
%\item[(f1')] $f\in C(\R, \R)$,
%\item[(f2')] $f(t)\equiv 0$ for $t\leq 0$,
%\item[(f3')] $\, \lim_{t \to 0} \frac{f(t)}{t}=0$,
%\item[(f4')] $\lim_{t \to +\infty} \frac{f(t)}{t^{2^*-1}}=a>0$, where $2^*:=\frac{2N}{N-2}$,
%\item[(f5')] for some $C>0$ and $\max\{2^*-2, 2\} <p< 2^*$, i.e. satisfying
%$$p \in \parag{&(4, 6)& \quad N =3, \\ &\Big(2, \frac{2N}{N-2}\Big)& \quad N\geq 4, }$$
%it results that
%$$f(t) \geq a t^{2^*-1} + C t^{p-1} \quad \hbox{for $t \geq 0$}.$$
%\end{itemize}
%
\begin{itemize}
	\item[(f1')] \label{(f1'l)}
Berestycki-Lions type assumptions with respect to $m_0>0$, that is
	\begin{itemize}
		\item[(f1.1)] $\, f\in C(\R, \R)$;
		\item[(f1.2)] $\, \lim_{t \to 0} \frac{f(t)}{t}=0$;
		\item[(f1.3')] $\, \lim_{t \to +\infty} \frac{f(t)}{t^{2^*-1}}=a>0$, where $2^*=\frac{2N}{N-2}$, and moreover
 for some $C>0$ and $\max\{2^*-2, 2\} <p< 2^*$, i.e. satisfying
$$p \in \parag{&(4, 6)& \quad N =3, \\ &\Big(2, \frac{2N}{N-2}\Big)& \quad N\geq 4, }$$
it results that
$$f(t) \geq a t^{2^*-1} + C t^{p-1} \quad \hbox{for $t \geq 0$};$$
		\item[(f1.4)] $\, F(t_0)> \frac{1}{2} m_0 t_0^2$ for some $t_0>0$. %, where $F(t):= \int_0^t f(s) ds$;
	\end{itemize}
%	\item[(f2)] $f(t)\equiv 0$ for $t\leq 0$.
\end{itemize}

See also Remark \ref{rem_ipotesi_reg} for some weakening and comments on the assumption \hyperref[(f1'l)]{\textnormal{(f1.3')}}. %(f5'). 

\begin{Theorem}\label{teo_main_locale}
Suppose $s=1$, $N\geq 3$ and that \hyperref[(V1)]{\textnormal{(V1)}}-\hyperref[(V2)]{\textnormal{(V2)}}, \hyperref[(f1'l)]{\textnormal{(f1')}}-\hyperref[(f2)]{\textnormal{(f2)}} %(f5')} 
hold. 		
Let $K$ be defined by \eqref{eq_def_K_2}.
Then, for small $\eps>0$ the equation 
$$-\eps^{2} \Delta v + V(x) v = f(v), \quad x \in \R^N$$
 has at least $\cupl(K)+1$ positive solutions, which belong to $C^{1, \sigma}(\R^N) \cap L^{\infty}(\R^N)$ for some $\sigma \in (0,1)$. Moreover, each of these sequences $v_{\eps}$ concentrates in $K$ as $\eps \to 0$: namely, 
there exist $x_{\eps} \in \R^N$ global maximum points of $v_{\eps}$, such that
$$\lim_{\eps\to 0}d(x_{\eps}, K) =0$$
and
\begin{equation}\label{eq_teo_stima_esp}
 v_{\eps}(x) \leq C' \textnormal{exp}\Big(-C''\Big|\frac{x-x_{\eps}}{\eps}\Big|\Big) \quad\hbox{ for $x \in \R^N$}
\end{equation}
where $C', C''>0$ are uniform in $\eps>0$.
Finally, for every sequence $\eps_n \to 0^+$ there exists a ground state solution $U$ of
$$-\Delta U + m_0 U = f(U), \quad x \in \R^N$$
and a point $x_0 \in K$
such that, up to a subsequence, $x_{\eps_n} \to x_0$ and
$$v_{\eps_n}(\eps_n \cdot + x_{\eps_n}) \to U \quad \hbox{as $n\to +\infty$}$$
in $H^s(\R^N)$ and locally on compact sets.
 \end{Theorem}

\claim Proof.
The arguments of the previous Sections apply mutatis mutandis. 
Indeed, we define in the same way the set of ground states $\widehat{S}$, which turns to be nonempty \cite{ZhZo1} and compact. Moreover to get the uniform $L^{\infty}(\R^N)$ bound, one can easily adapt the proof of Proposition \ref{prop_upperbound} after observing that by the chain rule it holds
$$|\nabla h(U)|^2 = \nabla U \cdot \nabla \tilde{h}(U), \quad U \in H^1(\R^N),$$
where we recall that $\tilde{h}' = (h')^2$. 
Then the truncation machinery can be set in motion, and one can prove $\widehat{S}_k = \widehat{S}$.
Existence, multiplicity and decay of solutions of the truncated problem are given by \cite[Theorem 1.1 and Remark 1.3]{CJT}; the regularity is instead a consequence of standard elliptic estimates \cite[Appendix B]{Stru0}.
\QED

\appendix

\pagestyle{fancy} %per iniziare a modificare
\fancyhf{} %azzera ciò che è predefinito delle testatine
%\nouppercase{} %fa sì che scriva minuscolo e non maiuscolo
% L=left (nella stessa pagina), R=right (nella stessa pagina), C=center (nella stessa pagina)
% E=even (page), O=odd (page)
% \leftmark, di default, contiene ciò che solitamente viene messo sulle pagine a sinistra, cioè il capitolo. Si può ridefinire.
% \rightmark, di default, contiene ciò che solitamente viene messo sulle pagine a destra, cioè la sezione. Si può ridefinire.
% \thepage inserisce il numero della pagina

%Per ridefinire:
% #1 è il titolo del capitolo (o sezione) corrente.
%\chaptername ti fa uscire "Chapter".
% \thechapter è il numero del capitolo corrente, stessa cosa \thesection
\renewcommand\sectionmark[1]{%
 \markboth{\MakeUppercase{Appendix – #1}}{}} %Fa comparire nella testatina Appendix - titolo della sezione

\fancyhead[LE, LO]{\nouppercase{\leftmark} }
\fancyhead[RE, RO]{\thepage}

%\phantomsection%

%\addchap{Appendix}

\chapter[Appendix]{}
\label{chap_app_alg_top}

\vspace{-2em}

\renewcommand{\thesection}{A}

%\automark{section}

%VERSIONE STAMPATA
%\fancyhead[LE, RO]{\thepage}
%\fancyhead[RE]{\nouppercase{\leftmark}}
%\fancyhead[LO]{\nouppercase{\leftmark}}

%
\fancyhead[L]{ \ifthenelse{\isodd{\value{page}}}{\nouppercase{\leftmark}}{\thepage} }
\fancyhead[R]{ \ifthenelse{\isodd{\value{page}}}{\thepage}{\nouppercase{\leftmark}} }

%%%%%%%%%%%%%%%%%%%%%%%%%%%%%%%%%%%%%%%%%%%%%%%%%%%%%%%%%%%%%%%%%%%%%
%%%%%%%%%%%%%%%%%%%%%%%%%%%%%%%%%%%%%%%%%%%%%%%%%%%%%%%%%%%%%%%%%%%%%%

\section{Some algebraic topology: the relative cup-length}

%\tor{Quasi tutta questa sezione andrebbe in arancione}

In order to estimate the number of critical points of certain functionals not bounded from below and above, it is useful to implement the algebraic-topological tool of the \emph{relative cup-length}, together with the more used \emph{relative category}. In this Appendix we briefly recall the basic notion of algebraic topology needed to define this object; afterwards we will highlight how it relates to the category and how they are exploited in order to gain multiple solutions of PDEs. Finally we will briefly recall also the definition of the genus.

\subsection{The singular cohomology}

%We briefly recall in this Section some notions from algebraic geometry which lead to the definition of cup-length. 

We start by defining the \emph{singular cohomology}. Here we essentially follow the self-contained description due to \cite{Cha0}, 
%based on the \emph{singular cohomology}, 
but we refer also to \cite{Hat0, Vic0, Mass2, Dol0, GH0,GH0, Spa0, Mass1, EiSt0, Bre0}. 

Fix $X$ a topological space (in our case it will be a subset of some Hilbert space, such as $\R^N$ or $H^s(\R^N)$, see Section \ref{sec_cup-length}), and fix an abelian group $G$: actually the choice of $G$ does not heavily influence the main properties of cohomology, and %to make things easy 
usually $G$ is chosen as a generic field $\mathbb{F}$ \cite{CJT}, or some specific ones like the real field $\R$ \cite{FoWi2,FLRW} or the $\Z_2$ field \cite{Szu0}.

Let $q\in \N$, and let $\Delta_q$ be the \emph{$q$-simplex} defined by
$$\Delta_q:=\left \{ \sum_{j=0}^q \lambda_j e_j \mid \lambda_j \geq 0, \; \sum_{j=0}^q \lambda_j=1 \right \} \equiv \left\{ (\lambda_0,\lambda_1, \dots, \lambda_q, 0, \dots) \mid \lambda_j \geq 0, \; \sum_{j=0}^q \lambda_j=1\right \} $$
where $e_0:=(0, 0, \dots)$, $e_1:=(0,1, \dots)$ and so on, are vectors in $\R^{\infty}$.
%, $e_j=(0, 0, \dots 1, \dots )$, the number $1$ in the $j$-th position.
We define the set of \emph{singular $q$-simplexes} by
$$\Sigma_q(X):=\{ \sigma: \Delta_q \to X \; \hbox{ continuous}\}.$$
Starting from $\Sigma_q(X)$ and $G$ we can build the free abelian group $C_q(X, G)$ with bases $\Sigma_q(X)$, that is
$$C_q(X, G):=\left\{ \sum_{i, \textnormal{ finite}} g_i \sigma_{i} \mid g_i\in G, \; \sigma_i \in \Sigma_q(X)\right\}$$
where the linear combination has to be intended in the formal sense\footnote{For example, if $G=R$ is a ring with unit $1_R$, we define the free abelian group in this way \cite[page 4]{Vic0}: 
start by identifying the elements $\sigma \in \sigma_q$ with the functions $f_{\sigma}:\Sigma_q(X)\to R$, $f_{\sigma}(\tau):= \parag{1_R \hbox{ if $\tau= \sigma$} \\ 0_R \hbox{ if $\tau\neq \sigma$}}$. Then set
$$C_q(X, G):=\{ f: \Sigma_q(X) \to \Z \mid f(\sigma)\neq 0 \hbox{ for a finite number of $\sigma \in \Sigma_q(X)$}\}$$
and observe that $\Sigma_q(X) \equiv \{ f_{\sigma}\}_{\sigma \in \Sigma_q(X)}$ is a basis for $C_q(X, G)$, that is, elements of $C_q(X,G)$ are of the form
$$f = \sum_{i, \textnormal{ finite}} g_i f_{\sigma_i}.$$
}.
We call $C_q(X, G)$ the set of \emph{singular $q$-chains}; here an inner summation and an external product (through $G$) can be easily defined. %\tr{It is a vectorial space with an inner sum and an external product through $G$. (?)}

We define now a boundary operator on $C_q(X,G)$, introducing it first on $\Sigma_q(X)$ and then extending it by linearity. Indeed, for any $q \geq 1$ and $\sigma \in \Sigma_q(X)$ and $j=0 \dots q$ we define $\sigma^{(j)} \in \Sigma_{q-1}$ by
%$$\sigma^{(j)}(x_0, x_1, \dots, x_{q-1}):= \sigma(x_0, x_1, \dots 0, \dots x_{q-1})$$
$$\sigma^{(j)}(x_0, x_1, \dots, x_{q-1}):= \sigma(x_0, x_1, \dots x_{j-1}, 0, x_j, \dots x_{q-1})$$
%\tr{CHECK definizione} %CCCOMMENT NOW
where $0$ is in the $j$-th position. Thus the \emph{boundary operator} is defined as
$$\partial \sigma:= \sum_{j=0}^q (-1)^j \sigma^{(j)}$$
and hence
$$\partial: C_q(X,G) \to C_{q-1}(X, G).$$
If $q=0$, the boundary operator $\partial: C_0(X,G)\to G$ is defined as $\partial(\sum g_i \sigma_i):= \sum g_i$ (we are formally setting $C_{-1}(X,G):=G$). We have that $\partial$ is a homomorphism.
Set
$$C_*(X,G):= \bigoplus_{q\geq 0} C_q(X,G)$$
we have %can imagine $\partial$ here defined, that is 
$\partial: C_*(X,G)\to C_*(X,G)$. It is a straightforward computation showing that
$$\partial^2 = 0$$
which is of key importance in the theory of homologies and cohomologies. With these ingredients it is possible to define a homology $H_*(X,G)$; we are anyway % but since we are 
interested in \emph{cohomologies}, and thus we need first to pass on homomorphisms and dualities. Thus we define the set of \emph{singular $q$-cochains}
$$C^q(X,G):= Hom(C_q(X,G), G);$$
by using the bracket notation 
$$[\sigma, c]:=c(\sigma)$$
 for every $c \in C^q(X,G)$ and $\sigma \in C_q(X,G)$, the definition of $C^q(X,G)$ rewrites as
$$[\sigma_1+\sigma_2, c]= [\sigma_1, c] + [\sigma_2, c] \quad \hbox{and} \quad \quad [g \sigma, c]=g[\sigma, c]$$
for every $c\in C^q(X,G)$, $\sigma,\sigma_1, \sigma_2 \in C_q(X,G)$ and $g \in G$. We can hence define the dual operator of $\partial$, named the \emph{coboundary operator}, by
$$[\sigma, \delta c]:= [\partial \sigma, c]$$
for every $c \in C^{q-1}(X, G)$, $\sigma \in C_q(X,G)$; thus
$$\delta: C^{q-1}(X,G)\to C^q(X,G),$$
which is a homomorphism. 
%Set
%$$C^*(X,G) := \bigoplus_{q \geq 0} C^q(X,G)$$
%we have %can think 
%$\delta: C^*(X,G) \to C^*(X,G)$, and we obtain
%$$\delta^2=0.$$
%In particular the last property easily implies that $\Ima(\delta) \triangleleft \Ker(\delta)$,
%% $\Ima(\delta) \subset \Ker(\delta)$ (actually $\Ima(\delta) \triangleleft \Ker(\delta)$); 
%thus we are allowed to define the \emph{singular $q$-cohomology group}
%$$H^q_{\Delta}(X,G):= \Ker(\delta) / \Ima(\delta)$$
%and hence the \emph{singular cohomology group} on $X$ with coefficients in $G$
%$$H^*_{\Delta}(X,G):= \bigoplus_{q\geq 0} H^q_{\Delta}(X,G).$$
Set
$$C^*(X,G) := \bigoplus_{q \geq 0} C^q(X,G)$$
we have %can think 
$\delta: C^*(X,G) \to C^*(X,G)$, and we obtain
$$\delta^2=0.$$
In particular the last property easily implies that $\Ima(\delta) \triangleleft \Ker(\delta)$,
% $\Ima(\delta) \subset \Ker(\delta)$ (actually $\Ima(\delta) \triangleleft \Ker(\delta)$); 
thus we are allowed to define the \emph{singular $q$-cohomology group}
$$H^q_{\Delta}(X,G):= \Ker(\delta_{|C^q(X,G)}) / \Ima(\delta_{|C^{q-1}(X,G)})$$
the sets $%Z:=
\Ker(\delta)$ and $%B:=
\Ima(\delta)$ are said, respectively, the sets of the \emph{cocycles} and of the \emph{coboundaries}. We highlight that $H^0(X,G)$ may be interpreted as the set of functions $X\to G$ constant on path-components of $X$ \cite[pages 198-199]{Hat0} (see also \cite[Proposition 3.11]{Vic0}, \cite[page 183]{GH0}, \cite[page 244]{Spa0}, \cite[Lemma 1.2]{Mass1}), %; \tr{thus it is trivial if and only if $G$ is trivial (?).}
while $H^q(\emptyset, G)$ is the trivial cohomology \cite[page 192]{Mass1}.

Moreover we define the \emph{singular cohomology group} on $X$ with coefficients in $G$
$$H^*_{\Delta}(X,G):= \bigoplus_{q\geq 0} H^q_{\Delta}(X,G).$$

Assume from now on $G=R$ to be a commutative ring with unit. On the cohomology $H^*_{\Delta}(X,R)$ (also called \emph{cohomology ring} of $X$ \cite[Remark 8.17]{Dol0}) we can define a \emph{cup product}: instead of introducing it in terms of \emph{cross product}, we give here a direct construction. % and a direct way to define it is the following. 
We start by defining it on $C^*(X,R)$, then by quotient we obtain it also on $H^*_{\Delta}(X,R)$. Indeed, fixed $p, q \geq 0$, we define
$$\phi_p: \Delta_{p}\to \Delta_{p+q}, \quad \beta_q: \Delta_{q} \to \Delta_{p+q}$$
the immersions in the first $p$ components and in the last $q$ components respectively, i.e.
$$\phi_p(\lambda_1, \dots, \lambda_p, 0, \dots) := (\lambda_1, \dots, \lambda_p, 0, \dots, 0, 0, \dots),$$
$$\beta_q(\lambda_1, \dots, \lambda_q, 0, \dots) := (0, \dots, 0, \lambda_1, \dots, \lambda_q, 0, \dots),$$
 so that, if $\sigma \in C_{p+q}(X,R)$, then $\sigma \phi_p \in C_p(X,R)$ and $\sigma \beta_q \in C_q(X,R)$. Thus we define, through the product in $R$, the cup product
$$[\sigma, c \smile d]:= [\sigma \phi_p, c] [\sigma \beta_q, d]$$
%\tr{per Hatcher dobbiamo essere in un anello e quello è il suo prodotto}
for any $c \in C^p(X,R)$, $d \in C^q(X,R)$ and $\sigma \in C_{p+q}(X,R)$, which implies
$$\smile: C^p(X,R) \times C^q(X,R) \to C^{p+q}(X,R)$$
and more generally, $\smile: C^*(X,R) \times C^*(X,R) \to C^*(X,R)$. 
Notice that multiplying $c \in C^p(X,R)$ with $d \in C^0(X,R)$ means multiplying by constant elements of the form $\sum_{i, \, \textnormal{finite}} g_i \sigma_i(e_0)$, with $g_i \in G$ and $\sigma_i \in \Sigma_{p}(X)$.
The cup product results bilinear, associative and unitary. % with unit. 
Moreover, it satisfies $c \smile d= (-1)^{pq} d \smile c$ (since $R$ is commutative %when $G=R$ is a commutative ring, 
 \cite[Theorem 3.11]{Hat0}), which implies that it is \emph{skew}-commutative: even if not properly commutative, it %(and not properly commutative), but 
nevertheless satisfies 
\begin{equation}\label{eq_zero_sse_zero}
c\smile d = 0 \iff d\smile c=0.
\end{equation}
Moreover, it holds 
$$\delta(c \smile d) = \delta c \smile d + (-1)^p \smile \delta d,$$
for $c \in C^p(X,R)$ and $d \in C^q(X,R)$; in particular, this easily implies that $\Ker(\delta)$ is a subalgebra of $C^*(X,R)$ and $\Ima(\delta)$ is an ideal of $\Ker(\delta)$. Thus, $\smile$ can be passed to the quotient and hence defined on 
$$\smile: H^*_{\Delta}(X,R) \times H^*_{\Delta}(X,R) \to H^*_{\Delta}(X,R).$$

\bigskip

Starting from a whatever cohomology $H^*(X,R)=H^*_{\Delta}(X,R)$ %(which we highlight \emph{need not} to be the singular cohomology, 
(see Section \ref{subsec_other_coho}), we can define the \emph{cup-length} as the length of the longest nontrivial chain of cup products in $H^*(X,R)$ (up to the constants in $H^0(X,R)$).
\begin{Definition}
Let $X$ be topological space and $R$ be a commutative ring with unit. We define the \emph{cup-length of $X$} as
%$$\cupl(X,R):= \max \{ l \in \N^* \mid \exists \alpha_i \in H^{q_i}(X,R), q_i \geq 1 \hbox{ for $i=1 \dots l$, s.t. } \alpha_1 \smile \dots \smile \alpha_l \neq 0\};$$
\begin{align*}
\cupl(X,R):= \max \{ l \in \N^* \mid \exists & \alpha_i \in H^{q_i}(X,R), q_i \geq 1 \hbox{ for $i=1 \dots l$,} \\
& \hbox{ s.t. } \alpha_1 \smile \dots \smile \alpha_l \neq 0\};
\end{align*}
%$$\cupl(X,R):= \max \{ l \in \N^* \mid \exists (\alpha_i)_{i=1\dots l} \subset H^*(X,R)\setminus H^0(X,R) \; s.t. \; \alpha_1 \smile \dots \smile \alpha_l \neq 0\};$$
if such $l$ does not exists, but $H^0(X,R)$ is not trivial, we set $\cupl(X,R):=0$, otherwise (if $H^0(X,R)=\{0\}$) we set $\cupl(X,R):=-1$.
%\tr{come fa a essere $l=0$ se parte da $1$? CAPISCI}
%while $\cupl(X,R):=-1$ if such $l\in \N$ does not exists. 
\end{Definition}
We notice that by \eqref{eq_zero_sse_zero}, the order in the cup product is of no importance. In the case $X$ is not connected, a slightly different definition (which makes the cup-length additive) can be found in \cite{Bar0}. 

For explicit computations of the cup-length we refer to \cite[Example 3.4 and page 19]{FLRW} 
 and to \cite{CLOT}: for instance if $B\subset \R^N$ is the closed unit ball, then $\cupl(\partial B)=1$ for $N \geq 2$; if $T^N$ is the $N$-dimensional torus, then $\cupl(T^N)=N$.

\subsubsection{Singular relative cohomology and cup-length}

We define now the cohomology and the cup-length relative to a \emph{topological pair} $(X,Y)$, that is $Y\subset X$ topological spaces. 
%All the objects defined in the previous Section are still valid. 
%In addition, we consider
Observed that
$$C_q(Y, G) \triangleleft C_q(X, G)$$
and %observe 
that $\partial$ conserves $C_q(Y, G)$, we can define the \emph{singular $q$-relative chain module}
$$C_q(X, Y, G):= C_q(X, G) / C_q(Y, G).$$
Notice that $C_q(X,\emptyset,G)\equiv C_q(X,G)$. 
%(?) %CCCOMMENT NOW
Considered the canonical projection $\pi_q: C_q(X, G) \to C_q(X, Y, G)$, we introduce
$$\tilde{\partial}: C_q(X, Y, G) \to C_{q-1}(X, Y, G)$$
the well defined function such that the canonical diagram commutes
$$\tilde{\partial} \circ \pi_q = \pi_{q-1} \circ \partial.$$
%(namely, some canonical diagram commutes). 
%From now on, t
The other definitions follows in the same way as before:
$$C^q(X, Y, G):= Hom(C_q(X, Y, G), G),$$
$$\tilde{\delta}: C^{q-1}(X, Y, G) \to C^q(X, Y, G),$$
$$H^q_{\Delta}(X, Y, G):= \Ker(\tilde{\delta})/\Ima(\tilde{\delta}),$$
and also $C_*(X, Y, G)$, $C^*(X, Y, G)$, $H^*_{\Delta}(X, Y, G)$ and $\smile$ (see also \cite[page 209]{Hat0}). 
Notice that, if $X$ is path-connected and $Y\neq \emptyset$, then $H^0(X,Y,R)$ is trivial \cite[page 183]{GH0}.

\bigskip

When $G=R$, we can define the \emph{relative cup-length} as the length of the longest chain of cup products in $H^*(X,R)$ multiplied with an element of $H^*(X,Y,R)$; see also \cite{Szu0, FoWi2, FLRW}.
\begin{Definition}
Let $(X,Y)$ be a topological pair and $R$ be a commutative ring with unit. We define the \emph{cup-length of $X$, relative to $Y$} as
%\begin{eqnarray*}
%\cupl(X, Y,R):= \max \{ l \in \N^* \mid &\exists \alpha_i \in H^{q_i}(X,R), q_i \geq 1 \hbox{ for $i=1 \dots l$,} \; \exists \alpha_0 \in H^*(X, Y,R)&\\
%& \hbox{s.t. } \alpha_0 \smile \alpha_1 \smile \dots \smile \alpha_l \neq 0\}.&
%\end{eqnarray*}
\begin{eqnarray*}
\lefteqn{\cupl(X, Y,R):=} \\
& \max \{ l \in \N^* \mid &\exists \alpha_i \in H^{q_i}(X,R), q_i \geq 1 \hbox{ for $i=1 \dots l$,} \; \exists \alpha_0 \in H^*(X, Y,R)\\
&& \hbox{s.t. } \alpha_0 \smile \alpha_1 \smile \dots \smile \alpha_l \neq 0\}.
\end{eqnarray*}
%\begin{eqnarray*}
%\cupl(X, Y,R):= \max \{ l \in \N^* \mid &\exists (\alpha_i)_{i=1\dots l} \subset H^*(X,R)\setminus H^0(X,R), \; \exists \alpha_0 \in H^*(X, Y,R)&\\
%& s.t. \; \alpha_0 \smile \alpha_1 \smile \dots \smile \alpha_l \neq 0\};&
%\end{eqnarray*}
if such $l\in \N$ does not exists, but $H^*(X,Y,R)$ is not trivial, we set $\cupl(X,Y,R):=0$, otherwise (if $H^*(X,Y,R) = \{0\}$) we set $\cupl(X,Y,R):=-1$. 
%(?) %CCCOMMENT NOW
 %è l'unica possibilità per essere $-1$, perché $id^*=id = 0$ quando c'è solo l'elemento nullo!
%\tr{come fa a essere $l=0$ se parte da $1$? CAPISCI}
%\tr{e $l=-1$?}
\end{Definition}
Notice that 
$$\cupl(X,R) = \cupl(X, \emptyset,R);$$
this is the same as taking $\alpha_0 \in H^0(X,R)$, since $H^*(X,\emptyset,R)$ is essentially $H^*(X,G)$ (see also \cite[page 256]{Mass1} and \cite[Proposition 12.3]{Bre0}).
Again, for explicit examples we refer to \cite[Example 3.4]{FLRW}: 
for instance, if $B\subset \R^N$ is the closed unit ball, then $\cupl(B, \partial B)=0$. % for $N\geq 1$.

%$H^*(X, Y,R)$ è un sottoinsieme di $H^*(X)$, così da poterci fare il cup product?

\subsubsection{Cup-length relative to a function}

Let us consider two topological pairs $(X, Y)$ and $(X', Y')$ and a continuous map $f:(X, Y) \to (X', Y')$, that is $f: X \to X'$ with $f(Y) \subseteq Y'$. 
It is possible to prove\footnote{
We show here the standard construction in the non-relative case \cite[Section VII.3]{Mass2}. 
Consider the induced function
$$f^{\#}:C^q(X') \to C^q(X)$$
such that
$$(f^{\#}(c))(\sigma) = c(f \circ \sigma) \quad \hbox{(i.e. $[\sigma, f^{\#}(c)]=[f \circ \sigma, c]$)}$$
for every $c \in C^q(X')$ and $\sigma \in C_q(X)$. %(or using bracket notation, $[\sigma, f^{\#}(c)]=[f \circ \sigma, c]$).
A straightforward computation shows that $ f^{\#} \circ \delta'=\delta \circ f^{\#}$, where $\delta': C^q(X') \to C^{q+1}(X')$; this easily implies that $f^{\#}(\Ker(\delta'))\subset \Ker(\delta)$ and $f^{\#}(\Ima(\delta'))\subset \Ima(\delta)$. This allows to pass to the quotient and define
$$f^*: \Ker(\delta')/\Ima(\delta') \to \Ker(\delta)/\Ima(\delta).$$
}
that $f$ induces and homomorphism of groups
$$f^*: H^*(X', Y', G) \to H^*(X, Y, G)$$
which is \emph{suitable functorial}, namely
$$(id)^*=id, \quad (gf)^*=f^* g^*, \quad f^*=g^* \hbox{ whenever $f, g$ homotopic}.$$
Moreover one can show that
$$\tilde{\partial} \circ f^* = f^* \circ \tilde{\partial}, \quad f^*(\alpha_1 \smile \alpha_2)= f^*(\alpha) \smile f^*(\alpha_2);$$
the second is said the \emph{naturality} of the cup product (see \cite[Section 7.8.6]{Dol0}). 
With this tool, when $G=R$, we can define the \emph{cup-length relative to $f$}, as the length of the longest chain of cup products in $f^*(H^*(X',R)) \subset H^*(X,R)$ multiplied with an element of $f^*(H^*(X',Y',R))\subset H^*(X,Y,R)$; see also \cite{BaWe0}. %, mediated through $f^*$.
\begin{Definition}
Let $(X, Y)$, $(X', Y')$ be two topological pairs, $R$ be a commutative ring with unit and $f:(X, Y) \to (X', Y')$ be continuous, with $f^*$ the induced homomorphism on the relative cohomolgies. We define the \emph{cup-length relative to $f$} as
%\begin{eqnarray*}
%\cupl(f,R):= \max \{ l \in \N^* \mid &\exists \alpha_i \in H^{q_i}(X',R), q_i \geq 1 \hbox{ for $i=1 \dots l$,} \; \exists \alpha_0 \in H^*(X', Y',R)&\\
%& \hbox{s.t. } f^*(\alpha_0) \smile f^*(\alpha_1) \smile \dots \smile f^*(\alpha_l) \neq 0\};&
%\end{eqnarray*}
\begin{eqnarray*}
\lefteqn{ \cupl(f,R):=} \\
& \max \{ l \in \N^* \mid &\exists \alpha_i \in H^{q_i}(X',R), q_i \geq 1 \hbox{ for $i=1 \dots l$,} \; \exists \alpha_0 \in H^*(X', Y',R)\\
&& \hbox{s.t. } f^*(\alpha_0) \smile f^*(\alpha_1) \smile \dots \smile f^*(\alpha_l) \neq 0\};
\end{eqnarray*}
%\begin{eqnarray*}
%\cupl(f,R):= \max \{ l \in \N^* \mid &\exists (\alpha_i)_{i=1\dots l} \subset H^*(X',R)\setminus H^0(X',R), \; \exists \alpha_0 \in H^*(X', Y',R)&\\
%& s.t. \; f^*(\alpha_0) \smile f^*(\alpha_1) \smile \dots \smile f^*(\alpha_l) \neq 0\};&
%\end{eqnarray*}
if such $l\in \N$ does not exist, but $f^*\nequiv 0$, it results $\cupl(f):=0$, otherwise (if $f^* \equiv 0$) we define $\cupl(f):=-1$.
\end{Definition}
Notice that 
$$\cupl(X, Y,R)=\cupl(id_{(X, Y)},R),$$
and in particular $\cupl(X,R)=\cupl(id_{(X, \emptyset)},R).$

\subsection{Other cohomologies}
\label{subsec_other_coho}

We highlight that other cohomologies could be used to define the cup-length: %other cohomologies could be used: 
for instance, the Alexander-Spanier cohomology \cite{CJT} and the \v{C}ech cohomology \cite{Szu0}. We sketch here how they are built, and then we point out how they are closely related to the singular cohomology.
\paragraph{Alexander-Spanier cohomology.}
We refer to \cite{Mass0, Mass1} (see also \cite{Spa0}). 
Let $X$ be a topological space and $G$ be a group. We define the \emph{abelian group of all the $q$-functions}
$$\Phi^q(X, G):=\{\varphi: X^{q+1} \to G\}$$
and its normal subgroup
$$\Phi_0^q(X, G):= \{ \varphi \in \Phi^q(X, G) \mid \varphi=0 \hbox{ in a neighborhood of the \emph{diagonal}} %\supp(\varphi)=\emptyset
 \}.$$
%where the support of a function is defined in the classical way for topological spaces. 
On $\Phi^q(X, G)$ we can define the \emph{coboundary operator} $\bar{\delta}: \Phi^q(X, G) \to \Phi^{q+1}(X,G)$ by
$$(\bar{\delta} \varphi)(x_0 \dots x_{q+1}):= \sum_{j=1}^{q+1} (-1)^j \varphi(x_0, \dots \hat{x}_j, \dots x_{q+1})$$
where $\hat{x}_j$ means that the variable is omitted; we have $\bar{\delta}^2=0$. Moreover, we can define the \emph{$q$-cochain}
$$\bar{C}^q(X, G):= \Phi^q(X,G)/\Phi^q_0(X, G)$$
and then the \emph{Alexander-Spanier cochain} $\bar{C}^*(X, G)$, on which we can %think to 
define $\bar{\delta}$ through quotients. Thanks to the property $\bar{\delta}^2=0$ we can pass to the quotient of $\Ker(\bar{\delta})$ over $\Ima(\bar{\delta})$ and obtain the \emph{%(\v{C}ech)-
Alexander-Spanier cohomology} $\bar{H}^*(X, G)$. Slightly different definitions, which focus on \emph{locally finitely valued} $q$-functions or which define $\Phi_0^q$ through \emph{supports}, %\emph{diagonal maps},
can be found in \cite{Mass0, Mass1}.
% could be also defined, see [Massey2].% and [Massey3] respectively.

\medskip

Once defined a \emph{relative Alexander-Spanier cohomology} $\bar{H}^*(X, Y, G)$, by exploiting \cite[Theorem 14.6.1 and Proposition 14.6.2]{Mass2}
% in [Massey1], 
one can show that actually, for $Y$ closed subset of $X$ paracompact Hausdorff space %\subset \R^N$ 
(for example a manifold, such as $\R^N$ or a more general Hilbert space, %which will be our case, 
see Section \ref{sec_cup-length}) and $G=R$ ring, it results that
\begin{equation}\label{eq_isom_Cech-AS_rel}
\bar{H}^*(X, Y, R) \approx H^*_{\Delta}(X, Y, R)
\end{equation}
that is, the Alexander-Spanier cohomology and the singular cohomology are isomorphic.

\paragraph{\v{C}ech cohomology.}

We refer to \cite{BT0, WP0} (see also \cite{Hat0, Dol0, Spa0, EiSt0}).
Let $X$ be a topological space and $G$ be a group (notice that we focus only on the case of a \emph{constant presheaf} with \emph{identical restrictions}). Let $\mathfrak{U}$ be a open covering of $X$ and define
$$\sigma:=(U_0, \dots U_q)$$
to be a \emph{$q$-simplex} if $U_i \in \mathfrak{U}$ and $|\sigma|:= \bigcap_{i=1}^q U_i \neq \emptyset $; $|\sigma|\subset X$ is called \emph{support} of $\sigma$. We thus define $\check{\Sigma}_q$ as the set of all the $q$-simplexes, and %. We thus define
$$\check{C}^q(\mathfrak{U}, G):= \{ \varphi : \check{\Sigma}_q \to G\}$$
the set of all the \emph{$q$-cochains}. On $\check{C}^q(\mathfrak{U}, G)$ we can define the \emph{coboundary operator} $\check{\delta}: \check{C}^q(\mathfrak{U}, G) \to \check{C}^{q+1}(\mathfrak{U}, G)$ as
$$(\check{\delta} \varphi)(U_0 \dots U_{q+1}):= \sum_{j=1}^{q+1} (-1)^j \varphi(U_0, \dots \hat{U}_j, \dots U_{q+1})$$
satisfying $\check{\delta}^2=0$. Thanks to this property we can define $\check{H}^*(\mathfrak{U}, G)$ by passing to the quotient the kernel and the image of $\check{\delta}$. Finally, considering the coverings of $X$ ordered by inclusion, we can define the \emph{\v{C}ech cohomology} as
$$\check{H}(X, G):= \lim_{\longrightarrow} \check{H}^*(\mathfrak{U}, G)$$
in the sense of the \emph{direct limits}. Notice that, if $X$ is an $n$-dimensional manifold and $\mathfrak{U}$ is a \emph{good cover}, i.e. every finite intersection of its elements is diffeomorphic to $\R^N$ (and there always exists such a good cover \cite[Theorem 5.1]{BT0}), then there is no need of passing to the direct limit, since it results that
$$\check{H}^*(X, G) \approx \check{H}^*(\mathfrak{U}, G);$$
in particular, the right-hand side does not depend on the particular good cover $\mathfrak{U}$. 

\medskip

By \cite[Proposition 15.8]{BT0} (see also \cite[page 257]{Hat0}) we have that
$$\check{H}^*(X,\Z)\approx H^*_{\Delta}(X, \Z)$$
whenever $X$ is a manifold. %(which will be our case, see Section \ref{...}). 
Moreover \cite[Corollary 6.9.9]{Spa0} 
$$\check{H}^*(X, G) \approx \bar{H}^*(X,G)$$
whenever $X$ is a closed subset of a manifold (or more generally, $X$ is a %(paracompact) 
Hausdorff space with coefficients in a module $G$ \cite[Corollary 6.8.8]{Spa0}). %; actually the paracompact assumption is not needed,

Once defined also the \emph{relative} \v{C}ech cohomology, one can prove \cite[pages 342 and 359]{Spa0} 
$$\check{H}^*(X, Y, G) \approx \bar{H}^*(X, Y, G)$$
whenever $X,Y$ are closed subset of a manifold. Thus, by combining this result with \eqref{eq_isom_Cech-AS_rel}, whenever $Y$ and $X$ are closed subsets of a manifold %$\R^N$ 
(such as $\R^N$ or a more general Hilbert space, %which will be our case, 
see Section \ref{sec_cup-length}) and $G=R$ ring, 
we have
$$\check{H}^*(X,Y,G)\approx H^*_{\Delta}(X, Y, G);$$
see also \cite[Proposition 8.6.12]{Dol0} for a direct proof in the case of a pair of \emph{ENR} (\emph{Euclidean Neighborhood Retracts}, which is the case for example of $Y\subset X\subset \R^N$ with $X$ retractible). 

\medskip

See also \cite{Bre0} for further relations on these three cohomologies.

\subsection{Properties of the cup-length}

Here we focus on the case $G:=\mathbb{F}$ for some field $\mathbb{F}$, and we drop the dependence on $G$ in the notations.
We collect some properties of the cup-length, see e.g. \cite[Lemma 2.6]{BaWe0}.

\begin{Lemma}\label{lem_propr_cupl}
We have the following properties.
\begin{itemize}
\item[(a)] For any $f:(A,B)\to (A', B')$ and $f':(A', B')\to (A'', B'')$ it results that
$$\cupl(f' \circ f) \leq \min\{ \cupl(f'), \cupl (f)\}.$$
As a consequence,
\begin{equation}\label{eq_propr_cupl}
\cupl(f' \circ f) \leq \cupl(A', B').
\end{equation}
\item[(b)] For any $f, g:(A,B) \to (A', B')$ homotopic, we have
$$\cupl(f)=\cupl(g).$$
\end{itemize}
\end{Lemma}

Finally, we cite the following key result \cite{Bar1} which can be found in \cite[Lemma 5.5]{CJT}.
\begin{Lemma}\label{lem_cupl_jK}
Consider the inclusion 
$$j:(I\times K, \partial I \times K) \to (I\times K_d, \partial I \times K_d)$$
for a whatever $K\subset \R^N$ compact, $K_d:=\{ x \in \R^N \mid d(x,K)\leq d\}$, and $I=[a,b]$. Then, for $d>0$ sufficiently small, we have
$$\cupl(j) \geq \cupl(K).$$
\end{Lemma}

\subsection{Relation with the Ljusternik-Schnirelmann category}
%Ljusternik and Schnirelmann

We recall here the definition of \emph{relative category}, by following \cite{Szu0, FLRW} and references therein (see also \cite{Bar0}).
% for a complete treatment on the topic of category.

\begin{Definition}
Let $X$ be a topological space and let $A, B$ be two closed subsets of $X$. We call the \emph{category of $A$ in $X$, relative to $B$}, and write 
$$k=\cat_{X, B}(A),$$
 the least integer $k\in \N$ such that there exist $A_0, A_1, \dots, A_k$ closed subsets of $X$ which verify
\begin{itemize}
\item $(A_i)_{i=0\dots k}$ cover $A$;
\item $(A_i)_{i=1 \dots k}$ are contractible in $X$, i.e. $id: A_i \to X$ is homotopic to a constant;
\item $A_0$ is deformable in $B$, i.e. there exists a continuous $h_0:[0,1]\times (A_0 \cup B) \to X$ such that $h_0(0, \cdot)=id$, $h_0(1, A_0)\subset B$ and $h_0(t, B)\subset B$ for each $t \in [0,1]$.
\end{itemize}
If such $k$ does not exists, we set $\cat_{X, B}(A):=+\infty$. %Notice that $\cat_X(A)=\cat_{X, \emptyset}(A)$.
\end{Definition}

Examples of computations can be found in \cite[Examples 2.2 and 3.7]{FLRW} and \cite[Remark 3.2]{FoWi1}.
For example, if $B$ is the unit ball in $\R^N$, then $\cat_{B, \partial B}(B)=1$ (while it is equal zero if $B$ is the unit ball in $H^s(\R^N)$); if $A$ is the annulus in $\R^N$ with $N\geq 2$, then $\cat_{A,\partial A}(A)=2$; moreover $\cat_{\R^2, \R}(\R^2)=\cat_{\R^2, (0,0)}(\R^2)=0$.

\begin{Remark}
$\,$
\begin{itemize}
\item If we drop the condition on $B$, $A_0$ and $h_0$, we have the classical definition of category, and simply write $\cat_X(A)$; more precisely
$$\cat_X(A)=\cat_{X, \emptyset}(A).$$
This definition can be given for a whatever $A$ (even not closed), and a posteriori one has $\cat_X(A)=\cat_X(\overline{A})$ (see \cite[Remark 1.12]{CLOT}).
\item We required the covering to be closed, but equivalently one can ask $A_0\dots A_k$ to be open (see \cite[Proposition 1.10]{CLOT}).
\item We do not require that $B\subset A_0 \subset A$, even if equivalent definitions could be given in this way (see e.g. \cite{FLRW}).
\item Some authors require the stronger condition $h_0(t, \cdot_{|B})=id_B$ (see e.g. \cite{CJT, FoWi1, FoWi2} and Remark 2.2 in \cite{Szu0}), and this modification would bring no differences in what follows. %(\tr{equivalente? VEDI}).
\item Considered a continuous map $f:(A,B) \to (A',B')$ one can define the category of $f$ by substituting, in the definition (with $A=X$), "$id: A_i \to A$" with "$f_{|A_i} :A_i \to A$", "$h_0:[0,1]\times (A_0 \cup B) \to A$" with "$h_0:[0,1]\times (A_0 \cup B) \to A$" and "$h_0(0, \cdot)=id$" with "$h_0(0, \cdot)=f$"; in this case $\cat_{A,B}(A)= \cat(id_{(A,B)})$. See \cite{BaWe0}. Anyway, we will not use this tool.
\end{itemize}
\end{Remark}

%Since we will work only with $X=\mc{X}_{\eps, \delta}$, to avoid cumbersome notation we will write $\cat(A):= \cat_{\mc{X}_{\eps, \delta}}(A)$ for the category and $\cat(A,B):=\cat_{\mc{X}_{\eps, \delta}, B}(A)$ for the relative category.

The following classical properties on category can be found, e.g., in \cite[Lemma 2.2]{BaWe0} and \cite[Lemma 1.13]{CLOT} (see also \cite[Proposition 2.9]{FLRW}).
\begin{Lemma}\label{prop_propriet_cat} 
Let $A$ be a closed subset of $X$. %\tr{be such that...} %We have the following
\begin{itemize}
\item $\#A \geq \cat_X(A)$;
%\item $\cat(\overline{A})=\cat(A)$;
\item If $A$ is compact, and every point in $A$ has an open neighborhood in $X$ contractible in $X$, then there exists an open neighborhood $N \subset X$ of $A$ such that $\cat_X(N)=\cat_X(A)$. In particular, if $A\subset X \subset X' \subset H$, with $A$ compact and $X$ open subset of the Hilbert space $H$, then the claim holds true for $\cat_{X'}(A)$.
%va applicato a $X=J^{\bar{c}+\delta}$, che è un chiuso in uno spazio di Hilbert
\end{itemize}
\end{Lemma}

Next proposition deals with some properties on relative category, and can be found, for instance, in \cite[Propositions 2.5 and 2.8]{Szu0} or \cite[Propositions 2.4 and 2.9]{FLRW} (see also \cite[Remark 3.2 and Propositions 3.4 and 3.5]{FoWi1}).
\begin{Lemma}\label{prop_proprieta_cat_rel}
Let $A,A', B, V$ be closed subsets of $X$. %We have the following
\begin{itemize}
\item %$\cat(A, A)=0$;
%$\cat_{X,A}(A)=0$; moreover, 
Then $\cat_{X,B}(A)=0$ if and only if $A$ can be deformed in $B$, i.e. there exists $h:(A\cup B) \times [0,1] \to X$ such that $h(0, \cdot)=id$, $h(t, B) \subset B$ %$h(t, \cdot_{|B})=id$ 
for each $t \in [0,1]$ and $h(1, A) \subset B$. As a consequence, if $A \subset B$, then $\cat_{X,B}(A)$. In particular, $\cat_{X,A}(A)=0$.
\item If $A\subset A'$, then %$\cat(A, B) \leq \cat (A', B)$;
$\cat_{X,B}(A) \leq \cat_{X,B}(A')$.
\item If $A \cup B \subset X \subset X'$, then $\cat_{X,B}(A) \geq \cat_{X',B}(A)$. In particular, if $B \subset A \subset X$, then $\cat_{A, B}(A) \geq \cat_{X,B}(A)$.
\item If %$\cat(V)<\infty$, then $\cat(\overline{A\setminus V}, B) \geq \cat(A, B) - \cat(V)$;
$\cat_X(V)<\infty$, then $\cat_{X,B}(\overline{A\setminus V}) \geq \cat_{X,B}(A) - \cat_X(V)$.
\item If there exists $\eta: [0,1] \times (A\cup B) \to X$ %\tr{\mc{X}_{\eps, \delta}}.
such that $\eta(1, A)\subset A'$ and $\eta([0,1], B) \subset B$, then %$\cat(A, B) \leq \cat(A', B)$.
$\cat_{X,B}(A) \leq \cat_{X,B}(A')$.
\end{itemize}
\end{Lemma}

The following lemma links the concepts of category (when $A=X$) and cup-length, and it can be found in \cite[Proposition 2.6 and Remark 2.7]{Szu0} (see also \cite[Theorem 1]{FoWi2} and \cite[Theorem 3.6]{FLRW}).
\begin{Lemma}\label{lem_coll_cat_cupl}
%For any $B\subset A \subset \R^N$ closed, we have
Let $B$ be a closed subset of a metric space $A$. Then %For any $B\subset A \subset \R^N$ closed, we have
%$$\cat(A,B) \geq \cupl(A,B)+1.$$
$$\cat_{A,B}(A) \geq \cupl(A,B)+1.$$
In particular, $\cat_{A}(A) \geq \cupl(A)+1$.
\end{Lemma}

%Since we will work only with $X=\mc{X}_{\eps, \delta}$, 
To avoid cumbersome notation we will write 
%$\cat(A):= \cat_{\mc{X}_{\eps, \delta}}(A)$ and $\cat(A,B):=\cat_{\mc{X}_{\eps, \delta}, B}(A)$.
$$\cat(A):= \cat_{A}(A), \quad \hbox{and} \quad \cat(A,B):=\cat_{A, B}(A).$$
%or, if $X$ is clear from the framework,
%$$\cat(A):= \cat_{X}(A), \quad \hbox{and} \quad \cat(A,B):=\cat_{X, B}(A).$$
Notice that, if $A \subset X$, then $\cat(A,B) \geq \cat_{X,B}(A)$ (and in particular $\cat(A) \geq \cat_X(A)$). 

\begin{Remark}\label{examples}
%Notice that the cup-length of a set $K$ is strictly related to the Lusternik-Schnirelmann category of $K$. 
We notice that in standard examples the inequality in Lemma \ref{lem_coll_cat_cupl} is actually an equality. 
Indeed, if $K$ is a contractible set or it is finite (e.g. a single point), then 
$$\cupl(K)+1=\cat_K(K)=1;$$
if $K=S^{N-1}$ is the $N-1$ dimensional sphere in $\R^N$, then $$\cupl(K) + 1= \cat_K(K) =2;$$ if $K=T^N$ is the $N$-dimensional torus, then $$\cupl(K) + 1 = \cat_K(K)= N + 1.$$
However in general the strict inequality may hold, % $$\cupl(K) +1 \leq \cat(K)$$ (
see \cite[Sections 2.8 and 9.23]{CLOT} for some examples. % where the strict inequality is attained).
\end{Remark}

\begin{Remark}\label{rem_monotonia}
When one deals with a functional which is not bounded from below, the tool of the relative category is needed. 
On the other hand, for any interval $I\subset \R$ and any neighborhood $K_d$ of $K$, considered the inclusion $$j:(I \times K, \partial I \times K) \to (I\times K_d, \partial I \times K_d)$$ the key relation
$$\cat(j) \geq \cat_K(K),$$
essential in the estimation of the relative category of two sublevels of the indefinite functional (see \cite[Remark 4.3]{CJT}) does not generally hold \cite[Remark 7.47]{CLOT}. Nevertheless, the same relation for the cup-length
$$\cupl(j) \geq \cupl(K)$$
holds true, as proved in \cite[Lemma 5.5]{CJT} (see also \cite[Proposition 3.5]{FLRW}). 
That is why we take advantage of the relative cup-length in order to get a bound on the number of solutions.
\end{Remark}

\subsection{Application to multiplicity of solutions}
\label{sec_mult_cat_cupl}

We sketch now how to obtain multiple solutions from the information on the category of a set.

Let indeed $J: X \to \R$ to be a $C^1$-functional on a function space $X$, and denote, for every $c\in \R$, $J^c:=\{J\leq c\}$ the sublevel at $c$ and $K_c:=\{J=c, \, J'=0\}$ the set of critical points at $c$. Assume the following:
\begin{itemize}
\item there exist $\bar{c} \in \R$ and $\delta>0$ such that 
%$J^{\bar{c}+\delta}$ is not deformable in $J^{\bar{c}-\delta}$ (in particular $[ \bar{c}-\delta < J \leq \bar{c}+\delta]\neq \emptyset$); 
%\item 
$K_c$ is compact for every $c \in [\bar{c}-\delta, \bar{c}+\delta]$ (for example, a Palais-Smale type condition holds at level $c$) and a Deformation Lemma holds around $K_c$; %$J^{\bar{c}\pm \delta}$ ...
\item $\bar{c}+\delta$ is a regular value; this is not restrictive, up to choosing properly $\delta$ (small), since otherwise we would have a sequence of critical values at levels $c+\delta_n$ with $\delta_1> ... >\delta_n \to 0$.
\item there exist a compact $K$ and two continuous maps $\phi_1, \phi_2$ such that $(I\times K, \partial I \times K) \stackrel{\phi_1} \to (J^{\bar{c}+\delta}, J^{\bar{c}-\delta}) \stackrel{\phi_2} \to (I\times K_d, \partial I \times K_d)$ is homotopic to the inclusion $j:(I\times K, \partial I \times K) \to (I\times K_d, \partial I \times K_d)$, where $K_d=\{ x \in \R^N \mid d(x,K)\leq d\}$ and $I=[a,b]$ for some $a,b \in \R$.
\end{itemize}
%By construction of the neighborhood $\mc{X}_{\eps, \delta}$ and Corollary \ref{corol_passaggio_critici} (recall that $\rho_0 < r_3 \leq r_2'$ and that $J_{\eps}(u) < E_{m_0} + R(\hat{\delta}, u) \leq E_{m_0} + \hat{\delta} \leq l_0'$ for $u \in \mc{X}_{\eps, \delta}$), we have
%$$\left\{ u \in (\mathcal{X}_{\eps, \delta})^{E_{m_0}+\hat{\delta}}_{E_{m_0}-\hat{\delta}} \mid J'_{\eps}(u)=0 \right\} \subset \left\{ u \in \mathcal{X}_{\eps, \delta} \mid J'_{\eps}(u)=0\right\} \subset \left\{ u \in H^s(\R^N) \mid I'_{\eps}(u)=0 \right\}.$$
%$$\left\{ u \in \mathcal{X}_{\eps, \delta}^{E_{m_0}+\hat{\delta}} \setminus \mathcal{X}_{\eps, \delta}^{E_{m_0}-\hat{\delta}} \mid J'_{\eps}(u)=0 \right\} \subset \left\{ u \in H^s(\R^N) \mid I'_{\eps}(u)=0 \right\}.$$

We want to show
\begin{align*}
%\# \{ u \textnormal{ solutions of \eqref{eq_princ_epsx}}\} &\geq& 
\# \left\{ u \in X \mid J(u) \in [\bar{c}-\delta, \bar{c}+\delta], \; J'(u)=0 \right\} %\left\{ u \in \mathcal{X}_{\eps, \delta}^{E_{m_0}+\hat{\delta}} \setminus \mathcal{X}_{\eps, \delta}^{E_{m_0}-\hat{\delta}} \mid J'_{\eps}(u)=0 \right\}
 &\stackrel{(i)}\geq \cat_{J^{\bar{c}+\delta}, J^{\bar{c}-\delta}} (J^{\bar{c}+\delta}) \\
&\stackrel{(ii)}\geq \cupl\left (J^{\bar{c}+\delta},\, J^{\bar{c}-\delta}\right) +1 \\
&\stackrel{(iii)}\geq \cupl(K) +1.
\end{align*}
%that is the claim, up to the proof of (i)--(iii).
which is an estimate from below on the number of critical points of $J$.

\smallskip

\textbf{Proof of (i).} This is a consequence of the Deformation Lemma and of the compactness of critical level sets, as done in \cite[Proposition 3.2]{Szu0} and \cite[Theorem 4.2]{FoWi1}
 %\ref{lem_def_lem} 
(see also \cite[Theorem 3]{FoWi2} and \cite[Theorem 6.1]{FLRW}). 
% We write it here for the sake of completeness and to highlight the variational formulation of the level sets of the critical points. 
Let thus define
$$k:= \cat_{J^{\bar{c}+\delta}, J^{\bar{c}-\delta}} (J^{\bar{c}+\delta}) \in \N \cup \{+\infty\};$$
if $k=0$ the claim is trivial, thus we can assume $k\geq 1$.
%notice that $k\neq 0$ since $J^{\bar{c}+\delta}$ is not deformable in $J^{\bar{c}-\delta}$. 
% $[ \bar{c}-\delta < J \leq \bar{c}+\delta]\neq \emptyset$ and thus $J^{\bar{c}+\delta} \not \subseteq J^{\bar{c}-\delta}$.
For each $j= 1 \dots k$ define %\tr{(VEDI ZERO)}
$$\Gamma_j:= \{ A\subset X \mid A \hbox{ closed}, \, \cat_{J^{\bar{c}+\delta}, J^{\bar{c}-\delta}} (A) \geq j\},$$
$$c_j :=\inf_{A \in \Gamma_j} \sup_{A} J.$$
Notice that, since $j\geq 1$, then each $A\in \Gamma_j$ cannot be included in $J^{\bar{c}-\delta}$, that is $c_j\geq \bar{c}-\delta$; moreover, since $j\leq k$, then $J^{\bar{c}+\delta}\in \Gamma_j$, which implies $c_j \leq \bar{c}+\delta$. Therefore
$$\bar{c}-\delta \leq c_1 \leq c_2 \leq \dots \leq \bar{c}+\delta.$$
Fix $j \in \{1 \dots k\}$ and let $p\in \N$ be such that
$$c_j = \dots = c_{j+p}=:c \in [\bar{c}-\delta, \bar{c}+\delta];$$ % \subset (\bar{c}-\delta', \bar{c}+\delta');$$
to reach the claim, it is sufficient to show that
\begin{equation}\label{eq_dim_Kc_p}
 \cat_{J^{\bar{c}+\delta}}(K_c) \geq p+1
\end{equation}
since $\# K_c \geq \cat_{J^{\bar{c}+\delta}}(K_c)$ and by combining the estimates for different values of $c_j$ (if $c_i \neq c_j$ we clearly have different critical points at the two levels).

We do some preliminary work.
We first exploit that $\bar{c}+\delta$ is a regular point to show that $c<\bar{c}+\delta$. 
%Indeed, if $c_j= \bar{c}+\delta$ for some $j$, then $K_{c_j}=\emptyset$ and thus by the deformation lemma there exists $\eta:[0,1]\times X \to X$ and an $\omega \in (0,\delta)$ such that $\eta(1,J^{c_j + \delta}) \subset J^{c_j-\omega}$ and $J(t,J^{\bar{c}-\delta}) \subset J^{\bar{c}-\delta}$ for each $t\in [0,1]$: by Lemma \ref{prop_proprieta_cat_rel} we have
%$$\cat_{J^{\bar{c}+\delta}, J^{\bar{c}-\delta}}(J^{c_j+\omega}) \leq \cat_{J^{\bar{c}+\delta}, J^{\bar{c}-\delta}} (J^{c_j-\omega}).$$
%On the other hand, since $c_j+\omega = \bar{c}+\delta+\omega \geq \bar{c}+\delta$, by the monotonicity of the category we have
%$$\cat_{J^{\bar{c}+\delta}, J^{\bar{c}-\delta}}(J^{c_j-\omega}) \geq \cat_{J^{\bar{c}+\delta}, J^{\bar{c}-\delta}} (J^{\bar{c}+\delta+\omega }) \geq \cat_{J^{\bar{c}+\delta}, J^{\bar{c}-\delta}} (J^{\bar{c}+\delta}) = k \geq j;$$
%thus $J^{c_j-\omega} \in \Gamma_j$, which implies $c_j \leq c_j-\omega$, absurd.
Indeed, since $K_{\bar{c}+\delta}=\emptyset$, by the Deformation Lemma there exist $\eta:[0,1]\times X \to X$ and an $\omega >0 $ % \in (0,\delta)$
 such that 
\begin{itemize}
\item $J(\eta(t,\cdot))\leq J(\cdot)$ for each $t\in [0,1]$, and thus $\eta:[0,1]\times J^{\bar{c}+\delta} \to J^{\bar{c}+\delta}$, 
\item $\eta(1,J^{\bar{c}+\delta+\omega}) \subset J^{\bar{c}+\delta-\omega}$, and thus $\eta(1,J^{\bar{c}+\delta}) \subset J^{\bar{c}+\delta-\omega}$,
\item $J(t,J^{\bar{c}-\delta}) \subset J^{\bar{c}-\delta}$ for each $t\in [0,1]$;
\end{itemize}
by Lemma \ref{prop_proprieta_cat_rel} we have %and the monotonicity of the category %we have
%$$\cat_{J^{\bar{c}+\delta}, J^{\bar{c}-\delta}}(J^{\bar{c}+\delta+\omega}) \leq \cat_{J^{\bar{c}+\delta}, J^{\bar{c}-\delta}} (J^{c_j-\omega}).$$
%On the other hand, since $c_j+\omega = \bar{c}+\delta+\omega \geq \bar{c}+\delta$, by the monotonicity of the category we have
$$\cat_{J^{\bar{c}+\delta}, J^{\bar{c}-\delta}}(J^{\bar{c}+\delta-\omega}) \geq 
%\cat_{J^{\bar{c}+\delta}, J^{\bar{c}-\delta}} (J^{\bar{c}+\delta+\omega }) \geq 
\cat_{J^{\bar{c}+\delta}, J^{\bar{c}-\delta}} (J^{\bar{c}+\delta}) = k \geq j;$$
thus $J^{\bar{c}+\delta-\omega} \in \Gamma_j$, which implies $c_j \leq \bar{c}+\delta-\omega < \bar{c}+\delta$ for each $j$, which is the claim.

Since $c<\bar{c}+\delta$, we have $K_c \subset \{ J < \bar{c}+\delta\} \subset J^{\bar{c}+\delta}$; moreover
%By the assumptions, %the Palais-Smale condition (see Proposition \ref{prop_Palais_Smale}) and definition of $\mc{X}_{\eps, \delta}$, 
%we have that 
$K_c$ is compact; thus by Lemma \ref{prop_propriet_cat} we have that there exists an open neighborhood $N$ of $K_c$ such that 
$$\cat_{J^{\bar{c}+\delta}}(N)=\cat_{J^{\bar{c}+\delta}}(K_c).$$
Corresponding to $N$, again by the Deformation Lemma there exist an $\omega \in (0, \bar{c}+\delta-c) $ and an
% $\eta:[0,1]\times X \to X$ such that $\eta(1, J^{c+\omega}\setminus N) \subset J^{c-\omega}$ and $\eta(t, J^{\bar{c}-\delta}) \subset J^{\bar{c}-\delta}$ % is the identity on $\mc{X}_{\eps, \delta}^{E_{m_0}-\hat{\delta}}$) 
%for each $t\in [0,1]$. % (actually $\eta(t, \cdot_{|J^{\bar{c}-\delta}})$ is the identity). 
 $\eta:[0,1]\times J^{\bar{c}+\delta} \to J^{\bar{c}+\delta}$ (notice that $J^{c+\omega} \cup J^{\bar{c}-\delta} \subset J^{\bar{c}+\delta}$) such that $\eta(1, J^{c+\omega}\setminus N) \subset J^{c-\omega}$ and $\eta(t, J^{\bar{c}-\delta}) \subset J^{\bar{c}-\delta}$ % is the identity on $\mc{X}_{\eps, \delta}^{E_{m_0}-\hat{\delta}}$) 
for each $t\in [0,1]$. % (actually $\eta(t, \cdot_{|J^{\bar{c}-\delta}})$ is the identity). 
By Lemma \ref{prop_proprieta_cat_rel} we have
\begin{equation}\label{eq_dim_esist_sol1}
\cat_{J^{\bar{c}+\delta}, J^{\bar{c}-\delta}}(J^{c+\omega}\setminus N) \leq \cat_{J^{\bar{c}+\delta}, J^{\bar{c}-\delta}}( J^{c-\omega} ).
\end{equation}
Corresponding to $\omega$, by definition of $c=c_{j+p}$ there exists an $A\in \Gamma_{j+p}$ such that $\sup_{A}J < c+\omega$, which means that $A\subset J^{c+\omega}$ and thus 
\begin{equation}\label{eq_dim_esist_sol2}
A \setminus N \subset J^{c+\omega} \setminus N.
\end{equation}
%Since the right-hand side is closed, we actually have
%$$\overline{A\setminus N} \subset \mc{X}_{\eps, \delta}^{c+\omega} \setminus N.$$

We prove now \eqref{eq_dim_Kc_p} by contradiction. Assume $\cat(K_c)\leq p <\infty$. Thus, by \eqref{eq_dim_esist_sol1}, \eqref{eq_dim_esist_sol2} and Lemma \ref{prop_proprieta_cat_rel} (notice that $A\setminus N$ is closed) we have %, since $\cat(N)=\cat(K_c)$,
\begin{align*}
 \cat_{J^{\bar{c}+\delta}, J^{\bar{c}-\delta}}(J^{c-\omega} ) &\geq \cat_{J^{\bar{c}+\delta}, J^{\bar{c}-\delta}}(J^{c+\omega}\setminus N) 
\geq \cat_{J^{\bar{c}+\delta}, J^{\bar{c}-\delta}}(A \setminus N) \\
&\geq \cat_{J^{\bar{c}+\delta}, J^{\bar{c}-\delta}}(A) - \cat_{J^{\bar{c}+\delta}}(N)
\geq (j+p)-p
= j.
\end{align*}
This means that $J^{c-\omega} \in \Gamma_j$, and thus
$$c_j \leq \sup_{J^{c-\omega}} J \leq c-\omega = c_j-\omega$$
which is an absurd.

\smallskip

\textbf{Proof of (ii).} This is a consequence of the property of algebraic topology given in Lemma \ref{lem_coll_cat_cupl}.

\smallskip

\textbf{Proof of (iii).} This is due to the existence of the homotopy %gained in Proposition \ref{prop_esist_homot} 
and properties of the cup-length. Indeed,
by \eqref{eq_propr_cupl} in Lemma \ref{lem_propr_cupl} (a), we have
$$\cupl\left (J^{\bar{c}+\delta},\, J^{\bar{c}-\delta}\right) \geq \cupl (\phi_2\circ \phi_1);$$
we highlight that the left-hand side deals with subsets of the function space $X$, while the right-hand side deals with subsets of $\R^N$.
%observed that
%$$\Psi_{\eps} \circ \Phi_{\eps} = \Psi_{\eps} \circ id_{\left(\mathcal{X}_{\eps, \delta}^{E_{m_0}+\hat{\delta}}, \, \mathcal{X}_{\eps, \delta}^{E_{m_0}-\hat{\delta}}\right)} \circ \Phi_{\eps}$$
%we have, by definition of $\cupl(A,B)$ and \ref{...}
%\begin{eqnarray*}
%\cupl\left (\mathcal{X}_{\eps, \delta}^{E_{m_0}+\hat{\delta}},\, \mathcal{X}_{\eps, \delta}^{E_{m_0}-\hat{\delta}}\right) &=& \cupl\left (id_{\left(\mathcal{X}_{\eps, \delta}^{E_{m_0}+\hat{\delta}}, \, \mathcal{X}_{\eps, \delta}^{E_{m_0}-\hat{\delta}}\right)}\right) \\
%&\geq & \min \left \{ \cupl(\Psi_{\eps}), \cupl\left (id_{\left(\mathcal{X}_{\eps, \delta}^{E_{m_0}+\hat{\delta}}, \, \mathcal{X}_{\eps, \delta}^{E_{m_0}-\hat{\delta}}\right)}\right), \cupl (\Phi_{\eps}) \right\} \\
%&\geq& \cupl \left( \Psi_{\eps} \circ id_{\left(\mathcal{X}_{\eps, \delta}^{E_{m_0}+\hat{\delta}}, \, \mathcal{X}_{\eps, \delta}^{E_{m_0}-\hat{\delta}}\right)} \circ \Phi_{\eps} \right) \\
%&=& \cupl (\Psi_{\eps} \circ \Phi_{\eps}).
%\end{eqnarray*}
Since $\phi_2 \circ \phi_1 $ is homotopic to the immersion $j$, we have by Lemma \ref{lem_propr_cupl} (b)
$$ \cupl (\phi_2\circ \phi_1) = \cupl(j).$$
Finally, we conclude thanks to Lemma \ref{lem_cupl_jK}.
\QED

\medskip

\subsection{The Krasnoselskii genus: a particular category}
\label{sec_app_genus}

In order to obtain existence of multiple solutions in the entire space $\R^N$, without any topology related to some potential $V$, it is useful to exploit some symmetry of the functionals, and some tool related to them.

In particular, we introduce the well known Krasnoselskii genus.

\begin{Definition}
For any $A$ closed subset of $\R^N \setminus \{0\}$, symmetric with respect to the origin (i.e. $A=-A$), the \emph{Krasnoselskii genus} is defined by
$$\genus(A):= \max \left\{ n \in \N \mid \exists \beta:A \to \R^n\setminus\{0\} \hbox{ continuous and odd}\right\};$$
if such $n$ does not exists, $\gamma(A):=+\infty$; moreover $\gamma(A)=0$ if (and only if) $A=\emptyset$.
\end{Definition}

The genus enjoys several standard properties \cite[Section 3]{Rab-1}.

\begin{Proposition}\label{prop_genus_g}
Let $A, B \subset \R^N\setminus \{0\}$ be closed and symmetric.
\begin{itemize}
%\item $\gamma(A)\leq N$;
\item if $A$ is finite, then $\genus(A)=1$;
\item $\genus(A\cup B) \leq \genus(A) + \genus(B)$;
\item if $\genus(B)<\infty$, then $\genus(\overline{A\setminus B}) \leq \genus(A) - \genus(B)$;
\item if $h: \R^N \to \R^N$ is continuous and odd, then $\genus(A)\leq \genus(\overline{h(A)})$;
\item if $A$ is compact, then there exists a closed, symmetric neighborhood $U \not \ni 0$ of $A$ such that $\genus(U)=\genus(A)<\infty$;
\item if $U$ is a symmetric neighborhood of the origin, then $\genus(\partial U)=N$.
\end{itemize}
\end{Proposition}

\begin{Example}
The genus describes, roughly, how a set is wrapped near the origin. Let $A$ be a closed subset of $\R^N\setminus \{0\}$, such that $A=-A$.
If $A=B \cup (-B)$, with $B \cap (-B)=\emptyset$, then $\genus(A)=1$. If $A$ is connected, then $\genus(A)\geq 2$. Moreover, $\genus(S^N)=N+1$.
\end{Example}

Actually, this tool reveals to be a subcase of the already introduced category. Indeed, considered the action of $\Z_2$ over $\R^N$ (which identifies $x$ with $-x$) we have the following relation \cite[Theorem 3.7]{Rab-1} (see also \cite{Fad0})
$$\genus(A) = \cat_{(\R^N\setminus \{0\})/\Z_2} (A/\Z_2).$$
This relation highlights the fact that the genus tool exploits not the topology of a particular subset of $\R^N$, but the topology induced by a symmetry relation.

%\bigskip

\bigskip

\bigskip

%\smallskip

\noindent
%\section*{{\normalsize \raggedleft{Acknowledgments}}}
\textbf{Acknowledgments.} The number of people to thank for this %wonderful 
three years journey is quite long, starting from my advisors, passing through collaborators and colleagues, up to family and friends. 
I have written a truly marvelous praise about this, which this margin unfortunately is too narrow to contain.

\smallskip

%\bigskip

\rightline{\emph{So long, and thanks for all the fish!}}

%\medskip
%
%\section*{\large{\raggedleft{Acknowledgments}}}
%The number of people to thank for this %wonderful 
%three years journey is quite long, starting from my advisors, passing through collaborators and colleagues, up to family and friends. 
%I have written a truly marvelous praise about this, which this margin unfortunately is too narrow to contain.
%
%\bigskip
%
%\bigskip
%
%\rightline{\emph{So long, and thanks for all the fish!}}

 %PER LA STAMPA

\newpage

\thispagestyle{empty}

\newpage

\small

\cleardoublepage
\phantomsection
\addcontentsline{toc}{chapter}{Bibliography}
\bibliographystyle{amsplain}

%VERSIONE STAMPATA
%\fancyhead[LE, RO]{\thepage}
%\fancyhead[RE]{\nouppercase{\leftmark}}
%\fancyhead[LO]{\nouppercase{\rightmark}}

\fancyhead[L]{ \ifthenelse{\isodd{\value{page}}}{\nouppercase{\leftmark}}{\thepage} }
\fancyhead[R]{ \ifthenelse{\isodd{\value{page}}}{\thepage}{\nouppercase{\leftmark}} }

\bigskip

\normalsize

%VERSIONE STAMPATA
%\newpage
%\thispagestyle{empty}
%\newgeometry{paperwidth=17cm, paperheight=24cm, textheight=225mm, textwidth=80mm, top=20mm, left=7mm, right=7mm}

%\section*{Fundings}
\textbf{Fundings:} 
The research in this thesis has been supported by INdAM-GNAMPA and by PRIN 2017JPCAPN ``Qualitative and quantitative aspects of nonlinear PDEs''.

%\newpage

\bigskip

\textbf{Mathematics Subject Classification:} 

\begin{tabular}{lllllllll}
30H40, % Zygmund spaces
&
33C05, %Classical hypergeometric functions, 2F1
&
35A01, %Existence problems for PDEs: global existence, local existence, non-existence
&
35A09, %Classical solutions to PDEs
&
35A15, %Variational methods applied to PDEs
&
35A21, %Singularity in context of PDEs
&
35B06, %Symmetries, invariants, etc. in context of PDEs
& 
35B09, %Positive solutions to PDEs
& 
35B25, %Singular perturbations in context of PDEs
\\
35B33, %Critical exponents in context of PDEs
&
35B38, %Critical points of functionals in context of PDEs (e.g., energy functionals)
& 
35B40, %Asymptotic behaviour of solutions to PDEs
& 
35B50, %Maximum principles in context of PDEs
&
35B65, %Smoothness and regularity of solutions to PDEs
&
35D30, %Weak solutions to PDEs
&
35D35, %Strong solutions to PDEs
&
35D40, %Viscosity solutions to PDEs
&
35G20, %Nonlinear higher-order PDEs
\\
35G50, %Systems of nonlinear higher-order PDEs
&
35J05, %Laplace operator, Helmholtz equation (reduced wave equation), Poisson equation
&
35J10, %Schrödinger operator, Schrödinger equation
&
35J15, %Second-order elliptic equations
&
35J20, %Variational methods for second-order elliptic equations
& 
35J47, %Second-order elliptic systems
&
35J50, %Variational methods for elliptic systems
&
35J60, %Nonlinear elliptic equations
&
35J61, %Semilinear elliptic equations
\\
35J91, %Semilinear elliptic equations with Laplacian, bi-Laplacian or poly-Laplacian
&
35Q40, %PDEs in connection with quantum mechanics
& 
35Q55, %NLS equations (nonlinear Schrödinger equations)
& 
35Q60, %PDEs in connection with optics and electromagnetic theory
& 
35Q70, %PDEs in connection with mechanics of particles and systems of particles
& 
35Q75, %PDEs in connection with relativity and gravitational theory
&
35Q85, %PDEs in connection with astronomy and astrophysics
&
35Q92, %PDEs in connection with biology, chemistry and other natural sciences
&
35R09, %Integro-partial differential equations
%& 
%35B09, %Positive solutions to PDEs
\\ 
35R11, %Fractional partial differential equations
& 
44A10, %Laplace transform
&
45G15, %Systems of nonlinear integral equations
&
45K05, %Integro-partial differential equations
&
45M05, %Asymptotics of solutions to integral equations
& 
45M20, %Positive solutions of integral equations
& 
46M20, %Methods of algebraic topology in functional analysis (cohomology, sheaf and bundle theory, etc.) 
&
47F05, %General theory of partial differential operators
&
47F10, %Elliptic operators and their generalizations
\\
47G10, %Integral operators
&
47H10, %Fixed-point theorems
&
47J30, %Variational methods involving nonlinear operators
& 
49J35, %Existence of solutions for minimax problems
&
55M30, %Lyusternik-Shnirel’man category of a space, topological complexity `a la Farber, topological robotics (topological aspects)
&
55N05, %Cech types
&
55N10, %Singular homology and cohomology theory
&
58E05, %Abstract critical point theory (Morse theory, Lyusternik-Shnirel’man theory, etc.) in infinite-dimensional spaces
&
58J05. %Elliptic equations on manifolds, general theory
\end{tabular}

\bigskip

\textbf{Key words:} 
%Asymptotic behaviour;
%Center of mass;
%Choquard equation;
%Concentration phenomena;
%Critical exponent;
%Relative cup-length;
%Double nonlocality;
%Existence and multiplicity; % of solutions;
%Even and odd nonlinearities,
%Fractional Laplacian;
%%Fractional nonlinear equation;
%Ground states;
%Hartree-type term;
%$L^2$-constraint; 
%Lagrang multiplier; %formulation and
%%Multidimensional odd 
%Mountain Pass paths;
%Nonlocal sources;
%Schrödinger equation;
%Nonlinear PDEs;
%Normalized solutions;
%Pohozaev identity;
%Polynomial decay;
%Positivity and sign;
%Prescribed mass;
%Qualitative properties; %of the solutions;
%Radial symmetry; %Radially symmetric solutions
%Regularity;
%Riesz potential;
%%Sign of the ground states;
%Singular perturbation;
%Spike solutions;
%Sublinear nonlinearities.
%%Symmetric solutions; %Radially symmetric solutions

\begin{tabular}{lll}
Asymptotic behaviour;
&
Center of mass;
&
Choquard-Pekar equation;
\\
Concentration phenomena;
&
Critical exponent;
&
Double nonlocality;
\\
Even and odd nonlinearities,
&
Existence and multiplicity;
&
Fractional Laplacian;
\\
Ground states;
&
Hartree-type term;
&
$L^2$-constraint; 
\\
Lagrange multiplier;
&
Mountain Pass paths;
&
Nonlocal sources;
\\
Nonlinear PDEs;
&
Normalized solutions;
&
Pohozaev identity;
\\
Polynomial decay;
&
Positivity and sign;
&
Prescribed mass;
\\
Qualitative properties;
&
Radial symmetry;
&
Regularity;
\\
Relative cup-length;
&
Riesz potential;
&
Schrödinger equation;
\\
Singular perturbation;
&
Spike solutions;
&
Sublinear nonlinearities.
\end{tabular}

%\bigskip

\bigskip

\noindent 
\textbf{Contacts:} \href{mailto: marco.gallo@uniba.it}{marco.gallo@uniba.it}, Dipartimento di Matematica, UNIBA (Bari).

$\qquad \quad \;$ \emph{Edit:} \href{mailto: marco.gallo1@unicatt.it}{marco.gallo1@unicatt.it}, %Dipartimento di Matematica e Fisica, 
Università Cattolica del Sacro Cuore (Brescia).
%
%
%
%\chapter*{Ringraziamenti}
%

%VERSIONE STAMPATA
%\restoregeometry

%\bigskip

%\bigskip

\bigskip

\bigskip

%\bigskip

%\newpage
%
%\cleardoublepage
%\thispagestyle{empty}

\noindent
%\emph{This thesis has been elaborated in the three years 2019--2022 and submitted first on 28th December 2022. The date of the defence is
%21st March 2023.}
\emph{This thesis has been elaborated in the three years 2019--2022, submitted first on 28th December 2022, and defended on 21st March 2023.}

\smallskip

%\noindent 
%\emph{\textbf{Internal Ph.D. commission:}
%%
%\begin{itemize}
%\item Prof. Anna Maria Candela (Università degli Studi di Bari Aldo Moro),
%\item Prof. Lorenzo D'Ambrosio (Università degli Studi di Bari Aldo Moro).
%\end{itemize}
%%
%\textbf{Ph.D. Thesis Evaluation commission:}
%%
%\begin{itemize}
%\item Prof. Jarosław Mederski (%Mathematical Institute of the 
%Polish Academy of Sciences, Warsaw),
%\item Prof. Tobias Weth (Goethe-Universität Frankfurt am Main).
%\end{itemize}
%%
%\textbf{Ph.D. Defense commission:}
%%
%\begin{itemize}
%\item Prof. Rossella Bartolo (Politecnico di Bari),
%\item Prof. Josè Carmona Tapia (Universidad de Almería),
%\item Prof. Pietro De Poi (Università degli Studi di Udine).
%\end{itemize}
%}

\noindent 
\emph{The contents are mainly based on the papers \cite{CG0, CG1, CGT1, CGT2, CGT3, CGT4, CGT5, Gal0, Gal1}, %, BoCiGa0, BG0},
and elaborated in collaboration with Prof. Silvia Cingolani (Università degli Studi di Bari Aldo Moro), Prof. Denis Bonheure (Université libre de Bruxelles) and Prof. Kazunaga Tanaka (Waseda University, Tokyo).
}

%\bigskip
%
%\emph{The thesis has been awarded with ``Enrico Jannelli'' Ph.D. thesis prize 2022/2023., whose examining board was composed by:
%%
%\begin{itemize}
%\item Prof. Piermarco Cannarsa (Università degli Studi di Roma Tor Vergata), 
%\item Prof. Alberto Farina (Université de Picardie Jules Verne, Amiens),
%\item Prof. Filomena Pacella (Sapienza Università di Roma).
%\end{itemize}
%}

%\bigskip
%
%\bigskip
%
%\noindent
%%\section*{{\normalsize \raggedleft{Acknowledgments}}}
%The number of people to thank for this %wonderful 
%three years journey is quite long, starting from my advisors, passing through collaborators and colleagues, up to family and friends. 
%I have written a truly marvelous praise about this, which this margin unfortunately is too narrow to contain.
%
%\medskip
%
%%\bigskip
%
%\rightline{\emph{So long, and thanks for all the fish!}}

\end{document}